%% file: mth.tex
\documentclass[10pt,russian]{book}
\usepackage[utf8]{inputenc}
\usepackage[T2A]{fontenc}
\usepackage{babel}
\usepackage{bm}       
\usepackage{amssymb}  
\usepackage{amsmath}  
\usepackage{stmaryrd} 
\usepackage{bbding}
\usepackage{calrsfs}
\usepackage{amsthm}
\usepackage[pic]{xy}  
\usepackage{srcltx}
\usepackage[dvipsnames]{xcolor}

\definecolor{MyLinkColour}{HTML}{0000B7}  
\definecolor{MyEmphColour}{HTML}{005781}
\definecolor{MyCiteColour}{HTML}{660099}  
\definecolor{MyFtNtColour}{HTML}{660099}  

\usepackage{fancyhdr}

\usepackage[nottoc]{tocbibind} 

\usepackage{tocloft}
\cftsetindents{chapter}{0em}{2em}  
\cftsetindents{section}{2em}{2.3em} 

\usepackage{chngcntr}
\counterwithout{section}{chapter}  

\usepackage{hyperref}
\hypersetup{
    final=true,  
    colorlinks=true,
    linkcolor=MyLinkColour,
    citecolor=MyCiteColour,
    linktoc=page,
    pdftitle={Представления и характеры конечных групп},
    pdfauthor={А. В. Заварницин; Д. О. Ревин},
    pdfsubject={Конспект лекций по теории представлений и характеров конечных групп},
    pdfkeywords={конечные группы, теория представлений, характеры, модулярные представления, брауэровы характеры, блоки характеров}
}

\usepackage[xindy={glsnumbers=false},
            toc,
            index,
            symbols,
            numberedsection=nameref
            ]{glossaries}

\usepackage{glossary-mcols} 

\newglossarystyle{mcolindexgroupfixed}{
  \setglossarystyle{mcolindexgroup}%
}

\newglossarystyle{symbolslong}{%
  \setglossarystyle{long}%
}

\makeglossaries

\newcommand*{\rusnewterm}[3][]{
      \newglossaryentry{#2}%
      {type={index},name={#3},description={\nopostdesc},#1}}

\newcommand*{\mem}[1]{\textcolor{MyEmphColour}{\textit{#1}}}

\renewcommand{\thechapter}{\Roman{chapter}}
\renewcommand{\thesection}{\arabic{section}}
\renewcommand{\thesubsection}{(\thesection.\arabic{subsection}\unskip )}

\makeatletter  
\renewcommand\@seccntformat[1]{\csname the#1\endcsname.\quad}
\makeatother

\makeatletter
\renewcommand\subsection{\@startsection{subsection}{2}{\z@}%
                                     {-3.25ex\@plus -1ex \@minus -.2ex}%
                                     {-1em \@minus -.1 em}%
                                     {\normalfont\normalsize\textbf}}
\makeatother

\newcommand{\la} {\langle}
\newcommand{\ra} {\rangle}
\newcommand{\se} {\subseteq}
\newcommand{\ot} {\otimes}
\newcommand{\st}{\star}

\renewcommand{\le}{\leqslant}
\renewcommand{\ge}{\geqslant}
\newcommand{\vf}{\varphi}
\newcommand{\vk}{\varkappa}
\newcommand{\ve}{\varepsilon}
\newcommand{\om}{\omega}
\renewcommand{\r}{\rho}
\newcommand{\s}{\sigma}
\renewcommand{\a}{\alpha}
\renewcommand{\b}{\beta}
\renewcommand{\d}{\delta}
\renewcommand{\l}{\lambda}
\renewcommand{\t}{\tau}
\renewcommand{\i}{\iota}
\newcommand{\z}{\zeta}
\newcommand{\Om}{\Omega}
\newcommand{\D}{\Delta}
\newcommand{\g}{\gamma}
\newcommand{\I}{\bm{1}}
\renewcommand{\O}{\bm{0}}

\newcommand{\x}{\chi}
\newcommand{\ov}{\overline}
\newcommand{\chk}{\check}

\newcommand{\FF} {\mathbb{F}}
\newcommand{\ZZ} {\mathbb{Z}}
\newcommand{\CC} {\mathbb{C}}
\newcommand{\QQ} {\mathbb{Q}}
\newcommand{\RR} {\mathbb{R}}
\newcommand{\nor}{\trianglelefteqslant}
\newcommand{\ba}[1] {\begin{array}{#1}}
\newcommand{\ea} {\end{array}}
\newcommand{\ld}{\ldots}

\def\clap#1{\hbox to 0pt{\hss#1\hss}}
\def\mathllap{\mathpalette\mathllapinternal}
\def\mathrlap{\mathpalette\mathrlapinternal}

\def\mathllapinternal#1#2{%
\llap{$\mathsurround=0pt#1{#2}$}}
\def\mathrlapinternal#1#2{%
\rlap{$\mathsurround=0pt#1{#2}$}}

\newcommand{\wtm}{\phantom{M}\mathllap{\widetilde{\phantom{R}\mkern2mu}\mathllap{M}}}
\newcommand{\whm}{\phantom{M}\mathllap{\widehat{\phantom{R}\mkern2mu}\mathllap{M}}}
\newcommand{\ti}{\tilde}
\newcommand{\wt}{\widetilde}
\newcommand{\wh}{\widehat}

\newcommand\Ker{\operatorname{Ker}}
\renewcommand\Im{\operatorname{Im}}
\newcommand\End{\operatorname{End}}
\newcommand\Aut{\operatorname{Aut}}
\newcommand\Inn{\operatorname{Inn}}
\newcommand\Hom{\operatorname{Hom}}
\newcommand\Ann{\operatorname{Ann}}
\newcommand\tr{\operatorname{tr}}
\newcommand\GL{\operatorname{GL}}
\newcommand\PSL{\operatorname{PSL}}
\newcommand\Syl{\operatorname{Syl}}
\newcommand\SL{\operatorname{SL}}
\newcommand\irr{\operatorname{irr}}
\newcommand\lin{\operatorname{lin}}
\newcommand\iBr{\operatorname{iBr}}
\newcommand\bl{\operatorname{bl}}
\newcommand\Gal{\operatorname{Gal}}
\newcommand\cf{\operatorname{cf}}
\newcommand\gch{\operatorname{gch}}
\renewcommand\f{\operatorname{f\mbox{}}}
\newcommand\diag{\operatorname{diag}}
\newcommand\supp{\operatorname{supp}}
\newcommand\vr{\varrho}
\newcommand\DD{\mathrm{D}}
\renewcommand\SS{\mathrm{S}}
\newcommand\N{\mathrm{N}}
\newcommand\Fr{\mathrm{Fr}}
\newcommand\St{\mathrm{St}}
\newcommand\Orb{\mathop\mathrm{Orb}\nolimits}
\renewcommand\C{\mathrm{C}}
\newcommand\OO{\mathrm{O}}
\newcommand\df{\mathrm{d}}
\newcommand\XX{\mathrm{X}}
\newcommand\W{\mathrm{W}}
\newcommand\II{\mathrm{I}}
\newcommand\n{\mathrm{n}}
\newcommand\Z {\mathrm{Z}}
\newcommand\J{\mathrm{J}}
\newcommand\M {\mathcal{M}}
\newcommand\A {\mathcal{A}}
\newcommand\B {\mathcal{B}}
\newcommand\CCC {\mathcal{C}}
\newcommand\R {\mathcal{R}}
\newcommand\K {\mathcal{K}}
\newcommand\X {\mathcal{X}}
\renewcommand\P {\mathcal{P}}
\renewcommand\H {\mathcal{H}}
\newcommand\E {\mathcal{E}}
\newcommand\Q {\mathcal{Q}}
\newcommand\Y {\mathcal{Y}}
\newcommand\ZZZ {\mathcal{Z}}
\newcommand\OOO {\mathcal{O}}
\newcommand\ZC {\mathcal{Z}}
\newcommand\SC {\mathcal{S}}
\newcommand\MM {\mathrm{M}}

\loadglsentries[index]{mth_indx}
\loadglsentries[symbols]{mth_symb}

\newtheoremstyle{thst}{\topsep}{\topsep}%
     {\itshape}
     {}
     {\bf}
     {. }
     { }
     {\thmnumber{#2\quad}\thmname{#1}\thmnote{ (#3)}}

\theoremstyle{thst}
\newtheorem{thm}[subsection]{Теорема}
\newtheorem{lem}[subsection]{Лемма}
\newtheorem{cor}[subsection]{Следствие}
\newtheorem{pre}[subsection]{Предложение}
\newtheorem{gip}[subsection]{Гипотеза}

\newtheoremstyle{opst}{\topsep}{\topsep}%
     {}
     {}
     {\bf}
     {. }
     { }
     {\thmnumber{#2\quad}\thmname{#1}\thmnote{ (#3)}}

\theoremstyle{opst}
\newtheorem{opr}[subsection]{Определение}
\newtheorem{prim}[subsection]{Пример}

\newcommand{\uk}[1]{$\big(${\it Указание.\ \ }#1$\big)$}

\newcommand{\myeqno}{\refstepcounter{subsection}\leqno{\bm{\thesubsection}}}

\newcommand{\mysubsection}{\refstepcounter{subsection}\noindent{\bf\thesubsection}\quad}

\newcommand{\textzam}[1]{\indent\llap{\lower.5ex\hbox{\Checkmark}\enspace\enspace}\ignorespaces#1}
\newcommand{\textupr}[1]{\medskip\noindent\llap{\PencilRightDown\quad}\ignorespaces#1\bigskip}
\newcommand{\textupl}[2]{\medskip\noindent\llap{\hyperlink{#1}{\PencilRightDown}\quad}\ignorespaces#2\bigskip}
\newcommand{\textuprn}[1]{\medskip\noindent\llap{\PencilRightDown\quad}\ignorespaces\mysubsection#1\bigskip}
\newcommand{\textupln}[2]{\medskip\noindent\llap{\hyperlink{#1}{\PencilRightDown}\quad}\ignorespaces\mysubsection#2\bigskip}
\newcommand{\textext}[1]{\smallskip\noindent\llap{\FourStarOpen\quad}\ignorespaces{\footnotesize #1\par}\medskip}

\title{\Huge \textsf{Представления и характеры\\конечных групп}\\[50pt]
             \large{конспект лекций}\\[10pt]
             \small{версия 1.1}\\[20pt]}

\author{\sc{А.\,В.\,Заварницин, \quad Д.\,О.\,Ревин}}

\date{\mbox{}\\[250pt]
      Новосибирск\\[10pt]
      2026 г.}

\begin{document}

\frontmatter

\maketitle

\renewcommand{\th}{\theta} 

\renewcommand{\ch}{\operatorname{ch}} 

\renewcommand{\proofname}{\sc{Доказательство}}

\setcounter{tocdepth}{1} 
\tableofcontents

\renewcommand{\sectionmark}[1]{\markright{\thesection.\ #1}}
\renewcommand{\chaptermark}[1]{%
\markboth{#1}{}}

\chapter{Предисловие}

Предлагаемый вниманию читателя текст представляет собой расширенный конспект лекций, в течение нес\-коль\-ких лет
читавшийся авторами для студентов и аспирантов Новосибирского государственного университета и Института
математики им. С.\,Л.\,Соболева СО РАН. Эти лекции были ориентированы на тех, кто специализируется
в области алгебры и знаком с основами линейной алгебры и теории групп. Публикуя данное пособие, мы также предполагаем знакомство читателя
с этими разделами математики. При необходимости можно обратиться к источникам \cite{w,vlm,km,l,m}.

Всякий, кто занимается конечными группами, знает, что для их изучения теория представлений
является незаменимым инструментом. Суть теории представлений можно кратко выразить как изучение групп
теоретико-кольцевыми методами. Групповые модули и представления над полями конечной характеристики естественно
возникают при исследовании подгрупповой структуры данной группы.

Уже в самом начале своего существования теория представлений и характеров продемонстрировала свою мощь,
позволив получить ряд исключительно важных результатов теории групп (теорема Бернсайда о разрешимости групп
порядка $p^a q^b$, теорема о строении групп Фробениуса, теорема Брауэра-Фаулера о конечности числа конечных
групп с данным централизатором инволюции и т.\,д.), причём некоторые из этих результатов до сих пор не имеют
чисто теоретико-группового доказательства. В XX веке теория конечных групп сделала огромный рывок.
Классификация простых конечных групп и многие другие глубокие результаты были получены с использованием теории
представлений и характеров. Предлагаемое пособие знакомит читателя с основами этой теории,
её связью с теорией обыкновенных
представлений и включает изучение таких понятий, как брауэровы характеры, $p$-блоки, дефектные группы.

Мы не претендуем на полноту и оригинальность. В изложении многих разделов мы следовали превосходным пособиям
М.\,Айзекса \cite{i} и Г.\,Наварро \cite{n}. Из русскоязычной литературы по данной тематике отметим
монографию В.\,А.\,Белоногова \cite{b}.
Некоторые важные результаты (третья основная теорема Брауэра,
теорема Дэйда о блоках с циклической дефектной группой, $Z^*$-теорема Глаубермана) приводятся
нами в заключительной части без доказательства.

В качестве приложений даны необходимые сведения из теории множеств, колец, полей, алгебраических чисел.
Доказательства многих утверждений оставляются читателю в качестве упражнений. Решение части упражнений
приведено в конце книги.
Возможно, полезным для ссылок будет также материал из приложения, содержащего таблицы характеров
(в том числе брауэровых), матрицы разложения и матрицы Картана некоторых отдельных конечных групп,
а также их серий.

При подготовке обновлённой версии мы внесли ряд изменений. 
Авторы хотели бы выразить искреннюю признательность А.\,А.\,Бутурлакину
за внимательное чтение первоначального текста, замечания и предложенные улучшения.

\mainmatter

\renewcommand{\chaptermark}[1]{%
{\markboth{\thechapter.\ #1}{}}}

\chapter{Алгебры и модули}

\section{Используемые обозначения и понятия\label{isp}}

Утверждения, помеченные символом \  \PencilRightDown \ , оставляются читателю в качестве упражнений.
В последнем разделе можно найти доказательства некоторых из них. Материал, напечатанный мелким шрифтом  и отмеченный
символом \ \FourStarOpen \ , носит вспомогательный характер.

Мы предполагаем, что читатель знаком с базовыми понятиями теории групп (см., например, \cite{km}). Используемые нами обозначения
в основном стандартны. Значительная их часть приведена в списке на стр. \pageref{symbols}, а некоторые мы поясним ниже.

Через \gls{dij} обозначается
\glsadd{idij}\mem{символ Кронекера}, который определяется равенством
$$
\d_{ij}=\left\{\ba{l}
1,\quad\mbox{при}\ i=j;\\
0,\quad\mbox{при}\ i\ne j.
\ea\right.
$$
Наибольший общий делитель целых чисел $m$ и $n$ будем обозначать через \gls{lmnr}. Натуральное число
вида $p^a$ или $p^aq^b$ для простых $p$ и $q$ называется \glsadd{iPrimNum}\mem{примарным}
или \glsadd{iBiprNum}\mem{бипримарным}, соответственно.

Для матрицы $M$ её транспонированную матрицу будем обозначать через \gls{atop}. Если  $M$ --- квадратная невырожденная матрица, то будем
обозначать через \gls{Mi}
матрицу $(M^{-1})^\top=(M^\top)^{-1}$ и называть её
\glsadd{iMatInvTr}\mem{обратно-транспонированной} к $M$.

Буквенные символы отображений и, в частности, морфизмов мы будем писать
как справа, так и слева от аргумента, мотивируя свой выбор соображениями удобства или
традиции. Вариантом правой записи будет экспоненциальная запись (например, $g^\a$ для образа элемента $g$ группы
под действием её автоморфизма $\a$).
Для композиции отображений также будут использоваться различные
виды записи. Запись $\vf\psi$ композиции отображений $\vf:X\to Y$ и $\psi:Y\to Z$ является <<правой>>, т.\,е.
$x(\vf\psi)=(x\vf)\psi$, а запись $\psi\circ\vf$
--- <<левой>>, т.\,е. $(\psi\circ\vf)(x)=\psi(\vf(x))$.

\mem{Тривиальной}\glsadd{iGrTrv} называется группа, состоящая только из нейтрального элемента.
При мультипликативной (аддитивной) записи операции нейтральный элемент группы $G$, а также иногда и тривиальную
группу, мы будем обозначать символом $1$ (соответственно, символом $0$). Порядок элемента $g$ группы $G$ будем
обозначать через
\gls{gmod}. Однако, в случае, когда $g$ является комплексным числом, мы сохраним
запись $|g|$ для классического обозначения модуля.
\mem{Центр}\glsadd{iCenGr} группы $G$, т.\,е. множество $\{z\in G\bigm| zg=gz\ \mbox{для всех}\ g\in G\}$,
обозначается через \gls{ZG}.
Подгруппа группы $G$, лежащая в $\Z(G)$, называется \glsadd{iSubGrCen}\mem{центральной}.
Если $\vf: G\to H$ --- гомоморфизм групп, то\ \ $\gls{kervf}=\{g\in G\mid g\vf=1\}$
--- его \glsadd{iKerGrHom}\mem{ядро}.  Для элементов $g,h\in G$ положим $h^g=g^{-1}hg$. Два элемента $x, y\in G$ называются
\glsadd{iCnjElms}\mem{сопряжёнными}, если существует $g\in G$ такой, что $y=x^g$. Сопряжённость
элементов является отношением эквивалентности в группе $G$.
\mem{Классом сопряжённости} \glsadd{iClsCnj} \gls{xG}, содержащим элемент $x\in G$ называется класс эквивалентности, которому
принадлежит $x$. Будем обозначать через
\gls{KlGr} множество всех классов сопряжённости группы $G$, а через
\gls{xsK} --- представитель класса $K\in \K(G)$.
Для элементов $x,y\in G$ определим \glsadd{ixy}
$\gls{xy}=x^{-1}y^{-1}xy$. \glsadd{iCmt}\mem{Коммутантом} \gls{Gpr} группы $G$
будем называть её подгруппу, порождённую коммутаторами $[x,y]$, для всех $x,y\in G$.

Пусть $X\se G$ --- произвольное подмножество группы $G$. Через
\gls{laXra} будем обозначать наименьшую подгруппу в $G$, содержащую $X$, и называть её
\glsadd{iSubGrGenX}\mem{подгруппой, порождённой множеством $X$}.
\glsadd{iCntrlGr}\mem{Централизатором} и
\glsadd{iNrmtr}\mem{нормализатором}
$X$ в $G$ называются, соответственно множества
\pagebreak[0] 
\begin{align*}
\gls{CGX}&=\{\,g\in G\bigm|g^{-1}xg=x \ \  \mbox{для всех}\ x\in X\,\}\\
\gls{NGX}&=\{\,g\in G\bigm|g^{-1}xg\in X \ \  \mbox{для всех}\ x\in X\,\}
\end{align*}
которые, как легко видеть, являются
подгруппами в $G$.
Если $g\in G$, то нормализатор $\N_G(\la g \ra)$ будем для краткости обозначать через
\gls{NGg}.

Термин <<собственный>> в отношении подгруппы (идеала, подкольца, подмножества и т.\,п.) $N$ группы (кольца,
множества и т.\,п.) $M$ означает, что $N\ne M$.

Под термином  \mem{максимальная}\glsadd{iSubGrMax}  (соответственно, \glsadd{iSubGrMin}\mem{минимальная})
подгруппа подразумевается максимальный (соответственно, минимальный)
по включению элемент среди всех собственных (соответственно, нетривиальных) подгрупп данной группы. Аналогичная
терминология будет использоваться также и в отношении колец, идеалов, модулей и т.\,п.

Пусть $G$ --- группа, $M$ --- множество. Будем говорить, что $G$ \mem{действует \glsadd{iActGr} $($справа\/$)$} на
$M$, если для любых $m\in M$, $g\in G$ однозначно определён элемент $mg\in M$ так, что выполнены условия $m1=m$
и $(mg)h=m(gh)$ для произвольных $m\in M$, $g,h\in G$. Множество $M$ однозначно представимо в виде объединения
непересекающихся \mem{$G$-орбит} \glsadd{iGOrb} так, что элементы $m_1,m_2\in M$ лежат в одной $G$-орбите
тогда и только тогда, когда существует элемент $g\in G$, для которого $m_1g=m_2$. Пусть $m\in M$.
Тогда $G$-орбиту, содержащую $m$,  будем называть \glsadd{iGOrbEtl}\mem{$G$-орбитой элемента} $m$
(или просто \mem{орбитой элемента} $m$ $)$ и обозначать через \gls{OrbGm}.
\mem{Стабилизатором \gls{StGm} элемента} $m\in M$ называется множество
$$
\St_G(m)=\{\,g\in G\mid mg=m\,\},
$$
которое оказывается подгруппой в $G$.
Если из контекста ясно, о действии какой группы $G$ идёт речь, то вместо $\St_G(m)$  и $\Orb_G(m)$ мы будем писать
\gls{Stm} и \gls{Orbm}, соответственно.
Действие, при котором $G$-орбита единственна (и, значит, совпадает с $M$), называется
\glsadd{iActGrTr}\mem{транзитивным}. Другими словами, $G$ действует на $M$ транзитивно, если
$M$ не пусто и для любых $m_1,m_2\in M$ существует элемент $g\in G$
такой, что $m_1g=m_2$.

\section{\texorpdfstring{$R$}{R}-модули}

Термин <<кольцо>> всегда будет обозначать ассоциативное кольцо
с единицей. Необходимые нам основные сведения о кольцах содержатся в приложении \ref{cv kol}.

\begin{opr} Пусть $R$ --- кольцо. Абелева группа $M$ с
аддитивной записью операции называется \glsadd{iRModL}\mem{левым $R$-модулем},
если для любых $m\in M$ и $r\in R$ однозначно определён элемент $rm\in M$, и для произвольных $m,n\in M$ и $r,s\in R$ выполнены условия

$$\ba{c}
r(m+n)=rm+rn, \\
(r+s)m=rm+sm, \\
(rs)m=r(sm), \\
1_R\,m=m.
\ea
$$
Аналогично, путём умножения справа на элементы из $R$, определяется \glsadd{iRModR}\mem{правый $R$-модуль}.
\end{opr}

\textzam{В случае, когда $R$ --- поле, определение левого $R$-модуля
становится обычным определением векторного пространства над этим полем.}

\textupr{У нулевого кольца существует только нулевой модуль.}

\textzam{Часто $R$-модули, $R$-подмодули,
и т.\,д. мы будем называть просто модулями,  подмодулями, и т.\,д., если из
контекста ясно, о чём идёт речь.}

Для кольца $R$, обладающего антиавтоморфизмом $\vf$, левый $R$-модуль $M$ можно превратить в правый $R$-модуль,
положив $mr=(r\vf)m$ для всех $r\in R$ и $m\in M$. Кроме того, произвольный левый $R$-модуль $M$ можно превратить в правый $R^{op}$-модуль,
положив $mr=rm$ для всех $r\in R^{op}$ и $m\in M$. В частности, если $R$ коммутативно, то
$M$ можно рассматривать и как левый, и как правый $R$-модуль.

\subsection{Примеры модулей} \label{ex mod}
\begin{list}{{\rm(}{\it\roman{enumi}\/}{\rm)}}
{\usecounter{enumi}\setlength{\parsep}{2pt}\setlength{\topsep}{5pt}\setlength{\labelwidth}{23pt}}

\item Левый (правый) идеал произвольного кольца $R$ является левым (правым) $R$-модулем.
В частности, само кольцо $R$ является левым (правым) $R$-модулем, который обозначается через
\gls{Rl0} (соответственно, через \gls{Rr0}) и называется
\glsadd{iRModReg}\mem{регулярным левым} $($\mem{правым}$\/)$ \mem{$R$-модулем}.

\item Произвольная абелева группа $A$ является одновременно левым и правым $\ZZ$-модулем, где произведение числа $n\in \ZZ$  на $a\in A$ равно $n$-й
(аддитивной) степени элемента $a$.

\item Произвольное кольцо является левым и правым модулем над любым своим подкольцом.

\item При правой (левой) записи отображений и композиций векторное пространство $V$ над полем $F$
является правым (левым) модулем над кольцом $\End_F(V)$ \big(соответственно, $\End_F(V)^{op}$, см. \ref{end op}\big) всех своих линейных преобразований.

\item Пространство строк (столбцов) длины $n$ над полем $F$ является правым (левым) $\MM_n(F)$-модулем.

\item Для произвольного множества $X$ абелева группа $\P(X)$ с операцией симметрической разности
является $\ZZ_2$-модулем, если для любого $A\in\P(X)$ положить $0A=\varnothing$ и $1A=A$.
\end{list}

\textupr{Проверить сформулированные в \ref{ex mod} утверждения.}

\textuprn{В случае, когда $X$ --- $n$-элементное множество,
найти базис и размерность $\P(X)$ как векторного пространства над полем $\ZZ_2$.}

\textzam{В дальнейшем определения, дающиеся только для левых модулей, аналогично формулируются и для правых модулей.}

\glsadd{iSubMod}\mem{Подмодулем} левого $R$-модуля $M$ называется такая
подгруппа $N\le M$, что для любых $n\in N$ и $r\in R$ выполнено $rn\in N$. Если $N$ --- подмодуль в $M$, то на
факторгруппе $M/N$ можно естественным образом задать структуру левого $R$-модуля (называемого
\glsadd{iQMod}\mem{фактормодулем} $M$ по $N$), полагая $r(m+N)=rm+N$ для любых
$m\in M$, $r\in R$.

Пусть $M,N$ --- левые $R$-модули. Гомоморфизм $\vf:M\to N$ абелевых групп такой, что для любых $m\in M$ и $r\in
R$ выполнено $(rm)\vf=r(m\vf)$ называется \glsadd{iHomRMod}\mem{гомоморфизмом $R$-модулей}
(или \glsadd{iMapRlin}\mem{$R$-линейным отображением}).
 Множество всех $R$-линейных
отображений $\vf:M\to N$ обозначается через
\gls{HomR}.
Для любых $\vf,\psi\in \Hom_R(M,N)$ определим $\vf+\psi$, положив $m(\vf+\psi)=m\vf+m\psi$ для любого $m\in M$.
Тем самым на $\Hom_R(M,N)$ задана структура абелевой группы.\footnote{На группе $\Hom_R(M,N)$ задать естественную структуру левого $R$-модуля,
вообще говоря, нельзя.}

\textuprn{Пусть $M$ --- левый (правый) $R$-модуль и $H=\Hom_R(\mbox{}^\circ\! R,M)$ \big(соответственно,
$H=\Hom_R(R^\circ,M)$\big). Зададим на абелевой группе $H$ умножение слева (справа)
на элементы из $R$ по правилу $s(r\vf)=(sr)\vf$ (соответственно, $(\vf r)(s)=\vf(rs)$ в левой записи) для любых $\vf\in H$, $r,s\in R$.
Показать, что тем самым на $H$ корректно задана структура левого (правого) $R$-модуля и имеет место
изоморфизм $R$-модулей $H\cong M$.}

\textuprn{\label{dv mod} Пусть $M$ --- левый (правый) $R$-модуль и
$\gls{Mst}=\Hom_R(M,\mbox{}^\circ\! R)$ \big(соответственно,
$\gls{Mst}=\Hom_R(M,R^\circ)$\big). Зададим на абелевой группе $M^*$ умножение справа (слева)
на элементы из $R$ по правилу $m(\vf r)=(m\vf)r$ (соответственно, $(r\vf)(m)=r\vf(m)$ в левой записи) для любых $\vf\in M^*$,
$m\in M$, $r\in R$. Показать, что тем самым на $M^*$ корректно задана структура правого (левого) $R$-модуля,
который будем называть \glsadd{iRModDl}\mem{двойственным} к $R$-модулю~$M$.}

Гомоморфизм модуля $M$ в себя называется \mem{эндоморфизмом}\glsadd{iEndRmod} $R$-модуля $M$.
Множество всех эндоморфизмов модуля $M$, обозначаемое
\gls{EndR}, образует кольцо, умножением в котором является
композиция отображений. Легко видеть, что $M$
является правым $\End_R(M)$-модулем относительно естественного действия.

\textupln{ndr prf}{\label{ndr}Пусть $R$ --- кольцо. Показать, что
\begin{list}{{\rm(}{\it\roman{enumi}\/}{\rm)}}
{\usecounter{enumi}\setlength{\parsep}{2pt}\setlength{\topsep}{5pt}\setlength{\labelwidth}{23pt}}
\item кольца $\End_R(R^\circ)$ и $R$ антиизоморфны;
\item для любого идемпотента $e\in R$ кольца $\End_R(eR)$ и $eRe$ антиизоморфны.
\end{list}
}

Если для $R$-модулей $M$ и $N$  существует взаимно однозначный гомоморфизм, то он называется
\glsadd{iIsoRMod}\mem{изоморфизмом}. В этом случае будем говорить, что модули $M$ и $N$
\mem{изоморфны} и писать $M\cong N$.


Пусть $N_1,\ld, N_s$ --- подмодули левого $R$-модуля $M$. Их
\glsadd{iSumRSubMods}\mem{суммой} $N_1 +\ld+ N_s$ называется
подмножество $\{\, n_1+\ld+n_s\bigm| n_i\in N_i,\ i=1,\ld,s\}$. Очевидно, что $N_1 +\ld+ N_s$ является
наименьшим подмодулем модуля $M$, содержащим подмодули $N_1,\ld, N_s$. Сумма подмодулей $N_1,\ld, N_s$
называется
\glsadd{iSumDirRSubMods}\mem{прямой}, если для любого $i=1,\ld,s$
пересечение подмодуля $N_i$ с суммой $\sum_{j\ne i}N_j$ тривиально. Прямая сумма подмодулей обозначается $N_1
\oplus \ld \oplus N_s$.

\mysubsection\label{prsm} Как и для колец,
для $R$-модулей $M_1,\ld, M_s$ можно определить внешним образом
\glsadd{iSumDirRMods}\mem{прямую сумму} $M=M_1\oplus\ld\oplus M_s$.
При этом каждый из модулей $M_1,\ld, M_s$ можно естественным образом отождествить с подмодулем модуля $M$ так,
что $M$ будет совпадать с прямой суммой таких подмодулей.

Если $X$ --- подмножество левого $R$-модуля $M$, то \glsadd{iSubModGenX}\mem{подмодулем,
порождённым} $X$, называется наименьший подмодуль из $M$, содержащий $X$.

\textzam{Введённые понятия, в частности, согласуются с данными ранее
определениями суммы и прямой суммы левых (правых) идеалов кольца, а также с определением идеала, порождённого
множеством.}

Легко видеть, что модуль $M$ порождается множеством $X$ тогда и только тогда, когда всякий элемент из $M$
представим в виде конечной суммы $\sum r_ix_i$, где $r_i\in R$, $x_i\in X$. В случае, когда для любого элемента
из $M$ такое представление однозначно, множество $X$ называется \mem{$R$-базисом} или
просто \glsadd{iBasMod}\mem{базисом} модуля $M$. Если $R$ --- поле, то введённое понятие базиса $R$-модуля совпадает с
обычным понятием базиса векторного пространства. В частности, в этом случае базис всегда существует, и
все базисы равномощны. Если же
$R$~--- не поле, то модуль $M$ может не обладать базисом, или же обладать базисами различной мощности.

\textuprn{Привести пример кольца $R$ и $R$-модуля $M$ таких, что $M$ не обладает $R$-базисом.\footnote{Пример свободного
$R$-модуля, обладающего базисами различной мощности, приведён в \cite[с. 61--62]{skor}}}

$R$-модуль $M$, обладающий $R$-базисом $X$ называется \glsadd{iRModFr}\mem{свободным}. В этом случае мы
будем говорить, что $M$ \mem{свободно порождается} множеством $X$.

\textuprn{\label{dv bas} Пусть $M$ --- свободный $R$-модуль с конечным $R$-базисом $m_1, \ld, m_k$.
Показать, что
\begin{list}{{\rm(}{\it\roman{enumi}\/}{\rm)}}
{\usecounter{enumi}\setlength{\parsep}{2pt}\setlength{\topsep}{5pt}\setlength{\labelwidth}{23pt}}
\item $M$ изоморфен прямой сумме $k$ копий регулярного $R$-модуля;
\item двойственный\footnote{см. \ref{dv mod}} модуль
$M^*$ свободный, и множество его элементов $m_1^*,\ld,m_k^*$ таких, что $m_im_j^*=\d_{ij}$, где $i,j=1,\ld,k$, образует
$R$-базис, который будем называть
базисом модуля $M^*$, \glsadd{iBasDl}\mem{двойственным} к базису $m_1,\ld, m_k$. В частности, в случае, когда $R$ --- поле, имеем
$\dim_RM=\dim_RM^*$.}

\end{list}

\begin{pre} \label{zri} Если кольцо $R$ порождается множеством $R^\times$ как $\Z(R)$-модуль, то
$\Z(R^\times)\le \Z(R)^\times$.
\end{pre}

\textupl{zri prf}{Доказать предложение \ref{zri}}

\begin{opr} Левый $R$-модуль $M$ называется \glsadd{iRModIndec}\mem{неразложимым},
если для любых подмодулей $M_1,\allowbreak M_2\le M$ из равенства $M=M_1\oplus M_2$ следует $M_1=0$ или
$M_2=0$. В противном случае $M$ называется \glsadd{iRModDec}\mem{разложимым}.
\end{opr}

\textupr{Какие векторные пространства над полем $F$ неразложимы как $F$-модули?}

\begin{pre} \label{idem osd} Пусть $R$ --- кольцо
и $e\in R$ --- идемпотент. Тогда следующие условия эквивалентны.
\begin{list}{{\rm(}{\it\roman{enumi}\/}{\rm)}}
{\usecounter{enumi}\setlength{\parsep}{2pt}\setlength{\topsep}{5pt}\setlength{\labelwidth}{23pt}}
\item $e$ --- примитивный идемпотент.
\item Левый идеал $Re$ (правый идеал $eR$) неразложим как левый (правый) $R$-модуль.
\item Единственными идемпотентами кольца $eRe$ являются $0$ и $e$.
\end{list}
\end{pre}

\textupl{idem osd prf}{Доказать предложение \ref{idem osd}.}

Двусторонний аналог предыдущего предложения также имеет место.

\begin{pre} \label{idem dsd} Пусть $R$ --- кольцо
и $e\in R$ --- центральный идемпотент. Тогда следующие условия эквивалентны.
\begin{list}{{\rm(}{\it\roman{enumi}\/}{\rm)}}
{\usecounter{enumi}\setlength{\parsep}{2pt}\setlength{\topsep}{5pt}\setlength{\labelwidth}{23pt}}
\item $e$ --- примитивный идемпотент кольца $\Z(R)$.
\item Идеал $Re$ не представим в виде прямой суммы двух ненулевых двусторонних идеалов из $R$.
\end{list}
\end{pre}

\textupl{idem dsd prf}{Доказать предложение \ref{idem dsd}.}

Следующее предложение обобщает \ref{ids}.

\begin{pre}\label{m id dec} Пусть $R$ --- кольцо,
$e_1,\ld,e_n\in R$ --- центральные попарно ортогональные идемпотенты и $M$ --- правый $R$-модуль. Справедливы
следующие утверждения.
\begin{list}{{\rm(}{\it\roman{enumi}\/}{\rm)}}
{\usecounter{enumi}\setlength{\parsep}{2pt}\setlength{\topsep}{5pt}\setlength{\labelwidth}{23pt}}
\item $Me_i$ --- $R$-подмодуль в $M$ для всех $i=1,\ld,n$.
\item Сумма $Me_1+\ld +Me_n$ является прямой.
\item Если $e_1+\ld+e_n=1_R$, то $M=Me_1\oplus\ld\oplus Me_n$.
\end{list}
\end{pre}
\textupl{m id dec prf}{Доказать предложение \ref{m id dec}.}

\section{Модули над алгебрами}

\begin{opr} Пусть $R$ --- коммутативное кольцо. Левый $R$-модуль $A$ называется
\glsadd{iRAlg}\mem{$R$-алгеброй} или \glsadd{iAlgRng}\mem{алгеброй над кольцом $R$}, если он является кольцом и
для произвольных $a,b\in A$ и $r\in R$ выполнено
$$
r(ab)=(ra)b=a(rb).
$$
\end{opr}

\textzam{Далее до конца раздела символом $R$ обозначается некоторое коммутативное кольцо, а символом $A$ ---
алгебра над $R$.}

\subsection{Примеры \texorpdfstring{$R$}{R}-алгебр}\label{ex r-alg}
\mbox{} 

\begin{list}{{\rm(}{\it\roman{enumi}\/}{\rm)}}
{\usecounter{enumi}\setlength{\parsep}{2pt}\setlength{\topsep}{5pt}\setlength{\labelwidth}{23pt}}
\item Кольцо $\MM_n(R)$ является $R$-алгеброй.
\item Для левого $R$-модуля $M$ кольцо $\End_R(M)$ является $R$-алгеброй.
Умножением в алгебре $\End_R(M)$ является композиция отображений.

\item Для конечной группы $G$ множество $RG$,
состоящее из всевозможных формальных\footnote{Слово <<формальных>> означает, что две такие
$R$-линейные комбинации равны тогда и только тогда, когда равны коэффициенты при соответствующих элементах
группы $G$, т.\,е. $RG$ свободно порождается множеством $G$ как $R$-модуль.} $R$-линейных комбинаций
элементов $g\in G$, т.\,е. выражений вида
$$
 \sum_{g\in G} a_g g,\quad a_g\in R
$$
является $R$-алгеброй, которая называется \mem{групповой алгеброй}\glsadd{iAlgGr} группы $G$.
Умножение в $RG$  продолжает по линейности умножение в $G$, т.\,е.
$$
 \Big(\sum_{g\in G} a_g g\Big) \Big(\sum_{g\in G} b_g g\Big)=
\sum_{g\in G}c_g g,\quad \text{где} \qquad c_g\,=\!\!\!\sum_{\substack{x,y\in G,\\xy=g}}\!\! a_x b_y.
$$

\item Произвольное кольцо $S$ является $R$-алгеброй над любым подкольцом $R\le \Z(S)$ своего центра.
В частности, любое поле является алгеброй над любым своим подполем.

\item Произвольное кольцо является $\ZZ$-алгеброй.

\item Для произвольного множества $X$ множество $\P(X)$ является $\ZZ_2$-алгеброй,
см. \ref{pr kol}, \ref{ex mod}.

\item Для произвольного множества $X$ обозначим через \gls{fRX} множество всех отображений $X\to R$.
Для двух отображений $\vf,\psi\in\f_R(X)$ естественно определяются сумма, произведение и
умножение на элементы кольца $R$:
$$(\vf+\psi)(x)=\vf(x)+\psi(x), \qquad(\vf\psi)(x)=\vf(x)\psi(x), \qquad(r\vf)(x)=r(\vf(x))$$
для всех $x\in X$,
$r\in R$. Относительно определённых операций множество $\f_R(X)$ становится
коммутативной $R$-алгеброй, которую мы будем называть \glsadd{iAlgRvFuncs}\mem{алгеброй $R$-значных функций на множестве} $X$.
\end{list}

\textupr{Проверить сформулированные в \ref{ex r-alg} утверждения.}

\mysubsection \label{1ot}
В случае, когда  $A$ --- либо ненулевая $R$-алгебра над полем $R$ либо групповая алгебра над кольцом, $R$ можно
отождествить c подкольцом $R\,1_A=\{r1_A\bigm|r\in R\}$ центра $\Z(A)$.

\textuprn{Показать, что для произвольных $R$-алгебр кольцо $R1_A$ может быть неизоморфно $R$.}

Из определения вытекает, что для $R$-алгебры $A$ противоположное кольцо $A^{op}$ также является 
$R$-алгеброй с сохранением прежней структуры $R$-модуля.

\textuprn{Показать, что групповая алгебра изоморфна своей противоположной, т.\,е. обладает антиавтоморфизмом.}


\textuprn{Пусть $G$ --- конечная группа и $R$ --- коммутативное кольцо. На алгебре $R$-значных функций $\f_R(G)$, рассматриваемой как $R$-модуль,
введём новую операцию умножения (так называемую <<свёртку>> функций), положив
$$
(\vf*\psi)(g)=\sum_{\substack{x,y\in G,\\xy=g}} \vf(x)\psi(y)
$$
для всех $\vf,\psi\in\f_R(G)$. Показать, что таким образом на $\f_R(G)$ задаётся структура $R$-алгебры, которая изоморфна
групповой алгебре $RG$. \uk{Рассмотреть отображение $\vf\mapsto\sum_{x\in G}\vf(x)x$.}}

Пусть $A$ --- $R$-алгебра. Подкольцо кольца $A$, являющееся подмодулем $R$-модуля $A$ называется \glsadd{iSubAlg}\mem{подалгеброй}
$R$-алгебры $A$. \glsadd{iIdlRAlg}\mem{Левым} (\mem{правым, двусторонним}) \mem{идеалом $R$-алгебры $A$}
называется произвольный левый (правый, двусторонний) идеал кольца $A$.

\textuprn{Показать, что левый (правый, двусторонний) идеал $R$-алгебры $A$ будет 
$R$-подмодулем $R$-модуля $A$, а подкольцо, вообще говоря, нет.}


Заметим, что центр $R$-алгебры $A$ всегда является подалгеброй. Исследуем более подробно центр $Z(A)$ в важном
для нас случае, когда $A$ --- групповая алгебра над коммутативным кольцом. Существенную роль при этом будут
играть классы сопряжённости группы.

Напомним, что через $\K(G)$ обозначается множество всех классов сопряжённости группы $G$, а через
$x_{\mbox{}_K}$~--- представитель класса $K\in \K(G)$.

Пусть $R$ --- коммутативное кольцо. Для $K\in \K(G)$ положим
$$\wh{K}=\sum_{x\in K} x\in RG.
$$
Будем называть \gls{Kwh} \glsadd{iSumClsK}\mem{суммой класса} $K$ или \glsadd{iClsSum}\mem{классовой суммой}.

\begin{pre} \label{kl sum}
Пусть $R$ --- коммутативное кольцо и $G$ --- конечная группа. Имеют место следующие утверждения.
\begin{list}{{\rm(}{\it\roman{enumi}\/}{\rm)}}
{\usecounter{enumi}\setlength{\parsep}{2pt}\setlength{\topsep}{5pt}\setlength{\labelwidth}{23pt}}
\item Если $K\in \K(G)$, то $\wh{K}\in \Z(RG)$.
\item Любой элемент из $\Z(RG)$ однозначно представим в виде $R$-линейной комбинации
элементов $\wh{K}$, $K\in \K(G)$.
\item Пусть $K,L,M\in \K(G)$  и $z\in M$. Тогда число $a_{\mbox{}_{KLM}}$ всевозможных пар
$(x,y)\in K\times L$ таких, что $xy=z$ не зависит от выбора $z$.
\item Если $K,L\in \K(G)$, то
$$
\wh{K}\wh{L}=\sum_{M\in \K(G)} a_{\mbox{}_{KLM}} \whm.
$$
\end{list}
\end{pre}

\textupl{kl sum prf}{Доказать предложение \ref{kl sum}.}

Как видно из \ref{kl sum}$(ii),(iv)$ классовые суммы $\wh{K}$, $K\in \K(G)$, образуют базис $R$-модуля
$\,\Z(RG)$, а числа \gls{aklm}
оказываются\footnote{Более строго, под структурными константами центра $\Z(RG)$ в этом
базисе следовало бы понимать образы чисел $a_{\mbox{}_{KLM}}$ относительно единственного кольцевого гомоморфизма $\ZZ\to R$, см. \ref{deg homs}.}
\glsadd{iStrConsZRG}\mem{структурными константами}
алгебры $\,\Z(RG)$ относительно этого базиса, т.\,е. однозначно определяют
произведение произвольных элементов из $\,\Z(RG)$.

\begin{cor} \label{zsr} Пусть $R$ --- коммутативное кольцо, $S$ --- его подкольцо и $G$ --- конечная группа. Тогда
$$
\Z(SG)=\Z(RG)\cap SG.
$$
\end{cor}
\textupl{zsr prf}{Доказать следствие \ref{zsr}.}

\mysubsection \label{pg}
Пусть $G$ --- конечная группа. Тогда любой гомоморфизм (эпиморфизм) $\vf:R\to S$ коммутативных колец естественным образом
поднимается до кольцевого гомоморфизма (эпиморфизма) $\ti\vf:RG\to SG$ соответствующих групповых алгебр и,
как следует из \ref{kl sum}$(ii)$, до гомоморфизма (эпиморфизма) $\ti\vf:\Z(RG)\to \Z(SG)$ их центров.

\begin{pre} \label{kvt} Пусть отображение $\ti\vf:RG\to SG$
$\big($соответственно, $\ti\vf:\Z(RG)\to \Z(SG)\big)$ определено как выше. Тогда  $\Ker\ti\vf$ состоит из
$\Ker\vf$-линейных комбинаций элементов $($соответственно, классовых сумм$\,)$ группы $G$.
\end{pre}

\textupl{kvt prf}{Доказать предложение \ref{kvt}.}

\begin{opr}
Пусть  $A$ и $B$ --- $R$-алгебры. Произвольное $R$-линейное отображение из $A$ в $B$, являющееся
гомоморфизмом колец, называется \glsadd{iHomRAlg}\mem{гомоморфизмом $R$-алгебр}.
\end{opr}

\textzam{Теория представлений группы $G$ изучает гомоморфизмы
$R$-алгебры $RG$ в алгебру $\MM_n(R)$.}

Пусть $I$  --- (двусторонний) идеал $R$-алгебры $A$. Легко проверить, что структуры факторкольца и
фактормодуля $A/I$ согласованы и определяют $R$-алгебру, называемую
\glsadd{iQuoAlg}\mem{факторалгеброй} $A$ по $I$. Отметим, что ядро $\Ker \vf$ любого гомоморфизма $R$-алгебр $\vf:A\to B$
является двусторонним идеалом в $A$, и имеет место изоморфизм $R$-алгебр $A/\Ker
\vf\cong \Im \vf$.

\subsection{Примеры гомоморфизмов \texorpdfstring{$R$}{R}-алгебр}\label{ex homr}
\begin{list}{{\rm(}{\it\roman{enumi}\/}{\rm)}}
{\usecounter{enumi}\setlength{\parsep}{2pt}\setlength{\topsep}{5pt}\setlength{\labelwidth}{23pt}}
\item Имеет место классический изоморфизм $F$-алгебр $\End_F(V)\cong\MM_n(F)$,
где $V$ --- $n$-мерное векторное пространство над полем $F$.

\item Пусть $Y$ --- подмножество множества $X$ и $R$ --- коммутативное кольцо. Тогда
отображение $\f_R(X)\to \f_R(Y)$, сопоставляющее каждой функции из $\f_R(X)$ её
ограничение на $Y$, является эпиморфизмом $R$-алгебр. Его ядро состоит из функций,
тождественно равных нулю на $Y$.

\item Для произвольного множества $X$ имеет место изоморфизм $\ZZ_2$-алгебр
$\P(X)\cong\f_{\ZZ_2}(X)$, см. \ref{ex r-alg}$(vi)$--$(vii)$.

\item Пусть $G$ --- группа, $H\nor G$ и $\vf: G \to G/H$ --- естественный
эпиморфизм групп. Тогда $\vf$ продолжается по  линейности до эпиморфизма групповых алгебр $RG\to
R(G/H)$. В случае, когда $R$ --- поле,
ядро этого эпиморфизма имеет размерность $(|H|-1)|G:H|$, а базисом этого ядра является множество
$$\big\{\,(1-h)x\bigm|x\in X, h\in H, h\ne 1\,\big\},$$
где $X$ --- полная система представителей смежных классов $G$ по $H$.
\end{list}

\textupr{Проверить сформулированные в \ref{ex homr} утверждения.}

\textuprn{Пусть $I\nor RG$ и $\vf: (RG)^\times \to (RG/I)^\times$ --- гомоморфизм групп,
являющийся ограничением  на $(RG)^\times$ естественного эпиморфизма $R$-алгебр $RG \to RG/I$. Показать, что
\begin{list}{{\rm(}{\it\roman{enumi}\/}{\rm)}}
{\usecounter{enumi}\setlength{\parsep}{2pt}\setlength{\topsep}{5pt}\setlength{\labelwidth}{23pt}}
\item $\ker \vf=\{a\in (RG)^\times\mid 1-a\in I\}$;
\item $\{a\in G\mid 1-a\in I\}\nor G$.
\end{list}
}

\begin{opr} Пусть
$A$ --- $R$-алгебра и $V$ --- правый $A$-модуль (где $A$ рассматривается как кольцо). Наделим $V$ структурой
левого $R$-модуля, положив $rv=v(r1_A)$ для произвольных $v\in V$ и $r\in R$ (это можно сделать в силу
коммутативности кольца $R$). Такой $A$-модуль $V$ с дополнительной структурой левого $R$-модуля будем называть
\glsadd{iModRAlg}\mem{модулем над алгеброй} $A$ или просто \glsadd{iAMod}\mem{$A$-модулем}.
Умножение элементов модуля $V$ на элементы алгебры $A$ мы будем называть
\glsadd{iActAlg}\mem{действием алгебры} $A$.

\end{opr}

\medskip\noindent\llap{\PencilRightDown\quad}\ignorespaces
\mysubsection{Пусть $A$ --- $R$-алгебра и $V$ --- $A$-модуль. Показать, что для любых
$v\in V$, $r\in R$, $a\in A$ выполнено
$$
r(va)=(rv)a=v(ra).
$$}

Из определения $A$-модуля $V$ следует существование естественного гомоморфизма $R$-алгебр $A\to \End_R(V)$, определяемого
правилом $a\mapsto a_V$, где $a_V\in \End_R(V)$ --- отображение $\gls{aV}:v\mapsto va$. Обратно, если $V$ --- левый $R$-модуль,
то для любого гомоморфизма $R$-алгебр $\vf:A\to \End_R(V)$ можно задать на $V$ структуру $A$-модуля, положив $va=v(a\vf)$
для любых $v\in V$, $a\in A$.

\textzam{В случае, когда $R$ --- поле, произвольный $A$-модуль $V$ является векторным пространством над $R$ и, в
частности, для него определены такие понятия как базис над $R$, размерность над $R$, и т.\,д.}

\subsection{Примеры \texorpdfstring{$A$}{A}-модулей} \label{pr am}
\begin{list}{{\rm(}{\it\roman{enumi}\/}{\rm)}}
{\usecounter{enumi}\setlength{\parsep}{2pt}\setlength{\topsep}{5pt}\setlength{\labelwidth}{23pt}}
\item \glsadd{iAModReg}\mem{Регулярный $A$-модуль} \gls{A0} --- сама алгебра $A$,
действующая на себе умножением справа.

\item Пусть $R$  --- коммутативное кольцо и $V$ --- левый $R$-модуль. Тогда $V$ является правым
модулем над $R$-алгеброй $\End_R(V)$ относительно естественного действия.

\item Пусть $S$ --- кольцо. Тогда любой правый $S$-модуль можно рассматривать как модуль над
$\ZZ$-алгеброй $S$ или как модуль над $R$-алгеброй $S$, где $R\le \Z(S)$.

\item Пусть $R$  --- коммутативное кольцо и $H$ --- некоторая подгруппа группы $G$. Поскольку алгебра $RH$ является
подалгеброй алгебры $RG$, то каждый $RG$-модуль $V$ будет также и $RH$-модулем, который мы будем обозначать
символом \gls{VH}.

\item Пусть $G$ --- группа, $F$ --- поле и $V$ --- векторное пространство конечной размерности $n$ над
$F$. Пусть действие групповой алгебры $FG$ на $V$ задаётся по правилу $vg=v$ для всех $g\in G$, $v\in V$ и
продолжается по линейности на всю алгебру $FG$. Таким образом на $V$ определяется структура $FG$-модуля. Этот
модуль мы будем называть \mem{тривиальным}\glsadd{iModRGTriv} $n$-мерным
\mem{модулем} алгебры $FG$. Одномерный тривиальный $FG$-модуль называется \glsadd{iModRGPrinc}\mem{главным}.

\item Пусть $G$ --- группа, $F$ --- поле и $V$ --- $FG$-модуль. На двойственном
\footnote{Т.е. на пространстве всех $F$-линейных отображений $\vf:V\to F$, см. \ref{dv mod}.} векторном пространстве
$V^*$ можно задать структуру $FG$-модуля, если положить для всех $\vf\in V^*$, $g\in G$
$$v(\vf g)=(vg^{-1})\vf,$$
где  $v\in V$, и продолжить это действие по линейности на всю
алгебру $FG$. Такой $FG$-модуль \gls{VCGrad} называется
\glsadd{iFGModCGrad}\mem{контрагредиентным}\footnote{Мы намеренно избегаем употребления для $V^*$
термина <<двойственный $FG$-модуль>>, зарезервированного нами для левого $FG$-модуля $\Hom(V,FG^\circ)$.
Обозначение $V^*$ будет в дальнейшем использоваться именно для контрагредиентного $FG$-модуля.}
к модулю $V$.

\item Пусть группа $G$ действует на конечном множестве $X$. Пусть $F$ --- поле,
 $V$ --- векторное пространство размерности $|X|$ над $F$ с базисом $\{e_x\mid x\in X\}$.
 Действие групповой алгебры $FG$ на $V$ задаётся действием групповых элементов на базисных
векторах по правилу $e_xg=e_{xg}$ для всех $g\in G$, $x\in X$, которое можно продолжить по линейности на всё
пространство $V$:
$$
\left(\sum_{\,x\in X}a_xe_x\right)g = \sum_{x\in X}a_xe_{xg}= \sum_{x\in X}a_{xg^{-1}}e_x,
$$
а затем на всю алгебру $FG$. Таким образом на $V$ определяется структура $FG$-модуля, называемого
\mem{подстановочным модулем}\glsadd{iModRGPerm}  групповой алгебры $FG$, соответствующим данному
действию группы на множестве. В частном случае, когда $G$~--- подгруппа симметрической группы $S_n$
подстановок множества $X=\{1,\ld,n\}$, определённый таким образом $FG$-модуль будем называть
\glsadd{iModRGPermNat}\mem{естественным подстановочным модулем}.
\end{list}

\textupr{Проверить сформулированные в \ref{pr am} утверждения.}

\textuprn{\label{reg perm}Пусть $F$ --- поле и $G$ --- группа. Показать, что регулярный модуль $FG^\circ$ изоморфен подстановочному
$FG$-модулю, соответствующему действию группы $G$ на себе правыми умножениями (так называемому
\glsadd{iActGrRReg}\mem{правому регулярному действию} группы $G$).}

\mysubsection Пусть $A$ --- $R$-алгебра и $V$ --- $A$-модуль. Тогда $A$-подмодуль  $W$ правого
$A$-модуля $V$ (здесь $A$ рассматривается как кольцо) автоматически
будет левым $R$-подмодулем, поскольку $rw=w(r1_A)\in W$ для всех
$r\in R$ и $w\in W$. Поэтому такие подмодули будут $A$-модулями над
алгеброй $A$. Это замечание устраняет двусмысленность термина
 <<\mem{$A$-подмодуль}>>\glsadd{iASubMod} для $A$-модуля $V$. Аналогично для $A$-подмодуля $W\le V$ фактормодуль $V/W$
будет $A$-модулем.

Легко видеть, что подмодули регулярного $A$-модуля $A^\circ$ --- это, в точности, правые идеалы алгебры $A$.
Более общо, правые идеалы алгебры $A$ позволяют строить подмодули в произвольном $A$-модуле $V$. Пусть $I\nor_r A$.
Как обобщение введённого в \ref{s p ids} произведения правых идеалов, обозначим через
$VI$ аддитивную подгруппу группы $V$,
порождённую всевозможными произведениями вида $vx$, где $v\in V$, $x\in I$. Ясно, что $VI$~---
$A$-подмодуль модуля $V$.

\textuprn{\label{diag aug}
Пусть группа $G$ действует на множестве $X$ и $V$ --- подстановочный $FG$-модуль,
определённый в примере \ref{pr am}$(vii)$. Введём подпространства
$$U=F\left(\sum_{\,x\in X}e_x\right) \quad \text{и} \quad W=\left\{\, \sum_{x\in X} a_xe_x\ \big|\  a_x\in F, \sum_{x\in X} a_x=0\, \right\}$$
пространства $V$. Показать, что справедливы утверждения.
\begin{list}{{\rm(}{\it\roman{enumi}\/}{\rm)}}
{\usecounter{enumi}\setlength{\parsep}{2pt}\setlength{\topsep}{5pt}\setlength{\labelwidth}{23pt}}
\item $U$ и $W$ являются $FG$-подмодулями модуля $V$.
\item В случае регулярного действия $G$ подмодули $U$ и $W$ являются идеалами алгебры $FG$.
\item Включение $U\le W$ справедливо тогда и только тогда, когда характеристика поля $F$ делит $|X|$.
\item Пусть $Y\se X$ --- объединение $G$-орбит. Тогда подпространства
$$
V_Y=\sum_{y\in Y}Fe_y \qquad \text{и} \qquad V_{X\setminus Y}=\sum_{x\in X\setminus Y}Fe_x
$$
являются $FG$-подмодулями модуля $V$ и справедливо разложение $V=V_Y\oplus V_{X\setminus Y}$.
\end{list}}

В обозначениях из \ref{diag aug} будем называть $U$ \glsadd{iSubModPermModDiag}\mem{диагональным},
а $W$ --- \glsadd{iSubModPermModDiff}\mem{разностным} подмодулями
подстановочного $FG$-модуля $V$. Очевидно, что $U$ изоморфен главному $FG$-модулю.

\begin{pre} \label{smods} Пусть $F$-поле.  Всякий
собственный ненулевой $FS_n$-подмодуль естественного подстановочного $FS_n$-модуля, см. {\rm \ref{pr am}}$(vii)$,
совпадает либо с диагональным, либо с разностным подмодулем.
\end{pre}
\textupr{Доказать предложение \ref{smods}.
\uk{Всякое подпространство векторного пространства задаётся однородной системой линейных уравнений.}}

\textupr{Показать, что утверждение, аналогичное \ref{smods}, справедливо для
естественного подстановочного $FA_n$-модуля при $n\ge 5$. Привести контрпримеры при $n=2,3,4$.}

Пусть $\vf: V\to W$ --- гомоморфизм $A$-модулей (рассматриваемых как правые модули над кольцом $A$). Тогда
$\vf$ будет также гомоморфизмом левых $R$-модулей (то есть $R$-линейным), поскольку
$$(rv)\vf=(v(r1_A))\vf=(v\vf)(r1_A)=r(v\vf)$$
для любых $r\in R$ и $v\in V$. Гомоморфизмы (эндоморфизмы, изоморфизмы) $A$-модулей мы будем называть также
\glsadd{iAHom}\mem{$A$-гомоморфизмами} (\glsadd{iAEnd}\mem{$A$-эндоморфизмами}, \glsadd{iAIso}\mem{$A$-изоморфизмами}).

Таким образом, $\Hom_A(V,W)$~--- множество всех $A$-гомоморфизмов из $V$ в $W$, а $\End_A(V)$~---
всех $A$-эн\-до\-мор\-физ\-мов модуля $V$.
 Если рассматривать $A$-модули $V$ и $W$ как левые $R$-модули, то $\Hom_A(V,W)$
состоит из тех элементов $\vf\in
\Hom_R(V,W)$, которые перестановочны с действием алгебры $A$, т.\,е. таких, что для любых $v\in V$ и $a\in A$
выполнено $(va)\vf=(v\vf)a$. В связи с этим отметим, что эндоморфизм регулярного $A$-модуля $A^\circ$ может не
быть эндоморфизмом алгебры $A$ и наоборот.

Группу $\Hom_A(V,W)$ легко превратить в левый $R$-модуль, положив $v(r\vf)=r(v\vf)$ для любых
$r\in R$ и $\vf\in \Hom_A(V,W)$. Композиция отображений задаёт на $\End_A(V)$ операцию
умножения.

\textuprn{\label{end mod} Пусть $A$ --- $R$-алгебра
и $V$ --- $A$-модуль. Показать, что $\End_A(V)$ является подалгеброй $R$-алгебры $\End_R(V)$. В частности, $V$
можно рассматривать как правый
$\End_A(V)$-модуль.}

\textuprn{Пусть $R$ --- коммутативное кольцо и $V$ --- левый $R$-модуль. Положим $A=\End_R(V)$. Показать, что
$\End_A(V)=\Z(A)$.}

\textuprn{\label{doub cont} Пусть $F$ --- поле и $G$ --- конечная группа. Показать, что для любого  $FG$-модуля $V$ конечной размерности над $F$
имеет место изоморфизм $V^{**}\cong V$.}

\normalmarginpar

В следующей теореме собраны некоторые факты об $A$-модулях, аналогичные теоремам о гомоморфизмах для групп.

\begin{thm}[о гомоморфизмах $A$-модулей] \label{thm o hom}\mbox{}
\begin{list}{{\rm(}{\it\roman{enumi}\/}{\rm)}}
{\usecounter{enumi}\setlength{\parsep}{2pt}\setlength{\topsep}{5pt}\setlength{\labelwidth}{23pt}}
\item Пусть $\vf:V\to W$ --- гомоморфизм $A$-модулей $V$ и $W$. Тогда $\Ker \vf\le V$, $\Im\vf\le W$ и
отображение $v+\Ker
\vf\mapsto v\vf$ для $v\in V$ определяет изоморфизм $A$-модулей $V/\Ker \vf\cong\Im \vf$.
\item Если имеются включения $A$-модулей $U\le W \le V$ и $M\le V$, то
$$
(W+M)/(U+M) \cong W/(U+\,W\cap M).
$$
В частности,
$$
(W+M)/M \cong W/(W\cap M).
$$
\item Пусть $U$ --- $A$-подмодуль $A$-модуля $V$. Тогда отображение $M\mapsto M/U$ задаёт
биекцию между $A$-подмодулями из $V$, содержащими $U$, и $A$-подмодуля\-ми $A$-модуля $V/U$. Эта биекция
сохраняет отношение включения. При этом, если $U\le M\le V$, то
$$
(V/U)/(M/U)\cong V/M.
$$
\end{list}
\end{thm}

\textupr{Доказать теорему \ref{thm o hom}.}

\begin{opr} Ненулевой\footnote{Нулевой модуль не является ни приводимым, ни неприводимым.}
$A$-модуль $V$ называется\footnote{Данное понятие неприводимости будет нами использоваться в основном для алгебр над полями.
В случае же, когда $R$ не является полем, это понятие может оказаться
не вполне удобным с точки зрения теории представлений, поскольку даже само кольцо $R$ не будет неприводимым как $R$-модуль.
Скажем, для изучения целочисленных
представлений (т.\,е. ситуации, когда $R=\ZZ$) более подходящие аналоги неприводимости определены в \cite[Глава XI]{cr}.}
\glsadd{iAModIrr}\mem{неприводимым} или \glsadd{iAmodSmp}\mem{простым},
если он не содержит подмодулей, отличных от $0$ и $V$. В противном случае $V$ называется \glsadd{iAModRed}\mem{приводимым}.
\end{opr}

Очевидно, что минимальные подмодули $A$-модуля $V$ --- это, в точности, неприводимые подмодули, а
максимальные --- это, те подмодули, фактор по которым неприводим.

\textuprn{Показать, что регулярный $A$-модуль $A^\circ$ неприводим тогда и только тогда, когда $A$ является телом.}

\textuprn{Пусть $R$-коммутативное кольцо. Показать, что следующие условия эквивалентны.
\begin{list}{{\rm(}{\it\roman{enumi}\/}{\rm)}}
{\usecounter{enumi}\setlength{\parsep}{2pt}\setlength{\topsep}{5pt}\setlength{\labelwidth}{23pt}}
\item Любой ненулевой левый $R$-модуль $V$ неприводим как правый $\End_R(V)$-модуль.
\item $R$ является полем.
\end{list}
\uk{Рассмотреть в качестве $V$ регулярный модуль $^\circ\! R$.}}

%

\begin{pre} \label{kon npr}
Пусть V --- FG-модуль и W --- его подмодуль. Доказать, что 
контрагредиентный модуль V* имеет подмодуль, изоморфный (V/W)* и фактормодуль, изоморфный W*.
В частности, V неприводим тогда и только тогда, когда V* неприводим.
\end{pre}
\textupr{Доказать предложение \ref{kon npr}. 
\uk{Рассмотреть множество 
$$
\{\vf\in V^*\,|\ v\vf=0 \ \mbox{для всех}\ v\in W \}.
$$
Воспользоваться упражнением \ref{doub cont}.
}}

Всякий неприводимый $A$-модуль, очевидно, неразложим. Обратное неверно, как, с учётом \ref{diag aug}, показывает следующее
упражнение.

\textuprn{\label{fsn}Пусть $V$ --- естественный подстановочный $FS_n$-модуль, $U$ --- его диагональный, $W$ --- разностный подмодули, и $n\ne 1,2$.
Доказать эквивалентность следующих утверждений.
\begin{list}{{\rm(}{\it\roman{enumi}\/}{\rm)}}
{\usecounter{enumi}\setlength{\parsep}{2pt}\setlength{\topsep}{5pt}\setlength{\labelwidth}{23pt}}
\item Модуль $V$ разложим.
\item Подмодуль $W$ неприводим.
\item Характеристика поля $F$ не делит $n$.
\end{list}
Показать, что  всякое разложение $V$ в прямую сумму собственных подмодулей совпадает с разложением
$V=U\oplus W$. Рассмотреть также случаи $n=1,2$.}

Конечный ряд $A$-модулей

$$
V=V_0\ge V_1\ge \ld \ge V_n=0 \myeqno\label{komr}
$$
называется \glsadd{iSerCompAMod}\mem{композиционным рядом} $A$-модуля $V$, если все
факторы $V_{i-1}/V_i$, $i=1,\ld,n$, называемые \glsadd{iQuoSerAMod}\mem{факторами ряда} \ref{komr},
неприводимы.

В случае, когда $R$
--- поле, примером $A$-модуля, обладающего композиционным рядом,
служит любой $A$-модуль $R$-алгебры $A$, имеющий как векторное пространство конечную размерность над $R$. В
общем же случае $A$-модуль  может не иметь композиционного ряда.

\textuprn{Привести пример $A$-модуля, не обладающего композиционным рядом.}

Два композиционных ряда некоторого $A$-модуля называются \glsadd{iSersCompEqv}\mem{эквивалентными}, если
они имеют одинаковое число факторов, и между этими факторами можно установить взаимно однозначное соответствие, при котором
соответствующие факторы будут изоморфны.

Следующая теорема, утверждающая единственность с точностью до эквивалентности композиционного ряда $A$-модуля
по существу вытекает из теоремы о гомоморфизмах \ref{thm o hom}.

\begin{thm}[Жордана--Гёльдера] \label{thm zh-g}\glsadd{iThmJorHol} Любые два
композиционных ряда $A$-модуля $V$ эквивалентны.
\end{thm}

\textupr{Доказать теорему \ref{thm zh-g}.}

 Теорема Жордана--Гёльдера
позволяет ввести важное понятие композиционных факторов для
$A$-модулей, обладающих композиционным рядом.

\begin{opr} \glsadd{iFctCmpAMod}\mem{Композиционными факторами} $A$-модуля $V$ называются факторы его композиционного ряда.
\end{opr}

Как вытекает из теоремы \ref{thm o hom}, наборы композиционных факторов изоморфных $A$-модулей по существу
совпадают. Следующее упражнение показывает, что есть примеры неизоморфных $A$-модулей с одинаковыми наборами
композиционных факторов.

\textuprn{Показать, что регулярный и тривиальный двумерный модули групповой алгебры $\FF_2S_2$ неизоморфны,
но имеют один и тот же набор композиционных факторов.}

\textuprn{Найти размерности композиционных факторов естественных подстановочных $FS_n$ модулей.}

Теорема Жордана-Гёльдера показывает важность изучения неприводимых модулей. Одним из ключевых фактов о
неприводимых модулях является следующее утверждение.

\begin{thm}[лемма Шура] \label{thm lsh}\glsadd{iLemSchr} Пусть $V,W$ --- неприводимые $A$-модули.
Тогда любой не\-нулевой элемент из $\Hom_A(V,W)$ является изоморфизмом.
\end{thm}

\textupr{Доказать теорему \ref{thm lsh}.}

Очевидным следствием из леммы Шура является следующее утверждение.

\begin{cor} \label{cor lsh}
Если $V$ --- неприводимый $A$-модуль, то алгебра $\End_A(V)$ является телом.
\end{cor}

Отметим важный частный случай.

\begin{cor} \label{cor sh} Пусть $F$ --- алгебраически замкнутое поле, $A$ --- $F$-алгебра,
$V$ --- неприводимый $A$-модуль и размерность $\dim_F V$ конечна. Тогда имеет место
изоморфизм $F$-алгебр $\ \End_A(V)\cong F$.
\end{cor}

\begin{proof} Очевидно, что $F1\le \End_A(V)$.
Рассмотрим произвольный элемент $\th\in \End_A(V)$. Тогда $\th$ является линейным преобразованием
конечномерного векторного пространства $V$ над алгебраически замкнутым полем $F$ и, значит, имеет
собственное значение $\l$. Следовательно эндоморфизм $\th-\l 1\in
\End_A(V)$ необратим. По лемме \ref{thm lsh} получаем $\th-\l 1=0$, т.\,е. $\th=\l 1\in F 1$.
Таким образом, $\End_A(V)=F1$, а так как $F$ --- поле, имеем $F1\cong F$.
\end{proof}

Покажем, что по существу источником всех неприводимых $A$-модулей является сама алгебра $A$. Для этого,
введём следующие понятия.

Пусть $V$ --- $A$-модуль. \mem{Аннулятором элемента}\glsadd{iAnnElmAMod}
$v\in V$ называется множество
$$\gls{Annlvr}=\big\{a\in A\bigm|va=0\big\}.$$

\mem{Аннулятором $A$-модуля}\glsadd{iAnnAmod} $V$ называется
$$
\gls{AnnlVr}=\bigcap_{v\in V} \Ann(v)=\big\{a\in A\bigm|va=0\ \ \mbox{для всех}\ \ v\in V\big\}.
$$

Отметим некоторые элементарные свойства аннуляторов.

\begin{pre} \label{ann pr}
Пусть
$V$ --- $A$-модуль и $v\in V$. Имеют место следующие утверждения.
\begin{list}{{\rm(}{\it\roman{enumi}\/}{\rm)}}
{\usecounter{enumi}\setlength{\parsep}{2pt}\setlength{\topsep}{5pt}\setlength{\labelwidth}{23pt}}
\item $\Ann(v)\nor_r A$.
\item $\Ann(V)\nor A$.
\item Из изоморфизма $A$-модулей $V\cong W$ следует равенство $\Ann(V)=\Ann(W)$.
\item Для любого идеала $I\nor A$ такого, что $I\le \Ann(V)$, модуль $V$ можно рассматривать
как $A/I$-модуль, если для $v\in V$ и $a\in A$ положить $v(a+I)=va$. При этом модуль $V$ неприводим
тогда и только тогда, когда он неприводим как $A/I$-модуль.
\end{list}
\end{pre}
\textupr{Доказать предложение \ref{ann pr}.}

Выше мы отмечали, что подмодуль максимален тогда и только тогда, когда фактор по нему неприводим. Сейчас
покажем, что любой неприводимый $A$-модуль может быть получен как фактормодуль регулярного модуля $A^\circ$ по
некоторому его  максимальному подмодулю.

\begin{pre}\label{ann el} Пусть $V$ --- неприводимый $A$-модуль $R$-алгебры $A$ и $v$ --- ненулевой элемент из $V$.
Тогда
\begin{list}{{\rm(}{\it\roman{enumi}\/}{\rm)}}
{\usecounter{enumi}\setlength{\parsep}{2pt}\setlength{\topsep}{5pt}\setlength{\labelwidth}{23pt}}
\item $\Ann(v)$ --- максимальный правый идеал в $A$;
\item $V\cong A^\circ/\Ann(v)$.
\end{list}
\end{pre}
\textupl{ann el prf}{Доказать предложение \ref{ann el}.}

\section{Радикал}

В дальнейшем мы зафиксируем некоторую систему представителей всех классов изоморфных неприводимых $A$-модулей и
будем обозначать её через \gls{MlAr}.

Теперь мы можем ввести следующее важное понятие.

\begin{opr}  \glsadd{iRadRAlg}\mem{Радикалом} $R$-алгебры $A$, обозначаемым через \gls{JlAr}, называется\footnote{В
литературе радикал алгебры также часто называется \glsadd{iRadJac}\mem{радикалом Джекобсона}.} пересечение аннуляторов всех
неприводимых $A$-модулей.
\end{opr}

Ввиду предложения \ref{ann pr}$(iii)$ можно записать
$$
\J(A)=\bigcap_{V\in \M(A)}\Ann(V).
$$
Из  \ref{ann pr}$(ii)$ следует также, что $\J(A)\nor A$.

\textext{Строго говоря, определённый таким образом радикал было бы правильнее назвать <<правым
радикалом>>, поскольку речь идёт о правых $A$-модулях. Однако, можно проверить, что если ввести аналогичное
понятие <<левого радикала>>, то правый и левый радикалы алгебры совпадут, см. \ref{lrad}.}

\begin{pre} \label{rad int max}
Радикал $\J(A)$ совпадает
 с пересечением всех максимальных правых идеалов алгебры~$A$.
\end{pre}

\begin{proof}
Пусть $J$ --- пересечение множества всех максимальных правых идеалов алгебры $A$. По
предложению \ref{ann el} множество $J$ содержится в аннуляторах $\Ann(v)$ всех ненулевых
элементов всех неприводимых $A$-модулей $V$. Поэтому
$$
J\se \bigcap_{V\in\, \M(A)}\ \bigcap_{0\ne v\in V} \Ann(v) = \bigcap_{V\in\, \M(A)} \Ann(V) =
\J(A).
$$
Для обратного включения достаточно показать, что любой максимальный правый идеал $I\nor_r A$ содержит аннулятор
некоторого неприводимого $A$-модуля. В качестве такого модуля можно взять фактормодуль $A^\circ/I$. В самом
деле, если $x\in \Ann(A^\circ/I)$, то, в частности, $x$ аннулирует элемент $1_A+I\in A^\circ/I$ и мы получаем
$I=(1_A+I)x=x+I$, т.\,е. $x\in I$.
\end{proof}

Доказанное предложение даёт <<внутреннюю>> и с практической точки зрения более удобную, чем
исходное определение, характеризацию радикала $\J(A)$ в терминах правых идеалов алгебры $A$.

\begin{pre} \label{dir sum} Пусть $A_1$, $A_2$ --- $R$-алгебры и $A=A_1\oplus A_2$.
Тогда $\J(A)=\J(A_1)\oplus\J(A_2)$.
\end{pre}

\textupl{dir sum prf}{Доказать предложение \ref{dir sum}.}

Зафиксируем кольцо $S$. Как видно из \ref{ex r-alg}$(iv),(v)$, это кольцо можно различными
способами превратить в алгебру над подходящим коммутативным кольцом $R$. Естественно возникает вопрос,
как соотносятся радикалы алгебр, возникающих из $S$. На самом деле эти радикалы совпадают. Это видно,
например, из предложения \ref{rad int max}, поскольку совпадают множества максимальных правых идеалов
таких алгебр (ведь по определению правые идеалы $R$-алгебры $S$ --- это, в точности, правые идеалы
кольца $S$). Таким образом, если определить \glsadd{iRadRng}\mem{радикал} \gls{JlSr} \mem{кольца} $S$  как пересечение  всех максимальных правых идеалов из $S$, то это
определение будет согласовано со определением радикала $R$-алгебры.

Следующие утверждения показывают связь между принадлежностью элемента $a$ алгебры $A$ её радикалу и обратимостью элемента $1-a$.

\begin{pre}\label{rad prop} Если $A$ --- $R$-алгебра и $a\in \J(A)$, то
элемент $1-a$ обратим.
\end{pre}

\begin{proof} Рассмотрим правый идеал $(1-a)A$ алгебры $A$. Если он собственный, то по предложению
\ref{max id}
содержится в некотором максимальном правом идеале $M$, и поэтому $1-a\in (1-a)A\le M$. Но с
другой стороны $a\in \J(A)\le M$ и, значит, $1\in M$. Противоречие. Поэтому $(1-a)A=A$ откуда,
очевидно, следует правая обратимость элемента $1-a$. По предложению \ref{obr}$(iii)$ элемент $1-a$
также обратим слева.
\end{proof}

\begin{pre}\label{id rad} Если $I\nor_r A$
и для любого $a\in I$ элемент $1-a$ обратим справа, то $I\se \J(A)$.
\end{pre}

\textupl{id rad prf}{Доказать предложение \ref{id rad}.}

\textuprn{\label{lrad}Показать, что <<левый радикал>> алгебры $A$ (т.\,е. пересечение аннуляторов всех неприводимых левых $A$-моду\-лей)
совпадает с $\J(A)$.}


\textuprn{Доказать, что $\J(A)$ не содержит ненулевых идемпотентов.}



Если $e$ --- идемпотент $R$-алгебры $A$, то $eAe$ --- $R$-алгебра с единицей $e$ в силу \ref{ids 2}$(ii)$.
Следующее утверждение описывает радикал алгебры $eAe$.

\begin{pre} \label{jre} Пусть $A$ --- $R$-алгебра и $e\in A$ --- идемпотент.
Тогда $$\J(eAe)=e\J(A)e=\J(A)\cap eAe.$$
\end{pre}
\textupl{jre prf}{Доказать предложение \ref{jre}.}

\begin{pre} \label{rad max}\mbox{}
\begin{list}{{\rm(}{\it\roman{enumi}\/}{\rm)}}
{\usecounter{enumi}\setlength{\parsep}{2pt}\setlength{\topsep}{5pt}\setlength{\labelwidth}{23pt}}
\item Любой правый идеал $R$-алгебры $A$, состоящий из нильпотентных элементов, лежит в радикале $\J(A)$.
\item Любой нильпотентный правый идеал $R$-алгебры $A$ лежит в радикале $\J(A)$.
\end{list}
\end{pre}

\textupl{rad max prf}{Доказать предложение \ref{rad max}.}

\begin{pre} \label{com nil}
В коммутативной $R$-алгебре множество всех нильпотентных элементов образует
идеал.
В частности, любой нильпотентный элемент лежит в радикале.
\end{pre}

\textupl{com nil prf}{Доказать предложение \ref{com nil}.}

$A$-модуль $V$ называется \glsadd{iModFinGen}\mem{конечно порождённым}, если он имеет конечное
число \glsadd{iGenAMod}\mem{порождающих}, т.\,е. таких элементов $v_1,\ld,v_n\in V$, что
$$
V=v_1A+\ld+v_nA.
$$

\begin{pre} \label{con por} Любой ненулевой конечно порождённый $A$-модуль содержит
максимальный подмодуль.
\end{pre}
\textupl{con por prf}{Доказать предложение \ref{con por}. \uk{Использовать лемму Цорна \ref{thm lz}.}}

\glsadd{iLemNak}
\begin{thm}[лемма Накаямы] \label{thm nak} Пусть $V$ --- $A$-модуль и $W\le V$.
Предположим, что фактормодуль $V/W$ конечно порождён.
Если $W+V\J(A)=V$, то $W=V$.
\end{thm}

\begin{proof} Достаточно доказать лемму для случая, когда $W=0$ и затем применить утверждение к
фактормодулю $V/W$. Итак, считаем, что $V$ --- конечно порождённый $A$-модуль и $V\J(A)=V$.
Требуется доказать, что $V=0$. Допустим, что $V\ne 0$. Тогда по предложению \ref{con por}
модуль $V$ содержит максимальный подмодуль $U$. Фактормодуль $V/U$ неприводим, поэтому
$(V/U)\J(A)=0$, откуда следует, что $V\J(A)\le U\ne V$. Противоречие. \end{proof}


Положив $W=0$ в теореме \ref{thm nak}, получаем важное следствие (также часто называемое леммой Накаямы).

\begin{cor} \label{nak kl} Если $V$ --- ненулевой конечно порождённый $A$-модуль, то $V\J(A)<V$.
\end{cor}

Отметим, что требование конечной порождённости модуля $V$ в лемме Накаямы является существенным. Например, в
качестве $V$ можно взять поле $\QQ$ рациональных чисел, а в качестве $A$ --- его подкольцо $\ZZ[\pi^{-1}]$, где
$\pi$ --- множество нечётных простых чисел, см. \ref{subrq}. Тогда $V$ не является конечно порождённым
$A$-модулем. Однако, $\J(A)=2A$ и имеет место равенство $V\J(A)=V$.

В качестве примера использования леммы Накаямы докажем следующее утверждение.

\begin{pre} \label{nrz} Пусть $S$ --- кольцо
и $R$ --- подкольцо центра $\Z(S)$. Предположим, что $S$ является конечно порождённым как $R$-модуль. Тогда
$\J(R)\le \J(S)\cap R$.

\end{pre}
\begin{proof} Рассмотрим $S$ как $R$-алгебру. Достаточно показать, что для произвольного
неприводимого $S$-модуля $V$ справедливо равенство $V\J(R)=0$. Поскольку $S$ конечно порождено как $R$-модуль,
можно записать $S=s_1R+\ld +s_nR$ для некоторых $s_i\in S$. Выберем ненулевой элемент $v\in V$. Тогда
$$V=vS=vs_1R+\ld +vs_nR,$$
откуда следует, что $V$ --- (ненулевой) конечно порождённый $R$-модуль. Из \ref{nak kl} следует, что
$V\J(R)<V$. Однако $V\J(R)$ является $S$-подмодулем модуля $V$, поскольку $\J(R)$ состоит из
центральных элементов кольца $S$. Значит, $V\J(R)=0$ в силу неприводимости $V$.
\end{proof}


Отметим, что в условиях предложения \ref{nrz} справедливо и обратное включение, т.\,е. имеет место равенство
$\J(R)=\J(S)\cap R$ (в случае, когда $S$ --- конечномерная алгебра над полем, см. \ref{rad cent}.

Следующее вспомогательное утверждение нам потребуется в дальнейшем.

\begin{pre} \label{l54} Пусть $G$ --- конечная группа и $\vf:R\to S$ --- гомоморфизм коммутативных колец такой,
что $\Ker \vf\se \J(R)$.
Рассмотрим естественное поднятие $\vf$ до кольцевого гомоморфизма
$\ti\vf:RG\to SG$ групповых алгебр. Пусть $e\in RG$ --- идемпотент и $x\in e(RG)e$. Если
$x\ti\vf=e\ti\vf$, то $x$ обратим в кольце $e(RG)e$.
\end{pre}
\begin{proof} Обозначим $A=e(RG)e$. Из \ref{ids 2}$(ii)$ вытекает, что  $A$ действительно является кольцом
и $e$ --- его единицей. По \ref{rad prop} достаточно показать, что $e-x\in \J(A)$. Из \ref{jre} следует, что
$\J(A)=\J(RG)\cap A$ и, значит,  достаточно проверить, что $e-x\in \J(RG)$. По условию $e-x\in \Ker\ti\vf$ и
мы покажем, что $\Ker\ti\vf\se \J(RG)$, откуда будет следовать требуемое.

В силу \ref{kvt} любой элемент $u\in \Ker\ti\vf$ имеет вид $u=\sum_{g\in G}u_gg$ для подходящих $u_g\in
\Ker\vf$. Значит, ввиду замечания \ref{1ot}, достаточно показать, что $\Ker\vf\se \J(RG)$. По условию $\Ker
\vf\se \J(R)$. Кроме того, $R\le \Z(RG)$ и $RG$ является конечно порождённым $R$-модулем. Поэтому из \ref{nrz}
следует, что  $\J(R)\se\J(RG)$ и, тем самым, утверждение доказано.
\end{proof}

Кольцо называется \glsadd{iRngLoc}\mem{локальным}, если оно содержит единственный
максимальный правый идеал.


\begin{pre} \label{loc eq} Пусть $S$ --- кольцо. Тогда следующие условия эквивалентны.
\begin{list}{{\rm(}{\it\roman{enumi}\/}{\rm)}}
{\usecounter{enumi}\setlength{\parsep}{2pt}\setlength{\topsep}{5pt}\setlength{\labelwidth}{23pt}}
\item $S$ содержит единственный максимальный правый  идеал $I_r$.
\item $S$ содержит единственный максимальный левый  идеал $I_l$.
\item Существует собственный идеал $I\nor S$ такой, что любой элемент из $S\setminus I$ обратим.
\end{list}
Кроме того, если выполнены условия $(i)$--$(iii)$, то $I_r=I_l=I=\J(S)$.
\end{pre}

\textupl{loc eq prf}{Доказать предложение \ref{loc eq}.}

\textuprn{Доказать, что единственными идемпотентами локального кольца являются $0$ и $1$.}



\textuprn{Пусть $p$ --- простое число, и $\pi$ --- множество всех простых чисел, отличных от $p$. Показать, что
кольцо $\ZZ[\pi^{-1}]$ является локальным.\footnote{См. \ref{subrq}.}}

Для локальных колец из леммы Накаямы вытекает следующее утверждение.

\begin{cor}[лемма Накаямы для локальных колец]\label{cor nak lok}
Пусть $S$ --- локальное кольцо с максимальным идеалом $I$,
$M$~--- конечно порождённый правый $S$-модуль, $N$ --- подмодуль модуля $M$ такой, что
$M=N+MI$. Тогда $M=N$.
\end{cor}

\begin{proof} Как мы замечали, $S$ можно рассматривать как $R$-алгебру для подходящего
коммутативного кольца $R$ и тогда $M$ и $N$ будут модулями над $R$-алгеброй $S$. Кроме того,
радикал $\J(S)$, очевидно, равен $I$. Значит, $M=N$ по лемме Накаямы, поскольку фактормодуль
$M/N$ конечно порождённого модуля $M$ также конечно порождён.
\end{proof}

\section{Конечномерные алгебры над полем}

Начиная с этого раздела, если не оговорено особо, мы будем предполагать, что $R=F$~--- поле, все
рассматриваемые $F$-модули (в том числе и $F$-алгебра $A$) конечномерны как векторные пространства над $F$, а алгебра $A$ ненулевая. В
частности, любой $A$-модуль обладает композиционным рядом, и, как мы уже отмечали, в этом случае можем
предполагать, что $F\le \Z(A)$, поскольку отображение $\a\mapsto \a1_A$ будет мономорфизмом колец.

\begin{pre} \label{fin irr num} Пусть $A$ --- $F$-алгебра.
Существует лишь конечное число попарно неизоморфных неприводимых $A$-модулей, т.\,е.
$|\M(A)|<\infty$.
\end{pre}

\textupr{Доказать предложение \ref{fin irr num}. \uk{Воспользоваться предложением \ref{ann el} и теоремой
Жордана-Гёльде\-ра
\ref{thm zh-g}.}}

\begin{pre} \label{rad nlp} Пусть $A$ --- $F$-алгебра. Тогда $\J(A)$ --- нильпотентный идеал.
\end{pre}

\textupl{rad nlp prf}{Доказать предложение \ref{rad nlp}.}

\textzam{Отметим, что радикал алгебры над произвольным кольцом, а также бесконечномерной алгебры над полем
может не быть нильпотентным.}

Из предложений \ref{rad max} и \ref{rad nlp} следует, что среди нильпотентных правых идеалов
$F$-алгебры $A$ существует единственный максимальный идеал, и этот идеал совпадает с радикалом
$\J(A)$. Также из предложения \ref{rad nlp} следует, что радикал $F$-алгебры состоит из
нильпотентных элементов. Значит, по предложению \ref{com nil} справедливо следующее
утверждение, дающее характеризацию радикала $F$-алгебр в коммутативном случае.

\begin{cor} \label{cor com rad}
Пусть $A$ --- коммутативная $F$-алгебра.
Тогда радикал $\J(A)$ совпадает с множеством нильпотентных элементов.
\end{cor}

\begin{pre} \label{rad cent}
Пусть $A$ --- $F$-алгебра.
Тогда для любой подалгебры $B\le \Z(A)$ имеет место равенство $\J(B)=\J(A)\cap B$.
\end{pre}

\textupl{rad cent prf}{Доказать предложение \ref{rad cent}.}

\textuprn{Найти радикалы следующих алгебр.
\begin{list}{{\rm(}{\it\roman{enumi}\/}{\rm)}}
{\usecounter{enumi}\setlength{\parsep}{2pt}\setlength{\topsep}{5pt}\setlength{\labelwidth}{23pt}}
\item $\MM_n(F)$, где $F$ --- поле.
\item $\ZZ_n$ как $\ZZ$-алгебра.
\item $\FF_p\ZZ_p$ для простого $p$.
\item $\FF_2S_3$.
\item $\FF_3S_3$.
\end{list}

\uk{$(iii)$~Рассмотреть идеал, порождённый элементом $1-a$, где $a$ --- порождающий группы $\ZZ_p$.
\quad $(iv)$~Сначала рассмотреть двусторонние идеалы $I_1$, $I_2$, порождённые, соответственно,
элементами $e_1=1+(123)+(132)$ и $e_2=(123)+(132)$. Воспользоваться предложениями \ref{inn dir sum} и
\ref{dir sum}. Рассмотреть дополнительно правые идеалы $J_1$, $J_2$, $J_3$, порождённые элементами
$v_1=e_2(1+(23))$, $v_2=e_2(1+(13))$, $v_3=e_2(1+(12))$.\quad
$(v)$~Рассмотреть правые идеалы $I_1$, $I_2$, порождённые, соответственно, элементами
$e_1=1+(12)$ и $e_2=1-(12)$. Показать, что модуль $I_i$ содержит единственный максимальный подмодуль
$J_i$, где $i=1,2$. Сравнить радикал с суммой $J_1+J_2$.}}



\section{Вполне приводимые модули}

Напомним, что $F$ всюду обозначает некоторое поле и $A$ --- конечномерную алгебру над $F$. Все
рассматриваемые $A$-модули считаются конечномерными пространствами над $F$.

Будем называть $A$-модуль $V$ \glsadd{iAModCompRed}\mem{вполне приводимым} или
\glsadd{iAModSemSmp}\mem{полупростым}, если он является прямой суммой неприводимых
$A$-подмодулей.





\textzam{Подчеркнём, что, в частности,  нулевой и неприводимый модуль вполне приводимы.}

\textuprn{\label{esn cr}Показать, что естественный подстановочный $FS_n$-модуль вполне приводим тогда и только тогда, когда
характеристика поля $F$ не делит $n$. \uk{Воспользоваться \ref{fsn}.}}

Следующее утверждение позволяет дать эквивалентное определение полной приводимости $A$-модулей.

\begin{pre} \label{equ ssm} Если $A$ --- $F$-алгебра и $V$ --- $A$-модуль,
 то следующие условия эквивалентны.
\begin{list}{{\rm(}{\it\roman{enumi}\/}{\rm)}}
{\usecounter{enumi}\setlength{\parsep}{2pt}\setlength{\topsep}{5pt}\setlength{\labelwidth}{23pt}}
\item Модуль $V$ вполне приводим.
\item Для любого $A$-подмодуля $U\le V$ существует $A$-подмодуль $W\le V$ такой, что $V=U\oplus W$.
\item Модуль $V$ является суммой $($необязательно прямой$\,)$ неприводимых подмодулей.
\end{list}
\end{pre}

\begin{proof}
$(i)\Rightarrow(ii)$ Пусть $V=\bigoplus V_i$ и $U\le V$. В силу конечномерности можно выбрать подмодуль
$W\le V$, максимальный среди всех, имеющих нулевое пересечение с $U$. Тогда $W+U=W\oplus U$ и достаточно
показать, что эта сумма совпадает с $V$. Пусть, напротив, $W+U<V$. Тогда существует $V_j$, не
содержащийся в $W+U$. В силу неприводимости $V_j$ получаем $(U+W)\cap V_j =0$. Но тогда $W+V_j>W$ и
подмодуль $W+V_j$ тривиально пересекается с $U$ вопреки выбору $W$.

$(ii)\Rightarrow(iii)$ Пусть $W$ --- максимальный подмодуль в $V$, являющийся суммой
неприводимых подмодулей. Если $W<V$, то $V=W\oplus U$ для ненулевого подмодуля $U\le V$. Но
тогда, ввиду конечномерности, $U$ содержит неприводимый подмодуль $M$. Значит, $W+M>W$ и
подмодуль $W+M$ --- сумма неприводимых подмодулей. Это противоречит выбору $W$.

$(iii)\Rightarrow(i)$ Пусть $V=\sum V_i$, где подмодули $V_i$ неприводимы. Пусть $W\le V$ ---
максимальный подмодуль, представимый как прямая сумма некоторых из подмодулей $V_i$. Если $W\ne V$, то
существует $V_j$, не лежащий в~$W$. Но тогда $W\cap V_j=0$ в силу неприводимости $V_j$ и сумма $W+V_j$
прямая вопреки максимальности~$W$.
\end{proof}

Отметим, что при доказательстве импликации $(iii)\Rightarrow(i)$ мы получили следующее более
сильное утверждение.

\begin{cor} \label{cor sum dir} Пусть $V$ --- $A$-модуль и $V=\sum V_i$,
где $V_i$ --- неприводимые подмодули. Тогда $V$ --- прямая сумма некоторых из подмодулей
$V_i$.
\end{cor}

\textuprn{Доказать, что
\begin{list}{{\rm(}{\it\roman{enumi}\/}{\rm)}}
{\usecounter{enumi}\setlength{\parsep}{2pt}\setlength{\topsep}{5pt}\setlength{\labelwidth}{23pt}}
\item вполне приводимый $A$-модуль неразложим тогда и только тогда, когда он неприводим;
\item любой подмодуль вполне приводимого модуля изоморфен некоторому его фактор модулю, и наоборот.
\end{list}}

\textuprn{\label{perm ncr}Пусть группа $G$ действует транзитивно на множестве $X$. Доказать, соответствующий
подстановочный $FG$-модуль не является вполне приводимым в случае, когда характеристика поля $F$ делит $|X|$.}


Очевидно, что прямая сумма вполне приводимых модулей вполне приводима. Следующее утверждение
показывает, что полная приводимость наследуется также при переходе к подмодулям и
фактормодулям.

\begin{pre} \label{ssm nsl} Пусть $V$ --- вполне приводимый $A$-модуль.
 Если $U$ ---
$A$-подмодуль из $V$, то $A$-модули $U$ и $V/U$ вполне приводимы.
\end{pre}

\textupl{ssm nsl prf}{Доказать предложение \ref{ssm nsl}.}

\textuprn{Доказать, что пересечение всех максимальных подмодулей вполне приводимого $A$-модуля равно нулю.}


%

Пусть $V$ --- вполне приводимый $A$-модуль и $M$ --- неприводимый $A$-модуль. \glsadd{iHmgCompAMod}\mem{$M$-однородной компонентой} модуля $V$
(или \glsadd{iCompWedAMod}\mem{компонентой Веддерберна}) называется сумма всех подмодулей из $V$, изоморфных $M$. Будем обозначать
\footnote{Если $V$ не содержит подмодулей, изоморфных $M$, то $\W_M(V)=0$.}
её через \gls{WMlVr}. Основные необходимые нам свойства компонент Веддерберна вполне приводимых модулей собраны в следующем утверждении.

\begin{pre} \label{ss wed}  Пусть $V$ --- вполне приводимый $A$-модуль.
Тогда
\begin{list}{{\rm(}{\it\roman{enumi}\/}{\rm)}}
{\usecounter{enumi}\setlength{\parsep}{2pt}\setlength{\topsep}{5pt}\setlength{\labelwidth}{23pt}}
\item справедливо разложение
$$
V=\bigoplus_{M\in \M(A)}\W_M(V).
$$
\end{list}

Зафиксируем произвольный неприводимый $A$-модуль $M$. Тогда
\begin{list}{{\rm(}{\it\roman{enumi}\/}{\rm)}}
{\usecounter{enumi}\setlength{\parsep}{2pt}\setlength{\topsep}{5pt}\setlength{\labelwidth}{23pt}}
\addtocounter{enumi}{1}
\item $\W_M(V)$ является $\End_A(V)$-подмодулем модуля $V$.
\end{list}

Зафиксируем также некоторое разложение $V=\bigoplus V_i$, где $V_i$ --- неприводимые $A$-модули. Тогда
\begin{list}{{\rm(}{\it\roman{enumi}\/}{\rm)}}
{\usecounter{enumi}\setlength{\parsep}{2pt}\setlength{\topsep}{5pt}\setlength{\labelwidth}{23pt}}
\addtocounter{enumi}{2}
\item имеет место равенство
$$
\W_M(V)=\displaystyle\bigoplus_{V_i\cong M} V_i;
$$
\item число подмодулей $V_i$, изоморфных $M$, является инвариантом модуля $V$, т.\,е. не зависит от данного
разложения $V=\bigoplus V_i$.
\end{list}
\end{pre}


\textupl{ss wed prf}{Доказать предложение \ref{ss wed}.}

Для вполне приводимого  $A$-модуля $V$ и неприводимого $A$-модуля $M$ будем обозначать через
\gls{nMlVr} число слагаемых, изоморфных $M$, в некотором (эквивалентно, в любом) разложении $V$ в прямую сумму неприводимых подмодулей.
В силу \ref{ss wed}$(iv)$ число
$$
\n_M(V)=\frac{\dim_F(\W_M(V))}{\dim_F M} \myeqno\label{nmv}
$$
зависит только от $M$ и $V$. Инварианты $\n_M(V)$ для всех $M\in \M(A)$ определяют вполне приводимый модуль $V$ с точностью до изоморфизма. А именно,
имеет место изоморфизм $A$-модулей
$$
V \cong \bigoplus_{M\in \M(A)} \left(\, \bigoplus_{i=1}^{\n_M(V)}M\right). \myeqno\label{vcong}
$$

\section{Полупростые \texorpdfstring{$F$}{F}-алгебры}

\begin{opr} $F$-алгебра $A$ называется \glsadd{iFAlgSemSmp}\mem{полупростой},
если $\J(A)=0$.
\end{opr}

\begin{pre}\label{fact ssm} Пусть $A$ --- $F$-алгебра.
Тогда факторалгебра $A/\J(A)$ полупроста.
\end{pre}

\textupl{fact ssm prf}{Доказать предложение \ref{fact ssm}.}

Полупростота $F$-алгебры и её модулей тесно связаны.

\begin{pre} \label{ssm crit}  Пусть $A$ --- $F$-алгебра. Тогда следующие условия
эквивалентны.
\begin{list}{{\rm(}{\it\roman{enumi}\/}{\rm)}}
{\usecounter{enumi}\setlength{\parsep}{2pt}\setlength{\topsep}{5pt}\setlength{\labelwidth}{23pt}}
\item Алгебра $A$ полупроста.
\item Регулярный модуль $A^\circ$ вполне приводим.
\item Каждый $A$-модуль вполне приводим.
\end{list}
\end{pre}

\begin{proof} $(i)\Rightarrow(ii)$ Достаточно показать, что существует конечный набор
максимальных правых идеалов $M_1,\ld,M_n$ алгебры $A$ такой, что $\cap_i M_i=0$. В самом деле, если это верно,
то $A^\circ$ можно изоморфно вложить во вполне приводимый модуль $A^\circ/M_1\oplus\ld\oplus A^\circ/M_n$ с
помощью отображения $a\mapsto(a+M_1,\ld,a+M_n)$. Тогда по предложению \ref{ssm nsl} модуль $A^\circ$ будет
вполне приводимым.


Пусть $N$ --- подмодуль из $A^\circ$, являющийся пересечением конечного числа максимальных
правых идеалов алгебры $A$ и имеющий минимальную размерность среди всевозможных таких
пересечений. Достаточно показать, что $N=0$. Пусть, напротив, $N\ne 0$. По условию $\J(A)=0$.
С другой стороны $\J(A)$ совпадает с пересечением всех максимальных правых идеалов из $A$ по
предложению \ref{rad int max}. Значит, существует максимальный правый идеал $M$ такой, что
$M\cap N<N$, вопреки выбору $N$.

$(ii)\Rightarrow(i)$ Поскольку модуль $A^\circ$ вполне приводим, то он равен сумме
неприводимых  подмодулей --- минимальных правых идеалов алгебры $A$. Поскольку $\J(A)$
аннулирует каждый простой модуль, то $0=A^\circ\J(A)=\J(A)$. Значит, алгебра $A$ полупроста.

$(ii)\Rightarrow(iii)$ Пусть $V$ --- $A$-модуль. В силу полной приводимости модуля $A^\circ$
мы можем записать $A^\circ =\sum_i I_i$, где $I_i$ --- некоторые минимальные правые идеалы
алгебры $A$. Выберем некоторый $F$-базис $v_1,\ld,v_n$ модуля $V$. Тогда
$$
V=VA=\sum_j v_jA= \sum_j\sum_i v_jI_i.
$$
Отображение $I_i\to v_jI_i$, действующее по правилу $a\mapsto v_ja$, $a\in I_i$,
является эпиморфизмом $A$-модулей. Поскольку модуль $I_i$ неприводим, то $v_jI_i\cong I_i$ или $v_jI_i\cong 0$.
Значит, $V$ --- сумма неприводимых модулей.

$(iii)\Rightarrow(ii)$ Очевидно.
\end{proof}

Поскольку регулярный модуль факторалгебры изоморфен фактормодулю регулярного модуля исходной
алгебры, из предложений \ref{ssm nsl} и \ref{ssm crit}  вытекает следующее утверждение.

\begin{cor} \label{cor fact ss}
Факторалгебра полупростой $F$-алгебры
полупроста.
\end{cor}

\textzam{Отметим, что, в отличие от факторалгебры, подалгебра полупростой алгебры может не быть полупростой.
В качестве примера можно рассмотреть подалгебру верхнетреугольных матриц в алгебре $\MM_n(F)$, $n\ge 2$. Эта
подалгебра содержит ненулевой нильпотентный идеал, состоящий из матриц с нулевой диагональю.}

\begin{cor} \label{cor min idl} Если $F$-алгебра $A$ полупроста,
то любой неприводимый $A$-модуль изоморфен минимальному правому идеалу из $A$.
\end{cor}

\textupl{cor min idl prf}{Доказать следствие \ref{cor min idl}.}

\textuprn{Пусть $V$ --- произвольный $A$-модуль.
\begin{list}{{\rm(}{\it\roman{enumi}\/}{\rm)}}
{\usecounter{enumi}\setlength{\parsep}{2pt}\setlength{\topsep}{5pt}\setlength{\labelwidth}{23pt}}
\item если $V\J(A)=0$, то $V$ вполне приводим,
\item фактормодуль $V/(V\J(A))$ вполне приводим,
\item $V\J(A)$ совпадает с пересечением максимальных подмодулей модуля $V$.
\end{list}
}



Исследуем теперь вопрос о полупростоте в важном для нас случае групповых алгебр.

\begin{pre} \label{mod ns}
Пусть $G$ --- конечная группа,
$p$ --- простой делитель порядка $|G|$, и $F$ --- поле характеристики $p$. Тогда групповая алгебра $FG$ не полупроста.
\end{pre}

\textupl{mod ns prf}{Доказать предложение \ref{mod ns}.} \uk{Показать, что $\sum_{g\in G}g\in \J(FG)$.}

Справедливо также утверждение, обратное к предложению \ref{mod ns}. Для его доказательства нам понадобится
следующее вспомогательное утверждение, являющееся несложным упражнением по линейной алгебре.

\begin{lem} \label{lem ff2} Пусть $\vf$ --- линейное преобразование
конечномерного векторного пространства $V$ над полем $F$. Показать, что $V=\Im\vf\oplus\Ker\vf$
тогда и только тогда, когда $V\vf^2=V\vf$.
\end{lem}

\textupl{lem ff2 prf}{Доказать лемму \ref{lem ff2}.}

\begin{thm}[Машке] \label{thm mash}\glsadd{iThmMasch} Пусть $G$ --- конечная группа
и $F$ --- поле, характеристика которого не делит\footnote{В частности, характеристика поля $F$ может быть равна нулю.}
порядок $|G|$. Тогда групповая алгебра $FG$ полупроста.
\end{thm}
\begin{proof} По предложению \ref{ssm crit} достаточно показать, что произвольный
$FG$-модуль $V$ вполне приводим. Пусть $U$ --- подмодуль из $V$. Найдём подмодуль $W\le V$ такой, что
$V=U\oplus W$. Векторное пространство $V$ можно представить в виде прямой суммы $U\oplus U_0$ для некоторого
подпространства $U_0\le V$. Пусть $\pi: V\to V$ --- линейное отображение, являющееся проекцией $V$ на $U$
параллельно $U_0$. Рассмотрим отображение $\ti{\pi}: V\to V$, которое строится по $\pi$ следующим образом:
$$
v\ti{\pi}=\frac{1}{|G|}\sum_{g\in G}vg\pi g^{-1}.
$$
Другими словами $\ti{\pi}$ является <<усреднением>> всех отображений, сопряжённых с $\pi$ элементами
группы $G$, причём деление на $|G|$ возможно ввиду ограничения на характеристику поля $F$. Заметим, что
$\ti{\pi}$ является $F$-линейным и $\Im \ti{\pi}\se U$, поскольку $U$ --- $FG$-подмодуль. Кроме
того, для любого $h\in G$ выполнено
$$
vh\ti{\pi}=\frac{1}{|G|}\sum_{g\in G}(vhg\pi g^{-1}h^{-1})h=v\ti{\pi}h,
$$
откуда следует, что $\ti\pi\in \End_{FG}(V)$. Более того, для любого $u\in U$ имеем
$u\ti\pi=u$. Поэтому $U=\Im \ti\pi$ и $\ti\pi^2=\ti\pi$. По лемме
\ref{lem ff2} получаем $V=U\oplus W$, где $W=\Ker \ti\pi$ --- $FG$-подмодуль модуля $V$.
\end{proof}

Объединяя утверждения \ref{mod ns} и \ref{thm mash}, получаем критерий полупростоты групповой алгебры.

\begin{cor} \label{cor cr pol} Пусть $G$ --- конечная группа и
$F$ --- поле. Групповая алгебра $FG$ полупроста тогда и только тогда, когда характеристика поля $F$ не делит
порядок группы $G$.
\end{cor}

\textext{Мы отмечали в \ref{ssm crit}, что полупростота алгебры $FG$ и её регулярного модуля эквивалентны,
а в \ref{reg perm} --- что регулярный $FG$-модуль является частным случаем подстановочного модуля. Сравнивая \ref{cor cr pol}
с критерием \ref{esn cr} полупростоты естественного подстановочного $FS_n$-модуля, а также с учётом \ref{perm ncr},
можно сформулировать следующий вопрос. Пусть группа $G$ транзитивно действует на множестве $X$.
Эквивалентна ли полная приводимость соответствующего подстановочного $FG$-модуля $V$ тому, что характеристика поля $F$ не
делит $|X|$? Ответ на этот вопрос, вообще говоря, отрицательный. Например, спорадическая группа Матьё
$M_{23}$ транзитивно действует на множестве из $23$-х элементов, но подстановочный $\FF_2M_{23}$-модуль для этого действия
не является вполне приводимым.}

\textuprn{Доказать следующее обобщение теоремы Машке. Пусть $G$ --- конечная группа, $H\le G$,
$F$ --- поле характеристики, не делящей индекс $|G:H|$. Если $V$ --- $FG$ модуль такой, что
$FH$-модуль $V_H$ вполне приводим, то $FG$-модуль $V$ также вполне приводим.
\uk{Следовать рассуждению из доказательства теоремы Машке.}}
%
%

\textuprn{Пусть $V$ --- конечная абелева группа и $G$ --- такая подгруппа группы
$\Aut(V)$, что $(|V|,|G|)=1$. Доказать следующие утверждения.

\begin{list}{{\rm(}{\it\roman{enumi}\/}{\rm)}}
{\usecounter{enumi}\setlength{\parsep}{2pt}\setlength{\topsep}{5pt}\setlength{\labelwidth}{23pt}}
\item \mbox{}\big(Аналог теоремы Машке для абелевых групп.\big) Если $V=U\oplus U_0$, где $U$ --- $G$-инвариантная подгруппа, то существует
$G$-инвариантная подгруппа $W\le V$ такая, что $V=U\oplus W$.
\item $V=\C_V(G)\oplus[V,G]$, где
\begin{align*}
  \C_V(G) & =\{v\in V\mid v^g=v \ \ \text{для всех}\ \ g\in G\}, \\
  [V,G] & = \langle v^g-v \ \mid \ v\in V, g\in G\rangle.
\end{align*}
Кроме того, проекция $\pi: V\to \C_V(G)$ находится\footnote{Ввиду того, что порядки групп $V$ и $G$ взаимно просты, уравнение $|G|x=w$ имеет
в группе $V$ единственное решение для любого $w\in V$.}
из соотношения
$$
|G|(v\pi)=\sum_{g\in G} v^g,
$$
\end{list}
где $v\in V$. \uk{Доказать аналог леммы \ref{lem ff2} для случая, когда $V$~--- конечная абелева группа и $\vf\in \End(V)$. Далее
следовать рассуждению из доказательства теоремы Машке.}}

Пусть $V$ --- $A$-модуль. Легко видеть, что аннулятор $\Ann(V)$ совпадает с ядром естественного
гомоморфизма $F$-алгебр $A\to \End_F(V)$, $a\mapsto a_V$. Обозначим через \gls{AV} образ этого гомоморфизма.
Таким образом, $A_V\cong A/\Ann(V)$.

\begin{pre} \label{mlt vp} Для вполне
приводимого $A$-модуля $V$ алгебра $A_V$ полупроста.
\end{pre}
\textupl{mlt vp prf}{Доказать предложение \ref{mlt vp}.}

Для $F$-алгебр полупростота тесно связано с другим важным понятием --- понятием простой алгебры.

\begin{opr}
Ненулевая $F$-алгебра называется \glsadd{iFAlgSmp}\mem{простой}, если она не
имеет собственных ненулевых двусторонних идеалов.
\end{opr}
Очевидно, что простая $F$-алгебра и даже прямая сумма простых $F$-алгебр будут полупростыми. Оказывается, верно и обратное:
всякая полупростая алгебра является прямой суммой простых, как показывает следующая теорема.

\begin{thm}[Веддерберна]\label{thm wedd}\glsadd{iThmWedd} Пусть $A$ --- полупростая $F$-алгебра.
\begin{list}{{\rm(}{\it\roman{enumi}\/}{\rm)}}
{\usecounter{enumi}\setlength{\parsep}{2pt}\setlength{\topsep}{5pt}\setlength{\labelwidth}{23pt}}
\item Справедливо разложение
$$
A=\bigoplus_{M\in \M(A)}\W_M(A^\circ),
$$
и его компоненты $\W_M(A^\circ)$ являются минимальными двусторонними идеалами в~$A$.
\end{list}

Зафиксируем  неприводимый $A$-модуль $M$.
\begin{list}{{\rm(}{\it\roman{enumi}\/}{\rm)}}
{\usecounter{enumi}\setlength{\parsep}{2pt}\setlength{\topsep}{5pt}\setlength{\labelwidth}{23pt}}
\addtocounter{enumi}{1}
\item Пусть $U$ --- неприводимый $A$-модуль. Если  $M\ncong U$,  то $\W_M(A^\circ)\le \Ann(U)$, а если
$M\cong U$, то
$$\W_M(A^\circ)\cap \Ann(U)=0.$$
В частности, имеет место разложение
\end{list}
\vspace{-2.5\topsep}
$$
A=\W_M(A^\circ)\oplus \Ann(M) \myeqno\label{wedd_dec}
$$

\mbox{}\quad алгебры $A$ в прямую сумму идеалов.
\begin{list}{{\rm(}{\it\roman{enumi}\/}{\rm)}}
{\usecounter{enumi}\addtocounter{enumi}{2}\setlength{\parsep}{2pt}\setlength{\topsep}{5pt}\setlength{\labelwidth}{23pt}}
\item $\W_M(A^\circ)$ является простой $F$-алгеброй, и отображение $a\mapsto a_M$ осуществляет
изоморфизм $F$-алгебр $\W_M(A^\circ)\cong A_M\le \End_F(M)$.
\end{list}
\end{thm}

\begin{proof} Поскольку, согласно \ref{ss wed}$(i)$, регулярный $A$-модуль допускает разложение
$$
A^\circ=\bigoplus_{M\in \M(A)}\W_M(A^\circ), \myeqno\label{reg dec}
$$
мы докажем $(i)$, если установим, что для фиксированного
неприводимого модуля $M$ компонента $\W_M(A^\circ)$ является минимальным двусторонним идеалом в $A$.
Заметим, что $\W_M(A^\circ)\ne 0$, т.\,к. по следствию  \ref{cor min idl}
модуль $A^\circ$ содержит подмодуль, изоморфный $M$. Сначала покажем, что $\W_M(A^\circ)$ --- идеал.
Его минимальность будет установлена чуть ниже. Как подмодуль регулярного модуля, $\W_M(A^\circ)$
является правым идеалом. Легко видеть также, что для любого $a\in A$
отображение $a_l:x\mapsto ax$ принадлежит алгебре $\End_A(A^\circ)$. По пункту $(ii)$ предложения
\ref{ss wed} выполнено включение $\W_M(A^\circ)a_l\se \W_M(A^\circ)$, т.\,е. правый идеал $\W_M(A^\circ)$
выдерживает умножения слева на элементы из $A$ и, значит, является двусторонним.

Если $U\ncong M$, то $\W_U(A^\circ)\cap\W_M(A^\circ)=0$ в силу \ref{ss wed}$(i)$. Поскольку $\W_U(A^\circ)$
и $\W_M(A^\circ)$ --- идеалы, имеет место равенство $\W_U(A^\circ)\W_M(A^\circ)=0$. По определению $\W_U(A^\circ)$ --- сумма
подмодулей из $A^\circ$, изоморфных $U$. Пусть $U_0\le \W_U(A^\circ)$ --- один из них. Тогда $U_0\W_M(A^\circ)=0$, и значит,
$\W_M(A^\circ)\se\Ann(U_0)=\Ann(U)$ ввиду изоморфизма $U_0\cong U$.

Пусть $U\cong M$. Предположим, что $I=\W_M(A^\circ)\cap \Ann(U)\ne 0$. Тогда любой
неприводимый подмодуль из $I$ является подмодулем $\W_M(A^\circ)$ и, значит, изоморфен $U$.
Пусть $U_0$ --- один из таких подмодулей. Тогда $U_0\le \Ann(U)=\Ann(U_0)$, откуда следует,
что  $U_0^2=0$, т.\,е. идеал $U_0$ нильпотентен. По предложению \ref{rad max} получаем $U_0\le
\J(A)$ вопреки полупростоте алгебры $A$. Отсюда следует, что сумма $\W_M(A^\circ)+\Ann(M)$ прямая.
Из доказанного следует, что эта сумма содержит все компоненты Веддерберна модуля $A^\circ$ и поэтому
совпадает с $A$ в силу разложения \ref{reg dec}. Тем самым $(ii)$ доказано.

Ввиду $(ii)$ всякий элемент $a\in A$ представим в виде $a=b+c$, где $b\in
\W_M(A^\circ)$ и $c\in \Ann(M)$. Поэтому $a_M=b_M$ и, значит, гомоморфизм $a\mapsto a_M$
сюръективно отображает $\W_M(A^\circ)$ на $A_M$. Если $a\in \W_M(A^\circ)$ и $a_M=0$, то в силу $(ii)$
получаем $a\in \W_M(A^\circ)\cap \Ann(M)=0$. Таким образом, отображение $a\mapsto a_M$ осуществляет
изоморфизм $F$-алгебр $\W_M(A^\circ)$ и $A_M$.

Теперь покажем, что идеал $\W_M(A^\circ)$ алгебры $A$ минимален. Пусть $I\nor A$ и $I<\W_M(A^\circ)$. 
Так как $\W_M(A^\circ)$~--- сумма $A$-подмодулей из $A^\circ$, изоморфных $M$, найдётся 
подмодуль $M_0\le\W_M(A^\circ)$ такой, что $M_0\cong M$ и $M_0\not\subseteq I$. Тогда $M_0\cap I=0$ 
ввиду неприводимости $M_0$. Значит, $M_0I\le M_0\cap I=0$, 
т.\,е. $I$ аннулирует $M_0\cong M$, что влечёт $I\leqslant\Ann(M)$. Но с другой стороны $I<\W_M(A^\circ)$.
Поэтому $I=0$ в силу \ref{wedd_dec}.

Простота $F$-алгебры $\W_M(A^\circ)$ теперь следует из её минимальности как идеала в $A$ ввиду того, что 
всякий её идеал будет также идеалом в $A$, как
вытекает из разложения \ref{reg dec} и доказанного выше равенства $\W_{M_1}(A^\circ)\W_{M_2}(A^\circ)=0$ для
неизоморфных неприводимых $A$-модулей $M_1$ и $M_2$ Отсюда получаем $(iii)$. Тем самым теорема 
полностью доказана.\end{proof}

По существу, теорема Веддерберна сводит изучение полупростых $F$-алгебр к простым.

\begin{pre} \label{mlt npr} Пусть $A$ --- $F$-алгебра.
\begin{list}{{\rm(}{\it\roman{enumi}\/}{\rm)}}
{\usecounter{enumi}\setlength{\parsep}{2pt}\setlength{\topsep}{5pt}\setlength{\labelwidth}{23pt}}
\item Если $A$ простая, то $|\M(A)|=1$.
\item Если $V$ --- неприводимый $A$-модуль, то алгебра $A_V$ простая.
\end{list}
\end{pre}

\textupl{mlt npr prf}{Доказать предложение \ref{mlt npr}.}

\textuprn{ \glsadd{iGrQuat}\mem{Группой кватернионов} называется группа
$$
\gls{Q8}=\la x,y\bigm| x^4=1,\ \  y^2=x^2, \ \ x^y=x^{-1} \ra.
$$
Сколько различных двусторонних идеалов в алгебре $\CC Q_8$?}

\mysubsection \label{cid}
Пусть $A$ --- полупростая $F$-алгебра. Из разложения \ref{thm wedd}$(i)$ алгебры $A$ в прямую сумму идеалов $\W_M(A^\circ)$,
а также из \ref{inn dir sum} следует,
что единица алгебры $A$ однозначно представима в виде суммы центральных идемпотентов
$$
1=\sum_{M\in \M(A)}e_{\mbox{}_M},
$$
где $\gls{esM}\in\W_M(A^\circ)$. Элемент $e_{\mbox{}_M}$ будем называть \glsadd{iIdmpCenFAlg}\mem{центральным
идемпотентом алгебры $A$, соответствующим неприводимому $A$-модулю} $M$. Этот идемпотент является единицей алгебры
$\W_M(A^\circ)$ и справедливо равенство $\W_M(A^\circ)=e_{\mbox{}_M}A$. Идемпотенты $e_{\mbox{}_M}$, где  $M$
пробегает $\M(A)$, попарно ортогональны и в силу простоты алгебр $\W_M(A^\circ)$ исчерпывают весь набор
примитивных идемпотентов алгебры $\Z(A)$.

\textzam{Отметим, что $e_{\mbox{}_M}$, вообще говоря, не является примитивным идемпотентом алгебры $A$.}

\begin{thm}[о двойном централизаторе] \label{thm dv cent}\glsadd{iThmDblCntr} Пусть $A$ --- простая $F$-алгебра
и $I\ne 0$--- правый идеал из $A$. Обозначим $B=\End_A(I)$. Тогда для любого $a\in A$
преобразование $a_{\mbox{}_I}$ является эндоморфизмом $B$-модуля $I$. Кроме того, $A_{\mbox{}_I}=\End_B(I)$
и отображение $a\mapsto a_{\mbox{}_I}$
является изоморфизмом $F$-алгебр $A\to \End_B(I)$.
\end{thm}

\begin{proof} Мы уже отмечали  (см. упражнение \ref{end mod}), что $B$ является
$F$-алгеброй, а $I$ --- $B$-модулем. Поэтому $\End_B(I)$ снова является $F$-алгеброй. Для
произвольного элемента $a\in A$ выполнены равенства
$$
(x\vf)a_{\mbox{}_I}=(x\vf)a=(xa)\vf=(xa_{\mbox{}_I})\vf
$$
при любых $x\in I$, $\vf\in B$. Поэтому $a_{\mbox{}_I}\in \End_B(I)$.

Покажем, что отображение $a\mapsto a_{\mbox{}_I}$ является изоморфизмом алгебр $A\to \End_B(I)$. Перед формулировкой
теоремы Веддерберна \ref{thm wedd} мы отметили, что оно является гомоморфизмом алгебр с ядром
$\Ann(I)$. Из простоты алгебры $A$ следует, что $\Ann(I)=0$ и, значит, это отображение
инъективно.

Прежде, чем доказывать сюръективность, сделаем одно замечание. Для произвольного $u\in I$ рассмотрим
отображение $u_l:I\to I$, действующее по правилу $x\mapsto ux$ для
любого $x\in I$. Поскольку для всех $x\in I$ и $a\in A$ выполнено
$$(xu_l)a=(ux)a=u(xa)=(xa)u_l,$$
имеем $u_l\in B$. В частности, действие $u_l$ на $I$ перестановочно с действием произвольного
элемента $\vf\in \End_B(I)$.


Итак, нам осталось показать, что $A_{\mbox{}_I}=\End_B(I)$.
Поскольку $1_{\mbox{}_I}\in A_{\mbox{}_I}$ --- единица в $\End_B(I)$, достаточно проверить, что $A_{\mbox{}_I}$ --- правый идеал в $\End_B(I)$.
Заметим, что $A=AI$, так как $AI$ --- ненулевой двусторонний идеал в $A$ и алгебра $A$ простая.
Поскольку элементы из $AI$ --- это суммы произведений вида $ax$, где $a\in A$ и $x\in I$,
то достаточно показать, что $(ax)_{\mbox{}_I}\vf\in (AI)_{\mbox{}_I}=A_{\mbox{}_I}$ для всякого $\vf\in \End_B(I)$.
Учитывая предыдущее замечание о перестановочности действия $\vf$ с операторами левого умножения на элементы из
$I$, получаем следующую цепочку равенств для произвольного $y\in I$
$$
y((ax)_{\mbox{}_I}\vf)=(y(ax)_{\mbox{}_I})\vf=(yax)\vf=(x(ya)_l)\vf=(x\vf)(ya)_l=ya(x\vf)=y(a(x\vf))_{\mbox{}_I},
$$
поскольку $ya\in I$. В силу произвольности $y$, имеет место требуемое включение
$(ax)_{\mbox{}_I}\vf=(a(x\vf))_{\mbox{}_I}\in (AI)_{\mbox{}_I}$. \end{proof}

\textext{Объясним, почему предыдущее утверждение носит название <<теорема о двойном централизаторе>>. Если $S$ ---
некоторая алгебра и $X\se S$ --- произвольное подмножество, то
\glsadd{iCntrlAlg}\mem{централизатором} $X$ в $S$ называется множество
$$\gls{CSX}=\{s\in S\bigm|xs=sx \ \  \mbox{для всех}\ x\in X\}$$ которое, как легко видеть, является
подалгеброй в $S$. В частности, $\Z(S)=\C_S(S)$. Теперь в обозначениях теоремы \ref{thm dv cent} мы можем записать
$$B=\End_A(I)=\C_{\End_F(I)}(A_{\mbox{}_I})\qquad\text{и}\qquad \End_B(I)=\C_{\End_F(I)}(B).$$
Поэтому основное утверждение теоремы состоит в том, что для простой
$F$-алгебры $A$ и её ненулевого правого идеала $I$ выполнено равенство
$$
A_{\mbox{}_I}=\C_{\End_F(I)}(\C_{\End_F(I)}(A_{\mbox{}_I})).
$$}

Важным следствием из теоремы \ref{thm dv cent} является описание строения простых алгебр над
алгебраически замкнутым полем.

\begin{cor} \label{cor sim azp}
Пусть $F$ --- алгебраически замкнутое поле и $A$ --- простая $F$-алгебра.
Если $V$ --- неприводимый $A$-модуль, то $A\cong \End_F(V)$. В частности,
$\dim_F A=n^2$, где  $n=\dim_F V$, и имеет место изоморфизм $A$-модулей
$$A^\circ\cong \underbrace{V\oplus\ld\oplus V}_n$$
\end{cor}
\begin{proof} Поскольку модуль $V$ изоморфен минимальному правому идеалу алгебры $A$, в силу следствия \ref{cor min idl},
можно считать, что $V\le A^\circ$. Положив $B=\End_A(V)$, по следствию
\ref{cor sh} из леммы Шура заключаем, что $B \cong F$. По теореме \ref{thm dv cent} имеем
$$A\cong\End_B(V)=\End_F(V)\cong \MM_n(F).
$$
Таким образом, $\dim_F A=n^2$. Из простоты алгебры $A$ и теоремы Веддерберна \ref{thm wedd}
следует, что $A^\circ=\W_V(A^\circ)$, а равенство \ref{nmv} влечёт, что $A^\circ$ --- прямая сумма
$$
\n_V(A^\circ)=\frac{\dim_F(\W_V(A^\circ))}{\dim_F V}=\frac{\dim_F(A^\circ)}{\dim_F V}=\frac{n^2}{n}=n
$$
подмодулей, изоморфных $V$.
\end{proof}

\begin{cor} \label{cor z ss}
Пусть $F$ --- алгебраически замкнутое поле и $A$ --- полупростая $F$-алгебра. Тогда
справедливы следующие утверждения.
\begin{list}{{\rm(}{\it\roman{enumi}\/}{\rm)}}
{\usecounter{enumi}\setlength{\parsep}{2pt}\setlength{\topsep}{5pt}\setlength{\labelwidth}{23pt}}
\item Для любого $M\in \M(A)$ имеем $\dim_F \W_M(A^\circ)=(\dim_F M)^2$, а также  $\n_M(A^\circ)=\dim_F M$.

\item Имеет место изоморфизм $A$-модулей
$$A^\circ \cong \bigoplus_{M\in \M(A)} \left(\, \bigoplus_{i=1}^{\dim_F M}M\right).$$

\item $\dim_F A=\sum_{M\in \M(A)} (\dim_F M)^2.$

\item Идемпотенты $e_{\mbox{}_M}$, где  $M$ пробегает $\M(A)$, образуют $F$-базис центра $\Z(A)$. В частности,
$$\dim_F Z(A)=|\M(A)|.$$
\end{list}
\end{cor}

\textupl{cor z ss prf}{Доказать следствие \ref{cor z ss}}

Для групповых алгебр, учитывая предложения \ref{kl sum} и \ref{thm mash}, следствие \ref{cor z ss} можно
переформулировать так:

\begin{cor} \label{cor z fg} Пусть $G$ --- конечная группа и $F$ --- алгебраически замкнутое поле,
характеристика которого не делит $|G|$. Тогда
\begin{list}{{\rm(}{\it\roman{enumi}\/}{\rm)}}
{\usecounter{enumi}\setlength{\parsep}{2pt}\setlength{\topsep}{5pt}\setlength{\labelwidth}{23pt}}
\item Имеет место изоморфизм $FG$-модулей
$$FG^\circ \cong \bigoplus_{M\in \M(FG)} \left(\, \bigoplus_{i=1}^{\dim_F M}M\right);$$
\item $|G|=\sum_{M\in \M(FG)} (\dim_F M)^2;$
\item $|\M(FG)|=|\K(G)|$.
\end{list}
\end{cor}

\chapter{Представления и характеры алгебр и групп}

\section{Представления алгебр}

\begin{opr} \glsadd{iRepAlgMat}\mem{Представлением}
$F$-алгебры $A$ называется гомоморфизм $F$-алгебр
$$
\X : A\to \MM_n(F)
$$
для некоторого $n$.
Число $n$ называется  \glsadd{iDegRep}\mem{степенью
представления} $\X$ и обозначается \gls{degX}. Представления степени $1$ мы будем часто называть
\glsadd{iRepLin}\mem{линейными}.
\end{opr}

\textzam{Символы представлений мы будем писать слева от аргумента.}

Легко видеть, что если $\X$ --- представление степени $n$ алгебры $A$ и $P$ --- невырожденная матрица из $\MM_n(F)$, то
отображение $a\mapsto P^{-1}\X(a)P$ определяет новое представление степени $n$ алгебры $A$. Однако это новое
представление, как мы сейчас увидим, мало отличается от исходного.

Пусть $\X$ и $\Y$ --- представления алгебры $A$, имеющие одну и ту же степень.
Будем называть $\X$ и $\Y$  \glsadd{iRepsEqv}\mem{эквивалентными},
 и записывать $\X\cong\Y$ если существует невырожденная матрица $P$ такая, что
$\Y(a)=P^{-1}\X(a)P$ для всех $a\in A$.

Из представлений легко строить модули, а из модулей --- представления. Пусть $V=F^n$ ---
пространство строк и $\X$ --- представление степени $n$. Если положить $va=v\X(a)$ для всех
$v\in V$, $a\in A$, то таким образом на $V$ задаётся структура $A$-модуля. Такой $A$-модуль
будем называть \glsadd{iModCorRep}\mem{модулем, соответствующим представлению} $\X$.

Обратно, пусть $V$ --- $A$-модуль, имеющий размерность $n$ над $F$. Выберем $F$-базис $v_1,\ld, v_n$
пространства $V$ и любому $a\in A$ поставим в соответствие матрицу\footnote{Другими словами,
$(\a_{ij})$ является матрицей преобразования $a_{V}$ в выбранном базисе $v_1,\ld, v_n$} $(\a_{ij})$, элементы которой находятся из разложений
$$v_ia=\a_{i1}v_1+\ld +\a_{in}v_n, \qquad i=1,\ldots,n.$$
Отображение $a\mapsto (\a_{ij})$ является представлением степени $n$ алгебры $A$. Такое представление мы будем
называть \glsadd{iRepCorMod}\mem{представлением, соответствующим} (в базисе $v_1,\ld, v_n$) \mem{модулю}
$V$.\footnote{Нулевому модулю соответствует \glsadd{iRepNul}\mem{нулевое представление}, т.\,е. представление нулевой степени.}

\begin{pre} \label{mod rep} Модули, соответствующие эквивалентным представлениям,
изоморфны. Обратно, представления, соответствующие изоморфным модулям, эквивалентны. В
частности, представления, соответствующие одному и тому же модулю в различных базисах,
эквивалентны.
\end{pre}

\textupr{Доказать предложение \ref{mod rep}.}


\textuprn{Пусть $\X$ --- представление $F$-алгебры $A$. Показать, что
\begin{list}{{\rm(}{\it\roman{enumi}\/}{\rm)}}
{\usecounter{enumi}\setlength{\parsep}{2pt}\setlength{\topsep}{5pt}\setlength{\labelwidth}{23pt}}
\item если $\Y$ --- представление алгебры $A$, эквивалентное $\X$, то $\Ker\X=\Ker\Y$;
\item если $V$ --- $A$-модуль, соответствующий $\X$, то $\Ker\X=\Ann(V)$.
\end{list}
}

Пусть $V$ --- векторное пространство размерности $n$ над полем $F$. Изоморфизм $F$-алгебр $\End_F(V)\cong
\MM_n(F)$ позволяет использовать термин <<представление>> также и для любого гомоморфизма $\X:A\to \End_F(V)$.
В частности, если  $V$ --- $A$-модуль, то отображение $A\to \End_F(V)$, задаваемое правилом $a\mapsto a_V$,
может рассматриваться как представление алгебры $A$.

\subsection{Примеры представлений алгебр} \label{pr pr}
Пусть $A$ --- $F$-алгебра.
\begin{list}{{\rm(}{\it\roman{enumi}\/}{\rm)}}
{\usecounter{enumi}\setlength{\parsep}{2pt}\setlength{\topsep}{5pt}\setlength{\labelwidth}{23pt}}
\item \glsadd{iRepFAlgReg}\mem{Регулярным представлением алгебры}
$A$
будем называть представление, соответствующее регулярному модулю $A^\circ$ в некотором базисе.
\item Если $\X$ --- представление $A$  и $B$ --- подалгебра в $A$, то ограничение \gls{XsB} представления $\X$ на $B$ является
представлением алгебры  $B$. При этом ясно, что $\deg \X_B =\deg \X$.
\item Если $I\nor A$, $\vf:A\to A/I$ --- естественный гомоморфизм алгебр
и $\wt\X$ --- представление факторалгебры $A/I$, то композиция $\X=\ov\X\circ\vf$
является представлением алгебры $A$.
При этом $\deg\X=\deg\wt\X$ и $I\se \Ker \X$. Обратно, если $\X$ --- представление алгебры $A$,
и $I\se \Ker \X$, то существует\footnote{Ср. предложение \ref{hom r}$(i)$.} единственное
представление $\wt\X$ факторалгебры $A/I$ такое, что $\X=\wt\X\circ\vf$.
\end{list}
\textupr{Проверить сформулированные в \ref{pr pr} утверждения.}

Представление $F$-алгебры $A$ называется \glsadd{iRepFAlgIrr}\mem{неприводимым}
(\glsadd{iRepFAlgRed}\mem{приводимым}, \glsadd{iRepFAlgCrd}\mem{вполне приводимым}, \glsadd{iRepFAlgInd}\mem{неразложимым}),
если таким является соответствующий ему $A$-модуль.

\textuprn{\label{nep f}Пусть $I\nor A$, $\vf:A\to A/I$ --- естественный гомоморфизм алгебр и $\X$ и $\wt\X$ --- представления
алгебр  $A$ и $A/I$, соответственно, такие, как в \ref{pr pr}$(iii)$. Показать, что $\X$ неприводимо
(приводимо, вполне приводимо, неразложимо) тогда и только тогда, когда таковым является $\wt\X$.}

Из сделанных выше замечаний о соответствии представлений и модулей, а также из предложения \ref{fin irr
num} получаем

\begin{cor} \label{fin irr rep} У любой $F$-алгебры $A$
существует лишь конечное число попарно неэквивалентных неприводимых представлений.
\end{cor}

\begin{pre} \label{int krs}
Пусть
$A$ --- $F$-алгебра и $\X_1,\ld,\X_s$ --- полная система её попарно неэквивалентных неприводимых
представлений. Тогда имеют место следующие утверждения.
\begin{list}{{\rm(}{\it\roman{enumi}\/}{\rm)}}
{\usecounter{enumi}\setlength{\parsep}{2pt}\setlength{\topsep}{5pt}\setlength{\labelwidth}{23pt}}
\item Справедливо равенство $$\J(A)=\bigcap_{i=1}^s\Ker \X_i.$$
\item Пусть $\vf:A\to A/J(A)$ --- естественный гомоморфизм $F$-алгебр и для любого $i=1,\ldots,s$ через
$\wt\X_i$ обозначено однозначно определённое представление алгебры $A/J(A)$ такое, что $\wt\X_i\circ\vf=\X_i$.
Тогда $\wt\X_1,\ldots,\wt\X_s$ --- полная система попарно неэквивалентных неприводимых
представлений алгебры $A/J(A)$.
\end{list}
\end{pre}

\textupl{int krs prf}{Доказать предложение \ref{int krs}.}

Такие свойства представлений, как неприводимость, полная приводимость и т.\,д. можно выразить на матричном языке,
не переходя к соответствующим модулям.

Пусть $\X$ --- приводимое представление алгебры $A$ и $V$ --- соответствующий ему модуль. Тогда
$V$ содержит нетривиальный собственный подмодуль $U$. Выберем в $V$ базис так, чтобы последние его
векторы образовывали базис $U$. Тогда в этом <<согласованном>> c рядом
$$
V>U>0
$$
базисе представление, соответствующее модулю $V$, имеет блочный вид
$$
a\mapsto
\left(\ba{cc}
\X_1(a) & * \\
\O & \X_2(a)
\ea\right), \quad a\in A,\myeqno\label{vpp}
$$
где $\X_1$ и $\X_2$ --- некоторые представления алгебры $A$. При этом $\X_2$ соответствует $A$-модулю $U$,
а $\X_1$ --- фактормодулю $V/U$. Таким образом, приводимое представление
$\X$ эквивалентно представлению \ref{vpp}, которое в дальнейшем для краткости мы будем записывать в виде
$$
\left(\ba{cc}
\X_1 & * \\
\O & \X_2
\ea\right).
$$
Обратно, представление алгебры $A$, эквивалентное представлению такого вида, очевидно, приводимо.

 Аналогично, представление $\X$ разложимо тогда и только тогда, когда оно эквивалентно представлению блочного
вида
$$\left(\ba{cc}\X_1 & \O \\
\O & \X_2
\ea\right),
$$
которое мы будем обозначать через $\X_1\oplus\X_2$ и называть
\glsadd{iSumDirRep}\mem{прямой суммой представлений}  $\X_1$ и $\X_2$.

Если $\X$ --- произвольное представление, то оно эквивалентно
блочно-верхнетреугольному представлению

$$\left(\ba{cccc}
\X_1 & * & \cdots & * \\
\O & \X_2 & \cdots & * \\
\cdots&\cdots & \cdots &\cdots \\
\O &  \O & \cdots & \X_k
\ea\right)
$$
с неприводимыми представлениями $\X_1,\ld,\X_k$ на диагонали. Чтобы показать это, достаточно, как и выше,
взять у модуля $V$, соответствующего $\X$, взять композиционный ряд
$$
V=V_0>V_1> \ld > V_k=0,
$$
и <<согласованный>> с ним базис. Неприводимые представления $\X_i$,
соответствующие факторам $V_{i-1}/V_i$, $i=1,\ld,k$, будем называть \mem{неприводимыми} (\mem{простыми})
\mem{компонентами}\glsadd{iCompIrrRep}  представления~$\X$. Из теоремы Жордана-Гёльдера
\ref{thm zh-g} следует, что набор неприводимых компонент представления $\X$ определён однозначно с точностью до
эквивалентности. В используемых обозначениях неприводимую компоненту $\X_1$ (компоненту $\X_k$) мы будем
называть \glsadd{iCompIrrRepUpp}\mem{верхней}
\glsadd{iCompIrrRepLow}(\mem{нижней}). Верхняя (нижняя) неприводимая
компонента представления соответствует неприводимому фактормодулю  (подмодулю) модуля $V$ и,
вообще говоря неверно, что она определена однозначно с точностью до эквивалентности, как можно
увидеть на примере регулярного представления полупростой, но не простой алгебры.

Пусть $\X$ --- вполне приводимое представление алгебры $A$. Поскольку модуль $V$, соответствующий $\X$,
изоморфен прямой сумме неприводимых модулей $V_1,\ld,V_k$, то $\X$ эквивалентно блочно-диагональному
представлению
$$\left(\ba{cccc}
\X_1 & \O & \cdots & \O \\
\O & \X_2 & \cdots & \O \\
\cdots&\cdots & \cdots &\cdots \\
\O &  \O & \cdots & \X_k
\ea\right),
$$
где $\X_i$ --- неприводимое представление, соответствующее модулю $V_i$, $i=1,\ld,k$. Обратно, если некоторое
представление $\X$ эквивалентно представлению такого вида  с неприводимыми компонентами $\X_i$, то оно вполне
приводимо.

\section{Характеры представлений алгебр}

\glsadd{iTrcMat}\mem{Следом} $n\times n$-матрицы $B=(\b_{ij})$ называется сумма
$\gls{trB}=\b_{11}+\ld+\b_{nn}$ её диагональных элементов.

\textuprn{\label{tr}Пусть $B,C$ --- $n\times n$-матрицы над полем $F$. Доказать, что $\tr (BC)=\tr (CB)$. В частности,
следы сопряжённых матриц совпадают.}

Пусть $A$ --- $F$-алгебра. \glsadd{iChrRepAlg}
\mem{Характером}
\mem{представления} $\X$
\mem{алгебры}
$A$ называется\footnote{Характером нулевого представления алгебры $A$ является нулевое отображение $A\to F$.}
отображение $\x: A\to F$, заданное соотношением $\x(a)=\tr\X(a)$,
для всех $a\in A$.


Соответствие между представлениями и модулями и предложение \ref{mod rep} позволяют говорить о
\glsadd{iChrAMod}\mem{характерах} произвольных \mem{$A$-модулей}.

\textzam{Отметим, что характер $\x$ некоторого представления $\X$ алгебры $A$ является $F$-линейным
отображением $A\to F$, но если $\deg \X>1$, то, вообще говоря, не является гомоморфизмом
$F$-алгебр.}

Если $\x_1$ и $\x_2$ --- характеры представлений $\X_1$ и $\X_2$, то их сумма $\x_1+\x_2$
также является характером, например, представления $\X_1\oplus\X_2$. Поэтому множество
характеров алгебры замкнуто относительно сложения.

\textuprn{Доказать, что характеры эквивалентных представлений совпадают.}

\textext{\label{neqch}
Вообще говоря, неэквивалентные представления $F$-алгебры $A$ могут обладать одинаковыми характерами. Например,
в качестве $A$ можно рассмотреть групповую алгебру $FG$, где $F$ --- поле характеристики $2$  и $G$ ---
циклическая группа порядка $2$, порождённая элементом $g$. Пусть представления $\X$ и $\Y$ алгебры $A$ задаются
на порождающем элементе $g$ равенствами
$$\X(g)=\left(
          \begin{array}{cc}
            1 & 0 \\
            0 & 1 \\
          \end{array}
        \right),\qquad
  \Y(g)=\left(
          \begin{array}{cc}
            1 & 1 \\
            0 & 1 \\
          \end{array}
        \right).
$$
Тогда, очевидно,  $\X\ncong\Y$, в то время как характеры этих представлений тождественно равны нулю. Однако,
имеет место следующее утверждение.}

\begin{pre} \label{rep trace} Пусть $A$ --- $F$-алгебра, где $F$ --- алгебраически замкнутое поле.
Пусть $\X_1,\ld,\X_s$ --- некоторые попарно неэквивалентные неприводимые представления алгебры $A$ и $\x_1,\ld,\x_s$
--- соответствующие им характеры. Тогда
\begin{list}{{\rm(}{\it\roman{enumi}\/}{\rm)}}
{\usecounter{enumi}\setlength{\parsep}{2pt}\setlength{\topsep}{5pt}\setlength{\labelwidth}{23pt}}
\item все представления $\X_i$ сюръективны,
\item существуют элементы $a_1,\ld,a_s\in A$ такие, что $\x_i(a_j)=\d_{ij}$,
\item характеры $\x_1,\ld,\x_s$ являются линейно независимыми. В частности, они все
ненулевые и попарно различны.
\end{list}
\end{pre}
\begin{proof} $(i)$ Пусть $M$ --- неприводимый $A$-модуль, соответствующий представлению $\X_i$.
Достаточно показать, что $A_M=\End_F(M)$. Заметим, что $M$ также является неприводимым
$A_M$-модулем, и $A_M$ --- простая алгебра в силу \ref{mlt npr}$(ii)$. Поэтому по следствию
\ref{cor sim azp} получаем $A_M=\End_F(M)$.

$(ii)$ Так как для всех неприводимых $A$-модулей $M_j$ имеет место включение $\J(A)\le \Ann(M_j)$,
из предложения \ref{ann pr}$(iv)$ следует,
что модули $M_j$ можно рассматривать как (неприводимые) модули полупростой алгебры $A/\J(A)$.

Пусть теперь $M_i$ --- неприводимый $A$-модуль, а $\x_i$ --- характер, соответствующие
представлению $\X_i$.
Выберем элемент $a_i\in A$ так, чтобы его образ при гомоморфизме $A\to
A/\J(A)$ лежал в $M_i$-однородной компоненте алгебры $A/\J(A)$ и действовал на
$M_i$ как какое-нибудь наперёд заданное преобразование из $\End_F(M_i)$ со следом $1$
(последнее возможно в силу $(i)$). Тогда $a_i$ аннулирует все неприводимые модули $M_j$, не изоморфные $M_i$,
что следует из  теоремы \ref{thm wedd}$(ii)$. Поэтому получаем
$\x_j(a_i)=\d_{ij}$.

$(iii)$ Если для некоторых $\a_j\in F$ выполнено равенство $\a_1\x_1+\ld+\a_s\x_s=0$, то, взяв значение обеих
частей на элементе $a_j$, определённом в $(ii)$, получим $\a_j=0$.
\end{proof}

\textext{Линейная независимость характеров неприводимых представлений групповых $F$-алгебр имеет место
и без предположения алгебраической замкнутости поля $F$, однако доказывается это
более сложно, см. \cite[Следствие 9.22]{i}.
Для произвольных же $F$-алгебр это, вообще говоря, неверно.
Пусть $F$ --- поле рациональных дробей $\FF_p(t)$, где $p$ --- простое число.
Пусть $A=F(\l)$, где $\l$ --- корень неприводимого многочлена $x^p-t\in F[x]$.
Т.\,к. алгебра $A$ --- поле, её регулярный модуль неприводим.
Можно проверить, что характер соответствующего (неприводимого) представления тождественно равен нулю.}

\begin{pre} \label{cent irr}
Пусть $\X:A\to \MM_n(F)$ --- неприводимое представление алгебры $A$ над алгебраически замкнутым полем $F$. Если
$M\in \MM_n(F)$ удовлетворяет $\X(a)M=M\X(a)$ для всех $a\in A$, то $M$ --- скалярная матрица. В частности, для
любого $z\in \Z(A)$ образ $\X(z)$ --- скалярная матрица. Следовательно, если алгебра $A$ коммутативна, то
$\deg\X=1$.
\end{pre}

\textupl{cent irr prf}{Доказать предложение \ref{cent irr}.}

Из предложения \ref{ann el} вытекает

\begin{pre} \label{irr rep} Всякое неприводимое представление $F$-алгебры $A$ является
верхней неприводимой компонентой её регулярного представления.
\end{pre}

\section{Представления и характеры групп}

Всюду далее символом $G$ будем обозначать некоторую конечную группу.

\begin{opr}
\glsadd{iRepGr}\mem{Представлением группы} $G$ над полем $F$ (или \glsadd{iFRepGr}\mem{$F$-представлени\-ем}) называется гомоморфизм групп
$$
\X : G\to \GL_n(F)
$$
для некоторого 
$n$.\end{opr} Представление $\X$ группы $G$ однозначно продолжается по
линейности до представления степени $n$ алгебры $FG$, которое мы, как правило, будем обозначать тем же символом
$\X$. Наоборот, если $\X$ --- представление алгебры $FG$, то его ограничение на $G$ будет
представлением группы. Поэтому все понятия, определённые для представлений алгебр (степень,
эквивалентность, неприводимость, и т.\,д.) можно перенести на представления групп.

\textuprn{Пусть $\X$ --- неприводимое представление группы $G$ над произвольным полем
$F$. Найти $\sum_{g\in G}\X(g)$.}

%


Ядро представления $\X$ группы $G$ будем обозначать через
\gls{kerX}, чтобы подчеркнуть его отличие от ядра $\Ker \X$ представления $\X$ алгебры $FG$.
Легко видеть, что ядра эквивалентных представлений группы совпадают.
Представление $\X$ называется \glsadd{iRepGrFth}\mem{точным}, если $\ker\X=1$.

\textext{Представляется естественным вопрос о связи ядер $\Ker\X$ и $\ker\X$. Ясно, что если $g\in \ker\X$, то $1-g\in \Ker\X$. Однако
в общем случае $\Ker\X$ как правило больше идеала алгебры $FG$, порождённого элементами  $1-g$ для $g\in \ker\X$.
Достаточно взять точное представление группы, квадрат степени которого меньше её порядка.
Например, для $F=\FF_2$ в качестве $\X$ можно рассмотреть изоморфизм $S_3\cong \GL_2(2)$.}

Пусть $\X$ --- представление группы $G$ над полем $F$.
Отображение $\x: G\to F$, задаваемое правилом
$\x(g)=\tr\X(g)$ для всех $g\in G$ называется \glsadd{iChrRepGr}\mem{характером представления} $\X$ \mem{группы}
$G$ или просто \glsadd{iFChrGr}\mem{$F$-характером группы} $G$.
Множество всех $F$-характеров группы $G$ будем обозначать через \gls{chFGr}.

Как следует из определения, характер представления $\X$ группы $G$ над полем $F$
совпадает с ограничением на $G$ характера
представления $\X$ алгебры $FG$.

Отметим, что, как и в случае алгебр, $F$-характеры групп, вообще говоря, не являются гомоморфизмами
групп.\footnote{Cм., однако, \ref{pr phg}$(ii)$.}

\subsection{Примеры представлений и характеров групп} \label{pr phg}
\begin{list}{{\rm(}{\it\roman{enumi}\/}{\rm)}}
{\usecounter{enumi}\setlength{\parsep}{2pt}\setlength{\topsep}{5pt}\setlength{\labelwidth}{23pt}}
\item У любой группы $G$ есть
\glsadd{iRepGrTrv}\mem{тривиальное $n$-мерное представление}, заданное
по правилу $\X(g)=\I_{n}$ для всех $g\in G$.  Это представление соответствует тривиальному
$n$-мерному $FG$-модулю.

\item Будем называть \glsadd{iFChrLin}\mem{линейным} $F$-характер $1$-мерного представления группы $G$.
Всякий линейный 
$F$-харак\-тер 
является гомоморфизмом групп $G\to F^\times$. Множество всех линейных $F$-характеров группы $G$ будем
обозначать \gls{LinFG}.

\item Тривиальное $1$-мерное $F$-представление группы $G$ будем называть
\glsadd{iRepGrPrnc}\mem{главным}, а его
характер \gls{1G} --- \glsadd{iFChrGrPrnc}\mem{главным} $F$-характером группы $G$. При этом $1_G(g)=1$ для всех $g\in G$.
Главное представление соответствует главному $FG$-модулю, см. \ref{pr am}.

\item \mem{Регулярное $F$-представление}\glsadd{iFRepGrReg}
группы $G$, которое по определению является ограничением на $G$ регулярного
представления алгебры $FG$. Обозначим через \gls{RF} регулярное $F$-представление,
соответствующее модулю $FG^\circ$ в естественном базисе из элементов группы $G$.
Тогда для любого $g\in G$ имеем
$$
\R_F(g)=(a_{xy})_{x,y\in G}, \qquad a_{xy}=\left\{\ba{ll}
1_F,& y=xg;\\
0_F,& y\ne xg.
\ea\right.
$$
Характер \gls{rhF} регулярного $F$-представления группы $G$ будем называть
\glsadd{iFChrReg}\mem{регулярным $F$-характером группы}~$G$. При этом для любого $g\in G$
$$
\r_{\mbox{}_F}(g)=\left\{\ba{ll}
|G|\cdot1_F,& g=1;\\
0_F,& g\ne 1.
\ea\right.
$$
\item Если  $\X$ --- $F$-представление группы $G$ и $H\le G$, то обозначим через
\gls{XsH} ограничение представления $\X$ на $H$. Тогда $\X_H$ является $F$-представлением группы $H$,
и если $V$ --- $FG$-модуль, соответствующий представлению $\X$, то $V_H$ --- $FH$-модуль, соответствующий
$\X_H$.
\item Если $N\nor G$, $\vf:G\to G/N$ --- естественный гомоморфизм групп
и $\wt\X$ --- $F$-представление факторгруппы $G/N$, то композиция $\X=\wt\X\circ\vf$
является $F$-представлением группы $G$.
При этом $\deg\X=\deg\wt\X$ и $N\se \ker\X$. Обратно, если $\X$ --- $F$-представление группы $G$,
и $N\se \ker\X$, то существует единственное
$F$-представление $\wt\X$ факторгруппы $G/N$ такое, что $\X=\wt\X\circ\vf$.
В обоих случаях $\X$ неприводимо
(приводимо, вполне приводимо, неразложимо) тогда и только тогда, когда таковым является $\wt\X$,
ср. \ref{pr pr}$(iii)$.
\item Пусть $\X$ --- $F$-представление группы $G$ и $V$ --- соответствующий ему $FG$-модуль.
Представление, соответствующее контрагредиентному модулю $V^*$ называется
\glsadd{iRepCntGrd}\mem{контрагредиентным представлением} для $\X$ и обозначается \gls{Xcon}.
В двойственном базисе\footnote{Cм. \ref{dv bas}.} модуля $V^*$ контрагредиентное представление $\X^*$
задаётся равенством\footnote{Напомним, что $M^\i=(M^{-1})^\top$
для квадратной невырожденной матрицы $M$.}
 $$\X^*(g)=\X(g)^\i,$$ а характеры $\x$ и $\x^*$ представлений $\X$ и $\X^*$  связаны соотношением
$$
\x^*(g)=\x(g^{-1})
$$
для всех $g\in G$.
\item Пусть группа $G$ действует на конечном множестве $X$ и $V$ --- соответствующий подстановочный
$FG$-модуль с базисом $\{e_x\mid x\in X\}$, см. \ref{pr am}$(vii)$. Характер $\x$ модуля $V$ мы будем
называть \glsadd{iChrPermFG}\mem{подстановочным характером групповой алгебры} $FG$,
соответствующим данному действию группы на множестве (или
\glsadd{iFChrAct}\mem{подстановочным $F$-характером действия} $G$ на $X$). Для любого $g\in G$ справедливо равенство
$$
\x(g)=\big|\{\,x\in X\mid xg=x\,\}\big|\cdot1_F.
$$
Частными случаями этой конструкции являются примеры $(iii)$, $(iv)$.

\end{list}
\textupl{pr phg prf}{Проверить сформулированные в \ref{pr phg} утверждения.}

\textupr{Показать, что группа кватернионов $Q_8$ имеет единственное с точностью до эквивалентности 
неприводимое $\CC$-представление степени 2.}

Обозначим через \gls{IrrFlGr} множество характеров всех неприводимых $F$-представлений группы $G$.
Ясно, что $\lin_F(G)\se\irr_F(G)$.

\textzam{Вообще говоря, характер неприводимого $F$-представления может одновременно быть
характером приводимого $F$-представления.\footnote{По этой причине мы в общем случае не называем элементы из $\irr_F(G)$
неприводимыми $F$-характерами.} Например, если $F$ имеет положительную характеристику $p$,
то таким будет главный $F$-характер, являющийся также характером тривиального $(p+1)$-мерного
представления.}

Из сделанных
выше замечаний о приведении произвольного представления к блочно-верхнетреугольному виду с неприводимыми
компонентами на диагонали следует, что всякий $F$-характер $\psi$ группы представляется в виде
$$
\psi=\sum_{\x\in \irr_F(G)}n_\x \x, \myeqno\label{c dec}
$$
для некоторых неотрицательных целых чисел $n_\x$.

Из \ref{tr} следует, что $F$-характеры группы являются примерами так называемых классовых функций.
\glsadd{iFuncCls}\mem{Классовой функцией} на группе $G$ со значениями в коммутативном кольце $R$ называется отображение
$\vf: G\to R$, принимающее для любого класса сопряжённости группы $G$ одинаковые значения на всех
элементах этого класса. Множество всех классовых функций на группе $G$ со
значениями в кольце $R$ будем обозначать через \gls{cfRlGr}.

Всякую  классовую функцию $\vf\in \cf_R(G)$ $\big($как и вообще любую
функцию из $\f_R(G)\,\big)$ можно естественным образом продолжить до $R$-линейного отображения из $RG$ в
$R$, которое мы без специальных оговорок будем обозначать тем же символом $\vf$. 
Подчеркнём, что для произвольного элемента $x\in RG$ 
равенство\footnote{Здесь $\vf\psi \in \cf_R(G)$ продолжается по $R$-линейности на $RG$. 
Мы не будем рассматривать произведение элементов из $\cf_R(RG)$ в смысле определения \ref{ex r-alg}$(vii)$.} 
$(\vf\psi)(x)=\vf(x)\psi(x)$,
вообще говоря, уже не имеет места.

Если $\vf\in \cf_R(G)$ и $H\le G$, то обозначим через \gls{vfsH} ограничение функции $\vf$ на подгруппу $H$.

\begin{pre} \label{cff} Пусть $G$ --- группа и
$R$ --- коммутативное кольцо. Тогда
\begin{list}{{\rm(}{\it\roman{enumi}\/}{\rm)}}
{\usecounter{enumi}\setlength{\parsep}{2pt}\setlength{\topsep}{5pt}\setlength{\labelwidth}{23pt}}
\item $\cf_R(G)$ является подалгеброй в $R$-алгебре $\f_R(G)$ всех $R$-значных функций на $G$;
\item если $H\le G$, то отображение,  $\vf\mapsto\vf_{H^{\vphantom{A^a}}}$, где $\vf\in\cf_R(G)$,
будет гомоморфизмом $R$-алгебр $\cf_R(G)\to\cf_R(H)$.
\item если $R$ --- поле, то $\dim_R\,\cf_R(G)=|\K(G)|$.
\end{list}
\end{pre}
\textupr{Доказать предложение \ref{cff}.}

В дальнейшем, как правило, мы будем рассматривать классовые функции со значениями в поле $F$.
Алгебру $\cf_F(G)$ будем называть \glsadd{iFAlgClsFncGr}\mem{$F$-алгеброй классовых функций на
группе} $G$.

\textuprn{Пусть $V$ --- $FG$-модуль с характером $\x$ и $H\le G$.
Показать, что характер $FH$-модуля $V_H$ совпадает с~$\x_{H^{\vphantom{A^a}}}$.}

Из предложения
\ref{rep trace}$(iii)$, следствия \ref{cor z fg}, теоремы \ref{thm mash} и замечаний после
предложения \ref{ss wed} вытекает

\begin{pre} \label{irr basis} Пусть $G$ --- конечная группа и $F$ --- алгебраически замкнутое
поле характеристики, не делящей $|G|$. Тогда справедливы следующие утверждения.
\begin{list}{{\rm(}{\it\roman{enumi}\/}{\rm)}}
{\usecounter{enumi}\setlength{\parsep}{2pt}\setlength{\topsep}{5pt}\setlength{\labelwidth}{23pt}}
\item $\r_{\mbox{}_F}=\sum_{\x\in \irr_F(G)}\x(1)\x$.
\item $|G|\cdot 1_F=\sum_{\x\in \irr_F(G)}\x(1)^2$.
\item $\irr_F(G)$ --- базис алгебры $\cf_F(G)$. В частности, в разложении {\rm\ref{c dec}}
произвольного $F$-характера $\psi$ образы коэффициентов $n_\x$ в поле $F$ $($т.\,е. элементы $n_\x\cdot1_F\in F)$
определяются однозначно.
\item Если характеристика поля $F$ равна нулю, то неэквивалентные $F$-представления
группы $G$ обладают различными $F$-характерами.
\end{list}
\end{pre}

Таким образом,  предложение \ref{irr basis} демонстрирует ключевую идею теории характеров.
Хотя значение $F$-характера представления $\X$ на каждом конкретном элементе $g\in G$
несёт мало информации о матрице $\X(g)$, совокупность этих значений по всем элементам группы $G$
может по существу целиком определять представление $\X$.

\textuprn{\label{sopr char}Пусть $G$ --- конечная группа и $F$ --- алгебраически замкнутое
поле характеристики, не делящей $|G|$. Показать, что два элемента $g,h\in G$ сопряжены тогда и только тогда, когда $\x(g)=\x(h)$ для всех $\x \in \irr_F(G)$. \uk{Если $g$ и $h$ не сопряжены, то существует классовая функция на $G$, принимающую различные значения на $g$ и $h$.}}

Как и для алгебр, множество $F$-характеров группы замкнуто относительно
сложения. Покажем, что оно также замкнуто относительно умножения. Для этого нам понадобится
конструкция тензорного произведения модулей групповой алгебры.

\section{Тензорное произведение \texorpdfstring{$FG$}{FG}-модулей}

Пусть $V$ и $W$ --- $FG$-модули. Мы построим новый $FG$-модуль
\gls{VotW}, называемый \glsadd{iProdTensFMod}\mem{тензорным произведением модулей} $V$ и
$W$. Выберем $F$-базисы $v_1,\ld,v_m$ и $w_1,\ld,w_n$ модулей $V$ и $W$, соответственно. Пусть $V\ot W$
обозначает векторное пространство над полем $F$ размерности $nm$, базисом которого являются векторы,
обозначаемые $v_i\ot w_j$. Другими словами, $V\ot W$
--- это множество всевозможных формальных линейных комбинаций $\sum
\a_{ij}(v_i\ot w_j)$, где $\a_{ij}\in F$. Если $v\in V$ и $w\in W$, то пусть $v=\sum\a_iv_i$
и $w=\sum\b_j w_j$ и положим по определению
$$
v\ot w=\sum_{i,j}\a_i\b_j(v_i\ot w_j)\in V\ot W.
$$
\textzam{Отметим, что, вообще говоря, не любой элемент из $V\ot W$ имеет вид $v\ot w$ для некоторых
$v\in V$ и $w\in W$.}
Определим действие
группы $G$ на $V\ot W$, полагая
$$
(v_i\ot w_j)g=v_ig\ot w_jg
$$
для всех $g\in G$ и продолжая это действие по линейности на всё пространство $V\ot W$ и всю
групповую алгебру $FG$. В частности, $(v\ot w)g=vg \ot wg$ для любых $v\in V$, $w\in W$ и
$g\in G$. Таким образом на $V\ot W$ задаётся структура $FG$-модуля.

\textuprn{Проверить корректность данного определения. Доказать, что $FG$-модуль $V\ot W$ определяется модулями
$V$ и $W$ с точностью до изоморфизма, т.\,е. не зависит от выбора базисов пространств $V$ и
$W$.}

\textzam{Важно отметить, что не для любого $a\in FG$ выполнено равенство $(v\ot w)a=va \ot wa$. Поэтому мы
сначала задали действие на $V\ot W$ для групповых элементов, и затем продолжили его по линейности. Если
же $A$ -- произвольная алгебра, а $V$ и $W$ --- её модули, то задать на пространстве $V\ot W$ структуру
$A$-модуля, которая в случае $A=FG$ совпадала бы с построенной выше, удаётся, вообще говоря, не всегда.}

Переведём понятие тензорного произведения $FG$-модулей на язык представлений.

Напомним, что \glsadd{iKronProdMat}\mem{кронекеровым произведением} матриц $A=(\a_{ij})\in
\MM_m(F)$ и $B=(\b_{kl})\in \MM_n(F)$ называется $mn\times mn$-матрица
$$
\gls{AoB}=
\left(
  \begin{array}{ccc}
    \a_{11}B & \ld & \a_{1m}B \\
    \ld & \ld & \ld \\
    \a_{m1}B & \ld & \a_{mm}B
  \end{array}
\right).
$$

\begin{pre} \label{tnz prop} Для любых $A,C\in \MM_m(F)$, $B,D\in \MM_n(F)$, $\a\in F$
\begin{list}{{\rm(}{\it\roman{enumi}\/}{\rm)}}
{\usecounter{enumi}\setlength{\parsep}{2pt}\setlength{\topsep}{5pt}\setlength{\labelwidth}{23pt}}
\item $(A\ot B)^\top=A^\top\ot B^\top$;
\item $\tr(A\ot B)=\tr A\cdot \tr B$;
\item $\det(A\ot B)=(\det A)^n\cdot (\det B)^m$;
\item $(A\ot B)(C\ot D)=AC\ot BD$;
\item $(A+C)\ot B=A\ot B+C\ot B$, \quad $B\ot (A+C)=B\ot A+B\ot C$;
\item $(\a A)\ot B=A\ot (\a B)=\a(A\ot B)$;
\item если  $\mu_1,\ld,\mu_m$ и $\nu_1,\ld,\nu_n$ --- характеристические
значения матриц $A$ и $B$, соответственно, то все характеристическими значения матрицы $A\ot B$ исчерпываются
величинами $\mu_i\nu_j$ при $i=1,\ld,m$ и $j=1,\ld,n$.
\end{list}
\end{pre}
\textupr{Доказать предложение \ref{tnz prop}.}

Пусть $\X$ и $\Y$ --- $F$-представления группы $G$ степеней $m$ и $n$, соответственно. Положим
$$
(\gls{XotY})(g)=\X(g)\ot\Y(g)
$$
для всех $g\in G$. Из \ref{tnz prop} следует, что таким образом определено отображение
$$\X\ot\Y:\,G\to \GL_{mn}(F),$$ и это отображение является $F$-представлением
группы $G$. Назовём $\X\ot\Y$ \glsadd{iTensProdRep}\mem{тензорным произведением}  представлений $\X$ и $\Y$.

\begin{pre} \label{tnz har} Пусть $V$ и $W$  ---  $FG$-модули c
базисами $v_1,\ld,v_m$ и $w_1,\ld,w_n$ и характерами $\x$ и $\psi$, соответственно. Пусть $\X$ и $\Y$
--- представления группы $G$, соответствующие $FG$-модулям $V$ и $W$ в выбранных базисах. Тогда
верны следующие утверждения.
\begin{list}{{\rm(}{\it\roman{enumi}\/}{\rm)}}
{\usecounter{enumi}\setlength{\parsep}{2pt}\setlength{\topsep}{5pt}\setlength{\labelwidth}{23pt}}
\item $F$-представление группы $G$, соответствующее модулю $V\ot W$ в базисе $$v_1\ot w_1,\ \ld,\ v_1\ot w_n,
\ \ld,\ v_m\ot w_1,\ \ld,\ v_m\ot w_n,$$ совпадает с $\X\ot\Y$.
\item Характер представления $\X\ot\Y$ равен $\x\psi$.
\end{list}
\end{pre}

\textupl{tnz har prf}{Доказать предложение \ref{tnz har}.}

Таким образом, предложение \ref{tnz har}$(ii)$ показывает, что
произведение $F$-характеров группы $G$ является\footnote{Вместе с тем множество $\irr_F(G)$, вообще говоря,
не замкнуто относительно умножения. Например, квадрат неприводимого $\CC$-характера степени $2$ группы $S_3$,
соответствующего разностному подмодулю её естественного
подстановочного модуля, не является характером никакого неприводимого представления.}
характером.

\textuprn{Показать, что
\begin{list}{{\rm(}{\it\roman{enumi}\/}{\rm)}}
{\usecounter{enumi}\setlength{\parsep}{2pt}\setlength{\topsep}{5pt}\setlength{\labelwidth}{23pt}}
\item $\lin_F(G)$ замкнуто относительно умножения и является группой;
\item если $\l\in\lin_F(G)$ и $\th\in \irr_F(G)$, то $\l\th\in\irr_F(G)$.
\end{list}}

Несмотря на то, что множество характеров группы $G$ замкнуто относительно сложения и
умножения, кольцом оно не является, поскольку разность характеров, вообще говоря, не является
характером. Иногда удобно рассматривать множество целочисленных линейных комбинаций
характеров неприводимых $F$-представлений группы $G$, которое будем обозначать\footnote{Сокращение от английского {\em generalised characters}.}
через $\gls{gchFG}$, а его элементы называть
\glsadd{iFChrGen}\mem{обобщёнными $F$-характерами}. Из сказанного выше следует, что
$\gch_F(G)$ является \glsadd{iRngGenFChr} кольцом.

\textuprn{Доказать, что классовая функция $\vf\in \cf_F(G)$ является обобщённым $F$-характером тогда и только тогда,
когда $\vf=\x-\th$ для некоторых $F$-характеров $\x$ и $\th$ группы $G$.}

\section{Индуцированные модули и классовые функции}

Пусть $V$ --- $FG$-модуль и $V=W_1\oplus\ld\oplus W_k$, где $W_1,\ld,W_k$ ---
подпространства пространства $V$. Предположим, что для любых $g\in G$ и $i=1,\ldots,k$ образ $W_ig$ совпадает с одним из подпространств 
$W_j$. Другими словами, группа $G$ действует на множестве $\{W_i\mid i=1,\ld,k\}$. Также предположим, что это действие транзитивно. 
В этом случае набор подпространств $W_1,\ld,W_k$ называется \glsadd{iImpSstMod}\mem{системой
импримитивности} степени $k$ модуля~$V$. Неприводимый $FG$-модуль $V$ называется
\glsadd{iModPrim}\mem{примитивным}, если он не имеет системы
импримитивности степени, большей $1$, в противном случае он называется
\glsadd{iMidImpr}\mem{импримитивным}.

\textzam{Подчеркнём, что если $W_1,\ld,W_k$  --- система импримитивности модуля $V$ степени, большей $1$, то
подпространства $W_i$ не являются $FG$-подмодулями.}

\begin{pre} \label{sst imp}
Пусть $V$ --- $FG$-модуль с системой импримитивности $W_1,\ld,W_k$. Обозначим через $W$
одно из подпространств $W_i$, и через $B$ --- некоторый его базис. Пусть $I=\{g\in G\, |\,
Wg=W\}$
--- стабилизатор подпространства $W$. Тогда имеют место следующие утверждения.
\begin{list}{{\rm(}{\it\roman{enumi}\/}{\rm)}}
{\usecounter{enumi}\setlength{\parsep}{2pt}\setlength{\topsep}{5pt}\setlength{\labelwidth}{23pt}}
\item $I\le G$ и $W$ является $FI$-подмодулем модуля $V_I$.
\item $k=|G:I|$.
\item Для любой системы $r_1,\ld, r_k$ представителей правых смежных классов $G$ по $I$
множество $Br_1\cup\ld\cup Br_k$ является базисом пространства $V$.
\item Обозначим через $\X$ представление группы $I$, соответствующее $FI$-подмодулю $W$ в
базисе $B$, и для любого $g\in G$ положим\footnote{Отметим, что $\X^\circ$ не является
представлением группы $G$, если $G\ne I$.}
$$
\X^\circ (g)=\left\{\ba{ll}
\X(g),& g\in I;\\
\O_{\mbox{}_{|B|}},& g\not\in I.
\ea\right.
$$
Выберем представители $r_1,\ld, r_k$ правых смежных классов $G$ по $I$.
Представление $\Y$ группы $G$, соответствующее $FG$-модулю $V$ в базисе $Br_1\cup\ld\cup Br_k$, задаётся
равенством
$$
\Y(g)=\left(
        \begin{array}{ccc}
          \X^\circ(r_1^{\mathstrut} gr_1^{-1}) & \ld & \X^\circ(r_1^{\mathstrut} gr_k^{-1}) \\
          \ld & \ld & \ld \\
          \X^\circ(r_k^{\mathstrut} gr_1^{-1}) & \ld & \X^\circ(r_k^{\mathstrut} gr_k^{-1})
        \end{array}
      \right).
$$
\item Всякий $FG$-модуль, обладающий системой импримитивности такой, что стабилизатор одной из её компонент совпадает с $I$, и соответствующее этой компоненте
представление подгруппы~$I$ эквивалентно $\X$, будет изоморфен~$V$.
\end{list}
\end{pre}

\textupl{sst imp prf}{Доказать предложение \ref{sst imp}}

Утверждение \ref{sst imp}$(v)$ в частности показывает, что импримитивный $FG$-модуль $V$ с точностью до изоморфизма
определяется представлением $\X$ подгруппы $I$. Из следующего предложения видно, что по  любому
представлению собственной подгруппы группы $G$ можно построить некоторый импримитивный $FG$-модуль.

\begin{pre} \label{ind mod} Пусть $H\le G$, $k=|G:H|$ и $W$ --- $FH$-модуль. Тогда существует
единственный с точностью до изоморфизма $FG$-модуль $V$, имеющий систему импримитивности $W_1,\ld,W_k$,
такую, что стабилизатор одного из подпространств этой системы $($которое мы обозначим через  $W_0)$ совпадает с $H$ и
имеет место изоморфизм $FH$-модулей $W_0\cong W$.
\end{pre}
\begin{proof}  Пусть $V$ --- внешняя прямая сумма $k$ копий пространства $W$. Мы
обозначим эти подпространства через $W\ot r$, где $r$ пробегает набор $R$ представителей правых смежных классов $G$ по
$H$, т.\,е. $V=\bigoplus_{r\in R}W\ot r$. Элементы пространства  $W\ot r$ будем записывать
$w\ot r$, т.\,е. $W\ot r= \{w\ot r\,|\, w\in W\}$.

Теперь для $g\in G$ и $w\in W$ определим $w\ot g \in V$ следующим образом. Пусть $g=hr$, где
$h\in H$ и $r\in R$. Тогда положим $w\ot g=wh\ot r$. Отметим, что $w\ot hg=wh \ot g$ для любых
$w\in W$, $h\in H$, $g\in G$. В частности, если $s\in Hr$, то $W\ot s=W\ot r$, и поэтому наша
конструкция не зависит от выбора набора представителей $R$ правых смежных классов.\footnote{На
тензорном языке $V=W\ot _{FH}FG$.}

Пусть $g\in G$. Положим $(w\ot r)g=w\ot rg$ для всех $w\in W$ и $r\in R$ и
продолжим это действие элемента $g$ по линейности на всё пространство $V$. Заметим, что $(w\ot g_1)g_2=w\ot
g_1g_2$ для всех $w\in W$, $g_1,g_2\in G$. Отсюда, как легко видеть, вытекает, что на $V$ задано действие группы $G$. Продолжив это действие
по линейности на всю алгебру $FG$, мы получаем структуру $FG$-модуля. При этом $\{W\ot r\,|\,r\in R\}$
является системой импримитивности модуля~$V$.

Ясно, что стабилизатор подпространства $W_0=W\ot 1$ совпадает с $H$, и изоморфизм
векторных пространств $\vf:W\to W_0$, действующий по правилу $w\vf=w\ot 1$, является
изоморфизмом $FH$-модулей, поскольку для всех $h\in H$, $w\in W$
$$
(w\vf)h=(w\ot 1)h=wh\ot 1=(wh)\vf.
$$

Всякий $FG$-модуль $U$, удовлетворяющий условию
предложения, изоморфен построенному модулю $V$. В самом деле, если $\X$ --- представление группы $H$,
соответствующее $FH$-модулю $W$, то в некоторых базисах представления, соответствующие
модулям $U$ и $V$ совпадут, поскольку будет иметь вид, указанный в \ref{sst imp}$(iv)$.
\end{proof}

%
%

Модуль $V$, существование которого утверждается в предложении \ref{ind mod}, называется $FG$-модулем,
\glsadd{iModInd}\mem{индуцированным} $FH$-модулем $W$, и обозначается \gls{VupG}.
Таким образом, мы получаем следующий критерий.

\begin{cor}\label{cr pr} Пусть $V$ --- $FG$-модуль.
Если $V$ индуцирован с $FH$-модуля собственной подгруппы $H<G$, то $V$ импримитивен. Обратно, если $V$ импримитивен и $W$ --- любая из компонент системы импримитивности, то
$V\cong W^G$, где $W$ рассматривается как 
$FH$-модуль, а $H=\St_G(W)$.
\end{cor}

\begin{proof} 
Как следует из определения, $FG$-модуль, индуцированный с собственной подгруппы $H$, обладает системой импримитивности степени $|G:H|>1$. Обратно, импримитивный модуль $V$ обладает системой импримитивности степени большей $1$. В силу транзитивности действия $G$ стабилизатор $H$ одного (произвольного) из подпространств $W$ этой системы будет собственной подгруппой в $G$. Из \ref{ind mod} следует, что $V\cong W^G$. 
\end{proof}

\begin{pre} \label{ed ind} Пусть $V$ --- подстановочный $FG$-модуль, соответствующий транзитивному
действию группы $G$, и пусть $H$ --- стабилизатор точки. Тогда имеет место изоморфизм $FG$-модулей
$V\cong W^G$, где $W$ --- главный $FH$-модуль. В частности, $FG^\circ\cong U^G$,
где $U$ --- главный $FE$-модуль единичной подгруппы $E$ группы $G$.
\end{pre}
\textupr{Доказать предложение \ref{ed ind}.}

\textuprn{Показать, что модуль, индуцированный неприводимым модулем, может не быть даже вполне приводимым.
\uk{Воспользоваться предложениями \ref{ed ind}, \ref{mod ns} и \ref{ssm crit}.}}

Следующее упражнение обобщает предложение \ref{ed ind}.

\textuprn{Пусть $H\le G$ и $F$ --- поле, характеристика которого не делит порядок $|H|$.
\begin{list}{{\rm(}{\it\roman{enumi}\/}{\rm)}}
{\usecounter{enumi}\setlength{\parsep}{2pt}\setlength{\topsep}{5pt}\setlength{\labelwidth}{23pt}}
\item Пусть $V$ --- главный $FH$-модуль и
$e=\displaystyle\frac{1}{|H|}\sum_{h\in H}h$. Показать, что $e$ --- идемпотент алгебры $FG$, и имеет место
изоморфизм $FG$-модулей $eFG\cong V^G$.
\item Пусть $V$ --- одномерный $FH$-модуль с $F$-характером $\l$ и
$e=\displaystyle\frac{1}{|H|}\sum_{h\in H}\l(h^{-1})h$. Чему изоморфен $FG$-модуль $eFG$?
\end{list}
}


\textuprn{Пусть $H \le G$ и $U$ -- $FH$-модуль. На множестве отображений
$$
\{f:G\to U\mid f(gh)=f(g)h\ \text{для всех}\ g\in G, h\in H\}
$$
можно задать структуру $FG$-модуля. $F$-линейность определяется покомпонентно, а действие элемента $g\in G$ на $f$
задаётся правилом $(fg)(x)=f(gx)$ для любых $x\in G$. Показать, что этот $FG$-модуль изоморфен $U^G$.}


\begin{pre}\label{prop ind} Имеют место следующие утверждения.
\begin{list}{{\rm(}{\it\roman{enumi}\/}{\rm)}}
{\usecounter{enumi}\setlength{\parsep}{2pt}\setlength{\topsep}{5pt}\setlength{\labelwidth}{23pt}}
\item Пусть $V$ --- $FG$-модуль. Тогда $V^G\cong V$.
\item Пусть $K\le H \le G$ и $U$ --- $FK$-модуль. Тогда имеет место изоморфизм $FG$-модулей
$(U^H)^G\cong U^G$.
\item Пусть $H\le G$ и $U\le V$ --- включение $FH$-модулей. Тогда $U^G\le V^G$,
причём $U<V$ если и только если $U^G< V^G$.
\item Пусть $H\le G$ и $U, V$ --- $FH$-модули. Тогда $(U\oplus V)^G\cong U^G\oplus V^G$.

\end{list}
\end{pre}

\textupl{prop ind prf}{Доказать предложение \ref{prop ind}.}

\begin{opr} Пусть $G$ --- группа, $H\le G$, $\X$ --- $F$-представление группы $H$, обладающее
характером $\x$ и $W$ --- $FH$-модуль, соответствующий представлению $\X$. Тогда представление
\gls{XuG} группы $G$, соответствующее индуцированному $FG$-модулю $W^G$ в некотором базисе, называется
\glsadd{iRepInd}\mem{представлением, индуцированным} представлением $\X$, а характер
\gls{chiuG} представления $\X^G$ называется \glsadd{iFChrInd}\mem{характером, индуцированным} $F$-характером
$\x$.
\end{opr}

\textzam{Отметим, что в силу предложения \ref{sst imp}$(v)$ представление $\X^G$, индуцированное
$F$-представлением $\X$ подгруппы $H$ группы $G$, эквивалентно
представлению $\Y$, приведённому в \ref{sst imp}$(iv)$. В частности, с точностью до эквивалентности, 
можно считать, что матрица $\X^G(g)$ для любого $g\in G$ 
является блочной $k\times k$-матрицей, где $k=|G:H|$, элементы которой --- $d\times d$-блоки,
где $d=\deg \X$, причём в каждой строке и каждом столбце этой матрицы ровно один блок является ненулевым.}

Строго говоря, определение индуцированного $F$-характера $\x^G$ зависит от представления $\X$,
характером которого является $\x$, а такое представление, как мы видели, не всегда
определяется характером $\x$ однозначно с точностью до эквивалентности. Однако из следующего
утверждения вытекает, что  характер $\x^G$ определяется только характером $\x$ и, значит,
данное выше определение корректно.

\begin{pre} \label{ind char} Пусть $G$ --- группа, $H\le G$,
и $\x$ --- $F$-характер группы $H$.
Для произвольного $g\in G$ определим
$$
\x^\circ (g)=\left\{\ba{ll}
\x(g),& g\in H;\\
0,& g\not\in H.
\ea\right.
$$
Тогда для любого набора $R$ представителей всех правых смежных классов $G$ по $H$
$$
\x^G(g)=\sum_{r\in R} \x^\circ(rgr^{-1}). \myeqno\label{fich}
$$
\end{pre}

\textupl{ind char prf}{Доказать предложение \ref{ind char}.}


По аналогии с формулой \ref{fich} для произвольной классовой функции $\vf\in \cf_F(H)$ подгруппы
$H$ группы $G$ можно определить понятие
индуцированной функции  \gls{vfG}, задав её
равенством
$$
\vf^G(g)=\sum_{r\in R} \vf^\circ(rgr^{-1}), \myeqno\label{fif}
$$
где $R$ --- набор представителей всех правых смежных классов $G$ по  $H$, а функция
$\vf^\circ: G\to F$ совпадает с $\vf$ на $H$ и тождественно равна нулю вне $H$.
Заметим, что значение $\vf^G(g)$  равно нулю, если $g$ не лежит ни в какой сопряжённой
c $H$ подгруппе. Поэтому для любого $g\in G$ определение \ref{fif} можно переписать в виде
$$
\vf^G(g)=\sum_{\{r\in R\,\mid\, g\in H^r\}} \vf(rgr^{-1}). \myeqno\label{fif ud}
$$
Такую форму определения индуцированной функции иногда удобнее использовать в доказательствах.

\begin{pre} \label{ind kl} Пусть $G$ --- группа, $H\le G$ и $\vf\in \cf_F(H)$.
\begin{list}{{\rm(}{\it\roman{enumi}\/}{\rm)}}
{\usecounter{enumi}\setlength{\parsep}{2pt}\setlength{\topsep}{5pt}\setlength{\labelwidth}{23pt}}
\item Функция  $\vf^G$, определяемая в {\rm\ref{fif}}, не зависит от выбора системы представителей
правых смежных классов и $\vf^G\in \cf_F(G)$.
\item Отображение $\cf_F(H)\to \cf_F(G)$, действующее то правилу $\vf\mapsto\vf^G$, является $F$-линейным.
\item $|H|\vf^G(g)=\sum_{x\in G} \vf^\circ(xgx^{-1})=\sum_{\{x\in G\,\mid\, g\in H^x\}} \vf(xgx^{-1})$.
\item Если характеристика поля $F$ не делит $|H|$, то для любого $g\in G$
$$
\vf^G(g)=|\C_G(g)|\sum_{x\in  R_g}\frac{\vf(x)}{|\C_H(x)|},
$$
где $R_g$ --- набор представителей всех классов сопряжённости группы $H$, содержащихся в классе $g^G$.
\item Если $H\le K \le G$, то $\vf^G=(\vf^K)^G$.
\item Если $\psi\in\cf_F(G)$, то $\vf^G\psi=(\vf\,\psi_{H^{\vphantom{A^a}}})^G$.
\item Если $K\le G$ и $G=HK$, то $(\vf^G)_{K^{\vphantom{A^a}}}=(\vf_{H\cap K^{\vphantom{A^a}}})^K$.
\end{list}
\end{pre}

\textupl{ind kl prf}{Доказать предложение \ref{ind kl}.}

Утверждение \ref{ind kl}$(i)$ позволяет называть $\vf^G$ \glsadd{iFuncClsInd}\mem{индуцированной классовой функцией}.

\section{Сопряжённые представления и классовые функции}

Между представлениями сопряжённых подгрупп, а также между их классовыми функциями можно установить
естественное взаимно однозначное соответствие.

Пусть $H\le G$ и $\vf\in \cf_F(H)$. Для любого $g\in G$ определим
отображение $\gls{vfg}:H^g\to F$ по правилу
$$
\vf^g(h^g)=\vf(h)
$$
для всех $h\in H$. Следующее предложение показывает, что отображение  $\vf^g$ принадлежит $\cf_F(H^g)$. Будем
называть его \mem{классовой функцией} подгруппы $H^g$, \glsadd{iFuncClsConjElm}\mem{сопряжённой}
с классовой функцией $\vf$ подгруппы $H$ \mem{элементом} $g$.

\begin{pre}\label{fcn pr} Пусть $H\le G$, $g,f\in G$ и $\vf\in \cf_F(H)$.
Тогда
\begin{list}{{\rm(}{\it\roman{enumi}\/}{\rm)}}
{\usecounter{enumi}\setlength{\parsep}{2pt}\setlength{\topsep}{5pt}\setlength{\labelwidth}{23pt}}
\item $\vf^g\in\cf_F(H^g)$.
\item $(\vf^g)^f=\vf^{gf}$.
\item Если $g\in HC_G(H)$, то $\vf^g=\vf$.
\item $(\vf^g)^G=\vf^G$.
\end{list}
\end{pre}

\textupl{fcn pr prf}{Доказать предложение \ref{fcn pr}.}

Пусть $H\le G$ и $\X$ --- представление группы $H$ над произвольным полем $F$. Для
элемента $g\in G$ определим \glsadd{iRepConjElm}\mem{сопряжённое представление \gls{Xg}} группы $H^g$ по правилу $\X^g(h^g)=\X(h)$.
Легко видеть, что $\X^g$ действительно является $F$-представлением группы $H^g$, причём $\deg\X^g=\deg\X$.

\begin{pre} \label{sop ek} Пусть $H\le G$, $g\in G$ и $\X,\Y$ --- $F$-представления подгруппы $H$.
Тогда
\begin{list}{{\rm(}{\it\roman{enumi}\/}{\rm)}}
{\usecounter{enumi}\setlength{\parsep}{2pt}\setlength{\topsep}{5pt}\setlength{\labelwidth}{23pt}}
\item представление $\X^g$ неприводимо тогда и только тогда, когда неприводимо $\X$;
\item $\X^g\cong\Y^g$ тогда и только тогда, когда $\X\cong\Y$;
\item если $\x$ --- $F$-характер представления $\X$, то $\x^g$ --- $F$-характер представления $\X^g$.
\end{list}
\end{pre}
\textupr{Доказать предложение \ref{sop ek}.}

Пусть $W$ --- $FH$-модуль и $U$ --- $FH^g$-модуль. Будем говорить, что \glsadd{iModConjElm}\mem{модули} $W$ и $U$
\mem{сопряжены}, если некоторые соответствующие
им представления групп $H$ и $H^g$ сопряжены.

\begin{pre} \label{spr mod} Пусть $V$ --- $FG$-модуль, $H\le G$, и $W$ --- подмодуль $FH$-модуля $V_H$.
\begin{list}{{\rm(}{\it\roman{enumi}\/}{\rm)}}
{\usecounter{enumi}\setlength{\parsep}{2pt}\setlength{\topsep}{5pt}\setlength{\labelwidth}{23pt}}
\item Если $g\in G$, то $Wg$ --- $FH^g$-подмодуль модуля $V_{H^g}$. Кроме того, модули $W$ и
$Wg$ сопряжены.
\item Если $g\in G$ и $M$ --- $FH^g$-подмодуль модуля $V_{H^g}$, сопряжённый с $W$, то $M\cong Wt$ для
некоторого элемента $t$ из смежного класса $\N_G(H)g$.
\item Если $M$ --- $FH$-подмодуль модуля $V_H$, изоморфный $W$, то для всех $g\in G$ имеет
место изоморфизм $FH^g$-модулей $Mg\cong Wg$.
\end{list}
\end{pre}

\textupl{spr mod prf}{Доказать предложение \ref{spr mod}.}

Пусть $H\le G$, $g\in G$, $\X$ --- $F$-представление группы $H$, и $W$ --- $FH$-модуль,
соответствующий $\X$. Символом $\gls{Wg}$  будем обозначать\footnote{Здесь $Wg$ понимается как единый символ, который в силу \ref{spr mod}$(i)$
хорошо согласуется со случаем, когда $W$ и $Wg$ содержатся в некотором $FG$-модуле.} $FH^g$-модуль, сопряжённый с $W$ и
соответствующий сопряжённому представлению $\X^g$.
Например, в качестве модуля $Wg$ можно взять само векторное пространство $W$, в котором образ
произвольного элемента $w\in W$ под действием $h^g\in H^g$ совпадает с $wh$. Или же $Wg$ можно рассматривать как множество
формальных произведений $wg$, где $w\in W$, с естественной структурой $F$-пространства и действием $(wg)h^g =(wh)g$ для всех $h^g\in H^g$.

\textuprn{Пусть $H\le G$ и $g\in G$. Показать, что
\begin{list}{{\rm(}{\it\roman{enumi}\/}{\rm)}}
{\usecounter{enumi}\setlength{\parsep}{2pt}\setlength{\topsep}{5pt}\setlength{\labelwidth}{23pt}}
\item если $W$ --- $FH$-модуль, то $W^G\cong (Wg)^G$,
\item если $V$ --- $FG$-модуль, то $V_{H^g}\cong V_H g$.
\end{list}
}

Пусть $H, K$ --- подгруппы группы $G$ и $t\in G$.
\glsadd{iDblCoset}\mem{Двойным смежным классом} группы $G$ по подгруппам $H$ и $K$ с представителем $t$
называется множество
$$HtK=\{htk\,|\,h\in H,\ k\in K\}.$$
Любые два двойных смежных класса либо не пересекаются, либо совпадают, и группа $G$ является
объединением своих двойных смежных классов.

\begin{pre}[Теорема Макки] \label{t mak}\glsadd{iThmMacK}
Пусть $H,K\le G$ и $T$ --- множество представителей всех двойных
смежных классов группы $G$ по $H$ и $K$. Пусть $W$--- $FH$-модуль. Тогда имеет место
изоморфизм $FK$-модулей
$$
(W^G)_K\cong \bigoplus_{t\in T}\ ((Wt)_{H^t\cap K})^K.
$$
\end{pre}

\textupl{t mak prf}{Доказать предложение \ref{t mak}.}

Читателю было бы полезно перевести формулировку теоремы Макки на язык представлений.

Отметим важные частные случаи теоремы \ref{t mak}.

\begin{cor} \label{cor mak}
\begin{list}{{\rm(}{\it\roman{enumi}\/}{\rm)}}
{\usecounter{enumi}\setlength{\parsep}{2pt}\setlength{\topsep}{5pt}\setlength{\labelwidth}{23pt}}
\item Пусть $V$ --- регулярный
$FG$-модуль и $K\le G$. Тогда $V_K$ --- прямая сумма $|G:K|$ регулярных $FK$-модулей.
\item Пусть $H\nor G$ и $T$ --- множество представителей всех смежных классов группы $G$ по
$H$. Пусть $W$~--- $FH$-модуль. Тогда
$$(W^G)_H\cong\bigoplus_{t\in T}\ Wt.$$
\end{list}
\end{cor}

Рассмотрим более подробно важный частный случай, когда $H$ --- нормальная подгруппа группы $G$.
Пусть $\X$ --- $F$-представление подгруппы $H$. В силу нормальности $H$ сопряжение любым
элементом $g\in G$ определяет новое представление $\X^g$. При этом, как вытекает из определения,
$$\X^g(h)=\X(ghg^{-1})$$
для всех $h\in H$. Таким образом, в силу  \ref{sop ek}
группа $G$ действует на множестве классов эквивалентности (неприводимых) $F$-представлений.

\textext{Следующий пример показывает, что это действие, вообще говоря, нетривиально.
Пусть $G=S_3$, $H=A_3\nor G$ и $\X,\Y$ --- одномерные $\CC$-представления подгруппы $H$, задающиеся
на её порождающем элементе $h=(123)$ равенствами $\X(h)=(\z)$ и  $\Y(h)=(\z^{-1})$, где
$\z=e^{\frac{2\pi i}{3}}$ ---
примитивный кубический комплексный корень из $1$. Очевидно, что $\X\ncong\Y$.
С другой стороны, для любого $g\in G\setminus H$ имеем $\X^g=\Y$, так как
$$
\X^g(h)=\X(ghg^{-1})=\X(h^{-1})=(\z^{-1})=\Y(h).
$$}

\mysubsection\label{pgin m} Пусть, как и выше, группа $G$ действует на множестве
классов эквивалентности $F$-представлений нормальной подгруппы $H$.
Стабилизатор класса, содержащего представление $\X$,
обозначим через  $\II_G(\X)$ и
назовём его  \glsadd{iGrInrRep}\mem{группой инерции представления} $\X$.
Другими словами,
$$\gls{IGX}=\{g\in G \mid \X^g\cong\X\}.$$

Аналогично, \glsadd{iGrInrMod}\mem{группой инерции} $FH$-\mem{модуля} $U$
назовём группу
$$\gls{IGU}=\{g\in G \mid Ug\cong U\}.$$
Легко видеть, что если модуль $U$ соответствует представлению  $\X$, то $\II_G(\X)=\II_G(U)$.

%


\mysubsection\label{conj aut} Определённое выше сопряжение классовых функций, представлений и модулей подгруппы $H$ элементами из $G$
можно обобщить на случай сопряжения элементами из $\Aut(G)$. Например, если $\vf\in \cf_F(H)$ и $\a\in \Aut(G)$, то, положив
$$
\vf^\a(h^\a)=\vf(h)\myeqno\label{spr aut}
$$
для всех $h\in H$, получим отображение $\gls{vfal}:H^\a\to F$, которое принадлежит $\cf_F(H^\a)$ и называется \glsadd{iFuncClsConjAut}\mem{классовой функцией,
сопряженной} с $\vf$ \mem{автоморфизмом} $\a$. Аналогичным образом для представления $\X$ подгруппы $H$ и соответствующего ему $FH$-модуля
$W$ определяются \glsadd{iRepConjAut}\mem{сопряжённое представление} $\gls{xal}$ подгруппы $H^\a$ и \glsadd{iModConjAut}\mem{сопряжённый} $FH^\a$-\mem{модуль}
$\gls{Wal}$. Большинство сформулированных в \ref{fcn pr}--\ref{spr mod}
свойств сопряжения элементами группы переносятся на случай сопряжения автоморфизмами. Легко видеть, что если $\a$ --- внутренний автоморфизм
группы $G$, являющийся сопряжением элементом $g\in G$, то $\vf^\a$, $\X^\a$, $W\a$ совпадают, соответственно, с $\vf^g$, $\X^g$, $Wg$.

\section{Ограничение на нормальные подгруппы}

Если $FG$-модуль $V$ неприводим, то о структуре $FH$-модуля $V_H$  для произвольной подгруппы
$H\le G$ в общем случае можно сказать немного. Однако, ситуация меняется в
случае, когда подгруппа $H$ нормальна.

\begin{thm}[Клиффорда] \label{thm cliff}\glsadd{iThmCliff}
Пусть $H\nor G$ и $V$ --- неприводимый $FG$-модуль. Пусть $U$ --- неприводимый $FH$-подмодуль
модуля $V_H$. Тогда справедливы следующие утверждения.
\begin{list}{{\rm(}{\it\roman{enumi}\/}{\rm)}}
{\usecounter{enumi}\setlength{\parsep}{2pt}\setlength{\topsep}{5pt}\setlength{\labelwidth}{23pt}}
\item $V_H=\sum_{g\in G} Ug$. В частности, модуль $V_H$ вполне приводим и любой неприводимый $FH$-подмодуль
модуля $V_H$ сопряжён с $U$.
\item $V_H=\oplus_M\W_M(V_H)$, где $M$ пробегает сопряжённые с $U$ попарно неизоморфные
подмодули из $V_H$. Группа $G$ действует транзитивно на однородных компонентах $\W_M(V_H)$. В частности, для
любого неприводимого $FH$-подмодуля $M\le V_H$ имеет место равенство $\n_M(V_H)=\n_U(V_H)$.
\item Обозначим $W=\W_U(V_H)$  и  $I=\II_G(U)$. Тогда $W$ --- неприводимый $FI$-модуль и $V\cong W^G$.
\end{list}
\end{thm}

\begin{proof} $(i)$ Поскольку $\sum_{g\in G} Ug$ --- ненулевое $G$-инвариантное подпространство в $V$, имеем
$V=\sum_{g\in G} Ug$ ввиду неприводимости модуля $V$. Значит, $V_H$ вполне приводим по
предложению \ref{equ ssm}$(iii)$. Из предложений \ref{cor sum dir} и \ref{ss wed}$(iii)$
следует, что всякий неприводимый подмодуль из $V_H$ изоморфен какому либо слагаемому $Ug$ и,
значит, сопряжён с $U$ по предложению \ref{spr mod}$(i)$.

$(ii)$ Разложение $V_H=\oplus_M\W_M(V_H)$ следует из $(i)$ и предложения \ref{ss wed}$(iii)$.
Установим транзитивность действия группы $G$. Пусть $M=Ug$ для подходящего $g\in G$.
Достаточно показать, что $\W_U(V_H)g=\W_M(V_H)$. Так как $\W_U(V_H)$ --- сумма
$FH$-подмодулей, изоморфных $U$, то из предложения \ref{spr mod}$(iii)$ следует, что
$\W_U(V_H)g$ --- сумма подмодулей, изоморфных $M$. Поэтому $\W_U(V_H)g\se \W_M(V_H)$. Обратное
включение получается аналогично.

В силу $(i)$ все неприводимые $FH$-подмодули $M$ из $V_H$ сопряжены с $U$ и поэтому имеют
одинаковую размерность. Кроме того, из транзитивности действия $G$ на однородных компонентах
следует, что все $\W_M(V_H)$ также имеют одинаковую размерность. Значит, число
$\n_M(V_H)=\dim_F(\W_M(V_H))/\dim_F M$ не зависит от модуля $M$.

$(iii)$ Из доказанного выше следует, что однородные компоненты $FH$-модуля $V_H$ образуют
систему импримитивности для $V$. Заметим, что  $I$ --- стабилизатор компоненты $W$,
поскольку произвольный элемент из $G$ принадлежит $I$ тогда и только тогда, когда он
переводит $U$ в изоморфный ему подмодуль, и значит, оставляет $W$ инвариантным.
Поэтому из предложения
\ref{ind mod} вытекает изоморфизм $V\cong W^G$. Неприводимость $FI$-модуля $W$ теперь следует из предложения
\ref{prop ind}$(iii)$ и неприводимости $FG$-модуля $V$.
\end{proof}

\nopagebreak
\section{Расширение основного поля}

В этом разделе мы будем преимущественно говорить о групповых $F$-алгебрах, хотя многие факты и
определения можно было бы перенести и на произвольные $F$-алгебры.

Пусть $E\ge F$ --- расширение полей и $\X$ --- $F$-представление группы $G$. Тогда $\X$ также
является $E$-представлением, и мы в этом случае будем обозначать его \gls{XE}. Если $F$-представления
$\X$ и $\Y$ эквивалентны, то также эквивалентны и $E$-представления $\X^E$ и $\Y^E$. Поэтому,
если $V$ --- $FG$-модуль, соответствующий представлению $\X$, то существует\footnote{На тензорном языке
имеют место изоморфизмы $EG\cong E\ot_F FG$ и $V^E\cong E\ot_F V$.}
однозначно определённый (с точностью до изоморфизма) $EG$-модуль $V^E$, соответствующий представлению $\X^E$.
Если представление $\X^E$ неприводимо, то также неприводимо и представление $\X$.
Однако, $\X^E$ может оказаться приводимым, даже если $\X$ неприводимо.

\textuprn{\label{pred z}Пусть $C=\la c\,|\, c^3=1\ra$
--- циклическая группа порядка $3$ и $\X$ --- её $\RR$-представление, определённое равенством
$$
\X(c)=\left(\ba{rr}
0 & 1 \\
-1 & -1
\ea\right).
$$
Показать, что $\X$ неприводимо, а $\X^\CC=\X_1\oplus\X_2$, где $\X_1$ и $\X_2$ --- неэквивалентные
неприводимые $\CC$-представления.}

\textuprn{\label{pred q}Пусть $\RR$-представление группы кватернионов
$$
Q_8=\la a,b\bigm| a^4=b^4=1,\ \  a^2=b^2, \ \ a^b=a^{-1} \ra
$$
определено равенствами
$$
\X(a)=\left(\ba{rrrr}
0 & \phantom{-}1 & 0 & 0 \\
-1 & 0 & 0 & 0 \\
0 & 0 & 0 & -1 \\
0 & 0 & \phantom{-}1 & 0
\ea\right)\qquad
\X(b)=
\left(\ba{rrrr}
0 & 0 & \phantom{-}1 & 0 \\
0 & 0 & 0 & \phantom{-}1 \\
-1 & 0 & 0 & 0 \\
0 & -1 & 0 & 0
\ea\right)
$$
Показать, что $\X$ неприводимо, а
$\X^\CC\cong \Y\oplus\Y$, где $\Y$ --- неприводимое $\CC$-представление степени $2$.}

Легко видеть,\footnote{См. \ref{pr phg}$(iv)$.}
что если $E \ge F$ --- расширение полей и $\R_F$ ---
регулярное $F$-представление группы $G$, то $(\R_F)^E=\R_E$ --- её регулярное $E$-представление.

\begin{pre} \label{ext npr} Пусть $E\ge F$ --- расширение полей и $\Y$ ---
неприводимое $E$-представление группы~$G$.
Тогда существует неприводимое $F$-представление $\X$ такое, что $\Y$ является неприводимой
компонентой представления $\X^E$.
\end{pre}
\textupl{ext npr prf}{Доказать предложение \ref{ext npr}. \uk{Использовать предложение \ref{irr rep}.}}

Чуть ниже\footnote{См. \ref{zam pp}.}  мы покажем единственность с точностью до эквивалентности $F$-представления
$\X$ из предложения \ref{ext npr}.

\begin{opr} $F$-представление $\X$ группы $G$ называется \glsadd{iRepGrAbsIrr}\mem{абсолютно неприводимым},
если представление $\X^E$ неприводимо для любого расширения $E$ поля
$F$.
\end{opr}

Всякое абсолютно неприводимое представление группы $G$, очевидно, неприводимо. Обратное неверно, как показывают
примеры \ref{pred z} и \ref{pred q} выше.

\begin{pre} \label{abs npr} Пусть $\X$ --- неприводимое $F$-представление степени $n$ группы $G$.
Тогда следующие условия эквивалентны.
\begin{list}{{\rm(}{\it\roman{enumi}\/}{\rm)}}
{\usecounter{enumi}\setlength{\parsep}{2pt}\setlength{\topsep}{5pt}\setlength{\labelwidth}{23pt}}
\item $\X$ абсолютно неприводимо.
\item $\X^E$ неприводимо для всякого конечного расширения $E \ge F$.
\item Централизатор образа $\X(G)$ в алгебре $\MM_n(F)$ состоит из скалярных матриц.
\item $\X(FG)=\MM_n(F)$.
\item Образ $\X(G)$ содержит $n^2$ линейно независимых над $F$ матриц.
\end{list}
\end{pre}


\begin{proof} $(i)\Rightarrow(ii)$ Очевидно.

$(ii)\Rightarrow(iii)$ Пусть матрица $M\in \MM_n(F)$ удовлетворяет равенству
$\X(g)M=M\X(g)$ для всех $g\in G$. Выберем
конечное расширение $E$ поля $F$, в котором $M$ имеет собственное значение $\l$. Тогда подпространство $U\le
E^n$ собственных векторов отвечающих собственному значению $\l$, ненулевое. Но $U$
инвариантно относительно любой матрицы $\X(g)$, поскольку $$u\X(g)M=uM\X(g)=\l u\X(g)$$ для всех $u\in U$, $g\in
G$. Из неприводимости представления $\X^E$ следует, что $U=E^n$, т.\,е. $M$ --- скалярная матрица и, в
частности, $\l\in F$.

$(iii)\Rightarrow(iv)$ Перейдём на язык модулей. Пусть $A=FG$ и $V=F^n$ --- $A$-модуль, соответствующий
представлению $\X$. Тогда при изоморфизме $\MM_n(F)\cong \End_F(V)$ образу $\X(FG)$ соответствует алгебра $A_V$, и нам достаточно
показать, что $A_V=\End_F(V)$.

По предложению \ref{mlt npr}$(ii)$ алгебра $A_V$ простая. Поскольку $V$ --- неприводимый $A_V$-модуль, то из
теоремы о двойном централизаторе \ref{thm dv cent} следует, что $A_V=\C_{\End_F(V)}(\C_{\End_F(V)}(A_V))$.
Однако из $(iii)$ следует, что $\C_{\End_F(V)}(A_V)$ состоит из скалярных преобразований пространства
$V$. Поэтому $A_V=\End_F(V)$.

$(iv)\Rightarrow(i)$ Пусть $E \ge F$ --- расширение полей. Поскольку любая матрица из $\MM_n(E)$ является
$E$-линейной комбинацией матриц из $\MM_n(F)$, имеем $\X^E(EG)=\MM_n(E)$. Но тогда представление $\X^E$ не может
быть эквивалентно представлению вида
$$\left(\ba{cc}
* & * \\
\O & *
\ea\right)
$$
и, значит, неприводимо.

Эквивалентность $(iv)$ и $(v)$ следует из того, что $\X$ является $F$-линейным и $\dim_F(\MM_n(F))=n^2$.
\end{proof}

Как видно из  предложения \ref{abs npr}$(iii)$--$(iv)$, абсолютную неприводимость
$F$-представления можно проверить, не рассматривая расширения поля $F$.

\begin{opr} Поле $F$ называется \glsadd{iFldDecGr}\mem{полем разложения группы} $G$, если любое неприводимое $F$-представление группы $G$ абсолютно
неприводимо.
\end{opr}

Поле, не имеющее собственных подполей, называется  \glsadd{iFldPrim}\mem{простым}. Легко
видеть, что простые поля исчерпываются полем рациональных чисел $\QQ$ и конечными полями $\FF_p$ для простого $p$.

\textuprn{Доказать, что для
группы $S_3$ любое (простое) поле является полем разложения.}

Поскольку любое конечное расширение алгебраически замкнутого поля $F$ совпадает с $F$, то из
предложения \ref{abs npr}$(ii)$ получаем

\begin{cor} \label{cor algz abs}
Алгебраически замкнутое поле $F$ является полем разложения для любой группы.
\end{cor}

Следующее утверждение оказывается полезным при проверке эквивалентности неприводимых
$F$-пред\-ставлений, 
а также в других ситуациях.

\begin{pre} \label{irr ab0} Пусть $\X$ --- неприводимое представление групповой алгебры $FG$ и
 $a\in FG$. Тогда существует элемент $b\in FG$ такой, что $\X(b)=\X(a)$ и $\Y(b)=0$ для всех
 неприводимых представлений $\Y$ алгебры $FG$, не эквивалентных $\X$.
\end{pre}
\begin{proof} Условимся об обозначениях. Пусть $A=FG/\J(FG)$ и $\vf:FG\to A$ --- естественный эпиморфизм $F$-алгебр.
Для любого неприводимого представления $\Y$ алгебры $FG$
через $\wt\Y$ обозначим (неприводимое) представление алгебры $A$ такое, что $\wt\Y\circ \vf=\Y$,
см. предложение \ref{int krs}$(ii)$.

Следующее рассуждение является видоизменением идеи из
доказательства пункта $(ii)$ предложения \ref{rep trace}.

Пусть $\X_1,\ld,\X_s$ --- полный набор представителей классов эквивалентности неприводимых
$FG$-представ\-лений. 
Из предложения \ref{int krs}$(ii)$ следует, что
$\wt\X_1,\ld,\wt\X_s$ --- полный набор неприводимых
представлений алгебры $A$, которая полупроста ввиду предложения
\ref{fact ssm}. В силу теоремы
\ref{thm wedd}$(ii)$ имеет место разложение $A=\W_M(A^\circ)\oplus \Ann(M)$, где $M$ ---
неприводимый $A$-модуль, соответствующий представлению $\wt\X$. Запишем $a\vf=b\vf+c\vf$,
для некоторых элементов $b,c\in FG$ таких, что $b\vf\in
\W_M(A^\circ)$ и $c\vf\in\Ann(M)$. Тогда имеем
$$\X(b)=\wt\X(b\vf)=\wt\X(a\vf)=\X(a)\qquad\mbox{и}\qquad\Y(b)=\wt\Y(b\vf)=0$$ для всех неприводимых
представлений $\Y$ алгебры $FG$, не эквивалентных $\X$.
\end{proof}

\begin{cor} \label{cor ext com} Пусть $E\ge F$ --- расширение полей, $\X$ и $\Y$ ---
неприводимые $F$-представления группы $G$ такие, что $\X^E$ и $\Y^E$ имеют общую неприводимую компоненту.
Тогда $\X\cong\Y$.
\end{cor}
\begin{proof}
Предположим, что $\X\ncong\Y$. По предложению \ref{irr ab0} существует элемент
$b\in FG$ такой, что $\X(b)=\X(1)=\I_{\deg \X}$ --- единичная матрица и $\Y(b)=\O_{\deg \Y}$
--- нулевая матрица. Тогда $\X^E(b)=\I_{\deg \X}$ и $\Y^E(b)=\O_{\deg \Y}$, и поэтому для любой общей неприводимой
компоненты $\ZC$ представлений $\X^E$ и $\Y^E$ должно быть выполнено, с одной стороны,
$\ZC(b)=\I_{\deg \ZC}$, а с другой стороны $\ZC(b)=\O_{\deg \ZC}$. Противоречие.
\end{proof}

\mysubsection \label{zam pp}  Отметим, что из этого следствия вытекает единственность
с точностью до эквивалентности неприводимого представления $\X$, существование которого утверждается в предложении
\ref{ext npr}.

\begin{cor} \label{lin ind fch} Пусть $F$ --- поле разложения группы $G$. Тогда характеры из $\irr_F(G)$ являются
линейно независимыми над $F$ (в частности, ненулевыми и попарно различными).
\end{cor}
\begin{proof} Пусть $\X_1,\ld,\X_s$
--- набор представителей всех классов эквивалентности неприводимых $F$-представлений группы $G$ и пусть $\x_i$ ---
характер представления $\X_i$. Продолжим $\X_i$ и $\x_i$ на всю групповую алгебру $FG$. Тогда $\X_i(FG)$ --- полная матричная алгебра ввиду \ref{abs npr}$(iv)$.
Значит, существует $a_i\in FG$ такой, что $\x_i(a_i)=1$. По предложению \ref{irr ab0} можем считать, что $\x_j(a_i)=0$ при $j\ne i$. Отсюда следует требуемое.
\end{proof}

\begin{opr} \label{moz zap} Пусть $E\ge F$ --- расширение полей и $\X$ --- $E$-представление степени $n$ группы $G$. Будем говорить, что $\X$
\glsadd{iRepGrZapF}\mem{записано над} $F$, если $\X(g)\in \MM_n(F)$ для всех $g\in G$. Также будем говорить, что $\X$ \mem{может
быть записано над} $F$, если оно эквивалентно некоторому $E$-представлению, записанному над $F$.
\end{opr}

\textuprn{Показать, что неприводимое $\CC$-представление группы кватернионов $Q_8$
степени~$2$ не может быть записано над $\RR$.}

\begin{pre} \label{ext abs} Пусть $F$ --- поле разложения группы $G$ и $\X_1,\ld,\X_s$
--- набор представителей всех классов эквивалентности неприводимых $F$-представлений группы $G$.
Если $E \ge F$ --- расширение полей, то $\X_1^E,\ld,\X_s^E$ --- набор представителей всех неприводимых
$E$-представлений группы $G$. В частности, любое неприводимое $E$-представление может быть записано над полем
$F$.
\end{pre}
\begin{proof} Из абсолютной неприводимости представлений $\X_i$ следует неприводимость
представлений $\X_i^E$, которые попарно неэквивалентны по следствию \ref{cor ext
com}. Пусть $\Y$ --- некоторое неприводимое $E$-представление. Из предложения \ref{ext npr} и
неприводимости представлений $\X_i^E$ вытекает, что $\Y\cong\X_i^E$ для некоторого $i$.
\end{proof}

\begin{cor} \label{dec ext}
Пусть $F$ --- поле разложения группы $G$ и $E\ge F$ --- расширение полей. Тогда $E$ также является
полем разложения группы $G$.
\end{cor}

\begin{proof} Пусть $\X$ --- неприводимое $E$-представление группы $G$. Из \ref{ext abs} следует,
что $\X\cong \Y^E$ для некоторого неприводимого $F$-представления $\Y$. Но $\Y$ абсолютно неприводимо,
и, значит, $\X$ тоже в силу \ref{abs npr}$(ii)$.
\end{proof}

\begin{cor} \label{cor us} В предложении {\rm \ref{irr basis}} вместо алгебраической замкнутости поля $F$
достаточно потребовать выполнения более слабого условия~--- чтобы $F$ было полем разложения
группы $G$.
\end{cor}

\begin{proof} Это следует из \ref{ext abs}, поскольку любое поле может быть вложено в своё алгебраическое замыкание.\end{proof}

В следующем разделе будут обсуждаться представления групп над полем  $\CC$ комплексных чисел.
Отметим здесь одно полезное следствие, касающееся таких представлений.

\begin{cor} \label{cor c q} Пусть $\ov{\QQ}$ --- алгебраическое замыкание поля $\QQ$
рациональных чисел $($т.\,е. поле всех алгебраических чисел, см. приложение \ref{ap alg}$\,)$.
Любое $\CC$-представление группы $G$ может быть записано над полем $\ov{\QQ}$.
\end{cor}

\textupl{cor c q prf}{Доказать следствие \ref{cor c q}.}

\chapter{Обыкновенные представления}

\section{Обыкновенные характеры}

В этом разделе мы остановимся более подробно на классическом частном случае --- когда основное
поле $F$ является полем $\CC$ комплексных чисел. В этом случае часто говорят об
\glsadd{iChrOrd}\glsadd{iRepOrd}\mem{обыкновенных} представлениях и характерах групп.
Множество $\ch_\CC(G)$ всех обыкновенных характеров группы $G$ будем для краткости обозначать~\gls{chGr}.

\begin{opr} \label{kerh} Пусть $G$ --- конечная группа,
$\x\in \ch(G)$ и $\X$ --- представление
с характером $\x$. \glsadd{iDegChrOrd}\mem{Степенью} \gls{degx} и \glsadd{iKerChrOrd}\mem{ядром} \gls{kerx}  характера $\x$  называются
$\deg\X$ и $\ker \X$, соответственно.
Характер $\x$ называется \glsadd{iChrOrdIrr}\mem{неприводимым}
  (\glsadd{iChrOrdFth}\mem{точным}), если таковым является представление $\X$.
\end{opr}

Определение \ref{kerh} корректно, поскольку из \ref{rep trace}$(iii)$ вытекает следующее утверждение
\begin{pre} \label{ob vz} Два
обыкновенных представления группы эквивалентны тогда и только тогда,
когда их характеры совпадают.
\end{pre}

\textzam{Отметим, что для произвольных $F$-характеров, аналогично определяемые понятия степени,
ядра, неприводимости, вообще говоря, некорректны,
поскольку могут существовать неэквивалентные $F$-представления с одинаковыми $F$-характерами.}

Ясно, что для любого $\x\in \ch(G)$ справедливо равенство
$$
\deg \x=\x(1).
$$

Обозначим $\gls{IrrlGr}=\irr_{\CC}(G)$. Каждому характеру $\x\in \irr(G)$ соответствует единственный с точностью до
изоморфизма неприводимый $\CC G$-модуль $\gls{Mchi}\in \M(\CC G)$.
Выберем некоторый базис в каждом модуле $M_\x$ и обозначим через \gls{Xchi} соответствующее представление.

Всякий обыкновенный характер $\th$ можно представить в виде
$$
\th=\sum_{\x\in \irr(G)}n_\x \x,\myeqno\label{x dec}
$$
см. \ref{c dec}, где $n_\x$ --- однозначно определённые неотрицательные целые числа. Характер $\x\in\irr(G)$,
для которого $n_\x>0$, мы будем называть
\glsadd{iCompIrrOrdChr}\mem{неприводимой компонентой} характера $\th$, а соответствующее число $n_\x$ ---
\glsadd{iMltIrrComp}\mem{кратностью неприводимой компоненты} $\x$. Таким образом, $\x\in\irr(G)$ является неприводимой компонентой
характера $\th$ тогда и только тогда, когда $\X_\x$ --- неприводимая компонента обыкновенного
представления с характером $\th$. Далее в \ref{kr vh} мы укажем способ нахождения кратностей неприводимых
компонент, использующий внутреннее произведение на алгебре классовых функций.

Для частного случая, когда в формуле \ref{x dec} характер $\th$ совпадает с регулярным $\CC$-характером $\gls{rh}=\r_{\mbox{}_\CC}$,
в силу \ref{irr basis}$(i)$ имеем $n_\x=\x(1)$ для любого $\x\in\irr(G)$, откуда получаем

\begin{pre} \label{fbern} Справедливы равенства
\begin{list}{{\rm(}{\it\roman{enumi}\/}{\rm)}}
{\usecounter{enumi}\setlength{\parsep}{2pt}\setlength{\topsep}{5pt}\setlength{\labelwidth}{23pt}}
\item $\r=\sum_{\x\in \irr(G)}\x(1)\x$;
\item $|G|=\sum_{\x\in \irr(G)} \x(1)^2$\quad (\mem{формула Бернсайда}\glsadd{iFormBurn}).
\end{list}
\end{pre}

Из \ref{cor z fg}$(iii)$ следует, что для конечной группы $G$  число неприводимых обыкновенных характеров
выражается только в групповых терминах и совпадает с количеством её классов сопряжённости. Другими словами
справедливо равенство $$|\irr(G)|=|\K(G)|.\myeqno\label{ir eq kl}$$

Подмножество $\lin_\CC(G)\se \irr(G)$ обыкновенных линейных характеров группы $G$ будем для краткости обозначать
$\gls{LinG}$.

\begin{pre} \label{ab npr}
Группа $G$ абелева тогда и только тогда, когда $\irr(G)=\lin(G)$.
\end{pre}
\textupl{ab npr prf}{Доказать предложение \ref{ab npr}}

\textuprn{Найти степени неприводимых обыкновенных характеров группы $S_3$.}

Для группы $G$ \glsadd{iTblOrdChr}\mem{таблицей $($обыкновенных$\,)$ характеров} \gls{XiG} называется квадратная матрица, строки которой индексированы неприводимыми
характерами $\x\in\irr(G)$, а столбцы --- классами (или представителями классов)
сопряжённости $K\in \K(G)$. При этом
элемент $\x$-й строки и $K$-го столбца таблицы $\XX(G)$ равен $\x(x_{\mbox{}_K})$, где
$x_{\mbox{}_K}$ обозначает представитель класса $K$.
Например, ниже приведена\footnote{См. \ref{ts3d}} таблица характеров группы $S_3$.
$$
  \begin{array}{c|ccc}
            \XX(S_3)              &\ 1\ & (12)           & (123)       \\
     \hline
    \x_1^{\vphantom{A^A}} &1      & \phantom{-}1  & \phantom{-}1\\
    \x_2                  &1      & -1            & \phantom{-}1\\
    \x_3                  &2      & \phantom{-}0  & -1
  \end{array}\myeqno\label{tchs3}
$$

Строка матрицы $\XX(G)$, индексированная главным характером $1_G$ обычно ставится на первое место.
Она находится по виду таблицы характеров однозначно, как единственная строка, все элементы которой равны~$1$. Также
на первом месте обычно ставится столбец, индексированный классом $1^G=\{1\}$. Он
находится однозначно, как единственный столбец, каждый элемент которого
является положительным вещественным и максимальным по модулю среди элементов своей строки.\footnote{См. \ref{perv st}}

\textuprn{\label{x cyc}Пусть $C=\la c\mid c^n=1\ra$ --- циклическая группа порядка $n$.
\begin{list}{{\rm(}{\it\roman{enumi}\/}{\rm)}}
{\usecounter{enumi}\setlength{\parsep}{2pt}\setlength{\topsep}{5pt}\setlength{\labelwidth}{23pt}}
\item Показать, что таблица $\XX(C)$ имеет вид
$$
\begin{array}{c|cllll}
            \XX(C)              & 1      & c  & c^2 & \ld & c^{n-1}\\
     \hline
    \x_1^{\vphantom{A^A}} & 1      & 1 & 1 & \ld  & 1\\
    \x_2                  & 1      & \z & \z^2 & \ld  & \z^{(n-1)} \\
    \x_3                  & 1      & \z^2 & \z^4 & \ld  & \z^{2(n-1)} \\
\ld & \ld   & \ld  & \ld  & \ld  & \ld \\
    \x_n                  & 1  & \z^{n-1} & \z^{(n-1)2} & \ld  & \z^{(n-1)^2}
  \end{array}
$$
где $\z=e^{2\pi i/n}$.
\item Доказать изоморфизм\footnote{Cр. \ref{lin ab}} групп $\lin(C)\cong C$.
\end{list}}

\textzam{Отметим, что задача вычисления таблицы характеров произвольной группы является достаточно сложной.
Для этого в каждом отдельном случае используются как правило свои специальные методы.
Особое значение имеет нахождение неприводимых характеров конечных простых групп. Для многих простых
и близких к ним групп <<небольшого>> порядка (в частности, для всех спорадических групп) таблицы характеров
и другая важная информация могут быть найдены в Атласе \cite{at}.
Таблицы характеров некоторых групп приведены также в приложении \ref{pril tch}.}

В силу \ref{rep trace}$(iii)$ получаем

\begin{cor} \label{tch nv} Для любой
группы $G$ матрица $\XX(G)$ является невырожденной.
\end{cor}

Для комплексного числа $\a\in \CC$ пусть $\ov\a$ обозначает число, комплексно-сопряжённое с
$\a$. Пусть $\ov{\ZZ}$ --- кольцо целых алгебраических чисел.\footnote{Определение и необходимые сведения о целых алгебраических числах приведены в приложении \ref{ap
alg}.} Примерами целых алгебраических чисел являются корни из $1$, т.\,е. комплексные корни
многочленов $x^m-1$ для различных натуральных $m$.

Через \gls{Qm} обозначим \glsadd{iFldCcl}\mem{$m$-е круговое поле}, т.\,е. подполе $\QQ(\z)$ поля $\CC$, где $\z$ ---
\glsadd{iPrmRt}\mem{примитивный корень степени $m$ из $1$}, т.\,е. элемент порядка $m$ из  группы $\CC^\times$.

\glsadd{iExpGr}\mem{Экспонентой} \gls{expG} группы $G$ называется наименьшее общее кратное порядков её элементов.

\begin{pre} \label{har prop} Пусть
$G$ --- группа экспоненты $m$, $g\in G$ и $k=|g|$. Пусть
$\X$ --- обыкновенное представление группы $G$ степени $n$ с  характером $\x$.
\begin{list}{{\rm(}{\it\roman{enumi}\/}{\rm)}}
{\usecounter{enumi}\setlength{\parsep}{2pt}\setlength{\topsep}{5pt}\setlength{\labelwidth}{23pt}}
\item Матрица $\X(g)$ подобна диагональной матрице $\diag(\z_1,\ld,\z_n)$ для некоторых $\z_i\in \CC^\times$.
\item $\z_i^k=1$. B частности, $|\z_i|=1$, \ $\ov{\z_i}=\z_i^{-1}$ и $\z_i\in \QQ_k$.
\item $\x(g)=\sum \z_i$. В частности, $|\x(g)|\le n$, $\x(g)\in \QQ_k\se \QQ_m$, и
$\x(g)\in \ov{\ZZ}$.
\item $|\x(g)|=n$ тогда и только тогда, когда $\X(g)$ --- скалярная матрица.
\item $\x(g)=n$ тогда и только тогда, когда $g\in \ker\x$.
\item $\x(g^{-1})=\ov{\x(g)}$.
\end{list}
\end{pre}
\textupl{har prop prf}{Доказать предложение \ref{har prop}}

\begin{pre} \label{reg har} Пусть $\r$ --- регулярный $\CC$-характер группы $G$. Тогда
\begin{list}{{\rm(}{\it\roman{enumi}\/}{\rm)}}
{\usecounter{enumi}\setlength{\parsep}{2pt}\setlength{\topsep}{5pt}\setlength{\labelwidth}{23pt}}
\item $\ker\r=1$;
\item для любого $g\in G$
$$
\r (g)=\left\{\ba{ll}
|G|,& g=1;\\
0,& g\ne 1;
\ea\right.
$$
\end{list}
\end{pre}

\textupl{reg har prf}{Доказать предложение \ref{reg har}.}

\begin{pre} \label{kpr} Пусть $\x\in \ch(G)$. Тогда
\begin{list}{{\rm(}{\it\roman{enumi}\/}{\rm)}}
{\usecounter{enumi}\setlength{\parsep}{2pt}\setlength{\topsep}{5pt}\setlength{\labelwidth}{23pt}}
\item $\ker\x=\{g\in G\mid\x(g)=\x(1)\}$.
\item Если $\x=\sum_{\t\in\irr(G)} n_\t\t$, где $n_\t$ --- неотрицательные целые числа, то
$$\ker\x=\bigcap_{\substack{\t\in\irr(G),\\n_\t>0}}\ker\t.$$
\item Имеет место равенство
$$\bigcap_{\t\in\irr(G)}\ker\t=1.$$
\item Если $N\nor G$, то
$$
N=\bigcap_{\substack{\t\in \irr(G),\\N\le   \ker\t}} \ker\t.
$$
\end{list}
\end{pre}
\textupl{kpr prf}{Доказать предложение \ref{kpr}.}

\textuprn{\label{perv st}Показать, что столбец таблицы характеров группы $G$, индексированный классом $1^G$,
находится однозначно как единственный столбец, каждый элемент которого
является положительным вещественным и максимальным по модулю среди элементов своей строки.\footnote{Ср. \ref{perv st1}.}}


\begin{pre} \label{ker indh} Пусть $H\le G$ и $\x\in \ch(H)$.
Тогда
$$
\ker(\x^G)=\bigcap_{x\in G}(\ker \x)^x.
$$
\end{pre}
\begin{proof} В силу \ref{kpr}$(i)$ и \ref{ind kl}$(iii)$ включение $g\in \ker(\x^G)$ эквивалентно равенству
$$
\sum_{x\in G}\x^\circ(xgx^{-1})=\sum_{x\in G}\x(1),
$$
где $\x^\circ$ совпадает с $\x$ на $H$ и тождественно равно нулю вне $H$. Из \ref{har prop}$(iii)$ следует, что
$|\x^\circ(xgx^{-1})|\le \x(1)$. Значит, $g\in \ker(\x^G)$ тогда и только тогда, когда $\x^\circ(xgx^{-1})= \x(1)$
для всех $x\in G$. Но это означает, что $xgx^{-1}\in \ker\x$, т.\,е. $g\in (\ker\x)^x$ для всех $x\in G$,
что и требовалось доказать. \end{proof}

\begin{pre} \label{ch gn} Пусть
$N\nor G$, $\wt G=G/N$ и $\vf:G\to \wt G$ --- естественный гомоморфизм.
\begin{list}{{\rm(}{\it\roman{enumi}\/}{\rm)}}
{\usecounter{enumi}\setlength{\parsep}{2pt}\setlength{\topsep}{5pt}\setlength{\labelwidth}{23pt}}
\item Отображение $\x\mapsto \wt\x$ такое, что $\x=\wt\x\circ\vf$, осуществляет взаимно однозначное соответствие между множествами
$\{\x\in \ch(G)\mid N\se\ker\x\}$ и $\ch(\wt G)$.
\item Имеем $\deg\x=\deg\wt\x$ и $\x\in \irr(G)$ тогда и только тогда, когда $\wt\x\in \irr(\wt G)$.
\item Характер $\wt \x$ будет точным тогда и только тогда, когда $N=\ker\x$.
\item Пусть $N\le H \le G$, $\wt H=H\vf$ и $\t$ --- характер группы $H$ такой, что $N\se\ker\t$.
Обозначим через $\wt\t$ такой характер группы $\wt H$, для которого $\t=\wt\t\circ\vf_{H^{\vphantom{A^a}}}$
и существование которого следует из~$(i)$. Тогда $N\se \ker \t^G$ и справедливо
равенство $\wt\t^{\,\wt G}=\wt{(\t^G)}$.
\end{list}
\end{pre}

\textupl{ch gn prf}{Доказать предложение \ref{ch gn}.}

Таблица характеров $\XX(G)$, вообще говоря, не определяет $G$ однозначно с точностью до изоморфизма.
Например, обе таблицы $\XX(D_{8})$ и $\XX(Q_8)$ можно отождествить (переставив, если необходимо, строки и столбцы)
со следующей таблицей.
$$
  \begin{array}{c|rrrrr}
                          & K_1 & K_2 & K_3 & K_4  & K_5 \\
     \hline
    \x_1^{\vphantom{A^A}} & 1\ & 1 \  & 1 \ & 1  \ & 1 \ \\
    \x_2                  & 1\ & 1 \  & -1\  & -1\  & 1\ \\
    \x_3                  & 1\ & 1 \  & 1 \ & -1 \ & -1\ \\
    \x_4                  & 1\ & 1 \  & -1\  & 1 \ & -1\ \\
    \x_5                  & 2\ & -2\  & 0 \ & 0  \ & 0 \
  \end{array}
$$
Тем не менее, в таблице характеров группы содержится много информации о её строении. Мы уже
видели в \ref{fbern}$(ii)$ и \ref{ab npr},
что по первому столбцу таблицы $\XX(G)$ можно найти порядок $|G|$ и выяснить, является ли группа $G$ абелевой.
Из предложения \ref{kpr}$(iv)$ следует, что по $\XX(G)$ также можно определить\footnote{Любая нормальная подгруппа $N\nor G$
является объединением классов сопряжённости, и под словом <<определить>> понимается, что
можно указать столбцы таблицы $\XX(G)$, индексированные классами, содержащимся в $N$.}
 все нормальные подгруппы группы $G$ и все отношения включения между ними.

Приведём алгоритм, с помощью которого
по таблице характеров группы можно найти таблицы характеров всех её факторгрупп. Пусть $N\nor G$.
Оставим в таблице $\XX(G)$ строки, индексированные только теми $\x\in\irr(G)$,
для которых $N\se \ker \x$. Из полученной прямоугольной таблицы удалим повторяющиеся столбцы. В результате получим квадратную матрицу, совпадающую с $\XX(G/N)$.

\textuprn{Доказать справедливость приведённого алгоритма нахождения таблицы характеров факторгруппы.
\uk{Воспользоваться \ref{ch gn} и \ref{sopr char}.}}

Напомним, что группа $G$ называется
\glsadd{iGrSol}\mem{разрешимой}, если ряд последовательных коммутантов
$$
G\ge G'\ge G''\ge\ld \myeqno{\label{der ser}}
$$
через конечное число доходит до единичной подгруппы. Другими словами, $G$ разрешима, если существует
число $n$ такое, что $G^{(n)}=1$, где $G^{(0)}=G$, $G^{(i)}=(G^{(i-1)})'$,  $i\ge 1$.
Подгруппа \gls{Glir} называется \glsadd{iCmtIth}\mem{$i$-м коммутантом} группы~$G$.

\textuprn{Доказать, что по таблице $\XX(G)$ можно выяснить, является ли группа $G$ разрешимой.
\uk{Конечная группа разрешима, тогда и только тогда, когда она обладает нормальным
рядом\footnote{\mem{Нормальным рядом}\glsadd{iSerNor} группы $G$ называется ряд подгрупп
$$
G=G_0\ge G_1\ge \ldots \ge G_n =1
$$
такой, что  $G_i\nor G$, $i=1,\ldots,n$.}
с факторами примарных порядков. При этом, как мы увидим далее в \ref{chr ccs}, по $\XX(G)$ можно также определить размеры классов сопряжённости группы $G$.}}

Следующее утверждение показывает, что по таблице характеров можно определить коммутант группы.

\begin{pre} \label{com cht} Пусть
$G$ --- группа. Тогда
\begin{list}{{\rm(}{\it\roman{enumi}\/}{\rm)}}
{\usecounter{enumi}\setlength{\parsep}{2pt}\setlength{\topsep}{5pt}\setlength{\labelwidth}{23pt}}
\item справедливо равенство
$$
G'=\bigcap_{\l\in \lin(G)} \ker\l\text{;}
$$
\item $\lin(G)\cong \lin(G/G')$;
\item $|\lin(G)|=|G/G'|$.\footnote{Cр. \ref{lin ab}}
\end{list}
\end{pre}
\textupl{com cht prf}{Доказать предложение \ref{com cht}.}


По таблице характеров можно определить центр группы. Для доказательства этого факта введём одно обозначение.
Пусть $\x\in \ch(G)$. Положим
$$
\gls{Zchi}=\{g\in G\bigm|\  |\x(g)|=\x(1)\}.
$$
В силу \ref{har prop}$(iv)$ можно эквивалентно определить
$$
\Z(\x)=\{g\in G\bigm|\  \X(g)\ \mbox{--- скалярная матрица}\},\myeqno\label{zh eq}
$$
где $\X$ --- обыкновенное представление с характером $\x$.

\begin{pre} \label{zh pr} Пусть
$\x\in \ch(G)$. Справедливы следующие утверждения.
\begin{list}{{\rm(}{\it\roman{enumi}\/}{\rm)}}
{\usecounter{enumi}\setlength{\parsep}{2pt}\setlength{\topsep}{5pt}\setlength{\labelwidth}{23pt}}
\item $\Z(\x)$ --- нормальная подгруппа в $G$,  содержащая $\ker\x$.
\item Группа $\Z(\x)/\ker\x$ циклическая. В частности, если $\x$ точный, то $\Z(\x)$ циклическая.
\item $\Z(\x)/\ker\x\le \Z(G/\ker\x)$.
\item Если $\x\in \irr(G)$, то $\Z(\x)/\ker\x=\Z(G/\ker\x)$. В частности, если $\x\in \irr(G)$ точный, то $\Z(\x)=\Z(G)$.
\end{list}
\end{pre}
\begin{proof}
Пусть $\X$ --- обыкновенное представление с характером $\x$.
Из \ref{zh eq} следует, что $\Z(\x)$
совпадает с полным прообразом относительно $\X$ пересечения $\X(G)$ с множеством скалярных матриц.
Поскольку это пересечение является центральной подгруппой в $\X(G)$,
в силу изоморфизма $\X(G)\cong G/\ker\x$ получаем $(i)$ и $(iii)$.
С другой стороны, это пересечение является
конечной подгруппой мультипликативной группы $\CC^\times$, и значит, будет циклическим в силу \ref{cyc fin}. Отсюда следует $(ii)$.

Предположим,  что $\x\in \irr(G)$ и для некоторого $g\in G$  элемент $g(\ker\x)$ лежит в $\Z(G/\ker\x)$.
Тогда матрица $\X(g)$  лежит в $\Z(\X(G))$.
В силу абсолютной неприводимости $\X$ из предложения \ref{abs npr}$(iii)$
следует, что $\X(g)$ --- скалярная матрица, т.\,е. $g\in \Z(\x)$. Отсюда получаем $(iv)$.
\end{proof}

\begin{cor} \label{zg ct} Пусть
$G$ --- конечная группа. Тогда
$$
\Z(G)=\bigcap_{\x\in \irr(G)}\Z(\x).
$$
\end{cor}
\begin{proof} Если $g\in \Z(G)$, то для любого $\x\in \irr(G)$ имеем $g(\ker\x)\in \Z(G/\ker\x)$.
Из предложения \ref{zh pr}$(iv)$ следует, что $g\in \Z(\x)$.

Обратно, пусть $g\in \Z(\x)$ для любого $\x\in \irr(G)$. Из \ref{zh pr}$(iii)$ следует,
что $g(\ker\x)$ коммутирует с $x(\ker\x)$ для любого $x\in G$. Но тогда $[g,x]$ лежит
в $\ker\x$ для всех $\x\in \irr(G)$. Из \ref{kpr}$(iii)$ следует, что $[g,x]=1$.
В силу произвольности $x$ получаем $g\in \Z(G)$, что и требовалось доказать.
\end{proof}

Из следствия \ref{zg ct} легко видеть, что центр $\Z(G)$ можно определить по таблице характеров $\XX(G)$.

Напомним, что конечная группа $G$ называется \glsadd{iGrNlp}\mem{нильпотентной}, если ряд её подгрупп
$$
\Z_0(G)\le \Z_1(G)\le \ldots
$$
через конечное число шагов достигает всей группы $G$,
где $\Z_0(G)=1$, а \gls{ZilGr}, $i\ge 1$, определяются равенством $\Z_i(G)/\Z_{i-1}(G)=\Z(G/\Z_{i-1}(G))$.
Для нильпотентной группы $G$ наименьшее число $n$, для которого $\Z_n(G)=G$ называется
\glsadd{iNlpCls}\mem{ступенью нильпотентности} группы $G$.

\textuprn{Доказать, что по таблице $\XX(G)$ можно выяснить, является ли группа $G$ нильпотентной,
а для нильпотентной группы $G$ определить её ступень нильпотентности.}

Далее в \ref{com opr} мы покажем, что по таблице характеров можно определить, является
ли данный элемент $g\in G$ коммутатором, т.\,е. элементом вида $a^{-1}b^{-1}ab$ для некоторых $a,b\in G$.

%

Пусть $G$ --- группа  экспоненты $m$, $\x$ --- её обыкновенный характер и
$\s\in\Gal(\QQ_m,\QQ)$. Обозначим через
\gls{chiov} и \gls{chisi} классовые функции, определённые равенствами
$$
\ov{\x}(g)=\ov{\x(g)}, \qquad \x^\s(g)=\x(g)^\s
$$
для любого $g\in G$.

\begin{pre} \label{gal conj} Пусть $G$ --- группа, $m=\exp G$ и $\x\in \ch(G)$. Тогда справедливы следующие утверждения.
\begin{list}{{\rm(}{\it\roman{enumi}\/}{\rm)}}
{\usecounter{enumi}\setlength{\parsep}{2pt}\setlength{\topsep}{5pt}\setlength{\labelwidth}{23pt}}
\item   $\ov{\x}\in\{\x^\t\,|\,\t\in\Gal(\QQ_m,\QQ)\}.$
\item  Если $\s\in \Gal(\QQ_m,\QQ)$, то $\x^\s\in \ch(G)$,
причём $\x^\s$ неприводим тогда и только тогда, когда $\x$ неприводим.
\item Пусть $\s\in \Gal(\QQ_m,\QQ)$. Тогда существует целое число $k$, взаимно простое с $m$, такое, что
\end{list}
\vspace{-2\topsep}
$$
\x^\s(g)=\x(g^k) \ \ \text{для любого}\ \  g\in G.\myeqno\label{xsg}
$$
\begin{list}{{\rm(}{\it\roman{enumi}\/}{\rm)}}
{\usecounter{enumi}\setlength{\parsep}{2pt}\setlength{\topsep}{5pt}\setlength{\labelwidth}{23pt}}
\item[] Обратно, пусть $k$ --- целое число, взаимно простое с $m$. Тогда
существует $\s\in \Gal(\QQ_m,\QQ)$ такое, что справедливо соотношение\footnote{Ср. \ref{har prop}$(vi)$} $\ref{xsg}$.
\end{list}
\end{pre}
\begin{proof}$(i)$ В силу \ref{har prop}$(iii)$ для любого $g\in G$ имеем $\x(g)\in \QQ_m$.
Поэтому ввиду \ref{ks f}$(ii)$ существует автоморфизм $\t\in \Gal(\QQ_m,\QQ)$ такой, что
$$\ov\x(g)=\ov{\x(g)}=\x(g)^\t=\x^\t(g).$$
для всех $g\in G$, т.\,е. $\ov\x=\x^\t$.

$(ii)$ Пусть $\X$ --- обыкновенное представление группы $G$ с характером $\x$. По \ref{cor c q}
существует 
$\ov{\QQ}$-представ\-ление 
$\Y$ такое, что  $\Y^\CC\cong\X$. Поэтому $\x$
будет характером $\ov{\QQ}$-представления $\Y$. В силу алгебраичности расширения $\QQ_m\ge \QQ$
из \ref{alg pro}$(ii)$ следует, что автоморфизм
$\s$ поля $\QQ_m$ продолжается до автоморфизма алгебраического замыкания $\ov{\QQ}$, который
также обозначим через $\s$. Пусть $g\in G$ и $\Y(g)=(\a_{ij})$, где $\a_{ij}\in \ov{\QQ}$.
Легко видеть, что отображение $\Y^\s:g\mapsto (\a_{ij}^\s)$ является $\ov{\QQ}$-представлением
группы $G$, и его характер совпадает с $\x^\s$. Поскольку $\s$ является автоморфизмом поля
$\ov{\QQ}$, представление $\Y^\s$ неприводимо тогда и только тогда, когда неприводимо
представление $\Y$. Таким образом, утверждения о неприводимости характеров $\x$ и $\x^\s$
равносильны в силу \ref{ext abs}.

\normalmarginpar

$(iii)$ Пусть $\z\in \CC$ --- примитивный корень степени $m$ из $1$. В силу \ref{ks f}$(iii)$
по данному $\s$ найдётся $k\in\ZZ$, взаимно простое с $m$, и обратно, для данного такого $k$
найдётся $\s$ такие, что справедливо равенство $\z^\s=\z^k$. Пусть $\x$
--- характер группы $G$ и $g\in G$. Рассмотрим обыкновенное представление $\X$ с
характером $\x$. Так как $|g|$ делит $m$, то из \ref{har prop} следует, что  матрица $\X(g)$ подобна
диагональной матрице $\diag(\z_1,\ld,\z_n)$, где $\z_i$
--- степени числа $\z$ и $n=\deg\x$. Тогда $\X(g^k)$ подобна $\diag(\z_1^k,\ld,\z_n^k)$ и,
значит,
$$
\x(g^k)=\tr \X(g^k)=\z_1^k+\ld+\z_n^k=(\z_1+\ld+\z_n)^\s=\x^\s(g).
$$
\end{proof}

В предыдущих обозначениях будем называть  $\ov{\x}$ и $\x^\s$ характерами,
\glsadd{iChrCompCnj}\mem{комплексно сопряжённым} и  \glsadd{iChrAlgCnj}\mem{алгебраически сопряжённым} с характером $\x$, соответственно.


\textuprn{Пусть $\X$ --- обыкновенное представление группы $G$ с характером $\x$. Показать, что $\ov{\x}$ является
характером контрагредиентного представления $\X^*$.}

\begin{pre}\label{alg spr}
Пусть  $G$ --- группа и  $g\in G$.
\begin{list}{{\rm(}{\it\roman{enumi}\/}{\rm)}}
{\usecounter{enumi}\setlength{\parsep}{2pt}\setlength{\topsep}{5pt}\setlength{\labelwidth}{23pt}}
\item Элементы $g$ и $g^{-1}$ сопряжены в $G$ тогда и только тогда, когда $\x(g)\in \RR$ для
всех $\x\in\irr(G)$.
\item Все элементы множества $\{\,g^k\bigm| k\in \ZZ, (k,|g|)=1\}$
сопряжены тогда и только тогда, когда $\x(g)\in \QQ$ для всех $\x\in\irr(G)$.
\end{list}
\end{pre}
\textupl{alg spr prf}{Доказать предложение \ref{alg spr}.}

\begin{cor} \label{sym int} Таблица характеров
симметрической группы $S_n$ целочисленна.
\end{cor}
\textupl{sym int prf}{Доказать следствие \ref{sym int}.}

\section{Соотношения ортогональности}

\mysubsection \label{cg} Пусть $\x\in\irr(G)$. Центральный идемпотент $e_{\mbox{}_{M_\x}}$ групповой алгебры
$\CC G$, соответствующий неприводимому модулю $M_\x$, будем обозначать \glsadd{iIdmpCenCG} через \gls{esx}. Как мы отмечали в
\ref{cid}, единица алгебры $\CC G$ представима в виде суммы
$$
1=\sum_{\x\in \irr(G)}e_\x.
$$
Кроме того, $e_\x$ --- единица простой алгебры $\W_{M_{\x}}(\CC G^\circ)$ и для любого $\x\in \irr(G)$ справедливо равенство
$$
\W_{M_{\x}}(\CC G^\circ)=e_\x\CC G.
$$
При этом идемпотенты $e_\x$ образуют базис $\Z(\CC G)$ в силу \ref{cor z ss}$(iv)$.


Нам потребуется явное выражение идемпотентов
$e_\x$ через естественный базис алгебры $\CC G$, состоящий из групповых элементов, а также
базис центра $\Z(\CC G)$, состоящий из классовых сумм.

\begin{pre} \label{id dec} Имеет место разложение
$$e_\x=\frac{\x(1)}{|G|}\sum_{g\in G}\x(g^{-1})g=\frac{\x(1)}{|G|}\sum_{K\in \K(G)}\x(x_{\mbox{}_K}^{-1})\wh{K}.$$
\end{pre}

\begin{proof} Достаточно доказать лишь первое равенство.
Пусть $e_{\x}=\sum_{g\in G}\a_g g$ для подходящих коэффициентов $\a_g\in \CC$. Из \ref{reg har}$(ii)$
следует, что для любого $g\in G$
$$
\r(e_{\x}g^{-1})=\sum_{h\in G}\a_h\r(hg^{-1})=\a_g|G|.
$$
Тогда
\ref{fbern}$(i)$ влечёт
$$
\a_g|G|=\sum_{\th\in \irr(G)}\th(1)\th(e_{\x}g^{-1}).
$$

Отметим, что ввиду \ref{thm wedd}$(ii)$ однородная компонента $\W_{M_{\x}}(\CC G)$ алгебры $\CC G$ аннулирует
любой модуль $M_\th$ при $\th\ne \x$ и поэтому $\X_\th(e_\x)=\O$, а в силу \ref{cg} также
$\X_\th(e_\th)=\X_\th(1)=\I$. Поэтому
$$
\X_\th(e_{\x}g^{-1})=\X_\th(e_{\x})\X_\th(g^{-1})=
\left\{\ba{ll}
\O,& \th\ne \x;\\
\X_\x(g^{-1}),& \th=\x,
\ea\right.
$$
и, значит, $\th(e_{\x}g^{-1})=\d_{\th,\x}\x(g^{-1})$. Следовательно $\a_g |G|=\x(1)\x(g^{-1})$, что
и требовалось показать.
\end{proof}

Условие ортогональности идемпотентов
$$
e_\x e_\th=\d_{\x,\th}\,e_\x,\myeqno\label{id ort}
$$
где $\x,\th\in\irr(G)$, см. \ref{cid}, позволяет вывести важное соотношение между
неприводимыми характерами.

\begin{pre}[Обобщённое соотношение ортогональности] \label{ob ort}\glsadd{iRelOrthGen}
Для произвольных $h\in G$, $\x,\th\in \irr(G)$
имеет место равенство
$$
\frac{1}{|G|}\sum_{g\in G}\x(gh)\th(g^{-1})=\d_{\x,\th}\frac{\x(h)}{\x(1).}
$$
\end{pre}
\begin{proof}
Подставим в соотношение \ref{id ort} выражение для идемпотентов $e_\x$ и $e_\th$
из предложения \ref{id dec} и сравним коэффициенты при групповых элементах в обеих частях равенства.

Зафиксируем $h\in G$. Коэффициент при $h$ справа равен
$$\d_{\x,\th}\frac{\x(1)}{|G|}\x(h^{-1}),$$
а слева ---
$$
\frac{\x(1)\th(1)}{|G|^2}\Big(\sum_{\substack{f,g\in G,\\ fg=h}}\x(f^{-1})\th(g^{-1})\Big)=
\frac{\x(1)\th(1)}{|G|^2}\Big(\sum_{g\in G}\x(gh^{-1})\th(g^{-1})\Big).
$$
Требуемое равенство получается приравниванием этих коэффициентов и заменой $h$ на $h^{-1}$.
\end{proof}

Подстановка $h=1$ в \ref{ob ort} даёт

\begin{cor}[Первое соотношение ортогональности] \label{cor perv ort}\glsadd{iRelOrthFrs}
Для произвольных $\x,\th\in \irr(G)$ имеет место равенство
$$
\frac{1}{|G|}\sum_{g\in G}\x(g)\th(g^{-1})=\d_{\x,\th}.
$$
\end{cor}

C учётом \ref{har prop}$(vi)$, первое соотношение ортогональности можно эквивалентно переписать
в виде

$$
\frac{1}{|G|}\sum_{g\in G}\x(g)\ov{\th(g)}=\d_{\x,\th}. \myeqno\label{perv ort2}
$$

\textuprn{\label{ts3d}Показать, что таблица характеров группы $S_3$ такая, как в \ref{tchs3}.}

Обозначим $\gls{cfG}=\cf_\CC(G)$. Для произвольных классовых
функций $\vf,\psi\in \cf(G)$ введём \glsadd{iSclPrdCompClsFnc}\mem{скалярное произведение}
$$
\gls{lvfpsrsG}=\frac{1}{|G|}\sum_{g\in G}\vf(g)\ov{\psi(g)}.
$$

\textuprn{Проверить, что для любых $\vf,\psi,\tau\in \cf(G)$, $\a\in \CC$ выполнены свойства
\begin{list}{{\rm(}{\it\roman{enumi}\/}{\rm)}}
{\usecounter{enumi}\setlength{\parsep}{2pt}\setlength{\topsep}{5pt}\setlength{\labelwidth}{23pt}}
\item $(\vf,\vf)_{\mbox{}_G}\in\RR$ и если $\vf\ne 0$, то $(\vf,\vf)_{\mbox{}_G}>0$;
\item $(\vf+\psi,\t)_{\mbox{}_G}=(\vf,\t)_{\mbox{}_G}+(\psi,\t)_{\mbox{}_G}$;
\item $(\a\vf,\psi)_{\mbox{}_G}=\a(\vf,\psi)_{\mbox{}_G}$;
\item $(\vf,\psi)_{\mbox{}_G}=\ov{(\psi,\vf)}_{\mbox{}_G}$.
\end{list}
Другими словами, $(\cdot,\cdot)_{\mbox{}_G}$ является обычным скалярным произведением, относительно которого
$\cf(G)$ становится унитарным пространством.}

Мы знаем из \ref{irr basis}$(iii)$, что $\irr(G)$ является базисом пространства $\cf(G)$. В
силу \ref{perv ort2}  для любых $\x,\th\in \irr(G)$ имеет место равенство
$$
(\x,\th)_{\mbox{}_G}=\d_{\x,\th},
$$
т.\,е. базис $\irr(G)$ ортонормированный. Отсюда получается простой метод нахождения
коэффициентов в разложении произвольной классовой функции $\vf\in \cf(G)$ по базису $\irr(G)$.
Если $(\vf,\x)_{\mbox{}_G}=\a_\x$, где $\x\in\irr(G)$, то $\vf=\sum_{\x\in \irr(G)}\a_\x \x$.
Другими словами для любой $\vf\in \cf(G)$ справедливо равенство
$$
\vf=\sum_{\x\in \irr(G)}(\vf,\x)_{\mbox{}_G}\x. \myeqno\label{kf raz}
$$
В частности, для кратности $n_\x$ неприводимой компоненты $\x$ обыкновенного характера $\th$ группы $G$
имеем соотношение $$
n_\x=(\th,\x)_{\mbox{}_G}. \myeqno\label{kr vh}
$$

\textuprn{Пусть $\th\in \ch(G)$  и $\x\in\irr(G)$.
Показать, что $(\th,\x)_{\mbox{}_G}\le\deg\th$.}

Используя введённое скалярное произведение, можно сформулировать простой критерий принадлежности
данной классовой функции множеству (неприводимых) характеров или кольцу обобщённых $\CC$-характеров $\gls{gchG}=\gch_\CC(G)$,
элементы
которого мы будем для краткости называть \glsadd{iChrGnrl}\mem{обобщёнными характерами} группы~$G$.

\begin{cor} \label{cor scl} Пусть $\vf,\psi\in \cf(G)$.
\begin{list}{{\rm(}{\it\roman{enumi}\/}{\rm)}}
{\usecounter{enumi}\setlength{\parsep}{2pt}\setlength{\topsep}{5pt}\setlength{\labelwidth}{23pt}}
\item $\vf\in \gch(G)$ тогда и только тогда, когда $(\vf,\x)_{\mbox{}_G}$
--- целое число для всех $\x\in\irr(G)$.
\item $\vf\in \ch(G)$ тогда и только тогда, когда $(\vf,\x)_{\mbox{}_G}$
--- неотрицательное целое число для всех $\x\in\irr(G)$.
\item Если $\vf,\psi\in \ch(G)$, то $(\vf,\psi)_{\mbox{}_G}=(\psi,\vf)_{\mbox{}_G}$
--- неотрицательное целое число.
\item Если $\vf\in \ch(G)$, то он неприводим тогда и только тогда, когда $(\vf,\vf)_{\mbox{}_G}=1$.
\end{list}
\end{cor}
\textupl{cor scl prf}{Доказать следствие \ref{cor scl}.}

Сопряжение сохраняет скалярное произведение классовых функций подгрупп и неприводимость их характеров как показывает
следующее утверждение.

\begin{pre} \label{ocn pr} Пусть $H\le G$, $x\in G$ и $\vf,\psi\in \cf(H)$.
Тогда
\begin{list}{{\rm(}{\it\roman{enumi}\/}{\rm)}}
{\usecounter{enumi}\setlength{\parsep}{2pt}\setlength{\topsep}{5pt}\setlength{\labelwidth}{23pt}}
\item $(\vf^x,\psi^x)_{\mbox{}_{H^x}}=(\vf,\psi)_{\mbox{}_H}$.
\item $(\th_{H^x},\vf^x)_{\mbox{}_{H^x}}=(\th_H,\vf)_{\mbox{}_H}$ для произвольной функции $\th\in \cf(G)$.
\item Если $\vf$ --- $($неприводимый$\,)$ характер группы $H$, то $\vf^x$ --- $($неприводимый$\,)$
      характер группы $H^x$.
\end{list}
\end{pre}

\textupl{ocn pr prf}{Доказать предложение \ref{ocn pr}.}

Скалярное произведение проясняет двойственность между операциями ограничения и индуцирования
классовых функций.

\begin{pre}[Закон взаимности Фробениуса] \label{vz fr}\glsadd{iZakVzFr} Пусть $H\le G$, $\vf\in \cf(H)$ и $\psi \in
\cf(G)$.
Тогда
$$
(\vf^G,\psi)_{\mbox{}_G}=(\vf,\psi_{H^{\vphantom{A^a}}})_{\mbox{}_H}
$$
\end{pre}

\textupl{vz fr prf}{Доказать предложение \ref{vz fr}.}

Напомним, что для матрицы $A$ её транспонированная матрица обозначается через $A^\top$. Если $A=(\a_{ij})\in
\MM_n(\CC)$, то положим $\gls{aov}=(\ov{\a}_{ij})$.

\begin{pre}[Второе соотношение ортогональности] \label{vtor ort}\glsadd{iRelOrthScd}
Для любых $g,h\in G$

$$
\sum_{\x\in \irr(G)}\x(g)\ov{\x(h)}=\left\{\ba{ll}
0,&\mbox{если $g$ и $h$ не сопряжены в $G$};\\
|\C_G(g)|& \mbox{в противном случае}.
\ea\right.
$$
\end{pre}
\begin{proof}
Поскольку характеры являются классовыми функциями, первое соотношение ортогональности \ref{perv ort2}
можно переписать в виде
$$
|G|\d_{\x,\th}=\sum_{K\in \K(G)} |K|\x(x_{\mbox{}_K})\ov{\th(x_{\mbox{}_K})}. \myeqno\label{gdxt}
$$
Пусть $X=\XX(G)$ и $D$ --- целочисленная диагональная матрица, индексированная классами
сопряжённости и состоящая из элементов $\d_{\mbox{}_{K,L}}|K|$. Тогда система равенств
\ref{gdxt} при $\x,\th\in \irr(G)$ эквивалентна матричному соотношению
$$
|G|\,\I_{|\K(G)|}=X D\ov{X}^\top.\myeqno\label{xdx}
$$
Согласно \ref{tch nv} матрица $X$ невырождена.
Домножив последнее равенство слева на $X^{-1}$, а справа на
$X$, получаем
$$
|G|\,\I_{|\K(G)|}=D\ov{X}^\top X,
$$
что эквивалентно системе
$$
|G|\d_{\mbox{}_{K,L}}=|K|\sum_{\x\in \irr(G)} \ov{\x(x_{\mbox{}_K})}\x(x_{\mbox{}_L})
$$
при $K,L\in \K(G)$. Отметим, что $|G|/|K|=|\C_G(x_{\mbox{}_K})|$.
Пусть теперь $K=g^G$ и $L=h^G$. Тогда с точностью до комплексного сопряжения доказанное равенство
совпадает с требуемым.
\end{proof}

\textuprn{\label{chr ccs} Показать, что по таблице характеров $\XX(G)$ можно определить размеры классов сопряжённости группы $G$.}

\textuprn{\label{perv st1} Показать, что столбец таблицы характеров группы $G$,
индексированный классом $1^G$ --- единственный столбец c неотрицательными вещественными элементами.}

\section{Определитель таблицы характеров}

Мы отмечали в \ref{tch nv}, что таблица характеров группы является невырожденной матрицей.
Доказательство предложения \ref{vtor ort} позволяет уточнить значение её определителя.
Это значение вычисляется с точностью до знака, поскольку порядок следования строк и столбцов таблицы характеров
не фиксирован.

\begin{pre} \label{detx} Справедливо равенство
$$
\det\XX(G) = \pm \,i^l\, \sqrt{\prod_{K\in \K(G)}|\C_G(x_{\mbox{}_K})|},
$$
где $i=\sqrt{-1}$, $l$ --- число пар комплексно сопряжённых характеров из $\irr(G)$, т.\,е.
двухэлементных множеств $\{\x,\ov\x\}$, $\x \in \irr(G)$.
\end{pre}
\begin{proof} Пусть $X=\XX(G)$. Комплексное сопряжение действует
на множестве $\irr(G)$ как подстановка, равная произведению $l$
независимых транспозиций. Значит, имеет место равенство
$$
\det \ov X = (-1)^l\det X.
$$
Поэтому, взяв определитель обеих частей равенства \ref{xdx}, имеем
$$
|G|^{|\K(G)|}=(-1)^l(\det X)^2 \prod_{K\in \K(G)}|K|.
$$
Учитывая, что $|G|/|K|=|\C_G(x_{\mbox{}_K})|$, получаем требуемое.
\end{proof}

\textupr{Пусть $C$ --- циклическая группа. Является ли определитель $\det \XX(C)$ вещественным? рациональным?}

\textupr{Доказать, что число пар комплексно сопряжённых характеров из $\irr(G)$ совпадает с числом пар <<взаимно обратных>> классов сопряжённости, т.\,е.
двухэлементных множеств $\{K,K^\i\}$, где $K\in \K(G)$ и \gls{Ki} обозначает класс, состоящий из элементов, обратных к элементам класса $K$. \uk{Воспользоваться
\ref{har prop}$(vi)$, а также более общим утверждением \ref{act br}, приведённым ниже.}}
\glsadd{iLemBra}
\begin{pre}[лемма Брауэра] \label{act br}
Пусть группа $A$ действует на множествах $\irr(G)$ и $\K(G)$ так, что
$$
\x^a(x_{\mbox{}_{K^a}})=\x(x_{\mbox{}_K}) \myeqno\label{comp act}
$$
для всех $a\in A$, $\x\in \irr(G)$, $K\in \K(G)$.
Тогда
\begin{list}{{\rm(}{\it\roman{enumi}\/}{\rm)}}
{\usecounter{enumi}\setlength{\parsep}{2pt}\setlength{\topsep}{5pt}\setlength{\labelwidth}{23pt}}
\item число неподвижных точек любого элемента $a\in A$ на множествах $\irr(G)$ и $\K(G)$ совпадает;
\item число орбит группы $A$ на $\irr(G)$ и $\K(G)$ совпадает.
\end{list}

\end{pre}
\begin{proof} Пусть $X=\XX(G)$ и $a\in A$. Обозначим через $P(a)$ матрицу $(p_{\x\eta})$, индексированную элементами $\irr(G)$, где
$p_{\x\eta}=1$, если $\eta=\x^a$, и $0$ в противном случае.
Тогда  $(\x,K)$-й элемент произведения $P(a)X$ равен $\sum_{\eta\in \irr(G)}p_{\x\eta}\,\eta(x_{\mbox{}_K})=\x^a(x_{\mbox{}_K})$. Аналогично, обозначим через
$Q(a)$ матрицу $(q_{\mbox{}_{LK}})$, индексированную элементами $\K(G)$, где $q_{\mbox{}_{LK}}=1$, если $K=L^a$, и $0$ в противном случае.
Тогда $(\x,K)$-й элемент произведения $XQ(a)$ равен $\sum_{L\in \K(G)}\x(x_{\mbox{}_L})\,q_{\mbox{}_{LK}}=\x(x_{\scriptscriptstyle{K^{a^{-1}}}})$.
По условию $\x^a(x_{\mbox{}_K})=\x(x_{\scriptscriptstyle{K^{a^{-1}}}})$, и значит, $P(a)X=XQ(a)$. Поскольку $X$ невырождена, имеем $Q(a)=X^{-1}P(a)X$ и, в частности,
$\tr P(a) = \tr Q(a)$. Однако $\tr P(a)$ совпадает с числом неподвижных точек при действии $a$ на $\irr(G)$, а $\tr Q(a)$~--- при действии на $\K(G)$.
Отсюда следует $(i)$. Утверждение $(ii)$ следует из $(i)$ и того, что число орбит группы $A$ выражается через число неподвижных
точек её элементов, см. \ref{nr go}$(ii)$.
\end{proof}

Одним из основных примеров, когда группа $A$ действует <<согласованно>> на множествах $\irr(G)$ и $\K(G)$, т.\,е. когда выполнено условие \ref{comp act},
является случай $A=\Aut(G)$ и действие $A$ на $\irr(G)$ --- сопряжение классовых функций, а на $\K(G)$ действие поэлементное. При этом
согласованность справедлива в силу определения \ref{spr aut}. Легко видеть, что в этом случае подгруппа внутренних
автоморфизмов $\Inn(G)$ лежит в ядре действия как на $\irr(G)$, так и на $\K(G)$.

\textupr{Пусть $m=\exp G$. Показать, что можно определить действие $A=\Gal(\QQ_m,\QQ)$ на $\K(G)$, которое будет согласовано с
алгебраическим сопряжением на $\irr(G)$. \uk{Воспользоваться \ref{gal conj}$(iii)$.}}

\textext{При всей схожести согласованное действие группы $A$ на множествах $\irr(G)$ и $\K(G)$ не всегда полностью
совпадает. Более строго, соответствующие подстановочные представления не
всегда подстановочно эквивалентны. Например, если $G$ --- неабелева группа порядка $27$ периода $9$ и
$A=\Aut(G)$, то размеры $A$-орбит на $\irr(G)$ равны $1,2,2,3,3$, а на $\K(G)$ --- $1,1,1,2,6$}

\section{Группы Фробениуса}

\begin{opr} Конечная группа $G$, обладающая такой подгруппой $H$, что $1<H<G$ и
$$H\cap H^x=1 \quad  \text{для всякого} \quad x\in G\setminus H,\myeqno\label{cf prop}$$
называется \glsadd{iGrFr}\mem{группой Фробениуса}, а подгруппа $H$ --- \glsadd{iCmplFr}\mem{дополнением Фробениуса} в $G$.
\glsadd{iKerFr}\mem{Ядром Фробениуса} называется множество элементов из $G$, не лежащих в объединении
$$
\bigcup_{x\in G} H^x,\myeqno\label{un cf}
$$
вместе с единичным элементом.
\end{opr}

\textupr{Показать, что всякая группа Фробениуса является транзитивной группой подстановок конечного множества,
в которой стабилизатор каждой точки нетривиален, а стабилизатор любых двух точек тривиален, и обратно,
всякая такая группа подстановок является группой Фробениуса.}

Несмотря на то, что формулировка следующего предложения является чисто теоретико-групповой,
все его известные доказательства не обходятся без использования теории характеров.

\begin{pre}[теорема Фробениуса] \label{gr fr}\glsadd{iThmFrob} Пусть $G$ --- группа Фробениуса с дополнением $H$, причём $|G:H|=n$.
Тогда ядро Фробениуса является нормальной подгруппой в $G$ порядка $n$.
\end{pre}
\begin{proof} Из \ref{cf prop} следует, что $H=\N_G(H)$, и значит, количество различных подгрупп вида $H^x$, где $x\in G$, равно $|G:\N_G(H)|=|G:H|=n$. А поскольку любые две различные подгруппы вида $H^x$ пересекаются тривиально, число неединичных элементов объединения \ref{un cf}
равно $n(|H|-1)=|G|-n$. Т.\,е. мощность множества элементов из $G$, не лежащих в \ref{un cf}, вместе с единичным элементом равна $|G|-(|G|-n)=n$. Очевидно также, что это множество выдерживает сопряжение в группе $G$. Основная сложность --- доказать, что это множество является подгруппой.

Для любого $\th\in \irr(H)\setminus\{1_H\}$ положим
$$
\psi_\th = d_\th 1_H-\th,
$$
где $d_\th=\deg\th$. Таким образом, $\psi_\th\in \gch(H)$, причём $\psi_\th(1)=0$ для всех $\th$. Покажем, что $((\psi_\th)^G)_H=\psi_\th$.

Напомним, что в силу \ref{fif ud} для любого $g\in G$ имеем
$$
(\psi_\th)^G(g)=\sum_{\{r\in R\,\mid\, g\in H^r\}} \psi_\th(rgr^{-1}), \myeqno\label{psthg}
$$
где $R$ --- набор представителей всех правых смежных классов $G$ по $H$. Если теперь  $1\ne h\in H$, то включение $h\in H^r$ возможно
лишь при $r\in H$ ввиду \ref{cf prop}, т.\,е. в этом случае \ref{psthg} влечёт $(\psi_\th)^G(h)=\psi_\th(rhr^{-1})=\psi_\th(h)$.
Если же $h=1$, то $(\psi_\th)^G(h)=n\psi_\th(1)=0=\psi_\th(h)$. Значит, действительно $((\psi_\th)^G)_H=\psi_\th$.

Учитывая доказанное, из закона взаимности \ref{vz fr} вытекает
$$
((\psi_\th)^G,(\psi_\th)^G)_{\mbox{}_G}=(\psi_\th,((\psi_\th)^G)_H)_{\mbox{}_H}=
(\psi_\th,\psi_\th)_{\mbox{}_H}=(d_\th 1_H-\th,d_\th 1_H-\th)_{\mbox{}_H}=d_\th^2+1,
$$
а также
$$
((\psi_\th)^G,1_G)_{\mbox{}_G}=(\psi_\th,1_H)_{\mbox{}_H}=
(d_\th 1_H-\th,1_H)_{\mbox{}_H}=d_\th.
$$
Следовательно, обобщённый характер $(\psi_\th)^G$ группы $G$ имеет вид $d_\th 1_G\pm \x_\th$ для некоторого $\x_\th\in \irr(G)\setminus\{1_G\}$. Если $(\psi_\th)^G = d_\th 1_G+\x_\th$, то это характер, и значит, его степень, равная $(\psi_\th)^G(1)$, положительна. Но это противоречит тому, что $(\psi_\th)^G(1)=0$ как мы видели выше. Значит, $(\psi_\th)^G = d_\th 1_G-\x_\th$.

Таким образом, мы имеем серию неприводимых неглавных характеров $\x_\th = d_\th 1_G - (\psi_\th)^G$
группы $G$, где $\th\in \irr(H)\setminus\{1_H\}$. Заметим, что $(\x_\th)_H=\th$, поскольку
$$
\x_\th(h)= d_\th  - (\psi_\th)^G(h) = d_\th - \psi_\th(h)=d_\th-(d_\th -\th(h))=\th(h)
$$
для всех $h\in H$. Рассмотрим пересечение ядер
$$
K=\bigcap_{\th\in \irr(H)\setminus\{1_H\}} \ker\x_\th.
$$
Это нормальная подгруппа в $G$. Пусть $h\in K\cap H$. Тогда $\th(h)=\x_\th(h)=\x_\th(1)=\th(1)$,
т.\,е. $h\in \ker\th$, для любого $\th\in \irr(H)\setminus\{1_H\}$. Очевидно также, что $h\in \ker 1_H$.
Из \ref{kpr}$(iii)$ следует, что $h=1$. Значит, $K$ тривиально пересекается с объединением \ref{un cf}.

Обратно, пусть элемент $g\in G$ не лежит в объединении \ref{un cf}.
Тогда $\x_\th(g)=d_\th-(\psi_\th)^G(g)=d_\th$ в силу \ref{psthg}, т.\,е. $g\in \ker\x_\th$, для любого $\th\in \irr(H)\setminus\{1_H\}$. Значит, $g\in K$. Таким образом, $K$ совпадает с рассматриваемым множеством элементов из $G$, не лежащих в \ref{un cf}, вместе с единичным элементом. Утверждение доказано.
\end{proof}

\textupr{Пусть $G$ --- группа Фробениуса с дополнением $H$ и ядром $K$. Показать, что
\begin{list}{{\rm(}{\it\roman{enumi}\/}{\rm)}}
{\usecounter{enumi}\setlength{\parsep}{2pt}\setlength{\topsep}{5pt}\setlength{\labelwidth}{23pt}}
\item $G=KH$ и $K\cap H=1$;
\item $\C_G(k)\se K$ для любого неединичного $k\in K$;
\item если $k\in K$ и $h\in H$, то $k^h=k$ тогда и только тогда, когда $k=1$ или $h=1$;
\item $|H|$ делит $|K|-1$ и $|\K(K)|-1$;
\item $|\K(G)|=|\K(H)|+(|\K(K)|-1)/|H|$.
\end{list}}

\section{Подстановочные характеры}

В этом разделе мы будем рассматривать действия групп на конечных множествах и
соответствующие этим действиям подстановочные $\CC$-характеры, см. \ref{pr phg}$(viii)$.

Пусть группа $G$ действует на множестве $X$ и $\x$ --- соответствующий характер.
Тогда для любого $g\in G$ имеем
$$
\x(g)=\big|\{\,x\in X\mid xg=x\,\}\big|, \myeqno\label{pod val}
$$
т.\,е. значение подстановочного характера на элементе $g$ равно числу его неподвижных точек при действии на $X$.
Следующее утверждение описывает характеры транзитивного действия
и, по существу, представляет собой переформулировку на языке характеров утверждения \ref{ed ind}.

\begin{pre} \label{pod ind} Пусть $G$ --- группа.
\begin{list}{{\rm(}{\it\roman{enumi}\/}{\rm)}}
{\usecounter{enumi}\setlength{\parsep}{2pt}\setlength{\topsep}{5pt}\setlength{\labelwidth}{23pt}}
\item Если $G$ действует транзитивно на множестве $X$,  $x\in X$ и $H=\St(x)$, то подстановочный характер
этого действия равен $(1_H)^G$.
\item Если $H\le G$, то $(1_H)^G$ совпадает с характером действия $G$
на правых смежных классах по $H$ согласно правилу $g:Hx\mapsto Hxg$ для всех $g,x\in G$.
\end{list}
\end{pre}
\begin{proof} $(i)$ Пусть $W$ --- соответствующий подстановочный $\CC G$-модуль, см. \ref{pr am}$(vii)$. Легко видеть, что одномерные
подпространства, порождённые базисными элементами $e_y$, $y\in X$, модуля $W$ образуют систему импримитивности. Стабилизатор
подпространства $V=\CC e_x$ совпадает с $H$, поскольку всякий элемент $g\in G$, оставляющий $V$ на месте, должен переводить $e_x$
в $ce_x$ для некоторого $c\in \CC$, но $e_xg=e_{xg}$, и значит, $xg=x$ и $c=1$ в силу линейной независимости базисных элементов,
т.\,е. $g\in H$. В частности, $V$ является
главным $\CC H$-модулем. Из предложения \ref{ind mod} следует, что $W\cong V^G$. Переходя к характерам, получаем
требуемое.

$(ii)$ При данном действии стабилизатор точки $H$, рассматриваемой как элемент множества правых смежных классов,
совпадает с подгруппой $H$. Поэтому требуемое следует из $(i)$.
\end{proof}

Транзитивное действие группы $G$ на множестве $X$ называется \glsadd{iActGrDblTrn}\mem{дважды транзитивным}, если $|X|\ge 2$ и для некоторого (для любого\footnote{Эквивалентность этих определений несложно проверить. Она также следует из \ref{nr go}$(iii)$.}) $x\in X$ стабилизатор $\St(x)$ действует транзитивно на $X\setminus\{x\}$.

\begin{pre} \label{nr go} Пусть $G$ действует на множестве $X$  и $\x$ --- соответствующий
подстановочный характер. Справедливы следующие утверждения.

\begin{list}{{\rm(}{\it\roman{enumi}\/}{\rm)}}
{\usecounter{enumi}\setlength{\parsep}{2pt}\setlength{\topsep}{5pt}\setlength{\labelwidth}{23pt}}
\item $\x=\sum_{O\in \Om}\x_{\mbox{}_O}$, где $\Om$ --- множество всех $G$-орбит на $X$, а $\x_{\mbox{}_O}$ ---
подстановочный характер, соответствующий действию $G$ на орбите $O$.

\item Число $G$-орбит на $X$ равно
$$
(\x,1_G)_{\mbox{}_G}=\frac{1}{|G|}\sum_{g\in G}{\x(g)}.
$$
В частности, действие $G$ на $X$ транзитивно тогда и только, когда $(\x,1_G)_{\mbox{}_G}=1$.
\item Если действие $G$ на $X$ транзитивно и $x\in X$, то число $\St(x)$-орбит на $X$ равно $(\x,\x)_{\mbox{}_G}$ и, в частности, не зависит от $x$.
\item Если $|X|\ge 2$, то действие $G$ на $X$ дважды транзитивно тогда и только, когда $(\x,\x)_{\mbox{}_G}=2$ или, эквивалентно,
когда $\x=1_G+\th$ для некоторого $\th\in \irr(G)$ такого, что $\th\ne 1_G$.
\end{list}
\end{pre}
\begin{proof} $(i)$ Это следует из равенства \ref{pod val}, поскольку $X=\cup_{O\in \Om} O$.

$(ii)$ Пусть $\Om$ --- множество всех $G$-орбит на $X$, $O\in \Om$ и $\x_{\mbox{}_O}$ ---
подстановочный характер, соответствующий действию $G$ на орбите $O$. Обозначим через $H$ стабилизатор
некоторого элемента из $O$. Поскольку действие $G$ на $O$ транзитивно,
из \ref{pod ind}$(i)$ и закона взаимности \ref{vz fr} следует, что
$$
(\x_{\mbox{}_O},1_G)_{\mbox{}_G}=((1_H)^G,1_G)_{\mbox{}_G}=(1_H,(1_G)_H)_{\mbox{}_H}=(1_H,1_H)_{\mbox{}_H}=1.
$$
Значит, в силу $(i)$ и произвольности $O$ получаем
$$
(\x,1_G,)_{\mbox{}_G}=\sum_{O\in \Om}(\x_{\mbox{}_O},1_G)_{\mbox{}_G}=|\Om|,
$$
что и требовалось доказать.

$(iii)$ Обозначим $H=\St(x)$. Ясно, что подстановочный характер действия $H$ на $X$ равен $\x_{H^{\vphantom{A^a}}}$. Поэтому
из $(ii)$ и \ref{vz fr} следует,
что число $H$-орбит на $X$ равно
$$
(\x_{H^{\vphantom{A^a}}},1_{H^{\vphantom{A^a}}})_{\mbox{}_G}=(\x,(1_{H^{\vphantom{A^a}}})^G)_{\mbox{}_G}=(\x,\x)_{\mbox{}_G},
$$
где последнее равенство следует из \ref{pod ind}$(i)$.

$(iv)$ Заметим, что для любого $x\in X$ существует как минимум две $\St(x)$-орбиты, одна и которых равна $\{x\}$.
Если $\x=1_G+\th$, где $\th\in \irr(G)$ и $\th\ne 1_G$, то $(\x,\x)_{\mbox{}_G}=2$
и действие дважды транзитивно в силу $(iii)$.

Обратно, пусть действие дважды транзитивно. Тогда для любого $x\in X$ существует ровно две $\St(x)$-орбиты и
$(\x,\x)_{\mbox{}_G}=2$ в силу $(iii)$.
Поскольку $G$-орбита единственна, из $(ii)$ следует, что $1_G$ является неприводимой
компонентой характера $\x$ кратности $1$, т.\,е. $\x=1_G+\th$, где $\th$ --- некоторый характер группы $G$ такой, что
$(\th,1_G)_{\mbox{}_G}=0$. Но тогда имеем
$$2=(\x,\x)_{\mbox{}_G}=(1_G+\th,1_G+\th)_{\mbox{}_G}=1+(\th,\th)_{\mbox{}_G},$$
откуда следует, что $(\th,\th)_{\mbox{}_G}=1$, т.\,е. $\th$ неприводим.
\end{proof}

\textupln{act pairs prf}{\label{act pairs}Пусть группа $G$ действует на множестве $X$ и  $\chi$ --- соответствующий подстановочный характер.
\begin{list}{{\rm(}{\it\roman{enumi}\/}{\rm)}}
{\usecounter{enumi}\setlength{\parsep}{2pt}\setlength{\topsep}{5pt}\setlength{\labelwidth}{23pt}}
\item Показать, что действие $G$ на $X$ дважды транзитивно
тогда и только тогда\footnote{Вытекающее отсюда эквивалентное определение является даже более классическим: группа $G$ действует на множестве $X$ дважды транзитивно, если $|X|\ge 2$ и для любых пар $(x_1,y_1),(x_2,y_2)\in X\times X$, где $x_1\ne y_1$ и $x_2\ne y_2$,  найдётся элемент $g\in G$ такой, что $(x_1,y_1)g=(x_2,y_2)$.}, когда транзитивно действие $G$ на множестве упорядоченных пар различных элементов 
$$\big\{ (x,y)\in X\times X \mid x\ne y \big\}$$ 
по правилу $g:(x,y)\mapsto (xg,yg)$ для всех $g\in G$. Найти соответствующий этому действию подстановочный характер.
\item Найти подстановочный характер действия $G$ на множестве неупорядоченных пар различных элементов 
$$\big\{ \{x,y\} \subseteq X \mid x\ne y \big\}$$ 
по правилу $g:\{x,y\}\mapsto \{xg,yg\}$ для всех $g\in G$.
\end{list}}

Отметим ряд свойств, которыми обладают подстановочные характеры транзитивного действия.

\begin{pre} \label{pod prop} Пусть $H\le G$ и $\x=(1_H)^G$. Тогда
\begin{list}{{\rm(}{\it\roman{enumi}\/}{\rm)}}
{\usecounter{enumi}\setlength{\parsep}{2pt}\setlength{\topsep}{5pt}\setlength{\labelwidth}{23pt}}
\item $\deg\x$ делит $|G|$;
\item $\x(g)$ --- неотрицательное целое число для любого $g\in G$;
\item $(\x,1_G)_{\mbox{}_G}=1$;
\item $(\x,\th)_{\mbox{}_G}\le\th(1)$ для всех $\th\in \irr(G)$;
\item $\x(g)\le \x(g^n)$ для любых $g\in G$ и $n\in\ZZ$;
\item $\x(g)=0$, если $|g|$ не делит $|G|/\x(1)$;
\item $|\N_G(g)|\deg \x$ делит $|G|\x(g)$ для любого $g\in G$. В частности, $\deg \x$ делит $|g^G|\x(g)$;
\item для любого $g\in G$ имеем
$$
|g^G\cap H|=\frac{\x(g)\,|H|}{|\C_G(g)|}.
$$
\end{list}
\end{pre}
\begin{proof} $(i)$ $\deg\x$ совпадает с индексом $|G:H|$,  который делит $|G|$.

$(ii)$--$(iii)$ Из \ref{pod ind}$(ii)$ следует, что $\x$ --- транзитивный подстановочный характер.
Поэтому его значения целые и неотрицательные в силу \ref{pod val} и кратность
компоненты $1_G$ равна $1$ в силу \ref{nr go}$(ii)$.

$(iv)$ Имеем
$$
(\x,\th)_{\mbox{}_G}=((1_H)^G,\th)_{\mbox{}_G}=(1_H,\th_H)_{\mbox{}_G}\le\th(1).
$$

$(v)$ Это следует из того, что $\x$ --- постановочный характер действия $G$ и
значения $\x$ на $g\in G$ совпадают с числом неподвижных точек $g$ при этом действии, которое всегда
не превосходит числа неподвижных точек любой степени $g^n$, $n\in \ZZ$.

$(vi)$ Поскольку $|g|$ не делит $|H|=|G|/\x(1)$, элемент $g$ не лежит
в подгруппах вида $H^s$, $s\in G$. Поэтому $\x(g)=0$ в силу \ref{fich}.

$(vii)$ Пусть $G$ действует транзитивно на множестве $X$ так, что
соответствующий подстановочный характер равен $\x$. Пусть $g\in G$. Обозначим $K=g^G$ и положим
$$
\A=\big\{(x,k)\in X\times K\bigm|\ xk=x\big\}.
$$
Мощность множества $\A$ можно подсчитать двумя способами: просуммировав либо по $x\in X$, либо по $k\in K$.
Соответственно, получаем
\begin{align*}
|\A|&=\sum_{x\in X} \big|\{k\in K\mid\ xk=x\}\big|=\sum_{x\in X} |K\cap\St(x)|;\\
|\A|&=\sum_{k\in K} \big|\{x\in X\mid\ xk=x\}\big|=\sum_{k\in K}\x(k)=|K|\x(g).
\end{align*}
В силу транзитивности действия, все стабилизаторы $\St(x)$,  $x\in X$, сопряжены в $G$. Поэтому мощность
$m=|K\cap\St(x)|$ не зависит от выбора $x\in X$. Значит, $|X|\,m=|K|\x(g)$.
Также отметим, что если $k\in K\cap\St(x)$, то для
любой степени $k^n$, сопряжённой с $k$, имеем $k^n\in K\cap\St(x)$. Число таких степеней
не зависит от $k$ и равно индексу $|\N_G(g):\C_G(g)|$, а потому является делителем $m=|K|\x(g)/\deg\x$
(поскольку множество $K\cap\St(x)$ разбивается на непересекающиеся множества мощности $m$, каждое из
которых состоит из сопряжённых степеней одного элемента).
Отсюда следует, что число $|\N_G(g)|\deg\x$ делит $|\C_G(g)|\,|K|\x(g)=|G|\x(g)$.

$(viii)$ Из \ref{ind kl}$(iii)$ следует, что произведение $|H|\x(g)$  равно
числу элементов $s\in G$ таких, что $sgs^{-1}\in H$. Это число совпадает с мощностью
$|g^G\cap H|$, умноженной на количество представлений фиксированного элемента класса $g^G$ в
виде $sgs^{-1}$, $s\in G$, которое равно $|\C_G(g)|$. Отсюда получаем требуемое.
\end{proof}

\textext{Иногда возникает вопрос, является данный характер $\x$ группы $G$ подстановочным характером транзитивного действия,
т.\,е. индуцированным с главного характера какой-то заранее неизвестной подгруппы из $G$.
В предложении \ref{pod prop} указаны некоторые необходимые для этого условия,\footnote{По этой причине в некоторых утверждениях
предложения \ref{pod prop} мы избегаем явного упоминания подгруппы $H$, хотя, например, ясно, что $|G|/\x(1)=|H|$.}
которые, вообще говоря, не являются достаточными.
Например, в \cite[стр. 70]{i} отмечается, что простая спорадическая группа Матьё $M_{22}$ обладает характером степени $56$, являющимся суммой главного и неприводимого характеров, для
которого выполнены утверждения \ref{pod prop}$(i)$--$(vii)$, но при этом $M_{22}$ не имеет подгрупп индекса $56$.

Кроме того, предложение \ref{pod prop} можно использовать для получения информации о
подгруппах групп с данной таблицей характеров.\footnote{Отметим, что для данных $g\in G$ и $\x\in \irr(G)$
значения $\x(g^n)$,  $|g|$  и  $|\N_G(g)|$ не всегда можно определить по таблице характеров группы $G$ (пример --- группы $D_8$ и $Q_8$).}
Скажем, рассматривая $\XX(A_4)$ (см. приложение \ref{pril tch}), можно заключить, что у группы $A_4$
не существует характеров степени $2$, которые удовлетворяли бы
условиям \ref{pod prop}$(ii)$--$(iii)$, и поэтому она не содержит подгрупп порядка $6$.}

\section{Характеры прямых произведений}

Пусть $H,K$ --- группы, $G=H\times K$, $\vf\in \cf(H)$,  $\psi\in \cf(K)$. Определим
отображение $\gls{fpsi}:G\to\CC$ по правилу
$$(\vf\times \psi)(hk)=\vf(h)\psi(k)$$
для всех $h\in H$ и
$k\in K$, которое будем называть \glsadd{iProdDirClsFuncs}\mem{прямым произведением}
функций $\vf$ и $\psi$.

\begin{pre} \label{dir ch} Пусть $G=H\times K$.
\begin{list}{{\rm(}{\it\roman{enumi}\/}{\rm)}}
{\usecounter{enumi}\setlength{\parsep}{2pt}\setlength{\topsep}{5pt}\setlength{\labelwidth}{23pt}}

\item Если $\vf\in \cf(H)$ и $\psi\in \cf(K)$, то $\vf\times \psi\in \cf(G)$.
\item Если $\x\in \ch(H)$ и $\th\in \ch(K)$, то $\x\times\th \in \ch(G)$.
\item Если $\vf_1,\vf_2\in \cf(H)$ и $\psi_1,\psi_2\in \cf(K)$, то
$(\vf_1\times\psi_1,\vf_2\times \psi_2)_{\mbox{}_G} = (\vf_1,\vf_2)_{\mbox{}_H}(\psi_1,\psi_2)_{\mbox{}_K}$.
\item $\irr(G)=\{\,\x\times\th\ |\ \x\in \irr(H),\ \th\in\irr(K)\}$.
\end{list}
\end{pre}
\begin{proof} $(i)$ Проверяется непосредственно.

$(ii)$ Используя естественные эпиморфизмы $G/K\cong H$ и $G/H\cong K$, характеры
$\x$ и $\th$ групп $H$ и $K$, соответственно, можно <<поднять>> до характеров $\wt{\x}$ и
$\wt{\th}$ группы $G$, положив $\wt{\x}(hk)=\x(h)$ и $\wt{\th}(hk)=\th(k)$ для всех $h\in H$,
$k\in K$. Тогда $\x\times\th=\wt{\x}\,\wt{\th}$ --- характер группы $G$.

$(iii)$ Имеем

\begin{align*}
(\phi_1\times\psi_1,\phi_2\times \psi_2)_{\mbox{}_G}=\frac{1}{|G|}\sum_{g\in
G}(\phi_1\times\psi_1)(g)\ov{(\phi_2\times\psi_2)(g)}=\frac{1}{|H||K|}\sum_{\substack{h\in H,\\k\in
K}}\phi_1(h)\psi_1(k)\ov{\phi_2(h)}\ov{\psi_2(k)}\\
\mbox{}=\left(\frac{1}{|H|}\sum_{h\in H}\phi_1(h)\ov{\phi_2(h)}\right)\left(
\frac{1}{|K|}\sum_{k\in K}\psi_1(k)\ov{\psi_2(k)}\right)= (\phi_1,\phi_2)_{\mbox{}_H}(\psi_1,\psi_2)_{\mbox{}_K}
\end{align*}

$(iv)$ Из $(ii)$ и $(iii)$ следует, что при различных выборах $\x\in \irr(H)$, $\th\in\irr(K)$ характеры $\x\times\th$ попарно различны и неприводимы.
Тот факт, что характеры $\x\times\th$ исчерпывают всё множество $\irr(G)$ можно доказать, воспользовавшись
\ref{ir eq kl} и заметив, что $|\K(G)|=|\K(H)|\cdot|\K(K)|$, либо можно применить
равенство \ref{fbern}$(ii)$ и соотношения
$$
\sum_{\substack{\x\in \irr(H),\\ \th\in \irr(K)}}(\x\times\th)(1)^2=
\sum_{\substack{\x\in \irr(H).\\ \th\in \irr(K)}}\x(1)^2\th(1)^2=\left(\sum_{\x\in \irr(H)}\x(1)^2\right)
\left(\sum_{\th\in \irr(K)}\th(1)^2\right)=|H||K|=|G|.
$$
\end{proof}

\textuprn{Пусть $\X$ --- обыкновенное представление группы $H$  с характером $\x$ и соответствующим $\CC H$-моду\-лем~$V$.
Аналогично, пусть $\Y$ --- представление группы $K$  с характером $\th$ и соответствующим $\CC K$-модулем~$W$.
Построить представление группы $H\times K$ и соответствующий ему $\CC(H\times K)$-модуль, характер которых совпадает с $\x\times\th$.}

Известно, что всякая конечная абелева группа изоморфна прямому произведению циклических групп, см. \cite[Теорема 8.1.2]{km}.
Поэтому из \ref{dir ch}, \ref{x cyc}$(ii)$ и \ref{com cht}$(ii)$ вытекает следующее утверждение.
\begin{pre} \label{lin ab} Для группы $G$ имеет место изоморфизм $\lin(G)\cong G/G'$. В частности, если $G$ абелева, то $\lin(G)\cong G$.
\end{pre}

\section{Характеры и нормальные подгруппы}

Пусть $H\nor G$, $x\in G$ и $\vf\in \cf(H)$. Определим отображение $\vf^x:H\to \CC$ по правилу
$$
\vf^x(h)=\vf(xhx^{-1})
$$
для всех $h\in H$.
Из \ref{fcn pr}$(i)$ следует, что $\vf^x$ --- классовая функция, сопряжённая с $\vf$, а из \ref{ocn pr}$(iii)$ --- что
группа $G$ действует на множестве $\irr(H)$, причём образ характера $\th\in
\irr(H)$ под действием $x\in G$ равен $\th^x$, и нормальная подгруппа $HC_G(H)\nor G$  лежит в ядре этого
действия. Стабилизатор элемента $\th\in \irr(H)$ относительно такого действия
назовём \glsadd{iGrInrChr}\mem{группой инерции характера} $\th$ и 
обозначим его через \gls{IGth}, ср. \ref{pgin m}. Другими словами, $$\II_G(\th)=\{x\in G \mid \th^x=\th\}.$$
Из сказанного следует, что $HC_G(H)\le \II_G(\th)$ для любого $\th\in \irr(H)$.

\textuprn{Пусть $H\nor G$.
\begin{list}{{\rm(}{\it\roman{enumi}\/}{\rm)}}
{\usecounter{enumi}\setlength{\parsep}{2pt}\setlength{\topsep}{5pt}\setlength{\labelwidth}{23pt}}
\item Найти $\II_G(1_H)$.
\item Показать, что для характера $\th\in\irr(H)$ справедливо равенство $\II_G(\th)=G$ в том и только в том
    случае, когда $\th$ принимает постоянные значения на всех классах сопряжённости группы $G$, содержащихся в $H$.
\end{list}}

Теорема Клиффорда \ref{thm cliff} об ограничении неприводимых модулей на нормальные подгруппы может
быть переформулирована на языке обыкновенных характеров в уточнённом виде. Мы приведём здесь независимое доказательство
этой важной теоремы, чтобы продемонстрировать эффективность теории характеров.

\begin{thm}[Клиффорда] \label{tcl char}\glsadd{iThmCliffChar} Пусть $H\nor G$ и $\x\in \irr(G)$. Пусть
$\th$ --- неприводимая компонента ограничения $\x_{H^{\vphantom{A^a}}}$ и
$e=(\x_{H^{\vphantom{A^a}}},\th)_{\mbox{}_H}$. Тогда справедливы следующие утверждения.
\begin{list}{{\rm(}{\it\roman{enumi}\/}{\rm)}}
{\usecounter{enumi}\setlength{\parsep}{2pt}\setlength{\topsep}{5pt}\setlength{\labelwidth}{23pt}}
\item  Если $\th_1,\th_2,\ld,\th_t$ --- полный набор
различных характеров вида $\th^x$, $x\in G$, то
$$
\x_{H^{\vphantom{A^a}}}=e(\th_1+\th_2+\ld+\th_t).
$$
\item Пусть $I=\II_G(\th)$. Существует единственный характер $\vf\in \irr(I)$ такой, что
$$(\vf_{H^{\vphantom{A^a}}},\th)_{\mbox{}_H}>0\qquad \text{и} \qquad (\vf^G,\x)_{\mbox{}_G}>0.$$
Для такого $\vf$ имеем $\vf_{H^{\vphantom{A^a}}}=e\,\th$ и $\vf^G=\x$.
\end{list}
\end{thm}
\begin{proof} $(i)$ Из закона взаимности Фробениуса \ref{vz fr} следует, что
$$
(\th^G,\x)_{\mbox{}_G}=(\th,\x_{H^{\vphantom{A^a}}})_{\mbox{}_H}>0,
$$
т.\,е. $\x$ является компонентой характера $\th^G$.
Поэтому всякая неприводимая компонента ограничения $\x_{H^{\vphantom{A^a}}}$ также будет
неприводимой компонентой характера $(\th^G)_H$. Изучим его более подробно.

В силу \ref{ind kl}$(iii)$ для произвольного $h\in H$ имеем
$$
\th^G(h)=\frac{1}{|H|}\sum_{x\in G}\th^\circ(xhx^{-1})=\frac{1}{|H|}\sum_{x\in G}\th^x(h),
$$
поскольку $xhx^{-1}\in H$ ввиду нормальности $H$. Значит, справедливо равенство классовых функций
$|H|(\th^G)_H=\sum_{x\in G}\th^x$. Поэтому все неприводимые компоненты характера $(\th^G)_H$ содержатся
среди $\th_1,\th_2,\ld,\th_t$. Из вышесказанного следует,
что это же справедливо и для характера $\x_{H^{\vphantom{A^a}}}$.
Таким образом, $\x_{H^{\vphantom{A^a}}}=\sum_{i=1}^t(\x_{H^{\vphantom{A^a}}},\th_i)_{\mbox{}_H}\th_i$.
Ввиду \ref{ocn pr}$(ii)$ произведение $(\x_{H^{\vphantom{A^a}}},\th_i)_{\mbox{}_H}=(\x_{H^{\vphantom{A^a}}},\th)_{\mbox{}_H}$ не зависит от $i$.
Отсюда следует требуемое.

$(ii)$ Заметим, что число $t$ различных сопряжённых с $\th$ характеров совпадает с мощностью
орбиты $\th$ при действии группы $G$ сопряжением и поэтому $t=|G:I|$.

Установим существование требуемого характера $\vf\in\irr(I)$.
По условию $\th$ --- неприводимая компонента ограничения
$\x_{H^{\vphantom{A^a}}}=(\x_{I^{\vphantom{A^a}}})_{H^{\vphantom{A^a}}}$. Поэтому существует неприводимая
компонента $\vf$ характера $\x_{I^{\vphantom{A^a}}}$ такая, что $\vf_{H^{\vphantom{A^a}}}$ имеет неприводимую
компоненту $\th$. В частности, $\vf$ удовлетворяет требуемым неравенствам
$(\vf_{H^{\vphantom{A^a}}},\th)_{\mbox{}_H}>0$ и $(\vf^G,\x)_{\mbox{}_G}=(\vf,\x_{I^{\vphantom{A^a}}})_{\mbox{}_G}>0$.

Теперь покажем, что для  любого $\vf\in\irr(I)$, удовлетворяющего неравенствам
$(\vf_{H^{\vphantom{A^a}}},\th)_{\mbox{}_H}>0$ и $(\vf^G,\x)_{\mbox{}_G}>0$, справедливы
равенства $\vf_{H^{\vphantom{A^a}}}=e\,\th$ и $\vf^G=\x$. В самом деле, в силу утверждения $(i)$,
применённого к группе $I$, её неприводимому характеру $\vf$ и нормальной подгруппе $H$, имеем
$\vf_{H^{\vphantom{A^a}}}=f\th$ для некоторого натурального $f$, поскольку $\th^x=\th$ для всех $x\in I$.
Так как $\vf$ --- неприводимая компонента характера $\x_{I^{\vphantom{A^a}}}$, рассмотрев
ограничения $\vf_{H^{\vphantom{A^a}}}$ и $\x_{H^{\vphantom{A^a}}}$, получим $f\le e$.
Поэтому выполнены соотношения
$$
et\deg\th=\deg\x\le\deg\vf^G=t\deg\vf=tf\deg\th\le te\deg\th
$$
и, значит, все неравенства в этой цепочке являются равенствами. В частности, $f=e$, откуда
получаем $\vf_{H^{\vphantom{A^a}}}=e\,\th$, а также $\deg\x=\deg\vf^G$, откуда $\vf^G=\x$.

Докажем единственность найденного характера $\vf$. Пусть существует
$\psi\in\irr(I)$ такой, что $\psi\ne\vf$
и $(\psi_{H^{\vphantom{A^a}}},\th)_{\mbox{}_H}>0$ и $(\psi^G,\x)_{\mbox{}_G}>0$.
Заметим, что $(\psi,\x_{I^{\vphantom{A^a}}})_{\mbox{}_I}=(\psi^G,\x)_{\mbox{}_G}>0$, т.\,е. $\psi$, равно как и $\vf$,
является неприводимой компонентой ограничения $\x_{I^{\vphantom{A^a}}}$. Поэтому
$$
e=(\x_{H^{\vphantom{A^a}}},\th)_{\mbox{}_H}=((\x_{I^{\vphantom{A^a}}})_{H^{\vphantom{A^a}}},\th)_{\mbox{}_H}\ge
((\vf+\psi)_{H^{\vphantom{A^a}}},\th)_{\mbox{}_H}=(\vf_{H^{\vphantom{A^a}}},\th)_{\mbox{}_H}+
(\psi_{H^{\vphantom{A^a}}},\th)_{\mbox{}_H}>(\vf_{H^{\vphantom{A^a}}},\th)_{\mbox{}_H}=e.
$$
Полученное противоречие завершает доказательство.
\end{proof}

Из теоремы \ref{tcl char} вытекает

\begin{pre}\label{ctc h} Пусть $H\nor G$, $\th\in\irr(H)$ и $I=\II_G(\th)$.
Положим
\begin{align*}
 \A&=\big\{\vf\in\irr(I) \bigm| (\vf_{H^{\vphantom{A^a}}},\th)_{\mbox{}_H}>0  \big\}, \\
  \B&=\big\{\x\in\irr(G) \bigm|  (\x_{H^{\vphantom{A^a}}},\th)_{\mbox{}_H}>0 \big\}.
\end{align*}
Отображение $\vf\mapsto\vf^G$ осуществляет биекцию $\A\to\B$.
\end{pre}
\begin{proof} Пусть $\vf\in \A$. Выберем $\x\in\irr(G)$ так, что $(\vf^G,\x)_{\mbox{}_G}>0$.
Поскольку $(\vf,\x_{I^{\vphantom{A^a}}})_{\mbox{}_I}=(\vf^G,\x)_{\mbox{}_G}>0$ и
$\th$ --- неприводимая компонента ограничения $\vf_{H^{\vphantom{A^a}}}$, отсюда следует, что
$\th$ также является неприводимой компонентой ограничения
$(\x_{I^{\vphantom{A^a}}})_{H^{\vphantom{A^a}}}=\x_{H^{\vphantom{A^a}}}$, т.\,е. $\x\in \B$.
Кроме того, мы находимся в условиях теоремы \ref{tcl char}, из пункта $(ii)$ которой следует, что $\vf^G=\x$,
т.\,е. соответствие $\vf\mapsto\vf^G$ действительно является
отображением $\A\to\B$. Инъективность этого отображения следует из единственности,
а сюръективность --- из существования характера $\vf\in\irr(I)$ в заключении утверждения \ref{tcl char}$(ii)$.
\end{proof}

\textupln{cor cl prf}{\label{cor cl}Пусть $H\nor G$ и $\th\in\irr(H)$. Запишем
$$
\th^G=\sum_{\x\in \irr(G)}n_\x \x,
$$
где $n_\x $ --- однозначно определённые неотрицательные целые числа.  Показать, что справедливо равенство
$$
|I:H|=\sum_{\x\in \irr(G)}n_\x^2,\myeqno\label{ihs}
$$
где $I=\II_G(\th)$.}

В случае, когда $H=1$ (а значит, $\th=1_H$ и $I=G$), характер $\th^G$ совпадает с регулярным и поэтому $n_\x=\x(1)$ в силу
\ref{fbern}$(i)$. Таким образом, формула Бернсайда \ref{fbern}$(ii)$ является частным случаем соотношения~\ref{ihs}.

\textext{В обозначениях упражнения \ref{cor cl} также справедлив тот факт, что все ненулевые
числа $n_\x$ являются делителями индекса $|I:H|$ (доказательство см. в \cite[Следствие 11.29]{i}).
Мы покажем далее в \ref{deg div}, что степени неприводимых обыкновенных характеров группы делят её порядок.
Возникающая таким образом аналогия между числами $n_\x$ и степенями неприводимых характеров факторгруппы $I/H$ не случайна:
ненулевые коэффициенты $n_\x$ будут степенями так называемых проективных неприводимых представлений группы~$I/H$.}

Характер $\x\in\ch(G)$ называется \glsadd{iChrOrdImpr}\mem{импримитивным},
если $\x=\th^G$, где $\th$ --- обыкновенный характер некоторой собственной
подгруппы $H\le G$. Неприводимый характер, не являющийся импримитивным,
называется \glsadd{iChrOrdPrm}\mem{примитивным}.

Из \ref{cr pr} следует, что характер примитивен (импримитивен) тогда и только тогда, когда он
является характером примитивного (импримитивного) $\CC G$-модуля.

Характер $\x\in\ch(G)$ называется \glsadd{iChrOrdHmg}\mem{однородным},
если он кратен неприводимому характеру, т.\,е. имеет вид $\x=n\th$ для натурального $n$ и $\th\in\irr(G)$.

\begin{pre} \label{prch n} Пусть $\x$ --- примитивный характер группы $G$ и $N\nor G$. Тогда
\begin{list}{{\rm(}{\it\roman{enumi}\/}{\rm)}}
{\usecounter{enumi}\setlength{\parsep}{2pt}\setlength{\topsep}{5pt}\setlength{\labelwidth}{23pt}}
\item ограничение $\x_{N^{\vphantom{A^a}}}$ будет
однородным характером;
\item если $\x$ точный и $N$ абелева, то $N\le\Z(G)$.
\end{list}
\end{pre}
\begin{proof} $(i)$ По определению $\x\in \irr(G)$. Пусть $\th$ --- неприводимая компонента
ограничения $\x_{N^{\vphantom{A^a}}}$ и $I=\II_G(\th)$.
По теореме \ref{tcl char} имеем $\x_{N^{\vphantom{A^a}}}=e(\th_1+\ld+\th_t)$, где $t=|G:I|$ и $\th_1,\ld,\th_t$ ---
различные сопряжённые с $\th$ характеры. Кроме того,
$\x=\vf^G$ для некоторого $\vf\in \irr(I)$. В силу примитивности $\x$ получаем, что $I=G$ и $t=1$. Отсюда
следует требуемое.

$(ii)$ Из $(i)$ следует,
что $\x_{N^{\vphantom{A^a}}}=e\l$ для некоторого неприводимого характера $\l$ группы $N$, который
линеен ввиду \ref{ab npr}. В частности, для любого $x\in N$
имеем $\x(x)=e\l(x)$ и, значит, $|\x(x)|=e=\x(1)$. Поэтому $N\le\Z(\x)$. Поскольку $\x$ точен
и неприводим, из \ref{zh pr}$(iv)$ следует, что $Z(\x)=\Z(G)$. Значит, $N\le\Z(G)$.
\end{proof}

\textzam{Отметим, что обращение утверждения \ref{prch n}$(i)$ неверно. Например,
ограничение неприводимого характера степени $5$ знакопеременной группы $G=A_5$ (см. приложение \ref{pril tch})
на любую её нормальную
подгруппу очевидным образом однородно, однако этот характер импримитивен, т.\,к. равен $\th^G$,
где $\th$ --- нетривиальный линейный характер группы $A_4$.}

\section{\texorpdfstring{$M$}{M}-группы}

Матрица, каждая строка и каждый столбец которой содержат ровно один ненулевой элемент,
называется \glsadd{iMatMnm}\mem{мономиальной}. Легко видеть, что мономиальная матрица
над полем невырождена.

Обыкновенный характер $\x$ группы $G$ называется \glsadd{iChrOrdMnm}\mem{мономиальным},
если $\x=\l^G$, где $\l$ --- линейный характер некоторой подгруппы $H\le G$.

В частности, всякий мономиальный характер степени больше $1$ импримитивен и, значит,
является характером импримитивного $\CC G$-модуля с системой импримитивности, состоящей из
одномерных подпространств.

\textuprn{Показать, что
\begin{list}{{\rm(}{\it\roman{enumi}\/}{\rm)}}
{\usecounter{enumi}\setlength{\parsep}{2pt}\setlength{\topsep}{5pt}\setlength{\labelwidth}{23pt}}
\item всякий мономиальный характер группы $G$ является характером представления $\X$
такого, что все матрицы $\X(g)$, $g\in G$, мономиальны;
\item характер неприводимого\footnote{Условие неприводимости существенно как видно
на примере прямой суммы двух линейных представлений.} представления $\X$ группы $G$ такого, что
все матрицы $\X(g)$, $g\in G$, мономиальны, будет мономиальным.
\end{list}
}

Класс групп, все неприводимые характеры которых мономиальны, представляет особый интерес.

\begin{opr}
Группа $G$ называется \glsadd{iMGr}\mem{$M$-группой}, если каждый характер $\x\in \irr(G)$ мономиален.
\end{opr}

\begin{pre} \label{mgr prop} Пусть $G$ --- $M$-группа и
$$
1=d_1<d_2<\ld <d_k
$$
--- все различные степени неприводимых характеров группы $G$. Если характер  $\x\in \irr(G)$ имеет степень $d_i$,
то $G^{(i)}\le \ker \x$, где $G^{(i)}$ обозначает $i$-й коммутант группы $G$.
\end{pre}
\begin{proof} Будем рассуждать индукцией по $i$.
Если $i=1$, то $\deg \x=1$, т.\,е. $\x$ --- линейный характер. Всякий линейный характер группы $G$ является
гомоморфизмом из $G$ в абелеву группу $\CC^\times$. Поэтому $G'\le \ker\x$.

Пусть $i>1$. Тогда $\deg\x>1$. Сначала заметим, что если $\th\in \irr(G)$ и $\deg\th < \deg\x$, то $\deg\th=d_j$
для некоторого $j<i$. Поэтому $G^{(i-1)}\le G^{(j)}\le \ker\th$.

Поскольку $G$ --- $M$-группа, найдётся подгруппа $H\le G$ и её линейный характер $\l$ такой, что $\x=\l^G$.
Подгруппа $H$ собственная, потому что $\deg\x>1$.
Заметим, что в силу \ref{vz fr} справедливы равенства
$$((1_H)^G,1_G)_{\mbox{}_G}=(1_H,1_H)_{\mbox{}_H}=1.$$
Значит, характер $(1_H)^G$ не является неприводимым, т.\,к. $H$ --- собственная подгруппа.
Если $\th$ --- его произвольная неприводимая компонента,
то
$$\deg\th<\deg (1_H)^G=|G:H|=\deg \l^G=\deg \x.$$
Ввиду вышесказанного имеем $G^{(i-1)}\le \ker\th$. По \ref{kpr}$(ii)$ в силу произвольности $\th$
справедливо включение $G^{(i-1)}\le \ker (1_H)^G$, а в силу \ref{ker indh} получаем
$$\ker (1_H)^G=\bigcap_{x\in G}H^x\le H.$$
Таким образом, $G^{(i-1)}\le H$.  Но тогда
$$
G^{(i)}=(G^{(i-1)})'\le H'\le \ker\l.
$$
В силу нормальности $G^{(i)}\nor G$, вновь применив \ref{ker indh}, получаем
$$
G^{(i)}\le \bigcap_{x\in G} (\ker\l)^x=\ker \l^G=\ker \x,
$$
что и требовалось доказать.
\end{proof}

\begin{cor} \label{mgr solv} Всякая $M$-группа разрешима.
\end{cor}
\begin{proof} Пусть $G$ --- $M$-группа. В обозначениях предложения \ref{mgr prop} получаем, что
$G^{(i)}$ содержится в пересечении ядер всех характеров $\x\in \irr(G)$, степень которых совпадает
с одним из чисел $d_1,d_2,\ld,d_i$. В частности,
$$G^{(k)}\le \bigcap_{\x\in \irr(G)} \ker\x=1$$
ввиду \ref{kpr}$(iii)$. Отсюда следует требуемое.
\end{proof}

Примером разрешимой группы, не являющейся $M$-группой, служит $\SL_2(3)$. Эта группа обладает
неприводимым обыкновенным характером степени $2$ (см. приложение \ref{pril tch}). Такой характер
не является индуцированным с собственной подгруппы, поскольку $\SL_2(3)$ не имеет подгрупп индекса $2$.

Нам понадобится следующий теоретико-групповой факт.
\normalmarginpar

\begin{pre} \label{nma} Пусть $G$ --- нильпотентная группа и $A\nor G$ --- максимальная
абелева нормальная подгруппа. Тогда $\C_G(A)=A$.
\end{pre}
\begin{proof} Cм. \cite[Теорема 16.2.6]{km}.
\end{proof}

Следующее предложение вместе со следствием \ref{mgr solv} показывает, что класс $M$-групп заключён между классами нильпотентных и разрешимых групп.

\begin{pre} \label{nil mgr} Всякая
нильпотентная группа является $M$-группой.
\end{pre}
\begin{proof} Пусть $G$ --- нильпотентная группа и $\x\in\irr(G)$. Пусть $H$ --- минимальная по включению
подгруппа в $G$, для которой существует $\th\in\irr(H)$ такой, что $\x=\th^G$. Поскольку индуцирование классовых функций
транзитивно \big(см. \ref{ind kl}$(v)$\big), получаем, что $\th$ --- примитивный характер группы $H$.
Пусть $\wt H=H/\ker\th$. Из предложения \ref{ch gn}$(iii)$-$(iv)$ следует, что характер $\wt\th$ группы $\wt H$,
соответствующий характеру $\th$, будет точным и примитивным. Поэтому всякая нормальная абелева подгруппа
группы $\wt H$ будет лежать в $\Z(\wt H)$ ввиду \ref{prch n}$(ii)$. Однако $\wt H$ нильпотентна
и, значит, в силу \ref{nma} содержит нормальную абелеву подгруппу, совпадающую со своим централизатором.
Следовательно, эта подгруппа совпадает с $\wt H$, т.\,е. группа $\wt H$ абелева.
Поэтому характер $\wt\th$, а, значит, и $\th$, линеен в силу неприводимости. Отсюда получаем требуемое.
\end{proof}

В качестве примера $M$-группы, которая не является нильпотентной, можно взять группу $G=S_3$.
Её единственный нелинейный неприводимый характер имеет вид $\l^G$, где $\l$ --- не главный
неприводимый характер знакопеременной подгруппы $A_3$.

\section{Центральные характеры}

Пусть $G$ --- группа и $\x\in \irr(G)$. Определим отображения $\om_\x:\Z(\CC G)\to \CC$.
Из предложения \ref{cent irr}
вытекает, что для любого $z\in\Z(\CC G )$ матрица $\X_\x(z)$ скалярна. Пусть $\X_\x(z)=\ve\I$ для
$\ve\in \CC$. Поскольку скалярная матрица совпадает с любой своей сопряжённой, число $\ve$ не
зависит от выбора представления $\X_\x$ с характером $\x$. Мы положим $\om_\x(z)=\ve$. Другими
словами, для всех $z\in \Z(\CC G)$
$$
\X_\x(z)=\om_\x(z)\I.
$$
Так как $\X_\x$ --- гомоморфизм алгебр, легко видеть, что $\om_\x : \Z(\CC G )\to \CC$ также
является гомоморфизмом, и следовательно линейным характером центра $\Z(\CC G )$.
Мы будем называть этот гомоморфизм
\glsadd{iChrCenCorIrr}\mem{центральным характером, соответствующим неприводимому характеру} $\x$.

Центральный характер $\om_\x$ полностью определяется своими значениями
на любом фиксированном базисе алгебры $\Z(\CC G)$.
Одним из таких базисов является набор центральных идемпотентов, как следует из \ref{cor z ss}$(iv)$.
Другим естественным базисом $\Z(\CC G)$, как мы знаем из
\ref{kl sum}$(ii)$, являются классовые суммы группы~$G$. Значения центральных характеров на элементах
этих базисов могут быть найдены явно.

\begin{pre} \label{cent dif}\mbox{}
\begin{list}{{\rm(}{\it\roman{enumi}\/}{\rm)}}
{\usecounter{enumi}\setlength{\parsep}{2pt}\setlength{\topsep}{5pt}\setlength{\labelwidth}{23pt}}
\item Пусть $\x\in \irr(G)$. Тогда для любого $\th\in \irr(G)$

$$\om_\x(e_\th)=\d_{\x,\th},$$
где $e_\th$ --- центральный идемпотент, соответствующий характеру $\th$. В частности, центральные характеры
$\om_\x$, $\x \in \irr(G)$, попарно различны.
\item Все гомоморфизмы $\CC$-алгебр $\Z(\CC G )\to \CC$ исчерпываются множеством $\{\om_\x\,|\,\x\in
\irr(G)\}$.
\item Для произвольного $z\in \Z(\CC G)$ имеет место представление
$$
z=\sum_{\x\in \irr(G)}\om_\x(z)e_\x.
$$
В частности, если $K\in \K(G)$, то
$$ \wh K=\sum_{\x\in \irr(G)} \om_\x(\wh K)e_\x.$$
\end{list}
\end{pre}
\vspace{-10pt}
\textupl{cent dif prf}{Доказать предложение \ref{cent dif}.}

Напомним, что значения характера $\x$ группы $G$ можно рассматривать на произвольном элементе алгебры~$\CC G$ как
след матрицы представления, соответствующего $\x$.

\begin{pre} \label{om ob} Пусть $\x\in \irr(G)$ и $z\in \Z(\CC G)$. Тогда справедливы следующие утверждения.
\begin{list}{{\rm(}{\it\roman{enumi}\/}{\rm)}}
{\usecounter{enumi}\setlength{\parsep}{2pt}\setlength{\topsep}{5pt}\setlength{\labelwidth}{23pt}}
\item $\om_\x(z)=\displaystyle\frac{\x(z)}{\x(1)}$. В частности, если  $K\in \K(G)$, то
\end{list}
\vspace{-1.5\topsep}
$$\om_\x(\wh{K})=\frac{|K|\x(x_{\mbox{}_K})}{\x(1)},\myeqno\label{ch val}$$
\begin{list}{{\rm(}{\it\roman{enumi}\/}{\rm)}}
{\usecounter{enumi}\addtocounter{enumi}{1}\setlength{\parsep}{2pt}\setlength{\topsep}{5pt}\setlength{\labelwidth}{23pt}}
\item[{}] где $x_{\mbox{}_K}$~--- представитель класса $K$.
\item Для произвольного $a\in \CC G$ имеем
$$
\x(za)=\om_\x(z)\x(a)=\frac{\x(z)\x(a)}{\x(1)}.
$$
\item Если $z\in \Z(G)$, то $\om_\x(z)^{|z|\strut}=1$. В частности, $|\om_\x(z)|=1$.
\end{list}
\end{pre}
\textupl{om ob prf}{Доказать предложение \ref{om ob}.}

Ввиду второго соотношения ортогональности \ref{vtor ort} по столбцам таблицы характеров $\XX(G)$,
индексированным классами $1^G$ и $K$ можно найти $|G|$ и $|\C_G(x_{\mbox{}_K})|$, соответственно,
а значит, и $|K|$. Поэтому из формулы \ref{ch val} вытекает, что
все значения $\om_\x$ на классовых суммах определяются по $\XX(G)$.

По аналогии с таблицей $\XX(G)$ определим
\mem{\glsadd{iTblOrdChr}таблицу \gls{OmG} центральных характеров}, соответствующих неприводимым обыкновенным характерам группы $G$
как квадратную матрицу, строки которой индексированы неприводимыми характерами $\x\in\irr(G)$,
а столбцы --- классами (или представителями классов) сопряжённости $K\in \K(G)$. При этом
элемент $\x$-й строки и $K$-го столбца матрицы $\Om(G)$
равен $\om_\x(\wh K)$. Например, ниже приведена таблица $\Om(S_3)$.
$$
  \begin{array}{c|ccc}
            \Om(S_3)              & 1      & (12)           & (123)       \\
     \hline
    \x_1^{\vphantom{A^A}} & 1      & \phantom{-}3  & \phantom{-}2\\
    \x_2                  & 1      & -3            & \phantom{-}2\\
    \x_3                  & 1      & \phantom{-}0  & -1
  \end{array}
$$
Из \ref{cent dif}$(iii)$ следует, что $\Om(G)$ является матрицей перехода от базиса $\wh{K}$, $K\in \K(G)$,
к базису $e_\x$, $\x\in \irr(G)$,  центра алгебры $\CC G$. Поэтому справедливо
\begin{cor} \label{ttsh nv} Для любой
группы $G$ матрица $\Om(G)$ является невырожденной.
\end{cor}

Следующее  утверждение показывает, что матрица $\Om(G)$ лежит в кольце $\MM_{|\irr(G)|}(\ov{\ZZ})$.

\begin{pre} \label{ch ac} Пусть $\x\in \irr(G)$ и $K\in \K(G)$. Тогда значение $\om_\x(\wh{K})$
является целым алгебраическим числом.
\end{pre}
\begin{proof} В \ref{kl sum}$(iv)$ показано, что для любых классов $K,L\in \K(G)$ имеет место равенство
$$\wh{K}\wh{L}=\sum_{M\in \K(G)} a_{\mbox{}_{KLM}} \whm,$$ 
где  $a_{\mbox{}_{KLM}}\in
\ZZ$. Поскольку $\om_\x$ --- гомоморфизм алгебр $\Z(\CC G)\to \CC$, отсюда следует, что
$$
\om_\x(\wh{K})\om_\x(\wh{L})=\sum_{M\in \K(G)} a_{\mbox{}_{KLM}} \om_\x(\whm). \myeqno\label{omhkl}
$$
Пусть $W$ --- подмодуль  $\ZZ$-модуля $\CC$, порождённый значениями $\{\om_\x(L)\,|\, L\in \K(G)\}$. Тогда $W$
--- ненулевой конечно порождённый $\ZZ$-модуль и в силу \ref{omhkl} для любого $K\in \K(G)$ выполнено
включение $\om_\x(\wh{K})W\se W$. Из предложения \ref{alg mlt} вытекает, что $\om_\x(\wh{K})$ --- целое
алгебраическое число.
\end{proof}

Следует подчеркнуть, что включение $\om_\x(\wh{K})\in \ov{\ZZ}$ не следует из того факта,
что $\x(x_{\mbox{}_K})\in \ov{\ZZ}$, поскольку произведение целого алгебраического числа и рационального числа,
вообще говоря, не лежит в $\ov{\ZZ}$, ср. \ref{lem alg cn}.

\section{Центр групповой алгебры \texorpdfstring{$\CC G$}{CG}}


Коэффициенты разложения произвольного элемента $z\in \Z(\CC G)$ по базису из центральных идемпотентов
можно найти, если известны значения $\x(z)$, где $\x\in\irr(G)$, воспользовавшись \ref{cent dif}$(iii)$ и \ref{om ob}$(i)$.
А именно,
$$
z=\sum_{\x\in \irr(G)}\om_\x(z)e_\x=\sum_{\x\in \irr(G)}\frac{\x(z)}{\x(1)}\,e_\x.
$$

В базисе из центральных идемпотентов структурные константы $\Z(\CC G)$ получаются из \ref{id ort}.
Значения структурных констант центра $\Z(\CC G)$ в базисе из классовых сумм, т.\,е. таких чисел  $a_{\mbox{}_{KLM}}$, что
$$\wh{K}\wh{L}=\sum_{M\in \K(G)} a_{\mbox{}_{KLM}} \whm$$
для любых $K,L\in \K(G)$,
были найдены в \ref{kl sum}$(iii)$--$(iv)$. Следующее утверждение позволяет выразить эти константы
через значения неприводимых характеров на элементах классов $K,L,M$.

\begin{pre} \label{str yav} Пусть $K,L,M\in \K(G)$. Тогда имеет место выражение
$$
a_{\mbox{}_{KLM}}=\frac{|K||L|}{|G|}\sum_{\x\in
\irr(G)}\frac{\x(x_{\mbox{}_K})\x(x_{\mbox{}_L})\ov{\x(x_{\mbox{}_M})}}{\x(1)}.
$$
\end{pre}

\textupl{str yav prf}{Доказать предложение \ref{str yav}. \uk{Воспользоваться \ref{cent dif}$(iii)$, \ref{ch val} и \ref{id dec}}.}

Таким образом, из \ref{str yav} и \ref{kl sum}$(iii)$ получаем

\begin{cor} \label{str chr} Пусть $K,L\in \K(G)$ и $z\in G$. Тогда
$$
\big|\{(x,y)\in K\times L\mid xy=z\}\big| = \frac{|K||L|}{|G|}\sum_{\x\in
\irr(G)}\frac{\x(x_{\mbox{}_K})\x(x_{\mbox{}_L})\ov{\x(z)}}{\x(1)}.
$$
\end{cor}

\textuprn{Пусть $s\ge 1$ и $K_1,\ld,K_s\in \K(G)$. Показать, что число
наборов $(x_1,\ld,x_s)$,  где $x_i\in K_i$, таких, что
$$
x_1\cdot\ld\cdot x_s=1
$$
равно $$\frac{|K_1|\cdot\ld\cdot|K_s|}{|G|}\sum_{\x\in
\irr(G)}\frac{\x(x_1)\cdot\ld\x(x_s)}{\x(1)^{s-2}}.$$
\uk{Применить индукцию по $s$.}}

\section{Степени неприводимых характеров}

Важным следствием предложения \ref{ch ac} является

\begin{pre} \label{deg div} Если $\x\in \irr(G)$, то $\deg\x$ делит $|G|$.
\end{pre}
\begin{proof} Из первого соотношения ортогональности \ref{cor perv ort} вытекает, что
$$
|G|=\sum_{g\in G}\x(g)\x(g^{-1}).
$$
Перепишем это равенство в терминах центрального характера $\om_\x$.
Из \ref{ch val} следует, что
$$
|G|=\sum_{K\in \K(G)}|K|\x(x_{\mbox{}_K})\x(x_{\mbox{}_K}^{-1})=\sum_{K\in
\K(G)}\x(1)\om_\x(\wh{K})\x(x_{\mbox{}_K}^{-1}).
$$
Поэтому
$$
\frac{|G|}{\x(1)}=\sum_{K\in \K(G)}\om_\x(\wh{K})\x(x_{\mbox{}_K}^{-1}).
$$
В силу \ref{ch ac}, \ref{ch val} и \ref{alg kol}$(i)$ выражение справа является целым
алгебраическим числом, а поскольку выражение слева --- рациональное число, по \ref{z cz}
обе части являются целыми рациональными числами. Отсюда следует требуемое.
\end{proof}

В действительности имеет место даже более сильный результат. Теорема Ито \ref{tito}
утверждает, что степень неприводимого обыкновенного характера группы $G$ делит индекс любой её нормальной
абелевой подгруппы. Для доказательства этого факта нам потребуются вспомогательные понятия и результаты,
которые представляют самостоятельный интерес.

Для произвольного элемента $g\in G$ обозначим через \gls{vklgr}
число способов (возможно, нулевое) представить $g$ в виде коммутатора
двух элементов из $G$. Другими словами,
$$
\vk(g)=\big|\big\{(a,b)\in G\times G\ |\ [a,b]=g\big\}\big|.
$$

\begin{pre} \label{cas} Пусть $G$ --- группа. Справедливы следующие утверждения.
\begin{list}{{\rm(}{\it\roman{enumi}\/}{\rm)}}
{\usecounter{enumi}\setlength{\parsep}{2pt}\setlength{\topsep}{5pt}\setlength{\labelwidth}{23pt}}
\item $\vk\in\cf(G)$.
\item Для любого $g\in G$ значение $\vk(g)$ делится на $|\Z(G)|^2$.
\item Имеет место разложение
$$
\vk=\sum_{\x\in \irr(G)}\frac{|G|}{\x(1)}\,\x.
$$
В частности, $\vk\in\ch(G)$.
\item Пусть $c=\sum_{a,b\in G}[a,b]$. Тогда $c\in \Z(\CC G)$ и справедливы равенства
$$
c=\sum_{K\in \K(G)}\vk(x_{\mbox{}_K})\wh{K}=\sum_{\x\in \irr(G)}\left(\frac{|G|}{\x(1)}\right)^2e_\x,
$$
где $x_{\mbox{}_K}$ --- представитель класса $K$.
\end{list}
\end{pre}
\begin{proof} $(i)$ Если $g=[a,b]$, то $g^x=[a^x,b^x]$ для любого $x\in G$. Поэтому $\vk$ является классовой
функцией на группе $G$.

$(ii)$ Пусть $z,w\in \Z(G)$. Отображение множества $\big\{(a,b)\in G\times G\ |\ [a,b]=g\big\}$ в себя,
при котором $(a,b)$ переходит в $(az,bw)$, определяет действие группы $\Z(G)\times \Z(G)$. Легко видеть, что все орбиты
этого действия имеют мощность $|\Z(G)|^2$. Отсюда следует требуемое.

$(iii)$--$(iv)$ Ввиду $(i)$  и \ref{id dec} имеем
\begin{align*}
c=\sum_{g\in G}\vk(g)g= \sum_{K\in \K(G)}\vk(x_{\mbox{}_K})\wh{K}
&=\sum_{K\in \K(G)}\vk(x_{\mbox{}_K})\left(\sum_{\x\in \irr(G)}\frac{|K|\x(x_{\mbox{}_K})}{\x(1)}\,e_\x\right)\\
&=\sum_{\x\in \irr(G)}\frac{|G|}{\x(1)}\left(\frac{1}{|G|}\sum_{K\in \K(G)}|K|\vk(x_{\mbox{}_K})\x(x_{\mbox{}_K})\right)e_\x=
\sum_{\x\in \irr(G)}\frac{|G|}{\x(1)}(\vk,\x)_{\mbox{}_G}e_\x,
\end{align*}
где в последнем равенстве мы воспользовались тем, что значения $\vk$ выдерживают комплексное сопряжение. В частности, $c\in \Z(\CC G)$.

Пусть $h\in G$ и $L=h^G$. Тогда в силу \ref{cent dif}$(iii)$ имеем
$$
\sum_{a\in G}a^{-1}ha=|\C_G(h)|\,\wh{L}=\frac{|G|}{|L|}\sum_{\x\in \irr(G)}\frac{|L|\x(h)}{\x(1)}\,e_\x
=\sum_{\x\in \irr(G)}\frac{|G|\x(h)}{\x(1)}\,e_\x.
$$
Значит, учитывая, что $e_\x$ --- центральные идемпотенты алгебры $\CC G$, получаем
\begin{align*}
c=\sum_{b\in G}\left(\ \sum_{a\in G}a^{-1}b^{-1}a\right)b=
\sum_{b\in G}\left(\sum_{\x\in \irr(G)}\frac{|G|\x(b^{-1})}{\x(1)}\,e_\x\right)b&=
\sum_{\x\in \irr(G)}\left(\frac{|G|}{\x(1)}\right)^2\left(\frac{\x(1)}{|G|}\sum_{b\in G}\x(b^{-1})b\right)e_\x\\
&=\sum_{\x\in \irr(G)}\left(\frac{|G|}{\x(1)}\right)^2e_\x^2
=\sum_{\x\in \irr(G)}\left(\frac{|G|}{\x(1)}\right)^2e_\x,
\end{align*}
где мы воспользовались выражением для $e_\x$ из \ref{id dec}. Отсюда следует $(iv)$.

Сравнивая коэффициенты при $e_\x$ в двух полученных выражениях для $c$, заключаем, что
$(\vk,\x)_{\mbox{}_G}=|G|/\x(1)$ для всех $\x\in\irr(G)$. Отсюда ввиду \ref{kf raz} получаем $(iii)$.
Значит, $\vk$ --- характер в силу \ref{deg div}.
\end{proof}

Элемент $c\in \Z(\CC G)$ из \ref{cas}$(iv)$, равный сумме всех коммутаторов группы $G$,
в литературе называют \glsadd{iElmKaz}\mem{элементом Казимира} алгебры $\CC G$.

\begin{cor}\label{com opr} По таблице
характеров $\XX(G)$ можно определить, представим ли элемент $g$ данного класса сопряжённости группы $G$
в виде $g=[a,b]$ для некоторых $a,b\in G$.
\end{cor}
\begin{proof} Элемент $g$ является коммутатором тогда и только тогда, когда $\vk(g)>0$.
В силу \ref{cas}$(iii)$ значение $\vk(g)$ можно найти по таблице $\XX(G)$.
\end{proof}

\textuprn{Пусть $G$ --- группа и $c=\sum_{a,b\in G}[a,b]$. Показать, что
$$
c=|G|\sum_{K\in \K(G)}\frac{\wh{K^{\phantom{\i}}}\wh{K^\i}}{|K|},
$$
где $K^\i$ обозначает класс, состоящий из элементов, обратных к элементам класса $K$.}

\textuprn{Пусть $G$ --- группа, $g\in G$ и $n$ --- целое число такое, что $(n,|g|)=1$. Показать,
что если $g$ является коммутатором, то элемент $g^n$ также будет коммутатором.
\uk{Воспользоваться \ref{gal conj}$(iii)$ и \ref{pr hr}$(i)$.}}

%
%

Теперь докажем следующее усиление предложения \ref{deg div}.

\begin{pre} \label{div pz} Пусть $\x\in\irr(G)$.
Тогда $\deg\x$ делит $|G:\Z(G)|$.
\end{pre}
\begin{proof} Обозначим $Z=\Z(G)$. Разделив на $|Z|^2$ обе суммы в равенстве из \ref{cas}$(iv)$,
получим
$$
\sum_{K\in \K(G)}\frac{\vk(x_{\mbox{}_K})}{|Z|^2}\,\wh{K}=
\sum_{\th\in \irr(G)}\left(\frac{|G:Z|}{\th(1)}\right)^2e_\th.
$$
Применив центральный характер $\om_\x$ и воспользовавшись \ref{cent dif}$(i)$, получим
$$
\sum_{K\in \K(G)}\frac{\vk(x_{\mbox{}_K})}{|Z|^2}\,\om_\x(\wh{K})=\left(\frac{|G:Z|}{\x(1)}\right)^2.
$$
В силу \ref{ch ac} и \ref{cas}$(ii)$ левая часть лежит в $\ov{\ZZ}$, поскольку является целочисленной линейной комбинацией
целых алгебраических чисел. Правая часть рациональна и, значит, является целым числом в силу \ref{z cz}. Отсюда получаем
требуемое. \end{proof}

\begin{cor} \label{cd kh} Пусть $\x\in\irr(G)$.
Тогда $\deg\x$ делит $|G:\Z(\x)|$.
\end{cor}
\begin{proof} Обозначим $N=\ker\x$ и $\wt G=G/N$. Из \ref{zh pr}$(iv)$ следует, что $|G:\Z(\x)|=|\wt{G}:\Z(\wt G)|$.
Пусть $\wt\x$ --- характер факторгруппы $\wt G$, соответствующий $\x$, см. \ref{ch gn}. Тогда $\deg\wt\x=\deg\x$.
Поскольку $\deg\wt\x$ делит $|\wt{G}:\Z(\wt G)|$, получаем требуемое.
\end{proof}

\begin{pre}[Теорема Ито] \label{tito}\glsadd{iThmIto} Пусть $A\nor G$ --- абелева нормальная подгруппа.
Тогда $\deg\x$ делит $|G:A|$ для любого $\x\in\irr(G)$.
\end{pre}
\begin{proof} Выберем характер $\l\in \irr(A)$ такой,
что $(\l,\x_{A^{\vphantom{A^a}}})_{\mbox{}_A}>0$ и обозначим $I=\II_G(\l)$. Из \ref{tcl char}$(ii)$
следует, что для некоторого $\vf\in\irr(I)$ и натурального числа $e$ справедливы равенства $\vf_{A^{\vphantom{A^a}}}=e\l$ и $\vf^G=\x$.
Поэтому для любого $a\in A$ имеем $\vf(a)=e\l(a)$, и, значит, $|\vf(a)|=e=\vf(1)$, поскольку характер $\l$ линейный.
Таким образом, $A\se \Z(\vf)$, и в силу \ref{cd kh} степень $\deg\vf$ делит индекс $|I:\Z(\vf)|$,
который в свою очередь делит $|I:A|$. Осталось заметить,
что $\deg \x=|G:I|\deg\vf$, откуда следует, что $\deg \x$ делит $|G:I|\,|I:A|=|G:A|$.
\end{proof}

\section{\texorpdfstring{$p^a q^b$}{paqb}-теорема Бернсайда}

Известная теорема Бернсайда о разрешимости конечных групп бипримарного порядка является
глубоким результатом, доказываемым с помощью теории характеров. Эта теорема опирается на
два вспомогательных утверждения, которые представляют самостоятельный интерес.

\begin{pre} \label{vch z} Пусть $\x\in\irr(G)$
и $K\in \K(G)$. Если $(|K|,\x(1))=1$, то либо $\x(x_{\mbox{}_K})=0$, либо $x_{\mbox{}_K}\in \Z(\x)$,
где $x_{\mbox{}_K}$ --- представитель класса $K$.
\end{pre}
\begin{proof} Поскольку $(|K|,\x(1))=1$, существуют целые числа $m,n$ такие, что $m|K|+n\x(1)=1$.
Поэтому в силу \ref{ch ac} величина
$$
m\frac{|K|\x(x_{\mbox{}_K})}{\x(1)}=\frac{(1-n\x(1))\x(x_{\mbox{}_K})}{\x(1)}=
\frac{\x(x_{\mbox{}_K})}{\x(1)}-n\x(x_{\mbox{}_K})
$$
является  целым алгебраическим числом. Положим $\a=\x(x_{\mbox{}_K})/\x(1)$.
Так как $n\x(x_{\mbox{}_K})\in \ov{\ZZ}$, имеем $\a\in \ov{\ZZ}$.

Предположим, что $x_{\mbox{}_K}\not\in \Z(\x)$, т.\,е. $|\x(x_{\mbox{}_K})|<\x(1)$. Тогда $|\a|<1$.

Обозначим $k=|x_{\mbox{}_K}|$. В силу \ref{har prop}$(iii)$ значение
$\x(x_{\mbox{}_K})$ лежит в $\QQ_k$  и является суммой $\x(1)$ корней из $1$ степени $k$. Поэтому для любого
автоморфизма $\s\in\Gal(\QQ_k,\QQ)$ образ $\x(x_{\mbox{}_K})^\s$ также лежит в $\QQ_k$ и является
суммой $\x(1)$ корней из $1$, откуда следует, что $|\x(x_{\mbox{}_K})^\s|\le\x(1)$, т.\,е. $|\a^\s|\le 1$.
Значит, положив
$$
\b=\prod_{\s\in\Gal(\QQ_k,\QQ)}\a^\s,
$$
имеем $|\b|<1$. Поскольку $\b^\s=\b$  для всех $\s\in\Gal(\QQ_k,\QQ)$, из \ref{nor gal}$(ii)$
следует, что $\b\in \QQ$. С другой стороны, все величины $\a^\s$ являются корнями того же многочлена
с целыми коэффициентами, что и $\a$, т.\,е. будут целыми алгебраическими числами.  Поэтому $\b\in \ov{\ZZ}$.
Из \ref{z cz} получаем, что $\b\in \ZZ$. Но тогда $\b=0$, поскольку $|\b|<1$. Значит, $\a^\s=0$
для некоторого $\s$, т.\,е. $\a=0$ и  $\x(x_{\mbox{}_K})=0$.
\end{proof}

\begin{pre}[$p^a$-лемма Бернсайда] \label{pa lb}\glsadd{iLemPBur} Пусть $G$ --- простая неабелева конечная группа и $K\in \K(G)$.
Если $|K|$ --- степень простого числа, то  $K=\{1\}$.
\end{pre}
\begin{proof} Пусть $K\in \K(G)$, $|K|=p^a$, где $p$ простое, и $x_{\mbox{}_K}$ --- представитель класса $K$.
Предположим, что $x_{\mbox{}_K}\ne 1$. Пусть $\x\in\irr(G)$ и $\x\ne 1_G$. Тогда $\ker\x=1$,
поскольку группа $G$ простая, и $\Z(\x)=\Z(G)=1$, поскольку $G$ неабелева, см. \ref{zh pr}$(iv)$.
Значит, если $p\nmid\x(1)$, то $\x(x_{\mbox{}_K})=0$ в силу предложения \ref{vch z}.
Поэтому из \ref{fbern}$(i)$ следует, что
$$
0=\r(x_{\mbox{}_K})=\sum_{\x\in \irr(G)} \x(1)\x(x_{\mbox{}_K})=
1+\sum_{\substack{\x\in \irr(G),\\p|\x(1)}} \x(1)\x(x_{\mbox{}_K}).
$$
Таким образом, $1+p\a=0$, где
$$
\a=\sum_{\substack{\x\in \irr(G),\\p|\x(1)}} \frac{\x(1)}{p}\x(x_{\mbox{}_K})
$$
--- целое алгебраическое число. Но из \ref{z cz} следует, что $\a=-1/p\not\in\ov{\ZZ}$, противоречие.
\end{proof}

\textext{Оказывается, что в произвольной конечной группе класс сопряжённых элементов, порядок которого равен степени простого числа,
порождает разрешимую (нормальную) подгруппу. Это усиление $p^\a$-леммы Бернсайда было доказано Л.\,C.\,Казариным с использованием теории
модулярных характеров, см. \cite{kaz}.}

\begin{pre}[$p^a q^b$-теорема Бернсайда] \label{pab b}\glsadd{iThmPQBur} Пусть $|G|=p^a q^b$, где $p$ и $q$ --- простые числа.
Тогда группа $G$ разрешима.
\end{pre}
\begin{proof} Будем рассуждать индукцией по порядку $|G|$. Можно считать, что $|G|>1$.
Если группа $G$ не проста, то она обладает собственной нетривиальной нормальной подгруппой $N$.
Группы $N$ и $G/N$ имеют бипримарные порядки и, значит, разрешимы по индукции. Следовательно,
$G$ также разрешима.

Пусть $G$ простая.  Выберем силовскую подгруппу $P$  группы $G$ так, что $P\ne 1$,
и элемент $x\in \Z(P)$ так, что $x\ne 1$.
Поскольку $P\le C_G(x)$, порядок $|x^G|=|G:C_G(x)|$ делит индекс $|G:P|$, т.\,е. является
степенью простого числа. Из предложения \ref{pa lb} следует, что группа $G$ абелева и, в частности,
разрешима. \end{proof}

\section{Брауэрова характеризация обобщённых характеров}

Обсудим методы определения принадлежности классовой функции $\vf\in \cf(G)$
множеству характеров $\ch(G)$ или кольцу обобщённых характеров $\gch(G)$.

Если известно множество $\irr(G)$, то из свойств скалярного произведения
характеров \ref{cor scl} вытекает, что $\vf\in\gch(G)$ тогда и
только тогда, когда $(\vf,\x)_{\mbox{}_G}\in \ZZ$ для всех $\x\in \irr(G)$. Более
типичной является ситуация, когда таблица характеров группы $G$ неизвестна, но имеется некоторый
класс $\H$ подгрупп группы $G$ такой, что $\vf_{H^{\vphantom{A^a}}}\in \gch(H)$
для всех $H\in \H$. Оказывается, это условие является достаточным для
включения $\vf\in  \gch(G)$, если $\H$ --- семейство так называемых элементарных брауэровых подгрупп, см. ниже.
В этом состоит суть характеризационной теоремы Брауэра \ref{thm char br}.

\begin{opr} Группа $E$ называется \glsadd{iGrPElmBr}\mem{$p$-элементарной (брауэровой)} для некоторого простого числа $p$, если она является прямым
произведением циклической группы и $p$-группы. Будем говорить, что $E$
\glsadd{iGrElmBr}\mem{элементарная (брауэрова)}, если она $p$-элементарна для некоторого
простого числа $p$.
\end{opr}

Заметим, что $p$-элементарную группу можно было бы эквивалентным образом определить как прямое произведение
$p$-группы и циклической $p'$-группы. Легко видеть, что элементарные группы нильпотентны и
всякая подгруппа и факторгруппа элементарной группы элементарна.

Пусть $G$ --- группа. На протяжении этого раздела символ $\E$ будет обозначать множество всех элементарных подгрупп
группы $G$.

\begin{thm}[Брауэра] \label{thm char br}\glsadd{iThmBra} Пусть $\vf\in \cf(G)$.
Тогда следующие условия эквивалентны.
\begin{list}{{\rm(}{\it\roman{enumi}\/}{\rm)}}
{\usecounter{enumi}\setlength{\parsep}{2pt}\setlength{\topsep}{5pt}\setlength{\labelwidth}{23pt}}
\item $\vf\in \gch(G)$
\item $\vf_{E^{\vphantom{A^a}}}\in \gch(E)$ для всякой подгруппы $E\in\E$.
\item $\vf$ является $\ZZ$-линейной комбинацией характеров вида $\l^G$, где $\l\in \lin(E)$ и $E\in \E$.
\end{list}
\end{thm}

Прежде чем доказывать теорему Брауэра нам потребуется ряд дополнительных результатов и определений.

Обозначим через $\A$, $\B$, $\CCC$ множества классовых функций  $\vf\in \cf(G)$, удовлетворяющих условиям $(i)$, $(ii)$, $(iii)$
теоремы \ref{thm char br}, соответственно. Теорема Брауэра утверждает, что  $\A=\B=\CCC$.
Мы знаем, что $\A=\gch(G)$ --- подкольцо в $\CC$-алгебре $\cf(G)$ с единицей $1_G$.
Между множествами $\A$, $\B$, $\CCC$ справедливы следующие соотношения.

\begin{pre} \label{pr abc} Имеем
\begin{list}{{\rm(}{\it\roman{enumi}\/}{\rm)}}
{\usecounter{enumi}\setlength{\parsep}{2pt}\setlength{\topsep}{5pt}\setlength{\labelwidth}{23pt}}
\item $\CCC\se\A\se \B$.
\item $\B$ --- подкольцо в $\,\cf(G)$.
\item $\CCC$ --- идеал в $\B$.
\end{list}
\end{pre}
\begin{proof} $(i)$--$(ii)$ Включение $\CCC\se\A$ следует из определения обобщённых характеров.
Включение $\A\se\B$ и тот факт, что $\B$ --- кольцо, справедливы, поскольку отображение
$\vf\mapsto\vf_{E^{\vphantom{A^a}}}$, где $\vf\in\cf(G)$,
является гомоморфизмом колец $\cf(G)\to\cf(E)$, см. \ref{cff}$(ii)$. Ясно также,
что $1_G\in\B$, т.\,е. $\B$ --- подкольцо  в $\cf(G)$.

$(iii)$ Пусть $\vf\in \CCC$ и $\psi\in \B$. Из определения множества $\CCC$ и линейности операции индуцирования
\big(см. \ref{ind kl}$(ii)$\big) следует, что $\vf$
является суммой обобщённых характеров вида $(\th_{\mbox{}_{(E)}})^G$, где
функции $\th_{\mbox{}_{(E)}}\in \gch(E)$ индексированы подгруппами $E\in\E$. Из \ref{ind kl}$(vi)$
следует, что произведение $\vf\,\psi$ является суммой обобщённых характеров
вида $(\th_{\mbox{}_{(E)}}\psi_{E^{\vphantom{A^a}}})^G$,
где  $E\in\E$, а поэтому и $\ZZ$-линейной комбинацией характеров вида $(\x_{\mbox{}_{(E)}})^G$,
где $\x_{\mbox{}_{(E)}}\in \irr(E)$ и $E\in\E$. Но всякая элементарная подгруппа $E$ нильпотентна. Поэтому
в силу \ref{nil mgr} любой характер из $\irr(E)$ имеет вид $\l^E$ для некоторого линейного характера $\l$
подгруппы из $E$, которая также будет элементарной. Из этих замечаний ввиду транзитивности
индуцирования следует требуемое.
\end{proof}

Таким образом, предложение \ref{pr abc} показывает, что для доказательства теоремы Брауэра достаточно
установить включение $1_G\in\CCC$.

Следующее утверждение носит вспомогательный характер.

\begin{pre} \label{l b} Пусть $X$ --- непустое
конечное множество и $R$ --- аддитивная подгруппа в $\f_{\ZZ}(X)$, замкнутая по умножению. Если $R$ не подкольцо $\big($т.\,е.
не содержит единицу кольца $\f_{\ZZ}(X)\big)$, то существуют элемент $x\in X$ и простое число $p$ такие, что
$p$ делит $\vf(x)$ для всех $\vf\in R$.
\end{pre}
\begin{proof} Для каждого $x\in X$ положим
$$
I_x=\{\,\vf(x)\bigm|\vf\in R\,\}.
$$
Легко видеть, что $I_x$ --- аддитивная подгруппа в $\ZZ$. Если для некоторого $x\in X$ имеет место
строгое включение $I_x\subset \ZZ$, то существует простое число $p$ такое, что $I_x\se p\,\ZZ$ и требуемое
выполнено. Предположим, что $I_x=\ZZ$ для всех $x\in X$. Тогда для любого $x$ найдётся функция $\vf_x\in R$
такая, что $\vf_x(x)=1$. Тогда $(\ve-\vf_x)(x)=0$, где $\ve$ --- единица кольца $\f_{\ZZ}(X)$.
Поэтому функция $\prod_{x\in X}(\ve-\vf_x)$ --- нулевой элемент из $\f_{\ZZ}(X)$.
Раскрывая скобки, получаем, что $\ve$ --- $\ZZ$-линейная комбинация произведений элементов $\vf_x\in R$, т.\,е.
$\ve\in R$ вопреки условию.
\end{proof}

Следующее утверждение, в частности, показывает, что $\ZZ$-линейные комбинации характеров вида $(1_H)^G$, где $H\le G$,
образуют подкольцо в $\gch(G)$.

\begin{pre} \label{per pod} Пусть $H,K\le G$. Тогда
$$
(1_H)^G(1_K)^G=\sum_{L\le H}a_{\mbox{}_L}(1_L)^G,
$$
где $a_{\mbox{}_L}$ --- неотрицательные целые числа.
\end{pre}
\begin{proof} Положим $\th=(1_K)^G$. Из \ref{ind kl}$(vi)$ следует, что
$$
(1_H)^G(1_K)^G=(1_H)^G\,\th=(\th_H)^G.
$$
Поскольку $\th_H$ --- подстановочный характер группы $H$, из \ref{pod ind}$(i)$ и \ref{nr go}$(i)$ следует, что $\th_H$
является суммой характеров вида $(1_L)^H$ где $L$ пробегает
некоторое множество подгрупп группы $H$. В силу аддитивности индуцирования $(\th_H)^G$
является суммой характеров вида $((1_L)^H)^G=(1_L)^G$, что и требовалось.
\end{proof}

Нам потребуется рассмотреть класс групп, более широкий, чем элементарные группы.

Пусть $p$ --- простое число. Группа $Q$ называется \glsadd{iGrPQsElm}\mem{$p$-квазиэлементарной},
если в $Q$ найдётся циклическая нормальная $p'$-подгруппа $A$, факторгруппа по которой будет $p$-группой
(в этом случае говорят, что $Q$ обладает циклическим
\glsadd{iNorPCmpl}\mem{нормальным $p$-дополнением} $A$).
Будем называть $Q$
\glsadd{iGrpQsElm}\mem{квазиэлементарной}, если она $p$-квазиэлементарна для некоторого
простого числа $p$.

Легко видеть, что всякая элементарная группа является квазиэлементарной, и любая подгруппа квазиэлементарной
группы будет квазиэлементарной.

В силу \ref{pod prop}$(ii)$ значения
характеров группы $G$ вида $(1_H)^G$, где $H\le G$, являются целочисленными. Справедливо следующее утверждение.

\begin{pre} \label{pzn qe} Пусть $g\in G$ и $p$ --- простое число.
Тогда существует $p$-квазиэлементарная подгруппа $Q\le G$ такая, что  $(1_Q)^G(g)$ не делится на $p$.
\end{pre}
\begin{proof} Рассмотрим подгруппу $\la g \ra\le G$. Она $p$-квазиэлементарна. Пусть $A$ ---
нормальное $p$-дополнение в $\la g \ra$ и $N=\N_G(A)$. Легко видеть, что $g\in N$. Поскольку
$\la g \ra/A$ --- $p$-группа, существует  подгруппа $Q\le N$ такая, что $Q/A\in \Syl_p(N/A)$ и
$\la g \ra\le Q$. Покажем, что $Q$ --- требуемая подгруппа. Она $p$-квазиэлементарна, так как $A$ является
в ней нормальным $p$-дополнением.

Заметим, что всякий элемент из $G$, централизующий $g$, также централизует $A$. Поэтому $\C_G(g)=\C_N(g)$.
Далее, если  $g^x\in Q$ для некоторого $x\in G$, то $A^x\le Q$. Но всякая $p'$-подгруппа из $Q$ содержится
в $A$, откуда следует, что $A^x=A$ и $x\in N$. Поэтому $g^G\cap Q=g^N\cap Q$. Из этих замечаний с
учётом \ref{pod prop}$(viii)$ получаем
$$
(1_Q)^G(g)=\frac{|g^G\cap Q|\,|\C_G(g)|}{|Q|}=\frac{|g^N\cap Q|\,|\C_N(g)|}{|Q|}=(1_Q)^N(g).
$$
Так как $(1_Q)^N$ --- подстановочный характер действия группы $N$ на правых смежных классах по $Q$,
нам достаточно установить, что число неподвижных точек элемента $g$ относительно этого действия не делится
на $p$. Подгруппа $A$ содержится в ядре этого действия, поскольку $A\nor N$ и $A\le Q$. Значит, искомое
число совпадает с числом неподвижных точек элемента $Ag$ при действии факторгруппы $N/A$ на правых смежных классах по $Q/A$.
Поскольку $Ag$~--- $p$-элемент, количество перемещаемых точек делится на $p$, откуда получаем,
что искомое число сравнимо по модулю $p$ с  индексом $|N/A:Q/A|$, который не делится на $p$ в силу выбора $Q$,
и требуемое доказано.
\end{proof}

\begin{pre} \label{od rep} Главный характер $1_G$
является $\ZZ$-линейной комбинацией характеров вида $(1_Q)^G$ для квазиэлементарных
подгрупп $Q$ группы~$G$.\end{pre}
\begin{proof}
Обозначим через $\Q$ множество всех квазиэлементарных подгрупп группы~$G$, а через
$R$ --- аддитивную группу $\ZZ$-линейных комбинаций характеров вида $(1_Q)^G$, где $Q\in \Q$.
В силу \ref{pod prop}$(ii)$ получаем, что $R$ --- подгруппа в $\cf_\ZZ(G)$.
Поскольку класс $\Q$ замкнут относительно взятия подгрупп, из \ref{per pod} вытекает,
что $R$ замкнута по умножению. Применим к $R$ предложение \ref{l b}.
Если $1_G\not\in R$, то найдутся простое число $p$ и элемент $g\in G$ такие, что $p$ делит $(1_Q)^G(g)$
для всякой $p$-квазиэлементарной подгруппы $Q\le G$. Однако, это противоречит предложению \ref{pzn qe}.
Поэтому имеет место требуемое включение $1_G\in R$.
\end{proof}

Нам потребуется ещё одно вспомогательное утверждение.

\begin{pre} \label{vs tb} Пусть $H=AP$,
где $A\nor H$, $P$ --- $p$-подгруппа в $H$ и $p\nmid|A|$. Если характер $\l\in \lin(A)$ такой, что $\l^x=\l$ для всех $x\in H$,
и $\C_A(P)\se \ker\l$, то $\l=1_A$.
\end{pre}
\begin{proof} Обозначим $K=\ker\l$. Поскольку $\l$ осуществляет изоморфизм между $A/K$ и конечной подгруппой в $\CC^\times$,
легко видеть, что $\l$ принимает постоянные значения на смежных классах $A$ по $K$ и все такие значения различны. По условию
$$\l(a)=\l^x(a)=\l(xax^{-1})$$
для всех $a\in A$ и $x\in P$. Поэтому $P$ оставляет инвариантным каждый смежный класс $A$ по $K$ при действии сопряжением.
Число перемещаемых точек в каждом смежном классе $Ka$, $a\in A$, при таком действии кратно $p$, а поскольку $p\nmid|K|$,
получаем, что $Ka\cap \C_A(P)\ne \varnothing$ и, значит, $A=K\C_A(P)$. Но по условию $\C_A(P)\se K$, откуда следует, что
$A=K$ и $\l=1_A$.
\end{proof}

Теперь мы можем доказать теорему Брауэра \ref{thm char br}.
\begin{proof} Как мы отметили выше, достаточно установить, что главный характер $1_G$ лежит в множестве $\CCC$, т.\,е.
является $\ZZ$-линейной комбинацией характеров вида $\l^G$, где $\l\in \lin(E)$ и $E\in \E$.

Можно считать, что $G$ не является элементарной, поскольку в противном случае утверждение очевидно.
Будем рассуждать индукцией по порядку $|G|$. Достаточно показать, что $1_G$ является $\ZZ$-линейной комбинацией
характеров, индуцированных с собственных подгрупп $H<G$. В самом деле, если это так, то,
учитывая, что по индукции для любой собственной подгруппы $H<G$ справедливо утверждение теоремы Брауэра \ref{thm char br} и, значит,
всякий характер $H$ представим в виде $\ZZ$-линейной комбинацией характеров вида $\l^H$,
где $\l$ --- линейный характер некоторой элементарной подгруппы из $H$,
применив транзитивность и линейность индуцирования, получим искомое представление для $1_G$.

Если $G$ не является квазиэлементарной, то требуемое следует из предложения \ref{od rep}. Поэтому можно считать, что
для некоторого простого числа $p$ группа $G$ обладает
нормальным циклическим $p$-дополнением $A$. Пусть $P\in \Syl_p(G)$ и $Z=\C_A(P)$. Поскольку $G$ не элементарна,
имеем $Z<A$ и $E=PZ<G$. Так как
$$((1_E)^G,1_G)_{\mbox{}_G}=(1_E,1_E)_{\mbox{}_E}=1,$$
получаем, что $(1_E)^G=1_G+\th$, где $\th$ --- характер
группы $G$. Если мы докажем, что всякая неприводимая компонента характера $\th$ индуцирована с собственной подгруппы,
то получим требуемое представление $1_G=(1_E)^G-\th$.

Пусть $\x$ --- неприводимая компонента $\th$. Покажем, что $1_A$ не является неприводимой компонентой
ограничения $\x_{A^{\vphantom{A^a}}}$. Из строения группы $G$ следует, что $G=AE$ и $A\cap E=Z$. Поэтому в силу
\ref{ind kl}$(vii)$ получаем цепочку равенств
$$
1=(1_Z,1_Z)_{\mbox{}_Z}=((1_Z)^A,1_A)_{\mbox{}_A}=\big(((1_E)^G)_A,1_A\big)_{\mbox{}_A}
=(1_A+\th_A,1_A)_{\mbox{}_A}=1+(\th_A,1_A)_{\mbox{}_A}.
$$
Значит, $(\th_{A^{\vphantom{A^a}}},1_{A^{\vphantom{A^a}}})_{\mbox{}_A}=0$ и $(\x_{A^{\vphantom{A^a}}},1_{A^{\vphantom{A^a}}})_{\mbox{}_A}=0$, как и утверждалось.

Пусть $\l$ --- неприводимая компонента ограничения $\x_{A^{\vphantom{A^a}}}$. Поскольку $A$ абелева, характер $\l$ линеен.
 По доказанному имеем $\l\ne 1_A$. Заметим,
что $Z\se\ker\l$. В самом деле, $Z\nor G$ и $Z\le E$, значит,
$$Z\se \bigcap_{x\in G} E^x=\ker((1_E)^G)$$
в силу \ref{ker indh}. Но тогда $Z$
содержится в ядре всякой неприводимой компоненты характера $(1_E)^G$, в частности $Z\se\ker\x$. Поскольку
$Z\le A$ и $\l$ --- компонента $\x_{A^{\vphantom{A^a}}}$, имеем $Z\se\ker\l$.

Теперь из предложения \ref{vs tb}, применённого к группе $G=AP$, следует, что существует элемент $x\in G$
такой, что $\l^x\ne \l$. В частности, группа инерции $I=\II_G(\l)$ является собственной подгруппой в $G$.
Из теоремы Клиффорда \ref{tcl char}$(ii)$ следует, что существует характер $\vf\in \irr(I)$ такой,
что $\x=\vf^G$. Таким образом требуемое установлено и теорема Брауэра доказана.
\end{proof}

\section{Обыкновенные характеры нулевого \texorpdfstring{$p$}{p}-дефекта}

Зафиксируем до конца этого раздела некоторое простое число $p$.

\begin{opr}\label{opch}
Пусть $n$ --- натуральное число.
Назовём \glsadd{iPPrtNum}\mem{$p$-частью}
числа $n$ максимальную степень $p$, делящую $n$, и обозначим её через \gls{np}, а \glsadd{iPPrPrtNum}\mem{$p'$-частью} числа $n$ ---
величину $n/n_p$, которую будем обозначать через \gls{npp}.
\end{opr}

По определению неотрицательное целое число $d$ является
\glsadd{iPDefOrdChr}\mem{$p$-дефектом} \label{pdef opr} характера $\x\in \irr(G)$, если
$$p^d=\left(\frac{|G|}{\deg\x}\right)_{\!p}.$$

Мы используем теорему Брауэра для изучения значений неприводимых обыкновенных характеров $p$-дефекта $0$,
т.\,е. таких $\x\in \irr(G)$, для которых $\x(1)_p=|G|_p$.

\begin{opr}\label{opr dot x}
Пусть $\x \in \irr(G)$ и $d$ --- $p$-дефект характера $\x$.
Определим отображение $\dot\x:G\to\CC$ по правилу
$$
\dot\x (g)=\left\{\ba{ll}
p^d\x(g),& \mbox{если}\ \ p\nmid |g|;\\
0,& \mbox{в противном случае}.
\ea\right.
$$
для любого $g\in G$.
\end{opr}
Легко видеть, что $\dot\x\in \cf(G)$. Покажем, что в действительности $\dot\x$ является
обобщённым характером.

\begin{pre} \label{dot x} Если $\x\in \irr(G)$, то $\dot\x\in \gch(G)$.
\end{pre}
\begin{proof} Пусть $E$ --- элементарная подгруппа из $G$ и $\theta\in \irr(E)$. По теореме \ref{thm char br}
достаточно проверить,
что $(\dot\x_{E^{\vphantom{A^a}}},\theta)_{\mbox{}_E}\in \ZZ$, где  $\dot\x_{E^{\vphantom{A^a}}}$ обозначает $(\dot\x)_{E^{\vphantom{A^a}}}$.

Поскольку $E$ нильпотентна, мы можем
записать $E=P\times Q$, где $P$ --- $p$-группа и $p\nmid |Q|$. Пусть $x\in E$. Если $p\nmid |x|$, то $x\in Q$.
Поэтому функция $\dot\x$ тождественно равна нулю на $E \setminus Q$ и $(\dot\x)_{Q^{\vphantom{A^a}}}=p^d\x_{Q^{\vphantom{A^a}}}$. Имеем
$$
(\dot\x_{E^{\vphantom{A^a}}},\theta)_{\mbox{}_E}=\frac{1}{|E|}\sum_{x\in E}\dot\x(x)\ov{\th(x)}
=\frac{p^d}{|E|}\sum_{x\in Q}\x(x)\ov{\th(x)}=\frac{p^d}{|P|}(\x_{Q^{\vphantom{A^a}}},\theta_Q)_{\mbox{}_Q},
\myeqno\label{xqt}
$$
поскольку $|E|=|P||Q|$.

Рассмотрим центральный характер $\om_\x$. Из \ref{ch val} следует, что если $g\in G$, то
$$\x(g)=\frac{\om_\x(\wh{g^G})\x(1)}{|g^G|}=\frac{\om_\x(\wh{g^G})\x(1)|\C_G(g)|}{|G|}.$$
 Тогда
$$
(\dot\x_{E^{\vphantom{A^a}}},\theta)_{\mbox{}_E}=\frac{p^d}{|E|}\sum_{x\in Q}\x(x)\ov{\th(x)}
=\frac{p^d\x(1)}{|E||G|}\sum_{x\in Q}\om_\x(\wh{x^G})|\C_G(x)|\ov{\th(x)}
=\frac{p^d\x(1)}{|Q||G|}\sum_{x\in Q}\om_\x(\wh{x^G})|\C_G(x):P|\ov{\th(x)},
\myeqno\label{xqt2}
$$
где последнее равенство следует из того, что $P\le \C_G(x)$ для любого $x\in Q$.

Из \ref{xqt} следует, во-первых, что $(\dot\x_{E^{\vphantom{A^a}}},\theta)_{\mbox{}_E}\in \QQ$, а во-вторых, что
$|P|(\dot\x_{E^{\vphantom{A^a}}},\theta)_{\mbox{}_E}\in \ZZ$. Из \ref{xqt2} следует, что величина
$$
\frac{|Q||G|}{p^d\x(1)}(\dot\x_{E^{\vphantom{A^a}}},\theta)_{\mbox{}_E}
$$
также является целой, т.\,к. с одной стороны она рациональна, а с другой --- целое алгебраическое число в силу
\ref{har prop}$(iii)$ и \ref{ch ac}. Однако по условию число  $a=|Q|\frac{|G|}{p^d\x(1)}$ является целым числом, не
делящимся на $p$, а $|P|$ --- степень $p$. Поэтому существуют числа $u,v\in \ZZ$ такие, что $ua+v|P|=1$. Умножив
обе части этого равенства на $(\dot\x_{E^{\vphantom{A^a}}},\theta)_{\mbox{}_E}$ и учитывая вышесказанное,
получаем $(\dot\x_{E^{\vphantom{A^a}}},\theta)_{\mbox{}_E}\in \ZZ$, что и требовалось доказать.
\end{proof}

\begin{pre} \label{pdef 0} Пусть $\x\in \irr(G)$ --- характер $p$-дефекта $0$. Тогда $\x(g)=0$ для любого
элемента $g\in G$ порядка, делящегося на $p$.
\end{pre}
\begin{proof} Из предложения \ref{dot x} следует, что  $\dot\x$ является обобщённым характером.
В частности, $(\dot\x,\x)_{\mbox{}_G}\in \ZZ.$ С другой стороны $\dot\x(g)=\x(g)$ при $p\nmid |g|$ т.\,к.
$p$-дефект характера $\x$ равен $0$, и поэтому
$$
(\dot\x,\x)_{\mbox{}_G}=\frac{1}{|G|}\sum_{g\in G,\, p\,\nmid\, |g|} |\x(g)|^2.
$$
Отсюда следует, что $0<(\dot\x,\x)_{\mbox{}_G}\le  (\x,\x)_{\mbox{}_G}=1$ и, значит,
$(\dot\x,\x)_{\mbox{}_G}=(\x,\x)_{\mbox{}_G}$. Но тогда
$$
0= (\x-\dot\x,\x)_{\mbox{}_G}=\frac{1}{|G|}\sum_{g\in G,\, p\,\mid\, |g|} |\x(g)|^2,
$$
откуда вытекает, что $\x(g)=0$ для всех $g\in G$ таких, что $p$ делит $|g|$.
\end{proof}

\chapter{Модулярные представления и брауэровы характеры}

Начиная с этого раздела, мы будем изучать модулярные представления конечных групп, т.\,е.
представления над полем простой характеристики.

\section{Брауэровы характеры}

Пусть $\ov{\ZZ}$ --- кольцо целых алгебраических чисел. Всюду далее мы фиксируем простое число $p$
и некоторый максимальный идеал $M$ кольца $\ov{\ZZ}$, содержащий идеал $p\ov{\ZZ}$.  Положим $F=\ov{\ZZ}/M$ и
пусть
$$
^{\gls{aaa}}:\ov{\ZZ}\to F
$$
--- естественный эпиморфизм колец.

\textuprn{Показать, что $F$ --- поле характеристики $p$.}

\begin{pre} \label{m int z} $M\cap \ZZ=p\ZZ$. В частности, $\ZZ^\st\cong \ZZ/p\ZZ=\FF_p$.
\end{pre}
\textupl{m int z prf}{Доказать предложение \ref{m int z}.}

Положим $$U=\{\z\in \CC\,|\,\z^m=1\ \ \mbox{для некоторого}\  \ m\in \ZZ\ \ \mbox{такого,
что}\
\ p\nmid m\}.$$ Другими словами, $U$ --- группа корней из $1$ степени, взаимно простой с $p$.
Ясно, что $U\se\ov{\ZZ}$.

\begin{pre} \label{mul gr azp} Имеют место следующие утверждения.
\begin{list}{{\rm(}{\it\roman{enumi}\/}{\rm)}}
{\usecounter{enumi}\setlength{\parsep}{2pt}\setlength{\topsep}{5pt}\setlength{\labelwidth}{23pt}}
\item Ограничение отображения $^\st$ на $U$ является изоморфизмом групп $U\to F^\times$.
\item Поле $F$ совпадает с алгебраическим замыканием своего простого подполя $\FF_p$.
\end{list}
\end{pre}
\begin{proof} Будем обозначать ограничение отображения $^\st$ на $U$ тем же символом.
Покажем инъективность гомоморфизма групп $^\st:U\to F^\times$. Пусть $1\ne \z\in U$. Тогда
$m=|\z|>1$ и $p\nmid m$. Допустим, что $\z^\st=1$. В кольце $\ov{\ZZ}[x]$ имеет место тождество
$$
1+x+\ld+x^{m-1}=\frac{x^m-1}{x-1}=\prod_{i=1}^{m-1}(x-\z^i).
$$
Сравнив правую и левую часть при $x=1$, получим, что в кольце $\ov{\ZZ}$ элемент $1-\z$ делит $m$.
Но раз $\z^\st=1$, то $(1-\z)^\st=0$ и, значит, $m^\st=0$. Тогда $m\in M$. В силу \ref{m int
z} число $m$ кратно $p$. Противоречие. Значит, отображение $^\st$ инъективно на $U$.

 Пусть $f\in F^\times$. Выберем элемент $\a\in \ov{\ZZ}$ такой, что
$\a^\st=f$. Так как $\a$ --- целое алгебраическое число, то
$$
\a^n+b_1\a^{n-1}+b_2\a^{n-2}+\ld+b_n=0
$$
для некоторых $b_1,b_2,\ld,b_n\in \ZZ$. Тогда
$$
f^n+b_1^\st f^{n-1}+b_2^\st f^{n-2}+\ld+b_n^\st=0
$$
и $f$ является алгебраическим над $\FF_p$ в силу \ref{m int z}.
То есть расширение $F \ge \FF_p$ является алгебраическим. Для завершения доказательства
достаточно показать, что если $K \ge F$
--- алгебраическое расширение полей, то $K^\times \se U^\st$. Отсюда будет следовать
сюръективность отображения $^\st$ и алгебраическая замкнутость поля $F$.

\normalmarginpar

Пусть $k\in K^\times$. Тогда элемент $k$ алгебраический над $\FF_p$ в силу \ref{alg trans}. То
есть $k$ является ненулевым элементом конечного поля $\FF_p(k)$ характеристики $p$. Поэтому
$k^s=1$, где $s=|\FF_p(k)|-1$ и, значит, $p\nmid s$. Но мы видели, что в $U^\st$ уже есть $s$
различных корней многочлена $x^s-1$ --- это степени элемента $\xi^\st$, где $\xi\in U$ ---
элемент порядка $s$. Значит,
$k$ совпадает с одним из этих корней, т.\,е. лежит в $U^\st$.
\end{proof}

Элемент $g\in G$  называется \glsadd{iElmPReg}\mem{$p$-регулярным} или
\glsadd{iPPrElm}\mem{$p'$-элементом}, если $p\nmid |g|$. Класс сопряжённости $g^G$
произвольного $p$-регулярного элемента $g$ мы будем также называть \glsadd{iClsCnjPReg}\mem{$p$-регулярным}. Ясно, что $p$-регулярный класс состоит из
$p$-регулярных элементов. Обозначим через
\gls{Gpp} множество всех $p$-регулярных элементов группы $G$, а через
\gls{KGpp} --- множество всех её $p$-регулярных классов.

\begin{lem} \label{lem ppt} Пусть $g\in G$. Тогда существуют однозначно определённые элементы
$g_{p'},g_p\in G$ такие, что $g_{p'}$ является $p$-регулярным, $g_p$ является $p$-элементом и
$\,g=g_pg_{p'}=g_{p'}g_p$. При этом $g_{p'},g_p\in \la g\ra$ и имеет место
равенство $\C_G(g)=\C_G(g_p)\cap\C_G(g_{p'})$.
\end{lem}

\textupl{lem ppt prf}{Доказать лемму \ref{lem ppt}.}

Элементы \gls{gp} и \gls{gpp} из леммы \ref{lem ppt} называются, соответственно,
\glsadd{iPPrtElm}\mem{$p$-частью} и \glsadd{iPPrPrtElm}\mem{$p'$-частью} элемента $g$.

Пусть $\pi$ --- подмножество множества всех простых чисел. Через $\pi'$ обозначается множество
простых чисел, не принадлежащих $\pi$. Элемент $g\in G$ называется  \glsadd{iPiElm}\mem{$\pi$-элементом},
если все простые делители его порядка принадлежат $\pi$. Множество
всех $\pi$-элементов группы $G$ обозначается через \gls{Gspi}.

Следующее упражнение обобщает лемму \ref{lem ppt}.

\textuprn{Пусть $g\in G$ и  $\pi$ --- произвольное множество простых чисел. Тогда существуют однозначно
определённые элементы $g_{\pi},g_{\pi'}\in G$ такие, что $g_{\pi}\in G_\pi$, $g_{\pi'}\in
G_{\pi'}$ и $g=g_{\pi}g_{\pi'}=g_{\pi'}g_{\pi}$. При этом $g_{\pi},g_{\pi'}\in \la g\ra$.}

\begin{opr} Пусть $\X:G\to \GL_n(F)$ --- $F$-представление группы $G$. Пусть $g\in G_{p'}$. По
\ref{mul gr azp}$(i)$ все характеристические значения матрицы $\X(g)$ имеют вид $\z_1^\st,\ld,\z_n^\st$
для однозначно определённых элементов $\z_1,\ld,\z_n\in U$. \glsadd{iChrBr}\mem{Брауэровым} или
\glsadd{iChrPMod}\mem{$p$-модулярным характером} $F$-представления $\X$ называется
отображение $\vf:G_{p'}\to\CC$, заданное для $g\in G_{p'}$ правилом
$$
\vf(g)=\z_1+\ld+\z_n.
$$
\end{opr}

\textzam{Значения $\vf(g)$, вообще говоря, зависят от выбранного идеала $M$.}

Функцию $G_{p'}\to \CC$, принимающую постоянное значение на каждом $p$-регулярном классе,
будем называть \glsadd{iFuncClsPReg}\mem{$p$-регулярной классовой
функцией}. Множество всех $p$-регулярных классовых функций группы $G$ будем обозначать через
\gls{cfGpp}. На этом множестве, как и на $\cf(G)$, задаётся структура $\CC$-алгебры
относительно естественных операций.

Определим для произвольной функции $\th\in\cf(G_{p'})$ \mem{комплексно-сопряжённую} функцию
$\ov{\th}$ по правилу $\ov{\th}(g)=\ov{\th(g)}$ для всех $g\in G_{p'}$. Если $H\le G$ и
$\th\in \cf(G_{p'})$, то обозначим через $\th_H$ ограничение $\th$ на $H_{p'}$. Ясно, что
$\th_H\in \cf(H_{p'})$.

\begin{pre}\label{br prop} Пусть $\vf$ --- брауэров характер $F$-представления $\X$ группы $G$.
\begin{list}{{\rm(}{\it\roman{enumi}\/}{\rm)}}
{\usecounter{enumi}\setlength{\parsep}{2pt}\setlength{\topsep}{5pt}\setlength{\labelwidth}{23pt}}
\item $\vf$ --- брауэров характер любого $F$-представления, эквивалентного  $\X$ (при
фиксированном выборе идеала $M$).
\item $\vf\in\cf(G_{p'})$.
\item $\vf(g^{-1})=\ov{\vf(g)}$ для любого $g\in G_{p'}$.
\item $\ov{\vf}$ является брауэровым характером контрагредиентного $F$-представления $\X^*$.
\item Если $H\le G$, то $\vf_{H^{\vphantom{A^a}}}$ --- брауэров характер группы $H$.
\end{list}
\end{pre}

\textupl{br prop prf}{Доказать предложение \ref{br prop}.}

\textzam{Следует подчеркнуть, что для произвольной $p$-регулярной классовой функции $\vf\in \cf
(G_{p'})$ в общем случае некорректно спрашивать, будет ли $\vf$ брауэровым характером группы
$G$ или нет, поскольку ответ на этот вопрос может существенно зависеть от выбора максимального
идеала $M$ кольца $\ov{\ZZ}$, содержащего простое число $p$. Есть примеры, показывающие, что
брауэров характер $\vf$ при одном выборе $M$ может перестать быть таковым, если идеал $M$
изменить. Кроме того, если $M$ фиксирован, $\vf$
--- брауэров характер и $\s$ --- автоморфизм поля $\CC$, то классовая функция $\vf^\s\in
\cf(G_{p'})$, заданная правилом, $\vf^\s(g)=\vf(g)^\s$ при всех $g\in G_{p'}$ может не
быть брауэровым характером группы $G$ (см. предложение \ref{mm}).}

Как и для $F$-характеров, \glsadd{iDegChrBr}\mem{степенью} брауэрова
характера $\vf$ назовём степень представления $\X$. Ясно, что она равна $\vf(1)$. Брауэров
характер главного $F$-представления называется \glsadd{iChrBrPrnc}\mem{главным}. Будем обозначать его через \gls{1Gpp}.

\begin{pre} \label{sum pr} Пусть $\X$ и $\Y$ --- $F$-представления группы $G$ с брауэровыми
характерами $\vf$ и $\psi$. Тогда брауэровы характеры $F$-представлений $\X\oplus\Y$ и
$\X\ot\Y$ равны, соответственно, $\vf+\psi$ и $\vf\psi$. В частности, множество брауэровых
характеров  группы $G$ замкнуто относительно сложения и умножения.
\end{pre}

\textupr{Доказать предложение \ref{sum pr}. \uk{Воспользоваться предложением \ref{tnz prop}.}}

Назовём брауэров характер $\vf$ представления $\X$
\glsadd{iChrBrIrr}\mem{неприводимым}, если неприводимо $\X$. Множество всех неприводимых брауэровых характеров
группы $G$ относительно максимального идеала $M$ будем обозначать через \gls{IBrMlGr}.
Если ясно, о каком идеале $M$ идёт речь, но нужно подчеркнуть, что характеристика поля $F$ равна $p$, то будем обозначать
это множество через \gls{IBrplGr}. Но в большинстве случаев будем просто использовать обозначение
\gls{IBrlGr}. Из \ref{fin irr rep} вытекает, что это множество конечно.

\begin{pre} \label{br sum} Справедливы следующие утверждения.
\begin{list}{{\rm(}{\it\roman{enumi}\/}{\rm)}}
{\usecounter{enumi}\setlength{\parsep}{2pt}\setlength{\topsep}{5pt}\setlength{\labelwidth}{23pt}}
\item Если два $F$-представления $\X$ и $\Y$ группы $G$ имеют один и тот же набор
неприводимых компонент (с учётом кратностей), то их брауэровы характеры совпадают.
\item Если $\psi\in \cf(G_{p'})$ является брауэровым характером, то
\end{list}
\vspace{-2.5\topsep}
$$\psi=\sum_{\vf\in \iBr(G)}a_\vf\vf, \myeqno\label{br dec}$$
\begin{list}{{\rm(}{\it\roman{enumi}\/}{\rm)}}
{\usecounter{enumi}\setlength{\parsep}{2pt}\setlength{\topsep}{5pt}\setlength{\labelwidth}{23pt}}
\item[] где $a_\vf$ --- неотрицательные целые числа.
\end{list}
\end{pre}
\begin{proof} Поскольку $\X$ и $\Y$ эквивалентны блочно-верхнетреугольным представлениям с
одинаковыми (с точностью до перестановки) неприводимыми компонентами на диагонали, то
брауэровы характеры $\X$ и $\Y$ равны сумме брауэровых характеров неприводимых компонент.
Отсюда также следует, что любой брауэров характер является суммой элементов из $\iBr(G)$ с
неотрицательными целыми коэффициентами.
\end{proof}

Как мы вскоре увидим (см. следствие \ref{cor br pr}), утверждение \ref{br sum}$(i)$ можно
обратить.

Следующее предложение объясняет, почему достаточно определять брауэровы характеры только на
множестве $p$-регулярных элементов: зная брауэров характер представления $\X$, можно
восстановить его $F$-характер.

\begin{pre} \label{br preg} Имеем
\begin{list}{{\rm(}{\it\roman{enumi}\/}{\rm)}}
{\usecounter{enumi}\setlength{\parsep}{2pt}\setlength{\topsep}{5pt}\setlength{\labelwidth}{23pt}}
\item пусть $\X$ --- $F$-представление группы $G$ с $F$-характером $\x$ и
брауэровым характером $\vf$. Тогда $\x(g)=\x(g_{p'})=\vf(g_{p'})^\st$ любого $g\in G$;
\item пусть $\x\in\ch(G)$. Тогда $\x(g)^\st=\x(g_{p'})^\st$ для любого $g\in G$;
\item $\irr_F(G)=\{\vf^\st\,|\,\vf\in
\iBr(G)\}$, где $\vf^\st(g)=\vf(g_{p'})^\st$ для всех $g\in G$.
\end{list}
\end{pre}
\begin{proof} $(i)$ Пусть $g\in G$. Ограничим представление $\X$ на циклическую группу $\la
g\ra$. Пусть $\X_i$ --- неприводимые компоненты этого ограничения. По \ref{cent irr} все
$\X_i$ одномерны. Так как $g_p\,,g_{p'}\in \la g\ra$, то $\X_i(g)=\X_i(g_p)\X_i(g_{p'})$. Но
$\X_i(g_p)\in F^\times$ --- $p$-элемент, откуда следует, что $\X_i(g_p)=1$. Поэтому
$\x(g)=\sum_i\X_i(g)=\sum_i\X_i(g_{p'})=\x(g_{p'})$.

$(ii)$ Пусть $\X$ --- обыкновенное представление с характером $\x$. Рассуждая, как в $(i)$, получаем
$\x(g)=\sum_i\X_i(g_p)\X_i(g_{p'})$, где $\X_i$ --- неприводимые компоненты ограничения $\X$ на $\la g \ra$.
Так как $\X_i(g_p)\in \ov{\ZZ}$ --- $p$-элемент, то $\X_i(g_p)^\st=1$. Значит,
$\x(g)^\st=\sum_i\X_i(g_{p'})^\st=\x(g_{p'})^\st$.

$(iii)$ Из $(i)$ вытекает, что $\vf^\st$ совпадает с $F$-характером любого $F$-представления,
брауэровым характером которого является $\vf$.
\end{proof}

\begin{cor} \label{br dif} Имеем
\begin{list}{{\rm(}{\it\roman{enumi}\/}{\rm)}}
{\usecounter{enumi}\setlength{\parsep}{2pt}\setlength{\topsep}{5pt}\setlength{\labelwidth}{23pt}}
\item брауэровы характеры неэквивалентных неприводимых $F$-представлений различны;
\item $|\iBr(G)|=|\M(FG)|$.
\end{list}
\end{cor}
\begin{proof} $(i)$ Если брауэровы характеры двух неприводимых $F$-представлений совпадают, то в силу \ref{br preg}$(i)$
совпадают и их $F$-характеры. Тогда сами представления эквивалентны по \ref{rep trace}$(iii)$.

$(ii)$ Следует из $(i)$.
\end{proof}

Имеет место даже более сильный результат, ср. \ref{lin ind fch}.

\begin{pre} \label{br lin ind} Пусть
$\vf_1,\ld,\vf_s\in\iBr(G)$ попарно различны.
\begin{list}{{\rm(}{\it\roman{enumi}\/}{\rm)}}
{\usecounter{enumi}\setlength{\parsep}{2pt}\setlength{\topsep}{5pt}\setlength{\labelwidth}{23pt}}
\item Если $\a_1,\ld,\a_s\in \ov{\ZZ}$ такие,
что $\sum_i \a_i\vf_i(x)\in M$ для всех $x\in G_{p'}$, то $\a_i\in M$ для всех $i$.
\item Характеры $\vf_1,\ld,\vf_s$  линейно независимо над $\CC$. В частности, всё множество $\iBr(G)$ линейно
независимо над $\CC$.
\end{list}
\end{pre}
\begin{proof}
$(i)$ Пусть $\X_i$ --- неприводимое $F$-представление с брауэровым характером $\vf_i$ и его
$F$-характер равен $\psi_i$. Так как $\sum_i \a_i\vf_i(x)\in M$ для всех $x\in G_{p'}$, то из
\ref{br preg}$(i)$ следует, что $\sum_i \a_i^\st\psi_i(x)=0$ для всех $x\in G$. Однако мы
знаем из \ref{rep trace}$(iii)$, что $F$-характеры $\psi_i$ линейно независимы, т.\,е.
$a_i^\star=0$ и, значит, $a_i\in M$.

$(ii)$ Достаточно доказать линейную независимость множества $\vf_1,\ld,\vf_s$ над $\ov{\QQ}$, поскольку отсюда
будет следовать\footnote{Если $K$ --- расширение $\ov{\QQ}$ с базисом $B$, то любую $K$-линейную комбинацию $\vf$
характеров $\vf_1,\ld,\vf_s$ можно переписать
как конечную линейную комбинацию элементов из $B$, коэффициенты которой --- $\ov{\QQ}$-линейные комбинации характеров $\vf_1,\ld,\vf_s$. Поскольку
значения этих характеров лежат в $\ov{\QQ}$, из того, что $\vf=0$, следует, что все такие коэффициенты равны нулю, а значит из
$\ov{\QQ}$-линейной независимости характеров $\vf_1,\ld,\vf_s$ будет следовать, что $K$-коэффициенты $\vf$ тоже равны нулю.}
линейная независимость над любым расширением поля $\ov{\QQ}$ и, в частности, над $\CC$.

\normalmarginpar

Пусть $\g_1,\ld,\g_s\in \ov{\QQ}$ такие, что $\sum_i \g_i\vf_i=0$. Предположим, что не все $\g_i$ равны нулю.
Из следствия \ref{req lem} вытекает существование элемента $\b\in \ov{\QQ}$ такого, что все произведения
$\b\g_i$ лежат в $\ov{\ZZ}$, но не все лежат в $M$. Поскольку $\sum_i (\b\g_i)\vf_i=0$, получаем противоречие с $(i)$.
\end{proof}

Из линейной независимости неприводимых брауэровых характеров вытекает, что коэффициенты
$a_\vf$ в разложении \ref{br dec} произвольного брауэрова характера $\psi$ по характерам из $\iBr(G)$ определяются
однозначно. Будем называть $a_\vf$ \mem{кратностью вхождения} неприводимого брауэрова
характера $\vf$ в характер $\psi$. Отметим, что если $V$ и $W$ --- $FG$-модули,
соответствующие представления которых имеют брауэровы характеры $\psi$ и $\vf$, то $a_\vf$
совпадает с числом композиционных факторов модуля $V$, изоморфных $W$. Переформулировав это на
языке представлений, получаем следующее обращение утверждения \ref{br sum}$(i)$.

\begin{cor} \label{cor br pr} Если два $F$-представления группы $G$ имеют одинаковые брауэровы
характеры, то у них один и тот же набор неприводимых компонент (с учётом кратностей).
\end{cor}

\textext{Равенство брауэровых характеров, вообще говоря, не влечёт
эквивалентность $F$-представлений. Неэквивалентные 
$F$-представле\-ния 
циклической группы
порядка $2$, приведённые на с. \pageref{neqch}, не только имеют равные $F$-характеры, но,
как легко видеть, и равные брауэровы характеры (при любом выборе идеала $M$).}

Ещё одним следствием линейной независимости является тот факт, что брауэров характер неприводим тогда и
только тогда, когда он не представим в виде суммы $\vf+\psi$ для некоторых брауэровых характеров $\vf$ и
$\psi$. В частности, в силу \ref{br prop}$(iv)$ получаем

\begin{pre} \label{spr br} $\vf\in \iBr(G)$ тогда и только тогда, когда $\ov{\vf}\in \iBr(G)$.
\end{pre}
Отметим, что этот факт также следует из \ref{kon npr}.

\section{Связь с обыкновенными характерами}

Брауэровы характеры важны тем, что осуществляют связь между обыкновенными характерами группы
$G$ и её $F$-характерами в характеристике $p$. Чтобы объяснить это подробнее, нам потребуется
ещё несколько определений и вспомогательных утверждений.

Пусть $K$ --- подкольцо поля $F$ и $\X$ --- $F$-представление группы $G$. По аналогии 
с определением \ref{moz zap} будем говорить, что $\X$ \glsadd{iRepGrZapR}\mem{может
быть записано над} $K$, если оно эквивалентно некоторому $F$-представлению с матричными элементами из $K$.

Определим подкольцо поля $\ov\QQ$ и идеал в этом подкольце следующим образом.

$$
\gls{Ztil}=\left\{\frac{\a}{\b}\,\Big|\ \a\in \ov{\ZZ},\ \b\in \ov{\ZZ}\setminus M\right\},\qquad
\gls{Mtil}=\left\{\frac{\a}{\b}\,\Big|\ \a\in M,\ \b\in \ov{\ZZ}\setminus M\right\}.
$$

\textuprn{Показать, что $\wt{Z}$ --- локальное кольцо с максимальным идеалом $\wtm$. В частности, $\J(\wt{Z})=\wtm$.}

\textuprn{Описать кольцо $\wt{Z}\cap \QQ$ и его идеал $\wtm\cap \QQ$.}

Отображение $^\st$ может быть естественно продолжено на кольцо $\wt{Z}$.
Если $\a\in \ov{\ZZ}$ и
$\b\in
\ov{\ZZ}\setminus M$, то
$$
\left(\frac{\a}{\b}\right)^\st=\frac{\a^\st}{\b^\st}.
$$

\textuprn{Показать, что введённое отображение $^\st:\wt{Z}\to F$ определено корректно и является гомоморфизмом
колец с ядром $\wtm$. В частности, $\ov{\ZZ}\cap\wtm=M$.}

\begin{lem} \label{lem qzm} Пусть $V$ --- векторное пространство над полем $\ov{\QQ}$ и $W$
--- конечно порождённый $\wt{Z}$-подмодуль $\wt{Z}$-модуля $V$. Тогда существуют элементы $w_1,\ld,w_s\in
W$, которые линейно независимы над $\ov{\QQ}$ и такие, что $W=\wt{Z} w_1\oplus\ld\oplus\wt{Z} w_s$.
\end{lem}
\begin{proof} Пусть
$$\t : W\to W/\wtm W$$
--- естественный эпиморфизм $\wt{Z}$-модулей. Фактормодуль $W\t$ можно рассматривать как векторное
пространство над полем $F=\wt{Z}/\wtm $, положив $\a^\st(w\t)=(\a w)\t$ для всех $\a\in
\wt{Z}$. Поскольку $\wt{Z}$-модуль $W$ конечно порождён, размерность
$\dim_F W\t$ конечна. Пусть $w_1\t,\ld,w_s\t$ --- базис пространства $W\t$, где $w_i\in W$. Положим $T=\wt{Z}
w_1+\ld +\wt{Z} w_s$. Тогда $T$ --- $\wt{Z}$-подмодуль модуля $W$, и имеет место равенство $W=T+\wtm W$. Из
следствия \ref{cor nak lok} вытекает, что $W=T$.

Остаётся проверить, что элементы $w_1,\ld,w_s$ линейно независимы над $\ov{\QQ}$. Если $\sum_i \g_i w_i=0$, где
$\g_i\in \ov{\QQ}$ не все равны нулю, то по следствию \ref{req lem} существует элемент $\b\in \ov{\QQ}$ такой,
что все произведения $\b\g_i$ лежат в $\ov{\ZZ}$, но не все лежат в $M$. Тогда $\sum_i \b\g_i w_i=0$. Применив
гомоморфизм $\t$, получим $\sum_i (\b\g_i)^\st(w_i\t)=0$, где не все коэффициенты $(\b\g_i)^\st$ равны нулю.
Это противоречит линейной независимости над $F$ элементов $w_i\t$.
\end{proof}

\begin{pre} \label{kom alg} Любое обыкновенное представление группы $G$
может быть записано над кольцом~$\wt{Z}$.
\end{pre}
\begin{proof} Пусть $\X$ --- обыкновенное представление группы $G$.
В силу следствия \ref{cor c q} можно считать, что $\X$ записано над $\ov{\QQ}$. Пусть $V$ ---
$\ov{\QQ}G$-модуль, соответствующий $\X$, и $v_1,\ld,v_n$
--- его базис. Рассмотрим $\wt{Z}$-подмодуль $W$ модуля $V$, порождённый всевозможными
произведениями $v_ig$, где $i=1,\ld,n$ и $g\in G$. Тогда $W$ является конечно порождённым $\wt{Z}$-модулем,
который инвариантен относительно действия группы $G$. По лемме \ref{lem qzm} существуют линейно независимые над
$\ov{\QQ}$ элементы $w_1,\ld,w_s\in W$ такие, что $W=\wt{Z}
w_1\oplus\ld\oplus\wt{Z} w_s$. В частности, $s\le n$. Обратно, $s\ge n$, поскольку $W$ содержит
$\ov{\QQ}$-базис $v_1,\ld,v_n$. Значит, $w_1,\ld,w_n$ --- базис пространства $V$, и $\ov\QQ$-представление
$\Y$, соответствующее модулю $V$ в этом базисе, записано над $\wt{Z}$ ввиду того, что $w_ig$ является
$\wt{Z}$-линейной комбинацией элементов $w_j$ для всех $g\in G$.
\end{proof}

\textext{Отметим, что если $\X$ --- представление группы $G$ над полем алгебраических
чисел $K$, то $\X$ не всегда может быть записано над кольцом целых величин поля $K$.
Однако для $\X$ справедливо следующее усиление предложения \ref{kom alg}.
Существует конечное расширение $E$ поля $K$ такое, что $\X$ может быть записано над кольцом целых величин поля $E$. Доказательство
см. в \cite[Satz V.12.5]{hup}.}

Продолжим отображение $^\st$ на кольца $\MM_n(\wt{Z})$ и $\wt{Z}[x]$ следующим образом. Пусть $A\in
\MM_n(\wt{Z})$ и $A=(\a_{ij})$. Положим $\gls{Ast}=(\a_{ij}^\st)$. Ясно, что отображение
$$
^\st: \MM_n(\wt{Z})\to \MM_n(F)
$$
является гомоморфизмом колец. Кроме того, $\det(A^\st)=\det(A)^\st$ для всех $A\in
\MM_n(\wt{Z})$. Мы также можем естественно продолжить отображение $^\st$ до гомоморфизма колец
многочленов
$$
^\st: \wt{Z}[x]\to F[x],
$$
причём заметим, что если корни $\z_i$ некоторого многочлена $f\in \wt{Z}[x]$ лежат в $\wt{Z}$,
то корнями многочлена $f^\st\in F[x]$ будут элементы $\z_i^\st$ поля $F$.

Пусть $\X$ --- обыкновенное представление группы $G$, записанное над $\wt{Z}$. Обозначим через
\gls{Xst} отображение $g\mapsto \X(g)^\st$. Тогда $\X^\st$ является
$F$-представлением группы $G$.

Для произвольной классовой функции $\th\in \cf(G)$ обозначим через
\gls{thh} ограничение $\th$ на $G_{p'}$. Легко видеть, что $\hat{\th}\in
\cf(G_{p'})$ для всех $\th\in \cf(G)$, и отображение
$\,\hat{\mbox{}}:\cf(G)\to\cf(G_{p'})$ является $\CC$-линейным и сюръективным.

\begin{pre} \label{ord br} Пусть $\x\in\ch(G)$.
\begin{list}{{\rm(}{\it\roman{enumi}\/}{\rm)}}
{\usecounter{enumi}\setlength{\parsep}{2pt}\setlength{\topsep}{5pt}\setlength{\labelwidth}{23pt}}
\item $\hat\x$ --- брауэров характер группы $G$ (для любого выбора максимального идеала $M$).
\item Если $\X$ --- обыкновенное представление
с характером $\x$, записанное над $\wt{Z}$, то $\hat \x$ является брауэровым характером $F$-представления $\X^\st$.
\end{list}
\end{pre}
\begin{proof} Докажем сразу оба утверждения.
Пусть $\X$ --- обыкновенное представление группы $G$ с характером $\x$. В силу
\ref{kom alg} можно считать, что $\X$ записано над $\wt{Z}$. Пусть  $\vf$ обозначает
брауэров характер $F$-представления $\X^\st$. Достаточно показать, что $\vf=\hat\x$ . Пусть
$g\in G_{p'}$ и $f\in
\wt{Z}[x]$
--- характеристический многочлен матрицы $\X(g)$. Имеем
$$f^\st(x)=\det(\X(g)-x\I)^\st=\det(\X(g)^\st-x\I)=\det(\X^\st(g)-x\I)$$
и, значит, $f^\st$ --- характеристический многочлен матрицы $\X^\st(g)$.

Заметим, что корни $\z_i$ многочлена $f$ являются корнями из $1$ степени $|g|$, т.\,е. лежат в
$U\se \ov{\ZZ}\se \wt{Z}$. Поэтому элементами $\z_i^\star$ исчерпываются все корни многочлена
$f^\star$ и, следовательно, характеристическими значениями матрицы $\X^\st(g)$. Таким образом,
по определению брауэрова характера $\vf(g)=\sum_i
\z_i=\x(g)=\hat{\x}(g)$. В силу произвольности $g$ получаем $\vf=\hat\x$.

Поскольку $p$-регулярная классовая функция $\hat\x$ не зависит от максимального идеала $M$ и
приведённое выше рассуждение верно для любого выбора $M$, то $\hat\x$ является брауэровым
характером группы $G$ при любом $M$.
\end{proof}

Следует подчеркнуть, что если обыкновенные представления $\X$ и $\Y$ эквивалентны и записаны над
$\wt{Z}$, то $F$-представления $\X^\st$ и $\Y^\st$, вообще говоря, могут не быть эквивалентны. Например,
пусть обыкновенные представления $\X$ и $\Y$ циклической группы $\la a\,|\, a^2=1\ra$ заданы по правилу
$$
\X(a)=\left(
             \begin{array}{rr}
               1 & 0 \\
               0 & -1 \\
             \end{array}
           \right),\qquad
\Y(a)=\left(
             \begin{array}{rr}
               1 & 1 \\
               0 & -1 \\
             \end{array}
           \right).
$$
Эти представления эквивалентны, поскольку их $\CC$-характеры равны. Однако при $p=2$
соответствующие $F$ представления
$$
\X^\st(a)=\left(
             \begin{array}{rr}
               1 & 0 \\
               0 & 1 \\
             \end{array}
           \right),\qquad
\Y^\st(a)=\left(
             \begin{array}{rr}
               1 & 1 \\
               0 & 1 \\
             \end{array}
           \right)
$$
неэквивалентны над $F$ (независимо от выбора идеала $M$). Отметим тем не менее, что $\X^\st$ и
$\Y^\st$ имеют одинаковый набор неприводимых компонент. Как видно из следующего утверждения
это верно и в общем случае.

\begin{cor} \label{cor br nesb} Пусть  $\X$ и $\Y$ --- эквивалентные обыкновенные представления группы $G$,
записан\-ные над $\wt{Z}$. Тогда $F$-представления $\X^\st$ и $\Y^\st$ имеют одинаковый набор
неприводимых компонент (с учётом кратностей).
\end{cor}
\begin{proof} Пусть $\x$ --- общий характер представлений $\X$ и $\Y$. Из \ref{ord br}$(ii)$ вытекает, что
  $F$-представ\-ления 
$\X^\st$  и $\Y^\st$ имеют один и тот же брауэров характер $\hat\x$.
Требуемое следует из \ref{cor br pr}.
\end{proof}

Сформулируем ещё два важных следствия.

\begin{cor} \label{cor col br} Справедливы следующие утверждения.
\begin{list}{{\rm(}{\it\roman{enumi}\/}{\rm)}}
{\usecounter{enumi}\setlength{\parsep}{2pt}\setlength{\topsep}{5pt}\setlength{\labelwidth}{23pt}}
\item Множество $\iBr(G)$ является базисом пространства $\cf(G_{p'})$.
\item $|\iBr(G)|=|\K(G_{p'})|$.
\item Пусть $E$ --- поле характеристики $p$, являющееся полем разложения группы $G$. Тогда
число попарно неэквивалентных неприводимых $E$-представлений группы $G$ равно $|\K(G_{p'})|$.
\end{list}
\end{cor}
\begin{proof} $(i)$ Ввиду \ref{br lin ind}$(ii)$ множество $\iBr(G)$ линейно независимо. Поэтому достаточно
выразить любую $p$-регулярную классовую функцию в виде $\CC$-линейной комбинации неприводимых брауэровых
характеров. Пусть $\th\in \cf(G_{p'})$. В силу сюръективности отображения
$\,\hat{\mbox{}}:\cf(G)\to\cf(G_{p'})$ существует функция $\eta\in \cf(G)$ такая, что $\hat\eta=\th$.
Тогда $\eta=\sum_{\x\in \irr(G)}\a_\x \x$ для некоторых $\a_\x\in \CC$. Применив $\,\hat{\mbox{}}$,
получим $\th=\sum_{\x\in \irr(G)}\a_\x \hat\x$. Но каждый $\hat\x$ является брауэровым характером по
\ref{ord br}$(i)$ и, значит, представим в виде линейной комбинации элементов из $\iBr(G)$ по \ref{br sum}$(ii)$.

$(ii)$ Это следует из $(i)$, т.\,к. размерность пространства $\cf(G_{p'})$, как легко видеть,
совпадает с $|\K(G_{p'})|$.

$(iii)$ По предложению \ref{ext abs} число попарно неэквивалентных неприводимых
$E$-представлений группы $G$ совпадает с числом попарно неэквивалентных неприводимых
$\ov{E}$-представлений, где $\ov{E}$ --- алгебраическое замыкание поля $E$, а поскольку $F$ можно отождествить с подполем поля
$\ov{E}$, это же число равно по \ref{ext abs} также и числу попарно
неэквивалентных неприводимых $F$-представлений. Значит, требуемое следует из $(ii)$ ввиду
\ref{br dif}$(i)$.
\end{proof}

Напомним, что через \gls{opG} обозначается наибольшая нормальная $p$-подгруппа группы $G$.

\begin{cor} \label{cor op ker} Справедливы
следующие утверждения.
\begin{list}{{\rm(}{\it\roman{enumi}\/}{\rm)}}
{\usecounter{enumi}\setlength{\parsep}{2pt}\setlength{\topsep}{5pt}\setlength{\labelwidth}{23pt}}
\item Единственным неприводимым $F$-представлением $p$-группы является её
главное $F$-представление.
\item Если $\X$ --- неприводимое $F$-представление группы $G$, то $\OO_p(G)\le \ker \X$.
\item Пусть $\X_1,\ld,\X_s$ --- полный набор попарно неэквивалентных неприводимых $F$-представлений  группы
$G$. Тогда
$$
\bigcap_{i=1}^s\ker \X_i=\OO_p(G).
$$
\end{list}
\end{cor}
\begin{proof} $(i)$ Любая $p$-группа $P$ имеет единственный $p$-регулярный класс. Значит, в силу \ref{cor col br}$(ii)$
единственным неприводимым $F$-представлением группы $P$ является её главное представление.

$(ii)$ Обозначим $P=\OO_p(G)$. По теореме Клиффорда \ref{thm cliff}$(i)$ ограничение $\X_P$
является вполне приводимым $F$-представлением группы $P$. Из $(i)$ вытекает, что $\X_P$ ---
прямая сумма главных $F$-представлений, т.\,е. $P\le \ker \X$.

$(iii)$ Обозначим $N=\bigcap_i\ker \X_i$. Включение $\OO_p(G)\le N$ следует из $(ii)$. Поскольку $N\nor G$, для
доказательства обратного включения достаточно установить, что $N$ состоит из $p$-элементов. Пусть $g\in N$.
Тогда элемент $1-g\in FG$ лежит в $\Ker \X_i$ для любого $i=1,\ld,s$, где $\X_i$ рассматривается как
представление алгебры $FG$. Поскольку $\X_1,\ld,\X_s$ образуют полный набор её неприводимых представлений, из
\ref{int krs}$(i)$ следует, что $1-g\in \J(FG)$. Поэтому $(1-g)^n=0$ для некоторого натурального $n$ в силу
\ref{rad nlp}. Пусть $k$ --- произвольное натуральное число такое, что $n\le p^k$.
Тогда $0=(1-g)^{p^k}=1-g^{p^k}$, т.\,е. $g^{p^k}=1$. Отсюда следует требуемое. \end{proof}

Следствие \ref{cor op ker}$(i)$ показывает, что, в отличие от обыкновенных представлений (см. \ref{kpr}$(ii)$),
пересечение ядер неприводимых компонент $F$-представления $\X$ группы
не всегда совпадает с ядром самого представления $\X$. Например, в качестве $\X$ можно взять любое точное (например, регулярное)
$F$-представление нетривиальной $p$-группы.

\section{Таблица брауэровых характеров. Числа разложения. Инварианты Картана}

\glsadd{iTblBrChr}\mem{Таблицей брауэровых характеров} \gls{PhiMG} группы $G$ (относительно максимального идеала $M$)
называется квадратная матрица, строки которой
индексированы неприводимыми брауэровыми характерами $\vf\in \iBr(G)$, а столбцы ---
$p$-регулярными классами сопряжённости $K\in \K(G_{p'})$.
Элемент $\vf$-й строки и $K$-го столбца равен $\vf(x_{\mbox{}_K})$. Если идеал $M$ фиксирован,
то вместо $\Phi_M(G)$ будем просто писать \gls{PhiG} или же \gls{PhipG}, если нужно подчеркнуть, что
характеристика основного поля равна $p$. Например, ниже приведена таблица брауэровых
характеров группы $S_3$ в характеристике $2$.
$$
  \begin{array}{c|cc}
    \Phi_2(S_3)                      & 1       &  (123)       \\
     \hline
    \vf_1^{\vphantom{A^A}} & 1      & \phantom{-}1 \\
    \vf_2                  & 2      & -1
  \end{array}
$$

Пусть $\Phi=\Phi(G)$. Из \ref{br lin ind}$(ii)$ следует, что $\Phi$ является
невырожденной матрицей. Так как элементы матрицы $\Phi$ лежат в $\ov{\ZZ}$, мы можем
рассмотреть её образ относительно отображения $^\st$.

\begin{pre} \label{phi st} Матрица $\Phi^\st$ является невырожденной.
\end{pre}
\begin{proof}
В силу  \ref{br preg}$(iii)$ множество $\{\vf^\st\,|\,\vf\in
\iBr(G)\}$, где $\vf^\st(g)=\vf(g_{p'})^\st$ для всех $g\in G$, совпадает с множеством всех
неприводимых $F$-характеров, которое линейно независимо над $F$ ввиду \ref{rep trace}$(iii)$.
Но любой $F$-характер однозначно определяются своими значениями на множестве
представителей $p$-регулярных классов, см. \ref{br preg}$(i)$.
Из этих замечаний следует, что матрица $\Phi^\st$, которая и состоит из таких значений, является невырожденной.
\end{proof}

\begin{opr}\label{opr dec} Пусть
$\x\in \irr(G)$. Ввиду \ref{ord br}$(i)$ мы можем записать

$$
\hat\x=\sum_{\vf\in \iBr(G)}d_{\x\vf}\vf. \myeqno\label{xhat}
$$
Однозначно определённые неотрицательные целые числа \gls{dxf} называются \glsadd{iNumsDec}\mem{числами разложения} группы
$G$ в характеристике $p$.
Матрица
$$\gls{D}=(d_{\x\vf})_{\x\in \irr(G),\vf\in \iBr(G)}$$
называется \glsadd{iMatDec}\mem{матрицей разложения} и, очевидно,
является матрицей отображения $\,\hat{\mbox{}}:\cf(G)\to\cf(G_{p'})$  в базисах $\irr(G)$ и
$\iBr(G)$.
\end{opr}

Если положить $\XX=\XX(G)$ и $\Phi=\Phi(G)$, а также обозначить через $\wh\XX$ матрицу,
состоящую из столбцов матрицы $\XX$, соответствующих $p$-регулярным классам, то равенство
\ref{xhat} можно переписать в матричном виде
$$
\wh\XX=\DD\Phi.
$$

\textzam{Несмотря на то, что формально определение матрицы разложения также зависит от выбора максимального
идеала $M$, мы покажем в \ref{mm}$(i.2)$, что при другом выборе $M$, содержащего $p\ov{\ZZ}$,
полученная матрица разложения будет отличаться от данной лишь перестановкой строк и столбцов.}

\begin{pre} \label{dek lin ind} Имеем
\begin{list}{{\rm(}{\it\roman{enumi}\/}{\rm)}}
{\usecounter{enumi}\setlength{\parsep}{2pt}\setlength{\topsep}{5pt}\setlength{\labelwidth}{23pt}}
\item ранг матрицы разложения равен $|\iBr(G)|$;
\item столбцы матрицы разложения линейно независимы.
\end{list}
\end{pre}
\begin{proof} Поскольку отображение $\,\hat{\mbox{}}:\cf(G)\to\cf(G_{p'})$ сюръективно,
$\dim \cf(G_{p'})=|\iBr(G)|$ и  $\irr(G)$ --- базис пространства $\cf(G)$, отсюда следует, что среди
$\{\hat\x\,|\,\x\in \irr(G)\}$ найдётся $|\iBr(G)|$ линейно независимых характеров. Поэтому у матрицы
разложения имеется $|\iBr(G)|$ линейно независимых строк. Так как число её столбцов также равно
$|\iBr(G)|$, то отсюда следуют оба утверждения.
\end{proof}

Следующее утверждение показывает, что все строки и столбцы матрицы разложения ненулевые.

\begin{pre} \label{br e} Имеем
\begin{list}{{\rm(}{\it\roman{enumi}\/}{\rm)}}
{\usecounter{enumi}\setlength{\parsep}{2pt}\setlength{\topsep}{5pt}\setlength{\labelwidth}{23pt}}
\item если $\vf\in \iBr(G)$,
то существует $\x\in \irr(G)$ такой, что $d_{\x\vf}\ne 0$;
\item если $\x\in \irr(G)$,
то существует $\vf\in \iBr(G)$ такой, что $d_{\x\vf}\ne 0$.
\end{list}
\end{pre}
\begin{proof} Пункт $(i)$ вытекает из утверждения \ref{dek lin ind}$(ii)$.  Пункт $(ii)$ очевиден,
поскольку $0\ne \x(1)=\sum_{\vf\in \iBr(G)} d_{\x\vf}\vf(1)$.
\end{proof}

Следующее утверждение показывает, что в теории брауэровых характеров интересен только случай,
когда простое число $p$ является делителем порядка $|G|$.

\begin{pre} \label{br nmod} Если $p\nmid |G|$, то $\iBr(G)=\irr(G)$.
\end{pre}
\begin{proof} Предположим, что  $p\nmid |G|$.  Из \ref{cor z fg}$(ii)$ следует, что в этом
случае $$|G|=\sum_{\vf\in \iBr(G)}\vf(1)^2.\myeqno\label{ogss}$$ С другой стороны
$$
|G|=\sum_{\x\in
\irr(G)}\x(1)^2=\sum_{\x\in\irr(G)}\left(\sum_{\vf\in\iBr(G)}d_{\x\vf}\vf(1)\right)^2=
\sum_{\x\in\irr(G)}\sum_{\vf,\psi\in\iBr(G)}d_{\x\vf}d_{\x\psi}\vf(1)\psi(1)
$$
\vspace{-10pt}
$$
\hspace{100pt}=\sum_{\vf,\psi\in\iBr(G)}a_{\vf,\psi}\,\vf(1)\psi(1) \ge  \sum_{\vf\in \iBr(G)}\vf(1)^2, \myeqno\label{ogg}
$$
где $a_{\vf,\psi}=\sum_{\x\in\irr(G)}d_{\x\vf}d_{\x\psi}$ совпадает с произведением $\vf$-го и
$\psi$-го столбцов матрицы разложения (в смысле стандартного скалярного произведения) и поэтому последнее неравенство вытекает из того, что
$a_{\vf,\psi}\ge 0$ при $\vf\ne \psi$ и $a_{\vf,\vf}\ge 1$ ввиду \ref{br e}$(i)$. Сравнивая
\ref{ogss} и \ref{ogg}, получаем, что в действительности $a_{\vf,\psi}=\d_{\vf,\psi}$, т.\,е.
каждому $\vf$ соответствует ровно один $\x$ такой что, $d_{\x\vf}\ne 0$ (и в этом случае
$d_{\x\vf}= 1$), причём разным $\vf$ соответствуют разные $\x$. Другими словами $\DD$ --- это
матрица подстановки, и $\x=\hat\x=\vf$ для любой пары соответствующих характеров $\vf$ и $\x$.
Отсюда следует требуемое.
\end{proof}

\begin{cor} \label{cor ss restr} Пусть $\vf$ --- брауэров характер группы $G$, $H\le G$ и $p\nmid |H|$.
Тогда $\vf_{H^{\vphantom{A^a}}}\in\ch(H)$.
\end{cor}
\begin{proof} Из \ref{br prop}$(v)$ следует, что $\vf_{H^{\vphantom{A^a}}}$ --- брауэров характер группы $H$. Поэтому он
является суммой элементов из $\iBr(H)$, и требуемое вытекает из \ref{br nmod} и условия
$p\nmid |H|$.
\end{proof}

\textuprn{Считая известной таблицы обыкновенных характеров знакопеременной группы $A_5$  и её подгрупп (см. приложение
\ref{pril tch}), найти её таблицы брауэровых характеров в характеристиках $2$, $3$, и $5$.  \uk{Использовать изоморфизм $A_5\cong \SL_2(4)$ в случае характеристики $2$.}}

 Для любой функции $\th\in \cf(G_{p'})$ определим отображение $\gls{chkth}:G\to \CC$ по правилу
$$
\chk\th(g)=\th(g_{p'})
$$
для всех $g\in G$. Легко видеть, что $\chk\th\in \cf(G)$. Ясно также, что для любой функции
$\th\in\cf(G_{p'})$ если $\eta=\chk\th$, то $\hat\eta=\th$.

Любой брауэров характер является $\CC$-линейной комбинацией брауэровых характеров вида
$\hat\x$, где $\x\in \irr(G)$, что следует из \ref{dek lin ind}$(i)$. В действительности имеет
место более сильный факт.

\begin{pre} \label{z lin} Пусть $\vf\in\iBr(G)$. Тогда
\begin{list}{{\rm(}{\it\roman{enumi}\/}{\rm)}}
{\usecounter{enumi}\setlength{\parsep}{2pt}\setlength{\topsep}{5pt}\setlength{\labelwidth}{23pt}}
\item $\chk\vf\in \gch(G)$;
\item $\vf$ является $\ZZ$-линейной комбинацией брауэровых характеров $\hat\x$, где $\x\in \irr(G)$.
\end{list}
\end{pre}
\begin{proof}  $(i)$ Воспользуемся характеризацией Брауэра обобщённых характеров.
Пусть $E\le G$ --- элементарная брауэрова подгруппа. Запишем $E=P\times Q$, где $P$ ---
$p$-группа и $Q$ --- $p'$-группа. Если $g\in E$, то существуют элементы $x\in P$ и $y\in Q$
такие, что $g=xy=yx$. По лемме \ref{lem ppt} имеем $x=g_p$ и $y=g_{p'}$. Значит,
$\chk\vf(g)=\vf(g_p')=\vf(y)=\vf_Q(y)$. Но $\vf_Q\in\ch(Q)$ в силу
\ref{cor ss restr}, и он продолжается до обыкновенного характера $1_P\times\vf_Q$ группы $E$ по
\ref{dir ch}$(ii)$. Другими словами, ограничение $\chk\vf$ на $E$ совпадает с $1_P\times\vf_Q$ и
требуемое вытекает из теоремы \ref{thm char br}.

$(ii)$ Это следует из $(i)$, поскольку на множестве $G_{p'}$ функции $\vf$ и $\chk\vf$
тождественно равны.
\end{proof}

Напомним, что для натурального числа $n$ и простого $p$ выражения $n_p$ и $n_{p'}$ обозначают, соответственно,
$p$-часть и $p'$-часть числа $n$.

Пусть $m=\exp G$. Поскольку значение произвольного брауэрова характера $\vf$ группы $G$
являются суммой комплексных корней степени $m_{p'}$ из $1$, то для любого $\s\in
\Gal(\QQ_{m_{p'}},\QQ)$ можно определить $p$-регулярную классовую функцию $\vf^\s$ по правилу
$$
\vf^\s(g)=\vf(g)^\s
$$
для всех $g\in G_{p'}$. В отличие от обыкновенных характеров, функция $\vf^\s$, вообще говоря,
может не быть брауэровым характером относительно того же выбора максимального идеала $M$, а
если $\vf\in \iBr(G)$ и $\vf^\s$ --- брауэров характер, то не всегда $\vf^\s\in\iBr(G)$.

Например, пусть  $G=J_1$ ---  простая спорадическая группа Янко и $p=11$. Рассмотрим элемент
$\s\in \Gal(\QQ_{m_{p'}},\QQ)$, который на примитивном корне $\xi=e^{2\pi i/5}$ действует по правилу
$\s:\xi\mapsto\xi^2$ и оставляет неподвижными все примитивные корни из $1$ степени, взаимно
простой с $5$. Существуют неприводимые брауэровы характеры $\vf_1,\vf_2,\vf_3\in
\iBr(G)$ степеней $7$, $49$ и $56$, соответственно, такие, что $\vf_1^\s$ и $\vf_2^\s$ не
являются брауэровыми характерами, в то время как $\vf_3^\s=\vf_1+\vf_2$.

В связи с этим представляет интерес следующее утверждение о поведении брауэровых характеров под
действием автоморфизмов поля $\ov{\QQ}$ или при смене максимального идеала $M$.

\begin{pre} \label{mm} Пусть
$m=\exp G$.
\begin{list}{{\rm(}{\it\roman{enumi}\/}{\rm)}}{\usecounter{enumi}\setlength{\parsep}{2pt}\setlength{\topsep}{5pt}\setlength{\labelwidth}{23pt}}
\item Предположим, что $M_1$ и $M_2$ --- максимальные идеалы кольца $\ov{\ZZ}$, содержащие
$p\ov{\ZZ}$.  Тогда
    \begin{list}{{\rm(}{\it\roman{enumi}\/}.{\rm\arabic{enumii})}}
    {\usecounter{enumii}\setlength{\parsep}{2pt}\setlength{\topsep}{2pt}\setlength{\labelwidth}{5pt}\setlength{\labelwidth}{50pt}}
    \item существует автоморфизм $\s\in
\Gal(\QQ_{m_{p'}},\QQ)$ такой, что отображение $\vf\mapsto\vf^\s$ является биекцией между
множествами $\iBr_{M_1}(G)$ и $\iBr_{M_2}(G)$;
    \item этот автоморфизм $\s$ однозначно поднимается до элемента из
$\Gal(\QQ_m,\QQ_{m_p})$ (который также обозначим через $\s$) так, что
$d_{\x^\s\vf^\s}=d_{\x\vf}$ для любых $\x\in\irr(G)$ и $\vf\in \iBr_{M_1}(G)$, и, в частности,
матрица разложения группы $G$ в характеристике $p$ определяется однозначно с точностью до
перестановки строк и столбцов.
    \end{list}
\item Предположим, что $\vf$ --- брауэров характер группы $G$ относительно некоторого выбора
максимального идеала $M_1$ кольца $\ov{\ZZ}$, содержащего $p\ov{\ZZ}$, и $\s\in
\Gal(\QQ_{m_{p'}},\QQ)$. Тогда существует максимальный идеал $M_2$ кольца $\ov{\ZZ}$, содержащий $p\ov{\ZZ}$,
относительно которого $\vf^\s$ является брауэровым характером. При этом $\vf$ неприводим
относительно $M_1$ тогда и только тогда, когда $\vf^\s$ неприводим относительно $M_2$.
\end{list}
\end{pre}
\begin{proof} $(i)$  Положим $F_1=\ov{\ZZ}/M_1$ и $F_2=\ov{\ZZ}/M_2$. Существует изоморфизм полей $\r:F_1\to F_2$
в силу \ref{mul gr azp}$(ii)$ и \ref{alg pro}$(i)$.

Пусть $^{\st_1}:\ov{\ZZ}\to F_1$ и $^{\st_2}:\ov{\ZZ}\to F_2$ --- естественные эпиморфизмы. Обозначим
через $\s$ автоморфизм группы $U$, который действует так, что
$(\xi^{\st_1})^\r=(\xi^\s)^{\st_2}$ для любого $\xi\in U$, делая коммутативной следующую
диаграмму.

$$
\xymatrix{
U\ar[r]^{ \textstyle{ ^{\st_1} }}\ar@{.>}[d]^{ \textstyle{\s} } & F_1^\times\ar[d]^{\textstyle{\r}}        \\
U\ar[r]^{\textstyle{^{\st_2}}}         & F_2^\times
}
$$
Такой автоморфизм существует ввиду \ref{mul gr azp}$(i)$. Отображение $\s$ можно поднять
до однозначно определённого элемента из группы $\Gal(\QQ_{m_{p'}},\QQ)$. В самом деле, так как
$\s$ --- автоморфизм группы $U$, то существует целое число $k$, взаимно простое с $m_{p'}$
такое, что $\z^\s=\z^k$, где $\z$ --- примитивный корень степени $m_{p'}$ из $1$. Поэтому
отображение $\s$ определяет автоморфизм из $\Gal(\QQ_{m_{p'}},\QQ)$, который мы будем также
обозначать через~$\s$.

Пусть $\X$ ---  $F_1$-представление группы $G$ с брауэровым характером $\vf$ относительно
$M_1$. Поскольку значения $\vf$ являются суммами корней из $1$ степени $m_{p'}$, они лежат в
$\QQ_{m_{p'}}$ и, значит, определена $p$-регулярная классовая функция $\vf^\s$. Покажем, что
брауэров характер $F_2$-представления $\X^\r$ группы $G$ относительно $M_2$ равен $\vf^\s$,
причём $\vf$ неприводим относительно $M_1$ тогда и только тогда, когда $\vf^\s$ неприводим
относительно $M_2$. Отсюда будет следовать утверждение $(i.1)$.

Пусть $g\in G_{p'}$ и $\z_1^{\st_1},\ld,\z_s^{\st_1}$ --- корни характеристического многочлена
$f(x)=\det(\X(g)-x\I_s)\in F_1[x]$ для подходящих $\z_1,\ld,\z_s\in U$. Тогда $\vf(g)=\z_1+\ld+\z_s$. Поскольку
$\r$ индуцирует изоморфизм колец многочленов $F_1[x]\to F_2[x]$, отсюда следует, что характеристическими
значениями матрицы $\X^\r(g)$ будут корни $(\z_1^{\st_1})^\r,\ld,(\z_s^{\st_1})^\r$ многочлена $f^\r(x)$. По
выбору отображения $\s$ имеем $(\z_i^{\st_1})^\r=(\z_i^\s)^{\st_2}$, и поэтому значение брауэрова характера
$F_2$-представления $X^\r$ относительно $M_2$ на элементе $g$ равно $\z_1^\s+\ld+\z_s^\s=\vf^\s(g)$, что и
требовалось показать.

Равносильность неприводимости брауэровых характеров $\vf$  и  $\vf^\s$ относительно идеалов
$M_1$ и $M_2$, соответственно, вытекает из того, что $F_1$ --- представление $\X$ неприводимо
тогда и только тогда, когда неприводимо $F_2$ --- представление $\X^\r$. В самом деле, если
$\X$ эквивалентно блочно-верхнетреугольному $F_1$-представлению $\Y$, то, применив изоморфизм
$\r$, получим, что $\X^\r$ эквивалентно блочно-верхнетреугольному 
$F_2$-представле\-нию 
$\Y^\r$, и наоборот.

Теперь докажем утверждение $(i.2)$. Заметим, что $\s$ однозначно продолжается\footnote{В самом деле, пусть $\xi$ -- примитивный корень
степени $m$ из $1$. Достаточно определить $\xi^\s$. Поскольку $\xi$ однозначно представ\'{и}м в виде произведения $\xi_{p'}\xi_p$ своих $p'$- и $p$-частей,
причём $\xi_{p'}\in \QQ_{m_{p'}}$, можно положить $\xi^\s=\xi_{p'}^\s\xi_p$. Тогда  $\s\in \Gal(\QQ_m,\QQ_{m_p})$ и это определение не зависит от выбора $\xi$,
т.\,к. любые два примитивных корня степени $m$ из $1$ отличаются возведением в степень, взаимно простую с $m$.}
до автоморфизма из группы $\Gal(\QQ_m,\QQ_{m_p})$.
Пусть $\x\in \irr(G)$. Тогда $(\hat\x)^\s=\wh{\x^\s}$. Однако
$$(\hat\x)^\s=\left(\ \sum_{\vf\in
\iBr_{M_1}(G)}d_{\x\vf}\vf\right)^\s
=\sum_{\vf\in \iBr_{M_1}(G)}d_{\x\vf}\vf^\s, $$
а также

$$\wh{\x^\s}=\sum_{\vf\in \iBr_{M_1}(G)}d_{\x^\s\vf^\s}\vf^\s.$$
Сравнив коэффициенты в правых частях, получим $d_{\x^\s\vf^\s}=d_{\x\vf}$, как и утверждалось
в $(i.2)$.

$(ii)$ Пусть $F_1=\ov{\ZZ}/M_1$, $^{\st_1}:\ov{\ZZ}\to F_1$ --- естественный эпиморфизм и $\X_1$ ---
$F_1$-представление группы $G$, брауэровым характером которого (относительно идеала $M_1$)
является $\vf$.

Из \ref{alg pro}$(ii)$ следует, что $\s$ допускает продолжение до элемента из
$\Gal(\ov\QQ,\QQ)$, который также обозначим через $\s$. Положим $M_2=M_1^\s$. Ясно, что
$\ov{\ZZ}^\s=\ov{\ZZ}$ и $M_2$ --- максимальный идеал кольца $\ov{\ZZ}$, содержащий $p\ov{\ZZ}$. Обозначим
$F_2=\ov{\ZZ}/M_2$ и пусть $^{\st_2}:\ov{\ZZ}\to F_2$ --- естественный гомоморфизм. Тогда $\s$
определяет отображение $\r: F_1\to F_2$, заданное правилом $(\a^{\st_1})^\r=(\a^\s)^{\st_2}$ для
всех $\a\in \ov{\ZZ}$. Рассуждая, как в начале доказательства пункта $(i)$, получим
совпадение отображения $\vf^\s$ с брауэровым характером $F_2$-представления $\X_2=\X_1^\r$
относительно $M_2$, а также равносильность неприводимости $\vf$ относительно $M_1$ и
неприводимости $\vf^\s$ относительно $M_2$.
\end{proof}

Тем не менее в группе $\Gal(\QQ_{m_{p'}},\QQ)$ существует циклическая подгруппа, элементы
которой действуют на множестве $\iBr(G)$.

\begin{pre} \label{aut br} Пусть $m=\exp G$ и автоморфизм $\s\in \Gal(\QQ_{m_{p'}},\QQ)$ переводит
примитивный корень  $\z$  степени $m_{p'}$ из $1$ в $\z^p$. Тогда $\vf\in \iBr(G)$ если и
только если $\vf^\s\in \iBr(G)$. Порядок подгруппы $\la\s\ra\le \Gal(\QQ_{m_{p'}},\QQ)$
совпадает с порядком $p$ по модулю $m_{p'}$.
\end{pre}
\textupl{aut br prf}{Доказать предложение \ref{aut br}. \uk{Использовать автоморфизм Фробениуса поля $F$.}}

Пусть $\DD$ --- матрица разложения группы $G$. Произведение
$$
\gls{C}=\DD^\top \DD
$$
называется \mem{матрицей Картана}\glsadd{iMatCart}  группы $G$ (в характеристике $p$).
Очевидно, что $\C$ является симметрической квадратной $|\iBr(G)|\times|\iBr(G)|$-матрицей с неотрицательными
целыми элементами. Поскольку $\DD$ имеет ранг $|\iBr(G)|$, из элементарной линейной алгебры вытекает,
что $\C$ положительно определена. Положим $\C=(c_{\vf\psi})_{\vf,\psi\in \iBr(G)}$. Тогда
$$
c_{\vf\psi}=\sum_{\x\in \irr(G)}d_{\x\vf}d_{\x\psi}
$$
--- скалярное произведение $\vf$-го и $\psi$-го столбцов матрицы разложения $\DD$.
Элементы $c_{\vf\psi}$ называют \mem{инвариантами Картана}\glsadd{iInvCart} группы $G$.

\begin{pre} \label{inv car} Матрица Картана группы $G$ определяется
однозначно с точностью до перестановки строк и столбцов. В частности, набор инвариантов Картана
группы $G$ определён однозначно.
\end{pre}

\textupr{Доказать предложение \ref{inv car}. \uk{Воспользоваться \ref{mm}$(i.2)$.}}

\begin{opr}\label{opr gl ner}
Для произвольного
$\vf\in \iBr(G)$ определим обыкновенный характер
$$
\gls{thsvf}=\sum_{\x\in \irr(G)}d_{\x\vf}\x,
$$
который называется \glsadd{iPrincIndecChr}\mem{главным неразложимым характером}, соответствующим брауэрову характеру $\vf$.
\end{opr}

Из определения
вытекает, что
$$
\wh{\th_\vf}=\sum_{\x\in \irr(G)}d_{\x\vf}\hat{\x}=
\sum_{\x\in \irr(G)}\sum_{\psi\in \iBr(G)} d_{\x\vf}d_{\x\psi}\psi=
\sum_{\psi\in \iBr(G)} c_{\vf\psi}\psi. \myeqno\label{thvfh}
$$
Пусть
\gls{Theta} --- \mem{таблица значений главных неразложимых характеров} на $G_{p'}$, т.\,е.
квадратная матрица, строки которой индексированы характерами
$\vf\in \iBr(G)$, а столбцы --- $p$-регулярными классами сопряжённости $K\in \K(G_{p'})$.
Элемент $\vf$-й строки и $K$-го столбца матрицы $\Theta$
равен $\th_\vf(x_{\mbox{}_K})$. Тогда из \ref{thvfh} следует, что
$$
\Theta=\C\Phi.
$$
Чуть ниже (предложение \ref{gl ner}$(i)$) мы покажем, что матрица $\Theta$ полностью определяет значения главных неразложимых характеров.

Пусть $\mu,\nu\in \cf(G)\cup\cf(G_{p'})$. Положим
$$
\gls{lvfpsrsGpp}=\frac{1}{|G|}\sum_{g\in G_{p'}}\mu(g)\ov{\nu(g)}.
$$
Обозначим через \gls{cfcirclGr} множество функций из $\cf(G)$,
тождественно равных нулю вне $G_{p'}$. Легко видеть, что $\cf^\circ(G)$ --- идеал алгебры $\cf(G)$. Заметим также, что
$\cf^\circ(G)$ является $\CC$-алгеброй,
и ограничение отображения $\,\hat{\mbox{}}{\mbox{}}:\cf(G)\to\cf(G_{p'})$ на $\cf^\circ(G)$ будет
изоморфизмом $\CC$-алгебр.

\textuprn{Проверить, что $\cf^\circ(G)$ и $\cf(G_{p'})$ являются унитарными пространствами
со скалярным произведением $(\cdot,\cdot)_{\mbox{}_{G_{p'}}}$.}

Основные свойства главных неразложимых характеров перечислены в следующем утверждении.

\begin{pre} \label{gl ner} Имеем
\begin{list}{{\rm(}{\it\roman{enumi}\/}{\rm)}}
{\usecounter{enumi}\setlength{\parsep}{2pt}\setlength{\topsep}{5pt}\setlength{\labelwidth}{23pt}}
\item $\th_\vf\in \cf^\circ(G)$ для любого $\vf\in \iBr(G)$.
В частности, значения главных неразложимых характеров определяются матрицей $\Theta$;
\item $(\vf,\th_\psi)_{\mbox{}_{G_{p'}}}=\d_{\vf,\psi}$ для любых $\vf,\psi\in \iBr(G)$;
\item множество $\{\th_\vf\,|\,\vf\in \iBr(G)\}$ является базисом пространства $\cf^\circ(G)$;
\item если $\x\in \cf^\circ(G)$ --- обыкновенный характер, то $\x$ является целочисленной
линейной комбинацией главных неразложимых характеров.
\item $\C^{-1}=\big(\,(\vf,\psi)_{\mbox{}_{G_{p'}}}\big)_{\vf,\psi\in \iBr(G)}$.
\item $(\th_\vf,\th_\psi)_{\mbox{}_G}=(\th_\vf,\th_\psi)_{\mbox{}_{G_{p'}}}=c_{\vf\psi}$ для
любых $\vf,\psi\in \iBr(G)$.
\item $|G|_p$ делит $\th_\vf(1)$ для любого $\vf\in \iBr(G)$.
\end{list}
\end{pre}
\begin{proof} $(i)$ Пусть $x\in G$ и $y\in G_{p'}$. Обозначим $K=x^G$, $L=y^G$.
Из второго соотношения ортогональности \ref{vtor ort}
вытекает
$$
\sum_{\x\in \irr(G)}\ov{\x(x)}\x(y)=\d_{\mbox{}_{K,L}}|\C_G(x)|.
$$
Поскольку $y$ --- $p$-регулярный элемент, имеем $\x(y)=\sum_{\vf\in\iBr(G)}d_{\x\vf}\vf(y)$. Значит,
$$
\sum_{\x\in \irr(G)}\ov{\x(x)}\x(y)=
\sum_{\x\in \irr(G)}\sum_{\vf\in \iBr(G)}d_{\x\vf}\ov{\x(x)}\vf(y)=
\sum_{\vf\in \iBr(G)}\ov{\th_\vf(x)}\vf(y),
$$
и, следовательно,
$$
\sum_{\vf\in \iBr(G)}\ov{\th_\vf(x)}\vf(y)=\d_{\mbox{}_{K,L}}|\C_G(x)|. \myeqno\label{dkl}
$$
Предположим теперь, что $x\not\in G_{p'}$. Тогда $K\ne L$ и в силу произвольности $y$ получаем
$$
\sum_{\vf\in \iBr(G)}\ov{\th_\vf(x)}\vf=0.
$$
По \ref{cor col br}$(i)$ неприводимые брауэровы характеры линейно независимы. Значит,
$\th_\vf(x)=0$ и $\th_\vf\in \cf^\circ(G)$.

$(ii)$ Обозначим через $Q$ целочисленную диагональную матрицу,
индексированную $p$-регулярными классами и состоящую
из элементов $\d_{\mbox{}_{K,L}}|\C_G(x_{\mbox{}_K})|$. В этих обозначениях равенство \ref{dkl}
для $x\in G_{p'}$ можно переписать в матричном виде
$$
\ov{\Theta}^\top \Phi=Q.
$$
Отсюда следует, что $$\I=\ov{\Theta}^\top (\Phi Q^{-1})=(\Phi Q^{-1})\ov{\Theta}^\top.$$
Выбрав произвольные $\vf,\psi\in \iBr(G)$, получаем
$$
\d_{\vf,\psi}=\sum_{K\in \K(G_{p'})}\vf(x_{\mbox{}_K})\,\frac{1}{|\C_G(x_{\mbox{}_K})|}\,\ov{\th_\psi(x_{\mbox{}_K})}=
\sum_{g\in G_{p'}}\frac{\vf(g)\ov{\th_\psi(g)}}{|G|}=(\vf,\th_\psi)_{\mbox{}_{G_{p'}}}.
$$

$(iii)$ Размерность $\dim \cf^\circ(G)$ очевидно равна $|\K(G_{p'})|$,
что совпадает c числом главных неразложимых характеров по \ref{cor col br}$(ii)$.
Линейная независимость этих характеров следует из $(ii)$ и, значит, они образуют базис пространства
$\cf^\circ(G)$.

$(iv)$ Из $(iii)$ следует, что $$\x=\sum_{\vf\in\iBr(G)} a_\vf\th_\vf$$
для некоторых $a_\vf\in \CC$. Тогда
$$
a_\vf=(\vf,\x)_{\mbox{}_{G_{p'}}}=(\chk\vf,\x)_{\mbox{}_{G}}\in \ZZ,
$$
где первое равенство следует из $(ii)$, второе --- из того, что $\x\in \cf^\circ(G)$,
а включение --- из \ref{z lin}$(i)$.

$(v)$ Пусть $\vf,\psi\in \iBr(G)$. Мы уже отмечали, что

$$
\wh{\th_\psi}=\sum_{\eta\in \iBr(G)} c_{\psi\eta}\eta.  \myeqno\label{thpsh}
$$
Поскольку $c_{\psi\eta}=c_{\eta\psi}$, из $(ii)$ следует, что
$$
\d_{\vf,\psi}=(\vf,\th_\psi)_{\mbox{}_{G_{p'}}}=(\vf,\wh{\th_\psi})_{\mbox{}_{G_{p'}}}
=\sum_{\eta\in \iBr(G)} c_{\eta\psi}(\vf,\eta)_{\mbox{}_{G_{p'}}},
$$
т.\,е. матрица $\big(\,(\vf,\psi)_{\mbox{}_{G_{p'}}}\big)_{\vf,\psi\in \iBr(G)}$
обратна к матрице Картана $\C$.

$(vi)$ В силу $(i)$, \ref{thpsh} и $(ii)$ получаем
$$
(\th_\vf,\th_\psi)_{\mbox{}_G}=(\th_\vf,\th_\psi)_{\mbox{}_{G_{p'}}}=
(\wh{\th_\vf},\th_\psi)_{\mbox{}_{G_{p'}}}=
\sum _{\eta\in \iBr(G)} c_{\eta\vf}(\eta,\th_\psi)_{\mbox{}_{G_{p'}}}=c_{\psi\vf}=c_{\vf\psi}
$$

$(vii)$ Пусть $\vf\in \iBr(G)$ и $P\in \Syl_p(G)$.\footnote{Напомним, что через
\gls{SylsplGr} обозначается множество $p$-силовских подгрупп группы $G$.}

Поскольку $(\th_\vf)_P$ и $1_P$ --- обыкновенные
характеры группы $P$, то число
$$
\big((\th_\vf)_P,1_P\big)_{\!\mbox{}_P}=\frac{1}{|P|}\,\th_\vf(1)
$$
является целым, а равенство следует из того, что  $\th_\vf(g)=0$ для любого неединичного $g\in P$
в силу $(i)$.  Отсюда вытекает $(vii)$.
\end{proof}


Пусть $\eta \in \cf(G_{p'})$. Определим отображение $\wt\eta:G\to\CC$ по правилу
$$
\wt\eta (g)=\left\{\ba{ll}
|G|_p\,\eta(g),& \mbox{если}\ \ g\in G_{p'};\\
0,& \mbox{иначе}
\ea\right.
$$
для любого $g\in G$ (cp. определение \ref{opr dot x}). Легко видеть, что $\wt\eta\in \cf^\circ(G)$
и отображение $\wt{}:\cf(G_{p'})\to \cf^\circ(G)$ является $\CC$-линейным и сюръективным.

\begin{pre} \label{ti} Если $\vf\in\iBr(G)$, то $\wt\vf\in \gch(G)$.
\end{pre}
\begin{proof}
В силу теоремы Брауэра \ref{thm char br} достаточно проверить, что если $E$ --- элементарная
брауэрова подгруппа, то $\wt\vf_{\mbox{}_E}$ --- её обобщённый характер. Запишем $E=P\times Q$,  где $P$ --- $p$-группа
и $Q$ --- $p'$-группа. Тогда
$$
\wt\vf_{\mbox{}_E} (g)=\left\{\ba{ll}
|G|_p\,\vf_{\mbox{}_Q}(g),& \mbox{если}\ \  g\in Q;\\
0,& \mbox{иначе}
\ea\right.
$$
для всех $g\in E$. Заметим, что $\vf_{\mbox{}_Q}\in\ch(Q)$ ввиду \ref{cor ss restr}.
Пусть $\r_{\mbox{}_P}$ --- регулярный характер группы $P$
и положим $\tau=\r_{\mbox{}_P}\times \vf_{\mbox{}_Q}$.
Из \ref{dir ch}$(ii)$ следует, что
$\tau\in\ch(E)$, и ввиду \ref{reg har}$(ii)$
получаем, что
$$
\tau(g)=\left\{\ba{ll}
|P|\,\vf_{\mbox{}_Q}(g),& \mbox{если}\ \ g\in Q;\\
0,& \mbox{иначе}
\ea\right.
$$
для всех $g\in E$.
Поэтому $\wt\vf_{\mbox{}_E}=\displaystyle\frac{|G|_p}{|P|}\,\tau\in\ch(E)$.
\end{proof}

Используя классовые функции вида $\wt\vf$, можно получить дополнительную информацию о матрице Картана.

\begin{pre} \label{cart det} Пусть $\C$ --- матрица Картана группы $G$. Тогда
 $$\det(\C)=\prod_{K\in \K(G_{p'})}|\C_G(x_{\mbox{}_K})|_p. \myeqno\label{detc}$$
В частности, определитель матрицы Картана является степенью числа $p$.
\end{pre}
\begin{proof} Сначала докажем, что $\det(\C)$ --- степень числа $p$.

Из определения отображения $\wt{}:\cf(G_{p'})\to \cf^\circ(G)$ следует, что для любых $\vf,\psi\in \iBr(G)$
имеет место соотношение
$$
(\wt\vf,\wt\psi)_{\mbox{}_{G}}=|G|_p^2(\vf,\psi)_{\mbox{}_{G_{p'}}}.
$$
Пусть $A=\big(\,(\wt\vf,\wt\psi)_{\mbox{}_{G}}\big)_{\vf,\psi\in \iBr(G)}$. Тогда, учитывая \ref{gl ner}$(v)$,
предыдущее соотношение можно переписать в матричном виде
$$
A=|G|_p^2\,\C^{-1}.
$$
Значит, $A\C=|G|_p^2\I$. Но в силу \ref{ti} матрица  $A$ целочисленная, и поэтому $\det(\C)$ --- степень числа $p$.

Теперь докажем равенство \ref{detc}. Как мы видели, имеет место матричное равенство
$$
\Theta=\C\Phi.
$$
Кроме того, в доказательстве \ref{gl ner}$(ii)$ было установлено, что
$$
\ov{\Theta}^\top \Phi=Q,
$$
где $Q$ --- диагональная матрица с элементами $\d_{\mbox{}_{K,L}}|\C_G(x_{\mbox{}_K})|$, индексированными
$p$-регулярными классами. Отсюда вытекает, что
$$
\ov{\Phi}^\top \C \Phi=Q,
$$
поскольку матрица $\C$ целочисленная и симметрическая. Заметим, что ввиду предложения \ref{spr br} матрица
$\ov\Phi$ отличается от $\Phi$ возможно лишь перестановкой строк, и поэтому
$\det(\ov{\Phi}^\top)=\det(\ov{\Phi})=\pm\det(\Phi)$. Значит,
$$
\det(\Phi)^2=\pm\frac{\prod_{K\in \K(G_{p'})}|\C_G(x_{\mbox{}_K})|}{\det(\C)}. \myeqno\label{dtfs}
$$
Поскольку значения брауэровых характеров лежат в $\ov{\ZZ}$, из предыдущего равенства следует, что
$$
\det(\Phi)^2\in \ov{\ZZ}\cap\QQ=\ZZ.
$$
Предположим, что $p$ делит  $\det(\Phi)^2$. Тогда
$$
0=(\det(\Phi)^2)^\st=(\det(\Phi)^\st)^2=(\det(\Phi^\st))^2
$$
и, значит, матрица $\Phi^\st$ вырожденная, вопреки \ref{phi st}.
Итак, $\det(\Phi)^2$ не делится на $p$, и поскольку $\det(\C)$ --- степень $p$, требуемое
следует из равенства \ref{dtfs}.
\end{proof}

\section{Блоки характеров}

В этом разделе мы определим понятие $p$-блока --- одного из ключевых объектов в теории модулярных
представлений, введённых Р.\,Брауэром.

По аналогии с центральными характерами $\om_\x$, определёнными ранее для любых $\x\in \irr(G)$,
введём также отображения $\gls{lsvf}: Z(FG)\to F$ для любого $\vf\in \iBr(G)$.
Пусть $\X$ --- $F$-представление с брауэровым характером $\vf$. В силу неприводимости $\X$ из \ref{cent irr}
следует, что для любого $z\in \Z(FG)$ матрица $\X(z)$ равна $a\I$ для некоторого $a \in F$. Мы положим
$\l_\vf(z)=a$.

В \ref{ch ac} мы показали, что для любых $\x\in \irr(G)$ и $K\in \K(G)$
значение центрального характера $\om_\x:\Z(\CC G)\to \CC$ на $\wh K$ является целым алгебраическим числом.
Как было отмечено в \ref{pg}, отображение $^\st:\ov{\ZZ}\to F$ можно
естественно поднять  до эпиморфизма колец $\Z(\ov{\ZZ} G)\to \Z(F G)$ и что
значения центрального характера $\om_\x$ на элементах из $\Z(\ov{\ZZ} G)$ лежат в $\ov{\ZZ}$.
Ясно также, что ограничение $\om_\x$ на $\Z(\ov{\ZZ} G)$ является гомоморфизмом $\ov{\ZZ}$-алгебр.
Из \ref{hom r}$(i)$ вытекает\footnote{Предложение \ref{hom r}$(i)$ применимо, поскольку
ядро гомоморфизма  $^\st$ на $\Z(\ov{\ZZ} G)$ содержится в $\Ker \om_\x^\st$.
В самом деле, ядро $^\st$ на $\Z(\ov{\ZZ} G)$ состоит из всевозможных линейных
комбинаций $\sum_{K\in \K(G)} \a_{\mbox{}_K} \wh K$,
где $\a_{\mbox{}_K}\in M$.  Однако образ такой линейной комбинации под действием $\om_\x$ равен
$\sum_{K\in \K(G)} \a_{\mbox{}_K} \om_\x(\wh K)\in M$ и, значит, обращается в $0$
под действием $^\st$.} существование и единственность гомоморфизма
колец $\gls{lsx}: \Z(FG)\to F$,  для которого коммутативна следующая диаграмма.

$$
\xymatrix{
\Z(\ov{\ZZ} G)\ar[dr]|-{ \textstyle{\om_\x^\st\,}}\ar[r]^-{ \textstyle{\om_\x}}\ar[d]_{ \textstyle{\st}}&
\ov{\ZZ} \ar[d]^{ \textstyle{\st}}       \\
\Z(FG)\ar@{.>}[r]_-{ \textstyle{\l_\x} }     &   F}
$$
Здесь $\om_\x^\st$ обозначает композицию отображений  $\om_\x:\Z(\ov{\ZZ} G)\to \ov{\ZZ}$ и $^\st: \ov{\ZZ}\to F$.

\begin{pre} \label{hom mod} Для любого $\eta\in \irr(G)\cup\iBr(G)$ отображение $\l_\eta:\Z(FG)\to F$ является
гомоморфизмом $F$-алгебр.
\end{pre}

\textupl{hom mod prf}{Доказать предложение \ref{hom mod}.}

Из предложения \ref{hom mod}, в частности, вытекает, что для любого $\x\in \irr(G)$ отображение
$\l_\x$ однозначно задаётся своими значениями
$$
\l_\x(\wh K)=\om_\x(\wh K)^\st \myeqno\label{dz}
$$
на базисных элементах $\wh K$ алгебры $\Z(FG)$.
Отметим, что $\wh K$ в левой части равенства \ref{dz}  является элементом из $FG$, а в правой --- из $\CC G$.
В дальнейшем мы будем преимущественно пользоваться равенством \ref{dz} для вычисления значений
центрального характера $\l_\x$.

\mysubsection\label{zwtr}
Напомним, что ранее мы продолжили отображение $^\st:\ov{\ZZ}\to F$ на локальное кольцо $\wt{Z}$.
Далее нам понадобится рассмотреть ограничение $\om_\x$ не только на $\Z(\ov{\ZZ} G)$,
но также и на $\Z(\wt{Z} G)$. Соображения, аналогичные приведённым выше,
показывают, что $\om_\x(\Z(\wt{Z} G))\se \wt{Z}$,
и позволяют поднять отображение $^\st$ до эпиморфизма колец $\Z(\wt{Z} G)\to \Z(F G)$.

\begin{pre} \label{com t} Во введённых обозначениях для любого $\x\in \irr(G)$ коммутативна следующая диаграмма
$$
\xymatrix{
\Z(\wt{Z} G)\ar[r]^-{ \textstyle{\om_\x}}\ar[d]^{ \textstyle{\st}}&
\wt{Z} \ar[d]^{ \textstyle{\st}}       \\
\Z(FG)\ar[r]^-{ \textstyle{\l_\x} }     &  F}
$$
\end{pre}
\textupl{com t prf}{Доказать предложение \ref{com t}.}

В отличие от центральных характеров $\om_\x$, введённые гомоморфизмы $\l_\x$ при $\x \in \irr(G)$,
равно как и $\l_\vf$ при $\vf \in \iBr(G)$, вообще говоря, не являются попарно различными.
Для характеров $\mu,\nu\in \irr(G)\cup\iBr(G)$ будем писать $\mu\sim\nu$ тогда и только тогда,
когда $\l_\mu=\l_\nu$. Очевидно, что $\sim$ является отношением эквивалентности.

\textuprn{Показать, что $1_G\sim 1_{G_{p'}}$.}

\begin{opr}
Классы эквивалентности множества $\irr(G)\cup\,\iBr(G)$
относительно отношения $\sim$ называются \glsadd{iPBlk}\mem{$\!p$-блоками} группы $G$.
Множество $p$-блоков обозначается через
\gls{BlplGr} или просто \gls{BllGr}.
\end{opr}

Пусть $B\in \bl(G)$. Из определения следует, что $p$-блок $B$ имеет <<обыкновенную часть>> и <<брауэрову часть>>.
Положим $\gls{IrrlBr}=B\cap\irr(G)$ --- \mem{$p$-блок обыкновенных характеров}\glsadd{iPBlkOrdChr}
и $\gls{IBrlBr}=B\cap\iBr(G)$ ---
\glsadd{iPBlkBrChr}\mem{$p$-блок брауэровых характеров}. Далее в \ref{bl ne} мы увидим, что для любого $p$-блока $B$
множества $\irr(B)$ и $\iBr(B)$ непусты.

Обозначим через \gls{lsB} однозначно определённый гомоморфизм $F$-алгебр $\Z(FG)\to F$, равный $\l_\eta$ для
произвольного $\eta\in B$ и будем называть его \glsadd{iChrCenCorBlk}\mem{центральным характером, соответствующим блоку} $B$.

\mem{Главным $p$-блоком}\glsadd{iPBlkPrnc} группы $G$ называется (единственный)
$p$-блок, содержащий характеры $1_G$ и $1_{G_{p'}}$.

Разбиение множества $\irr(G)$ на $p$-блоки обыкновенных характеров можно осуществить,
исходя из таблицы обыкновенных характеров группы $G$.
Из предложения \ref{hom mod} и формулы \ref{ch val} следует, что два характера $\x,\t\in \irr(G)$ лежат в одном $p$-блоке
тогда и только тогда, когда выполнены равенства
$$
\left(\frac{|K|\x(x_{\mbox{}_K})}{\x(1)}\right)^\st=\left(\frac{|K|\t(x_{\mbox{}_K})}{\t(1)}\right)^\st
$$
для всех $K\in \K(G)$. Строго говоря, эти равенства определяют эквивалентность $\x\sim\t$
при фиксированном выборе идеала $M$, содержащего число $p$. Следующее утверждение показывает, однако,
что эта эквивалентность на множестве $\irr(G)$ не зависит от $M$.

%

\begin{pre} \label{bl nez} Пусть $\x,\t\in \irr(G)$. Для того чтобы $\x$ и $\t$ лежали в одном $p$-блоке необходимо
и достаточно, чтобы для любого $K\in \K(G)$ величина $\om_\x(\wh K)-\om_\t(\wh K)$ лежала в каждом максимальном
идеале кольца $\ov{\ZZ}$, содержащем $p\ov{\ZZ}$.
В частности, разбиение множества $\irr(G)$ на $p$-блоки обыкновенных характеров не зависит от $M$.
\end{pre}
\begin{proof} Достаточность очевидна. Покажем необходимость. Допустим, что $\x\sim\t$ и $K\in \K(G)$.
Обозначим $\a=\om_\x(\wh K)-\om_\t(\wh K)$. Если мы покажем, что для некоторого натурального числа $n$ имеет место
включение $\a^n\in p\ov{\ZZ}$, то отсюда будет следовать требуемое. В самом деле, если $\a$ не содержится в некотором
максимальном идеале $I$ из $\ov{\ZZ}$, то и никакая степень элемента $\a$ не лежит в $I$ ввиду того, что поле $\ov{\ZZ}/I$
не имеет ненулевых нильпотентных элементов.

Пусть $m=\exp G$
и $\s\in \Gal(\QQ_m,\QQ)$.
В силу \ref{gal conj}$(iii)$ существует целое число $k$, взаимно простое с $m$, такое,
что $\eta(g)^\s=\eta(g^k)$ для любого обыкновенного характера $\eta$ и любого $g\in G$.
Заметим, что для такого $k$ имеет место равенство  $|g^G|=|(g^k)^G|$, где $g\in G$ --- произвольный элемент,
поскольку $g$ и $g^k$ порождают одну и ту же циклическую подгруппу и, значит, имеют одинаковые централизаторы в $G$.
Учитывая это, получаем
$$
\a^\s=\left(\frac{|K|\x(x_{\mbox{}_K})}{\x(1)}-\frac{|K|\t(x_{\mbox{}_K})}{\t(1)}\right)^\s=
\frac{|L|\x(x_{\mbox{}_K}^k)}{\x(1)}-\frac{|L|\t(x_{\mbox{}_K}^k)}{\t(1)}=\om_\x(\wh L)-\om_\t(\wh L)\in M,
$$
где $L$ --- класс сопряжённости, содержащий элемент $x_{\mbox{}_K}^k$,
а последнее включение вытекает из соотношения $\x\sim\t$.

Рассмотрим многочлен
$$
f(x)=\prod_{\s\in \Gal(\QQ_m,\QQ)} (x-\a^\s)\in \ov{\ZZ}[x].
$$
Его коэффициенты инвариантны относительно всех элементов группы $\Gal(\QQ_m,\QQ)$
и, следовательно, лежат в $\QQ$.
Поскольку $\ov{\ZZ}\cap \QQ=\ZZ$, имеем $f(x)\in \ZZ[x]$. С другой стороны мы показали, что
$\a^\s\in M$ для всех $\s\in \Gal(\QQ_m,\QQ)$. Поэтому все коэффициенты многочлена $f(x)$, кроме старшего,
лежат в $\ZZ\cap M=p\ZZ$. Подставив $\a$ в $f$, получим $0=f(\a)\equiv \a^n\mod {\,p\ov{\ZZ}}$,
где $n=|\Gal(\QQ_m,\QQ)|$, что и требовалось доказать.
\end{proof}

После того, как множество $\irr(G)$ разбито на $p$-блоки обыкновенных характеров, мы можем также
разбить множество $\iBr(G)$ на $p$-блоки брауэровых характеров, если известна матрица разложения.

\begin{pre} \label{mod bl} Имеем
\begin{list}{{\rm(}{\it\roman{enumi}\/}{\rm)}}
{\usecounter{enumi}\setlength{\parsep}{2pt}\setlength{\topsep}{5pt}\setlength{\labelwidth}{23pt}}
\item пусть $\x\in \irr(G)$ и $\vf\in \iBr(G)$. Если $d_{\x\vf}\ne 0$, то $\l_\x=\l_\vf$;
\item $\iBr(B)=\{\vf\in \iBr(G)\,|\,d_{\x\vf}\ne 0$ для некоторого $\x\in \irr(B)\}$.
\end{list}
\end{pre}
\begin{proof} Пусть $\x\in \irr(G)$ и $\X$ --- обыкновенное представление с характером $\x$. В силу
\ref{kom alg} можно считать, что $\X$ записано над $\wt{Z}$. Тогда $\hat\x$ является брауэровым
характером $F$-представления  $\X^\st$ ввиду \ref{ord br}$(ii)$.

Пусть $K\in \K(G)$.
По определению центрального характера $\om_\x$ имеем $\X(\wh K)=\om_\x(\wh K)\I$. Отсюда
$$\X^\st(\wh K)=\om_\x(\wh K)^\st\I=\l_\x(\wh K)\I.$$
Обозначим через $\Y$ блочно-верхнетреугольное $F$-представление с неприводимыми компонентами на главной
диагонали, эквивалентное представлению $\X^\st$.
Поскольку  матрица $\X^\st(\wh K)$ скалярная, получаем
$$\Y(\wh K)=\X^\st(\wh K)=\l_\x(\wh K)\I.$$
Так как любая неприводимая компонента $\ZZZ$ представления $\X^\st$ эквивалентна некоторой неприводимой
компоненте представления $\Y$, значение которой на $\wh K$ --- скалярная матрица, отсюда вытекает, что
$\ZZZ(\wh K)=\l_\x(\wh K)\I$. В частности, если $d_{\x\vf}\ne 0$ для некоторого $\vf\in \iBr(G)$,
то $\vf$ является брауэровым характером некоторой неприводимой компоненты представления $\X^\st$,
и поэтому $\l_\vf(\wh K)=\l_\x(\wh K)$. Отсюда, в силу произвольности $K$
и  того, что гомоморфизмы $\l_\vf$ и $\l_\x$ определяются своими значениями на классовых суммах,
получаем равенство $\l_\vf=\l_\x$,
а также включение
$$\{\,\vf\in \iBr(G)\,|\,d_{\x\vf}\ne 0\ \mbox{для некоторого}\ \x\in \irr(B)\,\}\se \iBr(B).$$
Обратно, если $\vf\in \iBr(B)$, то ввиду \ref{br e}$(i)$ существует
характер $\x\in \irr(G)$, для которого $d_{\x\vf}\ne 0$, и из только что доказанного вытекает, что
$p$-блок, которому принадлежит $\x$ должен совпадать с $B$.
\end{proof}

\begin{cor} \label{bl ne} Для любого $B\in \bl(G)$ имеем $\irr(B)\ne\varnothing$ и $\iBr(B)\ne\varnothing$.
\end{cor}
\begin{proof} Пусть $B\in \bl(G)$. Тогда либо $\x\in B$ для некоторого $\x \in \irr(G)$,
либо $\vf\in B$ для некоторого $\vf \in \iBr(G)$.  В силу \ref{br e} существуют, соответственно,
$\vf\in \iBr(G)$ и $\x\in \irr(G)$, такие, что $d_{\x\vf}\ne 0$, и тогда, соответственно,
 $\vf\in \iBr(B)$ и $\x\in \irr(B)$.
\end{proof}

\textuprn{Найти $2$- и  $3$-блоки группы $S_3$.}

В приложении \ref{pril tch} приведены таблицы обыкновенных характеров и матрицы разложения
для некоторых конечных групп по всем простым делителям $p$ их порядка. Исходя из этих данных,
можно определить $p$-блоки таких групп, используя предложения \ref{bl nez} и \ref{mod bl}.

\begin{opr} \mem{Графом Брауэра}\glsadd{iGrBr} группы $G$ называется граф с множеством вершин $\irr(G)$,
в котором два характера $\x,\t\in \irr(G)$ соединены ребром тогда и только тогда, когда существует характер
$\vf\in \iBr(G)$ такой,  что $d_{\x\vf}\ne 0$ и $d_{\t\vf}\ne 0$.
\end{opr}

Из предложения \ref{mod bl} вытекает, что если характеры $\x,\t\in \irr(G)$ соединены ребром в графе Брауэра,
то они лежат в одном $p$-блоке. Отсюда получаем

\begin{cor} \label{gr bl} Для любого $B\in \bl(G)$ множество
$\irr(B)$ является объединением компонент связности\footnote{Здесь и далее под компонентой связности графа Брауэра
мы подразумеваем множество вершин некоторой компоненты связности этого графа.} графа Брауэра.
\end{cor}

Существуют примеры, показывающие, что смежность вершин в графе Брауэра зависит от выбора максимального идеала $M$,
содержащего число $p$. Однако компоненты связности графа Брауэра, как мы увидим в \ref{bl ob}$(i)$,
в точности совпадают множествами $\irr(B)$ и поэтому не зависят от идеала $M$ ввиду \ref{bl nez}.
Чтобы показать это, сформулируем несколько вспомогательных утверждений.

\begin{pre} \label{os} Пусть $\A\se \irr(G)$ --- объединение компонент связности графа Брауэра. Положим
$$
\B=\{\vf\in \iBr(G)\,|\,d_{\x\vf}\ne 0\ \mbox{для некоторого}\ \x\in \A\}.
$$
Тогда для любых $x\in G_{p'}$ и $y\in G$ выполнено равенство
$$
\sum_{\x\in \A}\x(x)\ov{\x(y)}=\sum_{\vf\in \B}\vf(x)\ov{\th_\vf(y)}.
$$
\end{pre}
\begin{proof} Для любого $\x\in \A$ по условию имеем
$$
\x(x)=\sum_{\vf\in \iBr(G)}d_{\x\vf}\vf(x)=\sum_{\vf\in \B}d_{\x\vf}\vf(x).
$$
Также для любого $\vf\in \B$
$$
\th_\vf(y)=\sum_{\x\in\irr(G)}d_{\x\vf}\x(y)=\sum_{\x\in\A}d_{\x\vf}\x(y).
$$
Поэтому
$$
\sum_{\x\in \A}\x(x)\ov{\x(y)}=\sum_{\substack{\x\in\A,\\\vf\in \B}}d_{\x\vf}\vf(x)\ov{\x(y)}=
\sum_{\vf\in \B}\vf(x)\ov{\th_\vf(y)},
$$
что и требовалось доказать.
\end{proof}

Одним из следствий предложения \ref{os} является частичное обобщение второго соотношения
ортогональности для обыкновенных характеров.

\begin{cor} [Слабая ортогональность в блоке]\label{cor sl ort} Пусть $x\in G_{p'}$, $y\in G\setminus G_{p'}$.
Тогда
$$
\sum_{\x\in \irr(B)}\x(x)\ov{\x(y)}=0
$$
для любого $B\in \bl(G)$.
\end{cor}
\begin{proof} Поскольку $\irr(B)$ --- объединение компонент связности графа Брауэра, требуемое следует
из \ref{os}, ввиду того, что $\th_\vf(y)=0$ для всех $\vf$ по \ref{gl ner}$(i)$.
\end{proof}

Напомним, что через $e_\x$ обозначается определённый в \ref{cg} центральный идемпотент алгебры $\CC G$,
соответствующий неприводимому обыкновенному характеру $\x$.
Пусть $\A$ --- произвольное подмножество из $\irr(G)$. Определим центральный идемпотент алгебры $\CC G$
$$\gls{fsA}=\sum_{\x\in\A}e_\x.$$
Поскольку классовые суммы образуют базис центра $\Z(\CC G)$, для всякого $\A\se\irr(G)$
однозначно определены константы
$\gls{asAlKr}\in \CC$,
где $K\in \K(G)$, такие, что
$$
f_{\A}^{\phantom{A}}=\sum_{K\in\K(G)}\a_{\A}^{\phantom{A}}(K)\wh K.
$$

\begin{pre} \label{id con} Пусть $\A$ --- объединение компонент связности графа Брауэра.
Тогда для всех $K\in \K(G)$ имеем
\begin{list}{{\rm(}{\it\roman{enumi}\/}{\rm)}}
{\usecounter{enumi}\setlength{\parsep}{2pt}\setlength{\topsep}{5pt}\setlength{\labelwidth}{23pt}}
\item $\a_{\A}^{\phantom{A}}(K)=\displaystyle\frac{1}{|G|} \sum_{\x\in \A}\x(1)\ov{\x(x_{\mbox{}_K})}$;
\item $\a_{\A}^{\phantom{A}}(K)=0$, если $K\not\in \K(G_{p'})$;
\item $\a_{\A}^{\phantom{A}}(K)\in \wt{Z}$.
\end{list}
\end{pre}
\begin{proof} $(i)$ Это равенство является прямым следствием предложения \ref{id dec}.

$(ii)$ Из $(i)$ и предложения \ref{os} вытекает
$$
\a_{\A}(K)=\frac{1}{|G|} \sum_{\x\in \A}\x(1)\ov{\x(x_{\mbox{}_K})}=
\frac{1}{|G|}\sum_{\vf\in \B}\vf(1)\ov{\th_\vf(x_{\mbox{}_K})},
$$
где $\B=\{\vf\in \iBr(G)\,|\,d_{\x\vf}\ne 0\ \mbox{для некоторого}\ \x\in \A\}$.
Если $K\not\in \K(G_{p'})$, то в силу \ref{gl ner}$(i)$ имеем $\th_\vf(x_{\mbox{}_K})=0$ для всех $\vf\in \B$.

$(iii)$ Ввиду $(ii)$ можем считать, что $K\in \K(G_{p'})$. Тогда из $(i)$ и предложения \ref{os} следует, что
$$
\a_{\A}(K)=\frac{1}{|G|} \sum_{\x\in \A}\x(1)\ov{\x(x_{\mbox{}_K})}=
\frac{1}{|G|}\sum_{\vf\in \B}\ov{\vf(x_{\mbox{}_K})}\th_\vf(1).
$$
Заметим, что по \ref{gl ner}$(vii)$ для всех $\vf\in \B$ имеет место делимость $|G|_p\, \big| \,\th_\vf(1)$
и значения $\ov{\vf(x_{\mbox{}_K})}$ лежат в $\ov{\ZZ}$. Поэтому $\a_{\A}(K)$ является частным от деления
числа из $\ov{\ZZ}$ на $|G|_{p'}\in \ZZ\setminus M$ (см. предложение \ref{m int z}). Отсюда следует требуемое.
\end{proof}

\begin{pre} \label{min un} Если $\A\se\irr(G)$ и $f_{\A}^{\phantom{A}}\in \wt{Z} G$,
то $\A$ является объединением некоторых $p$-блоков обыкновенных характеров.
\end{pre}
\begin{proof}
Мы отмечали в \ref{zwtr}, что отображение $^\st$ можно поднять
до гомоморфизма колец $^\st:\Z(\wt{Z} G)\to \Z(FG)$.
Заметим, что в силу \ref{zsr} имеет место включение $f_{\A}^{\phantom{A}}\in \Z(\CC G)\cap \wt{Z} G=\Z(\wt{Z} G)$.
Рассмотрим значения на $f_\A^{\ \,\st_{\vphantom{A}}}$
центральных характеров $\l_{\x}$, $\x\in \irr(G)$.
Из \ref{com t} следует, что
$$
\l_\x(f_\A^{\ \,\st_{\vphantom{A}}})=\om_\x(f_\A^{\phantom{A}})^\st=
\sum_{\t\in \A}\om_\x(e_\t)^\st=\left\{\ba{ll}
1,& \x\in \A,\\
0,& \x\not\in \A,
\ea\right.\myeqno\label{lx}
$$
где последнее равенство вытекает из \ref{cent dif}$(i)$.

Предположим, что $\A$ не является объединением $p$-блоков обыкновенных характеров. Тогда существуют $p$-блок
$B\in \bl(G)$ и характеры $\x,\t\in \irr(B)$ такие, что $\x\in \A$ и $\t\not \in \A$. Из \ref{lx} следует,
что $\l_\x(f_\A^{\ \,\st_{\vphantom{A}}})=1$ и $\l_\t(f_\A^{\ \,\st_{\vphantom{A}}})=0$. Но $\l_\x=\l_\t$,
поскольку $\x$ и $\t$ лежат в одном $p$-блоке. Противоречие. Отсюда следует требуемое утверждение.
\end{proof}

Теперь для $p$-блоков обыкновенных характеров кроме определения мы получаем ещё две характеризации.

\begin{cor} \label{bl ob} Имеем
\begin{list}{{\rm(}{\it\roman{enumi}\/}{\rm)}}
{\usecounter{enumi}\setlength{\parsep}{2pt}\setlength{\topsep}{5pt}\setlength{\labelwidth}{23pt}}
\item множество $\{\irr(B)\,|\,B\in \bl(G)\}$ совпадает с множеством компонент связности графа Брауэра;
\item $p$-блоки обыкновенных характеров --- это, в точности, минимальные непустые подмножества $\A\se\irr(G)$,
для которых $f_{\A}\in \wt{Z} G$.
\end{list}
\end{cor}
\begin{proof} $(i)$ В силу \ref{gr bl} достаточно показать, что любая компонента связности $\A$
графа Брауэра является объединением $p$-блоков обыкновенных характеров. Из \ref{id con}$(iii)$
следует, что $f_\A^{\phantom{A}}\in \wt{Z} G$. Поэтому требуемое вытекает из \ref{min un}.

$(ii)$ Это вытекает из $(i)$, \ref{id con}$(iii)$ и \ref{min un}.
\end{proof}

\begin{pre} \label{l ch} Пусть $B\in \bl(G)$. Тогда
\begin{list}{{\rm(}{\it\roman{enumi}\/}{\rm)}}
{\usecounter{enumi}\setlength{\parsep}{2pt}\setlength{\topsep}{5pt}\setlength{\labelwidth}{23pt}}
\item для любого $x\in G_{p'}$
$$\sum_{\x\in\irr(B)}\ov{\x(x)}\,\x=\sum_{\vf\in\iBr(B)}\ov{\vf(x)}\,\th_\vf;$$
\item для любого $x\in G$
$$\sum_{\x\in\irr(B)}\ov{\x(x)}\,\hat\x=\sum_{\vf\in\iBr(B)}\ov{\th_\vf(x)}\,\vf;$$
\item для любого $x\in G_{p'}$
$$\sum_{\vf\in\iBr(B)}\ov{\vf(x)}\,\wh{\th_\vf}=\sum_{\vf\in\iBr(B)}\ov{\th_\vf(x)}\,\vf.$$
\end{list}
\end{pre}
\begin{proof} В силу \ref{bl ob}$(i)$ получаем, что $\irr(B)$ --- компонента связности графа Брауэра и
$$
\{\vf\in \iBr(G)\,|\,d_{\x\vf}\ne 0\ \mbox{для некоторого}\ \x\in B\}=\iBr(B)
$$
в силу \ref{mod bl}$(ii)$. Поэтому, положив $\A=B$ в \ref{os}, получаем $(i)$ и $(ii)$. Рассмотрев
образ относительно отображения $\,\hat{\mbox{}}:\cf(G)\to\cf(G_{p'})$ в $(i)$
и взяв $x\in G_{p'}$ в $(ii)$, получим $(iii)$.
\end{proof}

\textzam{Отметим, что равенство \ref{l ch}$(ii)$ обобщает слабую ортогональность
в блоке \ref{cor sl ort}, поскольку в случае $x\in G\setminus G_{p'}$ имеем $\th_\vf(x)=0$ в силу \ref{gl ner}$(i)$.}

Для $B\in \bl(G)$ положим $\gls{fsB}=f_{\irr(B)}$.
Иногда в литературе элементы $f_{\mbox{}_B}$ называют
\glsadd{iIdmpOsm}\mem{идемпотентами Осимы} алгебры $\CC G$. Из \ref{bl ob}$(ii)$ следует, что
$$
f_{\mbox{}_B}\in \Z( \wt{Z} G)\myeqno\label{fbz}
$$
для любого $B\in \bl(G)$.

\textupln{id os prim prf}{\label{id os prim}Показать, что идемпотенты $f_{\mbox{}_B}$, где $B\in \bl(G)$,
алгебры $\Z( \wt{Z} G)$ являются примитивными.}

Опишем строение групповой алгебры $\wt{Z} G$.

\begin{pre} \label{rtg} Имеем
\begin{list}{{\rm(}{\it\roman{enumi}\/}{\rm)}}
{\usecounter{enumi}\setlength{\parsep}{2pt}\setlength{\topsep}{5pt}\setlength{\labelwidth}{23pt}}
\item справедливо разложение
$$
\wt{Z} G= \bigoplus_{B\in \bl(G)}f_{\mbox{}_B} \wt{Z} G;
$$
\item $\Z(f_{\mbox{}_B} \wt{Z} G)=f_{\mbox{}_B}\Z(\wt{Z} G)$.
\end{list}
\end{pre}
\textupl{rtg prf}{Доказать предложение \ref{rtg}.}

Пусть $\bl(G)=\{B_1,\ld,B_k\}$. Упорядочив множества $\irr(G)$ и $\iBr(G)$ таким образом, чтобы вначале шли
характеры из $B_1$, затем из $B_2$, и т.\,д., в соответствии с предложением \ref{mod bl}  мы можем записать
матрицу  $\DD$ в виде
$$
\DD=
\left(
  \begin{array}{cccc}
    \DD_{B_1} & 0 & \ld & 0 \\
    0 & \DD_{B_2} & \ld & 0 \\
    \ld & \ld & \ld &\ld \\
    0 & 0 & \ld & \DD_{B_k}
  \end{array}
\right), \myeqno\label{zv22}
$$
где $\DD_{B_i}$ обозначает подматрицу, стоящую на пересечении строк и столбцов, соответствующих характерам из
блока $B_i$. Будем называть \gls{DsB}
\glsadd{iMatDecPBlk}\mem{матрицей разложения $p$-блока} $B_i$, а её элементы ---
\glsadd{iNumsDecPBlk}\mem{числами разложения $p$-блока} $B_i$.

\begin{pre} \label{bl dec} Для любого $B\in \bl(G)$
\begin{list}{{\rm(}{\it\roman{enumi}\/}{\rm)}}
{\usecounter{enumi}\setlength{\parsep}{2pt}\setlength{\topsep}{5pt}\setlength{\labelwidth}{23pt}}
\item ранг матрицы $\DD_B$ равен $|\iBr(B)|$;
\item $|\irr(B)|\ge|\iBr(B)|$.
\end{list}
\end{pre}

\textupl{bl dec prf}{Доказать предложение \ref{bl dec}.}

Пусть $\DD$ записана в виде \ref{zv22}. Тогда матрица Картана примет блочно-диагональный вид
$$
\C=\DD^\top \DD=\left(
  \begin{array}{cccc}
    \C_{B_1} & 0 & \ld & 0 \\
    0 & \C_{B_2} & \ld & 0 \\
    \ld & \ld & \ld &\ld \\
    0 & 0 & \ld & \C_{B_k}
  \end{array}
\right),
$$
где
$\gls{CsB}=(\DD_{B_i})^\top \DD_{B_i}^{\phantom\top}$ --- \glsadd{iMatCartPBlk}\mem{матрица Картана $p$-блока} $B_i$.
Покажем, что для любого $p$-блока его матрица разложения и матрица Картана не представимы в блочно-диагональном виде.

\begin{pre} \label{m idec}
Для любого $B\in\bl(G)$ матрица разложения $\DD_B$ и матрица Картана $\C_B$
не могут быть представлены в блочном виде
$$
\left(
  \begin{array}{cc}
    * & 0  \\
    0 & *
  \end{array}
\right)
$$
путём перестановки характеров из $\irr(B)$ и $\iBr(B)$, которыми индексируются строки и
столбцы матриц $\DD_B$ и $\C_B$.
\end{pre}
\begin{proof} Пусть для некоторого $B\in\bl(G)$ матрица $\DD_B$ имеет указанный вид и $\irr(B)=\A_1\cup\A_2$
--- соответствующее разбиение характеров, индексирующих строки $\DD_B$.
Тогда для любой пары характеров $\x_1\in \A_1$ и $\x_2\in \A_2$
не существует брауэрова характера $\vf\in \iBr(G)$ такого, что  $d_{\x_1\vf}\ne 0$ и $d_{\x_2\vf}\ne 0$.
Значит, $\A_1$ и $\A_2$ лежат в разных компонентах связности графа Брауэра вопреки \ref{bl ob}$(i)$.

Теперь предположим, что для некоторого $B\in\bl(G)$ матрица Картана  $\C_B$ имеет указанный вид
и $\iBr(B)=\B_1\cup\B_2$ --- соответствующее разбиение брауэровых характеров. Тогда для
любых $\vf_1\in \B_1$ и $\vf_2\in \B_2$ выполнено $0=c_{\vf_1\vf_2}=\sum_{\x\in \irr(B)}d_{\x\vf_1}d_{\x\vf_2}$,
т.\,е. $\vf_1$-й и $\vf_2$-й столбцы матрицы $\DD_B$ ортогональны. Рассмотрим столбцы матрицы $\DD_B$,
индексированные характерами из $\B_1$. Эти столбцы образуют матрицу, ненулевые строки которой
индексированы некоторым подмножеством $\A_1\se\irr(B)$, а нулевые --- подмножеством $\A_2\se\irr(B)$.
Заметим, что $\A_1\ne \varnothing$, т.\,к. иначе
$\DD_B$ имеет нулевой столбец, вопреки \ref{br e}$(i)$.
Можно считать, что матрица $\DD_B$ имеет блочный вид

$$
\ba{c@{}c}
 &
\scriptstyle\B_1 \qquad \scriptstyle\B_2 \\
\ba{c}
\scriptstyle\A_1 \\
\scriptstyle\A_2
\ea
&
\left(
\ba{c|c}
D_{11} &  D_{12} \\
\hline
\O^{\vphantom{A^A}} & D_{22}
\ea
\right)
\ea
$$
Покажем, что  $D_{12}$ ---
нулевая подматрица. В самом деле, если в какой-то её $\x$-й строке есть ненулевой элемент $d_{\x\vf_2}$, то
поскольку $\x$-я строка из $D_{11}$ ненулевая, существует $\vf_1\in \B_1$
такой, что $d_{\x\vf_1}\ne 0$ и, значит, $c_{\vf_1\vf_2}\ge d_{\x\vf_1}d_{\x\vf_2}>0$, вопреки
ортогональности столбцов.
Таким образом, в случае $\A_2\ne \varnothing$ матрица $\DD_B$ может быть приведена
к блочно-диагональному виду, что противоречит доказанному выше, а если $\A_2=\varnothing$, то $\DD_B$
имеет нулевой столбец из $D_{12}$ вопреки \ref{br e}$(i)$.
\end{proof}

Разбиение неприводимых характеров на блоки позволяет доказать
следующее усиление утверждения \ref{z lin}$(ii)$.

\begin{pre} \label{z lin bl} Пусть $B\in \bl(G)$ и $\vf\in \iBr(B)$.
Тогда $\vf$ является $\ZZ$-линейной
комбинацией брауэровых характеров $\hat\x$, где $\x\in \irr(B)$.
\end{pre}
\begin{proof} В силу \ref{z lin}$(ii)$ для некоторых $a_\x\in \ZZ$ мы можем записать
$$\vf=\sum_{\x\in \irr(G)}a_{\x}\hat\x=\vf_{\mbox{}_B}+\vf_{\mbox{}_0},$$
где
$$
\vf_{\mbox{}_B}=\sum_{\x\in \irr(B)}a_{\x}\hat\x,\qquad \vf_{\mbox{}_0}=\sum_{\x\in \irr(G)\setminus\irr(B)}a_{\x}\hat\x
$$
Выразив $\hat \x$ через брауэровы характеры, ввиду \ref{mod bl} получаем, что  $\vf_{\mbox{}_B}$ является
$\ZZ$-линейной комбинацией характеров из $\iBr(B)$, а $\vf_{\mbox{}_0}$ --- характеров из $\iBr(G)\setminus\iBr(B)$.
Значит, из линейной независимости брауэровых характеров следует, что равенство $\vf-\vf_{\mbox{}_B}=\vf_{\mbox{}_0}$
возможно лишь если обе части равны нулю. Отсюда получаем требуемое.
\end{proof}


\bigskip

\section{Блоки групповой алгебры}

Изучим подробнее связь примитивных идемпотентов алгебр $\Z(\CC G)$ и $\Z(F G)$.

Согласно \ref{id con}$(iii)$ для произвольного $B\in \bl(G)$ идемпотент $f_{ \mbox{}_B}$ лежит в $\Z(\wt{Z}
G)$. Назовём \glsadd{iIdmpCenFGCorr}\mem{центральным идемпотентом алгебры}
$FG$\mem{, соответствующим $p$-блоку} $B$, элемент $\gls{esB}=f_{\mbox{}_B}^{\ \st_{\vphantom{A}}}$.

\begin{pre} \label{e b} Справедливы следующие утверждения.
\begin{list}{{\rm(}{\it\roman{enumi}\/}{\rm)}}
{\usecounter{enumi}\setlength{\parsep}{2pt}\setlength{\topsep}{5pt}\setlength{\labelwidth}{23pt}}
\item Элементы $e_{\mbox{}_B}$, где $B\in \bl(G)$, являются попарно ортогональными центральными
идемпотентами алгебры $FG$.
\item $e_{\mbox{}_B}$ --- $F$-линейная комбинация классовых сумм $p$-регулярных классов.
\item $\l_{\mbox{}_{B\vphantom{'}}}(e_{\mbox{}_{B'}})=\d_{\mbox{}_{B,B'}}$ для любых $B,B'\in \bl(G)$.
\item Множество  $\{e_{\mbox{}_B}\,|\,B\in \bl(G)\}$ линейно независимо над $F$. В частности, все $e_{\mbox{}_B}$
ненулевые и попарно  различны.
\item $\sum_{\scriptscriptstyle{B\in \bl(G)}} e_{\mbox{}_B}=1$.
\item Пусть $z\in \Z(FG)$. Если $\l_{\mbox{}_B}(z)=0$ для всех $B\in \bl(G)$, то $z$ --- нильпотентный элемент.
\item Все гомоморфизмы $F$-алгебр $\Z(FG)\to F$ исчерпываются множеством $\{\l_{\mbox{}_B}\,|\, B\in \bl(G)\}$.
\item Любой идемпотент алгебры $\Z(FG)$ является суммой некоторых различных $e_{\mbox{}_B}$.
\item Все  $e_{\mbox{}_B}$ являются примитивными идемпотентами алгебры $\Z(FG)$.
\end{list}
\end{pre}
\begin{proof} $(i)$ Это следует из того, что $^\st:\Z(\wt{Z} G)\to \Z(F G)$ --- эпиморфизм колец и идемпотенты
$f_{\mbox{}_B}\in \Z(\wt{Z} G)$ попарно ортогональны.

$(ii)$ Это вытекает из \ref{id con}$(ii)$.

$(iii)$ Пусть $\x\in \irr(G)$. Из \ref{cent dif}$(i)$ следует, что значение $\om_{\x}(f_{\mbox{}_{B'}})$ равно $1$
если $\x\in \irr(B')$ и $0$ если $\x\not\in \irr(B')$. Поэтому, выбрав $\x$ из $\irr(B)$, ввиду \ref{com t} получим
$$\l_{\mbox{}_{B\vphantom{'}}}(e_{\mbox{}_{B'}})=\l_{\x}(f_{\mbox{}_{B'}}^{\ \,\st})
=\om_{\x}(f_{\mbox{}_{B'}})^\st=\d_{\mbox{}_{B,B'}}.$$

$(iv)$ Это вытекает из $(iii)$.

$(v)$ Это следует из равенств
$\sum_{B\in \bl(G)} f_{\mbox{}_B}=\sum_{\x\in \irr(G)} e_{\mbox{}_\x}=1$
ввиду того, что отображение $^\st:\Z(\wt{Z} G)\to \Z(F G)$ является гомоморфизмом колец.

$(vi)$ Покажем, что $\X(z)=\O$ для любого неприводимого представления $\X$ алгебры $FG$. Пусть
$\vf$ --- брауэров характер $\X$, рассматриваемого как (неприводимое) $F$-представление группы
$G$ и пусть $B$ --- $p$-блок, содержащий $\vf$. По условию имеем
$$
\X(z)=\l_\vf(z)\I=\l_{\mbox{}_B}(z)\I=0.
$$
Значит, $z$ лежит в ядре любого неприводимого представления $\X$ алгебры $FG$ и поэтому $z\in
\J(FG)$ ввиду \ref{int krs}$(i)$. Из предложения \ref{rad nlp} следует, что $z$ --- нильпотентный
элемент.

$(vii)$ Пусть $\l:\Z(FG)\to F$ --- гомоморфизм $F$-алгебр. Тогда ядро $\Ker\l$ является подпространством в
$\Z(FG)$ коразмерности $1$ и $\Z(FG)=\Ker \l \oplus F 1$. Поэтому $\l$ однозначно задаётся своим ядром.
Если $\l\ne\l_{\mbox{}_B}$, то $\Ker \l_{\mbox{}_B}\nsubseteq \Ker\l$, и существует элемент
$z_{\mbox{}_B}\in \Ker \l_{\mbox{}_B}\setminus \Ker\l$.

Предположим, что $\l\not\in\{\l_{\mbox{}_B}\,|\, B\in \bl(G)\}$. Рассмотрим элемент
$$z=\prod_{\mathllap{B}\in \bl(G\mathrlap{)}} z_{\mbox{}_B} \in \Z(FG).$$
Тогда $\l_{\mbox{}_B}(z)=0$ для всех $B\in \bl(G)$, и из $(vi)$  следует, что $z$ --- нильпотентный элемент.
Отсюда вытекает, что $\l(z)=0$ как нильпотентный элемент поля $F$. Но это противоречит тому, что
$\l(z)=\prod \l(z_{\mbox{}_B})\ne 0$. Отсюда следует требуемое.

$(viii)$ Пусть $e\in \Z(FG)$ --- идемпотент. Предположим, что $ee_{\mbox{}_B}\ne 0$ для некоторого $B\in \bl(G)$.
Покажем, что в этом случае $ee_{\mbox{}_B}=e_{\mbox{}_B}$.

Так как элемент $ee_{\mbox{}_B}$ --- ненулевой идемпотент, он не является нильпотентным. С другой стороны
для любого $p$-блока $B'$, отличного от $B$, имеем
$\l_{\mbox{}_{B'}}(ee_{\mbox{}_{B}})=\l_{\mbox{}_{B'}}(e)\l_{\mbox{}_{B'}}(e_{\mbox{}_{B}})=0$.
Поэтому из $(vi)$ следует, что $\l_{\mbox{}_{B}}(ee_{\mbox{}_B})\ne 0$.
Но $\l_{\mbox{}_{B}}(ee_{\mbox{}_B})$ --- идемпотент из $F$, значит, $\l_{\mbox{}_{B}}(ee_{\mbox{}_B})=1$.
Получаем, что $\l_{\mbox{}_{B}}(e_{\mbox{}_B}-ee_{\mbox{}_B})=1-1=0$ и
$\l_{\mbox{}_{B'}}(e_{\mbox{}_B}-ee_{\mbox{}_B})=\l_{\mbox{}_{B'}}(1-e)\l_{\mbox{}_{B'}}(e_{\mbox{}_B})=0$
для всех $B'\ne B$. Тогда из $(vi)$ следует, что элемент $e_{\mbox{}_B}-ee_{\mbox{}_B}$ нильпотентный.
Но он, очевидно, также идемпотент. Поэтому $e_{\mbox{}_B}-ee_{\mbox{}_B}=0$.

Итак, из доказанного следует, что $e=e(\sum e_{\mbox{}_B})=\sum ee_{\mbox{}_B}$ --- сумма
некоторых (различных) $e_{\mbox{}_B}$.

$(ix)$ Пусть идемпотент $e_{\mbox{}_B}$ не является примитивным. Тогда $e_{\mbox{}_B}=e_1+e_2$,
где $e_1,e_2\in \Z(FG)$ --- ненулевые ортогональные идемпотенты и мы получаем
$$e_1=e_1^2=e_1(e_{\mbox{}_B}-e_2)=e_1e_{\mbox{}_B}=e_{\mbox{}_B},$$
где последнее равенство следует из доказательства предыдущего пункта. Но тогда $e_2=e_{\mbox{}_B}-e_1=0$. Противоречие.
\end{proof}

\textzam{Отметим, что в алгебре $FG$ идемпотенты $e_{\mbox{}_B}$ могут не быть примитивными.}

\begin{pre} \label{jzfg}Имеет
место равенство
$$\J(\Z(FG))=\bigcap_{B\in \bl(G)}\Ker \l_{\mbox{}_B}.$$
\end{pre}
\textupr{Доказать предложение \ref{jzfg}.}

Алгебра $FG$ допускает в силу \ref{ids}$(iii)$ разложение в прямую сумму своих двусторонних идеалов
$$
FG=\bigoplus_{B\in \bl(G)}e_{\mbox{}_B}FG, \myeqno\label{zv21}
$$
причём из примитивности $e_{\mbox{}_B}$ как идемпотентов алгебры $\Z(FG)$ и
предложения \ref{inn dir sum}$(iii)$ следует, что идеал $e_{\mbox{}_B}FG$
не представим в виде прямой суммы ненулевых двусторонних идеалов алгебры $FG$.
Мы можем рассматривать \ref{zv21} как обобщение разложения \ref{thm wedd}$(i)$ из теоремы Веддерберна
на случай, когда алгебра $FG$ не обязательно является полупростой. Для $B\in \bl(G)$ будем обозначать соответствующую
неразложимую компоненту $e_{\mbox{}_B}FG$ через
\gls{WsBlFGr} или просто \gls{WsB} и называть её
\glsadd{iBlkAlg}\mem{блоком алгебры} $FG$\mem{, соответствующим $p$-блоку} $B$. Ясно, что блок $\W_B$ является $F$-алгеброй
с единицей $e_{\mbox{}_B}$.

\textzam{Отметим, что в некоторой литературе часто именно компоненты $\W_B$ называют $p$-блоками
группы $G$. В силу \ref{dec uniq} набор блоков $\{\W_B\,|\, B\in \bl(G)\}$ определяется однозначно
алгеброй $FG$ как набор неразложимых двусторонних идеалов.}

\begin{pre} \label{bl z} Имеем
\begin{list}{{\rm(}{\it\roman{enumi}\/}{\rm)}}
{\usecounter{enumi}\setlength{\parsep}{2pt}\setlength{\topsep}{5pt}\setlength{\labelwidth}{23pt}}
\item Алгебра $\Z(FG)$ допускает разложение в прямую сумму своих идеалов
$$
\Z(FG)=\bigoplus_{B\in \bl(G)}\Z(\W_B).
$$
\end{list}
Пусть $B\in \bl(G)$. Тогда выполнены следующие утверждения.
\begin{list}{{\rm(}{\it\roman{enumi}\/}{\rm)}}
{\usecounter{enumi}\setlength{\parsep}{2pt}\setlength{\topsep}{5pt}\setlength{\labelwidth}{23pt}}
\addtocounter{enumi}{1}
\item $\Z(\W_B)=e_{\mbox{}_B}\Z(FG)$.
\item Единственным неприводимым представлением $F$-алгебры $\Z(\W_B)$ является ограничение на неё гомоморфизма
$\l_{\mbox{}_B}$.
\item $\Z(\W_B)=Fe_{\mbox{}_B}\oplus \J(\Z(\W_B))$.
\item $\Z(\W_B)$ --- локальное кольцо.
\item $\J(\Z(\W_B))=e_{\mbox{}_B}\J(\Z(FG))$.
\item Если $I,J$ --- идеалы алгебры $\Z(FG)$ и $e_{\mbox{}_B}\in I+J$, то либо $e_{\mbox{}_B}\in I$
либо $e_{\mbox{}_B}\in J$.
\end{list}
\end{pre}
\begin{proof} $(i)$ Это следует из \ref{z pr}$(i)$.

$(ii)$ Это следует из \ref{ids}$(iv)$.

$(iii)$ Из \ref{cent irr} следует, что всякое неприводимое представление алгебры $\Z(\W_B)$ одномерно.
Пусть $\l:\Z(\W_B)\to F$ --- одно из таких представлений. Ввиду разложения $(i)$ можно продолжить $\l$ до
неприводимого представления $\Z(FG)\to F$. Тогда из  \ref{e b}$(vii)$ следует, что $\l=\l_{\mbox{}_{B'}}$
для некоторого $B'\in \bl(G)$. Но если $B'\ne B$, то $\l_{\mbox{}_{B'}}(\W_B)=0$ в силу \ref{e b}$(iii)$.
Значит, $B'=B$, откуда следует требуемое.

$(iv)$ Пусть $\l$ --- ограничение на $\Z(\W_B)$ гомоморфизма $\l_{\mbox{}_B}$. Так как $\l$ --- единственное
неприводимое представление алгебры $\Z(\W_B)$, то $\J(\Z(\W_B))=\Ker\l$ --- $F$-подпространство в $\Z(\W_B)$
коразмерности
$1$. Поскольку $\l(e_{\mbox{}_B})=1\ne 0$, имеем требуемое разложение $\Z(\W_B)=Fe_{\mbox{}_B}\oplus \Ker\l$.

$(v)$ В силу $(iv)$ радикал $\J(\Z(\W_B))$ имеет в  $\Z(\W_B)$ коразмерность $1$ и содержится в любом
максимальном идеале. Поэтому любой максимальный идеал совпадает с $\J(\Z(\W_B))$.
Следовательно, $\Z(\W_B)$ --- локальное кольцо.

$(vi)$ Из $(ii)$ и \ref{jre} получаем
$$\J(\Z(\W_B))=\J(e_{\mbox{}_B}\Z(FG))=e_{\mbox{}_B}\J(\Z(FG)).$$

$(vii)$ Так как $e_{\mbox{}_B}I$ --- идеал в $\Z(FG)$ и $e_{\mbox{}_B}\in \Z(\W_B)$, то
$e_{\mbox{}_B}I\nor \Z(\W_B)$. Если $e_{\mbox{}_B}I=\Z(\W_B)$, то
$$
e_{\mbox{}_B}\in\Z(\W_B)=e_{\mbox{}_B}I\se I.
$$
Поэтому можем считать, что $e_{\mbox{}_B}I$ --- собственный идеал в $\Z(\W_B)$. Аналогично,
$e_{\mbox{}_B}J$ --- собственный идеал в $\Z(\W_B)$. В силу локальности кольца $\Z(\W_B)$
получаем $e_{\mbox{}_B}I\se\J(\Z(W_B))$ и $e_{\mbox{}_B}J\se\J(\Z(W_B))$. По условию
$e_{\mbox{}_B}=x+y$ для некоторых $x\in I$ и $y\in J$. Поэтому
$$
e_{\mbox{}_B}=e_{\mbox{}_B}^2=e_{\mbox{}_B}x+e_{\mbox{}_B}y\in e_{\mbox{}_B}I+e_{\mbox{}_B}J\se\J(\Z(W_B)).
$$
Но $e_{\mbox{}_B}$ --- единица алгебры. Противоречие.
\end{proof}

Подобно неприводимым обыкновенным и брауэровым характерам группы $G$ неприводимые $\CC G$- и $FG$-модули
естественно разбиваются на классы эквивалентности, соответствующие $p$-блокам группы $G$.
Следующее утверждение даёт критерий принадлежности модуля такому классу.

\begin{pre} \label{ir mod bl} Пусть $B\in \bl(G)$.

$(i)$ Если $M_\x$ --- неприводимый $\CC G$-модуль с обыкновенным характером $\x\in \irr(G)$, то

$$
M_\x f_{\mbox{}_B}=\left\{\ba{ll}
M_\x,& \mbox{если}\ \ \x\in B;\\
0,& \mbox{если}\ \ \x\not\in B.
\ea\right.
$$

$(ii)$ Если $M_\vf$ --- неприводимый $FG$-модуль с брауэровым характером $\vf\in \iBr(G)$, то

$$
M_\vf e_{\mbox{}_B}=\left\{\ba{ll}
M_\vf,& \mbox{если}\ \ \vf\in B;\\
0,& \mbox{если}\ \ \vf\not\in B.
\ea\right.
$$
\end{pre}
\begin{proof} $(i)$ Идемпотент $e_\x$ действует тождественно на $M_\x$ в силу \ref{thm wedd}$(iii)$, поскольку
является единицей компоненты Веддерберна $\W_{M_\x}(\CC G^\circ)$. Значит,
$M_\x=M_\x e_\x$. Осталось заметить, что произведение $e_\x f_{\mbox{}_B}$ равно $e_\x$, если $\x\in B$,
и $0$, если $\x\not\in B$. Отсюда следует требуемое.

$(ii)$ Так как $\sum_{B\in \bl(G)} e_{\mbox{}_B}=1$,  имеем $M_\vf=\bigoplus_{B\in \bl(G)}M_\vf
e_{\mbox{}_B}$ в силу \ref{m id dec}$(iii)$. Из неприводимости модуля $M_\vf$ следует, что существует
единственный $B\in \bl(G)$ для которого $M_\vf=M_\vf e_{\mbox{}_B}$ и $M_\vf e_{\mbox{}_{B'}}=0$ для
всех $B'\ne B$. В частности, $e_{\mbox{}_B}$ действует тождественно на $M_\vf$, т.\,е. если $\X$ ---
$F$-представление, соответствующее $M_\vf$, то $\X(e_{\mbox{}_B})=\I$. Но $e_{\mbox{}_B}\in \Z(FG)$ и,
значит, $\X(e_{\mbox{}_B})=\l_\vf(e_{\mbox{}_B})\I$. Поэтому
$1=\l_\vf(e_{\mbox{}_B})=\l_{\mbox{}_{B'}}(e_{\mbox{}_B})$, где $B'$ --- $p$-блок, содержащий $\vf$. Из
\ref{e b}$(iii)$ следует, что $B'=B$.
\end{proof}

\begin{cor} \label{bl alg} Для всякого $B\in \bl(G)$
брауэров характер любого композиционного фактора блока $\W_B$ принадлежит $B$.
\end{cor}
\begin{proof} Поскольку $\W_B=e_{\mbox{}_B} FG$, идемпотент $e_{\mbox{}_B}$ действует как тождественное
преобразование на блоке $\W_B$, а, значит, и на любом его композиционном факторе $V$. Из \ref{ir mod
bl}$(ii)$ следует, что брауэров характер $FG$-модуля $V$ лежит в $B$.
\end{proof}

Естественно спросить, каковы брауэров характер и размерность блока $W_B$.

\begin{pre} \label{br har bl} Пусть $B\in \bl(G)$.

$(i)$ Брауэров характер $FG$-модуля $\W_B$ равен
$$
\sum_{\x\in \irr(B)}\x(1)\hat\x\,=\sum_{\vf\in \iBr(B)}\th_\vf(1)\vf\,=\sum_{\vf\in \iBr(B)}\vf(1)\wh{\th_\vf}.
$$

$(ii)$ $\dim_F W_B=\sum_{\x\in \irr(B)}\x(1)^2=\sum_{\vf\in \iBr(B)}\th_\vf(1)\vf(1)$.

\end{pre}
\begin{proof} $(i)$ Пусть  $\eta_{\mbox{}_B}$ --- брауэров характер блока $\W_B$. Поскольку
$$
FG=\bigoplus_{B\in \bl(G)}\W_B,
$$
брауэров характер регулярного модуля $FG^\circ$ равен $\sum_{B\in \bl(G)}\eta_{\mbox{}_B}$.

Пусть $\R$ --- регулярное $\CC$-представление группы $G$ относительно базиса из групповых элементов и
$\r$ --- его характер. Поскольку $\R$ записано над $\ZZ\se \wt{Z}$, из \ref{ord br}$(ii)$ получаем,
что брауэров характер регулярного $F$-представления $\R^\st$ группы $G$ равен $\hat\r$.

Итак, в силу \ref{fbern}$(i)$ имеем

$$\sum_{B\in \bl(G)}\eta_{\mbox{}_B}=\hat\r=\sum_{\x\in \irr(G)} \x(1)\hat\x=
\sum_{B\in \bl(G)}\left(\sum_{\x\in \irr(B)} \x(1)\hat\x\right).$$
Поскольку неприводимые компоненты брауэрова характера $\eta_{\mbox{}_B}$ лежат в $\iBr(B)$ по \ref{bl alg},
и неприводимые компоненты брауэрова характера $\sum_{\x\in \irr(B)} \x(1)\hat\x$ лежат в $\iBr(B)$ по \ref{mod bl}$(i)$,
получаем
$$\eta_{\mbox{}_B}=\sum_{\x\in \irr(B)} \x(1)\hat\x$$
для всех $B\in \bl(G)$. Далее, положив $x=1$ в \ref{l ch}$(i),(ii)$, получаем $(i)$.

$(ii)$ Это следует из $(i)$.
\end{proof}

Пусть $B\in \bl(G)$ и $\eta\in \cf(G)$.
Тогда функция
$$
\gls{etasB}=\sum_{\x\in \irr(B)}(\eta,\x)_{\mbox{}_G}\,\x
$$
называется  \glsadd{iBPrtClsFunc}\mem{$B$-частью} классовой функции $\eta$.
Легко видеть, что $\eta=\sum_{B\in \bl(G)}\eta_{\mbox{}_B}$, и
для любого $B\in \bl(G)$  отображение $\eta\mapsto \eta_{\mbox{}_B}$
является $\CC$-линейным преобразованием алгебры $\cf(G)$.

\begin{pre} \label{bch} Пусть $B\in \bl(G)$, $\eta\in \cf(G)$.
Тогда
$$
\eta_{\mbox{}_B}(x)=\eta(f_{\mbox{}_B}x)
$$
для всех $x\in \CC G$.
\end{pre}
\begin{proof}
В силу линейности можно считать, что $\eta\in\irr(G)$.  Тогда
$$
\eta_{\mbox{}_B}=\left\{\ba{ll}
\eta, & \mbox{если}\ \ \eta\in B;\\
0, & \mbox{если}\ \ \eta\not\in B.
\ea\right. \myeqno\label{etb}
$$
Пусть $M=M_\eta$ и $\X=\X_\eta$. Из \ref{ir mod bl}$(i)$ следует, что если $\eta \in B$, то $M=Mf_{\mbox{}_B}$
и, значит, $f_{\mbox{}_B}$ действует на $M$ тождественно, т.\,е. $\X(f_{\mbox{}_B})=\I$,
а если $\eta \not\in B$, то $Mf_{\mbox{}_B}=0$, т.\,е. $\X(f_{\mbox{}_B})=\O$. Поэтому
$$
\eta(f_{\mbox{}_B}x)=\tr \X(f_{\mbox{}_B}x) =\tr \big(\X(f_{\mbox{}_B})\X(x)\big)=\left\{\ba{ll}
\tr \X(x)=\eta(x), & \mbox{если}\ \ \eta\in B;\\
\tr \O=0, & \mbox{если}\ \ \eta\not\in B.
\ea\right.
$$
Сравнив с \ref{etb}, получаем $\eta_{\mbox{}_B}(x)=\eta(f_{\mbox{}_B}x)$.
\end{proof}

\textuprn{Доказать, что отображение $\eta\mapsto\eta_{\mbox{}_B}$ является эндоморфизмом алгебры $\cf(G)$.}

\section{Дефект блока. Дефектная группа}

Существует тесная связь между $p$-блоками и $p$-подгруппами группы $G$, к описанию которой мы приступаем.

Для произвольных подмножеств $X,Y$ группы $G$ определим отношение
\glsadd{iGIncl}\mem{$G$-включения}. Будем говорить, что $X$\ \mem{$G$-содержится} в $Y$ и записывать
\gls{XleGY}, если некоторое $G$-сопряжённое с $X$ подмножество
содержится в $Y$ или, что эквивалентно, $X$ содержится в некотором $G$-сопряжённом с $Y$ подмножестве.

Пусть $K\in\K(G)$ --- произвольный класс сопряжённости. Положим
$$\gls{DplKr}=\bigcup_{x\in K}\Syl_p(\C_G(x)).$$
Как правило мы будем писать \gls{DlKr} вместо $\D_p(K)$. Легко видеть, что множество $\D(K)$ является классом
сопряжённых $p$-подгрупп группы $G$. Его элементы будем называть \glsadd{iGrPDefCls}\mem{$p$-дефектными группами класса}
$K$. Для любого $K\in\K(G)$ зафиксируем некоторый представитель
$\gls{dellKr}\in \D(K)$. Будем также иногда писать \gls{delplKr} вместо $\d(K)$.
\glsadd{iDefCls}\mem{Дефектом класса} $K$  (или, точнее, его \glsadd{iPDefCls}\mem{$p$-дефектом})
  называется число \gls{dK}, определяемое равенством $|\d(K)|=p^{\,\df(K)}$.

Пусть $P$ --- произвольная $p$-подгруппа группы $G$. Мы можем рассмотреть совокупность всех
классов $K\in\K(G)$, $p$-дефектные группы которых $G$-содержатся в $P$. Символом
$\Z_P(FG)$ обозначим  $F$-линейную оболочку классовых сумм таких классов. Другими словами,
$$
\Z_P(FG)=\bigoplus_{\substack{K\in \K(G),\\ \d(K)\le_G P}}F\wh K.
$$

Из данных определений сразу вытекает

\begin{pre} \label{och cor} Справедливы
следующие утверждения.
\begin{list}{{\rm(}{\it\roman{enumi}\/}{\rm)}}
{\usecounter{enumi}\setlength{\parsep}{2pt}\setlength{\topsep}{5pt}\setlength{\labelwidth}{23pt}}
\item Пусть $P,Q\le G$ --- $p$-подгруппы и $P\le_G Q$. Тогда $\Z_P(FG)\se\Z_Q(FG)$.
\item Если $p$-подгруппы $P$ и $Q$ сопряжены в $G$, то $\Z_P(FG)=\Z_Q(FG)$.
\item Для любого $K\in \K(G)$ имеем $\wh K\in \Z_{\d(K)}(FG)$.
\end{list}
\end{pre}

\begin{pre} \label{zp id} Справедливы
следующие утверждения.
\begin{list}{{\rm(}{\it\roman{enumi}\/}{\rm)}}
{\usecounter{enumi}\setlength{\parsep}{2pt}\setlength{\topsep}{5pt}\setlength{\labelwidth}{23pt}}
\item Для любой $p$-подгруппы $P\le G$ имеем $\Z_P(FG)\nor\Z(FG)$.
\item Для любых $p$-подгрупп $P,Q\le G$ выполнено равенство
$$
\Z_P(FG)\cap\Z_Q(FG)=\sum_{x,y\in G}\Z_{P^x\cap Q^y}(FG).
$$
\end{list}
\end{pre}
\begin{proof} $(i)$ Из \ref{kl sum}$(iv)$ следует, что для произвольных классов $K,L\in \K(G)$ соответствующие
классовые суммы из $\Z(FG)$ перемножаются по правилу
$$
\wh{K}\wh{L}=\sum_{M\in \K(G)} a_{\mbox{}_{KLM}}^{\,\st} \whm,
$$
где $a_{\mbox{}_{KLM}}$ --- число всевозможных пар
$(x,y)\in K\times L$ таких, что $xy=z$ для фиксированного (но произвольного) представителя $z\in M$.
Если $z$ выбран так, что $\d(M)\le \C_G(z)$, то на множестве
таких пар $(x,y)$ можно определить действие группы $\d(M)$, положив $(x,y)^d=(x^d,y^d)$ для любого $d\in \d(M)$.
Поскольку $\d(M)$ --- $p$-группа, её орбиты при этом действии будут иметь мощность, равную степени числа $p$.

Предположим, что $a_{\mbox{}_{KLM}}^{\,\st}\ne 0$ для некоторого класса $M$. Из предыдущего замечания
следует, что $\d(M)$ имеет одноэлементную орбиту, т.\,е. существует пара элементов $x\in K$ и $y\in L$,
централизуемых группой $\d(M)$. Другими словами, при нашем предположении имеем
$\d(M)\le_G \d(K)$ и $\d(M)\le_G \d(L)$.

Пусть теперь $K$ --- класс, для которого $\d(K)\le_G P$, и $L$ --- произвольный класс. Из доказанного выше
следует, что произведение $\wh{K}\wh{L}$ является линейной комбинацией классовых сумм $\wh{M}$,
соответствующих классам $M$, для которых, в частности, $\d(M)\le_G\d(K)\le_G P$, т.\,е. $\wh{K}\wh{L}\in \Z_P(FG)$.
Поэтому $\Z_P(FG)$ --- идеал в $\Z(FG)$.

$(ii)$ Включение $\sum_{x,y\in G}\Z_{P^x\cap Q^y}(FG)\se \Z_P(FG)\cap\Z_Q(FG)$ следует из соотношений
$P^x\cap Q^y\le_G P$ и $P^x\cap Q^y\le_G Q$ для всех $x,y\in G$. Докажем обратное включение. Пусть
$z\in \Z_P(FG)\cap\Z_Q(FG)$. Тогда
$$z=\sum_{K\in \K(G)} a_{\mbox{}_K}\wh K$$
для некоторых $a_{\mbox{}_K}\in F$. Поскольку $z\in\Z_P(FG)$, получаем, что если
$a_{\mbox{}_K}\ne 0$,  то $\d(K)\le_G P$. Аналогично, если $a_{\mbox{}_K}\ne 0$,  то $\d(K)\le_G Q$.
То есть для каждого $K\in \K(G)$ такого, что $a_{\mbox{}_K}\ne 0$, найдутся $x,y\in G$, для которых
$\d(K)\se P^x\cap Q^y$. Но тогда $\wh K\in\Z_{P^x\cap Q^y}(FG)$ и, значит, $z\in \sum_{x,y\in G}\Z_{P^x\cap Q^y}(FG)$,
что и требовалось доказать.
\end{proof}

Как соотносятся идеалы алгебры $\Z(FG)$ вида $\Z_P(FG)$ с идеалами $\Z(W_B)$?
 Оказывается, что для каждого $B\in \bl(G)$
существует единственный класс сопряжённых минимальных по включению $p$-подгрупп $P$, для которых $\Z(W_B)$ лежит в
$\Z_P(FG)$.

\begin{pre} \label{def ex} Пусть $B\in \bl(G)$. Существует единственный класс сопряжённых $p$-подгрупп $P$,
для которых выполнены следующие два условия.
\begin{list}{{\rm(}{\it\roman{enumi}\/}{\rm)}}
{\usecounter{enumi}\setlength{\parsep}{2pt}\setlength{\topsep}{5pt}\setlength{\labelwidth}{23pt}}
\item $\Z(W_B)\se \Z_P(FG)$;
\item если $\Z(W_B)\se \Z_Q(FG)$ для некоторой $p$-подгруппы $Q$, то $P\le_G Q$.
\end{list}
\end{pre}

\textupl{def ex prf}{Доказать предложение \ref{def ex}. \uk{Воспользоваться \ref{bl z}$(vii)$.}}

\begin{opr} \label{opr def gr}
Пусть $B\in \bl(G)$.
Класс сопряжённых $p$-подгрупп, существование и единственность которого
утверждается в \ref{def ex}, будем обозначать через
\gls{DlBr}, а его элементы называть \glsadd{iGrDefPBlk}\mem{дефектными группами $p$-блока} $B$.
Для любого $B\in\bl(G)$ зафиксируем некоторый представитель $\gls{dellBr}\in \D(B)$.
\end{opr}

Из определения следует

\begin{pre} \label{dg cor} Для любого $B\in \bl(G)$
 имеем $e_{\mbox{}_B}\in \Z_{\d(B)}(FG)$.
\end{pre}

Ввиду важности введённого понятия, мы дадим другое, эквивалентное определение дефектной группы $p$-блока.

Пусть $B\in \bl(G)$. Запишем
$$
e_{\mbox{}_B}=\sum_{K\in\K(G)}a_{\mbox{}_B}(K)\wh K, \myeqno\label{zv6}
$$
где $a_{\mbox{}_B}:\K(G)\to F$ --- однозначно определённое отображение. Из \ref{id con}$(i),(iii)$ следует, что
$$
a_{\mbox{}_B}(K)=\left(\frac{1}{|G|}
\sum_{\x\in \irr(B)}\x(1)\ov{\x(x_{\mbox{}_K})}\right)^\st. \myeqno\label{abk}
$$
Кроме того, в силу \ref{id con}$(ii)$ имеем
$a_{\mbox{}_B}(K)=0$ если $K\not\in \K(G_{p'})$. Поэтому
$$
e_{\mbox{}_B}=\sum_{K\in\K(G_{p'})}a_{\mbox{}_B}(K)\wh K. \myeqno\label{zzv7}
$$
Поскольку
$$
1=\l_{\mbox{}_B}(e_{\mbox{}_B})=\sum_{K\in\K(G_{p'})}a_{\mbox{}_B}(K)\l_{\mbox{}_B}(\wh K),
$$
существует по меньшей мере один класс $K\in \K(G)$ для которого $a_{\mbox{}_B}(K)\ne 0$ и $\l_{\mbox{}_B}(\wh K)\ne 0$.
Любой такой $K$ мы назовём \glsadd{iClsDefPBlk}\mem{дефектным классом $p$-блока} $B$.
Из определения сразу вытекает

\begin{pre} \label{dk cor} Дефектный
класс любого $p$-блока является $p$-регулярным.
\end{pre}

Мы покажем, что $p$-дефектные группы дефектных классов $p$-блока $B$ --- это в точности его дефектные группы.

\begin{thm}[min\;\!-\!-max] \label{thm min max}\glsadd{iThmMinMax} Пусть $B\in \bl(G)$ и $K\in \K(G)$.
\begin{list}{{\rm(}{\it\roman{enumi}\/}{\rm)}}
{\usecounter{enumi}\setlength{\parsep}{2pt}\setlength{\topsep}{5pt}\setlength{\labelwidth}{23pt}}
\item Если $\l_{\mbox{}_B}(\wh K)\ne 0$, то $\d(B)\le_G\d(K)$.
\item Если $a_{\mbox{}_B}(K)\ne 0$, то $\d(K)\le_G\d(B)$.
\end{list}
\end{thm}
\begin{proof} $(i)$ Рассмотрим идеалы $\Z_{\d(K)}(FG)$ и $\Ker \l_{\mbox{}_B}$
алгебры $\Z(FG)$. Имеем $\wh K \in \Z_{\d(K)}(FG)$, но с другой стороны
$\wh K\not \in \Ker \l_{\mbox{}_B}$, поскольку $\l_{\mbox{}_B}(\wh K)\ne 0$.
Значит, $\Z_{\d(K)}(FG)\nsubseteq\Ker\l_{\mbox{}_B}$. Однако идеал $\Ker\l_{\mbox{}_B}$ имеет
коразмерность $1$ в $\Z(FG)$. Поэтому
$$
\Z(FG)=\Z_{\d(K)}(FG)+\Ker\l_{\mbox{}_B}.
$$
В частности, из \ref{bl z}$(vii)$ следует, что идемпотент $e_{\mbox{}_B}$ лежит либо в $\Z_{\d(K)}(FG)$,
либо в $\Ker\l_{\mbox{}_B}$. Но последнее невозможно, поскольку $\l_{\mbox{}_B}(e_{\mbox{}_B})=1\ne 0$.
Значит, $e_{\mbox{}_B}\in \Z_{\d(K)}(FG)$, откуда получаем $\W_B=e_{\mbox{}_B}\Z(FG)\se \Z_{\d(K)}(FG)$.
Поэтому $\d(B)\le_G\d(K)$.

$(ii)$ Из условия следует, что классовая сумма $\wh K$ входит в разложение \ref{zzv7} с ненулевым коэффициентом.
Однако $e_{\mbox{}_B}$ лежит в идеале $\Z_{\d(B)}(FG)$. По определению этого идеала получаем, что $\d(K)\le_G\d(B)$.
\end{proof}

\begin{cor} \label{cor bl kl gr} Для любого $B\in \bl(G)$
и произвольного дефектного класса $K$ блока $B$ имеем $\D(B)=\D(K)$.
\end{cor}
\begin{proof} Пусть $K$ --- дефектный класс $p$-блока $B$. Тогда $a_{\mbox{}_B}(K)\ne 0$
и $\l_{\mbox{}_B}(\wh K)\ne 0$. Из теоремы \ref{thm min max} следует, что группы $\d(K)$ и $\d(B)$ сопряжены,
т.\,е. $\d(K)\in \D(B)$ и $\d(B)\in \D(K)$. \end{proof}

Какие $p$-подгруппы группы $G$ могут быть дефектными группами её $p$-блоков?
Оказывается,\footnote{См.\ref{tgr}} что любая
дефектная группа имеет вид $P\cap Q$ для некоторых $P,Q\in \Syl_p(G)$.
Нашей ближайшей целью является установить более
слабый результат, который утверждает, что дефектная группа произвольного $p$-блока содержит любую
нормальную $p$-подгруппу группы $G$.

\normalmarginpar

Докажем вспомогательное утверждение.

\begin{pre} \label{k nil} Пусть $K\in \K(G)$ и $K\nsubseteq \C_G(\OO_p(G))$.
Тогда $\wh K$ --- нильпотентный элемент  алгебры $FG$.
\end{pre}
\begin{proof} Обозначим $P=\OO_p(G)$. Группа $P$ действует сопряжением на $K$. Пусть $O$ --- одна из
орбит этого действия. Заметим, что $O$ содержится в одном смежном классе $G$ по $P$, поскольку
любые два её представителя сопряжены элементом из $P$ и, значит, имеют одинаковые образы в факторгруппе $G/P$.
Кроме того, из условия следует, что  $K\cap \C_G(\OO_p(G))=\varnothing$. Поэтому
$|O|>1$ и, в частности, $p\bigm| |O|$.

Пусть $\X$ --- неприводимое $F$-представление группы $G$. Ввиду \ref{cor op ker}$(ii)$ имеем $P\se \ker \X$,
т.\,е. $\X$ принимает постоянное значение на каждом смежном классе $G$ по $P$.
Из сделанных выше замечаний получаем
$$
\sum_{y\in O}\X(y)=|O|^\st\X(x)=0,
$$
поскольку $p\bigm| |O|$, где $x\in O$ --- произвольный представитель. Просуммировав по всем орбитам, имеем $\X(\wh K)=0$ для любого неприводимого
представления $\X$. Поэтому элемент $\wh K$ лежит в $\J(FG)$ в силу предложения \ref{int krs}$(i)$ и, значит, является
нильпотентным.
\end{proof}

\begin{pre} \label{op def} Для любого $B\in \bl(G)$ имеем $\OO_p(G)\le\d(B)$.
\end{pre}
\begin{proof} Пусть $K$ --- дефектный класс блока $B$. Тогда $\l_{\mbox{}_B}(\wh K)\ne 0$ и, значит, элемент
$\wh K$ не является нильпотентным. Тогда  из \ref{k nil} следует, что $K\se \C_G(\OO_p(G))$. Поэтому
$\OO_p(G)\se \d(K)$. Но $\d(K)\in \D(B)$ ввиду \ref{cor bl kl gr}, и мы получаем требуемое.
\end{proof}

\begin{opr} \label{opr def} Пусть $B\in \bl(G)$ и $|\d(B)|=p^d$.
Число $d$, которое мы обозначим через
\gls{dD}, называется \glsadd{iDefPBlk}\mem{дефектом $p$-блока} $B$.
\end{opr}

Это определение корректно, поскольку все дефектные группы $p$-блока $B$ сопряжены и, значит,
имеют одинаковые порядки. Мы покажем, что дефект $\df(B)$ можно определить, зная лишь степени
характеров из $\irr(B)$ или $\iBr(B)$.

\textuprn{Доказать, что если $\OO_p(G)\ne 1$, то $G$ не имеет $p$-блоков дефекта $0$.}

\textuprn{Каковы дефекты блоков абелевой группы?}


\begin{pre} \label{def deg} Пусть $B\in \bl(G)$ и $d=\df(B)$.
Запишем $|G|_p=p^a$. Тогда
\begin{list}{{\rm(}{\it\roman{enumi}\/}{\rm)}}
{\usecounter{enumi}\setlength{\parsep}{2pt}\setlength{\topsep}{5pt}\setlength{\labelwidth}{23pt}}
\item $p^{a-d}$ --- максимальная степень числа $p$, делящая $\x(1)$ для
всех $\x\in \irr(B)$;
\item $p^{a-d}$ --- максимальная степень числа $p$, делящая $\vf(1)$ для
всех $\vf\in \iBr(B)$.
\end{list}
\end{pre}
\begin{proof} $(i)$ Пусть $K$ --- дефектный класс $p$-блока $B$.
В силу \ref{cor bl kl gr} имеем $\d(K)\in \D(B)$, и поэтому из условия получаем
$$p^d=|\d(K)|=|\C_G(x_{\mbox{}_K})|_p=\frac{|G|_p}{|K|_p}=\frac{p^a}{|K|_p},$$
т.\,е. $|K|_p=p^{a-d}$.

Пусть $\x\in \irr(B)$. Так как $K$ --- дефектный класс $p$-блока $B$, то
$$
0\ne \l_{\mbox{}_B}(\wh K)=\l_\x(\wh K)=\om_\x(\wh K)^\st=
\left(\frac{|K|\x(x_{\mbox{}_K})}{\x(1)}\right)^\st. \eqno{(*)}
$$
Но $\x(x_{\mbox{}_K})\in \ov{\ZZ}$ и, значит,
$\displaystyle\frac{|K|}{\x(1)}\not\in \wtm$ в силу $(*)$.  Отсюда следует, что $\x(1)_p\ge |K|_p=p^{a-d}$, т.\,е. $p^{a-d}\bigm| \x(1)$.

С другой стороны
$$
0\ne a_{\mbox{}_B}(K)=\left(\frac{1}{|G|} \sum_{\x\in \irr(B)}\x(1)\ov{\x(x_{\mbox{}_K})}\right)^\st=
\left(\frac{1}{|G|} \sum_{\vf\in \iBr(B)}\th_\vf(1)\ov{\vf(x_{\mbox{}_K})}\right)^\st, \eqno{(**)}
$$
где последнее равенство следует из \ref{os}, поскольку $x_{\mbox{}_K}\in G_{p'}$ в силу \ref{dk cor}.
Так как $|G|_p\bigm | \th_\vf(1)$ ввиду \ref{gl ner}$(vii)$, имеем
$\displaystyle\frac{\th_\vf(1)}{|G|}\in \wt{Z}$ для всех
$\vf$ и, значит, из $(**)$ следует, что найдётся характер $\vf\in \iBr(B)$ для которого
$\ov{\vf(x_{\mbox{}_K})}\not\in \wtm$.

Вследствие \ref{z lin bl} характер $\vf$ является $\ZZ$-линейной комбинацией брауэровых характеров
$\{\hat\t\,|\,\t\in \irr(B)\}$. Поэтому $\ov{\hat\t(x_{\mbox{}_K})}=\ov{\t(x_{\mbox{}_K})}\not\in \wtm$
для некоторого $\t\in \irr(B)$. Однако число
$\a=\displaystyle\frac{|K|\ov{\t(x_{\mbox{}_K})}}{\t(1)}$ лежит в $\ov{\ZZ}$
и $\ov{ \t(x_{\mbox{}_K})}=\a\displaystyle\frac{\t(1)}{|K|}$, откуда следует, что
$\displaystyle\frac{\t(1)}{|K|}\not\in\wtm$. Поэтому $\t(1)_p\le |K|_p=p^{a-d}$.

$(ii)$ В действительности имеет место более сильный факт: множества натуральных чисел
$$\big\{\x(1)\bigm|\,\x\in \irr(B)\big\},\qquad \big\{\vf(1)\bigm|\, \vf\in \iBr(B)\big\}$$
имеют одинаковый наибольший общий делитель. Это следует из того, что $\x(1)=\hat\x(1)$ --- целочисленная
линейная комбинация значений $\vf(1)$ с коэффициентами, являющимися числами разложения $p$-блока $B$,
а в свою очередь $\vf(1)$ --- целочисленная  линейная комбинация значений $\x(1)$ в силу \ref{z lin bl}.
\end{proof}

\textext{В связи с предложением \ref{def deg} отметим, что в общем случае степени неприводимых брауэровых
характеров не делят порядок группы $G$. Более того,
 если $\vf\in \iBr(G)$, то $p$-часть $\vf(1)_p$, вообще говоря, может
превосходить $|G|_p$. Так, например, простая спорадическая группа Маклафлина $McL$
порядка $2^7\cdot 3^6\cdot 5^3 \cdot 7\cdot 11$ обладает неприводимым $2$-модулярным
характером степени $2^9\cdot 7$.}


В условиях предложения \ref{def deg} для произвольного характера $\eta\in B$ можно записать
$$
\eta(1)_p=p^{a-d+h},
$$
где $h\ge 0$ --- однозначно определённое целое число, которое называется
\mem{$p$-высотой} (или просто \glsadd{iHtIrrOrdChr}\mem{высотой})   характера $\eta$.

Напомним, что ранее (см. стр. \pageref{pdef opr}) мы определили $p$-дефект произвольного
неприводимого обыкновенного характера. Очевидно, что
если $B\in\bl(G)$, то для любого характера  из $\irr(B)$ сумма  его $p$-дефекта и $p$-высоты совпадает с $\df(B)$.

Из предложения \ref{def deg} следует, что каждый $p$-блок содержит по крайней мере один обыкновенный
и один брауэров характер высоты $0$. В связи с этим отметим следующее утверждение.

\begin{thm}[Гипотеза Брауэра
\glsadd{iCnjBrZerHt}
о нулевой высоте\footnote{Brauer's height-zero conjecture.}] \label{gbr nv}
Если
$B\in\bl_p(G)$, то все характеры из $\irr(B)$ имеют $p$-высоту $0$ тогда и только тогда,
когда $\d(B)$ --- абелева группа.
\end{thm}
Доказательство гипотезы было завершено в работе \cite{mnst}.

\textuprn{Найти $2$-высоту всех неприводимых обыкновенных и  $2$-модулярных
характеров группы $A_5$. Проиллюстрировать для этого случая справедливость теоремы \ref{gbr nv}.}

\section{Блоки дефекта \texorpdfstring{$0$}{0}}

Наиболее простое описание имеют $p$-блоки нулевого дефекта.\footnote{Ср. предложение \ref{pdef 0}.}

\begin{pre} \label{nul def} Пусть $B\in\bl(G)$. Тогда следующие условия эквивалентны.

\begin{list}{{\rm(}{\it\roman{enumi}\/}{\rm)}}
{\usecounter{enumi}\setlength{\parsep}{2pt}\setlength{\topsep}{5pt}\setlength{\labelwidth}{23pt}}
\item $|\irr(B)|=|\iBr(B)|$.
\item $\x(g)=0$ для любого $\x\in\irr(B)$ и любого $g\in G\setminus G_{p'}$.
\item $\x(g)=0$ для любого $\x\in\irr(B)$ и любого неединичного $p$-элемента $g$.
\item $\df(B)=0$.
\item Существует $\x\in \irr(B)$ такой, что $|\x(1)|_p=|G|_p$.
\item $|\irr(B)|=1$.
\end{list}
\end{pre}
\begin{proof} $(i)\Rightarrow(ii)$ Из $(i)$ и \ref{bl dec} следует, что матрица разложения $\DD_B$ блока $B$
обратима. Положим
$$\DD_B^{-1}=(a_{\vf\x})_{\vf\in \iBr(B),\,\x\in \irr(B)}.
$$
Тогда для любого $\x\in \irr(B)$ получаем
$$
\sum_{\vf\in\iBr(B)}a_{\vf\x}\th_\vf=\sum_{\vf\in\iBr(B)}\sum_{\t\in\irr(B)}a_{\vf\x}d_{\t\vf}\t=
\sum_{\t\in\irr(B)}\d_{\x,\t}\t=\x.
$$
Поэтому, в силу \ref{gl ner}$(i)$, $\x$ обращается в ноль на элементах из $G\setminus G_{p'}$.

$(ii)\Rightarrow(iii)$ Очевидно.

$(iii)\Rightarrow(iv)$ Пусть $P\in \Syl_p(G)$. Из $(iii)$ следует, что для произвольного $\x\in \irr(B)$
выполнено равенство
$$
(\x_P,1_P)_{\mbox{}_P}=\frac{\x(1)}{|P|}.
$$
Так как левая часть --- целое число и $\x(1)$ делит $|G|$, отсюда следует, что $\x(1)_p=|G|_p$.
В силу произвольности~$\x$, из \ref{def deg}$(i)$ вытекает, что $\df(B)=0$.

$(iv)\Rightarrow(v)$ Из $(iv)$ и \ref{def deg}$(i)$ следует, что $\x(1)_p=|G|_p$ для любого
$\x\in \irr(B)$.

$(v)\Rightarrow(vi)$ Пусть $\x\in \irr(B)$ такой, что $\x(1)_p=|G|_p$. В частности, $\x(1)/|G|\in \wt{Z}$.
По \ref{id dec}, для центрального идемпотента $e_\x$ имеет место разложение
$$e_\x=\frac{\x(1)}{|G|}\sum_{g\in G}\x(g^{-1})g.$$
Поэтому $e_\x\in \wt{Z} G$ и из \ref{min un} получаем, что одноэлементное множество $\{\x\}$
является объединением $p$-блоков обыкновенных характеров. Значит, $\irr(B)=\{\x\}$.

$(vi)\Rightarrow(i)$ Ввиду \ref{bl ne} и \ref{bl dec}$(ii)$, имеем $1\le |\iBr(B)|\le |\irr(B)|=1$,
откуда следует $(i)$.
\end{proof}

Для натурального числа $n$ и произвольного множества $\pi$
простых чисел определим \glsadd{iPiPrtNum}\mem{$\pi$-часть} \gls{npi}
числа $n$ как наибольший делитель $n$, все простые делители которого принадлежат $\pi$
(ср. определение \ref{opch}). Напомним, что $\pi'$ обозначает множество простых чисел, не принадлежащих $\pi$.

Следующее упражнение частично обобщает предложение \ref{nul def}.

\textuprn{Пусть $\x\in \irr(G)$ и  $\pi$ --- множество простых чисел. Доказать равносильность
следующих утверждений.
\begin{list}{{\rm(}{\it\roman{enumi}\/}{\rm)}}
{\usecounter{enumi}\setlength{\parsep}{2pt}\setlength{\topsep}{5pt}\setlength{\labelwidth}{23pt}}
\item  $|G|_\pi$ делит $\x(1)$.
\item $\x(g)=0$ для любого неединичного $g\in G_\pi$.
\item $\x(g)=0$ для любого $g\in G\setminus G_{\pi'}$.
\end{list}}

\textuprn{Доказать справедливость теоремы \ref{gbr nv} для $p$-блоков дефекта $0$ и $1$.}

В связи с блоками дефекта $0$ упомянем ещё одну важную гипотезу. Пусть $P$ --- некоторая $p$-подгруппа
группы $G$. Назовём \mem{$p$-весом}\glsadd{iPWght}  пару $(P,\x)$, где $\x\in \irr(\N_G(P))$ --- такой характер,
что $P\le\ker\x$ и $\x(1)_p=|\N_G(P)/P|_p\,$, другими словами $\x$, рассматриваемый как характер  группы $\N_G(P)/P$,
принадлежит $p$-блоку нулевого дефекта. Если $(P,\x)$ --- $p$-вес, то для любого $g\in G$ пара $(P^g,\x^g)$,
как легко видеть, также является $p$-весом, т.\,е. $G$ действует сопряжением на множестве $p$-весов.

\normalmarginpar
\begin{gip}[Альперина
\glsadd{iConjAlpWt}
о весах\footnote{Alperin's weight conjecture.}] \label{galv}
Число классов сопряжённости $p$-весов группы $G$ совпадает с $|\iBr(G)|$.
\end{gip}
Более точную формулировку этой гипотезы мы приведём чуть ниже,\footnote{См. \ref{alb}.} когда сформулируем
понятие индуцированного блока.

\section{Гомоморфизм Брауэра. Индуцированные блоки}

Остановимся более подробно на связи $p$-блоков группы и её подгрупп.

Пусть $P$ --- произвольная $p$-подгруппа группы $G$. Положим $C=\C_G(P)$ и $N=\N_G(P)$.
Определим $F$-линейное отображение
$\b_{\mbox{}_P}:\Z(FG)\to \Z(FN)$, положив для произвольного $K\in \K(G)$
$$\b_{\mbox{}_P}(\wh K)=\sum_{\mathllap{c}\in K\cap \mathrlap{C}}\ c$$
и продолжив это отображение по линейности на $\Z(FG)$.  Поскольку $K\cap C$ является
объединением классов сопряжённости
группы $N$, образ $\b_{\mbox{}_P}$ действительно лежит в $\Z(FN)$.

\begin{pre} \label{br hom} Для любой $p$-подгруппы $P$ отображение $\b_{\mbox{}_P}:\Z(FG)\to \Z(FN)$,
где $N=\N_G(P)$, является гомоморфизмом $F$-алгебр.
\end{pre}
\begin{proof} Достаточно проверить, что
$$\b_{\mbox{}_P}(\wh K\wh L)=\b_{\mbox{}_P}(\wh K)\b_{\mbox{}_P}(\wh L) \myeqno\label{zv11}$$
для любых классов $K,L\in \K(G)$.  Зафиксируем $c\in C=\C_G(P)$ и обозначим
\begin{align*}
 \A&=\big\{\ (x,y)\in K\times L\bigm|\  xy=c\, \big\}, \\
  \A_0&=\big\{\ (x,y)\in (K\cap C)\times(L\cap C)\bigm|\  xy=c\, \big\}.
\end{align*}
Из правила умножения классовых сумм \ref{kl sum}$(iv)$ следует, что коэффициенты при $c$ в левой и правой частях
соотношения \ref{zv11} равны $|\A|^\st$ и $|\A_0|^\st$, соответственно.

Достаточно показать, что $|\A|\equiv|\A_0|\mod{p}$. Заметим, что $P$ действует на множестве $\A$
по правилу $(x,y)^u=(x^u,y^u)$ для любого $u\in P$. Легко видеть, что $\A_0$ --- это, в точности, множество
неподвижных точек группы $P$ относительно такого действия. Отсюда следует требуемое.
\end{proof}

\begin{opr} Отображение \gls{bsP} называется
\mem{гомоморфизмом Брауэра}.
\end{opr}

\begin{pre} \label{esg} Пусть
$P$ --- $p$-подгруппа группы $G$, $C=\C_G(P)$ и $K\in\K(G)$.
Следующие утверждения эквивалентны.
\begin{list}{{\rm(}{\it\roman{enumi}\/}{\rm)}}
{\usecounter{enumi}\setlength{\parsep}{2pt}\setlength{\topsep}{5pt}\setlength{\labelwidth}{23pt}}
\item $\b_{\mbox{}_P}(\wh K)\ne 0$.
\item $K\cap C\ne\varnothing$.
\item $P\le_G\d(K)$.
\end{list}
\end{pre}
\textupl{esg prf}{Доказать предложение \ref{esg}.}

Изучим ядро гомоморфизма Брауэра.

\begin{pre} \label{ker br} Пусть $P$ --- $p$-подгруппа группы $G$.
\begin{list}{{\rm(}{\it\roman{enumi}\/}{\rm)}}
{\usecounter{enumi}\setlength{\parsep}{2pt}\setlength{\topsep}{5pt}\setlength{\labelwidth}{23pt}}
\item Выполнено равенство
$$
\Ker \b_{\mbox{}_P} =\sum_{\substack{K\in \K(G),\\   P\ \nleqslant_G\ \d(K)}} F\wh K.
$$
\item Для произвольного $B\in \bl(G)$ имеем $\b_{\mbox{}_P}(e_{\mbox{}_B})=0$ тогда и только тогда,
когда $P\nleqslant_G \d(B)$.
\end{list}
\end{pre}
\begin{proof} Положим $C=\C_G(P)$.

$(i)$ Пусть элемент $z\in \Z(FG)$ лежит в $\Ker \b_{\mbox{}_P}$. Запишем
$$z=\sum_{K\in \K(G)} a_{\mbox{}_K}\wh K,$$
где $a_{\mbox{}_K}\in F$. Поскольку для различных $K\in \K(G)$ множества $K\cap C$ попарно не пересекаются,
равенство $\b_{\mbox{}_P}(z)=0$ возможно тогда и только тогда, когда $\b_{\mbox{}_P}(\wh K)=0$ для
всех $K$ таких, что $a_{\mbox{}_K}\ne 0$. Поэтому требуемое следует из \ref{esg}.

$(ii)$ Запишем
$$
e_{\mbox{}_B}=\sum_{K\in\K(G)}a_{\mbox{}_B}(K)\wh K.
$$
Из $(i)$ следует, что $\b_{\mbox{}_P}(e_{\mbox{}_B})=0$ тогда и только тогда, когда для любого класса
$K\in \K(G)$ с условием $a_{\mbox{}_B}(K)\ne 0$ выполнено  $P\nleqslant_G \d(K)$.

Но если это так и $P\le_G \d(B)$, то взяв в качестве $K$ дефектный класс блока $B$, получим
с одной стороны $a_{\mbox{}_B}(K)\ne 0$, а с другой $P\le_G \d(K)$, т.\,к. $\D(B)=\D(K)$ в силу \ref{cor bl kl gr}.
Противоречие.

Обратно, если $P\nleqslant_G \d(B)$ и $a_{\mbox{}_B}(K)\ne 0$, то $\d(K)\le_G\d(B)$ в силу \ref{thm min max}$(ii)$
и, значит, $P\nleqslant_G \d(K)$.\end{proof}

Таким образом, из \ref{ker br}$(ii)$ мы получаем ещё одну характеризацию дефектных групп $p$-блоков.

\begin{cor} \label{cor db ch} Пусть $B\in \bl(G)$.
Тогда $\D(B)$ состоит из максимальных по включению
$p$-подгрупп $P$ группы $G$, для которых $\b_{\mbox{}_P}(e_{\mbox{}_B})\ne 0$.
\end{cor}

Теперь определим операцию, в некотором смысле двойственную к гомоморфизму Брауэра. Эта двойственность
станет более ясной чуть ниже, когда мы сформулируем первую основную теорему Брауэра.

Пусть $H\le G$ и $b\in \bl(H)$. Используя центральный характер $\l_b:\Z(FH)\to F$, определим
$F$-линейное отображение $\l_b^{\, G}:\Z(FG)\to F$, положив
$$
\l_b^{\, G}(\wh K)=\l_b\left(\sum_{x\in K\cap H}x\right)
$$
для произвольного $K\in \K(G)$. Если $\l_b^{\, G}$ является
гомоморфизмом алгебр (что, вообще говоря, бывает не всегда), то в силу \ref{e b}$(ii)$ существует единственный
$p$-блок $B\in \bl(G)$ такой, что $\l_b^{\, G}=\l_{\mbox{}_B}$. В этом случае будем говорить, что определён
\glsadd{iPBlkInd}\mem{индуцированный $p$-блок} $B$ и обозначать его через
\gls{buG}.

\textuprn{Доказать, что для любого $B\in \bl(G)$ индуцированный блок $B^G$ определён и совпадает с $B$.}

Приведём пример, показывающий, что для блока $b$ подгруппы
$H$ группы $G$ индуцированный блок $b^G$ может быть не определён.
Отметим, что в этом примере блок $b$ является главным и одновременно имеет нулевой дефект.

\begin{prim}
Пусть $p=2$, $G=S_3$, $H=A_3\le G$ и $b$ ---
главный $2$-блок (он же блок дефекта $0$) группы $H$.
Значения гомоморфизма $\l_b$ на классовых суммах $\wh L$, где $L\in \K(H)$ приведены в следующей таблице.
$$
  \begin{array}{c|ccc}
       L                     & L_1   &  L_2 &  L_3    \\
 \hline
       x_{\mbox{}_L}^{\vphantom{A^{A^A}}}         & 1    &  (123) &  (132)  \\
     \om_{1_H}(\wh L)=|L|^{\vphantom{A^{A^A}}} & 1   &  1  &  1  \\

    \l_b(\wh L)=|L|^{\st^{\vphantom{A^A}}}   & 1    & 1 & 1
  \end{array}
$$

\noindent
Легко найти значения отображения $\l_b^G$ на классовых суммах $\wh K$, где $K\in \K(G)$.
$$
  \begin{array}{c|ccc}
     K                     & K_1   &  K_2 &  K_3\\
   \hline
       x_{\mbox{}_K}^{\vphantom{A^{A^A}}}       & 1   &  (12) &  (123)   \\
       \sum_{x\in K\cap H}x^{\vphantom{A^{A^A}}}         & \wh L_1    &  0 & \wh L_2+\wh L_3 \\
     \l_b^G(\wh K)^{\vphantom{A^{A^{A^A}}}} & 1   &  0  &  0
  \end{array}
$$

\noindent
Группа $G$ имеет два $2$-блока, один из которых главный (обозначим его $B_1$), а другой (обозначим его $B_2$)
имеет нулевой дефект (см. приложение \ref{pril tch}). Пусть $\irr(B_2)=\{\x\}$.
Определим значения гомоморфизмов $\l_{\mbox{}_{B_1}}$ и $\l_{\mbox{}_{B_2}}$.
$$
  \begin{array}{c|ccc}
     K                     & K_1   &  K_2 &  K_3  \\
   \hline
        1_G(x_{\mbox{}_K})^{\vphantom{A^{A^A}}} & 1 & 1 & 1 \\
   \om_{1_G}(\wh K)=|K|  & 1    &  3 & 2   \\
     \l_{\mbox{}_{B_1}}(\wh K)=|K|^{\st^{\vphantom{A^A}}}_{\vphantom{A_{A_A}}} & 1   &  1  &  0 \\
        \hline
        \x(x_{\mbox{}_K})^{\vphantom{A^{A^A}}} & 2 & 0 & -1 \\
   \om_{\x}(\wh K)^{\vphantom{A^{A^A}}}  & 1    &  0 & -1    \\
     \l_{\mbox{}_{B_2}}(\wh K)^{\vphantom{A^{A^A}}} & 1   &  0  &  1 \\

  \end{array}
$$

\noindent
Как видно,  $\l_b^G$ не совпадает ни $\l_{\mbox{}_{B_1}}$, ни с $\l_{\mbox{}_{B_2}}$. Таким образом, блок $b^G$
не определён.
\end{prim}

\begin{pre} \label{ind def} Пусть $b\in \bl(H)$, где $H\le G$, и предположим, что определён $p$-блок $b^G$.
Тогда $\d(b)\le_G\d(b^G)$.
\end{pre}
\begin{proof} Пусть $K$ --- дефектный класс $p$-блока $b^G$. Тогда
$$
0\ne \l_{b^G}(\wh K)=\l_b^{\, G}(\wh K)=\l_b\left(\sum_{x\in K\cap H}x\right).
$$
В частности, $K\cap H\ne \varnothing$ и существует класс $L\in \K(H)$ такой, что $L\se K$ и $\l_b(\wh L)\ne 0$.
Из теоремы min--max \ref{thm min max}$(i)$ следует, что $\d(b)\le_H\d(L)$. Поскольку $\d(L)\in \Syl_p(\C_H(x))$
для некоторого $x\in L\se K$, имеем $\d(L)\le P$ для некоторой подгруппы $P\in \Syl_p(\C_G(x))$.
Так как $P\in \D(K)=\D(b^G)$, получаем $\d(b)\le_H\d(L)\le_G \d(b^G)$, т.\,е. $\d(b)\le_G \d(b^G)$.
\end{proof}

Гомоморфизм Брауэра можно использовать для получения достаточного
условия существования индуцированных $p$-блоков.

\begin{pre} \label{comp ind} Пусть $P\le G$ --- $p$-подгруппа и $P\C_G(P)\le H\le \N_G(P)$.
Тогда справедливы следующие утверждения.
\begin{list}{{\rm(}{\it\roman{enumi}\/}{\rm)}}
{\usecounter{enumi}\setlength{\parsep}{2pt}\setlength{\topsep}{5pt}\setlength{\labelwidth}{23pt}}
\item  Образ гомоморфизма Брауэра $\b_{\mbox{}_P}$ лежит в $\Z(FH)$, и для любого $b\in \bl(H)$ диаграмма
$$
\xymatrix{
\Z(FG)\ar[dr]^{ \textstyle{\l_b^{\, G}}}\ar[d]_{\textstyle{\b_{\mbox{}_P}}}&  \\
\Z(FH)\ar[r]_{ \textstyle{\l_b}}&F
}
$$
коммутативна. В частности, $\l_b^{\, G}$ --- гомоморфизм $F$-алгебр.
\item Для любого $b\in \bl(H)$ определён индуцированный $p$-блок $b^G$.
\item Для любого $B\in \bl(G)$ имеет место равенство
$$
\b_{\mbox{}_P}(e_{\mbox{}_B})=\sum_{\substack{b\in \bl(H),\\B\,=\,b^G}}e_b.
$$
\item Пусть $B\in \bl(G)$. Тогда $B=b^G$ для некоторого $b\in \bl(H)$ в том и только в том случае,
когда $P\le_G \d(B)$.
\end{list}
\end{pre}
\begin{proof} $(i)$ Пусть $C=\C_G(P)$. Так как $H$ нормализует $C$, образ гомоморфизма
Брауэра $\b_{\mbox{}_P}$ действительно лежит
в $\Z(FH)$. Проверим, что $\l_b^{\vphantom{G}}\circ\b_{\mbox{}_P}=\l_b^{\, G}$. Для этого
достаточно показать, что
$$
\l_b\left(\,\sum_{x\in K\cap H}x\right)=\l_b\left(\,\sum_{x\in K\cap C}x\right)
$$
для всякого $K\in \K(G)$, т.\,е. что $\l_b$ принимает нулевое значение на сумме элементов из множества
$K\cap (H\setminus C)$, которое является объединением некоторых классов $L\in \K(H)$, не пересекающихся с $C$.
Поскольку $P\nor H$, имеем $\C_H(\OO_p(H))\le C$, т.\,е. для каждого такого класса $L$  сумма $\wh L$ ---
нильпотентный элемент из $FH$ в силу \ref{k nil}. Поэтому $\l_b(\wh L)=0$, откуда следует требуемое.

$(ii)$ Это следует из $(i)$, поскольку $\l_b^G$ --- гомоморфизм $F$-алгебр.

$(iii)$ Пусть $B\in \bl(G)$. Рассмотрим значение $\b_{\mbox{}_P}(e_{\mbox{}_B})$. Поскольку
$\b_{\mbox{}_P}$ --- гомоморфизм $F$-алгебр, элемент $\b_{\mbox{}_P}(e_{\mbox{}_B})\in\Z(FH)$
является идемпотентом. В силу \ref{e b}$(viii)$ можно записать
$$
\b_{\mbox{}_P}(e_{\mbox{}_B})=e_{b_1}+\ld +e_{b_k}\myeqno\label{zv10}
$$
для некоторых различных $b_1,\ld,b_k\in \bl(H)$. Заметим, что ввиду $(i)$ равенство $B=b^G$ для некоторого
$b\in \bl(H)$ имеет место тогда и только тогда, когда гомоморфизм $F$-алгебр
$\l_b^{\vphantom{G}}\circ\b_{\mbox{}_P}:\Z(FG)\to F$ совпадает с $\l_{\mbox{}_B}$, т.\,е. когда значение
$\l_b^{\vphantom{G}}\circ\b_{\mbox{}_P}$ на $e_{\mbox{}_B}$ равно $1$. В силу \ref{zv10} это выполнено тогда и
только тогда, когда $b$ совпадает в одним из $b_i$. Отсюда получаем требуемое.

$(iv)$ Пусть $B\in \bl(G)$.

Если $B=b^G$ для некоторого $b\in \bl(H)$, то из \ref{ind def} следует, что
$\d(b)\le_G\d(B)$, а из \ref{op def} получаем, что $P\le \OO_p(H)\le \d(B)$. Поэтому $P\le_G\d(B)$.

Обратно, предположим, что $P\le_G\d(B)$. Тогда $\b_{\mbox{}_P}(e_{\mbox{}_B})\ne 0$ ввиду \ref{ker br}$(ii)$.
Из $(iii)$ вытекает, что существует $b\in \bl(H)$ для которого $B=b^G$.
\end{proof}

\begin{cor} \label{lndef} Пусть блок $B\in \bl(G)$ такой, что  $\d(B)$ --- нормальная подгруппа группы $G$.
Тогда
$$
e_{\mbox{}_B}=\sum_{\substack{K\in\K(G_{p'}),\\\D(K)=\D(B)}}a_{\mbox{}_B}(K)\wh K.
$$
\end{cor}
\begin{proof} Напомним, что имеет место представление
$$
e_{\mbox{}_B}=\sum_{K\in\K(G_{p'})}a_{\mbox{}_B}(K)\wh K, \myeqno\label{zv8}
$$
см. \ref{zzv7}.
В силу нормальности $\d(B)$ множество $\D(B)$ одноэлементно. Поэтому достаточно доказать, что
если $a_{\mbox{}_B}(K)\ne 0$, то $\d(K)=\d(B)$.
Из теоремы min--max \ref{thm min max} следует, что если $a_{\mbox{}_B}(K)\ne 0$, то $\d(K)\le \d(B)$.

Для доказательства обратного включения применим предложение \ref{comp ind}, положив $P=\d(B)$ и $H=G$.
Тогда, ввиду того, что $B=B^G$, из \ref{comp ind}$(iii)$ следует равенство
$\b_{\mbox{}_P}(e_{\mbox{}_B})=e_{\mbox{}_B}$. Применив гомоморфизм Брауэра
$\b_{\mbox{}_P}$ к обеим частям соотношения \ref{zv8}, получим
$$
\b_{\mbox{}_P}(e_{\mbox{}_B})=\sum_{K\in\K(G_{p'})}a_{\mbox{}_B}(K)\b_{\mbox{}_P}(\wh K).
$$
Поскольку для любого $K$ значение $\b_{\mbox{}_P}(\wh K)$ является суммой некоторых различных элементов из
$K$, в силу линейной независимости получаем, что равенство $\b_{\mbox{}_P}(e_{\mbox{}_B})=e_{\mbox{}_B}$
возможно лишь тогда, когда $\b_{\mbox{}_P}(\wh K)=\wh K$ для всех $K$ таких, что $a_{\mbox{}_B}(K)\ne 0$.
В частности, если $a_{\mbox{}_B}(K)\ne 0$, то $K\cap C\ne \varnothing$,  где $C=\C_G(\d(B))$
и, значит, $\d(B)\le\d(K)$, что и требовалось доказать.
\end{proof}

Для доказательства ещё одного достаточного условия существования индуцированного блока
нам потребуется следующее утверждение.

\begin{pre} \label{trn} Пусть $H\le N\le G$ и $b\in \bl(H)$. Пусть определён индуцированный блок $b^N\in \bl(N)$.
Тогда блок $b^G$ определён если и только если определён блок $(b^N)^G$, и в этом случае $b^G=(b^N)^G$.
\end{pre}
\textupl{trn prf}{Доказать предложение \ref{trn}.}

\begin{pre} \label{dgi} Пусть $H\le G$ и $b\in \bl(H)$.
Если $\C_G(\d(b))\le H$, то определён индуцированный блок $b^G$.
\end{pre}
\begin{proof} Пусть $P=\d(b)$ и $T=P\C_G(P)$. По условию имеем $T\le H$. Из \ref{comp ind}$(iv)$
вытекает, что существует блок $b_0\in \bl(T)$ такой, что $b_0^H=b$. С другой стороны,
в силу \ref{comp ind}$(ii)$ индуцированный блок $b_0^G$ определён. По \ref{trn} получаем,
что блок $b^G=(b_0^H)^G$ определён и совпадает с $b_0^G$.
\end{proof}

Как легко видеть, условие существования индуцированного блока $b^G$, приведённое в \ref{dgi},
не зависит от выбора дефектной группы $\d(b)$. Также отметим, что
\ref{comp ind}$(ii)$ является следствием  \ref{dgi}, так как если $P\C_G(P)\le H\le \N_G(P)$,
то для любого блока $b\in \bl(H)$ выполнено включение $\C_G(\d(b))\le H$, поскольку $P$ содержится в $\d(b)$
ввиду \ref{op def}.

Пусть $P$ --- $p$-подгруппа группы $G$. Введём следующие обозначения
\begin{align*}
\gls{KlGPr}&=\big\{K\in\K(G)\bigm|\,P\in \D(K)\big\},\\
\gls{BllGPr}&=\big\{B\in\bl(G)\bigm|\,P\in \D(B)\big\}.
\end{align*}

\begin{pre} \label{p np} Пусть $P$ --- $p$-подгруппа группы $G$. Положим $C=\C_G(P)$, $N=\N_G(P)$.
Тогда отображение $K\mapsto K\cap C$ осуществляет биекцию множеств $\K(G|P)\to\K(N|P)$.
\end{pre}
\begin{proof} Пусть $K\in \K(G|P)$ и $P\in \Syl_p(\C_G(x))$ для подходящего $x\in K$. Тогда $x\in K\cap C$.
Пусть $y\in K\cap C$. Тогда $P\in \Syl_p(\C_G(y))$, и поскольку $y=x^g$ для некоторого $g\in G$,
также имеем $P^g\in \Syl_p(\C_G(y))$. Значит, $P=P^{gc}$ для некоторого $c\in \C_G(y)$. Тогда $x^{gc}=y^c=y$,
т.\,е. $x$ и $y$ сопряжены элементом из $N$. Таким образом, $K\cap C\in \K(N)$. Очевидно также,
что $P\in \D(K\cap C)$ и что отображение $K\mapsto K\cap C$ инъективно.

Докажем сюръективность.
Пусть $L\in \K(N|P)$ и $K\in \K(G)$ --- класс, содержащий $L$. Проверим, что $P\in \D(K)$.
Пусть $x\in L$ --- элемент, для которого $P\in \Syl_p(\C_N(x))$ и пусть $P\le S\in \Syl_p(\C_G(x))$.
Если $P<S$, то из условия $N=\N_G(P)$ и того, что $S$ --- $p$-группа следует, что $P<\N_S(P)=S\cap N\le \C_N(x)$.
Это противоречит тому, что $P\in \D(L)$.  Значит, $P=S$ и $P\in \D(K)$.
\end{proof}

\begin{thm}[Первая основная теорема Брауэра] \label{thm br 1}\glsadd{iThmBraFstMain} Пусть $P$ --- $p$-подгруппа группы $G$ и $N=\N_G(P)$.
\begin{list}{{\rm(}{\it\roman{enumi}\/}{\rm)}}
{\usecounter{enumi}\setlength{\parsep}{2pt}\setlength{\topsep}{5pt}\setlength{\labelwidth}{23pt}}
\item Отображение $b\mapsto b^G$ является биекцией множеств $\bl(N|P)\to \bl(G|P)$.
\item Для любого $b\in \bl(N|P)$ имеем $\b_{\mbox{}_P}(e_{\mbox{}_B})=e_b^{\vphantom{G}}$, где $B=b^G$.
\end{list}
\end{thm}
\begin{proof} Обозначим $C=\C_G(P)$.

Сначала покажем, что образ отображения $b\mapsto b^G$, где $b\in \bl(N|P)$, действительно лежит в $\bl(G|P)$.
Пусть $b\in \bl(N|P)$. В силу \ref{comp ind} определён индуцированный $p$-блок $B=b^G$ и
$\l_b^{\vphantom{G}}\circ\b_{\mbox{}_P}=\l_{\mbox{}_B}$.
Требуется проверить, что $B\in \bl(G|P)$. Пусть $L$ --- дефектный класс $p$-блока $b$. Тогда $L\in \K(N|P)$
ввиду \ref{cor bl kl gr}. Пусть $K\in \K(G)$ --- класс, содержащий $L$. Из \ref{p np} следует, что $K\in \K(G|P)$
и $L=K\cap C$. Поэтому
$$
\l_{\mbox{}_B}(\wh K)=\l_b^{\vphantom{G}}(\b_{\mbox{}_P}(\wh K))=\l_b^{\vphantom{G}}(\wh L)\ne 0.
$$
По теореме min--max \ref{thm min max}$(i)$ получаем $\d(B)\le_G P$. С другой стороны, из \ref{ind def}
следует, что $P\le_G\d(B)$. Поэтому $P\in \D(B)$.

Теперь проверим, что отображение $b\mapsto b^G$ сюръективно. Пусть $B\in \bl(G|P)$. Из \ref{comp ind}$(iii)$
следует, что существуют $p$-блоки $b_1,\ld,b_k\in \bl(N)$ такие, что $b_i^G=B$ и
$$
\b_{\mbox{}_P}(e_{\mbox{}_B})=e_{b_1}+\ld+e_{b_k}.
$$
Так как $P\nor N$, из \ref{op def} вытекает, что $P\le \d(b_i)$ для всех $i=1,\ld,k$. Из \ref{ind def} также
получаем $\d(b_i)\le_G\d(B)=P$ для всех $i=1,\ld,k$. Отсюда следует, что $P\in \D(b_i)$, и тем самым
сюръективность отображения $b\mapsto b^G$ доказана.  Более того, если докажем инъективность, то получим, что
$k=1$, т.\,е. выполнено $(ii)$.

Допустим, что $b_1^G=B=b_2^G$ для некоторых $b_1,b_2\in \bl(N|P)$. Тогда
$\l_{b_1}^{\vphantom{G}}\circ \b_{\mbox{}_P}=\l_{\mbox{}_B}=\l_{b_2}^{\vphantom{G}}\circ\b_{\mbox{}_P}$
в силу \ref{comp ind}$(i)$. Таким образом, для любого $K\in \K(G)$ имеем $\l_{b_1}(\wh L)=\l_{b_2}(\wh L)$,
где $L=K\cap C$. В частности, из \ref{p np} вытекает, что значения $\l_{b_1}$ и $\l_{b_2}$ совпадают
на классовых суммах для всех классов из $\K(N|P)$.  Из соотношения \ref{zv6}
получаем
$$
\ba{rl}
1=&\l_{b_1}(e_{b_1})=\displaystyle\sum_{L\in\K(N)}a_{b_1}(L)\l_{b_1}(\wh L)
=\sum_{L\in\K(N|P)}a_{b_1}(L)\l_{b_1}(\wh L), \\ [20pt]
&\l_{b_2}(e_{b_1})=\displaystyle\sum_{L\in\K(N)}a_{b_1}(L)\l_{b_2}(\wh L)
=\sum_{L\in\K(N|P)}a_{b_1}(L)\l_{b_2}(\wh L),
\ea
\myeqno\label{zzv5}
$$
где последние равенства в каждой цепочке выполнены, поскольку ввиду теоремы min--max \ref{thm min max} произведение
$a_{b_1}(L)\l_{b_i}(\wh L)$, $i=1,2$, ненулевое только если $\d(L)\le_G P$ и $P\le_G \d(L)$,
т.\,е. если $L\in \K(N|P)$. Из сделанных выше замечаний вытекает, что правые части в \ref{zzv5}
совпадают. Поэтому $\l_{b_2}(e_{b_1})=1$ и $b_1=b_2$ в силу \ref{e b}$(ii)$.
\end{proof}

Говорят, что подгруппа  $P$ группы $G$ является
\glsadd{iRadPSubGr}\mem{радикальной $p$-подгруппой}, если $P=\OO_p(\N_G(P))$.

В качестве одного из следствий теоремы \ref{thm br 1} мы получаем следующее усиление утверждения \ref{op def},
дающее необходимое условие для того, чтобы $p$-подгруппа была дефектной группой некоторого $p$-блока.

\begin{cor} \label{cor p rad} Для любого $B\in \bl(G)$
дефектная группа $\d(B)$ является радикальной $p$-подгруппой.
\end{cor}

\textupl{cor p rad prf}{Доказать следствие \ref{cor p rad}.}

Пусть $B\in \bl(G)$.  Говорят, что \mem{$p$-вес} $(P,\x)$ группы $G$
\glsadd{iPWghtBlgBlk}\mem{принадлежит блоку} $B$, если
$B=b^G$, где $b\in \bl(\N_G(P))$ и $\x\in \irr(b)$.

\textuprn{Доказать следующие утверждения.
\begin{list}{{\rm(}{\it\roman{enumi}\/}{\rm)}}
{\usecounter{enumi}\setlength{\parsep}{2pt}\setlength{\topsep}{5pt}\setlength{\labelwidth}{23pt}}
\item Если  $(P,\x)$ --- $p$-вес, принадлежащий блоку $B$, то $P\in\D(B)$. В частности, $P$ --- радикальная $p$-подгруппа.
\item Сопряжённые $p$-веса принадлежат одному блоку.
\end{list}}

Гипотеза Альперина о весах \ref{galv} в уточнённой формулировке звучит следующим образом.

\begin{gip}[Альперина о весах, блочная версия] \label{alb}
Число классов сопряжённости $p$-весов группы $G$, принадлежащих
блоку $B\in \bl(G)$, совпадает с $|\iBr(B)|$.
\end{gip}

Отметим, что имеется ряд усилений гипотезы Альперина, принадлежащих различным авторам, среди которых интенсивно изучаемые в настоящее время 
\glsadd{iCnjDade}\mem{гипотеза Дэйда}
и 
\glsadd{iCnjCharTr}\mem{гипотеза о тройках характеров}
\footnote{Character Triple conjecture.} см.~\cite{rossi,spath}. Их доказательство сведено к случаю 
\glsadd{iGrpQsSimp}\mem{квазипростых}
 групп, т.\,е. таких групп $G$, что $G=G'$ и факторгруппа $G/\Z(G)$ проста.

\section{Теоремы Робинсона и Грина}

В этом разделе мы изложим важные результаты Робинсона о числе блоков с данной дефектной группой,
которые затем применим для доказательства теоремы Грина о свойстве дефектных групп.

Определим \mem{отображение Робинсона} \glsadd{iMapRob} $\gls{Rb}:\Z(FG)\to \Z(FG)$
по правилу
$$
\vr(x)=\sum_{B\in \bl(G)}\l_{\mbox{}_B}(x)e_{\mbox{}_B}
$$

\begin{pre} \label{rbp} Имеют место следующие утверждения.
\begin{list}{{\rm(}{\it\roman{enumi}\/}{\rm)}}
{\usecounter{enumi}\setlength{\parsep}{2pt}\setlength{\topsep}{5pt}\setlength{\labelwidth}{23pt}}
\item $\vr$ является эндоморфизмом $F$-алгебры $\Z(FG)$.
\item $\vr(e_{\mbox{}_B})=e_{\mbox{}_B}$ для любого $B\in \bl(G)$.
\item $\vr^2=\vr$.
\item $\Ker \vr=\J(\Z(FG))$.
\end{list}
\end{pre}

\textupl{rbp prf}{Доказать предложение \ref{rbp}.}

Рассмотрим действие $\vr$ на классовых суммах. Пусть $L\in \K(G)$. Тогда
$$
\vr(\wh L)=\sum_{B\in \bl(G)}\l_{\mbox{}_B}(\wh L)e_{\mbox{}_B}=
\sum_{B\in \bl(G)}\l_{\mbox{}_B}(\wh L)\left(\sum_{K\in\K(G)}a_{\mbox{}_B}(K)\wh K\right)=
\sum_{K\in\K(G)}\left(\sum_{B\in \bl(G)}\l_{\mbox{}_B}(\wh L)a_{\mbox{}_B}(K)\right)\wh K.
$$
Таким образом, матрица отображения $\vr$ в базисе из классовых сумм равна
$$
\left(\sum_{B\in \bl(G)}\l_{\mbox{}_B}(\wh L)a_{\mbox{}_B}(K)\right)_{K,L\in \K(G)}.
$$
Наша ближайшая цель --- найти удобную формулу для вычисления сумм вида
$$
\sum_{B\in \bl(G)}\l_{\mbox{}_B}(\wh L)a_{\mbox{}_B}(K),
$$
где $K,L\in \K(G)$.

Пусть $P\in \Syl_p(G)$. Для $K,L\in \K(G)$ положим
$$
\gls{Akl}=\{(x,y)\in K\times L\mid Px=Py)\}. \myeqno\label{zv3}
$$

\begin{pre} \label{omkl} Пусть $K,L\in \K(G)$.
\begin{list}{{\rm(}{\it\roman{enumi}\/}{\rm)}}
{\usecounter{enumi}\setlength{\parsep}{2pt}\setlength{\topsep}{5pt}\setlength{\labelwidth}{23pt}}
\item $|\A_{K,L}|=|\A_{L,K}|$ и эта величина не зависит от выбора $P\in \Syl_p(G)$ в определении {\rm \ref{zv3}}.
\item $\displaystyle\frac{|\A_{K,L}|}{|K|}\in \wt{Z}$.
\item Если $K\in \K(G_{p'})$, то
$$
\sum_{B\in\bl(G)}\l_{\mbox{}_B}(\wh L)a_{\mbox{}_B}(K)=\left(\frac{|\A_{K,L}|}{|K|}\right)^\st
$$
\end{list}
\end{pre}
\begin{proof} Утверждение $(i)$ напрямую следует из определения и сопряжённости силовских
подгрупп.

%

В силу \ref{str chr} для фиксированного $z\in G$ имеем
$$
\big|\{(x,y)\in K\times L\mid xy=z\}\big|=\frac{|K|\,|L|}{|G|}
\sum_{\x\in
\irr(G)}\frac{\x(x_{\mbox{}_K})\x(x_{\mbox{}_L})\ov{\x(z)}}{\x(1)}.
$$
Обозначим $K^{-1}=(x_{\mbox{}_K}^{-1})^G=\{y^{-1}\mid y\in K\}$. Тогда
$|K|=|K^{-1}|$ и
\begin{align*}
|\A_{K,L}|=\sum_{z\in P}\big|\{(x,y)\in K\times L\mid yx^{-1}=z\}\big|=&
\sum_{z\in P}\big|\{(x,y)\in K^{-1}\times L\mid yx=z\}\big|\\
=\sum_{z\in P}\frac{|K^{-1}|\,|L|}{|G|}
\sum_{\x\in\irr(G)}\frac{\ov{\x(x_{\mbox{}_K})}\x(x_{\mbox{}_L})\ov{\x(z)}}{\x(1)}
=&\frac{|K|\,|L|}{|G|}
\sum_{\x\in\irr(G)}\left(\ \sum_{z\in P}\ov{\x(z)}\right)\frac{\ov{\x(x_{\mbox{}_K})}\x(x_{\mbox{}_L})}{\x(1)}\\
=&\frac{|K|}{|G|}
\sum_{\x\in\irr(G)}|P|(1_P,\x_P)_{\mbox{}_P}\ov{\x(x_{\mbox{}_K})}\ \frac{|L|\x(x_{\mbox{}_L})}{\x(1)}.
\end{align*}
Поэтому
$$
\frac{|\A_{K,L}|}{|K|}=
\frac{|P|}{|G|}
\sum_{\x\in\irr(G)}(\x_P,1_P)_{\mbox{}_P}\ov{\x(x_{\mbox{}_K})}\om_\x(\wh{L}).\myeqno\label{zv4}
$$
Поскольку $\om_\x(\wh{L}),\x(x_{\mbox{}_K})\in \ov{\ZZ}$, $(\x_P,1_P)_{\mbox{}_P}\in \ZZ$, $|P|/|G|\in \wt{Z}$,
отсюда следует  $(ii)$.

Применив $^\st$ к обеим частям равенства \ref{zv4}, получим
\begin{align*}
\left(\frac{|\A_{K,L}|}{|K|}\right)^\st=
\left(\frac{|P|}{|G|}\right)^\st
\sum_{B\in\bl(G)}\ \left(\sum_{\x\in\irr(B)}
(\x_P,1_P)_{\mbox{}_P}\ov{\x(x_{\mbox{}_K})}\right)^\st
\l_{\mbox{}_B}(\wh{L})\\
=\sum_{B\in\bl(G)}
\left( \frac{1}{|G|}
\sum_{\x\in\irr(B)}|P|(\x_P,1_P)_{\mbox{}_P}\ov{\x(x_{\mbox{}_K})}\right)^\st\l_{\mbox{}_B}(\wh{L}).
\end{align*}
Предположим теперь, что $K\in \K(G_{p'})$. Из \ref{cor sl ort} следует, что
$$
\sum_{\x\in\irr(B)}|P|(\x_P,1_P)_{\mbox{}_P}\ov{\x(x_{\mbox{}_K})}=
\sum_{z\in P}\sum_{\x\in\irr(B)}\x(z)\ov{\x(x_{\mbox{}_K})}=
\sum_{\x\in\irr(B)}\x(1)\ov{\x(x_{\mbox{}_K})}.
$$
Поэтому, в силу \ref{abk} получаем, что
$$
\left(\frac{|\A_{K,L}|}{|K|}\right)^\st=\sum_{B\in\bl(G)}
\left( \frac{1}{|G|}\sum_{\x\in\irr(B)}\x(1)\ov{\x(x_{\mbox{}_K})}\right)^\st\l_{\mbox{}_B}(\wh{L})=
\sum_{B\in\bl(G)}a_{\mbox{}_B}(K)\l_{\mbox{}_B}(\wh L).
$$
\end{proof}

Из \ref{omkl}$(ii)$ вытекает следующее утверждение.
\begin{cor} \label{cmd} Пусть
$|G|_p=p^a$, $K\in \K(G)$, $|\d(K)|=p^d$. Тогда $\displaystyle\frac{|\A_{K,L}|}{p^{a-d}}\in \ZZ$
для любого $L\in \K(G)$.
\end{cor}

Для произвольной $p$-подгруппы $P$ группы $G$ по аналогии с введённым ранее обозначением $\K(G|P)$
положим
$$\gls{KlGpPr}=\big\{K\in\K(G_{p'})\bigm|\,P\in \D(K)\big\}.$$

Пусть $D$ --- нормальная $p$-подгруппа группы $G$.
Запишем $|G|_p=p^a$, $|D|=p^d$. Определим
$$
\gls{AlDr}=\left(\frac{|\A_{K,L}|}{p^{a-d}}\right)_{K,L\in \K(G_{p'}|D)}.
$$
Ввиду \ref{cmd} и \ref{omkl}$(i)$ матрица $A(D)$
является целочисленной и симметрической. Отметим, что в случае, когда $D$ не является дефектной группой
никакого $p$-регулярного класса, матрица $A(D)$ является пустой.

Теперь мы можем приступить к определению числа блоков группы $G$ с данной дефектной группой $D$.
Отметим, что в силу первой основной теоремы Брауэра можно считать, что подгруппа $D$ нормальна в $G$.

\begin{thm}[Робинсона] \label{trob}\glsadd{iThmRobn} Пусть $D$ --- нормальная $p$-подгруппа группы $G$.
Тогда число $p$-блоков с дефектной группой $D$ равно
рангу матрицы $A(D)^\st$.
\end{thm}
\begin{proof} Пусть $U$ обозначает $F$-линейную оболочку классовых сумм $\big\{\wh L\bigm| L\in \K(G_{p'}|D)\big\}$.
Покажем, что $\vr(U)\se U$. Если $L\in \K(G_{p'}|D)$ и $\l_{\mbox{}_B}(\wh L)\ne 0$ для
некоторого $B\in \bl(G)$, то из теоремы min--max \ref{thm min max} вытекает, что $\d(B)\le D$.
Поскольку $D\le\OO_p(G)$ и $\OO_p(G)$ содержится в дефектной группе любого блока группы $G$, отсюда
следует, что $D=\d(B)$. Значит, если $B\not\in \bl(G|D)$, то $\l_{\mbox{}_B}(\wh L)=0$.
Поэтому для $L\in \K(G_{p'}|D)$ имеем
$$
\vr(\wh L)=\sum_{B\in \bl(G)}\l_{\mbox{}_B}(\wh L)e_{\mbox{}_B}=
\sum_{B\in \bl(G|D)}\l_{\mbox{}_B}(\wh L)e_{\mbox{}_B}. \myeqno\label{zv2}
$$
В силу \ref{lndef} для любого $B\in \bl(G|D)$ имеет место равенство
$$
e_{\mbox{}_B}=\sum_{K\in\K(G_{p'}|D)}a_{\mbox{}_B}(K)\wh K. \myeqno\label{ebo}
$$
Поэтому
$$
\vr(\wh L)=\ \sum_{K\in\K(G_{p'}|D)}\left(\sum_{B\in \bl(G|D)}\l_{\mbox{}_B}(\wh L)a_{\mbox{}_B}(K)\right)\wh K
$$
и, значит, $\vr(U)\se U$.

Заметим, что ранг $F$-линейного отображения $\vr_{\mbox{}_U}:U\to U$ совпадает с $|\bl(G|D)|$. Это следует
из того, что идемпотенты $\{\,e_{\mbox{}_B}\mid B\in \bl(G|D)\,\}$ образуют базис образа $\vr(U)$:
они линейно независимы, каждый из них лежит в образе $\vr(U)$ по \ref{rbp}$(ii)$ и \ref{ebo},
и образ $\vr(U)$ содержится в их линейной оболочке в силу \ref{zv2}.

С другой стороны ранг отображения $\vr_{\mbox{}_U}$ совпадает с рангом его матрицы, которая
в базисе из классовых сумм $\{\,\wh K\mid K\in\K(G_{p'}|D)\,\}$ равна
$$
\left(\sum_{B\in \bl(G|D)}\l_{\mbox{}_B}(\wh L)a_{\mbox{}_B}(K)\right)_{L,K\in\K(G_{p'}|D)}.
$$
Как было отмечено выше, $\l_{\mbox{}_B}(\wh L)=0$ при $B\not\in \bl(G|D)$, поэтому суммирование
можно распространить на все блоки, т.\,е. ранг $\vr_{\mbox{}_U}$ равен рангу матрицы

$$
\left(\sum_{B\in \bl(G)}\l_{\mbox{}_B}(\wh L)a_{\mbox{}_B}(K)\right)_{L,K\in\K(G_{p'}|D)}=
\left(\frac{|\A_{K,L}|}{|K|}\right)^\st_{L,K\in\K(G_{p'}|D)}
$$
по \ref{omkl}$(iii)$, которая, как легко видеть, отличается от матрицы $A(D)^\st$ транспонированием и
умножением на невырожденную диагональную матрицу с диагональными элементами
$\displaystyle\left(\frac{1}{|K|_{p'}} \right)^\st$
по всем $K\in\K(G_{p'}|D)$.
Из этих рассуждений следует, что число $|\bl(G|D)|$ равно рангу матрицы $A(D)^\st$.
\end{proof}

\textuprn{Пусть $P\in \Syl_p(G)$. Показать, что

$(i)$ если $P\nor G$, то $|\bl(G)|=\big|\{K\in \K(G_{p'})\mid K\se \C_G(P)\}\big|$;

$(ii)$ $\big|\bl(G|P)\big|=\big|\{K\in \K(G_{p'})\mid K\cap \C_G(P)\ne \varnothing\}\big|$.}

Напомним, что группа называется \mem{$p$-разрешимой}, \glsadd{iGrPSol} если каждый её композиционный фактор является $p$- или $p'$-группой. Через \gls{opprG} обозначается наибольшая нормальная $p'$-подгруппа группы $G$.

\textuprn{Пусть $G$ --- $p$-разрешимая группа и $\OO_{p'}(G)=1$. Доказать, что $|\bl(G)|$=1.
\uk{Использовать тот факт, что $\C_{G}(\OO_p(G))\le \OO_p(G)$.}}

%

\glsadd{iSubGrFrt}\mem{Подгруппой Фраттини} \gls{FrlGr}  группы $G$ называется пересечение всех её максимальных подгрупп.\footnote{Обычно
подгруппа Фраттини обозначается символом $\Phi(G)$, однако мы используем его для таблицы модулярных характеров.} В случае,
когда $G$ --- $p$-группа, $\Fr(G)$ совпадает с наименьшей нормальной подгруппой, факторгруппа по которой
элементарная абелева.

Следующее  упражнение показывает, что нахождение числа $p$-блоков группы $G$ с нормальной дефектной
группой $D$ сводится к случаю, когда $D$ элементарная абелева.

\textuprn{Пусть $D$ --- нормальная $p$-подгруппа группы $G$. Показать, что
$$|\bl(G|D)|=\big|\bl\big(\,G/\Fr(D)\mid D/\Fr(D)\,\big)\big|.$$
\uk{Воспользоваться тем, что если $p$-регулярный элемент $x$ действует тождественно на факторах
нормального ряда $p$-группы $P$ или на факторгруппе $P/\Fr(P)$, то $x$ действует тождественно на $P$.}}

\begin{thm}[Грина] \label{tgr}\glsadd{iThmGrn} Пусть $D$ --- дефектная группа некоторого блока группы $G$ и $P$ --- $p$-силовская
подгруппа, содержащая $D$. Тогда существует $p$-регулярный элемент $x\in G$ такой, что $D\in \D(x^G)$ и $P\cap P^x=D$.
\end{thm}
\begin{proof} Пусть $|D|=p^d$ и $|P|=p^a$.

Сначала мы докажем утверждение в случае, когда $D\nor G$. По условию $|\bl(G|D)|>0$
и из теоремы Робинсона \ref{trob} следует, что
существуют $K,L\in \{G_{p'}|D\}$ такие, что
$$\frac{|\A_{K,L}|}{p^{a-d}}\not\equiv 0 \mod{p}.$$
Обозначим $T=\{Pg\mid g\in G\}$. Тогда
$$
\A_{K,L}=\{(y,z)\in K\times L\mid Py=Pz\}=\bigcup_{Pg\in T}(Pg\cap K)\times(Pg\cap L),
$$
и значит,
$$
|\A_{K,L}|=\sum_{Pg\in T}|Pg\cap K||Pg\cap L|.
$$
Заметим, что $P$ действует на множестве $T$ правыми умножениями.  Рассмотрим мощность орбиты класса $Pg\in T$
относительно этого действия. Поскольку для $u\in P$ равенство $Pgu=Pg$ выполнено тогда и только тогда,
когда $u\in P^g$, мощность орбиты класса $Pg$ равна $|P:\,P\cap P^g|$.

Заметим также, что для $u\in P$ справедливы равенства
$$
(Pg\cap K)^u=(Pg)^u\cap K^u=Pgu\cap K,\myeqno\label{zv1}
$$
и значит, $|Pg\cap K|=|Pgu\cap K|$. Аналогично, $|Pg\cap L|=|Pgu\cap L|$.

Обозначим через $\OOO$ набор представителей орбит множества $T$ относительно
действия подгруппы $P$ правыми умножениями. Тогда из сделанных замечаний получаем
$$
|\A_{K,L}|=\sum_{Pg\in \OOO}|P:\,P\cap P^g||Pg\cap K||Pg\cap L|.
$$

Поскольку, как мы заметили, подгруппа $P\cap P^g$, действуя правыми умножениями, оставляет
неподвижным смежный класс $Pg$, из \ref{zv1} следует, что она же, действуя сопряжением, оставляет
неподвижным пересечение $Pg\cap K$. При этом мощность орбиты произвольного
элемента $x\in Pg\cap K$ равна $|P\cap P^g:\C_{P\cap P^g}(x)|$. В силу того, что $D\in \D(K)$ и
$D\nor G$, справедливо включение $\C_{P\cap P^g}(x)\le D$. Однако, $D\le P$ и $D\le \C_G(x)$,
поэтому справедливо и обратное включение, т.\,е. $\C_{P\cap P^g}(x)=D$.

Отсюда следует, что $|P\cap P^g:D|$ делит $|Pg\cap K|$ и, аналогично, $|Pg\cap L|$.
А поскольку
$$p^{a-d}=|P:D|=|P:P\cap P^g||P\cap P^g:D|,$$ число
$$
\frac{|P:\,P\cap P^g||Pg\cap K||Pg\cap L|}{p^{a-d}}
$$
делится на $|P\cap P^g:D|$.
 Итак, ввиду того, что
$$
\frac{|\A_{K,L}|}{p^{a-d}}=
\sum_{Pg\in \OOO}\frac{|P:\,P\cap P^g||Pg\cap K||Pg\cap L|}{p^{a-d}}\not\equiv 0 \mod{p},
$$
должен существовать $g\in G$, для которого $Pg\in \OOO$, $|P\cap P^g:D|=1$ и $|Pg\cap K|\ne 0$.
Пусть $x\in Pg\cap K$. Тогда $P\cap P^x=P\cap P^g=D$,  $D\in \D(x^G)$, и значит, $x$ --- искомый
$p$-регулярный элемент.

Теперь рассмотрим общий случай. Обозначим $N=\N_G(D)$. Выберем $S\in \Syl_p(N)$ так, чтобы $P\cap N\le S$.
Из первой основной теоремы Брауэра \ref{thm br 1} следует,
что $D$ является дефектной группой некоторого блока группы $N$. Поэтому из доказанного выше следует,
что существует $p$-регулярный элемент $x\in N$ такой, что $D\in \D(x^N)$ и $S\cap S^x=D$.
В силу \ref{p np}, справедливо равенство $x^N=x^G\cap C$, где $C=\C_G(D)$, причём $D\in \D(x^G)$.

Осталось показать, что $P\cap P^x=D$. По условию имеем $D\le P$. Поскольку $x\in N$, также имеем $D\le P^x$,
т.\,е. $D\le P\cap P^x$.
Допустим, что $D<P\cap P^x$. Тогда по известному свойству $p$-групп
$$D<\N_{P\cap P^x}(D)=P\cap P^x\cap N=(P\cap N)\cap (P\cap N)^x\le S\cap S^x.$$
Это противоречие завершает доказательство теоремы.
\end{proof}

\section{Высшие числа разложения. Вторая основная теорема Брауэра}

Прежде, чем определять высшие числа разложения и формулировать теорему Брауэра,
докажем ряд технических утверждений.

Напомним, что для любого $B\in \bl(G)$ элемент $f_{\mbox{}_B}$ является идемпотентом алгебры $\Z(\wt{Z} G)$,
см. \ref{bl ob}$(ii)$.
Кроме того, как мы уже отмечали, отображение $^\st:\wt{Z}\to F$
можно естественно поднять до эпиморфизма колец $\Z(\wt{Z} G)\to \Z(F G)$,
также обозначаемого через $^\st$.

\begin{pre} \label{blob} Пусть
$B\in \bl(G)$ и $x\in \Z(\wt{Z} G)$. Если $\l_{\mbox{}_B}(x^\st)=1$,
то существует $y\in f_{\mbox{}_B}\Z(\wt{Z} G)$ такой, что $xy=f_{\mbox{}_B}$.
\end{pre}
\begin{proof} Предположим сначала, что $x$ лежит в кольце $f_{\mbox{}_B}\Z(\wt{Z} G)$.
Тогда требуется показать, что $x$ обратим в этом кольце. Имеем $f_{\mbox{}_B}x=x$ и,
значит, $e_{\mbox{}_B}x^\st=x^\st$. Поэтому в силу \ref{e b}$(iii)$
$$
\l_{\mbox{}_{B'}}(x^\st-e_{\mbox{}_B})=\l_{\mbox{}_{B'}}(x^\st)=\l_{\mbox{}_{B'}}(e_{\mbox{}_B}x^\st)=0
$$
для всех блоков $B'\ne B$. Кроме того $\l_{\mbox{}_B}(x^\st-e_{\mbox{}_B})=1-1=0$ по условию. Из \ref{e b}$(vi)$
следует, что элемент $e_{\mbox{}_B}-x^\st$, лежащий в  $e_{\mbox{}_B}\Z(FG)$, является нильпотентным.
В силу \ref{com nil} и \ref{rad prop} элемент $x^\st$ обратим в кольце $e_{\mbox{}_B}\Z(FG)$. Значит,
существует $u\in f_{\mbox{}_B}\Z(\wt{Z} G)$ такой, что
$$(xu)^\st=e_{\mbox{}_B}=f_{\mbox{}_B}^\st.$$
Поэтому мы можем применить
утверждение \ref{l54} к гомоморфизму коммутативных колец $^\st:\wt{Z}\to F$, идемпотенту
$f_{\mbox{}_B}\in\wt{Z} G$ и
элементу $xu\in f_{\mbox{}_B}\Z(\wt{Z} G)\le f_{\mbox{}_B}\wt{Z} G$ \big(т.\,к. мы отмечали, что
ядро отображения $^\st:\wt{Z}\to F$ совпадает с $\J(\wt{Z})=\wtm$\big).
Из \ref{l54} следует, что существует $v\in f_{\mbox{}_B}\wt{Z} G$ для
которого $xuv=f_{\mbox{}_B}$. Таким образом, $uv$ --- обратный к $x$ в кольце $\wt{Z} G$.
Имеем $x\in f_{\mbox{}_B}\Z(\wt{Z} G)=\Z(f_{\mbox{}_B}\wt{Z} G)$ в силу \ref{rtg}$(ii)$.
Поэтому $x$ обратим в $\Z(f_{\mbox{}_B}\wt{Z} G)$ по \ref{zr}$(i)$ и,
значит, $uv\in f_{\mbox{}_B}\Z(\wt{Z} G)$, что и требовалось доказать.

Теперь рассмотрим общий случай. Пусть $x\in \Z(\wt{Z} G)$ удовлетворяет
условию. Тогда $f_{\mbox{}_B}x\in f_{\mbox{}_B}\Z(\wt{Z} G)$ и
$$
\l_{\mbox{}_B}\big(\,(f_{\mbox{}_B}x)^\st\big)=\l_{\mbox{}_B}(e_{\mbox{}_B})\l_{\mbox{}_B}(x^\st)=1.
$$
По доказанному выше существует $y\in f_{\mbox{}_B}\Z(\wt{Z} G)$ такой, что $f_{\mbox{}_B}x y=f_{\mbox{}_B}$,
и, значит, $f_{\mbox{}_B}y$ --- искомый элемент из $f_{\mbox{}_B}\Z(\wt{Z} G)$.
\end{proof}

Пусть $R$ --- коммутативное кольцо. \mem{Носителем}\glsadd{iSuppElmGrAlg}
\gls{supp} элемента
$$
x= \sum_{g\in G} a_g g,\quad a_g\in R
$$
групповой алгебры $RG$ называется подмножество $\{g\in G\mid a_g\ne 0\}$ группы $G$.

\begin{pre} \label{exw} Пусть $b\in \bl(H)$, где $H\le G$, и предположим,
что определён $p$-блок $B=b^G\in \bl(G)$.
Тогда справедливы следующие утверждения.
\begin{list}{{\rm(}{\it\roman{enumi}\/}{\rm)}}{\usecounter{enumi}\setlength{\parsep}{2pt}\setlength{\topsep}{5pt}\setlength{\labelwidth}{23pt}}
\item Существует элемент $w\in f_b\wt{Z} G f_b$, для которого
    \begin{list}{{\rm(}{\it\roman{enumi}\/}.{\rm\arabic{enumii})}}
    {\usecounter{enumii}\setlength{\parsep}{2pt}\setlength{\topsep}{2pt}\setlength{\labelwidth}{5pt}\setlength{\labelwidth}{50pt}}
    \item $\supp w\se G\setminus H$;
    \item $H$ централизует $w$;
    \item $(1-f_{\mbox{}_B})f_b^{\vphantom A}=(1-f_{\mbox{}_B})w$.
    \end{list}
\item Для любого $w\in \wt{Z} G$, обладающего свойствами $(i.1)$--$\,(i.3)$,
и любого элемента $h\in H$ такого, что $\C_G(h_p)\le H$, где $h_p$ --- $p$-часть элемента $h$, существуют
элементы $w_0,\ld,w_{p-1}\in f_b\wt{Z} G f_b$ такие, что $w=w_0+\ld +w_{p-1}$ и $w_i^h=w_{i+1}$, $i=0,\ld,p-1$,
где индексы берутся по модулю $p$.
\end{list}
\end{pre}
\begin{proof} $(i)$ Запишем $f_{\mbox{}_B}=u-v$, где $u,v$ --- однозначно определённые
элементы из $\wt{Z} G$ такие, что $\supp u\se H$ и $\supp v\se G\setminus H$. Заметим,
что $u\in \Z(\wt{Z} H)$. Действительно,
$f_{\mbox{}_B}\in \Z(\wt{Z} G)$ и, значит, $u$ является линейной комбинацией сумм
вида $\sum_{x\in K\cap H}$, где $K\in \K(G)$, а множества $K\cap H$ являются объединениями классов сопряжённости
группы $H$.

Поскольку $B=b^G$, имеем $\l_{\mbox{}_B}=\l_b^G$. Значит,
$$
1=\l_{\mbox{}_B}(e_{\mbox{}_B})=\l_{\mbox{}_B}(f_{\mbox{}_B}^\st)=\l_b^G(u^\st-v^\st)=\l_b(u^\st),
$$
ввиду того, что $\supp v\cap H=\varnothing$. Из \ref{blob} следует, что найдётся элемент
$y\in f_b\Z(\wt{Z} H)$ такой, что $uy=f_b$. Покажем, что элемент $w=vy$ удовлетворяет $(i.1)$--$(i.3)$
и лежит в $f_b\wt{Z} G f_b$, т.\,е. является искомым.

Так как $\supp v\se G\setminus H$ и $\supp y\se H$, имеем $(i.1)$.
Ввиду того, что $H$ централизует элементы $y,u,f_{\mbox{}_B}$ и
выполнено равенство $w=(u-f_{\mbox{}_B})y$, получаем $(i.2)$.
Так как $f_b$ --- единица кольца $f_b\Z(\wt{Z} H)$ и $y$ --- его элемент, то
$$wf_b=vyf_b=vy=w,$$
а из $(i.2)$ получаем $f_bw=wf_b$, т.\,е. $w\in f_b\wt{Z} G f_b$. Наконец,
$$
(1-f_{\mbox{}_B})(f_b^{\vphantom A}-w)=(1-f_{\mbox{}_B})(uy-vy)=
(1-f_{\mbox{}_B})(u-v)y=(1-f_{\mbox{}_B})f_{\mbox{}_B}y=0,
$$
откуда вытекает $(i.3)$.

$(ii)$ Поскольку группа $H$  централизует $w$, она действует сопряжением на носителе $\supp w$ и, кроме того,
коэффициент элемента $w$ при произвольном $g\in G$ равен коэффициенту при $g^u$ для любого $u\in H$. Заметим,
что мощность любой орбиты группы $\la h \ra$ на $\supp w$ делится на $p$. В самом деле, если
это не так, то найдётся $g\in \supp w$ такой, что $|\la h \ra : \C_{\la h \ra}(g)|$ не делится на $p$,
т.\,е. $\la h_p\ra \le \C_{\la h \ra}(g)$. Но тогда по условию $g\in \C_G(h_p)\le H$ вопреки тому,
что $g\in \supp w\se G\setminus H$ по $(i.1)$.

Пусть
$\OOO$ --- полное множество
представителей $\la h \ra$-орбит множества $\supp w$. Тогда
$$
\supp w=\dot\bigcup_{g\in \OOO}\Orb(g).
$$
Из сказанного выше следует, что для любого $g\in \supp w$ индекс $\big|\la h \ra:\C_{\la h \ra}(g)\big|$
делится на $p$, и поэтому $\C_{\la h \ra}(g)\le \la h^p\ra$. Значит, обозначив $\la h^p\ra$-орбиту элемента $g$
через $T(g)$, получим, что все множества $T(g)^{h^i}$, $i=0,\ld,p-1$, попарно не пересекаются, и
$$
\Orb(g)=T(g)\dot\cup T(g)^h\dot\cup \ld \dot\cup T(g)^{h^{p-1}}.
$$
Следовательно, $w$ допускает однозначное представление в виде $w=w_0'+\ld+w_{p-1}'$, где
$$\supp w_i'=\dot\bigcup_{g\in \OOO}T(g)^{h^i},\quad i=0,\ld,p-1.$$
В частности, $(w_i')^h=w_{i+1}'$, $i=0,\ld,p-1$, где индексы берутся по модулю $p$. Так как $w\in f_b\wt{Z} G
f_b$, то
$$w=f_bwf_b=f_bw_0'f_b+\ld+f_bw_{p-1}'f_b=w_0+\ld+w_{p-1}$$
--- требуемое представление, где $w_i=f_bw_i'f_b\in f_b \wt{Z} G f_b$, $i=0,\ld,p-1$,
поскольку
$$w_i^h=(f_bw_i'f_b)^h=f_b(w_i')^hf_b=f_bw_{i+1}'f_b=w_{i+1}$$
для всех $i$.
\end{proof}

\begin{pre} \label{is} Пусть $b\in \bl(H)$, где $H\le G$, и предположим,
что определён $p$-блок $B=b^G\in \bl(G)$.
Пусть $h\in H$ такой, что $\C_G(h_p)\le H$, где $h_p$ --- $p$-часть элемента $h$.
Тогда $\x(f_bh)=0$ для любого $\x\in \irr(G)\setminus \iBr(B)$.
\end{pre}

\begin{proof} Пусть $M=M_\x$ --- неприводимый $\CC G$-модуль с характером $\x$.
В силу \ref{fbz} имеем $f_b\in \Z( \wt{Z} H)$ и, значит, из \ref{m id dec}$(iii)$ следует, что
модуль $M$, рассматриваемый как $\CC H$-модуль, допускает разложение
$$M=Mf_b\oplus M(1-f_b).$$
Поскольку  элемент $f_bh$ аннулирует $M(1-f_b)$ и действует на
$Mf_b$ также, как $h$, величина $\x(f_bh)$ равна сумме характеристических
значений элемента $h$ (с учётом кратности) на $Mf_b$.

Положим $V=Mf_b$. Пусть $\z=e^{\frac{2\pi i}{p}}$. Достаточно показать, что для любого $\a\in \CC$ кратность
$\a$ как характеристического значения элемента $h$ на $V$ совпадает с кратностью $\z\a$. В самом деле, если это
верно, то кратности значений $\a,\z\a,\ld,\z^{p-1}\a$ совпадают. А поскольку $1+\z+\ld+\z^{p-1}=0$, отсюда
получим требуемое.

Поскольку $V$, рассматриваемый как $\CC\la h\ra$-модуль, вполне приводим и неприводимые $\CC\la h\ra$-модули
одномерны, $V$ обладает базисом из собственных векторов элемента $h$, в частности, каждое характеристическое
значение элемента $h$ на $V$ является собственным. Обозначим через $U_\a$ подпространство собственных векторов
элемента $h$, соответствующее собственному значению $\a$. Тогда кратность $\a$ равна $\dim_{\CC}(U_\a)$. Пусть
$w,w_0,\ld,w_{p-1}\in f_b \wt{Z} G f_b$ --- элементы, существование которых утверждается в предложении
\ref{exw} для данных подгруппы $H$, $p$-блока $B=b^G$ и элемента $h\in H$. Заметим, что, поскольку $f_b$ ---
единица кольца $f_b \wt{Z} G f_b$, справедливы соотношения
$$
Vw_i=Mf_b w_i=Mw_if_b\se Mf_b=V.
$$
Поэтому для элемента $s\in \wt{Z} G$, определённого равенством
$$
s=w_0+\z^{-1}w_1+\ld +\z^{-(p-1)}w_{p-1}, \myeqno\label{spr}
$$
имеем $Vs\se V$. Кроме того
$$
s^h=w_0^h+\z^{-1}w_1^h+\ld +\z^{-(p-1)}w_{p-1}^h=w_1+\z^{-1}w_2+\ld +\z^{-(p-1)}w_0=\z s.
$$
Поэтому для любого $u\in U_\a$ получаем
$$
(us)h=uhs^h=\z\a(us).
$$
Другими словами действие элемента $s$ переводит собственное подпространство $U_\a$ в
подпространство $U_{\z\a}$. Если мы покажем, что отображение $v\mapsto vs$ инъективно на $V$,
то получим, что
$$
\dim_{\CC}(U_\a)=\dim_{\CC}(U_\a s)\le\dim_{\CC}(U_{\z\a})
$$
и, аналогично,
$$
\dim_{\CC}(U_\a)\le\dim_{\CC}(U_{\z\a})\le\dim_{\CC}(U_{\z^2\a})\le\ld
\le\dim_{\CC}(U_{\z^p\a})=\dim_{\CC}(U_\a),
$$
откуда будет следовать требуемое.

Обозначим $f=(1-f_{\mbox{}_B})f_b^{\vphantom A}$. Тогда $f$ --- идемпотент алгебры $\wt{Z} G$,
поскольку является произведением перестановочных идемпотентов.

Заметим, что $(1-f_{\mbox{}_B})s\in f \wt{Z} G f$. Действительно, с одной стороны имеем $f_b s f_b=s$, что
следует из \ref{spr} и того, что $f_b w_i f_b=w_i$, $i=0,\ld, p-1$, а с другой стороны $1-f_{\mbox{}_B}\in
\Z(\wt{Z} G)$ и, значит,
$$
(1-f_{\mbox{}_B})s=(1-f_{\mbox{}_B})f_b^{\vphantom A}sf_b^{\vphantom A}=
(1-f_{\mbox{}_B})f_b^{\vphantom A}\big((1-f_{\mbox{}_B})s\big)(1-f_{\mbox{}_B})f_b^{\vphantom A}=f(1-f_{\mbox{}_B})sf.
$$

Также заметим, что $\z^\st=1$, поскольку в группе $F^\times$ нет нетривиальных $p$-элементов.
Поэтому $s^\st=w^\st$. Тогда из \ref{exw}$(i.3)$ получаем
$$
f^\st=\big((1-f_{\mbox{}_B})w\big)^\st=(1-f_{\mbox{}_B})^\st w^\st=
(1-f_{\mbox{}_B})^\st s^\st=\big((1-f_{\mbox{}_B}) s\big)^\st.
$$
Поэтому, применяя утверждение \ref{l54} к гомоморфизму коммутативных
колец $^\st:\wt{Z}\to F$, идемпотенту
$f\in\wt{Z} G$ и элементу $(1-f_{\mbox{}_B})s\in f \wt{Z} G f$, заключаем, что существует $y\in \wt{Z} G$,
для которого $(1-f_{\mbox{}_B})sy=f$.

Элемент $f$ действует тождественно на $V$. В самом деле, поскольку $\x\not \in B$, из
\ref{ir mod bl}$(i)$ следует, что идемпотент $1-f_{\mbox{}_B}$
действует тождественно даже на всём модуле $M$,
а $f_b$ действует тождественно на $V$ в силу равенства $V=Mf_b$.

Из сделанных замечаний следует, что для произвольного $v\in V$
$$
vsy=v(1-f_{\mbox{}_B})sy=vf=v,
$$
и поэтому если $vs=0$, то $v=0$. Тем самым действие $s$ на $V$ инъективно и утверждение доказано.
\end{proof}

Теперь приступим к формулировке основных понятий и результатов этого раздела.
Для определения обобщения чисел разложения докажем следующее утверждение.

\begin{pre} \label{dnp} Пусть $x$ --- $p$-элемент группы $G$, $H=\C_G(x)$ и $\x\in \irr(G)$.
Тогда

$(i)$ для произвольного $\vf\in \iBr(H)$ существует единственное число $d_{\x\vf}^{\,x}\in \CC$
такое, что равенство
$$
\x(xy)=\sum_{\vf\in \iBr(H)}d_{\x\vf}^{\,x}\vf(y)
$$
справедливо для всех $p$-регулярных $y\in H$;

$(ii)$ $d_{\x\vf}^{\,x}\in \QQ_{|x|}\cap \ov{\ZZ}$ для всех $\vf\in \iBr(H)$.

\end{pre}
\begin{proof} $(i)$ Для ограничения $\x_{H^{\vphantom{A^a}}}$ запишем
$$
\x_{H^{\vphantom{A^a}}}=\sum_{\t\in \irr(H)}a_{\t}\t,
$$
где $a_\t=(\x_{H^{\vphantom{A^a}}},\t)_{\mbox{}_H}\in \ZZ$. Поскольку $x\in \Z(H)$, из \ref{om ob}$(ii)$ следует, что
$$
\x(xy)=\sum_{\t\in \irr(H)}a_{\t}\t(xy)=\sum_{\t\in \irr(H)}a_{\t}\om_\t(x)\t(y)
=\sum_{\t\in \irr(H)}\sum_{\vf\in \iBr(H)}a_{\t}\om_\t(x)d_{\t\vf}\vf(y)
=\sum_{\vf\in \iBr(H)}d_{\x\vf}^{\,x}\vf(y),
$$
где
$$
d_{\x\vf}^{\,x} = \sum_{\t\in \irr(H)}a_\t\, d_{\t\vf}\,\om_\t(x)=
\sum_{\t\in \irr(H)} \frac{(\x_{H^{\vphantom{A^a}}},\t)_{\mbox{}_H}d_{\t\vf}}{\t(1)}\,\t(x).
\myeqno\label{vyse}
$$
Единственность чисел $d_{\x\vf}^{\,x}$ следует из линейной независимости над $\CC$
множества $\iBr(H)$, см. \ref{cor col br}$(i)$.

$(ii)$ следует из выражения \ref{vyse} и \ref{om ob}$(iii)$.
\end{proof}

\textuprn{Показать, что $d^{\,x}_{\x\vf}\in \ZZ[\z]$, где $\z=e^{\frac{2\pi i}{|x|}}$, т.\,е.
$d^{\,x}_{\x\vf}$ является $\ZZ$-линейной комбинацией чисел $\z^i$, при $i=0,\ld,|x|-1$.}

\begin{opr} \label{hidec} Для $p$-элемента $x\in G$ числа \gls{dxxvf}, где  $\x\in \irr(G)$ и $\vf\in \iBr(\C_G(x))$,
определённые в \ref{dnp} называются \glsadd{iNumsDecHghr}\mem{высшими числами разложения}.
\end{opr}

Легко видеть, что обычные числа разложения $d_{\x\vf}$ соответствуют случаю $x=1$.

\begin{thm}[Вторая основная теорема Брауэра] \label{bro2}\glsadd{iThmBraScdMain} Пусть $\x\in \irr(G)$ и $\vf\in \iBr(H)$,
где $H=\C_G(x)$ для некоторого $p$-элемента $x \in G$. Пусть $B\in \bl(G)$ и $b\in \bl(H)$ такие, что
$\x\in B$ и $\vf \in b$. Тогда если $b^G\ne B$, то $d_{\x\vf}^{\,x}=0$.
\end{thm}
\begin{proof}
Заметим, что блок $b^G$ определён в силу \ref{comp ind}. В силу линейной
независимости брауэровых достаточно доказать, что если $b^G\ne B$, то
линейная комбинация $\sum_{\psi\in\iBr(b)}d_{\x\psi}^{\,x}\psi$ является нулевым элементом
алгебры $\cf(H_{p'})$.  По \ref{vyse} имеем
\begin{align*}
\sum_{\psi\in\iBr(b)}d_{\x\psi}^{\,x}\psi=
\sum_{\psi\in\iBr(b)} \left(\sum_{\t\in \irr(H)}  a_\t\, d_{\t\psi}\,\om_\t(x)\right)\psi
=\sum_{\psi\in\iBr(b)} \left(\sum_{\t\in \irr(b)}  a_\t\, d_{\t\psi}\,\om_\t(x)\right)\psi\\
=\sum_{\t\in \irr(b)} a_\t\, \om_\t(x)\left(\sum_{\psi\in\iBr(b)}d_{\t\psi}\psi\right)
=\sum_{\t\in \irr(b)} a_\t\, \om_\t(x)\,\t,
\end{align*}
где $a_\t=(\x_{H^{\vphantom{A^a}}},\t)_{\mbox{}_H}$.
Достаточно проверить, что значение последнего выражения на произвольном элементе $y\in H_{p'}$ равно нулю.
Поскольку $x\in\Z(H)$, из \ref{om ob}$(iii)$ и \ref{bch} следует, что
\begin{align*}
\sum_{\t\in \irr(b)} a_\t\,\om_\t(x)\,\t(y)&=
\sum_{\t\in \irr(b)} a_\t\,\t(xy)\\
&=\sum_{\t\in \irr(H)} a_\t\,\t_b(xy)
=\sum_{\t\in \irr(H)} (\x_{H^{\vphantom{A^a}}},\t)_{\mbox{}_H}\t(f_b\,xy)=\x(f_b\,xy).
\end{align*}
Теперь применим предложение \ref{is} к данной группе $H$  и элементу $h=xy\in H$.
Поскольку $x$ и $y$ перестановочны, имеем $h_p=x$. Поэтому, если $b^G\ne B$,
то $\x\not\in b^G$, и из \ref{is} вытекает, что $\x(f_b\,xy)=\x(f_bh)=0$,
как и требовалось. \end{proof}

\section{Следствия из второй основной теоремы Брауэра}

Вторая основная теорема Брауэра является мощным и чрезвычайно полезным инструментом.
Выведем из неё некоторые следствия.

\begin{pre} \label{chnz} Пусть $\x\in \irr(G)$ и $g\in G$. Если $\x(g)\ne 0$,  то $g_p$ принадлежит
некоторой дефектной группе $p$-блока, содержащего $\x$.
\end{pre}
\begin{proof} Пусть $x=g_p$ и запишем $g=xy$ для подходящего $p$-регулярного элемента $y\in C=\C_G(x)$.
Тогда
$$
\x(g)=\x(xy)=\sum_{\vf\in \iBr(C)}d_{\x\vf}^{\,x}\vf(y).
$$
По условию существует характер $\vf\in \iBr(C)$ такой, что $d_{\x\vf}^{\,x}\ne 0$. Из \ref{bro2} следует, что
если $\vf\in b\in \bl(C)$, то $\x\in b^G$. В силу \ref{op def} имеем $x\in \OO_p(C)\le \d(b)$,
а из \ref{ind def} получаем $\d(b)\le_G\d(b^G)$. Поэтому $x$ лежит в некоторой дефектной группе блока $b^G$.
\end{proof}

Отметим, что прямым следствием из \ref{chnz} является доказанное в \ref{nul def} утверждение о том, что
обыкновенный характер из $p$-блока дефекта $0$ обращается в $0$ на всех не-$p$-регулярных элементах.

\begin{opr} \label{psech} Пусть $g\in G$. Множество
$$
\gls{SsGlgr}=\{x\in G\mid (x_p)^G=(g_p)^G\}
$$
назовём \mem{$p$-сечением}\glsadd{iPSecGr} группы $G$, содержащем элемент $g$. Другими словами, $\SS_G(g)$ является
множеством тех элементов из $G$, $p$-часть которых сопряжена с $p$-частью элемента $g$.
\end{opr}

\begin{pre} \label{svp} Справедливы следующие утверждения.
\begin{list}{{\rm(}{\it\roman{enumi}\/}{\rm)}}
            {\usecounter{enumi}\setlength{\parsep}{2pt}\setlength{\topsep}{5pt}\setlength{\labelwidth}{20pt}
            }
\item $\SS_G(1)=G_{p'}$.
\item $g\in \SS_G(g)=\SS_G({g_p})$ для любого $g\in G$.
\item Для любого $g\in G$ множество $\SS_G(g)$ является объединением
некоторых классов сопряжённости группы $G$.
\item Пусть $U$ --- полное множество представителей классов сопряжённости $p$-элементов
группы $G$. Тогда
$$
G=\dot\bigcup_{u\in U}\SS_G(u).
$$
\item Пусть $g\in G$ --- $p$-элемент и $T$ --- полное множество представителей классов сопряжённости $p$-регулярных
элементов группы $C=\C_G(g)$. Если $C\le H\le G$, то
 $$\SS_H(g)=\dot\bigcup_{t\in T}(gt)^H.$$
\end{list}
\end{pre}
\textupl{svp prf}{Доказать предложение \ref{svp}.}

\begin{pre} \label{psv} Пусть $x\in G$ и $\eta\in \cf(G)$. Если $\eta$ тождественно равно нулю на $\SS_G(x)$,
то $B$-часть $\eta_{\mbox{}_B}$ для любого $B\in \bl(G)$ также тождественно равна нулю на $\SS_G(x)$.
\end{pre}
\begin{proof} Пусть сначала $B$ --- произвольный блок группы $G$. В силу \ref{svp}$(ii)$ можно считать,
что $x$ --- $p$-элемент.
Достаточно показать, что $\eta_{\mbox{}_B}(xt)=0$ для любого $p$-регулярного элемента $t\in C=\C_G(x)$.
Справедливы соотношения
$$
\eta_{\mbox{}_B}(xt)=\sum_{\x\in \irr(B)}(\eta,\x)_{\mbox{}_G}\x(xt)=\sum_{\x\in \irr(B)}(\eta,\x)_{\mbox{}_G}
\!\!\sum_{\vf\in \iBr(C)}d_{\x\vf}^{\,x}\vf(t)=
\sum_{\vf\in \iBr(C)}\left(\sum_{\x\in \irr(B)}(\eta,\x)_{\mbox{}_G}d_{\x\vf}^{\,x}\right)\vf(t).
\myeqno\label{etxt}
$$
Просуммировав по всем $B\in \bl(G)$, получим
$$
\sum_{\vf\in \iBr(C)}\left(\sum_{\x\in
\irr(G)}(\eta,\x)_{\mbox{}_G}d_{\x\vf}^{\,x}\right)\vf(t)=\sum_{B\in\bl(G)}\eta_{\mbox{}_B}(xt)=\eta(xt)=0,
\myeqno\label{lnez}
$$
где последнее равенство следует из условия.

Теперь фиксируем блок $B\in \bl(G)$.
Тогда в силу линейной независимости множества $\iBr(C)$
из \ref{lnez} следует, что
для любого $\vf\in \iBr(C)$ выполнены равенства
$$
0=\sum_{\x\in \irr(G)}(\eta,\x)_{\mbox{}_G}d_{\x\vf}^{\,x}=
\sum_{\x\in \iBr(B)}(\eta,\x)_{\mbox{}_G}d_{\x\vf}^{\,x}+
\sum_{\x\in \irr(G)\setminus\iBr(B)}(\eta,\x)_{\mbox{}_G}d_{\x\vf}^{\,x}. \myeqno\label{sraz}
$$
Заметим, что в случае, когда $b^G=B$, где $\vf\in b\in\bl(C)$, последнее слагаемое в \ref{sraz} равно
нулю в силу теоремы \ref{bro2}, т.\,е. в этому случае имеем
$$
\sum_{\x\in \iBr(B)}(\eta,\x)_{\mbox{}_G}d_{\x\vf}^{\,x}=0.
$$
Однако, это равенство также справедливо и в случае, когда $b^G\ne B$, опять-таки в силу теоремы
\ref{bro2}. Поэтому из \ref{etxt} получаем
$\eta_{\mbox{}_B}(xt)=0$, как и требовалось.
\end{proof}

Теперь докажем усиление утверждения \ref{cor sl ort}.

\begin{pre}[Ортогональность в блоке] \label{ortbl}\glsadd{iOrthBlk} Пусть $x,y\in G$ такие, что $p$-части $x_p$ и $y_p$
не сопряжены в $G$.
Тогда
$$
\sum_{\x\in \irr(B)}\x(x)\ov{\x(y)}=0
$$
для любого $B\in \bl(G)$.
\end{pre}
\begin{proof} Определим классовую функцию
$$\eta=\sum_{\x\in \irr(G)}\ov{\x(y)}\x$$
и положим $u=x_p$. Если $t\in \C_G(u)$ --- $p$-регулярный элемент, то элементы $ut$ и $y$ не сопряжены в $G$,
поскольку не сопряжены их $p$-части $(ut)_p=u=x_p$ и $y_p$. Поэтому из второго соотношения ортогональности
\ref{vtor ort} следует, что $\eta(ut)=0$. Значит, $\eta$ обращается в ноль на $\SS_G(x)$. Из \ref{psv}
вытекает, что
$$
0=\eta_{\mbox{}_B}(x)=\sum_{\x\in \irr(B)}\x(x)\ov{\x(y)},
$$
как и требовалось.\end{proof}

Для дальнейшего нам потребуется ещё один вспомогательный результат, который, по сути, является усилением
пункта $(i)$ теоремы min--max \ref{thm min max}.

\begin{pre} \label{xpz} Пусть $P$ --- $p$-подгруппа группы $G$, $B\in \bl(G|P)$ и $K\in\K(G)$.
Если $\l_{\mbox{}_B}(\wh K)\ne 0$, то существует $x\in K\cap\C_G(P)$ такой, что $x_p\in \Z(P)$.
\end{pre}
\begin{proof} Обозначим $C=\C_G(P)$, $N=\N_G(P)$. Из первой основной теоремы Брауэра \ref{thm br 1}
следует, что существует блок $b\in \bl(N|P)$ такой, что $B=b^G$. Если в \ref{comp ind}$(i)$ положить $H=N$, то
получим $\l_b^{\vphantom{G}}\circ\b_{\mbox{}_P}=\l_{\mbox{}_B}$ и, значит,
$$
0\ne \l_{\mbox{}_B}(\wh K) = \l_b^{\vphantom{G}}\left(\,\sum_{x\in K\cap C}x\right).
$$
Поэтому существует класс $L\in \K(N)$ такой, что $L\se K\cap C$ и $\l_b^{\vphantom{G}}(\wh L)\ne 0$.
Выберем произвольный элемент  $x\in L$. Заметим, что $x_p\in C$, поскольку $x\in C$. Кроме того,
для любого $\t\in \irr(b)$ имеем $\t(x)\ne 0$, поскольку
$$
0 \ne \l_b^{\vphantom{G}}(\wh L)=\l_\t^{\vphantom{G}}(\wh L)=\om_\t(\wh L)^\st
=\left(\frac{|L|\t(x)}{\t(1)}\right)^\st
$$
в силу \ref{dz} и \ref{ch val}. Значит, из \ref{chnz} следует, что $x_p$ лежит в некотором элементе из $\D(b)$.
Однако, $\D(b)=\{P\}$ в силу нормальности $P$ в $N$. Поэтому $x_p\in P\cap C =\Z(P)$, т.\,е.
$x$ --- искомый элемент из $K\cap C$.
\end{proof}

В \ref{nul def} было установлено, что всякий неприводимый обыкновенный характер, обращающийся в ноль на всех
неединичных $p$-элементах, принадлежит блоку дефекта $0$. Усилением этого является
следующее утверждение.

\begin{thm}[Кнёрра] \label{tkno}\glsadd{iThmKnr} Пусть $\x\in \irr(G)$ и $\x(x)=0$ для любого $x\in G$ порядка $p$.
Тогда $\x$ принадлежит $p$-блоку дефекта $0$.
\end{thm}
\begin{proof} Пусть $\x\in B\in \bl(G)$ и пусть $P=\d(B)$. Требуется показать, что $P=1$.
Положим
$$U=\big\{x\in G\bigm| |x|=p\big\}.$$
Тогда $U$ --- объединение классов сопряжённости группы $G$. Если $K$ --- один из таких классов,
то по условию имеем
$$
\l_{\mbox{}_B}(\wh K)=\om_\x(\wh K)^\st=\left(\frac{|K|\x(x_{\mbox{}_K})}{\x(1)}\right)^\st=0.
$$
Просуммировав по всем таким классам $K$, получим $\l_{\mbox{}_B}\left(\,\sum_{x\in U}x\right)=0$.

Положим $N=\N_G(P)$ и $C=\C_G(P)$. В силу первой основной теоремы Брауэра \ref{thm br 1} существует блок
$b\in \bl(N|P)$ такой, что $B=b^G$. Кроме того, по \ref{comp ind}$(i)$ имеем
$\l_b^{\vphantom{G}}\circ\b_{\mbox{}_P}=\l_{\mbox{}_B}$. Поэтому
$$
0=\l_{\mbox{}_B}\left(\,\sum_{x\in U}x\right)=\l_b\left(\,\sum_{x\in U\cap C}x\right).
$$
Заметим, что $U\cap C$ --- объединение некоторых классов сопряжённости группы $N$.
Если $L$ --- один из таких классов, причём $\l_b(\wh L)\ne 0$, то, применив
\ref{xpz} к группе $N$, её блоку $b$ и классу $L$, получим, что существует элемент $x\in L$
такой, что $x_p\in \Z(P)$. Однако, $x=x_p$ и $\Z(P)\nor N$, поскольку $P\nor N$. Поэтому $L\se \Z(P)$.
Значит,
$$
0=\l_b\left(\,\sum_{x\in U\cap C}x\right)=\l_b\left(\,\sum_{x\in U\cap \Z(P)}x\right).
$$
Обозначим
$$
H=\{x\in \Z(P)\mid x^p=1\}
$$
и положим $\wh H=\sum_{x\in H}x$. Тогда $H\le \Z(P)$ и $\wh H\in \Z(FN)$. Кроме того,
$U\cap \Z(P)=H\setminus\{1\}$. В частности,
$$
\l_b(\wh H )=\l_b\left(\,\sum_{x\in H\setminus\{1\}}x\right)+\l_b(1)=
\l_b\left(\,\sum_{x\in U\cap \Z(P)}x\right)+1=1.
$$
Кроме того, так как $H$ --- $p$-подгруппа, имеем
$$
\wh H \wh H=\sum_{x\in H} \sum_{y\in H} xy =\sum_{x\in H} \wh H = |H|^\st \wh H=\left\{\ba{ll}
0, & \mbox{если}\ \ H\ne 1;\\
1, & \mbox{если}\ \ H=1.
\ea\right.
$$
Значит, поскольку
$1=\l_b(\wh H)^2=\l_b(\,\wh H\wh H)$, получаем $H=1$, а это возможно лишь когда $P=1$.
Отсюда вытекает требуемое.
\end{proof}

Докажем ещё ряд вспомогательных утверждений.

\begin{pre} \label{chvz} Пусть $P$ --- $p$-подгруппа группы $G$ и $B\in \bl(G|P)$.
Пусть $K$ --- дефектный класс блока $B$ и $x\in K$ --- такой элемент, что $P\in \Syl_p(\C_G(x))$.
Тогда справедливы следующие утверждения.
\begin{list}{{\rm(}{\it\roman{enumi}\/}{\rm)}}
{\usecounter{enumi}\setlength{\parsep}{2pt}\setlength{\topsep}{5pt}\setlength{\labelwidth}{23pt}}
\item Если $\x\in \irr(B)$ --- характер высоты $0$, то $\x(xu)^\st\ne 0$ для любого $u\in P$.
\item Если $u\in \Z(P)$ и $L=(xu)^G$, то $\l_{\mbox{}_B}(\wh L)\ne 0$.
\end{list}
\end{pre}
\begin{proof} $(i)$ По определению дефектного класса $\l_{\mbox{}_B}(\wh K)\ne 0$.
Поскольку $\x$ --- характер высоты $0$, имеем $|\x(1)|_p=|K|_p$. Значит,
$\displaystyle\frac{|K|}{\x(1)}$ --- обратимый элемент кольца $\wt{Z}$ и поэтому
$$
0\ne \l_{\mbox{}_B}(\wh K)=\om_\x(\wh K)^\st=\left(\frac{|K|\x(x)}{\x(1)}\right)^\st=
\left(\frac{|K|}{\x(1)}\right)^\st\x(x)^\st.
$$
Отсюда следует, что $\x(x)^\st\ne 0$.

В силу \ref{dk cor} имеем $x\in G_{p'}$. Поэтому $(xu)_{p'}=x$ для любого $u\in P$.
Из \ref{br preg}$(ii)$ получаем $\x(xu)^\st=\x(x)^\st\ne 0$.

$(ii)$ Поскольку $(xu)_p=u$ и $(xu)_{p'}=x$, в силу \ref{lem ppt} имеем $\C_G(xu)=\C_G(x)\cap \C_G(u)$.
Однако, по условию $P$ централизует $x$ и $u$. Поэтому $P\le \C_G(xu)$. Заметим, что $|\C_G(xu):P|$ делит
индекс $|\C_G(x):P|$, который взаимно прост с $p$ в силу условия $P\in \Syl_p(\C_G(x))$. Поэтому
индекс $|\C_G(xu):P|$ также взаимно прост с $p$ и, значит, $P\in \Syl_p(\C_G(xu))$.
Отсюда следует, что $|L|_p=|K|_p$.

Выберем характер $\x\in B$ высоты $0$. Тогда $|\x(1)|_p=|K|_p=|L|_p$.
Поэтому $\displaystyle\frac{|L|}{\x(1)}$ --- обратимый элемент кольца $\wt{Z}$.
Отсюда, применив $(i)$, получаем
$$
\l_{\mbox{}_B}(\wh L)=\left(\frac{|L|\x(xu)}{\x(1)}\right)^\st=\left(\frac{|L|}{\x(1)}\right)^\st\x(xu)^\st\ne 0,
$$
как и утверждалось.
\end{proof}

Следующее утверждение устанавливает новую связь между дефектными группами блоков $b$ и $b^G$ (ср. \ref{ind def}).

\begin{pre} \label{zds} Пусть $b\in \bl(H)$, где $H\le G$, и предположим, что определён $p$-блок $B=b^G\in \bl(G)$.
Тогда любой элемент из $\Z(\d(B))$ сопряжён в группе $G$ с некоторым элементом из $\Z(\d(b))$.
\end{pre}
\begin{proof} Выберем произвольный элемент $u\in\Z(\d(B))$.
Пусть $K=x^G$ --- дефектный класс блока $B$, причём $\d(B)\in \Syl_p(\C_G(x))$. Обозначим $L=(xu)^G$.
Тогда  $\l_{\mbox{}_B}(\wh L)\ne 0$ в силу \ref{chvz}$(ii)$. По определению индуцированного блока имеем
$$
0\ne \l_{\mbox{}_B}(\wh L)=\l_b\left(\sum_{y\in L\cap H}y\right).
$$
Значит, существует класс $M\in \K(H)$ такой, что $M\se L$ и $\l_b(\whm)\ne 0$. Применив \ref{xpz} к блоку $b$
группы $H$  и классу $M$, получим, что существует $z\in M$, для которого $z_p\in \Z(\d(b))$.
Поскольку элементы $z$ и $xu$ сопряжены в $G$, то сопряжены и их $p$-части $z_p$ и $(xu)_p=u$. Отсюда следует
требуемое.
\end{proof}

Теперь мы можем охарактеризовать блоки, индуцирующие блок дефекта $0$.

\begin{pre} \label{bli} Пусть $b\in \bl(H)$, где $H\le G$, и предположим, что определён $p$-блок $B=b^G\in \bl(G)$.
Тогда $\df(B)=0$ в том и только в том случае, когда $\df(b)=0$.
\end{pre}
\begin{proof} В силу \ref{ind def} имеем $\d(b)\le_G\d(B)$. Поэтому, если $\df(B)=0$, то
$\df(b)=0$. Обратно, предположим, что $\df(b)=0$. Тогда из \ref{zds} следует,
что $\Z(\d(B))=1$, а значит, $\d(B)=1$, откуда следует требуемое.
\end{proof}

\section{Обзор некоторых дальнейших результатов}

В этом разделе результаты приводятся без доказательств.

Аналог утверждения \ref{bli} для главных блоков также имеет место, однако доказывается более сложно.

\begin{thm}[Третья основная теорема Брауэра] \label{brt3}\glsadd{iThmBraThdMain} Пусть $b\in \bl(H)$, где $H\le G$, и предположим,
что определён $p$-блок $B=b^G\in \bl(G)$.
Тогда $B$ --- главный блок группы $G$ в том и только в том случае, когда $b$ --- главный блок группы $H$.
\end{thm}


Рассмотрим пример, иллюстрирующий третью основную теорему Брауэра.

\begin{prim}
Пусть $p=2$, $G=A_5$ --- знакопеременная группа
подстановок множества $\{1,2,3,4,5\}$, $H\cong A_4$ --- подгруппа группы $G$,
состоящая из элементов, для которых $5$ является неподвижным символом. Тогда $H$ имеет единственный
$2$-блок $b$ (главный, см. приложение \ref{pril tch}). Выясним, определён ли индуцированный
блок $b^G$ и совпадает ли он с главным $2$-блоком
группы $G$. Значения гомоморфизма $\l_b$ на классовых суммах $\wh L$, где $L\in \K(H)$,
приведены в следующей таблице.

$$
  \begin{array}{c|cccc}
       L                     & L_1   &  L_2 &  L_3   & L_4 \\
 \hline
       x_{\mbox{}_L}^{\vphantom{A^{A^A}}}         & 1    &  (12)(34) &  (123)   & (132)   \\
     \om_{1_H}(\wh L)=|L|^{\vphantom{A^{A^A}}} & 1   &  3  &  4    & 4    \\

    \l_b(\wh L)=|L|^{\st^{\vphantom{A^A}}}   & 1    & 1 & 0   & 0
  \end{array}
$$

\noindent
Найдём значения отображения $\l_b^G$ на классовых суммах $\wh K$, где $K\in \K(G)$.
$$
  \begin{array}{c|ccccc}
     K                     & K_1   &  K_2 &  K_3   & K_4  & K_5\\
   \hline
       x_{\mbox{}_K}^{\vphantom{A^{A^A}}}       & 1   &  (12)(34) &  (123)   & (12345) & (13524) \\
       \sum_{x\in K\cap H}x^{\vphantom{A^{A^A}}}         & \wh L_1    &  \wh L_2 & \wh L_3+\wh L_4   & 0 & 0   \\
     \l_b^G(\wh K)^{\vphantom{A^{A^A}}} & 1   &  1  &  0    & 0 & 0
  \end{array}
$$
\noindent
Также вычислим значения $\l_{\mbox{}_B}$, где $B$ --- главный $2$-блок группы $G$.
$$
  \begin{array}{c|ccccc}
     K                     & K_1   &  K_2 &  K_3   & K_4  & K_5\\
   \hline
   \om_{1_G}(\wh K)=|K|^{\vphantom{A^{A^{A^A}}}}  & 1    &  15 & 20   & 12 & 12   \\
     \l_{\mbox{}_B}(\wh K)=|K|^{\st^{\vphantom{A^A}}} & 1   &  1  &  0    & 0 & 0
  \end{array}
$$
\noindent
Как видно, значения отображений  $\l_b^G$ и $\l_{\mbox{}_B}$ совпадают на всех классовых суммах группы $G$.
Поэтому $\l_b^G=\l_{\mbox{}_B}$ и, значит, индуцированный блок $b^G$ определён и совпадает с $B$.
\end{prim}

Описание блоков с циклическими дефектными группами (иногда называемых \mem{циклическими блоками})
было завершено Дэйдом и является глубоким и полезным результатом.

\begin{thm}[Дэйда] \label{bcd}\glsadd{iThmDade} Пусть $B\in \bl(G)$ и $\d(B)$ --- циклическая группа.
Пусть $|\d(B)|=p^d>1$ и $|\iBr(B)|=e$. Тогда справедливы следующие утверждения.
\begin{list}{{\rm(}{\it\roman{enumi}\/}{\rm)}}{\usecounter{enumi}\setlength{\parsep}{2pt}\setlength{\topsep}{5pt}\setlength{\labelwidth}{23pt}}
\item $e\mid p-1$.
\item $|\irr(B)|=\displaystyle\frac{p^d-1}{e}+e$.
\item Если $d=1$ и $e=p-1$, то положим $\irr(B)=\{\x_0,\x_1,\ld,\x_e\}$. Если либо $d>1$, либо $e<p-1$,
то выполнены следующие утверждения.
      \begin{list}{{\rm(}{\it\roman{enumi}\/}.{\rm\arabic{enumii})}}
          {\usecounter{enumii}\setlength{\parsep}{2pt}\setlength{\topsep}{2pt}\setlength{\labelwidth}{5pt}\setlength{\labelwidth}{50pt}}
\item Множество $\irr(B)$ однозначно разбивается на два непересекающихся подмножества
$$\{\t_1,\t_2,\ld,\t_t\}, \qquad \{\x_1,\x_2,\ld,\x_e\},$$
где $t=\displaystyle\frac{p^d-1}{e}$. В этом случае
характеры из первого подмножества называются \mem{исключительными}. Любой (фиксированный) исключительный характер
обозначим через $\x_0$.
    \item $\hat\t_1=\hat\t_2=\ld =\hat\t_t$. В частности, $d_{\t_i,\vf}=d_{\x_0,\vf}$ для
          любого $\vf\in \iBr(B)$ и любого $i=1,\ld,t$.
    \end{list}
\item Для любого $\vf\in \iBr(B)$ ровно два из чисел разложения
$d_{\x_0,\vf}$, $d_{\x_1,\vf}$, \ld, $d_{\x_e,\vf}$ равны $1$, а остальные равны $0$. В частности, все числа
разложения блока $B$ равны $0$ или $1$.
\item Все характеры из $\irr(B)$ имеют высоту $0$.
\end{list}
\end{thm}


Теорема \ref{bcd} позволяет построить \glsadd{iBrTr}\mem{дерево Брауэра}, которое описывает строение
блока $B$ с циклической дефектной группой. Оно имеет $e+1$ вершин, взаимно однозначно соответствующих
характерам $\x_0$, $\x_1$, \ld, $\x_e$. Причём в случае, когда либо  $d>1$, либо $e<p-1$, вершина,
соответствующая характеру $\x_0$ называется \mem{исключительной} и обычно выделяется особо. Количество рёбер
равно $e$, и рёбра взаимно однозначно соответствуют характерам из $\iBr(B)$. Две вершины, соответствующие
$\x_i$ и $\x_j$, соединены ребром тогда и только тогда, когда существует характер $\vf\in \iBr(B)$, для
которого $d_{\x_i,\vf}=d_{\x_j,\vf}=1$. Из свойств чисел разложения и пункта $(iv)$ теоремы
\ref{bcd} легко следует, что полученный граф связен и является деревом.

\begin{figure}[htb]
\begin{center}
\setlength{\unitlength}{0.8mm}
\begin{picture}(100,35) (-25,0)
 \put (5,15){\circle*{1.2}}\put (5,15){\line(1,0){15}}
 \put (20,15){\circle*{1.2}}\put (20,15){\line(1,0){13.5}}
 \put (35,15){\circle*{1.2}}\put (50,15){\line(-1,0){13.5}}
  \put (35,15){\circle{3.0}}
 \put (50,15){\circle*{1.2}}\put (50,15){\line(1,0){15}}
 \put (65,15){\circle*{1.2}}
 \put (50,15){\line(0,-1){15}} \put (50,0){\circle*{1.2}}
 \put (50,15){\line(0,1){15}} \put (50,30){\circle*{1.2}}

 \put (6,17){$\scriptstyle 1$}
 \put (21,17){$\scriptstyle 675$}
 \put (37,17){$\scriptstyle 2048$}
 \put (51,17){$\scriptstyle 1728$}
 \put (66,17){$\scriptstyle 300$}
 \put (52,1){$\scriptstyle 27$}
 \put (52,32){$\scriptstyle 27$}
\end{picture}
\caption{Главный $13$-блок группы $\mbox{}^2F_4(2)'$}\label{t13}
\end{center}
\end{figure}

В качестве примера рассмотрим главный $13$-блок $B$ простой группы Титса $\mbox{}^2F_4(2)'$ порядка
$2^{11}\cdot 3^3\cdot 5^2\cdot 13$. В этом случае $|\iBr(B)|=6$, имеется два исключительных
характера, и дерево Брауэра выглядит, как показано на рис. \ref{t13}. Вершины дерева
помечены степенями соответствующих характеров, а исключительная вершина  обведена.

В заключении курса сформулируем несколько глубоких теоретико-групповых теорем, вытекающих из теории
обыкновенных и брауэровых характеров.

Напомним, что \mem{обобщённой группой кватернионов} \glsadd{iGrQuatGen}
называется группа
$$
\gls{Q2n}=\big\la x,y\bigm| x^{2^{m-1}}=1,\ \  y^2=x^{2^{m-2}}, \ \ x^y=x^{-1} \big\ra,
$$
где $m\ge 3$. Легко видеть, что $|Q_{2^m}|=2^m$ и $Q_8$ совпадает с обычной группой кватернионов.

Элемент порядка $2$ группы $G$ называется \glsadd{iInvl}\mem{инволюцией}.

\textuprn{Доказать, что обобщённая группа кватернионов содержит единственную инволюцию.}

Из следующего утверждения следует, что конечная группа не может быть простой в случае, когда
её $2$-силовская подгруппа является обобщённой группой кватернионов.

\begin{thm}[Брауэра\;\!-\!-Сузуки] \label{tbs}\glsadd{iThmBraSuz} Пусть $2$-силовская подгруппа группы $G$
является обобщённой группой кватернионов. Если $\OO_{2'}(G)=1$, то $|\Z(G)|=2$.
\end{thm}

Если $2$-силовская подгруппа $P$ группы $G$ является обобщённой группой кватернионов порядка, большего, чем $8$,
то теорема Брауэра-Сузуки может быть доказана с использованием только теории обыкновенных характеров. Наиболее сложная
часть теоремы --- случай, когда  $P\cong Q_8$ --- первоначально была доказана с помощью
теории блоков. Однако позже, в 1974 году, Глауберман показал, что и здесь можно обойтись
только результатами из теории обыкновенных характеров. Чисто теоретико-группового доказательства теоремы
Брауэра-Сузуки до сих пор не известно.

\begin{opr} \label{izi} Пусть $u\in G$ --- инволюция и $u\in P\in \Syl_2(G)$.
Инволюция $u$ называется \glsadd{iInvlIsl}\mem{изолированной}, если она не сопряжена в $G$
ни с одним элементом из $P\setminus\{u\}$.
\end{opr}

\textuprn{Доказать, что

$(i)$ изолированность инволюции $u$ не зависит от выбора $2$-силовской подгруппы $P$ в определении \ref{izi};

$(ii)$ изолированная инволюция лежит в центре любой содержащей её $2$-силовской подгруппы.}

Для группы $G$ обозначим через \gls{Zst} полный прообраз в $G$ центра факторгруппы $G/\OO_{2'}(G)$.

Следующий фундаментальный результат является одним из наиболее важных инструментов локального анализа.

\begin{thm}[$\Z^*$-теорема Глаубермана] \label{zsg}\glsadd{iThmZGla} Пусть $G$ обладает изолированной инволюцией $u$. Тогда
$u\in \Z^*(G)$.
\end{thm}

Отметим, что $\Z^*$-теорема является обобщением теоремы Брауэра-Сузуки, однако её доказательство
при рассмотрении некоторых частных случаев опирается на эту теорему. Также отметим, что на сегодняшний день
не известно доказательства $\Z^*$-теоремы, которое не  использовало бы теорию брауэровых характеров.

%
%

\backmatter

\setcounter{section}{0}
\renewcommand{\thesection}{\Asbuk{section}}
\renewcommand{\thesubsection}{(\mbox{\thesection}.\arabic{subsection}\unskip )}

\makeatletter  
\renewcommand\@seccntformat[1]{\csname the#1\endcsname.\quad}
\makeatother

\renewcommand{\chaptermark}[1]{%
\markboth{#1}{}}

\chapter{Приложения}

\section{Сведения из теории множеств\label{sv tmn}}
Нам потребуется лемма Цорна. Прежде чем её сформулировать, напомним некоторые определения из теории множеств.

Множество  $X$ называется \glsadd{iMnChU}\mem{частично упорядоченным}, если между любой парой его элементов $x$ и $y$
задано соотношение $x\le y$ такое, что для произвольных $x,y,z\in X$ справедливо
\begin{list}{{\rm(}{\it\roman{enumi}\/}{\rm)}}
{\usecounter{enumi}\setlength{\parsep}{2pt}\setlength{\topsep}{5pt}\setlength{\labelwidth}{23pt}}
\item $x\le x$;
\item если $x\le y$ и $y\le x$, то $x=y$;
\item если $x\le y$ и $y \le z$, то $x\le z$.
\end{list}
Подмножество $Y\se X$ называется \glsadd{iChn}\mem{цепью}, если для любых $x,y\in Y$ имеем либо $x\le y$, либо $y\le x$.
Подмножество $Z\se X$ \glsadd{iSbBnd}\mem{ограничено сверху}, если существует элемент $x\in X$ такой, что $z\le x$ для всех $z\in Z$.
Элемент $x\in X$ называется \glsadd{iElmMax}\mem{максимальным}, если для любого $y\in X$ из того, что $x\le y$ следует, что $x=y$.

\glsadd{iLemZorn}
\begin{thm}[леммa Цорна]\label{thm lz}
Любое\footnote{Для пустого множества лемма Цорна также верна, т.\,к. в нём пустая цепь не ограничена сверху, и значит, посылка утверждения ложна.}
частично упорядоченное множество, в котором всякая цепь ограничена сверху, содержит максимальный элемент.
\end{thm}
Известно, что лемма Цорна эквивалентна аксиоме выбора, принимаемой в теории множеств. Как правило мы будем использовать лемму Цорна применительно
к семейству подмножеств некоторого множества, частично упорядоченному относительно включения.

\section{Сведения из теории колец\label{cv kol}}

В этом приложении вводятся некоторые понятия и обозначения, связанные с кольцами.
Напомним, что термин <<кольцо>> всегда обозначает ассоциативное кольцо с единицей.

Пусть $R$ --- кольцо. Единицу и ноль этого кольца будем обозначать символами \gls{1R} и \gls{0R}, соответственно,
или просто $1$ и $0$, если из
контекста ясно, о каком кольце идёт речь.\footnote{Кольцо $R$, в котором $1_R= 0_R$, одноэлементно
и называется \glsadd{iRngZer}\mem{нулевым}.} Аддитивная подгруппа
 $S$ в $R$ называется \mem{подкольцом} \glsadd{iSubRng}
(обозначается $S\le R$), если $1_R\in S$ и для любых $a,b\in S$ выполнено $ab\in S$. В частности, $S$ является
кольцом и $1_S=1_R$. Примером подкольца в произвольном кольце $R$ является его \glsadd{iCenRng}\mem{центр}
$$\gls{ZlRr}=\{z\in R\bigm|  zr=rz\ \mbox{для всех}\ r\in R\}.$$
Аддитивная подгруппа $I$ в $R$ называется \glsadd{iIdlRngL}\glsadd{iIdlRngR}\mem{левым $($правым$\,)$ идеалом}
и обозначается $I\gls{IdRl}R$ ($I\gls{IdRr}R$), если для любых $a\in I$, $r\in R$ выполнено $ra\in I$ (соответственно $ar\in
I$). Подмножество $I\se R$, одновременно являющееся левым и правым идеалом, называется \mem{двусторонним
идеалом}\glsadd{iIdlRngTs}  или просто \mem{идеалом}\glsadd{iIdlRng} и обозначается $I\gls{IdR}R$.

\textuprn{Когда подкольцо является идеалом кольца? Когда центр является идеалом? Когда идеал является подкольцом?
Когда $\{0_R\}$ является подкольцом в $R$?}

\subsection{Примеры колец, подколец и идеалов}\label{pr kol}
\begin{list}{{\rm(}{\it\roman{enumi}\/}{\rm)}}
{\usecounter{enumi}\setlength{\parsep}{2pt}\setlength{\topsep}{5pt}\setlength{\labelwidth}{23pt}}

\item Кольцо целых чисел $\ZZ$. Идеалы кольца
$\ZZ$ исчерпываются множествами $n\ZZ$ целых чисел, кратных произвольному числу $n$.

\item Кольцо $\ZZ_n$ вычетов по модулю $n$. Идеалы $\ZZ_n$ исчерпываются
множествами вычетов по модулю $n$ элементов из $m\ZZ$ для произвольного делителя $m$ числа $n$.

\item Пусть $V$ --- абелева группа. Множество
\gls{EndV} её эндоморфизмов является кольцом относительно операций поэлементного сложения и (правой) композиции \big(т.\,е.
$v(\vf+\psi)=v\vf+v\psi$ и $v(\vf\psi)=(v\vf)\psi$ для произвольных $\vf,\psi\in \End(V)$, $v\in V$\big).
Если одновременно  $V$ является векторным пространством над полем $F$, то кольцо линейных преобразований
\gls{EndFV} этого пространства будет подкольцом\footnote{Аналогичные примеры подколец в $\End(V)$  можно
получить, если рассмотреть все эндоморфизмы, сохраняющие на $V$ какую-либо другую дополнительную структуру
(модуля, кольца, и т.\,п.).}
 в $\End(V)$.

\item $R[x_1,\ld,x_n]$ --- кольцо многочленов от $n$ переменных над произвольным коммутативным кольцом
$R$. Многочлены без свободного члена образуют идеал в $R[x_1,\ld,x_n]$.


\item \gls{MnR} --- кольцо матриц
размера\footnote{Кольцо $M_0(R)$ является нулевым и состоит из \glsadd{iMatEmp}\mem{пустой матрицы} размера $0\times 0$.} $n\times n$ над
коммутативным кольцом $R$. Центр $\Z(\MM_n(R))$
совпадает с множеством скалярных матриц с элементами  из $\Z(R)$. Если $R$ --- поле, то единственным ненулевым
идеалом кольца $\MM_n(R)$ является само это кольцо. Единицу этого кольца будем часто обозначать символом
\gls{In} или просто $\I$, а ноль --- символом
\gls{On} или просто $\O$.

\item Множество \gls{PX} всех подмножеств
некоторого множества $X$ является кольцом относительно операций симметрической разности и пересечения.
Напомним, что \mem{симметрической разностью} \glsadd{iAtrA} множеств $A$ и $B$ называется множество
$$\gls{AtrB}=(A\cup B)\setminus(A\cap B)= (A\setminus B)\cup(B\setminus A).$$
Пусть $A\subseteq X$. Тогда $\{\varnothing,X,A,X\!\setminus\! A\}\le \P(X)$ и $\P(A)\nor\P(X)$. Другим примером
идеала в $\P(X)$ служит множество всех конечных подмножеств множества $A$.
\end{list}

\textupr{Проверить сформулированные в \ref{pr kol} утверждения.}

\begin{pre}\label{max id}
Любой собственный правый $($левый, двусторонний$\,)$ идеал кольца $R$ содержится в
некотором максимальном правом $($левом, двустороннем$\,)$ идеале.
\end{pre}

\textupl{max id prf}{Доказать предложение \ref{max id}.
\uk{Использовать лемму Цорна \ref{thm lz} и тот факт, что $R$ содержит единицу.}}

Для любого кольца $R$ можно определить \glsadd{iRngOp}\mem{противоположное кольцо}
\gls{Rop}, аддитивная группа которого та же, что и у $R$, а произведение $r$ на $s$ равно $sr$ для любых $r,s\in
R^{op}$. Ясно, что  $(R^{op})^{op}=R$ для произвольного кольца $R$, и $R^{op}=R$ для коммутативного кольца $R$.
Кроме того, левые (правые) идеалы кольца $R$ являются правыми (левыми) идеалами кольца $R^{op}$.

\mysubsection\label{end op}
Если  в примере \ref{pr kol}$(iii)$ для кольца  $\End(V)$ использовать не правую, а левую запись эндоморфизмов группы $V$ и их композиций,
то получится кольцо $\End(V)^{op}$.

Легко видеть, что если $S$ --- подкольцо, а $I$ --- правый (левый, двусторонний) идеал некоторого кольца, то
$S\cap I$ --- правый (левый, двусторонний) идеал кольца $S$. Отметим также, что пересечение произвольного
семейства подколец (левых, правых, двусторонних идеалов) кольца снова является подкольцом (соответственно, левым, правым,
двусторонним идеалом).

Если $S$  --- подкольцо кольца $R$ и $X\se R$, то обозначим через \gls{Sx} наименьшее подкольцо в $R$, содержащее $S$ и $X$. Будем называть $S[X]$
\mem{подкольцом} кольца $R$, \glsadd{iSubRngGenX}\mem{порождённым} $S$ и $X$ или
\mem{подкольцом, полученным присоединением} к $S$ \mem{элементов} множества $X$. В случае, когда множество
$X=\{r_1,\ld,r_n\}$ конечно, будем для краткости писать $S[r_1,\ld,r_n]$ вместо $S[\{r_1,\ld,r_n\}]$. Легко
показать, что $S[r_1,\ld,r_n]=S[r_1,\ld,r_{n-1}][r_n]$.

 Аналогично, если $F\ge K$ --- расширение полей и $X\se F$, то \gls{Kx} --- наименьшее
подполе поля $F$, содержащее $K$ и $X$. Будем называть $K(X)$ \mem{подполем} поля $F$, \mem{порождённым} $K$
и $X$ или \glsadd{iSubFldGenX}  \mem{подполем} поля $F$, \mem{полученным присоединением} к $K$ \mem{элементов} множества $X$.
Очевидно, что выполнено включение подколец $K[X]\se K(X)$, которое, вообще говоря, является строгим.
Для $a_1,\ld,a_n \in F$ обозначим $K(\{a_1,\ld,a_n\})=K(a_1,\ld,a_n)$. Как и в
случае колец, легко проверить, что $K(a_1,\ld,a_n)=K(a_1,\ld,a_{n-1})(a_n)$.

\textuprn{Проверить следующие утверждения.
\begin{list}{{\rm(}{\it\roman{enumi}\/}{\rm)}}
{\usecounter{enumi}\setlength{\parsep}{2pt}\setlength{\topsep}{5pt}\setlength{\labelwidth}{23pt}}
\item Пусть $R$ --- коммутативное кольцо, $S\le R$ и $r_1,\ld,r_n\in R$. Тогда подкольцо
$S[r_1,\ld,r_n]$ состоит из всевозможных значений многочленов $f(r_1,\ld,r_n)$, где $f$ пробегает
$S[x_1,\ld,x_n]$.
\item  Пусть $F\ge K$ --- расширение полей  и $a_1,\ld,a_n \in F$. Тогда подполе $K(a_1,\ld,a_n)$ состоит из значений выражения
$f(a_1,\ld,a_n)/g(a_1,\ld,a_n)$, где $f,g\in K[x_1,\ld,x_n]$ и $g(a_1,\ld,a_n)\ne 0$.
\end{list}}

\textuprn{\label{subrq}Описать подкольца поля $\QQ$. \uk{
Рассмотреть кольца
$\gls{Zpim1}=\ZZ[\,\frac{1}{p}\mid p\in\pi\,]$, где $\pi$ --- некоторое
подмножество множества всех простых целых чисел.}}

\begin{opr}\label{s p ids}
Пусть $X$ --- подмножество кольца $R$. Правым (левым, двусторонним) идеалом из $R$,
\glsadd{iIdlGenX}\mem{порождённым множеством} $X$, называется наименьший правый (левый,
двусторонний) идеал из $R$, содержащий $X$.\footnote{В частности, из определения вытекает, что правый (левый,
двусторонний) идеал, порождённый пустым множеством, нулевой.} Легко видеть, что такой идеал совпадает с
множеством всевозможных конечных сумм\footnote{Здесь и далее мы естественно предполагаем, что сумма нулевого
числа слагаемых равна нулю.} вида $\sum x_ir_i$ в случае правого идеала,  вида $\sum r_ix_i$ в случае левого
идеала, и вида $\sum r_ix_is_i$ в случае двустороннего идеала, где  $r_i,s_i\in R$, $x_i\in X$.

Если $X=\{x_1,\ld,x_n\}$ --- конечное множество, то правый (левый) идеал, порождённый $X$, обозначается через
$x_1R+\ld+x_nR$ (соответственно, $Rx_1+\ld+Rx_n$), а двусторонний идеал, порождённый $X$, --- через
$Rx_1R+\ld+Rx_nR$. Здесь выражение $RxR$ обозначает множество конечных сумм вида $\sum r_ixs_i$, где
$r_i,s_i\in R$.

Пусть $I_1,\ld,I_n$ --
правые (левые, двусторонние) идеалы кольца $R$. \glsadd{iSumIds}\mem{Суммой идеалов} $I_1,\ld,I_n$ называется множество $$
I_1+\ld+I_n= \{\, r_1+\ld +r_n\bigm| r_i\in I_i,\ i=1,\ld,n
\},$$
а \glsadd{iProdIds}\mem{произведением} $I_1\cdot \ld \cdot I_n$ называется множество всевозможных
конечных сумм произведений $r_1\cdot\ld\cdot r_n$, где $r_i\in I_i$, $i=1,\ld,n$.
\end{opr}

Легко видеть, что сумма, произведение, а также пересечение правых (левых, двусторонних)
идеалов снова будет правым (левым, двусторонним) идеалом. В частности, для правого (левого,
двустороннего) идеала $I$ из $R$ определены степени $I^n$ для произвольного $n>0$. Если же $n=0$,
то положим по определению $I^0=R$.

Заметим, что для произвольных двусторонних идеалов $I_1,\ld, I_n$ из $R$ имеет место включение
$$
I_1\cdot\ld\cdot I_n\le I_1\cap\ld\cap I_n.
$$

Правый (левый) идеал $I\nor_r R$ называется \glsadd{iIdlNlp}\mem{нильпотентным}, если
$I^n=0$ для некоторого $n$. Отметим, что $I^n=0$ тогда и только тогда, когда любое произведение $n$
элементов из $I$ равно нулю.

Элемент $r\in R$ называется \glsadd{iInvElm}\mem{обратимым справа},
если у него есть \glsadd{iObrElm}\mem{правый обратный}, т.\,е. такой элемент $s\in R$, что $rs=1_R$. \mem{Обратимый
слева} элемент определяется аналогичным образом. Обратимый справа и слева элемент называется просто
\mem{обратимым}.

\begin{pre} \label{pro} Собственный
(правый, левый) идеал кольца не содержит обратимых (справа, слева) элементов и, в частности, единицу кольца.
\end{pre}
\textupl{pro prf}{Доказать предложение \ref{pro}.}

\textzam{Существуют примеры колец, в которых собственный правый идеал содержит обратимый слева элемент.}





\textuprn{Доказать, что у обратимого элемента $r\in R$ любой правый обратный совпадает с любым левым обратным.
Этот однозначно определённый элемент называется просто \mem{обратным} к $r$ и обозначается $r^{-1}$.}

Для кольца $R$ положим
$$
\gls{Rx} = \{r\in R \mid r\ \mbox{обратим}\}.
$$

\begin{pre} \label{zr} Пусть $R$ --- кольцо. Тогда
\begin{list}{{\rm(}{\it\roman{enumi}\/}{\rm)}}
{\usecounter{enumi}\setlength{\parsep}{2pt}\setlength{\topsep}{5pt}\setlength{\labelwidth}{23pt}}
\item $R^\times$ --- группа;
\item $\Z(R)^\times=\Z(R)\cap R^\times$;
\item $\Z(R)^\times\le \Z(R^\times)$.
\end{list}
\end{pre}

\textupl{zr prf}{Доказать предложение \ref{zr}.}

Множество $R^\times$, являющееся группой ввиду \ref{zr}$(i)$,
называется \glsadd{iGrInvElms}\mem{группой обратимых элементов} или
\glsadd{iGrMltRng}\mem{мультипликативной группой} кольца~$R$.

В связи с \ref{zr}$(iii)$ отметим, что для произвольного кольца $R$
включение $\Z(R^\times)\le \Z(R)^\times$, вообще говоря, не имеет места. См., однако, \ref{zri}.


\subsection{Примеры}
\begin{list}{{\rm(}{\it\roman{enumi}\/}{\rm)}}
{\usecounter{enumi}\setlength{\parsep}{2pt}\setlength{\topsep}{5pt}\setlength{\labelwidth}{23pt}}
\item Имеем $\ZZ^\times =\{1,-1\}$.
\item Группа $\ZZ_n^\times$ состоит из вычетов по модулю $n$ целых чисел, взаимно простых с $n$. Она имеет порядок $\phi(n)$,
где \gls{phi} --- \glsadd{iFuncEul}\mem{функция Эйлера}. Если $n=p_1^{a_1}\ldots p_k^{a_k}$ --- разложение числа $n$ на простые множители,
то
$$
\phi(n)=(p_1-1)p_1^{a_1-1}\ldots (p_k-1)p_k^{a_k-1}.
$$
\item Группа $\MM_n(R)^\times$ обратимых матриц над коммутативным кольцом $R$ обозначается через
\gls{GLnR} и называется \glsadd{iGLnR}\mem{полной линейной группой степени $n$ над $R$}. Для матрицы $A\in\MM_n(R)$ включение $A\in\GL_n(R)$
имеет место тогда и только тогда, когда $\det A\in R^\times$.
\end{list}

Кольцо $R$ называется \glsadd{iSqrng} \mem{телом}, если $R^\times=R\setminus\{0\}$.\footnote{В частности, нулевое кольцо не является телом.}

\begin{pre}\label{obr}\mbox{}
\begin{list}{{\rm(}{\it\roman{enumi}\/}{\rm)}}
{\usecounter{enumi}\setlength{\parsep}{2pt}\setlength{\topsep}{5pt}\setlength{\labelwidth}{23pt}}
\item Если все ненулевые элементы ненулевого кольца $R$ обратимы справа (слева), то $R$ является телом.
\item Если $I\nor_l R$ и в $R\setminus I$ любой элемент обратим справа,
то в $R\setminus I$ любой элемент также обратим и слева.
\item Если $I\nor_r R$
и для любого $r\in I$ элемент $1-r$ обратим справа, то для любого $r\in I$ элемент
$1-r$ также обратим и слева.
\end{list}
\end{pre}

\textupl{obr prf}{Доказать предложение \ref{obr}.}

\textupln{mpl prf}{\label{mpl}Показать, что если идеал кольца является максимальным как правый идеал, то он
также является максимальным как левый идеал, и факторкольцо по нему является телом.}

Элемент $r\in R$ называется \glsadd{iNlpElm}\mem{нильпотентным}, если $r^n=0$ для
некоторого $n$. Из сказанного выше вытекает, что любой правый (левый) нильпотентный идеал кольца $R$
состоит из нильпотентных элементов. Заметим также, что для нильпотентного элемента $r$ элемент $1-r$ обратим,
поскольку для подходящего $n$

$$
(1-r)(1+r+\ld+r^{n-1})=(1+r+\ld+r^{n-1})(1-r)=1.
$$

\textuprn{Найти нильпотентные элементы кольца $\ZZ_n$.}

Пусть $R,S$ -- кольца. Гомоморфизм $\vf:R\to S$ аддитивных групп этих колец будем называть \glsadd{iHomRng}\mem{гомоморфизмом колец},
если он сохраняет умножение, т.\,е.\footnote{при правой записи отображений}
$(xy)\vf=(x\vf)(y\vf)$ для любых $x,y\in R$, а также сохраняет
единицу: $(1_R)\vf=1_S$.

\textuprn{\label{deg homs}Из нулевого кольца в ненулевое нет гомоморфизмов. Также нет гомоморфизмов, например, из $\ZZ_2$ в $\ZZ_3$.
Из произвольного кольца в нулевое существует единственный гомоморфизм. Также существует единственный гомоморфизм
из кольца $\ZZ$ в произвольное кольцо.}

Отметим, что для гомоморфизма колец $\vf:R\to S$ \glsadd{iKerRng}\mem{ядро}, обозначаемое\footnote{Напомним,
что для ядер групповых гомоморфизмов, в отличие от кольцевых,
мы, как правило, используем обозначение <<$\ker$>>, а не <<$\Ker$>>.} \gls{Ker}, является идеалом в $R$, а
образ $\Im\vf$ --- подкольцом в $S$. Сюръективные гомоморфизмы называются
\glsadd{iEpiRng}\mem{эпиморфизмами}, инъективные~--- \glsadd{iMonoRng}\mem{мономорфизмами}, а биективные \glsadd{iIsoRng}\mem{изоморфизмами}. Легко проверить, что
гомоморфизм $\vf$ инъективен тогда и только тогда, когда $\Ker\vf=\{0\}$. Любой изоморфизм $R\to R$ называется
\glsadd{iAutRng}\mem{автоморфизмом кольца} $R$.


\subsection{Примеры}\label{pr hr}
\begin{list}{{\rm(}{\it\roman{enumi}\/}{\rm)}}
{\usecounter{enumi}\setlength{\parsep}{2pt}\setlength{\topsep}{5pt}\setlength{\labelwidth}{23pt}}
\item Если $m$ и $n$ --- натуральные числа такие, что $m$ делит $n$, то существует естественный эпиморфизм колец
$\ZZ_n\to\ZZ_m$, действующий по правилу $a+n\ZZ\,\mapsto\, a+m\ZZ$ для любого $a\in\ZZ$. Его ограничение на
группу обратимых элементов будет эпиморфизмом\footnote{Отметим, что всякий гомоморфизм колец $\vf:R\to S$ индуцирует
гомоморфизм мультипликативных групп $\vf^\times:R^\times\to S^\times$. Однако, $\vf^\times$ может не быть
сюръективным даже если $\vf$ сюръективен.
}
групп $\ZZ_n^\times\to\ZZ_m^\times$.
\item Если $R$ --- коммутативное кольцо и $r\in R$, то отображение $f(x)\mapsto f(r)$ будет эпиморфизмом
колец $R[x]\to R$.
\item Пусть $X$ --- множество и $Y\se X$. Тогда отображение $A\mapsto A\cap Y$ будет эпиморфизмом
колец $\P(X)\to\P(Y)$.
\item Если $r\mapsto\wt{r}$ --- гомоморфизм (эпиморфизм, мономорфизм) коммутативных колец $ R\to S$, то отображения
$$(a_{ij})\mapsto (\wt{a}_{ij}) \qquad \text{и} \qquad\textstyle\sum a_ix^i\mapsto \sum \wt{a}_ix^i$$
будут гомоморфизмами (эпиморфизмами, мономорфизмами) колец $\MM_n(R)\to \MM_n(S)$ и $R[x]\to S[x]$, соответственно.
\end{list}

\begin{pre} \label{hom r} Пусть $R,S,T$ --- кольца.
\begin{list}{{\rm(}{\it\roman{enumi}\/}{\rm)}}
{\usecounter{enumi}\setlength{\parsep}{2pt}\setlength{\topsep}{5pt}\setlength{\labelwidth}{23pt}}
\item Если $\vf: R\to S$ и $\psi: R\to T$ --- гомоморфизмы, причём $\vf$ --- эпиморфизм и $\Ker \vf\se \Ker \psi$,
то существует единственный гомоморфизм $\t: S\to T$ такой, что $\psi=\vf\t$.
$$
\xymatrix{
R\ar[r]^{ \textstyle{\vf}}\ar[d]_{ \textstyle{\psi}}& S\ar@{.>}[dl]^{ \textstyle{\t} }        \\
T       &}
$$
\item Если $\vf: S\to R$ и $\psi: T\to R$ --- гомоморфизмы, причём $\vf$ --- мономорфизм и $\Im \vf\supseteq \Im
\psi$, то существует единственный гомоморфизм $\t: T\to S$ такой, что $\psi=\t\vf$.
$$
\xymatrix{
R& S\ar[l]_{ \textstyle{\vf}}         \\
T\ar@{.>}[ur]_{ \textstyle{\t} }\ar[u]^{ \textstyle{\psi}}  &}
$$
\end{list}
\end{pre}
\vspace{-10pt}
\textupl{hom r prf}{Доказать предложение \ref{hom r}.}

Для колец $R,S$ гомоморфизм (изоморфизм) $\vf:R\to S$ их аддитивных групп называется \glsadd{iAntHomRng}\mem{антигомоморфизмом}
\glsadd{iAntIsoRng}$($\mem{антиизоморфизмом}$\,)$ колец, если $(1_R)\vf=1_S$ и $(xy)\vf=(y\vf)(x\vf)$ для любых $x,y\in R$.
Антиизоморфизм $R\to R$ называется \glsadd{iAntAutRng}\mem{антиавтоморфизмом} кольца $R$.

Например, антиавтоморфизмом матричного кольца $\MM_n(R)$ является отображение $A\to A^\top$ для всех $A\in
\MM_n(R)$.

\textuprn{Проверить следующие утверждения.
\begin{list}{{\rm(}{\it\roman{enumi}\/}{\rm)}}
{\usecounter{enumi}\setlength{\parsep}{2pt}\setlength{\topsep}{5pt}\setlength{\labelwidth}{23pt}}
\item Композиция двух антигомоморфизмов колец является гомоморфизмом, а композиция гомоморфизма и
антигомоморфизма (в произвольном порядке) является антигомоморфизмом.
\item Гомоморфизм (антигомоморфизм) $\vf:R\to S$ осуществляет антигомоморфизмы (гомоморфизмы)
$R\to S^{op}$ и $R^{op}\to S$.
\item Отображение, обратное к  антиизоморфизму, является антиизоморфизмом.

\item Тождественное отображение $R\to R$ является антиавтоморфизмом тогда и только тогда, когда кольцо $R$
коммутативно.
\end{list}
}

Для двустороннего идеала $I\nor R$ на факторгруппе $R/I$ вводится умножение $(x+I)(y+I)=xy+I$.
Корректность этого умножения и аксиомы кольца для $R/I$ легко проверяются. Кольцо $R/I$ называется
\mem{факторкольцом} кольца $R$ по идеалу $I$. Отображение $x\mapsto x+I$  из $R$ в $R/I$ является эпиморфизмом колец,
ядро которого совпадает с $I$. Таким образом, ядра кольцевых гомоморфизмов --- это, в точности, двусторонние идеалы.

Существует конструкция, позволяющая строить новые кольца из уже имеющихся.

\begin{opr}\label{opr dir kol} Пусть $R_1,\ld, R_n$ --- кольца.
Их \glsadd{iSumDirRng}\mem{прямой суммой} $R_1\oplus\ld\oplus R_n$ называется множество
$$
\{ (r_1,\ld,r_n)\bigm|  r_i\in R_i,\ \ i=1,\ld,n \}
$$
с покомпонентными операциями сложения и умножения элементов. Ясно, что прямая сумма колец
будет кольцом.
\end{opr}

\textuprn{Пусть $R_1,R_2,R_3$ --- кольца.
Показать, что
\begin{list}{{\rm(}{\it\roman{enumi}\/}{\rm)}}
{\usecounter{enumi}\setlength{\parsep}{2pt}\setlength{\topsep}{5pt}\setlength{\labelwidth}{23pt}}
\item  $R_1 \oplus R_2 \cong R_2 \oplus R_1$;
\item $(R_1 \oplus R_2)\oplus R_3\cong R_1 \oplus R_2\oplus R_3\cong R_1 \oplus (R_2\oplus R_3)$;
\item $(R_1 \oplus R_2)^\times \cong R_1^\times\times R_2^\times$.
\end{list}

Показать также, что если кольца $R_1$ и $R_2$ коммутативны и $R=R_1\oplus R_2$, то

\begin{list}{{\rm(}{\it\roman{enumi}\/}{\rm)}}
{\usecounter{enumi}\setlength{\parsep}{2pt}\setlength{\topsep}{5pt}\setlength{\labelwidth}{23pt}}
\addtocounter{enumi}{3}
\item $\MM_n(R)\cong\MM_n(R_1)\oplus\MM_n(R_2)$;
\item $R[x]\cong R_1[x]\oplus R_2[x]$.
\end{list}
}

Положим $R=R_1\oplus\ld\oplus R_n$. Для каждого $i=1,\ld,n$ рассмотрим множество $I_i$ элементов из $R$ вида
$(0,\ld,0,r,0,\ld,0)$, где все компоненты с номерами, отличными от $i$, нулевые, а $r$ пробегает все элементы
кольца $R_i$. Легко видеть, что $I_i\nor R$. Кроме того, каждый идеал $I_i$ является кольцом (но, вообще
говоря, не подкольцом в $R$), изоморфным $R_i$. Поэтому в дальнейшем мы часто будем отождествлять $R_i$ с
$I_i$. Отметим, что кольцо $R$ совпадает с суммой $I_1+\ld+I_n$. Кроме того, для любого $i=1,\ld,n$ пересечение
$I_i$ с суммой $\sum\limits_{j\ne i} I_j$ тривиально. Последнее свойство позволяет наряду с данным выше
определением <<внешней>> прямой суммы колец ввести понятие <<внутренней>> прямой суммы идеалов кольца.

Пусть $I_1,\ld,I_n$ -- правые (левые, двусторонние) идеалы кольца $R$. Сумма $I_1+\ld +I_n$ называется
\mem{прямой}, \glsadd{iSumDirIds} если для любого $i=1,\ld,n$ пересечение идеала $I_i$ с суммой
$\sum\limits_{j\ne i} I_j$ тривиально.


\begin{pre} \label{inn dir sum}
Пусть $R$ --- кольцо, совпадающее с прямой суммой
своих правых (левых, двусторонних) идеалов $I_1,\ld,I_n$. Тогда справедливы следующие утверждения.
\begin{list}{{\rm(}{\it\roman{enumi}\/}{\rm)}}
{\usecounter{enumi}\setlength{\parsep}{2pt}\setlength{\topsep}{5pt}\setlength{\labelwidth}{23pt}}
\item Каждый элемент $r\in R$ однозначно представляется в виде суммы $r=r_1+\ld+r_n$, где $r_i\in
I_i$. В частности, $1_R=e_1+\ld+e_n$ для некоторых $e_i\in I_i$, $i=1,\ld, n$.
\item Для каждого $i,j=1,\ld, n$ имеем $e_i^2=e_i$ и $e_ie_j=0$ при $i\ne j$.

\item Каждый идеал $I_i$,  $i=1,\ld, n$, совпадает с идеалом $e_i R$ (соответственно,  $Re_i$, $Re_iR$).
\end{list}
Если все идеалы  $I_1,\ld,I_n$ двусторонние, то дополнительно справедливы следующие утверждения.
\begin{list}{{\rm(}{\it\roman{enumi}\/}{\rm)}}
{\usecounter{enumi}\setlength{\parsep}{2pt}\setlength{\topsep}{5pt}\setlength{\labelwidth}{23pt}}
\addtocounter{enumi}{3}
\item $e_1,\ld, e_n\in\Z(R)$.
\item Каждый идеал $I_i$ является кольцом, а $e_i$ --- его единицей.
\item Имеет место изоморфизм колец $R\cong I_1\oplus\ld\oplus I_n$.
\end{list}
\end{pre}
\textupr{Доказать предложение \ref{inn dir sum}.}

Прямую сумму правых (левых, двусторонних) идеалов $I_1,\ld,I_n$ будем обозначать $I_1\oplus\ld\oplus I_n$. Такая запись не
вызовет недоразумений: в случае двусторонних идеалов понятия прямой суммы колец и прямой суммы идеалов по существу эквивалентны
в силу предложения \ref{inn dir sum}$(v)$, а в случае односторонних идеалов эта запись согласуется с вводимым далее
определением прямой суммы $R$-модулей, см. \ref{prsm}.

Элемент $e$ кольца $R$ называется \glsadd{iIdmp}\mem{идемпотентом}, если $e^2=e$. Идемпотент называется
\glsadd{iIdmpCen}\mem{центральным}, если он лежит в центре $\Z(R)$. Два идемпотента $e$ и $f$
называются \glsadd{iIdmpsOrt}\mem{ортогональными}, если $ef=fe=0$. Легко видеть, что сумма
двух ортогональных идемпотентов снова будет идемпотентом. Идемпотент кольца называется
\glsadd{iIdmpPrim}\mem{примитивным}, если он не представим в виде суммы двух ненулевых
ортогональных идемпотентов.\footnote{В частности, нулевой идемпотент примитивен.}

Из предложения \ref{inn dir sum} видно, что разложение кольца $R$ в прямую сумму идеалов $I_1\oplus\ld\oplus
I_n$ однозначно определяет разложение единицы $1_R$ в сумму ортогональных идемпотентов $e_1+\ld+e_n$, которые
являются центральными в случае, если идеалы $I_1,\ld,I_n$ двусторонние. С другой стороны, как показывают
следующие два предложения, наличие кольце $R$ (центрального) идемпотента, отличного от нуля и единицы,
определяет разложение $R$ в прямую сумму собственных (двусторонних) идеалов.

\begin{pre}\label{ids} Если $R$ --- кольцо и
$e_1,\ld,e_n\in R$ --- центральные попарно ортогональные идемпотенты, то справедливы следующие утверждения.
\begin{list}{{\rm(}{\it\roman{enumi}\/}{\rm)}}
{\usecounter{enumi}\setlength{\parsep}{2pt}\setlength{\topsep}{5pt}\setlength{\labelwidth}{23pt}}
\item $e_i R\nor R$, $i=1,\ld,n$.
\item $e_1+\ld+e_n$ --- центральный идемпотент.
\item Если $e_1+\ld+e_n=1_R$, то $R=e_1R\oplus\ld\oplus e_nR$.
\item $\Z(e_iR)=e_i\Z(R)$, $i=1,\ld,n$.
\end{list}
\end{pre}
\textupr{Доказать предложение \ref{ids}.}


\begin{pre}\label{ids 2} Если $R$ --- кольцо и $e\in R$ --- идемпотент, то
\begin{list}{{\rm(}{\it\roman{enumi}\/}{\rm)}}
{\usecounter{enumi}\setlength{\parsep}{2pt}\setlength{\topsep}{5pt}\setlength{\labelwidth}{23pt}}
\item элемент $1-e$ также является идемпотентом и $\,e(1-e)=(1-e)e=0$,
\item $R=eR\oplus(1-e)R$,
\item подмножество $eRe=\{ere\bigm|r\in R\}$ является кольцом, a идемпотент $e$ --- единицей этого кольца.
\end{list}
\end{pre}

\textupr{Доказать предложение \ref{ids 2}.}


\begin{pre} \label{subdir} Пусть $R$ --- кольцо. Если $R=I_1\oplus\ld\oplus I_n$, где $I_i\nor_l R$  (соответственно, $I_i\nor_r R$), $i=1,\ld,n$, 
и $I\nor_r R$ (соответственно, $I\nor_l R$), то $I=(I\cap I_1)\oplus\ld\oplus(I\cap I_n)$.

%
\end{pre}

\textupl{subdir prf}{Доказать предложение \ref{subdir}.}


\begin{pre} \label{z pr} Пусть $R_1,\ld, R_n$ --- кольца
и $R=R_1\oplus\ld\oplus R_n$. Тогда
\begin{list}{{\rm(}{\it\roman{enumi}\/}{\rm)}}
{\usecounter{enumi}\setlength{\parsep}{2pt}\setlength{\topsep}{5pt}\setlength{\labelwidth}{23pt}}
\item $\Z(R_i)\nor\Z(R)$ для всех $i=1,\ld,n$.
\item $\Z(R)=\Z(R_1)\oplus\ld\oplus\Z(R_n)$;
\end{list}
\end{pre}

\textupr{Доказать предложение \ref{z pr}}

\begin{pre} \label{dec uniq} Пусть
$R$ --- кольцо и $$R=I_1\oplus\ld\oplus I_m=J_1\oplus\ld\oplus J_n,$$ где  $I_i,J_j$ --- ненулевые двусторонние
идеалы из $R$, не представимые в виде прямой суммы меньших ненулевых двусторонних идеалов. Тогда $m=n$ и
существует подстановка $\s$  на множестве индексов $\{1,\ld, m\}$ такая, что $I_i=J_{i\s}$ для всех
$i=1,\ld,m$.
\end{pre}
\textupl{dec uniq prf}{Доказать предложение \ref{dec uniq}.}

\section{Элементы теории полей} \label{polya}

В этом приложении мы приведём факты из теории полей, используемые в основном тексте.

\begin{pre}\label{cyc fin} Всякая конечная подгруппа мультипликативной группы поля является циклической.
\end{pre}
\begin{proof} Пусть $F$ поле и $A\le F^\times$ --- конечная подгруппа. Обозначим через $n$ наибольший порядок элементов группы $A$, и
пусть $a\in A$ --- элемент порядка $n$. Тогда
порядок любого элемента из $A$ делит $n$. $\big($Это следует из того, что если в конечной абелевой группе есть элементы порядков $n_1$ и $n_2$, то
в ней также есть элемент порядка, равного наименьшему общему кратному $n_1$ и $n_2$.$\big)$
Поэтому все элементы из $A$, включая $n$ степеней элемента $a$, являются корнями многочлена $x^n-1$ над $F$. Но многочлен над полем не может иметь
корней больше, чем его степень. Значит, степенями $a$ исчерпываются все элементы группы $A$, т.\,е. $A$ циклическая.
\end{proof}

Пусть $F$ --- подполе поля $E$. В этом случае будем говорить, что $E$ --- \mem{расширение поля} $F$ или,
что $E \ge F$~--- \glsadd{iExtFld}\mem{расширение полей}.

Пусть $E\ge F$ --- расширение полей. Элемент $a\in E$ называется \glsadd{iElmAlg}\mem{алгебраическим} над $F$, если он
является корнем некоторого многочлена
$$
a_0x^n+ a_1 x^{n-1}+ a_2 x^{n-2} +\ld +a_n,
$$
где $a_0,a_1,\ld,a_n\in F$ и $a_0\ne 0$. Расширение $E\ge F$ называется \glsadd{iExtFldAlg}\mem{алгебраическим}, если
всякий элемент из $E$ алгебраический над $F$.

Поле $F$ называется \glsadd{iFldAlgCl}\mem{алгебраически замкнутым}, если любой
многочлен $f\in F[x]$ такой, что $\deg f\ge 1$, имеет в $F$ хотя бы один корень.
\glsadd{iClAlgFld}\mem{Алгебраическим замыканием} $\ov{F}$ поля $F$ называется его алгебраически замкнутое алгебраическое расширение.

\begin{pre} \label{alg pro}\mbox{}
\begin{list}{{\rm(}{\it\roman{enumi}\/}{\rm)}}
{\usecounter{enumi}\setlength{\parsep}{2pt}\setlength{\topsep}{5pt}\setlength{\labelwidth}{23pt}}
\item У каждого поля $F$ существует алгебраическое замыкание, и между любыми двумя алгебраическими замыканиями поля $F$ существует изоморфизм,
действующий тождественно на $F$.
\item Пусть $E\ge F$ --- алгебраическое расширение полей. Тогда
существует вложение поля $E$ в алгебраическое замыкание $\ov F$ тождественное на $F$, причём
всякий автоморфизм поля $E$, централизующий $F$, продолжается до автоморфизма $\ov F$.
\end{list}
\end{pre}
\begin{proof} $(i)$ См. \cite[\S 72]{w}.

$(ii)$ См. \cite[гл. VII, \S 2, Теорема 2]{l}.
\end{proof}

Расширение полей $E \ge F$ называется \glsadd{iExtFldFin}\mem{конечным},
если размерность поля $E$, рассматриваемого как векторное пространство над $F$, является конечным.
Эта размерность называется \glsadd{iDegExt}\mem{степенью расширения} $E\ge F$.

\textuprn{\label{fin alg}Доказать следующие утверждения.
\begin{list}{{\rm(}{\it\roman{enumi}\/}{\rm)}}
{\usecounter{enumi}\setlength{\parsep}{2pt}\setlength{\topsep}{5pt}\setlength{\labelwidth}{23pt}}
\item Всякое конечное расширение полей является алгебраическим.
\item Если элемент $a$ алгебраический над полем $F$, то расширение $F(a)\ge F$ конечно.
\end{list}}

\begin{pre} \label{kon algf} Расширение полей $E\ge F$ конечно тогда и только тогда, когда
$E=F(a_1,\ldots,a_n)$ для некоторых алгебраических над $F$ элементов $a_1,\ldots,a_n\in E$ и некоторого $n\ge 0$.
\end{pre}
\begin{proof} Пусть $F\leqslant E$ конечно и $a_1,\ldots,a_n$ --- базис $E$ над $F$. Тогда
все элементы $a_i$ алгебраические над $F$ ввиду \ref{fin alg}$(i)$ и, очевидно, $E=F(a_1,\ldots,a_n)$. Обратно, пусть $E=F(a_1,\ldots,a_n)$
и все $a_i$ алгебраические над $F$. По индукции поле $K=F(a_1,\ldots,a_{n-1})$ --- конечное расширение поля $F$ и $E=K(a_n)$.
Заметим, что степень $E$ над $F$ равна произведению степеней расширений $E\ge K$ и $K\ge F$, а расширение $E\ge K$ конечно в силу \ref{fin alg}$(ii)$,
поскольку элемент $a_n$ алгебраический над $F$ и, значит, над $K$ тоже. Отсюда следует требуемое.
\end{proof}

\begin{pre} \label{alg trans} Пусть
$K\ge E\ge F$ --- цепочка расширений полей. Если $K$ алгебраическое над $E$ и $E$ алгебраическое над $F$, то $K$
алгебраическое над $F$.
\end{pre}
\begin{proof} Пусть $a\in K$ --- корень многочлена $a_0x^n+\ld +a_n\in E[x]$. Рассмотрим поля $L=F(a_0,\ldots,a_n)$ и $P=L(a)$.
Элемент $a$ --- алгебраический над $L$, значит, расширение $P\ge L$ конечно. Кроме того, расширение $L\ge F$ также конечно в силу \ref{kon algf},
поскольку элементы $a_0,\ldots, a_n \in E$ алгебраические над $F$. Значит, расширение $P\ge F$ конечно, и поэтому $a$  алгебраический над $F$
ввиду \ref{fin alg}$(i)$.
\end{proof}

Пусть $F$ --- поле и $f\in F[x]$. \glsadd{iFldDecPol}\mem{Полем разложения многочлена} $f$ называется подполе алгебраического замыкания $\ov F$,
полученное присоединением к $F$ всех корней многочлена $f$.

Пусть $E$ --- алгебраическое расширение поля $F$. Тогда $E$ --- подполе алгебраического замыкания $\ov F$ ввиду \ref{alg pro}$(ii)$.
Расширение $E\ge F$ называется \glsadd{iExtFldNor}\mem{нормальным}, если всякий автоморфизм поля $\ov F$,
централизующий~$F$, оставляет $E$ инвариантным.

\begin{pre} \label{not raz}\mbox{}
\begin{list}{{\rm(}{\it\roman{enumi}\/}{\rm)}}
{\usecounter{enumi}\setlength{\parsep}{2pt}\setlength{\topsep}{5pt}\setlength{\labelwidth}{23pt}}
\item Расширение полей $E \ge F$ нормально тогда и только тогда, когда $E$ является полем разложения некоторого
семейства многочленов из $F[x]$.
\item Пусть $K\ge E\ge F$ --- цепочка расширений полей. Если расширение $K\ge F$ нормально, то и $K\ge E$ нормально.
\end{list}
\end{pre}
\begin{proof} $(i)$ См. \cite[\S 41]{w}.

$(ii)$ Рассмотрим $K$ как подполе алгебраического замыкания $\ov E$ поля $E$. Если некоторый автоморфизм $\a$ поля $\ov E$ централизует $E$, то он также централизует $F$.
Заметим, что $\ov E=\ov F$, т.\,к. $\ov E$ алгебраическое над $F$ ввиду предложения \ref{alg trans},
применённого к цепочке $\ov E\ge E\ge F$.
Значит, из нормальности $K$ над $F$ следует, что $\a$ оставляет поле $K$  инвариантным. Поэтому $K$ нормально над $E$.
\end{proof}

Пусть $E\ge F$ --- нормальное расширение полей характеристики $0$. \glsadd{iGrGalExt}\mem{Группой Галуа расширения}
$E\ge F$ называется подгруппа $\Aut E$, состоящая из автоморфизмов, централизующих подполе $F$. Эта подгруппа
обозначается \gls{GalEF}.

\begin{pre} \label{nor gal} Пусть $E\ge F$ --- нормальное конечное расширение полей характеристики $0$. Тогда
\begin{list}{{\rm(}{\it\roman{enumi}\/}{\rm)}}
{\usecounter{enumi}\setlength{\parsep}{2pt}\setlength{\topsep}{5pt}\setlength{\labelwidth}{23pt}}
\item группа $\Gal(E,F)$ конечна и степень расширения $E\ge F$ совпадает с $|\Gal(E,F)|$;
\item поле $F$ совпадает с множеством неподвижных точек всех элементов из $\Gal(E,F)$.
\end{list}
\end{pre}
\begin{proof} $(i)$ См. \cite[\S 57]{w}.

$(ii)$ Пусть $K$ --- множество неподвижных точек всех элементов из $\Gal(E,F)$. Имеем цепочку расширений
$E\ge K\ge F$. В силу \ref{not raz}$(ii)$ расширение $E\ge K$ нормально. Значит, из $(i)$ следует, что степень
$E$ над $K$ равна $|\Gal(E,K)|$. По предположению, всякий элемент из $\Gal(E,F)$ лежит в $\Gal(E,K)$.
В частности, $|\Gal(E,K)|\ge |\Gal(E,F)|$. Но $|\Gal(E,F)|$ совпадает со степенью $E$ над $F$ по $(i)$. То есть,
с одной стороны степень $E$ над $K$ не меньше степени $E$ над $F$, а с другой, очевидно, не больше её.
Значит, эти степени равны и, следовательно, $K=F$.
\end{proof}

Напомним, что $m$-е круговое поле $\QQ_m$ --- это расширение $\QQ(\z)$ поля $\QQ$, где $\z\in\CC$ ---
примитивный корень степени $m$ из $1$.

\begin{pre} \label{ks f} Пусть
$m$ --- натуральное число и  $\z$ --- примитивный корень степени $m$ из $1$.
\begin{list}{{\rm(}{\it\roman{enumi}\/}{\rm)}}
{\usecounter{enumi}\setlength{\parsep}{2pt}\setlength{\topsep}{5pt}\setlength{\labelwidth}{23pt}}
\item Комплексное сопряжение $\ov{\phantom{\a}}:\a\mapsto\ov\a$, $\a\in \CC$, оставляет инвариантным
поля $\ov{\QQ}$ и $\QQ_m$.
\item Ограничение отображения $\ov{\phantom{\a}}$ на $\QQ_m$ является элементом группы $\Gal(\QQ_m,\QQ)$.
\item Имеет место изоморфизм
$$
\Gal(\QQ_m,\QQ)\cong \ZZ_m^\times,
$$
сопоставляющий каждому $\s\in \Gal(\QQ_m,\QQ)$ класс вычетов $k+m\ZZ$ для некоторого целого числа $k$, взаимно простого с $m$, такой, что
$\z^\s=\z^k$.
\end{list}
\end{pre}
\begin{proof} Поскольку отображение $\ov{\phantom{\a}}$ переставляет корни любого многочлена с рациональными
коэффициентами, множество $\ov\QQ$ алгебраических чисел  инвариантно относительно комплексного сопряжения.
Так как $\QQ_m$ --- поле разложения многочлена $x^m-1$, расширение полей $\QQ_m\ge \QQ$ нормально
в силу \ref{not raz}$(i)$. Поэтому комплексное сопряжение, являясь автоморфизмом поля $\ov\QQ$, оставляет
инвариантным поле $\QQ_m$.  Значит, существует автоморфизм из $\Gal(\QQ_m,\QQ)$, действие которого на $\QQ_m$
совпадает с комплексным сопряжением.
Изоморфизм $\Gal(\QQ_m,\QQ)\cong \ZZ_m^\times$ из $(iii)$ доказан в \cite[\S 60]{w}.
Для $\s\in \Gal(\QQ_m,\QQ)$ требуемое число $k$ существует, поскольку $\z^\s$ также является примитивным корнем степени $m$ из $1$,
а все такие корни являются порождающими группы корней степени $m$ из $1$, которая циклическая в силу \ref{cyc fin}.
Обратно, если $k$ взаимно просто с $m$, то отображение $\z\mapsto\z^k$ продолжается до единственного автоморфизма из $\Gal(\QQ_m,\QQ)$.
Это следует из того, что элементы группы $\Gal(\QQ_m,\QQ)$ однозначно определяются своим действием
на $\z$ и того, что порядок $|\Gal(\QQ_m,\QQ)|$ совпадает порядком $|\ZZ_m^\times|$, т.\,е. с числом примитивных корней степени $m$ из $1$.
\end{proof}

\section{Алгебраические числа} \label{ap alg}

Комплексное число называется
\glsadd{iNumAlg}\mem{алгебраическим}, если оно является алгебраическим элементом расширения полей
$\CC\ge \QQ$. Комплексное число называется \glsadd{iNumIntAlg}\mem{целым алгебраическим},
если оно является корнем некоторого многочлена\footnote{Существенным здесь
является тот факт, что старший коэффициент многочлена равен $1$. В частности, степень $n\ge 1$. }
$$
x^n+ b_1 x^{n-1} + b_2 x^{n-2} +\ld +b_n,
$$
где $b_1,b_2,\ld,b_n\in \ZZ$. Очевидно, что всякое целое алгебраическое число будет алгебраическим.
Обратное утверждение, вообще говоря, неверно.

\begin{pre}\label{z cz}
Рациональное число $r\in \QQ$ будет целым алгебраическим
тогда и только тогда, когда $r\in \ZZ$.
\end{pre}
\textupr{Доказать предложение \ref{z cz}.}

Чтобы подчеркнуть разницу с целыми алгебраическими числами, элементы кольца $\ZZ$ целых чисел мы будем
иногда называть \glsadd{iNumIntRat}\mem{целыми рациональными}.

Оказывается, что целые алгебраические числа образуют кольцо, а алгебраические числа --- поле. Для
доказательства этого нам потребуются вспомогательные результаты.

\begin{lem} \label{lem alg cn} Справедливы следующие утверждения.
\begin{list}{{\rm(}{\it\roman{enumi}\/}{\rm)}}
{\usecounter{enumi}\setlength{\parsep}{2pt}\setlength{\topsep}{5pt}\setlength{\labelwidth}{23pt}}
\item Пусть $\a$ --- алгебраическое число. Тогда
существует $m\in \ZZ$ такое, что число $m\a$ целое алгебраическое.
\item Пусть $\om$ --- целое алгебраическое число. Тогда для любого $0\ne m\in \ZZ$ число $\om/m$
алгебраическое.
\end{list}
\end{lem}

\textupl{lem alg cn prf}{Доказать лемму \ref{lem alg cn}.}

\begin{pre} \label{alg mlt}
Пусть $W$ --- ненулевой конечно порождённый подмодуль $\ZZ$-модуля $\CC$.
Если элемент $\om\in \CC$ таков, что $\om W\se W$, то $\om$ будет целым алгебраическим числом.
\end{pre}
\begin{proof} Пусть $W$ порождается как $\ZZ$-модуль элементами $\g_1,\ld,\g_s\in \CC$.
Имеем $\om\g_i\in W$ для $i=1,\ld,s$. Таким образом, $\om\g_i=\sum_{j=1}^s a_{ij}\g_j$, где $a_{ij}\in \ZZ$.
Отсюда следует, что $\sum_{j=1}^s(a_{ij}-\d_{ij}\om)\g_j=0$. Поскольку не все $\g_j$ равны нулю, получаем
$\det (a_{ij}-\d_{ij}\om)=0$. Выписывая определитель, мы видим, что $\om$ является корнем многочлена с целыми
коэффициентами, старший коэффициент которого (с точностью до знака) равен $1$. Таким образом, $\om$ --- целое
алгебраическое число.
\end{proof}

\begin{pre} \label{alg kol}  Имеем
\begin{list}{{\rm(}{\it\roman{enumi}\/}{\rm)}}
{\usecounter{enumi}\setlength{\parsep}{2pt}\setlength{\topsep}{5pt}\setlength{\labelwidth}{23pt}}
\item множество целых алгебраических чисел образует кольцо;
\item множество алгебраических чисел образует поле.
\end{list}
\end{pre}
\begin{proof}  $(i)$ Пусть $\mu$ и $\nu$  -- целые алгебраические числа. Легко видеть, что число $-\mu$ также целое
алгебраическое, поэтому достаточно показать, что $\mu+\nu$ и $\mu\nu$ --- целые алгебраические числа. Пусть
$\mu$ и $\nu$ являются корнями многочленов $x^m+ a_1 x^{m-1} +\ld +a_m$ и $x^n+ b_1 x^{n-1} +\ld +b_n$,
соответственно, где $a_i,b_j\in \ZZ$. Пусть $W$ --- $\ZZ$-модуль, порождённый произведениями $\mu^i\nu^j$, где
$0\le i<m$, $0\le j<n$. Поскольку $m,n\geqslant 1$, имеем $W\ne 0$.
Так как $\mu^m$ выражается в виде целочисленной линейной комбинации степеней $\mu^i$, $i<m$, и, аналогично,
$\nu^n$ выражается через $\nu^j$, $j<n$, отсюда следует, что $\mu W\se W$ и $\nu W\se W$,
а, значит, и $(\mu +\nu)W\se W$ и $(\mu\nu)W\se W$. В силу предложения \ref{alg mlt} отсюда следует,
что $\mu+\nu$ и $\mu\nu$ --- целые алгебраические числа.

$(ii)$ Пусть $\a$ и $\b$ --- алгебраические числа. Из \ref{lem alg cn}$(i)$ следует, что существуют целые
алгебраические числа $\mu,\nu$  и $m,n\in \ZZ$ такие, что $\a=\mu/m$ и $\b=\nu/n$. Поэтому $-\a$, $\a+\b$,
$\a\b$ являются алгебраическими в силу $(i)$ и \ref{lem alg cn}$(ii)$. Наконец, пусть $\a\ne 0$ и $a_0\a^n+ a_1
\a^{n-1}+\ld +a_n=0$, где $a_i\in
\QQ$. Можно считать, что $a_n\ne 0$. Тогда $a_n\a^{-n}+ a_{n-1}
\a^{-(n-1)}+\ld +a_0=0$, т.\,е. $\a^{-1}$ также является алгебраическим, откуда следует требуемое. \end{proof}

Таким образом множество всех алгебраических чисел совпадает с алгебраическим замыканием $\ov{\QQ}$ поля
рациональных чисел.

Кольцо всех целых алгебраических чисел мы будем обозначать через \gls{ZZZ}.

\mem{Полем алгебраических чисел} называется конечное расширение $F$ поля $\QQ$ рациональных
чисел. Не ограничивая общности, мы будем считать, что $F$ является подполем поля комплексных
чисел $\CC$. Из конечности расширения $F\ge \QQ$ легко следует, что элементы поля $F$ являются
алгебраическими числами.

Кольцо $D=F\cap\ov{\ZZ}$, состоящее из всех целых алгебраических чисел, содержащихся в поле
алгебраических $F$, будем называть \glsadd{iRngIntFld} \mem{кольцом целых величин поля} $F$.

Нашей главной целью будет изучение свойств идеалов колец целых величин полей алгебраических
чисел.

Напомним, что для конечно порождённой свободной абелевой группы или, что то же, для свободного $\ZZ$-модуля $A$
имеет место изоморфизм $A\cong \ZZ\oplus\ld\oplus
\ZZ$, где число слагаемых определено однозначно и называется
\glsadd{iZZRnk}\mem{$\ZZ$-рангом} модуля $A$. Поскольку $\ZZ$-ранг является обобщением понятия размерности
векторного пространства, мы будем обозначать его
\gls{dimZ}.
Известно \cite[\S 16]{cr}, что всякий подмодуль свободного $\ZZ$-модуля $A$ ранга
$n$ также является свободным $\ZZ$-модулем ранга, не больше $n$.

\normalmarginpar

\begin{pre} \label{rank} Пусть $F$ --- поле алгебраических чисел
с кольцом целых величин $D$. Любой идеал $A$ кольца $D$ является конечно порождённым
свободным $\ZZ$-модулем. Если $A\ne 0$, то $\dim_\ZZ A=\dim_\QQ F$.
\end{pre}
\begin{proof} См. \cite[Предложение 12.2.2]{ir}.
\end{proof}


Коммутативное кольцо $R$ называется \glsadd{iIntDmn}\mem{областью целостности},
если $R\ne 0$ и для любых $rs\in R$ из равенства $rs=0$ следует, что либо $r=0$, либо $s=0$.
Например, любое ненулевое подкольцо поля будет областью целостности. Верно и обратное: всякая
область целостности
вкладывается в некоторое поле, например в своё \mem{поле частных}.

\normalmarginpar

\begin{pre} \label{kon oc} Всякая конечная область целостности будет полем.
\end{pre}

\textupl{kon oc prf}{Доказать предложение \ref{kon oc}.}

Назовём область целостности $R$ \mem{целозамкнутой}\glsadd{iIntClsdDmn}
в своём поле частных $K$, если любой элемент из $K$, являющийся корнем некоторого многочлена
$$
x^n+ b_1 x^{n-1} + b_2 x^{n-2} +\ld +b_n,
$$
где $b_1,b_2,\ld,b_n\in R$, лежит в $R$. Например, предложение \ref{z cz} утверждает, что кольцо $\ZZ$
целозамкнуто в $\QQ$.

\textzam{Всюду далее мы будем предполагать, что $F$ --- поле алгебраических чисел с кольцом целых величин
$D$.}

\begin{pre} \label{ma} Справедливы следующие утверждения.
\begin{list}{{\rm(}{\it\roman{enumi}\/}{\rm)}}
{\usecounter{enumi}\setlength{\parsep}{2pt}\setlength{\topsep}{5pt}\setlength{\labelwidth}{23pt}}
\item Для любого $\a\in F$
существует число $m\in \ZZ$ такое, что $m\a\in D$. В частности, $F$ совпадает с полем частных
области целостности $D$.
\item Кольцо $D$ целозамкнуто в $F$.
\end{list}
\end{pre}
\begin{proof} $(i)$ Поскольку $\a$ --- алгебраическое число, то из \ref{lem alg cn}$(i)$ следует
существование числа $m\in \ZZ$ такого, что $m\a\in \ov{\ZZ}$. Так как $m\a\in F$, то $m\a\in D$.

$(ii)$ Если  элемент $\a\in F$ удовлетворяет равенству $\a^n+ b_1 \a^{n-1}+\ld +b_n=0$, где $b_i \in D$, то
кольцо $D[\a]$ является конечно порождённым $D$-модулем. В силу \ref{rank} кольцо $D$ является
конечно-порождённым $\ZZ$-модулем, откуда следует, что и $D[a]$ является (ненулевым) конечно-порождённым
$\ZZ$-модулем. Тогда $\a\in R$ по \ref{alg mlt}, откуда следует, что $\a\in D$. \end{proof}

\begin{lem} \label{lem int z} Если $A$ --- ненулевой идеал в $D$, то $A\cap \ZZ\ne 0$.
\end{lem}
\textupl{lem int z prf}{Доказать лемму \ref{lem int z}.}

\begin{pre} \label{fak fin} Если $A$ --- ненулевой идеал в $D$,
то факторкольцо $D/A$ конечно.
\end{pre}
\begin{proof} В силу \ref{lem int z} существует $0\ne a\in A\cap\ZZ$. Тогда идеал $aD$
содержится в $A$. Поэтому достаточно показать, что факторкольцо $D/aD$ конечно. На самом деле
мы покажем, что оно состоит из $a^n$ элементов, где $n=\dim_\QQ F$.

Согласно предложению \ref{rank} существуют элементы $\om_1,\ld,\om_n\in D$ такие, что
$D=\om_1\ZZ\oplus\ld\oplus\om_n\ZZ$, где $n=\dim_\QQ F$. Отметим, что отсюда следует, что если $\om\in aD$ и
$\om=\sum_i b_i\om_i$, $b_i\in \ZZ$, то все коэффициенты $b_i$ делятся на $a$.

Пусть $\Gamma=\{\sum_i a_i\om_i\,|\,0\le a_i<a \}$. Покажем, что $\Gamma$ является множеством
представителей для смежных классов $D$ по $aD$. Пусть $\om\in D$ тогда существуют $b_i\in \ZZ$
такие, что $\om =\sum_i b_i\om_i$. Запишем $b_i=q_ia+a_i$, где $0\le a_i <a$. Тогда очевидно,
что $\om - \g \in aD$, где $\g =\sum_i a_i
\om_i\in \Gamma$.

Таким образом, каждый смежный класс $D$ по $aD$ содержит некоторый элемент из $\Gamma$. Если
$\om=\sum_i a_i\om_i$ и $\om'=\sum_i a_i'\om_i$ --- элементы из $\Gamma$, лежащие в одном классе
$D$ по $aD$, то $\om-\om'\in aD$ и, значит, целые числа $a_i-a_i'$ делятся на $a$. Так как
$0\le a_i,a_i'<a$, отсюда следует, что $a_i=a_i'$ и $\om=\om'$. Таким образом, $\Gamma$ есть
множество представителей смежных классов $D$ по $aD$ и факторкольцо $D/aD$ содержит $a^n$
элементов, как и утверждалось.
\end{proof}

Кольцо $R$ называется \mem{нётеровым (слева)}, если любая возрастающая цепочка его (левых) идеалов $I_1\se
I_2\se\ld$ стабилизируется. Другими словами, существует натуральное число $s$ такое, что $I_s=I_{s+1}=\ld$

\begin{cor} \label{cor neter} Кольцо $D$ нётерово.
\end{cor}

\textupl{cor neter prf}{Доказать следствие \ref{cor neter}.}

Идеал $I$ коммутативного кольца $R$ называется \glsadd{iIdlSmp}\mem{простым}, если
$R/I$ --- область целостности. Это эквивалентно тому, что $I\ne R$ и множество $R\setminus I$
замкнуто относительно умножения.

\textuprn{Показать, что
\begin{list}{{\rm(}{\it\roman{enumi}\/}{\rm)}}
{\usecounter{enumi}\setlength{\parsep}{2pt}\setlength{\topsep}{5pt}\setlength{\labelwidth}{23pt}}
\item всякий ненулевой простой идеал кольца $\ZZ$ имеет вид $p\ZZ$ для некоторого простого числа $p$;
\item всякий максимальный идеал коммутативного кольца будет простым. Верно ли обратное?
\end{list}}

\begin{cor} \label{cor pr max} Каждый ненулевой простой идеал из $D$ максимален.
\end{cor}
\textupl{cor pr max prf}{Доказать следствие \ref{cor pr max}.}

Мы установили три важных свойства кольца $D$, которые легли в основу следующего определения.

\begin{opr} Область целостности, которая является нётеровой, целозамкнутой в своём поле частных и в
которой всякий ненулевой простой идеал максимален называется \glsadd{iDdkDmn}\mem{дедекиндовой областью}.
\end{opr}

Таким образом, кольца целых величин полей алгебраических чисел являются дедекиндовыми
областями. Замечательным свойством дедекиндовых областей является справедливость аналога
основной теоремы арифметики: всякий ненулевой идеал дедекиндовой области допускает однозначное
разложение в произведение простых идеалов. В этом идеалы схожи с целыми числами. Кроме того,
как и для целых чисел, множество идеалов дедекиндовой области можно вложить в более широкий
класс <<дробных идеалов>>, в котором всякий ненулевой идеал будет обратим.

%

\begin{opr} \glsadd{iIdlFrc}\mem{Дробным идеалом} поля $F$ называется конечно порождённый $D$-подмодуль
$D$-моду\-ля $F$. Дробный идеал поля $F$ называется
\glsadd{iIdlInt}\mem{целым}, если он содержится в $D$.
\end{opr}

\begin{pre} \label{cel ids} Всякий идеал кольца $D$ является целым идеалом поля $F$.
\end{pre}
\begin{proof} В силу \ref{rank} всякий идеал кольца $D$ конечно порождён как $\ZZ$-модуль, а
значит, и конечно порождён как $D$-модуль, т.\,е. является дробным идеалом в $F$.
\end{proof}

\begin{pre} \label{inv cel} Пусть $A$ --- дробный идеал в $F$. Тогда в $F$ существует дробный
идеал $B$ такой, что $AB=D$.
\end{pre}
\begin{proof}
См. \cite[Теорема 18.13]{cr}
\end{proof}

\begin{cor} \label{req lem} Пусть $I$ --- собственный идеал из $\ov{\ZZ}$ и пусть
$\g_1,\ld,\g_n\in \CC$ --- алгебраические числа, не все равные нулю. Тогда в поле $\QQ(\g_1,\ld,\g_n)$
существует элемент $\b$ такой, что $\g_i\b\in \ov{\ZZ}$ для всех $\,i$, но не все элементы $\g_i\b$ лежат в $I$.
\end{cor}
\begin{proof} Легко видеть, что поле $F=\QQ(\g_1,\ld,\g_n)$ является полем алгебраических
чисел. Пусть $D=F\cap \ov{\ZZ}$ --- кольцо целых величин поля $F$. Рассмотрим $D$-модуль $A$, порождённый
элементами $\g_1,\ld,\g_n$. Тогда $A$ --- дробный идеал поля $F$. Из \ref{inv cel} следует, что в $F$
существует дробный идеал $B$ такой, что $AB=D$. В частности, $AB\nsubseteq I$ и, значит, существует элемент $\b
\in B$ такой, что $A\b \nsubseteq I$. Но тогда $\g_i\b\not\in I$ для  некоторого $i$. Отсюда следует требуемый
результат.
\end{proof}

\section{Доказательства некоторых утверждений}

{\sc \hypertarget{max id prf}{Доказательство}} \ref{max id}. Пусть $I$ --- собственный правый (левый, двусторонний) идеал кольца $R$ и пусть $\SC$
--- множество собственных правых (левых, двусторонних) идеалов кольца $R$, содержащих $I$.
Если показать, что $\SC$ содержит максимальный элемент, то он, очевидно, будет максимальным правым (левым,
двусторонним)  идеалом в $R$, содержащим $I$. По лемме Цорна \ref{thm lz} достаточно проверить, что любая цепь
$$
I_1\le I_2\le\ld\myeqno\label{cep}
$$
элементов из $\SC$ ограничена сверху. Пусть $L=\bigcup_{i\ge1} I_i$.
Очевидно, что $L$ --- правый (левый, двусторонний) идеал в $R$. Кроме того,  $L$ --- верхняя грань
для цепи $\ref{cep}$
и нужно лишь проверить, что $L\in \SC$. Если это не так, то $L=R$, поскольку $L$ содержит $I$. Значит, $1_R$
лежит в некотором идеале $I_{i_0}$. Но тогда $R=I_{i_0}$, что противоречит включению $I_{i_0}\in\SC$.\qed

\medskip\noindent
{\sc \hypertarget{pro prf}{Доказательство}} \ref{pro}. Пусть $I\nor_r R$, $I\ne R$ и элемент $r\in I$ обратим справа, т.\,е. существует $s\in R$
такой, что $rs=1$. Тогда $1\in I$ и, значит, для любого $t\in R$ имеем $t=1t\in I$. Следовательно, $I=R$
вопреки предположению. Случай левого идеала рассматривается аналогично.

\medskip\noindent
{\sc \hypertarget{zr prf}{Доказательство}} \ref{zr}. $(ii)$ Включение  $\Z(R)^\times\le \Z(R)\cap R^\times$ очевидно. Обратно, если элемент $z\in
Z(R)$ обратим в $R$ и $w$ --- его обратный, то для произвольного $r\in R$, умножив равенство $zr=rz$ справа и
слева на $w$, получим $rw=wr$. Поэтому $w\in \Z(R)$ и, значит $z\in Z(R)^\times$.

$(iii)$ Всякий элемент $z\in \Z(R)^\times$ лежит в $R^\times$ и перестановочен со всеми
элементами из $R$. В частности, $z\in \Z(R^\times)$.
\qed

\medskip\noindent
{\sc \hypertarget{obr prf}{Доказательство}} \ref{obr}. Поскольку $(i)$ является частным случаем $(ii)$ при $I=\{0\}$, то докажем $(ii)$. Пусть $r\in
R\setminus I$. Существует элемент $s\in R$ такой, что $rs=1$. Ясно, что $s\in R\setminus I$.
Также существует элемент $t\in R\setminus I$ такой, что $st=1$. Теперь $t=rst=r$ и, значит, $sr=1$, т.\,е. элемент $r$ обратим слева.

$(iii)$ Пусть $r\in I$ и $1-s$ --- правый обратный для $1-r$. Тогда
$$
(1-r)(1-s)=1, \myeqno\label{1a1b}
$$
откуда следует, что $s=r(s-1)\in I$. Но тогда по предположению $1-s$ также обратим справа.
Пусть $t$ --- его правый обратный. Тогда, умножив справа на $t$ обе части равенства \ref{1a1b} получим $1-r=t$.
То есть $(1-s)(1-r)=1$, что и требовалось доказать.
\qed

\medskip\noindent
{\sc \hypertarget{mpl prf}{Доказательство}} \ref{mpl}. Пусть $I$ --- идеал кольца $R$, являющийся максимальным правым идеалом. Тогда $R$ ненулевое и в
факторкольце $R/I$ каждый ненулевой элемент обратим справа. По \ref{obr}$(i)$ в $R/I$ каждый ненулевой элемент обратим слева.
Поэтому $R/I$ --- тело и $I$ --- максимальный левый идеал в $R$.

\medskip\noindent
{\sc \hypertarget{hom r prf}{Доказательство}} \ref{hom r}. $(i)$ Определим отображение $\t: S\to T$ следующим образом. Пусть $s\in S$. В
силу сюръективности гомоморфизма $\vf$ существует элемент $r\in R$ такой, что $r\vf=s$. Тогда мы положим $s\t=r\psi$.
Проверим корректность этого определения. Пусть $r_0\in R$ также является прообразом элемента $s$ относительно
$\vf$. Тогда существует элемент $k\in \Ker\vf$ такой, что $r_0-r=k$. По условию $\Ker\vf\se\Ker\psi$ и, значит,
$k\in \Ker\psi$. Поэтому $(r_0-r)\psi=k\psi=0$, т.\,е. $r_0\psi=r\psi$.

Если $s_1,s_2\in S$ и элементы $r_1,r_2\in R$ удовлетворяют равенствам $r_1\vf=s_1$ и $r_2\vf=s_2$, то $(r_1+r_2)\vf=s_1+s_2$
и $(r_1r_2)\vf=s_1s_2$. Поэтому $(s_1+s_2)\t=(r_1+r_2)\psi=r_1\psi+r_2\psi=s_1\t+s_2\t$ и, аналогично,
$(s_1s_2)\t=s_1\t s_2\t$. Ясно также, что $(1_S)\t=1_T$, поскольку $(1_R)\vf=1_S$ и $(1_R)\psi=1_T$.
Итак, $\t$ является требуемым гомоморфизмом колец.

Пусть $\t_0: S\to T$ --- другой гомоморфизм с условием $\psi=\vf\t_0$. Возьмём $s\in S$.
Существует элемент $r\in R$ такой, что $r\vf=s$. Тогда $s\t_0=r\vf\t_0=r\psi=r\vf\t=s\t$.
Значит, $\t_0=\t$.

Утверждение $(ii)$ доказывается аналогично.
\qed

\medskip\noindent
{\sc \hypertarget{subdir prf}{Доказательство}} \ref{subdir}. Достаточно доказать включение $I\se(I\cap I_1)\oplus\ld\oplus(I\cap I_n)$. Пусть $r\in
I$. Запишем $1=e_1+\ld+e_n$, где $e_i\in I_i$. 
Пусть $I_i\nor_l R$ для всех $i$ и $I\nor_r R$. Тогда имеем $r=r(e_1+\ld +e_n)=re_1+\ld +re_n$, и $re_i\in I_i$, так
как $I_i\nor_l R$. С другой стороны, $re_i\in I$, поскольку $I\nor_r R$. Следовательно, $re_i\in I\cap I_i$,
откуда следует требуемое. Доказательство для случая $I_i\nor_r R$ для всех $i$ и $I\nor_l R$ аналогично, 
если записать $r=(e_1+\ld +e_n)r$.
\qed

\medskip\noindent
{\sc \hypertarget{dec uniq prf}{Доказательство}} \ref{dec uniq}. Так как $I_i\nor R$ из \ref{subdir} следует, что $I_i=(I_i\cap
J_1)\oplus\ld\oplus(I_i\cap J_n)$. Поскольку $I_i\cap J_j\nor R$, по условию получаем, что существует индекс
$i\s\in \{1,\ld,n\}$ такой, что $I_i\cap J_{i\s}=I_i$ и $I_i\cap J_j=0$ при $j\ne i\s$, т.\,е. $I_i\se J_{i\s}$
для всех $i=1,\ld,m$. Так как $R=\bigoplus_i I_i=\bigoplus_i J_{i\s}$ и все идеалы $J_j$ ненулевые, отсюда
следует, что отображение $\s:\{1,\ld,m\}\to \{1,\ld,n\}$ сюръективно и, в частности, $m\ge n$.

Аналогично, существует сюръективное отображение $\t:\{1,\ld,n\}\to \{1,\ld,m\}$ такое что $J_j\se I_{j\t}$ для
всех $j=1,\ld,n$. В частности, $n\ge m$. Таким образом, $m=n$,
т.\,е. отображение $\s$ биективно. \qed

\medskip\noindent
{\sc \hypertarget{ndr prf}{Доказательство}} \ref{ndr}. Поскольку $(i)$ следует из $(ii)$, докажем $(ii)$.
Пусть отображение $\a:\End_R(eR)\to
eRe$ действует по правилу $\a(\vf)=e\vf$ для любого $\vf \in \End_R(eR)$. Так как $e\vf\in eR$, имеем
$e\vf=e(e\vf)$. Кроме того, $e\vf=(ee)\vf=(e\vf)e$ ввиду того, что $\vf \in \End_R(eR)$.
Значит, $e\vf=e(e\vf)e$, т.\,е. образ $\a$ действительно лежит в $eRe$.

Пусть $\vf,\psi \in \End_R(eR)$. Тогда
$$
\a(\vf+\psi)=e(\vf+\psi)=e\vf+e\psi=\a(\vf)+\a(\psi),
$$
т.\,е. отображение $\a$ аддитивно. Покажем инъективность. Пусть $\a(\vf)=0$. Тогда для любого $r\in R$ имеем
$$
(er)\vf=(e\vf)r=\a(\vf)r=0r=0,
$$
откуда следует, что $\vf$ --- нулевой элемент из $\End_R(eR)$.

Покажем, что $\a$ сюръективно. Пусть $r\in R$. Определим $\vf\in \End_R(eR)$ по правилу $(es)\vf=eres$ для
любого $s\in R$. Легко видеть, что $\vf$ является $R$-линейным, т.\,е. действительно лежит в $\End_R(eR)$.
Из равенства $\a(\vf)=e\vf=ere$ вытекает сюръективность $\a$.

Пусть, наконец, $\vf,\psi \in \End_R(eR)$. Тогда
$$
\a(\vf\psi)=e\vf\psi=(e(e\vf))\psi=e\psi\,e\vf=\a(\psi)\,\a(\vf),
$$
и значит, $\a$ является антиизоморфизмом колец.
\qed

\medskip\noindent
{\sc \hypertarget{zri prf}{Доказательство}} \ref{zri}. Пусть $z\in \Z(R^\times)$. В силу \ref{zr}$(ii)$ достаточно показать, что $z\in \Z(R)$. Пусть
$r\in R$. По условию  $r=r_1z_1+\ld +r_k z_k$ для подходящих $r_i\in R^\times$, $z_i\in \Z(R)$, $i=1,\ld,k$.
Поэтому $zr=rz$, откуда следует требуемое.
\qed

\medskip\noindent
{\sc \hypertarget{idem osd prf}{Доказательство}} \ref{idem osd}. Утверждение $(ii)$ будем рассматривать в случае правых идеалов. Рассуждение для левых
аналогично.

$(i)\Rightarrow(ii)$ Предположим, что $eR=I_1\oplus I_2$, где $I_1,I_2\le eR$ --- правые идеалы из $R$.
Запишем $e=e_1+e_2$, где $e_i\in I_i$, $i=1,2$.
Поскольку умножение слева на $e$ действует тождественно на все элементы из $eR$, мы получаем
$e_i=ee_i=e_1e_i+e_2e_i$. Но $e_je_i\in I_j$, $j=1,2$, и, ввиду однозначности разложения,
отсюда следует $e_i=e_i^2$ и $e_je_i=0$ при $j\ne i$. В силу примитивности $e$, один из идемпотентов $e_1$, $e_2$ равен нулю.
Обозначим его $e_i$. Тогда $e=e_j\in I_j$, где $j\ne i$, и $eR\le I_j$, откуда следует $I_i=0$.

$(ii)\Rightarrow(iii)$ Пусть, напротив, $e_1$ --- идемпотент из $eRe$, отличный от $0$ и $e$.
Обозначим $e_2=e-e_1$. Тогда $e_2$ также идемпотент из $eRe$, ортогональный $e_1$ и отличный от $0$ и $e$.
Поскольку $e_1=ee_1\in eR$ и $e_2=ee_2\in eR$, то $e_1R\le eR$ и $e_1R\le eR$.
Значит, $e_1R+e_2R\le eR$. Обратно, для любого $er\in eR$ имеем $er=e_1r+e_2r\in e_1R+e_2R$, т.\,е. $e_1R+e_2R=eR$.
Осталось показать, что сумма $e_1R+e_2R$ прямая. Пусть $s\in e_1R\cap e_2R$. Тогда существуют $r_1,r_2\in R$
такие, что $s=e_1r_1=e_2r_2$. Но тогда $s=e_1r_1=e_1(e_1r_1)=e_1(e_2r_1)=0$ в силу ортогональности идемпотентов
$e_1$ и $e_2$.

$(iii)\Rightarrow(i)$ Допустим, что $e=e_1+e_2$, где $e_i\in R$, $i=1,2$, --- ортогональные идемпотенты.
Тогда $ee_i=(e_1+e_2)e_i=e_i^2=e_i$ и, аналогично, $e_ie=e_i$. Значит, $e_i=ee_ie\in eRe$. По
условию $e_i=0$ или $e$.\qed

\medskip\noindent
{\sc \hypertarget{idem dsd prf}{Доказательство}} \ref{idem dsd}. $(i)\Rightarrow(ii)$ Пусть $eR=I_1\oplus I_2$, где $I_1,I_2\le eR$ ---
двусторонние идеалы из $R$. Запишем $e=e_1+e_2$, где $e_i\in I_i$, $i=1,2$. Поскольку $e$ ---
единица кольца $Re$, то $e_i$ --- центральные идемпотенты кольца $eR$ по предложению
\ref{inn dir sum}. Покажем, что $e_i\in \Z(R)$. Пусть $r\in R$. Учитывая, что
$re_i,e_ir\in Re$ и $e_i\in \Z(Re)$, мы получаем $e_ir=e_i(e_ir)=e_ire_i=(re_i)e_i=re_i$. Из
примитивности идемпотента $e$ следует, что $e_i=0$ для $i=1$ или $2$, то есть $e=e_j$ где
$i\ne j\in\{1,2\}$. Поэтому $eR\le I_j$ и $I_i=0$.

$(ii)\Rightarrow(i)$ Предположим, что $e=e_1+e_2$, где $e_i\in Z(R)$, $i=1,2$, ---
ортогональные идемпотенты. Тогда $e_i=e_ie\in Re$ и, значит, идеалы $I_i=Re_i$ лежат в $Re$.
Любой элемент $re\in Re$ можно записать как $re=re_1+re_2$, то есть $Re=I_1+I_2$. Покажем, что
эта сумма прямая. Если $s\in I_1\cap I_2$, то существуют $r_1,r_2\in R$ такие, что
$s=r_1e_1=r_2e_2$. Но тогда $s=r_1e_1=(r_1e_1)e_1=(r_2e_2)e_1=0$ в силу ортогональности
идемпотентов $e_1$ и $e_2$. По условию $I_i=0$ для $i=1$ или $2$ и, значит, $e_i=0$.
\qed

\medskip\noindent
{\sc \hypertarget{m id dec prf}{Доказательство}} \ref{m id dec}. $(i)$ $Me_i$ --- подгруппа аддитивной группы $M$, которая выдерживает умножение на любой
элемент $r\in R$, поскольку $Me_ir=Mre_i\se Me_i$ в силу того, что $e_i\in \Z(R)$. Следовательно,
$Me_i$ --- $R$-подмодуль.

$(ii)$ Пусть $m_1e_1+\ld+m_ne_n=0$ для некоторых $m_i\in M$. Умножив справа на $e_i$ и учитывая ортогональность
идемпотентов, получим $m_ie_i=0$ для любого $i$. Отсюда вытекает, что сумма $Me_1+\ld+Me_n$ прямая.

$(iii)$ Для любого $m\in M$ имеем $m=m1_R=me_1+\ld+me_n\in Me_1+\ld+Me_n$. Учитывая $(ii)$, получаем требуемое.
\qed

\medskip\noindent
{\sc \hypertarget{kl sum prf}{Доказательство}} \ref{kl sum}. $(i)$ Пусть $K\in \K(G)$. Достаточно показать, что любой элемент $g\in G$ коммутирует
с $\wh{K}$. Это следует из соотношений
$$
g^{-1}\wh{K}g=\sum_{x\in K}g^{-1}xg=\wh{K}.
$$

$(ii)$ Пусть элемент $a=\sum_{x\in G} \a_x x$ лежит в
$\Z(RG)$. Тогда для каждого $g\in G$
$$
\sum_{x\in G}\a_x x=a=g^{-1}ag=\sum_{x\in G}\a_x (g^{-1}xg).
$$
Сравнивая коэффициенты при элементах $x\in G$ в правой и левой частях, получим, что $\a_x=\a_{gxg^{-1}}$
для всех $g\in G$. Таким образом, коэффициенты при любых элементах одного и того же класса сопряжённости
совпадают и, значит, $a$ является $R$-линейной комбинацией элементов $\wh{K}$, $K\in
\K(G)$.

Так как классовые суммы $\wh{K}$, где $K\in \K(G)$, являются суммами
элементов непересекающихся подмножеств группы $G$, то коэффициенты в разложении любого элемента из $\Z(RG)$
в виде $R$-линейной комбинации сумм $\wh{K}$ определяются однозначно.

$(iii)$ Если $z_1, z_2\in M$, то существует элемент $g\in G$ такой, что $z_2=z_1^g$. Поэтому между
множеством пар $(x_1,y_1)$, где $x_1\in K$, $y_1\in L$, таких, что $x_1y_1=z_1$, и множеством пар
$(x_2,y_2)$, где $x_2\in K$, $y_2\in L$, таких, что $x_2y_2=z_2$, можно установить взаимно однозначное
соответствие $(x_1,y_1)\mapsto (x_1^g,y_1^g)$. Поэтому мощности этих множеств совпадают.

$(iv)$ Фиксируем $z\in M$. Коэффициент при $\whm$ в произведении $\wh{K}\wh{L}$ равен коэффициенту при
$z$ и, значит, совпадает с числом $a_{\mbox{}_{KLM}}$ пар $(x,y)$, где $x\in K$, $y\in L$, таких,
что $xy=z$.
\qed

\medskip\noindent
{\sc \hypertarget{zsr prf}{Доказательство}} \ref{zsr}. Включение $\Z(SG)\supseteq \Z(RG)\cap SG$  очевидно. Пусть $z\in \Z(SG)$
из \ref{kl sum}$(ii)$ следует, что $z$ является $S$-линейной комбинацией классовых сумм $\wh K$,
а, значит, и их $R$-линейной комбинацией. Поэтому $z\in\Z(RG)$ опять-таки в силу \ref{kl sum}$(ii)$.
\qed

\medskip\noindent
{\sc \hypertarget{kvt prf}{Доказательство}} \ref{kvt}. Элемент $\sum_{g\in G} a_g g$ из $RG$ лежит в $\Ker\ti\vf$ тогда и только тогда, когда
$$0=\Big(\sum_{g\in G} a_g g\Big)\ti\vf=\sum_{g\in G} a_g\vf\ g.$$
Ввиду однозначности представления элемента $0\in SG$ в виде $S$-линейной комбинации элементов группы $G$ имеем
$a_g\vf=0$ для всех $g\in G$. Это означает, что $a_g\in \Ker\vf$ для всех $g\in G$.

Аналогично, элемент $\sum_{K\in \K(G)} a_{\mbox{}_K}\wh{K}$ из $\Z(RG)$ лежит в $\Ker\ti\vf$ тогда и только
тогда, когда $a_{\mbox{}_K}\in \Ker\vf$ для всех $K\in \K(G)$.
\qed

\medskip\noindent
{\sc \hypertarget{ann el prf}{Доказательство}} \ref{ann el}. Пусть $v\in V$ и $v\ne 0$. Отображение $\tau:a\mapsto va$ является гомоморфизмом $A$-модулей
$A^\circ\to V$. Его образ ненулевой, поскольку содержит $1_A\tau=v$. Из простоты $V$ следует, что $\tau$ ---
эпиморфизм. Положим $I=\Ker\tau$. Очевидно, что $I=\Ann(v)$. Но тогда $I$ --- максимальный правый идеал алгебры
$A$ и  $A^\circ/I\cong V$ по теореме \ref{thm o hom}.

\medskip\noindent
{\sc \hypertarget{dir sum prf}{Доказательство}} \ref{dir sum}. Как следует из замечаний после определения \ref{opr dir kol}, алгебры $A_1$ и $A_2$
можно рассматривать как идеалы в алгебре $A$. По предложению \ref{subdir} любой максимальный правый
идеал из $A$ имеет вид $I_1\oplus A_2$ или $A_1\oplus I_2$, где $I_1$ и $I_2$
--- максимальные правые идеалы в $A_1$ и $A_2$, соответственно. Поэтому равенство
$\J(A)=\J(A_1)\oplus\J(A_2)$ следует из предложения
\ref{rad int max}. \qed

\medskip\noindent
{\sc \hypertarget{id rad prf}{Доказательство}} \ref{id rad}. Приведём два доказательства этого предложения.

Пусть, напротив, $I\nsubseteq \J(A)$. Из предложения \ref{rad int max} следует, что $I$ не содержится в
некотором максимальном правом идеале $M$. Тогда $R=I+M$ и $1=a+m$ для некоторых $a\in I$, $m\in M$. Но
по условию элемент $m=1-a$ обратим справа. Это противоречит тому, что он содержится в собственном правом
идеале $M$.

Второе доказательство. Пусть, напротив, $I\nleqslant \J(A)$. Тогда существует неприводимый $A$-модуль
$V$ такой, что $VI\ne 0$, а, значит, существует элемент $0\ne v\in V$, для которого $vI\ne 0$. В силу
простоты $V$ получаем $vI=V$. В частности, существует элемент $a\in I$, для которого $va=v$, т.\,е.
$v(1-a)=0$. По условию элемент $1-a$ обратим справа, откуда следует $v=0$, противоречие.
\qed

\medskip\noindent
{\sc \hypertarget{jre prf}{Доказательство}} \ref{jre}. $(i)$ Включение $e\J(A)e\se \J(A)\cap eAe$ очевидно.

$(ii)$ Докажем включение $\J(A)\cap eAe\se \J(eAe)$. Легко видеть, что $\J(A)\cap eAe$ --- идеал алгебры $eAe$.
В силу \ref{id rad} достаточно показать, что для всякого элемента $a\in \J(A)\cap eAe$ элемент $e-a$ обратим
справа в $eAe$. Поскольку $a\in \J(A)$, из \ref{rad prop} следует, что существует
элемент $b\in A$ такой, что $(1-a)b=1$. Тогда
$$
(e-a)(ebe)=e^2be-aebe=ebe-eabe=e(b-ab)e=e^2=e,
$$
где $ae=ea$, поскольку $a\in eAe$ и $e$ --- единица кольца $eAe$.
Итак, $ebe$ --- правый обратный к $e-a$ в $eAe$, что и требовалось показать.

$(iii)$ Докажем включение $\J(eAe)\se e\J(A)e$. Пусть $a\in \J(eAe)$.
Тогда $a=eae$ опять-таки поскольку $a\in eAe$ и $e$ --- единица кольца $eAe$.
Достаточно показать, что $a\in \J(A)$, т.\,к. отсюда будет следовать, что $a=eae\in e\J(A)e$.

Рассмотрим правый идеал $aA$ из $A$. Если мы покажем, что для всякого элемента $ax\in aA$, где $x\in A$,
элемент $1-ax$ обратим справа в $A$, то в силу \ref{id rad} получим $a\in aA\se\J(A)$, как и требуется.

Так как $a\in \J(eAe)$ и $exe\in eAe$, то $axe=aexe\in \J(eAe)$, и в силу \ref{rad prop}
существует элемент $b\in eAe$ такой, что $(e-axe)b=e$ и, значит,
$$
e=(e-axe)b=(1-ax)eb=(1-ax)b,
$$
поскольку $eb=b$. Отсюда получаем
$$
(1-ax)(1+bax)=(1-ax)+(1-ax)bax=1-ax+eax=1,
$$
поскольку $ea=a$.
Итак, элемент $1+bax$ --- правый обратный к $1-ax$ в $A$, что и требовалось показать.
\qed

\medskip\noindent
{\sc \hypertarget{rad max prf}{Доказательство}} \ref{rad max}. $(i)$ Пусть идеал  $I\nor_r A$ состоит из нильпотентных элементов. Поскольку для
любого $a\in I$ элемент $1-a$ обратим, то требуемое следует из леммы
\ref{id rad}.

\medskip\noindent
$(ii)$ Утверждение следует из пункта $(i)$, поскольку в нильпотентном идеале всякий элемент
нильпотентен.

Приведём ещё одно доказательство. Пусть идеал  $I\nor_r A$ нильпотентен и $I^n=0$. Для
включения $I\le\J(A)$ достаточно проверить, что $I$ аннулирует любой неприводимый $A$-модуль
$V$. Пусть, напротив, $VI\ne 0$. Тогда $VI=V$ в силу простоты $V$, и мы получаем
$$
V=VI=(VI)I=VI^2=\ld=VI^n=0,
$$
противоречие. \qed

\medskip\noindent
{\sc \hypertarget{com nil prf}{Доказательство}} \ref{com nil}. Если $a^n=0$ и $b^m=0$, то, используя коммутативность, получим
$(a+b)^{m+n}=0$ и $(ta)^n=0$ для любого $t\in A$. Второе утверждение следует из обратимости
элемента $1-a$ и предложения \ref{id rad}. \qed

\medskip\noindent
{\sc \hypertarget{con por prf}{Доказательство}} \ref{con por}. Пусть $V$ --- ненулевой конечно порождённый $A$-модуль, $v_1,\ld, v_n$
--- его порождающие и  $\SC$ --- множество собственных подмодулей модуля $V$.
По лемме Цорна \ref{thm lz} достаточно показать, что любая цепь
$$
V_1\le V_2\le\ld\myeqno\label{tsep}
$$
элементов из $\SC$ ограничена сверху. Пусть $W=\bigcup_{i\ge1} V_i$. Тогда $W$ --- верхняя грань для цепи \ref{tsep}
и нужно лишь проверить, что $W\in \SC$. Если это не так, то $W=V$ и, значит, каждый порождающий $v_j$ лежит в
некотором подмодуле $V_{i_j}$. Но тогда все порождающие содержатся в $V_{i_0}$, где $i_0=\max\{i_1,\ld,i_n\}$.
Следовательно, $V=V_{i_0}$, что противоречит включению $V_{i_0}\in \SC$.

\medskip\noindent
{\sc \hypertarget{loc eq prf}{Доказательство}} \ref{loc eq}. Докажем, что $(iii)$ эквивалентно $(i)$ и $(ii)$.

Пусть верно $(i)$.
Поскольку максимальный правый идеал $I_r$ единственный, то он совпадает
с $\J(S)$ по предложению \ref{rad int max}. Пусть $a\not\in I_r$.
Тогда, ввиду правой обратимости элемента $a+I_r$ из $S/I_r$ существует элемент $s\in S$ такой, что
$as=1-c$ для некоторого $c\in \J(S)$. Но тогда $as$ обратим справа в $S$ по предложению \ref{rad prop}.
Следовательно, и элемент $a$ обратим справа. Поэтому в силу \ref{obr}$(ii)$ имеет место $(iii)$.

Доказательство того, что из $(ii)$ следует $(iii)$ проводится <<симметрично>>.

Пусть верно $(iii)$. Тогда любой максимальный правый (левый) идеал содержится в $I$,
т.\,к. в нём нет обратимых справа (слева) элементов. Значит, имеют место пункты $(i)$ и $(ii)$.
\qed

\medskip\noindent
{\sc \hypertarget{rad nlp prf}{Доказательство}} \ref{rad nlp}. Заметим, что идеал $\J(A)^n$ является $F$-подпро\-странством в $A$ и,
значит, конечномерен. Поскольку $A$ содержит $F$, $\J(A)^n$ конечно порождён как $A$-модуль.
Если $\J(A)^n\ne 0$, то по теореме \ref{thm nak} имеет место строгое включение
$\J(A)^{n+1}<\J(A)^n$, и поэтому из соображения размерности ясно, что идеал $\J(A)$
нильпотентен. \qed

\medskip\noindent
{\sc \hypertarget{rad cent prf}{Доказательство}} \ref{rad cent}. Из предложения \ref{rad nlp} следует,
что $\J(A)\cap B$ --- нильпотентный идеал
алгебры $B$. Значит, по предложению \ref{rad max} имеется включение $\J(A)\cap B\le \J(B)$.
Обратно, пусть $z\in\J(B)$. Достаточно показать, что $z\in
\J(A)$. Правый идеал $zA$ нильпотентен, поскольку $z$ коммутирует со всеми
элементами из $A$ и $z^n=0$ для некоторого $n$. Поэтому $z\in zA\le \J(A)$.

Отметим, что включение $\J(B)\le\J(A)\cap B$ также вытекает из предложения \ref{nrz}, поскольку $A$
конечно порождена как $F$-модуль, а значит, и как $B$-модуль ввиду включения $F\le B$
\qed

\medskip\noindent
{\sc \hypertarget{ssm nsl prf}{Доказательство}} \ref{ssm nsl}. По пункту $(ii)$ предложения \ref{equ ssm} фактормодуль $V/U$ изоморфен
некоторому подмодулю из $V$. Поэтому достаточно доказать полную приводимость подмодуля $U$. Пусть $W\le
U$. Тогда по пункту $(ii)$ предложения \ref{equ ssm} существует подмодуль $W_0\le V$ такой, что
$V=W\oplus W_0$. Тогда $U=W\oplus (U\cap W_0)$. \qed

\medskip\noindent
{\sc \hypertarget{ss wed prf}{Доказательство}} \ref{ss wed}. $(i)$ Модуль $V$ является суммой своих неприводимых подмодулей в силу
\ref{equ ssm}$(iii)$, а каждый неприводимый подмодуль $M\le V$ содержится в компоненте $\W_M(V)$. Значит,
$V=\sum_{M\in \M(A)}\W_M(V)$. Эта сумма является прямой, поскольку
из теоремы Жордана--Гёльдера \ref{thm zh-g} и следствия \ref{cor sum dir} вытекает, что всякий композиционный фактор модуля $\W_M(V)$ изоморфен
$M$, а модуль $\sum_{\{U\in \M(A)\mid U\ncong M\}}\W_U(V)$ не имеет композиционных факторов изоморфных $M$, поэтому эти
модули пересекаются по нулю для любого $M\in \M(A)$.

$(ii)$ Пусть $\vf\in\End_A(V)$. Достаточно показать, что
$\W_M(V)\vf\se\W_M(V)$. Заметим, что если $U\le V$ и $U\cong M$, то $U\vf$ --- подмодуль из $V$,
являющийся гомоморфным образом подмодуля $U$. В силу неприводимости $U\vf$ изоморфен $0$ или $M$, т.\,е.
$U\vf \se \W_M(V)$. Отсюда следует требуемое включение.

$(iii)$ Очевидно, что правая часть содержится в $\W_M(V)$. Обратно, пусть $U\le V$ и $U\cong M$. Пусть
$\pi_i\in \End_A(V)$ --- проекция $V$ на прямое слагаемое $V_i$. Очевидно, что $U\se \bigoplus_iU\pi_i$.
С другой стороны, если $U\pi_i\ne 0$, то $U\pi_i\cong M$ и, значит, $U\se \bigoplus_{V_i\cong M} V_i$.

$(iv)$ Число подмодулей $V_i$, изоморфных $M$ равно $\dim_F(\W_M(V))/\dim_F M$ в силу $(ii)$.
\qed

\medskip\noindent
{\sc \hypertarget{fact ssm prf}{Доказательство}} \ref{fact ssm}. Достаточно заметить, что у алгебр $A$ и $A/\J(A)$ одинаковое множество
неприводимых модулей.
\qed

\medskip\noindent
{\sc \hypertarget{cor min idl prf}{Доказательство}} \ref{cor min idl}. Если $V$ --- неприводимый $A$-модуль, то по предложению \ref{ann el} он изоморфен
фактор модулю $A^\circ/I$, где $I$ --- некоторый максимальный правый идеал алгебры $A$. По предложению
\ref{ssm crit} модуль $A^\circ$ вполне приводим и, значит, существует правый идеал $J$ такой, что $A^\circ
=I\oplus J$. Поскольку $J\cong A^\circ/I\cong V$ и модуль $V$ простой, то $J$ --- минимальный правый идеал.
\qed

\medskip\noindent
{\sc \hypertarget{lem ff2 prf}{Доказательство}} \ref{lem ff2}.
Пусть $V=\Im\vf\oplus\Ker\vf$ и $v\in V$. Можно записать $v=u\vf+w$ для подходящих $u\in V$, $w\in \Ker\vf$.
Тогда $v\vf=u\vf^2$, и значит, $V\vf\se V\vf^2$. Обратное включение очевидно.

Пусть теперь $V\vf=V\vf^2$. Сначала докажем, что $V=\Im\vf+\Ker\vf$. Если $v\in V$, то $v\vf=u\vf^2$ для некоторого $u\in V$.
Запишем $v=u\vf +(v-u\vf)$. Ясно, что $u\vf\in \Im\vf$, и  поскольку $(v-u\vf)\vf=v\vf-u\vf^2=0$, имеем $v-u\vf \in \Ker\vf$.

Осталось показать, что $\Im\vf\cap\Ker\vf=0$. Заметим, что из условия $V\vf=V\vf^2$  и конечномерности $V$ следует, что $\vf$
действует взаимно однозначно на $V\vf$. В частности, из равенства $u\vf^2=0$ для $u\in V$ вытекает $u\vf=0$.

Пусть $v\in \Im\vf\cap\Ker\vf$. Тогда имеем $v=u\vf$ для подходящего $u\in V$, а также $0=v\vf=u\vf^2$. По замечанию выше
получаем $0=u\vf=v$, что и требовалось.\qed

\medskip\noindent
{\sc \hypertarget{mod ns prf}{Доказательство}} \ref{mod ns}. Обозначим $t=\sum_{g\in G}g$. Легко видеть, что $t\in \Z(FG)$ и
$t^2=|G|t=0$, т.\,е. элемент $t$ нильпотентный. Из следствия \ref{cor com rad} вытекает, что
$t\in\J(\Z(FG))$. Тогда  $t\in\J(FG)$ по предложению \ref{rad cent}.
\qed

\medskip\noindent
{\sc \hypertarget{mlt vp prf}{Доказательство}} \ref{mlt vp}. По условию $V=\sum V_i$, где $V_i$ --- неприводимые $A$-модули. Поскольку $\J(A)\le
\Ann(V_i)$ для всех $i$, имеем включение $\J(A)\le \Ann(V)$. С другой стороны, $A_V\cong A/\Ann(V)$
и, значит, алгебра $A_V$ изоморфна гомоморфному образу полупростой алгебры $A/\J(A)$, т.\,е. сама является
полупростой в силу следствия \ref{cor fact ss}.
\qed

\medskip\noindent
{\sc \hypertarget{mlt npr prf}{Доказательство}} \ref{mlt npr}. $(i)$ Из разложения \ref{thm wedd}$(i)$ и простоты алгебры
$A$ вытекает, что ненулевой является в точности одна компонента Веддерберна $\W_M(A^\circ)$, где $M\in \M(A)$.
А поскольку для любого модуля $M\in \M(A)$ справедливо $\W_M(A^\circ)\ne 0$ (как мы уже отмечали в начале доказательства теоремы \ref{thm wedd}),
имеем $|\M(A)|=1$.

$(ii)$ Обозначим $B=A_V$. Поскольку $B\cong A/\Ann(V)$ и $\J(A)\se\Ann(V)$, то алгебра
$B$ полупроста. Заметим, что $V$ --- неприводимый $B$-модуль и по теореме Веддерберна \ref{thm wedd}$(ii)$
имеет место разложение. $B=\W_V(B^\circ)\oplus \Ann_{B}(V)$. Но тогда $\Ann_{B}(V)=0$, т.\,к. $B$
является факторалгеброй алгебры $A$ по аннулятору $\Ann(V)$. Значит, $B$~--- простая алгебра. \qed

\medskip\noindent
{\sc \hypertarget{cor z ss prf}{Доказательство}} \ref{cor z ss}. $(i)$ Из теоремы Веддерберна \ref{thm wedd} следует, что
$$\dim_F A=\sum_{M\in
\M(A)} \dim_F \W_M(A^\circ),$$ где $\W_M(A^\circ)$ --- простые $F$-алгебры, и, кроме того, $M$ является
единственным, с точностью до изоморфизма, неприводимым $\W_M(A^\circ)$-модулем. По следствию \ref{cor
sim azp} получаем $\dim_F \W_M(A^\circ)=(\dim_F M)^2$. Значит, из равенства \ref{nmv} имеем $\n_M(A^\circ)=\dim_F M$.

$(ii)$ Прямое следствие $(i)$ и формулы \ref{vcong}.

$(iii)$ Получается путём вычисления размерностей над $F$ обеих частей соотношения $(ii)$.

$(iv)$ Из теоремы \ref{thm wedd} и предложения \ref{z pr} следует, что
$\Z(A)=\bigoplus_{M\in\M(A)}\Z(\W_M(A^\circ))$. По следствию \ref{cor sim azp} простая $F$-алгебра
$\W_M(A^\circ)$ изоморфна $\End_F(M)\cong \MM_n(F)$, где $n=\dim_F M$. Однако центр $\Z(\MM_n(F))$ одномерен,
поскольку, как мы отмечали в \ref{pr kol}, он состоит из скалярных матриц.
Значит, единица $e_{\mbox{}_M}$ алгебры $\W_M(A^\circ)$ является единственным, с точностью до умножения на скаляр,
элементом из $\Z(\W_M(A^\circ))$. Поэтому идемпотенты $e_{\mbox{}_M}$, $M\in\M(A)$, образуют базис $Z(A)$ и
размерность $\dim_F\Z(A)$ равна $|\M(A)|$.
\qed

\medskip\noindent
{\sc \hypertarget{int krs prf}{Доказательство}} \ref{int krs}. $(i)$ Пусть $\deg\X_i=n_i$ и $V_i=F^{n_i}$. Тогда из замечаний о
соответствии представлений и модулей следует, что $A$-модули $V_i$, $i=1,\ld,s$, с действием $va=v\X_i(a)$,
образуют полный набор попарно неизоморфных неприводимых $A$-модулей. Легко видеть, что $\Ker\X_i=\Ann(V_i)$.
Отсюда следует требуемое равенство.

$(ii)$ Требуемое  следует из \ref{pr pr}$(iii)$, и \ref{nep f}.
\qed

\medskip\noindent
{\sc \hypertarget{cent irr prf}{Доказательство}} \ref{cent irr}. Чтобы показать, что $M$ --- скалярная матрица, можно воспользоваться предложением \ref{rep
trace}, утверждающим, что $\Im \X =\MM_n(F)$. Приведём, однако, другое доказательство, использующее, по
существу, только лемму Шура.

Пусть $V=F^n$. Тогда $V$ является $A$-модулем с действием $va=v\X(a)$. Если $\vf$ --- линейное
преобразование пространства $V$, действующее по правилу $v\vf=vM$, то $\vf \in \End_A(V)$. Из
следствия
\ref{cor sh} вытекает, что  $\vf$ --- скалярное преобразование и $M$ --- скалярная
матрица.

Если $z\in Z(A)$, то матрица $\X(z)$ перестановочна с $\X(a)$ для всех $a\in A$ и поэтому
скалярна по только что доказанному. Если $A$ коммутативна, то $Z(A)=A$. Значит, любое
подпространство из $V$ будет $A$-инвариантным. Из неприводимости $V$ следует, что $\deg\X=1$.
\qed

\medskip\noindent
{\sc \hypertarget{pr phg prf}{Доказательство}} \ref{pr phg}$(iv)$.
Выберем в качестве базиса регулярного модуля $FG^\circ$ множество
элементов группы $G$ и пусть $\R_F$ --- соответствующее регулярное представление. Пусть $g\in G$ и
$\R_F(g)=(\a_{x,y})$, $x,y\in G$. Тогда $\a_{x,y}\ne 0$ если и только если $xg=y$, и в этом случае
$\a_{x,y}=1_F$. Поэтому
$$\r(g)=\big|\{x\in G\mid xg=x\}\big|\cdot1_F=\left\{\ba{ll}
|G|\cdot1_F,& g=1;\\
0_F,& g\ne 1.
\ea\right.$$

\medskip\noindent
{\sc \hypertarget{tnz har prf}{Доказательство}} \ref{tnz har}.  Пусть $g\in G$, $\X(g)=(\a_{ij})\in \MM_n(F)$, $\Y(g)=(\b_{rs})\in \MM_m(F)$.

$(i)$ Имеем
$$
(v_i\ot w_r)g=v_ig\ot w_rg=\sum_{j=1}^m\sum_{s=1}^n\a_{ij}\b_{rs}(v_j\otimes w_s)
$$
и, значит, матрица $(\X\ot\Y)(g)$ равна $\X(g)\ot \Y(g)$.

$(ii)$ Имеем
$$
\x(g)=\sum_{i=1}^m \a_{ii}, \qquad \psi(g)=\sum_{r=1}^n \b_{rr}.
$$
Пусть $\th$ --- характер представления $\X\ot \Y$. Ввиду $(i)$
$$
\th(g)=\tr (\X(g)\ot\Y(g))= \sum_{i=1}^m\sum_{r=1}^n
\a_{ii}\b_{rr}=\big(\sum_{i=1}^m\a_{ii}\big)\big(\sum_{r=1}^n\b_{rr}\big)=\x(g)\psi(g).
$$
Отсюда следует требуемое. \qed

\medskip\noindent
{\sc \hypertarget{sst imp prf}{Доказательство}} \ref{sst imp}. $(i)$ Очевидно.

$(ii)$ Это следует из транзитивности действия группы $G$ на компонентах системы
импримитивности и того, что $I$ --- стабилизатор одной из компонент.

$(iii)$ Множество $\{Wr_1,\ld,Wr_k\}$ совпадает с множеством $\{W_1,\ld,W_k\}$. Поскольку $Br_i$
--- базис подпространства $Wr_i$ и $W$ --- прямая сумма подпространств $W_i$, то
базис $W$ --- объединение базисов этих подпространств.

$(iv)$ Пусть $g\in G$ и блочная $k\times k$-матрица $\Y(g)$ состоит из клеток $A_{ij}$.
Покажем, что $A_{ij}=\X^\circ(r_i g r_j^{-1})$.

Для любого $r_i$ существует единственный $r_{i_0}$ такой, что $r_i g r_{i_0}^{-1}\in I$.
Поэтому  $\X^\circ(r_i g r_{i_0}^{-1})=\X(r_i g r_{i_0}^{-1})$, а $\X^\circ(r_i g
r_{i_0}^{-1})=\O$ для всех $j\ne i_0$. С другой стороны, для любого $b\in B$

$$(br_i)g=b(r_igr_{i_0}^{-1})r_{i_0}=b\X(r_igr_{i_0}^{-1})r_{i_0}.$$
Поэтому $A_{i,i_0}=\X(r_igr_{i_0}^{-1})=\X^\circ(r_igr_{i_0}^{-1})$ и $A_{ij}=\O=\X^\circ(r_i
g r_j^{-1})$ при $j\ne i_0$.

$(v)$ Из $(iii)$ и $(iv)$ следует, что у модуля $V$ имеется базис, в котором соответствующее представление
группы $G$ однозначно определяются представлением $\X$ группы $I$. Отсюда следует требуемое.
\qed

\medskip\noindent
{\sc \hypertarget{prop ind prf}{Доказательство}} \ref{prop ind}. $(i)$ Очевидно.

$(ii)$ Если $X=\{x_1,\ld,x_m\}$ --- набор представителей правых смежных классов $H$ по $K$ и
$Y=\{y_1,\ld,y_n\}$ --- набор представителей правых смежных классов $G$ по $H$, то $x_iy_j$ --- набор
представителей правых смежных классов $G$ по $K$. Запишем $U^H=\bigoplus_{i}U\ot x_i$. Тогда
$$
(U^H)^G=\bigoplus_j\left(\bigoplus_i\, U\ot x_i\right)\ot y_j=\bigoplus_{ij}\, U\ot x_i\ot y_j.
$$ С другой стороны,
$U^G=\bigoplus_{ij}U\ot x_iy_j$. Покажем, что отображение $\vf:(U^H)^G\to U^G$, действующее по правилу $(u\ot
x_i\ot y_j)\vf=u\ot x_iy_j$ для всех $u\in U$, $x_i\in X$, $y_j\in Y$, и продолженное по линейности на весь
модуль $(U^H)^G$, является изоморфизмом $FG$-модулей. Легко видеть, что $\vf$
--- изоморфизм векторных пространств. Пусть $g\in G$. Запишем $y_jg=hy_{j'}$ и $x_ih=kx_{i'}$
для подходящих $h\in H$, $k\in K$. Тогда
$$
(u\ot x_i\ot y_j)g\vf=(u\ot x_ih\ot y_{j'})\vf=(uk\ot x_{i'}\ot y_{j'})\vf= uk\ot x_{i'}y_{j'} =u\ot x_i
y_jg=(u\ot x_i\ot y_j)\vf g.
$$
Отсюда следует требуемое.

$(iii)$ Пусть $R$ --- набор представителей правых смежных классов $G$ по $H$. Тогда
$V^G=\bigoplus_{r\in R} V\ot r$. Очевидно, что $U^G=\bigoplus_{r\in R} U\ot r$ ---
подпространство в $V^G$, выдерживающее действие группы $G$. Если $U<V$ и $v\in V\setminus U$,
то $v\ot r\in V^G\setminus U^G$ для любого $r\in R$. Обратная импликация очевидна.

$(iv)$ Пусть $R$ как выше. Имеем
$$
(U\oplus V)^G=\bigoplus_{r\in R}\, (U\oplus V)\ot r= \bigoplus_{r\in R}\, (U\oplus 0)\ot r \ \oplus \ \bigoplus_{r\in R}\, (0\oplus V)\ot r.
$$
Оба прямых слагаемых в правой части инвариантны относительно $G$ и как $FG$-модули
естественно изоморфны $U^G$ и $V^G$, соответственно.
\qed

\medskip\noindent
{\sc \hypertarget{ind char prf}{Доказательство}} \ref{ind char}. Пусть $\X$ --- некоторое $F$-представление группы $H$ с характером $\x$.
Тогда $\x^G$ --- характер индуцированного представления $\X^G$. Если для любого $g\in G$ определить
$$
\X^\circ (g)=\left\{\ba{ll}
\X(g),& g\in H;\\
\O,& g\not\in H,
\ea\right.
$$
то $\x^\circ(g)=\tr\X^\circ(g)$ и из \ref{sst imp}$(iv)$ следует, что значение характера представления
$\X^G$ на элементе $g$ равно $\sum_{r\in R} \x^\circ(rgr^{-1})$, где $R$ --- произвольная система
представителей всех правых смежных классов $G$ по $H$.

\medskip\noindent
{\sc \hypertarget{ind kl prf}{Доказательство}} \ref{ind kl}. $(i)$ Отметим сначала, что для произвольных $x\in G$ и $h\in H$ имеет место равенство
$$
\vf^\circ(hxh^{-1})=\vf^\circ(x),\myeqno\label{fcir}
$$
поскольку $\vf\in \cf_F(H)$, и включение $hxh^{-1}\in H$ выполнено тогда и только тогда, когда $x\in H$.

Обозначим $\vf_R^G(g)$ значение правой части равенства \ref{fif}. Пусть $R_0$ --- другой
набор представителей всех правых смежных классов $G$ по $H$. Тогда найдутся такие элементы $h_r\in H$, где $r\in R$,
что $R_0=\{h_rr\mid r\in R\}$. В силу \ref{fcir} имеем

$$
\vf_{R_0}^G(g)=\sum_{r_0\in R_0} \vf^\circ(r_0^{\vphantom{-1}}gr_0^{-1})=\sum_{r\in R} \vf^\circ(h_rrgr^{-1}h_r^{-1})=
\sum_{r\in R} \vf^\circ(rgr^{-1})=\vf_{R}^G(g).
$$

$(ii)$ Это следует из того, что значения функции $\vf^G$ аддитивно зависят от значений функции $\vf$.

$(iii)$Поскольку группа $G$ является объединением $|H|$ попарно непересекающихся наборов $R_1,\ld,R_{|H|}$
представителей правых смежных классов $G$ по $H$, ввиду $(i)$ получаем
$$
\sum_{x\in G} \vf^\circ(xgx^{-1})=\sum_{i=1}^{|H|}\sum_{x\in R_i}\vf^\circ(xgx^{-1})=|H|\vf^G(g).
$$
Второе равенство следует из того, что значение $\vf^\circ(xgx^{-1})$ совпадает с $\vf(xgx^{-1})$, если $g\in H^x$,
и равно нулю в противном случае.

$(iv)$ Любой элемент из класса сопряжённости $g^G$ представ\'{и}м в виде $ygy^{-1}$, где $y\in G$, ровно
$|\C_G(g)|$ способами. Поэтому в силу $(ii)$
$$
|H|\vf^G(g)=\sum_{y\in G} \vf^\circ(ygy^{-1})= |\C_G(g)|\sum_{x\in g^G} \vf^\circ(x)= |\C_G(g)|\sum_{x\in
g^G\cap H} \vf(x). \myeqno\label{hfgg}
$$
Пусть $R_g$ --- набор представителей всех классов сопряжённости группы $H$, содержащихся в классе $g^G$.
Так как $\vf\in\cf_F(H)$, то правая часть в \ref{hfgg} равна
$$
|\C_G(g)|\sum_{x\in R_g} |x^H| \vf(x).
$$
По условию $|H|\ne 0$  в поле $F$. Поэтому  $|x^H|/|H|=1/|\C_H(x)|$. Разделив обе части на $|H|$, получим
требуемое.

$(v)$ Сначала заметим, что если $X$ и $Y$ --- полные  наборы представителей правых смежных классов $K$ по
$H$ и $G$ по $K$, соответственно, то множество
$$
R=\{\,xy\mid x\in X,\ y\in Y\,\}
$$
является полным  набором представителей правых смежных классов $G$ по $H$.

Из \ref{fif ud} следует, что
$$
\vf^G(g)=\sum_{\{r\in R\,\mid\, g\in H^r\}} \vf(rgr^{-1})=\sum_{\{x\in X,\,y\in Y\,\mid\, g\in H^{xy}\}} \vf(xygy^{-1}x^{-1}).
$$
Аналогично имеем
$$
(\vf^K)^G(g)=\sum_{\{y\in Y\,\mid\, g\in K^y\}} \vf^K(ygy^{-1})=
\sum_{\{y\in Y\,\mid\, g\in K^y\}}\ \sum_{\{x\in X\,\mid\, ygy^{-1}\in H^x\}} \vf(xygy^{-1}x^{-1}).
$$
Правые части полученных выражений совпадают, поскольку из включений $g\in H^{xy}$ и $x\in K$ вытекает, что $g\in K^y$.
Отсюда следует требуемое.

$(vi)$ Для произвольного $g\in G$ имеем
\begin{align*}
(\vf\,\psi_{H^{\vphantom{A^a}}})^G(g)=\sum_{\{r\in R\,\mid\, g\in H^r\}} (\vf\,\psi_{H^{\vphantom{A^a}}})(rgr^{-1})
&=\sum_{\{r\in R\,\mid\, g\in H^r\}} \vf(rgr^{-1})\psi_{H^{\vphantom{A^a}}}(rgr^{-1})\\
&=\Big(\sum_{\{r\in R\,\mid\, g\in H^r\}} \vf(rgr^{-1})\Big)\psi(g)=\vf^G(g)\psi(g)=(\vf^G\psi)(g),
\end{align*}
где мы воспользовались тем, что $\psi_{H^{\vphantom{A^a}}}(rgr^{-1})=\psi(g)$, поскольку $\psi\in\cf_F(G)$.

$(vii)$ Поскольку $G=HK$, множество $R$ представителей всех правых смежных классов $G$ по $H$ можно выбрать так, что
$R\se K$. При этом, как легко видеть,  $R$ будет также множеством представителей всех правых смежных классов $K$ по $H\cap K$.
Поэтому для любого $k\in K$ имеем
$$
(\vf^G)_{K^{\vphantom{A^a}}}(k)=\sum_{\{r\in R\,\mid\, k\in H^r\}} \vf(rkr^{-1})=
\sum_{\{r\in R\,\mid\, k\in (H\cap K)^r\}} \vf_{H\cap K^{\vphantom{A^a}}}(rkr^{-1})=(\vf_{H\cap K^{\vphantom{A^a}}})^K(k),
$$
так как для $r\in R$  включение $k\in H^r$ справедливо тогда и только тогда, когда $k\in (H\cap K)^r$.
\qed

\medskip\noindent
{\sc \hypertarget{fcn pr prf}{Доказательство}} \ref{fcn pr}. $(i)$ Пусть $h,k\in H$.
Если элементы $h^g,k^g\in H^g$ сопряжены элементом $l^g\in H^g$,
где $l\in H$, то $h$ и $k$ сопряжены с помощью $l$. Поэтому $\vf^g(h^g)=\vf(h)=\vf(k)=\vf^g(k^g)$.

$(ii)$ Имеем $$(\vf^g)^f(h^{gf})=\vf^g(h^g)=\vf(h)=\vf^{gf}(h^{gf})$$
для всех $h\in H$. Отсюда следует требуемое.

$(iii)$ Если $g\in H$, то $\vf^g=\vf$, поскольку $\vf$ --- классовая функция на $H$. Если $g\in C_G(H)$,
то $H^g=H$ и для любого $h\in H$ получаем $\vf^g(h^g)=\vf(h)=\vf(h^g)$, т.\,е. $\vf^g=\vf$.

$(iv)$ Пусть $R$ --- набор представителей всех правых смежных классов $G$ по $H$.
Тогда $R^x=\{r^x\mid r\in R\}$ --- набор представителей всех правых смежных классов $G$ по $H^x$.
Достаточно показать, что для любого $g\in G$ справедливо $(\vf^x)^G(g^x)=\vf^G(g)$.
В силу равенства \ref{fif ud} имеем
$$
(\vf^x)^G(g^x)=\sum_{\{r^x\in R^x\,\mid\, g^x\in (H^x)^{r^x}\}} \vf^x(r^xg^x(r^x)^{-1})=
\sum_{\{r\in R\,\mid\, g\in H^r\}} \vf(rgr^{-1})=\vf^G(g).
$$
\qed

\medskip\noindent
{\sc \hypertarget{spr mod prf}{Доказательство}} \ref{spr mod}. $(i)$ Выберем произвольный элемент $h^g\in H^g$, где $h\in H$. Тогда
$$(Wg)h^g=Wg(g^{-1}hg)=(Wh)g=Wg,$$ т.\,е. $Wg$ является $FH^g$-модулем. Выберем базис $w_1,\ld,w_n$
подпространства $W$. Тогда $w_1g,\ld,w_ng$ --- базис $Wg$. Пусть $\X$ и $\Y$ --- представления, соответствующие
$FH$-модулю $W$ и  $FH^g$-модулю $Wg$, соответственно, относительно выбранных базисов. Покажем, что $\Y=\X^g$.

Пусть $h\in H$ и $\X(h)=(\a_{ij})$.  Тогда $\X^g(h^g)=\X(h)=(\a_{ij})$ и
$$
w_i \X(h)=w_i h=\sum_j\a_{ij}w_j.
$$
Поэтому
$$
w_i g\Y(h^g)=(w_i g)h^g=(w_i h)g=(\sum_j\a_{ij}w_j)g=\sum_j\a_{ij}w_jg
$$
и, значит, $\Y(h^g)=(\a_{ij})=\X^g(h^g)$.

$(ii)$ Пусть $\X$ --- представление, соответствующее модулю $W$. По условию $FH^g$-модуль $M$
соответствует представлению $\X^t$ для некоторого $t\in G$ такого, что $H^g=H^t$ Однако, мы только что
показали, что $FH^t$-модуль $Wt$ также соответствует представлению $\X^t$. Значит, $M\cong Wt$ и
$tg^{-1}\in \N_G(H)$.

$(iii)$ Поскольку модули $M$ и $W$ изоморфны, то в подходящих базисах они соответствуют одному и тому же
представлению $\X$. Тогда $FH^g$-модули $Mg$ и $Wg$ соответствуют представлению $\X^g$ группы $H^g$ и,
значит, также изоморфны.
\qed

\medskip\noindent
{\sc \hypertarget{t mak prf}{Доказательство}} \ref{t mak}. Каждый двойной смежный класс $HtK$ является объединением правых смежных классов по
$H$. Пусть $K_t$ --- набор элементов из $K$ такой, что смежные классы $Htk$ попарно различны при различных $k\in K_t$, и для
которого $HtK=\cup_{k\in K_t}Htk$. Покажем, что $K_t$
является полным набором представителей правых смежных классов группы $K$ по подгруппе $H^t\cap K$.
Заметим, что два элемента $k_1,k_2\in K_t$ совпадают тогда и только тогда, когда $tk_1k_2^{-1}t^{-1}\in
H$, т.\,е. $k_1k_2^{-1}\in H^t$. С другой стороны $k_1k_2^{-1}\in K$. Значит, $k_1=k_2$ тогда и только
тогда, когда $k_1$ и $k_2$ лежат в одном правом смежном классе группы $K$ по $H^t\cap K$. Обратно, для
любого $k\in K$ существуют $h\in H$ и $k_0\in K_t$ такие, что $tk=htk_0$. Поэтому $kk_0^{-1}=h^t\in
H^t\cap K$, т.\,е. $K_t$ --- это полная система представителей $K$ по $H^t\cap K$.

Как и в доказательстве предложения \ref{ind mod}, представим $FG$-модуль $W^G$ в виде прямой
суммы
$$W^G=\bigoplus_{t\in T}\bigoplus_{k\in K_t}W\ot tk.$$
Легко видеть, что для любого $t\in T$ слагаемое
$$\bigoplus_{k\in K_t}W\ot tk\myeqno\label{opwtk}$$
является
подмодулем $FK$-модуля $(W^G)_K$.

С другой стороны, из \ref{spr mod}$(i)$ следует, что для любого $t\in T$ подпространство $Wt$ является
$FH^t$-модулем, а значит и $F(H^t\cap K)$-модулем $(Wt)_{H^t\cap K}$.
Поскольку $K_t$ --- это полная система представителей $K$ по $H^t\cap K$, имеет место равенство подпространств
$$((Wt)_{H^t\cap K})^K=\bigoplus_{k\in K_t}Wt\ot k. \myeqno\label{wtht}$$

Достаточно показать, что отображение $\vf$ из $FK$-модуля \ref{opwtk} в $FK$-модуль \ref{wtht}, определённое по
правилу $(w\ot tk)\vf=wt\ot k$ для всех $w\in W$, $k\in K_t$,  является изоморфизмом.

Пусть $x\in K$ и $k\in K_t$. Тогда существуют $h\in H$ и $k_0\in K_t$ такие, что $kx=h^tk_0$ и, значит,
$$\ba{l}
(w\ot tk)\vf x=(wt\ot k)x=wt\ot h^tk_0=wht\ot k_0, \\
(w\ot tk)x \vf=(w\ot htk_0)\vf=(wh\ot tk_0)\vf=wht\ot k_0,
\ea
$$
откуда следует требуемое. \qed

\medskip\noindent
{\sc \hypertarget{ext npr prf}{Доказательство}} \ref{ext npr}. Если $\R_F$ --- регулярное $F$-представление, то $(\R_F)^E=\R_E$ --- регулярное
$E$-представ\-ле\-ние. Для нахождения неприводимых компонент представления $\R_E$ нужно взять неприводимые
компоненты $\X_i$ представления $\R_F$ и затем найти неприводимые компоненты представлений $\X_i^E$ \big(это
следует из теоремы Жордана-Гёльдера \ref{thm zh-g}\big). Осталось заметить, что по предложению
\ref{irr rep} представление $\Y$ является неприводимой компонентой регулярного $E$-представления
$\R_E$. \qed

\medskip\noindent
{\sc \hypertarget{cor c q prf}{Доказательство}} \ref{cor c q}.
По следствию \ref{cor algz abs} поле $\ov{\QQ}$ является полем разложения для группы $G$.
Из теоремы \ref{thm mash} и предложения \ref{ssm crit} следует, что любое $\CC$-представление $\X$
группы $G$ вполне приводимо, т.\,е. $\X$ эквивалентно прямой сумме
неприводимых $\CC$-представлений, каждое из которых по предложению
\ref{ext abs} может быть записано над $\ov{\QQ}$.
Отсюда следует, что и само $\X$  может быть записано над $\ov{\QQ}$.
\qed

\medskip\noindent
{\sc \hypertarget{ab npr prf}{Доказательство}}  \ref{ab npr}. Воспользуемся соотношением
$$
|G|=\sum_{\x\in \irr(G)} \x(1)^2, \myeqno\label{auxg}
$$
где число слагаемых совпадает с $|\K(G)|$ ввиду \ref{ir eq kl}.

Если $G$ абелева, то $|\K(G)|=|G|$ и, значит, все слагаемые в \ref{auxg} равны $1$. Обратно, пусть
$\x(1)=1$ для всех $\x\in \irr(G)$. Тогда, чтобы правая часть в \ref{auxg} равнялась $|G|$
необходимо, чтобы число слагаемых было равно $|G|$, откуда следует $|\K(G)|=|G|$ и, значит,
$G$ абелева.
\qed

\medskip\noindent
{\sc \hypertarget{har prop prf}{Доказательство}} \ref{har prop}. $(i)$ Ограничение представления $\X$ на циклическую подгруппу $\la g\ra$ является вполне
приводимым по теореме Машке \ref{thm mash} и поэтому эквивалентно блочно-диагональному представлению с
неприводимыми компонентами по диагонали. Из \ref{ab npr} следует, что неприводимые компоненты одномерны и, в
частности, $\X(g)$ подобна диагональной матрице $\diag(\z_1,\ld,\z_n)$.

$(ii)$ Равенство $\z_i^k=1$ следует из $(i)$ и того, что $g^k=1$. Значит, $|\z_i|=1$ и
$\z_i^{-1}=\ov{\z_i}$.

$(iii)$ Поскольку следы подобных матриц равны, то из $(i)$ вытекает, что $\x(g)=\sum \z_i$.
Поэтому из $(ii)$ следует, что $|\x(g)|\le \sum |\z_i|=n$, а также включение $\x(g)\in \QQ_k$.
Так как $k$ делит $m$, имеем $\QQ_k\se \QQ_m$. Поскольку $\z_i\in \ov{\ZZ}$ и $\ov{\ZZ}$ --- кольцо,
получаем $\x(g)\in \ov{\ZZ}$.

$(iv)$ Если $\X(g)=\z\I_n$, то $\x(g)=n\z$. Но $|\z|=1$ в силу $(ii)$. Значит, $|\x(g)|=n$.
Обратно, пусть $|\x(g)|=n$, т.\,е. $\x(g)$ лежит на комплексной окружности радиуса $n$ с центром в нуле.
Так как $\x(g)=\sum \z_i$ --- сумма $n$ комплексных чисел с модулем $1$,
из геометрических соображений ясно, что равенство $|\x(g)|=n$ возможно лишь когда все $\z_i$ равны, т.\,е. $\X(g)$ --- скалярная матрица.

$(v)$ Если $\X(g)\in \ker\X$, то $\X(g)=\I_n$ и $\x(g)=n$. Обратно, пусть $\x(g)=n$. Из $(iv)$ следует, что
$\X(g)=\z\I_n$ для некоторого $\z\in \CC^\times$. Тогда $n=\x(g)=n\z$, т.\,е. $\z=1$ и $\X(g)\in \ker\X$.

$(vi)$ Так как $\X(g^{-1})=\X(g)^{-1}$, а $\X(g)$ подобна $\diag(\z_1,\ld,\z_n)$ в силу $(i)$, то матрица
$\X(g^{-1})$ подобна $\diag(\z_1^{-1},\ld,\z_n^{-1})$. Значит, из $(ii)$ следует, что $\x(g^{-1})=\ov{\x(g)}$.
\qed

\medskip\noindent
{\sc \hypertarget{reg har prf}{Доказательство}} \ref{reg har}. $(i)$ Обозначим через $\R$ регулярное представление группы $G$,
соответствующее модулю $\CC G^\circ$ в естественном базисе.\footnote{см. \ref{pr phg}$(iv)$.}
Тогда $\r$ --- характер представления $\R$ и для любого $g\in G$ имеем
$$
\R(g)=(a_{xy})_{x,y\in G}, \qquad a_{x,y}=\left\{\ba{ll}
1,& y=xg;\\
0,& y\ne xg;
\ea\right.
$$
В частности, $\R(g)$ --- единичная матрица  тогда и только тогда, когда
 $g=1$. Значит, $\ker\r=1$.

$(ii)$ См. \ref{pr phg}$(iv)$.\qed

\medskip\noindent
{\sc \hypertarget{kpr prf}{Доказательство}} \ref{kpr}. $(i)$ Это следует из \ref{har prop}$(v)$.

$(ii)$ Если $g\in \ker \t$ для всех $\t\in \irr(G)$ таких, что $n_\t>0$, то
$$\x(g)=\sum n_\t\t(g)=\sum n_\t\t(1)=\sum n_\t\deg \t=\deg\x=\x(1)$$
в силу $(i)$, и, значит,  $g\in \ker\x$. Обратно, пусть $g\in \ker\x$.
Тогда $$\sum n_\t \t(g) = \x(g)=\x(1)= \sum n_\t \t(1).$$
Но $|\t(g)|\le \t(1)$ в силу \ref{har prop}$(iii)$, откуда следует, что  $\t(g)=\t(1)$, т.\,е.  $g\in \ker \t$ для всех
$\t$ таких, что $n_\t>0$.

$(iii)$ Это следует из $(ii)$, если положить $\x=\r$ и воспользоваться утверждением \ref{fbern}$(i)$.

$(iv)$ В силу \ref{pr phg}$(vi)$, существует взаимно однозначное соответствие
между неприводимыми представлениями группы $G$, содержащими $N$ в своём ядре, и всеми неприводимыми представлениями
факторгруппы $G/N$. Из $(iii)$ вытекает, что $\cap_{\t \in\irr(G/N)} \ker\t=1$. Из этих замечаний следует требуемое.
\qed

\medskip\noindent
{\sc \hypertarget{ch gn prf}{Доказательство}} \ref{ch gn}. Утверждения $(i)$--$(iii)$ следуют из \ref{pr phg}$(vi)$,
определения точного характера и
взаимно однозначного соответствия между обыкновенными характерами и классами эквивалентности
обыкновенных представлений группы, см. \ref{ob vz}.

Докажем $(iv)$. Включение  $N\se\ker \t^G$ следует из \ref{ker indh} и поэтому характер
$\phantom{(\t^G)}\mathllap{\widetilde{\phantom{(t^G}\mkern10mu}\mathllap{(\t^G)}}$
определён. Для $g\in G$ обозначим $\wt g=g\vf$. Если $R$ --- полный набор представителей
правых смежных классов $G$ по $H$, то $\{\wt r\mid r\in R\}$ --- полный набор представителей правых смежных классов $\wt G$ по $\wt H$.
В силу \ref{fif ud} для любого $g\in G$ имеем
$$
\phantom{(\t^G)}\mathllap{\widetilde{\phantom{(t^G}\mkern10mu}\mathllap{(\t^G)}}
(\wt g)=\t^G(g)=\sum_{\{r\in R\,\mid\, g\in H^r\}} \t(rgr^{-1})=
\sum_{\{r\in R\,\mid\, \wt g\in \wt H^{\wt r}\}} \wt\t(\wt r\,\wt g\, \wt r^{\,-1})=\wt\t^{\,\wt G}(\wt g).
$$
Отсюда следует требуемое.
\qed

\medskip\noindent
{\sc \hypertarget{com cht prf}{Доказательство}} \ref{com cht}. $(i)$ Проверим, что
для любого $\l\in \irr(G)$ равенство $\deg\l=1$ справедливо тогда и только тогда, когда $G'\le \ker\l$.
В силу \ref{kpr}$(iv)$ отсюда будет следовать $(i)$.

Всякий линейный характер $\l$ группы $G$ является
гомоморфизмом из $G$ в абелеву группу $\CC^\times$. Поэтому $G'\le \ker\l$. Обратно,
если $\l\in \irr(G)$ и $G'\le \ker\l$, то в силу \ref{ch gn}$(ii)$, характер $\l$ однозначно
соответствует некоторому характеру $\wt\l\in \irr(G/G')$ той же степени.
Однако группа $G/G'$ абелева и все её неприводимые обыкновенные характеры линейны
ввиду \ref{ab npr}. Отсюда получаем, что $\deg\l=\deg\wt\l=1$.

$(ii)$ Мы показали в $(i)$, что множество характеров $\x\in \irr(G)$, содержащих $G'$ в своём ядре, совпадает с $\lin(G)$.
В частности, $\irr(G/G')=\lin(G/G')$ ввиду абелевости $G/G'$ (или можно сослаться на \ref{ab npr}). Поэтому
из \ref{ch gn}$(vi)$ следует, что отображение $\lin(G)\to \lin(G/G')$, действующее по правилу $\l \mapsto \wt\l$,
где $\l=\wt\l\circ\vf$ для естественного гомоморфизма $\vf:G\to G/G'$ является взаимно однозначным
соответствием. Это соответствие также является групповым гомоморфизмом, поскольку для
любых $\l_1,\l_2\in \lin(G)$ и $g\in G$ имеем
$$
\wt{\l_1\l_2}(g\vf)= (\l_1\l_2)(g) = \l_1(g)\l_2(g) = \wt\l_1(g\vf)\wt\l_2(g\vf).
$$
Отсюда следует $(ii)$.

$(iii)$ Из $(ii)$ ввиду абелевости факторгруппы $G/G'$ и утверждения \ref{ab npr} следует, что
$$
|\lin(G)|=|\lin(G/G')|=|\irr(G/G')|=|G/G'|,
$$
где последнее равенство вытекает из \ref{cor z fg}$(iii)$. Отсюда следует $(iii)$.
\qed

\medskip\noindent
{\sc \hypertarget{alg spr prf}{Доказательство}} \ref{alg spr}.
$(i)$ Если $g$ и $g^{-1}$ сопряжены, то $\x(g)=\x(g^{-1})=\ov{\x(g)}$ в силу \ref{har prop}$(iv)$, т.\,е.
$\x(g)\in \RR$. Обратно, пусть $g\in K$, и $g^{-1}\in K^\i$ где $K,K^\i\in
\K(G)$. Поскольку для любого $\x\in \irr(G)$ по условию $\x(g)=\ov{\x(g)}=\x(g^{-1})$, отсюда
следует, что столбцы таблицы характеров группы $G$, соответствующие классам $K$ и $K^\i$, совпадают. Из
невырожденности таблицы характеров вытекает, что $K=K^\i$ и, значит, элементы $g$ и $g^{-1}$ сопряжены.

$(ii)$ Пусть  $g$ и $g^k$ сопряжены для всех $k$ таких, что $(k,|g|)=1$. Возьмём $\x\in
\irr(G)$. По \ref{har prop}$(iii)$ имеем $\x(g)\in \QQ_m$, где $m=\exp G$.
В силу \ref{gal conj}$(iii)$ существует целое число $k$ взаимно простое с $m$ такое, что для любого
$\s\in
\Gal(\QQ_m,\QQ)$ выполнены равенства $(\x(g))^\s=\x^\s(g)=\x(g^k)$. Поскольку $k$ также
взаимно просто с $|g|$, то предположению отсюда вытекает, что $(\x(g))^\s=\x(g)$ и, значит, число
$\x(g)$ лежит в неподвижном подполе поля $\QQ_m$ относительно группы $\Gal(\QQ_m,\QQ)$. Из \ref{nor gal}$(ii)$
следует, что это подполе равно $\QQ$.
\normalmarginpar

Докажем обратное утверждение. Пусть $k\in\ZZ$, $(k,|g|)=1$ и элементы $g$ и $g^k$ лежат в
классах сопряжённости $K$ и $K'$, соответственно. В силу \ref{ks f}$(iii)$ существует (единственный)
автоморфизм $\s\in\Gal(\QQ_{|g|},\QQ)$, такой, что $\z^\s=\z^k$, где $\z\in\QQ_{|g|}$ ---
примитивный корень степени $|g|$ из $1$. Из
\ref{har prop} следует, что для любого $\x\in \irr(g)$ значение $\x(g)$ --- сумма степеней
числа $\z$. По условию $\x(g)\in \QQ$, т.\,е. $\x(g)=\x(g)^\s=\x(g^k)$. Как и выше отсюда
следует, что столбцы таблицы характеров группы $G$, соответствующие классам $K$ и $K'$,
совпадают и, значит, элементы $g$ и $g^k$ сопряжены.
%
%
%
%
%
\qed

\medskip\noindent
{\sc \hypertarget{sym int prf}{Доказательство}} \ref{sym int}. Если $g\in S_n$ и $(k,|g|)=1$,
то элементы $g$ и $g^k$ имеют одинаковое цикловое строение и,
значит, сопряжены в $S_n$. Из \ref{alg spr}$(ii)$ вытекает, что значения всех неприводимых характеров рациональны. C
другой стороны по \ref{har prop}$(iii)$ они являются целыми алгебраическими числами, и поэтому целые
рациональные в силу \ref{z cz}.
\qed

\medskip\noindent
{\sc \hypertarget{cor scl prf}{Доказательство}} \ref{cor scl}. $(i)$ -- $(ii)$ Утверждения следуют из того, что произведение
$(\vf,\x)_{\mbox{}_G}$ совпадает с коэффициентом при $\x$ в разложении $\vf$ в виде линейной комбинации
неприводимых характеров.

$(iii)$ -- $(iv)$ Пусть $\vf=\sum_{\x\in \irr(G)} a_\x \x$ и $\psi=\sum_{\x\in \irr(G)} b_\x
\x$ для неотрицательных целых чисел $a_\x$ и $b_\x$. Тогда
$(\vf,\psi)_{\mbox{}_G}=(\psi,\vf)_{\mbox{}_G}=\sum_{\x} a_\x b_\x$ --- неотрицательное целое число.
Поэтому $(\vf,\vf)_{\mbox{}_G}=\sum_\x \a_\x^2=1$ тогда и только тогда, когда один из коэффициентов
$a_\x$ равен $1$, а остальные --- $0$.
\qed

\medskip\noindent
{\sc \hypertarget{ocn pr prf}{Доказательство}} \ref{ocn pr}. $(i)$ Если $h$ пробегает всю
группу $H$, то  $xhx^{-1}$ тоже. Поэтому
$$(\vf^x,\psi^x)_{\mbox{}_{H^x}}=\frac{1}{|H^x|}\sum_{h^x\in H^x}\vf^x(h^x)\ov{\psi^x(h^x)}
                                =\frac{1}{|H|}\sum_{h\in H}\vf(h)\ov{\psi(h)}
                                =(\vf,\psi)_{\mbox{}_H}.$$

$(ii)$ Поскольку $(\th_H)^x=\th_{H^x}$, требуемое следует из $(i)$.

$(iii)$ Пусть $\X$ --- представление группы $H$ с характером $\vf$. Тогда $\X^x$, определённое
правилом $\X^x(h^x)=\X(h)$ для всех $h\in H$ также будет представлением $H^x$, характер которого
совпадает с $\vf^x$. Заметим, что образ группы $H$ относительно представления $\X$
совпадает с образом $H^x$ относительно $\X^x$ как множества матриц. Поэтому эти представления одновременно
приводимы или неприводимы (т.\,е. эквивалентны или неэквивалентны блочно-верхнетреугольному представлению).

Также отметим, что одновременная неприводимость (или приводимость) характеров $\vf$ и $\vf^x$ также следует
из $(i)$ и \ref{cor scl}$(iv)$.
\qed

\medskip\noindent
{\sc \hypertarget{vz fr prf}{Доказательство}} \ref{vz fr}.
В случае поля $\CC$ из \ref{ind kl}$(ii)$ вытекает, что
$$
\vf^G(g)=\frac{1}{|H|}\sum_{x\in G} \vf^\circ(xgx^{-1}).
$$
Поэтому
$$
(\vf^G,\psi)_{\mbox{}_G}=\frac{1}{|G|}\sum_{g\in
G}\vf^G(g)\ov{\psi(g)}=\frac{1}{|G|}\frac{1}{|H|}\sum_{g\in G}\sum_{x\in
G}\vf^\circ(xgx^{-1})\ov{\psi(g)}.
$$
Переобозначив $y=xgx^{-1}$ и заметив, что $\psi(g)=\psi(y)$, получим
$$
(\vf^G,\psi)_{\mbox{}_G}=\frac{1}{|G|}\frac{1}{|H|}\sum_{x\in G}\sum_{y\in
G}\vf^\circ(y)\ov{\psi(y)}=
\frac{1}{|H|}\sum_{y\in G}\vf^\circ(y)\ov{\psi(y)}=
\frac{1}{|H|}\sum_{y\in H}\vf(y)\ov{\psi(y)}=(\vf,\psi_{H^{\vphantom{A^a}}})_{\mbox{}_H}.
$$
\qed

\medskip\noindent
{\sc \hypertarget{act pairs prf}{Доказательство}} \ref{act pairs}. 
$(i)$ Предположим, что действие $G$ дважды транзитивно и пусть $(x_1,y_1), (x_2,y_2)\in X\times X$, где $x_1\ne y_1$ и $x_2\ne y_2$. В силу транзитивности,
$x_2=x_1g$ для некоторого $g\in G$. Заметим, что $y_2,y_1g\in X\setminus\{x_2\}$. Ввиду дважды транзитивности найдётся $h\in \St(x_2)$ такой, что $y_1gh=y_2$. Значит, 
$$
(x_1,y_1)gh = (x_2,y_1g)h = (x_2,y_2)
$$
и действие $G$ на парах транзитивно.

Обратно, пусть $G$ действует транзитивно на упорядоченных парах $(x,y)\in X\times X$, где $x\ne y$. Пусть $x\in X$ и $y_1,y_2\in X\setminus \{x\}$. Тогда найдётся $g\in G$ такой, что $(x,y_1)g=(x,y_2)$, т.\,е. $xg=x$ и $y_1g=y_2$. В частности, $g\in \St(x)$, действие $\St(x)$ на $X\setminus \{x\}$ транзитивно и, значит, $G$ действует дважды транзитивно.

Пусть $\chi^{(2)}$ --- подстановочный характер действия $G$ на упорядоченных парах $(x,y)\in X\times X$, где $x\ne y$. Такая пара $(x,y)$ 
неподвижна под действием $g\in G$ тогда и только тогда, когда $xg=x$ и $yg=y$. Значит,  число таких пар равно $\chi(g)^2$ минус число пар вида $(x,x)$ где $xg=x$, которых ровно $\chi(g)$. Таким образом,  $\chi^{(2)}=\chi^2-\chi$.

$(ii)$ Пусть $\chi^{\{2\}}$ --- подстановочный характер действия $G$ на множестве неупорядоченных пар $\{x,y\}\subseteq X$, где $x\ne y$.
Пара $\{x,y\}$ неподвижна под действием $g\in G$ в двух случаях.
Во-первых, когда $xg=x$ и $yg=y$. Число таких пар равно $\chi(g)(\chi(g)-1)/2$.
Во-вторых, когда $xg=y$ и $yg=x$, то есть $x$ и $y$ неподвижны под действием $g^2$, но не являются неподвижными под действием $g$. Число таких пар равно $(\chi(g^2)-\chi(g))/2$. Таким образом, 
$$
\chi^{\{2\}}(g)=\frac{1}{2}\chi(g)(\chi(g)-1)+\frac{1}{2}(\chi(g^2)-\chi(g))=
\frac{1}{2}(\chi^2(g)+\chi(g^2))-\chi(g).
$$

\medskip\noindent
{\sc \hypertarget{cor cl prf}{Доказательство}} \ref{cor cl}.
Для любого $\x\in \irr(G)$ из закона взаимности \ref{vz fr} получаем $n_\x=(\th^G,\x)=(\th,\x_{H^{\vphantom{A^a}}})$,
а по теореме Клиффорда \ref{tcl char} имеем $\x_{H^{\vphantom{A^a}}}=n_\x(\th_1+\ldots+\th_t)$, где
$\th_1,\th_2,\ld,\th_t$ --- полный набор различных характеров, сопряжённых с $\th$ элементами из $G$. Значит,
$$
(\th^G)_{H^{\vphantom{A^a}}}=\sum_{\x\in \irr(G)}n_\x \x_{H^{\vphantom{A^a}}}
= \left(\sum_{\x\in \irr(G)}n_\x^2\right) (\th_1+\ldots+\th_t).
$$
Взяв степень обеих частей этого равенства, получим
$$
|G:H|\,\deg\th =\left(\sum_{\x\in \irr(G)}n_\x^2\right)t\,\deg\th.
$$
Учитывая, что $|G:H|=|G:I||I:H|$ и $t=|G:I|$, получаем требуемое.
\qed

\medskip\noindent
{\sc \hypertarget{cent dif prf}{Доказательство}} \ref{cent dif}. $(i)$ В доказательстве предложения \ref{id dec} мы отмечали,
что $\X_\th(e_\x)=\O$ при $\x\ne\th$  и  $\X_\th(e_\th)=\I$. Отсюда следует требуемое.

$(ii)$ Пусть $\t:\Z(\CC G)\to \CC$ --- гомоморфизм $\CC$-алгебр. Поскольку центральные идемпотенты $e_\x$
образуют базис алгебры $\Z(\CC G)$, гомоморфизм $\t$ однозначно задаётся значениями $e_\x\t$, которые
являются идемпотентами из $\CC$, т.\,е. равны $0$ или $1$. В силу $(i)$ достаточно показать, что
существует $\x\in \irr(G)$ такой, что $e_\x \t=1$ и $e_\th\t=0$ для всех  $\th\ne\x$.
Так как $1_{\Z(\CC G)}\t=1_\CC$, то $e_\x \t=1$ хотя бы для одного $\x\in \irr(G)$.  Если существуют
различные $\x_1,\x_2\in \irr(G)$ такие, что $e_{\x_1}\t=e_{\x_2}\t=1$, то $(e_{\x_1}+e_{\x_2})\t=2$, вопреки
тому, что образ идемпотента $e_{\x_1}+e_{\x_2}$ под действием $\t$ должен быть идемпотентом из $\CC$.

$(iii)$ Пусть $z\in Z(\CC G)$. Поскольку центральные идемпотенты $e_\x$ образуют базис алгебры $Z(\CC G)$,
существуют константы $\a_\x\in \CC$ такие, что
$$
z=\sum_{\x\in \irr(G)}\a_\x e_\x.
$$
Подействовав гомоморфизмом $\om_\th$, где $\th\in\irr(G)$,  на обе части этого равенства, в силу $(i)$ получим
требуемое равенство $\om_\th(z)=\a_{\th}$.\qed

\medskip\noindent
{\sc \hypertarget{om ob prf}{Доказательство}} \ref{om ob}.
$(i)$ Взяв след обеих частей равенства $\X_\x(z)=\om_\x(z)\I$, получим $\x(z)=\om_\x(z)\x(1)$. В частности,
для $K\in \K(G)$ имеем
$$
\om_\x(\wh{K})\x(1)=\x(\wh{K})=\sum_{x\in K}\x(x)=|K|\x(x_{\mbox{}_K}).
$$

$(ii)$ Имеем $$\X_\x(za)=\X_\x(z)\X_\x(a)=\om_\x(z)\X_\x(a).$$
Переходя к следам и воспользовавшись $(i)$, получим требуемые равенства.

$(iii)$ Это следует из равенств
$$
\I=\X_\x(1)= \X_\x(z^{|z|})=\X_\x(z)^{|z|}=\om_\x(z)^{|z|}\I.
$$

\qed

\medskip\noindent
{\sc \hypertarget{str yav prf}{Доказательство}} \ref{str yav}.
Пусть $\x\in \irr(G)$. Подставим в равенство
$$
\om_\x(\wh{K})\om_\x(\wh{L})=\sum_{M\in \K(G)} a_{\mbox{}_{KLM}} \om_\x(\whm)
$$
выражения для $\om_\x(\wh{K})$, $\om_\x(\wh{L})$, $\om_\x(\whm)$ из \ref{ch val}. Тогда
$$
\frac{|K||L|}{\x(1)}\x(x_K)\x(x_L)=\sum_{M\in \K(G)}  a_{\mbox{}_{KLM}}|M| \x(x_M).
$$
Пусть $N\in \K(G)$. Умножим обе части
последнего равенства на $\ov{\x(x_{\mbox{}_N})}$, просуммируем по $\x$ и воспользуемся вторым
соотношением ортогональности \ref{vtor ort}. Имеем

\begin{align*}
|K||L|\sum_{\x\in
\irr(G)}\frac{\x(x_{\mbox{}_K})\x(x_{\mbox{}_L})\ov{\x(x_{\mbox{}_N})}}{\x(1)}=&
\sum_{M\in \K(G)}a_{\mbox{}_{KLM}}|M| \sum_{\x\in \irr(G)}  \x(x_{\mbox{}_M})\ov{\x(x_{\mbox{}_N})}\\
\mbox{}=\sum_{M\in \K(G)}& a_{\mbox{}_{KLM}} |M|\cdot \d_{\mbox{}_{M,N}} |\C_G(x_{\mbox{}_M})|=
a_{\mbox{}_{KLN}}|N||\C_G(x_{\mbox{}_N})|.
\end{align*}
Учитывая, что $|N||\C_G(x_{\mbox{}_N})|=|G|$ и переобозначая $N$ через $M$, получаем требуемый
результат.
\qed

\medskip\noindent
{\sc \hypertarget{m int z prf}{Доказательство}} \ref{m int z}.
Включение $p\ZZ\se M\cap \ZZ$ очевидно. Покажем обратное включение.
Допустим, напротив, что $m\in M\cap \ZZ$ и $p\nmid m$. Тогда $ap+bm=1$ для некоторых $a,b\in
\ZZ$. Отсюда следует, что $1\in M$. Противоречие.
\qed

\medskip\noindent
{\sc \hypertarget{lem ppt prf}{Доказательство}} \ref{lem ppt}.
Пусть $|g|=n=n_pn_{p'}$, где $n_p$ --- максимальная степень $p$, делящая $n$. Поскольку $n_p$ и $n_{p'}$ взаимно просты,
существуют $a,b\in \ZZ$ такие, что $1=an_p+bn_{p'}$. Тогда $g=g^{an_p+bn_{p'}}=g_{p'}g_p$, где $g_{p'}=g^{an_p}$ и $g_p=g^{bn_{p'}}$.
Поскольку $g^{n_p}$ является $p$-регулярным, $g_{p'}$ тоже $p$-регулярен. Поскольку $g^{n_{p'}}$ является $p$-элементом, $g_p$ --- тоже $p$-элемент.
Элементы $g_p$ и $g_{p'}$ перестановочны, т.\,к. являются степенями $g$. По этой же причине всякий элемент, перестановочный с $g$ перестановочен с
$g_p$ и $g_{p'}$. Если имеется другое представление $g=\ov{g}_{p'}\ov{g}_p$, где $\ov{g}_{p'}$ и $\ov{g}_p$ --- перестановочные $p'$- и $p$-элементы,
то эти элементы также являются степенями $g$ и поэтому элемент $\ov{g}_{p'}g_{p'}^{-1}=\ov{g}_pg_p^{-1}$ единичный,
поскольку одновременно является $p'$- и $p$-элементом. Отсюда следует единственность $g_p$ и $g_{p'}$.
\qed

\medskip\noindent
{\sc \hypertarget{br prop prf}{Доказательство}} \ref{br prop}. $(i)$ -- $(ii)$ Оба утверждения следуют из того, что для всякого $g\in
G_{p'}$ матрица $\X(g)$ имеет тот же набор характеристических корней, что и любая сопряжённая
с ней матрица.

$(iii)$ Мы уже отмечали, что матрица $\X(g)$ подобна $\diag(\z_1^\st,\ld,\z_n^\st)$, где $\z_1,\ld,\z_n\in U$.
Поэтому $\X(g^{-1})=\X(g)^{-1}$ подобна $\diag((\z_1^{-1})^\st,\ld,(\z_n^{-1})^\st)$. Поскольку
$\z_i^{-1}=\ov{\z_i}$, отсюда следует, что
$$\vf(g^{-1})=\sum\ov{\z_i}=\ov{\sum\z_i}=\ov{\vf(g)}.$$

$(iv)$ По \ref{pr phg}$(vii)$, мы можем считать, что $\X^*(g)=\X(g)^\i$ для всякого $g\in G$. Пусть $g\in G_{p'}$ и
$\z_1^\st,\ld,\z_n^\st\in F$ --- характеристические корни матрицы $\X(g)$ для подходящих $\z_1,\ld,\z_n\in U$.
Тогда  характеристические корни для $\X(g)^\i$ равны $(\z_1^{-1})^\st,\ld,(\z_n^{-1})^\st$, откуда следует, что
значение брауэрова характера представления $\X^*$ на элементе $g$ равно
$$\sum\z_i^{-1}=\sum\ov{\z_i}=\ov{\sum\z_i}=\ov{\vf(g)}=\ov{\vf}(g).$$

$(v)$ Это следует из того, что $H_{p'}\se G_{p'}$ и ограничение представления $\X$ на $H$
будет $F$-представлением группы $H$.
\qed

\medskip\noindent
{\sc \hypertarget{aut br prf}{Доказательство}} \ref{aut br}.
Пусть $\r$ --- автоморфизм Фробениуса поля $F$,
действующий по правилу $\a^\r=\a^p$ для любого $\a\in F$. Пусть $\s$ --- автоморфизм группы
$U$, являющийся <<поднятием>> автоморфизма $\r$, т.\,е. для произвольного $\xi\in U$ выполнено
соотношение $(\xi^\st)^\r=(\xi^\s)^\st$. Отсюда следует, что $\xi^\s=\xi^p$ и поэтому $\s$
однозначно продолжается до автоморфизма из $\Gal(\QQ(\xi),\QQ)$, который также обозначим через
$\s$. Таким образом $\s$ определён на всех элементах поля $\QQ(U)$. В частности,
$\a_1^\s+\a_2^\s=(\a_1+\a_2)^\s$ для любых $\a_1,\a_2\in \QQ(U)$.

Дальнейшее рассуждение повторяет идею из доказательства \ref{mm}$(i.1)$. Пусть $\vf\in
\iBr(G)$ и $\X$ --- $F$-представление с брауэровым характером $\vf$. Если
$\z_1^\st,\ld,\z_s^\st$ --- характеристические значения матрицы $\X(g)$, где $\z_i\in U$, то
$\vf(g)=\z_1+\ld+\z_s$. С другой стороны характеристическими значениями матрицы $\X^\r(g)$ являются
$(\z_1^\st)^\r,\ld,(\z_s^\st)^\r$. По определению автоморфизма $\s$ имеем $(\z_i^\st)^\r=(\z_i^\s)^\st$, откуда
следует, что значение брауэрова характера $F$-представления $\X^\r$ равно
$$\z_1^\s+\ld+\z_s^\s=(\z_1+\ld+\z_s)^\s=\vf(g)^\s=\vf^\s(g).$$

Поскольку $\z^{\s^i}=\z^{p^i}$ для любого $i$, где $\z\in U$ --- примитивный корень степени
$m_{p'}$ из $1$, то порядок $\s$ как элемента из $\Gal(\QQ_{m_{p'}},\QQ)$ равен порядку $p$ по
модулю $m_{p'}$.
\qed

\medskip\noindent
{\sc \hypertarget{hom mod prf}{Доказательство}} \ref{hom mod}.
В случае, когда $\eta\in\iBr(G)$, утверждение следует из того факта, что $F$-представле\-ние
$\X:FG\to \MM_n(F)$ с брауэровым характером $\eta$ является гомоморфизмом алгебр.

Пусть $\eta\in \irr(G)$. Поскольку отображение $\l_\eta$ является гомоморфизмом колец, достаточно показать,
что его действие на базисных элементах из $\Z(FG)$ перестановочно с умножением на скаляры из $F$. Для
любых $K\in \K(G)$ и $\a^\st\in F$, где $\a\in R$, имеем
$$
\l_\eta(\a^\st\wh K)=\om_\x(\a\wh K)^\st=\a^\st\om_\x(\wh K)^\st.
$$
\qed

\medskip\noindent
{\sc \hypertarget{com t prf}{Доказательство}} \ref{com t}. Пусть $z\in \Z(\wt{Z} G)$. Тогда $z=\sum_{K\in \K(G)} \a_{\mbox{}_K}\wh K$ для
некоторых $\a_K\in \wt{Z}$. Из \ref{dz} получаем
$$
(\om_\x(z))^\st=\sum_{K\in \K(G)} (\a_{\mbox{}_K})^\st\om_\x(\wh K)^\st=
\sum_{K\in \K(G)} (\a_{\mbox{}_K})^\st\l_\x(\wh K)=
\l_\x\Big(\sum_{K\in \K(G)} (\a_{\mbox{}_K})^\st\wh K\Big)=\l_\x(z^\st)
$$
в силу того, что $\om_\x:\Z(\wt{Z} G)\to \wt{Z}$ --- гомоморфизм $\wt{Z}$-алгебр,
$\l_\x:\Z(FG)\to F$ --- гомоморфизм $F$-алгебр и $^\st$~--- гомоморфизм колец.
\qed

\medskip\noindent
{\sc \hypertarget{id os prim prf}{Доказательство}} \ref{id os prim}.
Достаточно показать, что всякий идемпотент $f$ алгебры $\Z( \wt{Z} G)$ является суммой некоторых $f_{\mbox{}_B}$.
Поскольку $\wt{Z}$ --- подкольцо в $\CC$, всякий идемпотент из $\Z( \wt{Z} G)$ лежит в
$\Z( \CC G)$, см. \ref{zsr}, и значит, является суммой некоторых $e_\x$. (Заметим при этом, что
сами $e_\x$ могут не лежать в $\Z( \wt{Z} G)$.) Из предложения \ref{min un} следует, что такой
идемпотент является суммой некоторых $f_{\mbox{}_B}$.
\qed

\medskip\noindent
{\sc \hypertarget{rtg prf}{Доказательство}} \ref{rtg}. $(i)$ Это следует из \ref{ids}$(iii)$ и \ref{fbz}, поскольку
$$\sum_{B\in \bl(G)}f_{\mbox{}_B}=\sum_{\x\in \irr(G)}e_{\x}=1$$
и идемпотенты $f_{\mbox{}_B}$ попарно ортогональны.

$(ii)$ Это следует из \ref{ids}$(iv)$.
\qed

\medskip\noindent
{\sc \hypertarget{bl dec prf}{Доказательство}} \ref{bl dec}. В силу \ref{dek lin ind}$(ii)$ столбцы матрицы $D$ линейно независимы.
Поэтому из равенства \ref{zv22}
видно, что линейно независимы столбцы всех матриц разложения $D_{B_i}$.
Так как число столбцов у $D_{B_i}$ совпадает с $|\iBr(B_i)|$, а число строк --- с $|\irr(B_i)|$, то
отсюда следуют оба утверждения.
\qed

\medskip\noindent
{\sc \hypertarget{def ex prf}{Доказательство}} \ref{def ex}.
Легко видеть, что для любого $B\in \bl(G)$ идеал $\Z(W_B)$ содержится в $\Z_P(FG)$ для
$P\in \Syl_p(P)$, т.\,к. в этом случае $\Z_P(FG)=\Z(FG)$. Отсюда следует существование требуемого
класса $p$-подгрупп.

Пусть теперь $P$ и $Q$ минимальные по включению $p$-подгруппы, удовлетворяющие  условиям $\Z(W_B)\se \Z_P(FG)$
и  $\Z(W_B)\se \Z_Q(FG)$, соответственно. Тогда в силу \ref{zp id}$(ii)$
$$e_{\mbox{}_B}\in \Z(W_B)\se \Z_P(FG)\cap\Z_Q(FG)=\sum_{x,y\in G}\Z_{P^x\cap Q^y}(FG).$$
Из  \ref{bl z}$(vii)$ следует, что найдутся $x,y\in G$ такие, что $e_{\mbox{}_B}\in \Z_{P^x\cap Q^y}(FG)$.
Значит,
$$\Z(W_B)=e_{\mbox{}_B}\Z(FG)\se \Z_{P^x\cap Q^y}(FG)=\Z_{P^{xy^{-1}}\cap Q}(FG)=\Z_{P\cap Q^{yx^{-1}}}(FG).$$
В силу минимальности $p$-подгрупп $P$ и $Q$ имеем $Q\se P^{xy^{-1}}$ и $P\se Q^{yx^{-1}}$, откуда следует,
что $P$ и $Q$ сопряжены.
\qed

\medskip\noindent
{\sc \hypertarget{esg prf}{Доказательство}} \ref{esg}.
Легко видеть, что $\b_{\mbox{}_P}(\wh K)=0$ тогда и только тогда,
когда $K\cap C=\varnothing$, что эквивалентно существованию в $K$ элемента $x$, централизующего $P$, т.\,е.
включению $P\le S$ для некоторой $S\in \Syl_p(\C_G(x))\se \D(K)$, которое означает, что $P\le_G\d(K)$.
\qed

\medskip\noindent
{\sc \hypertarget{trn prf}{Доказательство}} \ref{trn}.
По определению отображения $\l_{b^N}^{\, G}$ для любого $K\in \K(G)$ имеем
$$
\l_{b^N}^{\, G}(\wh K)=\l_{b^N}\left(\sum_{x\in K\cap N}x\right)=\l_b^{\,N}\left(\sum_{x\in K\cap N}x\right)
=\l_b\left(\sum_{x\in K\cap N\cap H}x\right)=\l_b\left(\sum_{x\in K\cap H}x\right)=\l_b^{\, G}(\wh K).
$$
Значит, отображения  $\l_{b^N}^{\, G}$ и $ \l_b^{\, G}$ совпадают и поэтому либо одновременно являются, либо
не являются гомоморфизмами $F$-алгебр. В первом случае блоки $b^G$ и $(b^N)^G$ определены  и совпадают.
\qed

\medskip\noindent
{\sc \hypertarget{cor p rad prf}{Доказательство}} \ref{cor p rad}.
 Пусть $P=\d(B)$ и $N=\N_G(P)$. Поскольку $B\in \bl(G|P)$, из теоремы \ref{thm br 1}$(i)$
следует, что множество $\bl(N|P)$ непусто. В силу \ref{op def} получаем $\OO_p(\N_G(P))\le P$.
С другой стороны $P\nor N$ и, значит, $P\le\OO_p(\N_G(P))$. Поэтому $\OO_p(\N_G(P))=P$.
\qed

\medskip\noindent
{\sc \hypertarget{rbp prf}{Доказательство}} \ref{rbp}.
$(i)$ В силу \ref{kl sum}$(ii)$ базис алгебры $\Z(FG)$ состоит из
классовых сумм. Поэтому достаточно проверить, что для произвольных $K,L\in \K(K)$ справедливо
равенство $\vr(\wh K\wh L)=\vr(\wh K)\vr(\wh L)$. Это легко установить, используя
\ref{kl sum}$(iv)$, \ref{e b}$(i)$, и тот факт, что $\l_{\mbox{}_B}$ --- гомоморфизм $F$-алгебр.

$(ii)$ Это следует из \ref{e b}$(iii)$.

$(iii)$ Это следует из $(ii)$.

$(iv)$ В силу линейной независимости идемпотентов $e_{\mbox{}_B}$, элемент $x\in\Z(FG)$
лежит в ядре $\vr$ тогда и только тогда, когда $\l_{\mbox{}_B}(x)=0$ для всех $B\in \bl(G)$ и,
значит, лежит в $\J(\Z(FG))$ ввиду \ref{jzfg}.
\qed

\medskip\noindent
{\sc \hypertarget{svp prf}{Доказательство}} \ref{svp}.
$(i)$ -- $(iv)$ следуют из определения $p$-сечения.

$(v)$ Если $x\in \SS_H(g)$, существует элемент $h\in H$ такой, что $(x_p)^h=g_p=g$, т.\,е.
$x^h=(x_px_{p'})^h=gc$, где $c=(x_{p'})^h\in C$ --- $p$-регулярный элемент. Тогда найдётся элемент $k\in
H$ такой, что $c^k=t$ для подходящего $t\in T$. Значит, $x^{hk}=gc^k=gt$, т.\,е. $x\in (gt)^H$.
Обратно, для любого $t\in T$ имеем $gt\in\SS_H(g)$, поскольку $(gt)_p=g$. Тогда из $(iii)$ следует, что
$(gt)^H\se \SS_H(g)$ и, значит, имеет место обратное включение.
\qed


\medskip\noindent
{\sc \hypertarget{lem alg cn prf}{Доказательство}} \ref{lem alg cn}.
$(i)$ Пусть алгебраическое число $\a$ является корнем многочлена $f(x)=a_0x^n+ a_1
x^{n-1}+ \ld +a_n$, где $a_0,a_1,\ld,a_n\in \QQ$ и $a_0\ne 0$. Умножив, если необходимо, многочлен $f(x)$ на
наименьшее общее кратное знаменателей своих коэффициентов, можно считать, что $a_0,a_1,\ld,a_n\in
\ZZ$. Тогда
$$a_0^{n-1}f(x)=(a_0x)^n+a_1(a_0x)^{n-1}+a_2a_0(a_0x)^{n-2}+\ld+a_na_0^{n-1}.$$
Поскольку $a_0^{n-1}f(\a)=0$, мы видим, что число $a_0\a$ является целым алгебраическим для
целого рационального числа $a_0$, что и требовалось доказать.

$(ii)$ Пусть $\om$ --- целое алгебраическое число, являющееся корнем многочлена $f(x)=x^n+ b_1 x^{n-1}+\ld
+b_n, $ где $b_1,b_2,\ld,b_n\in \ZZ$, и пусть $0\ne m\in \ZZ$. Разделив обе части равенства $f(\om)=0$ на $m^n$
легко проверить, что $\om/m$ является корнем многочлена
$$y^n+\frac{b_1}{m}y^{n-1}+\frac{b_2}{m^2}y^{n-2}+\ld+\frac{b_n}{m^n},$$
и, значит, является алгебраическим.
\qed

\medskip\noindent
{\sc \hypertarget{kon oc prf}{Доказательство}} \ref{kon oc}.
Пусть $R$ --- конечная область целостности и $0\ne r\in R$. Рассмотрим последовательные степени
$r,r^2,r^3,\ld $ В силу конечности $R$ существуют $m\le n$ такие, что $r^m=r^n$. Тогда $r^m(1-r^{n-m})=0$. Так
как $r^m\ne 0$, отсюда вытекает $r^{n-m}=1$, т.\,е. элемент $r$ обратим.
\qed

\medskip\noindent
{\sc \hypertarget{lem int z prf}{Доказательство}} \ref{lem int z}. Пусть $\a\in A$ и $\a\ne 0$. Существуют $a_i\in \ZZ$ такие, что
$\a^m+a_1\a^{m-1}+\ld+a_m=0$, причём можно считать, что $a_m\ne 0$. Но тогда $a_m$ --- ненулевой элемент из
$A\cap\ZZ$. \qed

\medskip\noindent
{\sc \hypertarget{cor neter prf}{Доказательство}} \ref{cor neter}. Пусть $A_1\se A_2\se\ld$ --- цепочка идеалов из $D$. Так как факторкольцо $D/A_1$ конечно
в силу \ref{fak fin}, то существует лишь конечное число различных идеалов, содержащих $A_1$. \qed

\medskip\noindent
{\sc \hypertarget{cor pr max prf}{Доказательство}} \ref{cor pr max}. Если $P$ --- простой идеал в $D$, то в силу \ref{fak fin} факторкольцо
$D/P$ --- конечная область целостности. Из \ref{kon oc} вытекает, что $D/P$ --- поле, а потому
идеал $P$ максимален. \qed

\newpage

\section{Таблицы характеров некоторых групп} \label{pril tch}

\bigskip
\noindent
$G=S_3$

$|G|=6=2\cdot 3$

$$
\XX(S_3)\qquad\qquad
  \begin{array}{c|rrr}
       K                   & 1a      & 2a           & 3a       \\
        |K|                   & 1      & 3           & 2       \\
       |\C_G(x_{\mbox{}_K})|    & 6      & 2           & 3       \\
     \hline
    \x_1^{\vphantom{A^A}} & 1      & 1  & 1\\
    \x_2                  & 1      & -1            & 1\\
    \x_3                  & 2      & .  & -1
  \end{array}
$$

$$
\Phi_2(S_3)\qquad\qquad
  \begin{array}{c|rr}
     K                    & 1a       &  3a       \\
    |K|                   & 1        & 2       \\
  |\C_G(x_{\mbox{}_K})|    & 6        & 3       \\

     \hline
    \vf_1^{\vphantom{A^A}} & 1      & 1\\
    \vf_2                  & 2      & -1
  \end{array}\qquad\qquad
  \begin{array}{c|rr}
  D & \vf_1        & \vf_2       \\
     \hline
    \x_1^{\vphantom{A^A}} & 1      & .\\
    \x_2   & 1      & .\\
    \x_3   & .      & 1
  \end{array}\qquad\qquad
  C=\left(
      \begin{array}{cc}
        2& . \\
        . & 1
      \end{array}
    \right)
$$

$$
\Phi_3(S_3)\qquad\qquad
  \begin{array}{c|rr}
     K                    & 1a       &  2a       \\
    |K|                   & 1        & 3       \\
  |\C_G(x_{\mbox{}_K})|    & 6        & 2       \\

     \hline
    \vf_1^{\vphantom{A^A}} & 1      & 1\\
    \vf_2                  & 1      & -1
  \end{array}\qquad\qquad
  \begin{array}{c|rr}
  D & \vf_1        & \vf_2       \\
     \hline
    \x_1^{\vphantom{A^A}} & 1      & .\\
    \x_2   & .      & 1\\
    \x_3   & 1      & 1
  \end{array}\qquad\qquad
  C=\left(
      \begin{array}{cc}
        2& 1 \\
        1 & 2
      \end{array}
    \right)
$$
\hrulefill

\bigskip
\noindent
$G=A_4$

$|G|=12=2^2\cdot 3$

$$
\XX(A_4)\qquad\qquad
  \begin{array}{c|rrrr}
       K                   & 1a      & 2a           & 3a      & 3b \\
        |K|                   & 1      & 3           & 4      & 4 \\
       |\C_G(x_{\mbox{}_K})|    & 12      & 4           & 3     & 3  \\
     \hline
    \x_1^{\vphantom{A^A}} & 1      & 1  & 1  & 1\\[2pt]
    \x_2                  & 1      & 1 & \z  & \ov\z\\[2pt]
    \x_3                  & 1      & 1  & \ov\z & \z\\[2pt]
    \x_4                  & 3      & -1  & . & .
  \end{array}\qquad\qquad
  \begin{array}{c}
    \mbox{Иррациональные величины} \\[10pt]
    \z=\displaystyle\frac{-1+i\sqrt3}{2^{\vphantom A}}=e^{2\pi i/3}
  \end{array}
$$

$$
\Phi_2(A_4)
\qquad\qquad
  \begin{array}{c|rrr}
     K                    & 1a       &  3a  & 3b      \\
    |K|                   & 1        & 4    & 4     \\
  |\C_G(x_{\mbox{}_K})|    & 12        & 3     & 3    \\

     \hline
    \vf_1^{\vphantom{A^A}} & 1      & 1    &1 \\[2pt]
    \vf_2                  & 1      & \z   &\ov\z\\[2pt]
    \vf_3                  & 1      & \ov\z &\z\\[2pt]
  \end{array}\qquad\qquad
  \begin{array}{c|rrr}
  D & \vf_1        & \vf_2  &\vf_3     \\
     \hline
    \x_1^{\vphantom{A^A}} & 1   & . & .\\
    \x_2 & .   & 1 & .\\
    \x_3 & .   & . & 1\\
    \x_4 & 1   & 1 & 1
\end{array}\qquad\qquad
  C=\left(
      \begin{array}{ccc}
        2& 1 & 1 \\
        1 & 2 & 1 \\
        1 & 1 & 2
      \end{array}
    \right)
$$

$$
\Phi_3(A_4)\qquad\qquad
  \begin{array}{c|rr}
     K                    & 1a       &  2a       \\
    |K|                   & 1        & 3       \\
  |\C_G(x_{\mbox{}_K})|    & 12        & 4       \\

     \hline
    \vf_1^{\vphantom{A^A}} & 1      & 1\\
    \vf_2                  & 3      & -1\\
  \end{array}\qquad\qquad
  \begin{array}{c|rrr}
  D & \vf_1        & \vf_2      \\
     \hline
    \x_1^{\vphantom{A^A}} & 1   & . \\
    \x_2 & 1   & . \\
    \x_3 & 1   & . \\
    \x_4 & .   & 1
\end{array}\qquad\qquad
  C=\left(
      \begin{array}{cc}
        3& . \\
        . & 1
      \end{array}
    \right)
$$
\newpage

\noindent
$G=S_4$

$|G|=24=2^3\cdot 3$

$$
\XX(S_4)\qquad\qquad
  \begin{array}{c|rrrrr}
       K                   & 1a      & 2a    &2b   & 3a     & 4a  \\
        |K|                   & 1      & 6    & 3     & 8     & 6   \\
       |\C_G(x_{\mbox{}_K})|    & 24      & 4   & 8     & 3     & 4   \\
     \hline
    \x_1^{\vphantom{A^A}} & 1   & 1  & 1  & 1 & 1\\
    \x_2                  & 1   & -1 & 1  & 1 & -1\\
    \x_3                  & 2   &  .  & 2  & -1  & . \\
    \x_4                  & 3   &  1 & -1  & . & -1 \\
    \x_5                  & 3   &  -1 & -1 & .  & 1
  \end{array}
$$

$$
\Phi_2(S_4)
\qquad\qquad
  \begin{array}{c|rr}
     K                    & 1a       &  3a     \\
    |K|                   & 1        & 8      \\
  |\C_G(x_{\mbox{}_K})|    & 24        & 3     \\

     \hline
    \vf_1^{\vphantom{A^A}} & 1      & 1     \\
    \vf_2                  &  2     & -1
  \end{array}\qquad\qquad
  \begin{array}{c|rr}
  D & \vf_1        & \vf_2      \\
     \hline
    \x_1^{\vphantom{A^A}} &  1  & .  \\
    \x_2 &  1  & . \\
    \x_3 &  .  & 1 \\
    \x_4 &  1  & 1 \\
    \x_5 &  1  & 1
\end{array}\qquad\qquad
  C=\left(
      \begin{array}{cc}
         4 & 2   \\
         2 & 3
      \end{array}
    \right)
$$

$$
\Phi_3(S_4)\qquad\qquad
  \begin{array}{c|rrrr}
     K                    & 1a   &  2a   & 2b & 4a   \\
    |K|                   & 1    & 6     &  3 & 6    \\
  |\C_G(x_{\mbox{}_K})|    & 12   & 4    &  8 &  4   \\

     \hline
    \vf_1^{\vphantom{A^A}} & 1   & 1  & 1 & 1\\
    \vf_2                  &  1   & -1   & 1  & -1 \\
    \vf_3                  &  3   &  1  & -1  & -1\\
    \vf_4                  &  3   &  -1  & -1  & 1 \\
  \end{array}\qquad\qquad
  \begin{array}{c|rrrr}
  D & \vf_1        & \vf_2  & \vf_3 & \vf_4    \\
     \hline
    \x_1^{\vphantom{A^A}} & 1   & . & . & . \\
    \x_2 & . & 1 & . & .\\
    \x_3 & 1 & 1 & . & .\\
    \x_4 & . & . & 1 & .\\
    \x_5 & . & . & . & 1
\end{array}\qquad\qquad
  C=\left(
      \begin{array}{cccc}
       2 & 1 & . & . \\
       1 & 2 & . & . \\
       . & . & 1 & . \\
       . & . & . & 1
      \end{array}
    \right)
$$
\hrulefill

\bigskip
\noindent
$G=\SL_2(3)$

$|G|=24=2^3\cdot 3$

$$
\XX(\SL_2(3))\qquad\qquad
  \begin{array}{c|rrrrrrr}
       K                     & 1a & 2a &  3a  &  3b  & 4a & 6a & 6b \\
        |K|                  & 1  & 1  &  4   &   4  &  6 & 4 &  4 \\
       |\C_G(x_{\mbox{}_K})| & 24 & 24 &  6   &   6  &  4 & 6 & 6  \\
     \hline
    \x_1^{\vphantom{A^A}}    & 1  &  1 &  1   &   1  &  1 & 1 & 1\\[2pt]
    \x_2                     & 1  &  1 &  \z  & \ov\z&  1 & \z &\ov\z\\[2pt]
    \x_3                     & 1  &  1 & \ov\z&  \z  &  1 &\ov\z& \z\\[2pt]
    \x_4                     & 2  & -2 & -1   &  -1  &  . &  1 & 1 \\[2pt]
    \x_5                     & 2  & -2 &-\z   &-\ov\z&  . & \z &\ov\z \\[2pt]
    \x_6                     & 2  & -2 &-\ov\z&  \z  &  . &\ov\z& \z \\[2pt]
    \x_7                     & 3  &  3 &  .   &   .  & -1 & . & .
  \end{array}\qquad\qquad\qquad
  \begin{array}{c}
    \mbox{Иррациональные величины} \\[10pt]
    \z=\displaystyle\frac{-1+i\sqrt3}{2^{\vphantom A}}=e^{2\pi i/3}
  \end{array}
$$
$$
\Phi_2(\SL_2(3))
\qquad\qquad
  \begin{array}{c|rrr}
     K                    & 1a       &  3a   & 3b \\
    |K|                   & 1        & 4     & 4 \\
  |\C_G(x_{\mbox{}_K})|    & 24        & 6   & 6 \\

     \hline
    \vf_1^{\vphantom{A^A}} & 1      & 1    & 1\\[2pt]
    \vf_2                  & 1     &  \z   & \ov\z \\[2pt]
    \vf_3                  & 1     & \ov\z   & \z
  \end{array}\qquad\qquad
  \begin{array}{c|rrr}
  D & \vf_1   & \vf_2  & \vf_3 \\
     \hline
    \x_1^{\vphantom{A^A}} &  1  & . & .\\
    \x_2 &  .  & 1 & .\\
    \x_3 &  .  & . & 1\\
    \x_4 &  .  & 1 & 1\\
    \x_5 &  1  & . & 1\\
    \x_6 &  1  & 1 & .\\
    \x_7 &  1  & 1 & 1
\end{array}\qquad\qquad
  C=\left(
      \begin{array}{ccc}
         4 & 2 & 2 \\
         2 & 4 & 2 \\
         2 & 2 & 4
      \end{array}
    \right)
$$

$$
\Phi_3(\SL_2(3))\qquad\qquad
  \begin{array}{c|rrr}
     K                    & 1a   & 2a   & 4a  \\
    |K|                   & 1    & 1    & 6   \\
  |\C_G(x_{\mbox{}_K})|   & 24   & 24   & 4   \\

     \hline
    \vf_1^{\vphantom{A^A}} & 1 & 1 & 1  \\
    \vf_2                  & 2 & -2 & . \\
    \vf_3                  & 3 & 3 & -1
  \end{array}\qquad\qquad
  \begin{array}{c|rrrr}
  D & \vf_1        & \vf_2  & \vf_3     \\
     \hline
    \x_1^{\vphantom{A^A}} & 1   & . & .  \\
    \x_2 & 1 & . & . \\
    \x_3 & 1 & . & . \\
    \x_4 & . & 1 & . \\
    \x_5 & . & 1 & . \\
    \x_6 & . & 1 & . \\
    \x_7 & . & . & 1
\end{array}\qquad\qquad
  C=\left(
      \begin{array}{ccc}
       3 & . & .  \\
       . & 3 & .  \\
       . & . & 1
      \end{array}
    \right)
$$
\hrulefill

\bigskip
\noindent
$G=A_5\cong \PSL_2(5)\cong \SL_2(4)$

$|G|=60=2^2\cdot 3\cdot 5$

$$
\begin{array}{l}
\XX(A_5)
\end{array}\qquad\qquad
  \begin{array}{c|rrrrr}
       K                     & 1a     & 2a    & 3a   & 5a   &  5b \\
        |K|                  & 1      & 15    & 20   & 12   &  12 \\
       |\C_G(x_{\mbox{}_K})|  & 60     & 4     & 3    & 5    &  5  \\
     \hline
    \x_1^{\vphantom{A^A}}    & 1      & 1     & 1    & 1    &  1  \\
    \x_2                     & 3      & -1    & .    & \a   & \a^{\mathrlap*}\\
    \x_3                     & 3      & -1    & .    & \a^{\mathrlap*} & \a\\
    \x_4                     & 4      & .     & 1    & -1   &  -1 \\
    \x_5                     & 5      & 1     & -1   & .    &   .
  \end{array}\qquad\qquad\qquad
  \begin{array}{c}
    \mbox{Иррациональные величины} \\[10pt]
    \a=\displaystyle\frac{1+\sqrt5}{2}=-\z^2-\z^3,\\
    \a^*=\displaystyle\frac{1-\sqrt5}{2}=-\z-\z^4,\\[7pt]
    \z=e^{2\pi i/5}
  \end{array}
$$

$$
\Phi_2(A_5)\qquad\qquad
  \begin{array}{c|rrrrr}
       K                     & 1a   & 3a   & 5a   &  5b \\
        |K|                  & 1    & 20   & 12   &  12 \\
       |\C_G(x_{\mbox{}_K})|  & 60   & 3    & 5    &  5  \\
     \hline
    \vf_1^{\vphantom{A^A}}    & 1    & 1    & 1    &  1  \\
    \vf_2                     & 2    & -1    & -\a^{\mathrlap*} & -\a\\
    \vf_3                     & 2    & -1    & -\a   & -\a^{\mathrlap*}\\
    \vf_4                     & 4    & 1    & -1   &  -1 \\
  \end{array}\qquad\qquad
    \begin{array}{c|rrrr}
  D & \vf_1        & \vf_2  &\vf_3  &\vf_4 \\
     \hline
    \x_1^{\vphantom{A^A}} & 1   & . & . & .\\
    \x_2 & 1  & 1 & . & . \\
    \x_3 & 1  & . & 1 & . \\
    \x_4 & .  & . & . & 1 \\
    \x_5 & 1 & 1 & 1 & .
\end{array}\qquad\qquad
  C=\left(
      \begin{array}{cccc}
        4 & 2 & 2 & .\\
        2 & 2 & 1 & .\\
        2 & 1 & 2& . \\
        . & . & .& 1
      \end{array}
    \right)
$$

$$
\Phi_3(A_5)\qquad\qquad
  \begin{array}{c|rrrr}
       K                     & 1a   &  2a   & 5a   &  5b \\
        |K|                  & 1    &  15   & 12   &  12 \\
       |\C_G(x_{\mbox{}_K})|  & 60   &  4    & 5    &  5  \\
     \hline
    \vf_1^{\vphantom{A^A}}    & 1    & 1    & 1    &  1  \\
    \vf_2                     & 3    & -1   & \a   & \a^{\mathrlap*}\\
    \vf_3                     & 3    & -1   & \a^{\mathrlap*} & \a\\
    \vf_4                     & 4    & .    & -1   &  -1 \\
  \end{array}\qquad\qquad
  \begin{array}{c|rrrr}
  D & \vf_1        & \vf_2  &\vf_3  &\vf_4 \\
     \hline
    \x_1^{\vphantom{A^A}} & 1   & . & . & .\\
    \x_2 & . & 1 & . & . \\
    \x_3 & . & . & 1 & . \\
    \x_4 & . & . & . & 1 \\
    \x_5 & 1 & . & . & 1
\end{array}\qquad\qquad
  C=\left(
      \begin{array}{cccc}
        2 & . & . & 1\\
        . & 1 & . & .\\
        . & . & 1 & . \\
        1 & . & . & 2
      \end{array}
    \right)
$$

$$
\Phi_5(A_5)\qquad\qquad
  \begin{array}{c|rrrr}
       K                     & 1a   &  2a   & 3a   \\
        |K|                  & 1    &  15   & 20   \\
       |\C_G(x_{\mbox{}_K})|  & 60   &  4    & 3    \\
     \hline
    \vf_1^{\vphantom{A^A}}    & 1    & 1    &  1   \\
    \vf_2                     & 3    & -1   &  .  \\
    \vf_3                     & 5   & 1   &  -1  \\
  \end{array}\qquad\qquad
  \begin{array}{c|rrr}
  D & \vf_1        & \vf_2  &\vf_3  \\
     \hline
    \x_1^{\vphantom{A^A}} & 1 & . & . \\
    \x_2 & . & 1 & .  \\
    \x_3 & . & 1 & .  \\
    \x_4 & 1 & 1 & .  \\
    \x_5 & . & . & 1
\end{array}\qquad\qquad
  C=\left(
      \begin{array}{ccc}
        2 & 1 & . \\
        1 & 3 & . \\
        . & . & 1
      \end{array}
    \right)
$$
\newpage

\noindent
$G=\PSL_2(7)\cong \GL_3(2)$

$|G|=168=2^3\cdot 3\cdot 7$

$$
\begin{array}{l}
\XX\big(PSL_2(7)\big)
\end{array}\qquad\qquad
  \begin{array}{c|rrrrrr}
       K                     & 1a     & 2a    & 3a  & 4a  & 7a   &  7b \\
        |K|                  & 1      &  21   & 56  & 42  & 24   & 24    \\
       |\C_G(x_{\mbox{}_K})|  & 168    & 8     & 3   & 4   & 7    & 7    \\
     \hline
    \x_1^{\vphantom{A^A}}    & 1      &  1    &  1  &  1  & 1    &  1  \\
    \x_2                     & 3      & -1    &  .  &  1  & \a   & \ov\a\\
    \x_3                     & 3      & -1    &  .  &  1  & \ov\a & \a\\
    \x_4                     & 6      &  2    &  .  &  .  & -1   &  -1 \\
    \x_5                     & 7      & -1    &  1  & -1  & .    &   .\\
    \x_6                     & 8      &  .    & -1  &  .  & 1    &   1
  \end{array}\qquad\qquad\qquad
  \begin{array}{c}
    \mbox{Иррациональные величины} \\[10pt]
    \a=\displaystyle\frac{-1+i\sqrt7}{2^{\vphantom A}}=\z+\z^2+\z^4,\\[7pt]
    \z=e^{2\pi i/7}
  \end{array}
$$

$$
\Phi_2\big(PSL_2(7)\big)\qquad\qquad
  \begin{array}{c|rrrrr}
       K                     & 1a   & 3a   & 7a   &  7b \\
        |K|                  & 1    & 56   & 24   &  24 \\
       |\C_G(x_{\mbox{}_K})|  & 168   & 3    & 7    &  7  \\
     \hline
    \vf_1^{\vphantom{A^A}}    & 1    & 1    & 1    &  1  \\
    \vf_2                     & 3    & .    & \a &  \ov\a\\
    \vf_3                     & 3    & .    & \ov\a   & \a\\
    \vf_4                     & 8    & -1    & 1   &  1 \\
  \end{array}\qquad\quad
    \begin{array}{c|rrrr}
  D & \vf_1        & \vf_2  &\vf_3  &\vf_4 \\
     \hline
    \x_1^{\vphantom{A^A}} & 1   & . & . & .\\
    \x_2 & .  & 1 & . & . \\
    \x_3 & .  & . & 1 & . \\
    \x_4 & .  & 1 & 1 & . \\
    \x_5 & 1 & 1 & 1 & . \\
    \x_6 & . & . & . & 1
\end{array}\qquad
  C=\left(
      \begin{array}{cccc}
        2 & 1 & 1 & .\\
        1 & 3 & 2 & .\\
        1 & 2 & 3& . \\
        . & . & .& 1
      \end{array}
    \right)
$$

$$
\Phi_3\big(PSL_2(7)\big)\quad
  \begin{array}{c|rrrrr}
       K                     & 1a   &  2a &  4a   & 7a   &  7b \\
        |K|                  & 1    &  21 &  42   & 24   &  24 \\
       |\C_G(x_{\mbox{}_K})|  & 168   &  8  &  4    & 7    &  7  \\
     \hline
    \vf_1^{\vphantom{A^A}}    & 1    &  1  &  1   & 1    &  1  \\
    \vf_2                     & 3    & -1  &  1   & \a   & \ov\a\\
    \vf_3                     & 3    & -1  &  1   & \ov\a & \a\\
    \vf_4                     & 6    &  2  &  .   & -1   &  -1 \\
    \vf_5                     & 7    & -1  & -1   & .   &  . \\
  \end{array}\quad
  \begin{array}{c|rrrrr}
  D & \vf_1        & \vf_2  &\vf_3  &\vf_4  &\vf_5\\
     \hline
    \x_1^{\vphantom{A^A}} & 1   & . & . & . &.\\
    \x_2 & . & 1 & . & . & . \\
    \x_3 & . & . & 1 & . & .\\
    \x_4 & . & . & . & 1 & . \\
    \x_5 & . & . & . & . & 1 \\
    \x_6 & 1 & . & . & . & 1
\end{array}\quad
  C=\left(
      \begin{array}{ccccc}
        2 & . & . & . & 1\\
        . & 1 & . & . & .\\
        . & . & 1 & . & .\\
        . & . & . & 1 & .\\
        1 & . & . & . & 2
      \end{array}
    \right)
$$

$$
\Phi_7\big(PSL_2(7)\big)\qquad\quad
  \begin{array}{c|rrrr}
       K                     & 1a   &  2a   & 3a    & 4a   \\
        |K|                  & 1    &  21   & 56    & 42  \\
       |\C_G(x_{\mbox{}_K})|  & 168   &  8    & 3     & 4   \\
     \hline
    \vf_1^{\vphantom{A^A}}    & 1    & 1    &  1    &  1  \\
    \vf_2                     & 3    & -1   &  .    &  1 \\
    \vf_3                     & 5   & 1   &  -1   &  -1 \\
    \vf_4                     & 7   & -1   &  1   &  -1 \\
  \end{array}\qquad\quad
  \begin{array}{c|rrrr}
  D & \vf_1  & \vf_2  &\vf_3 &\vf_4  \\
     \hline
    \x_1^{\vphantom{A^A}} & 1 & . & . & .\\
    \x_2 & . & 1 & . & . \\
    \x_3 & . & 1 & . & . \\
    \x_4 & 1 & . & 1 & . \\
    \x_5 & . & . & . & 1 \\
    \x_6 & . & 1 & 1 & .
\end{array}\qquad\quad
  C=\left(
      \begin{array}{cccc}
        2 & . & 1 & . \\
        . & 3 & 1 & . \\
        1 & 1 & 2 & . \\
        . & . & . & 1
      \end{array}
    \right)
$$
\hrulefill

\newpage
\noindent
$G=\la c\mid c^n=1\ra$

$|G|=n$

$$
\XX(G)\qquad
  \begin{array}{c|ccccc}
            \XX(C)              & 1      & c  & c^2 & \ld & c^{n-1}\\
     \hline
    \x_1^{\vphantom{A^A}} & 1      & 1 & 1 & \ld  & 1\\
    \x_2                  & 1      & \z & \z^2 & \ld  & \z^{(n-1)} \\
    \x_3                  & 1      & \z^2 & \z^4 & \ld  & \z^{2(n-1)} \\
\ld & \ld   & \ld  & \ld  & \ld  & \ld \\
    \x_n                  & 1  & \z^{n-1} & \z^{(n-1)2} & \ld  & \z^{(n-1)^2}
  \end{array}\qquad\qquad\qquad
  \begin{array}{c}
    \mbox{Иррациональные величины} \\[5pt]
    \z=e^{2\pi i/n}
  \end{array}
$$
\hrulefill

\bigskip
\noindent
$G=D_{2n}=\la x,y\mid x^n=y^2=1,\ x^y=x^{-1}\ra,\ \ n\ge 3$

$|G|=2n$

$$
\begin{array}{c}
\XX(D_{2n})\\
n\ \ \text{нечётно}
\end{array}
\qquad\quad
  \begin{array}{c|rrc}
       K                     & 1a   &  2a  & (x^r)^G   \\
        |K|                  &  1   &  n   &  2     \\
       |\C_G(x_{\mbox{}_K})|  & 2n   &  2    & n   \\
           x_{\mbox{}_K}   &    1    &  y   & x^r     \\
     \hline
    \x_1^{\vphantom{A^A}}    & 1    & 1    &  1   \\
    \x_2                     & 1    & -1   &  1  \\
    \th_s                     & 2   & .   &  \z^{rs}+\z^{-rs}\\
  \end{array}\qquad\qquad\qquad
  \begin{array}{c}
      \mbox{Параметры} \\[5pt]
    r,s \ =\ 1,\ld,(n-1)/2\\[10pt]

    \mbox{Иррациональные величины} \\[5pt]
    \z=e^{2\pi i/n}
  \end{array}
$$

$$
\begin{array}{c}
\XX(D_{2n})\\
n=2m
\end{array}
\qquad\quad
  \begin{array}{c|rcrrc}
       K                     & 1a   &  2a  & 2b & 2c & (x^r)^G   \\
        |K|                  &  1   &  1   &  m & m & 2 \\
       |\C_G(x_{\mbox{}_K})|  & 2n   &  2n   & 4 & 4 & n \\
           x_{\mbox{}_K}   &    1    &  x^m   & y & xy & x^r     \\
     \hline
    \x_1^{\vphantom{A^A}}    & 1    & 1    &  1 & 1 & 1  \\
    \x_2                     & 1    & 1   &  -1 & -1 & 1 \\
    \x_3                     & 1    & (-1)^m  &  1 & -1 & (-1)^r \\
    \x_4                     & 1    & (-1)^m  &  -1 & 1 & (-1)^r \\
    \th_s                     & 2   & (-1)^s\,2  & . & . & \z^{rs}+\z^{-rs}\\
  \end{array}\qquad
  \begin{array}{c}
      \mbox{Параметры} \\[5pt]
    r,s\ =\ 1,\ld,m-1\\[10pt]

    \mbox{Иррациональные величины} \\[5pt]
    \z=e^{2\pi i/n}
  \end{array}
$$
\hrulefill

\bigskip
\noindent
$G=Q_{2^m}=\la x,y\mid x^{2^{m-1}}=1,\ y^2=x^{2^{m-2}},\ x^y=x^{-1}\ra,\ \ m\ge3$

$|G|=2^m$

$$
\XX(Q_{2^m})
\qquad\quad
  \begin{array}{c|ccccc}
       K                     & 1a   &  2a  & 4a & 4b & (x^r)^G  \\
        |K|                  &  1   &  1   & 2^{m-2} & 2^{m-2} & 2 \\
       |\C_G(x_{\mbox{}_K})|  & 2^m   &  2^m   & 4 & 4 & 2^{m-1} \\
           x_{\mbox{}_K}   &    1    &  y^2   & y & xy & x^r     \\
     \hline
    \x_1^{\vphantom{A^A}}    & 1    & 1    &  1 & 1 & 1  \\
    \x_2                     & 1    & 1   &  \mathllap{-}1 & \mathllap{-}1 & 1 \\
    \x_3                     & 1    & 1  &  1 & \mathllap{-}1 & (-1)^r \\
    \x_4                     & 1    & 1  &  \mathllap{-}1 & 1 & (-1)^r \\
    \th_s                     & 2   & \mathllap{-}2  & . & . & \z^{rs}+\z^{-rs}\\
  \end{array}\qquad\qquad
  \begin{array}{c}
      \mbox{Параметры} \\[5pt]
    r,s\ =\ 1,\ld,2^{m-2}-1\\[10pt]

    \mbox{Иррациональные величины} \\[5pt]
    \z=e^{2\pi i/2^{m-1}}
  \end{array}
$$
\hrulefill

\newpage
\noindent
$G=\PSL_2(q)$

\medskip
$|G|=\frac{1}{(2,q-1)}q(q^2-1)$

$$
\begin{array}{l}
\XX\big(PSL_2(q)\big)\\
q=2^m
\end{array}\qquad
  \begin{array}{c|cccc}
       K                     & 1a      & 2a    & (x^r)^G  & (y^t)^G    \\
        |K|                  & 1       & q^2-1   & q(q+1)  & q(q-1)  \\
       |\C_G(x_{\mbox{}_K})| & q(q^2\!-\!1)& q     & q-1 &  q+1  \\
        x_{\mbox{}_K}        &  1      & z     &  x^r  & y^t \\
     \hline
    \x_1^{\vphantom{A^A}}    & 1      &  1    &  1  &  1  \\
    \x_2                     & q      & .    &  1  &  \mathllap{-}1  \\
    \th_s                    & q+1    & 1    &  \z^{rs}\!+\!\z^{-rs}  &  .  \\
    \t_h                     & q-1     &  \mathllap{-}1  &  .  & -\xi^{th}\!-\xi^{-th}
  \end{array}\qquad
  \begin{array}{c}
      \mbox{Элементы} \\[3pt]
    |z|=2,\ |x|=q-1,\ |y|=q+1\\[10pt]

      \mbox{Параметры} \\[3pt]
    r,s\ = \ 1,\ld,(q-2)/2;\\
    t,h\ = \ 1,\ld,q/2\\[10pt]

    \mbox{Иррациональные величины} \\[3pt]
    \z=e^{2\pi i/(q-1)},\ \xi=e^{2\pi i/(q+1)}
  \end{array}
$$

\medskip

$$
\begin{array}{l}
\Phi_p\big(PSL_2(q)\big)\\
q=2^m,\\
 p\mid(q-1)
\end{array}
\qquad
  \begin{array}{c|cccc}
       K                     & 1a   & 2a   & (x_0^r)^G   & (y^t)^G \\
        |K|                  & 1    & q^2-1 & q(q+1)   &  q(q-1) \\
       |\C_G(x_{\mbox{}_K})|  & q(q^2-1) & q & q-1    &  q+1  \\
       x_{\mbox{}_K}        &  1      & z   &  x_0^r  & y^t \\
     \hline
    \vf_1^{\vphantom{A^A}}    & 1    & 1    & 1    &  1  \\
    \vf_2                     & q    & .    & 1 &  \mathllap{-}1\\
    \psi_j                     & q+1    & 1    & \z_0^{rj}\!+\!\z_0^{-rj}   & .\\
    \eta_f                     & q-1    & \mathllap{-}1    & .   & -\xi^{tf}\!-\xi^{-tf} \\
  \end{array}\qquad
      \begin{array}{c|rrrr}
  D & \vf_1   & \vf_2  &\psi_j  &\eta_f \\
     \hline
    \x_1^{\vphantom{A^A}} & 1   & . & . & .\\
    \x_2 & .  & 1 & . & . \\
    \th_{s_0} & 1  & 1 & . & . \\
    \th_{s_i} & .  & . & \delta_{ij} & . \\
    \t_h & . & . & . & \delta_{hf}
    \end{array}
$$

$$
  C=\left(
      \begin{array}{cccc}
        \frac{p^d+1}{2} & \frac{p^d-1}{2} & . & .\\
        \frac{p^d-1}{2} & \frac{p^d+1}{2} & . & .\\
        . & . & p^d \I_{(n-1)/2}& . \\
        . & . & .& \I_{q/2}
      \end{array}
    \right)
\qquad
  \begin{array}{c}
\mbox{Обозначения} \\[2pt]
 p^d=(q-1)_p,\ n=(q-1)_{p'},\\
 x_0=x^{p^d},\ \z_0=\z^{p^d}=e^{2\pi i/n}\\[10pt]

    \mbox{Параметры} \\[3pt]
    r,i,j\ = \ 1,\ld,\frac{n-1}{2};\\
    f,h,t\ = \ 1,\ld,\frac{q}{2};\\
    s_0=n,2n,\ld,\frac{p^d-1}{2}n;\\
    s_i=i,n\pm i,2n\pm i,\ld,\frac{p^d-1}{2}n\pm i
  \end{array}
$$

\medskip

$$
\begin{array}{l}
\Phi_p\big(PSL_2(q)\big)\\
q=2^m,\\
 p\mid(q+1)
\end{array}
\qquad
  \begin{array}{c|cccc}
       K                     & 1a   & 2a   & (x^r)^G   & (y_0^t)^G \\
        |K|                  & 1    & q^2-1 & q(q+1)   &  q(q-1) \\
       |\C_G(x_{\mbox{}_K})|  & q(q^2-1) & q & q-1    &  q+1  \\
       x_{\mbox{}_K}        &  1      & z   &  x^r  & y_0^t \\
     \hline
    \vf_1^{\vphantom{A^A}}    & 1    & 1    & 1    &  1  \\
    \vf_2                     & q-1    & -1    & . &  \mathllap{-}2\\
    \psi_j                     & q+1    & 1    & \z^{rj}\!+\!\z^{-rj}   & .\\
    \eta_f                     & q-1    & \mathllap{-}1    & .   & -\xi_0^{tf}\!-\xi_0^{-tf} \\
  \end{array}\qquad
      \begin{array}{c|rrrr}
  D & \vf_1   & \vf_2  &\psi_j  &\eta_f \\
     \hline
    \x_1^{\vphantom{A^A}} & 1   & . & . & .\\
    \x_2 & 1  & 1 & . & . \\
    \th_s & .  & . & \delta_{sj} & . \\
    \t_{h_0} & .  & 1 & . & . \\
    \t_{h_i} & . & . & . & \delta_{if}
    \end{array}
$$

$$
  C=\left(
      \begin{array}{cccc}
        2 & 1 & . & .\\
        1 & \frac{p^d+1}{2} & . & .\\
        . & . &  \I_{(q/2)-1}& . \\
        . & . & .& p^d \I_{(n-1)/2}
      \end{array}
    \right)
\qquad
  \begin{array}{c}
\mbox{Обозначения} \\[2pt]
 p^d=(q+1)_p,\ n=(q+1)_{p'},\\
 y_0=y^{p^d},\ \xi_0=\xi^{p^d}=e^{2\pi i/n}\\[10pt]

    \mbox{Параметры} \\[3pt]
    r,s,j\ = \ 1,\ld,\frac{q}{2}-1;\\
    i,f,t\ = \ 1,\ld,\frac{n-1}{2};\\
    h_0=n,2n,\ld,\frac{p^d-1}{2}n;\\
    h_i=i,n\pm i,2n\pm i,\ld,\frac{p^d-1}{2}n\pm i
  \end{array}
$$

\newpage
\bigskip
$$
\begin{array}{l}
\XX\big(PSL_2(q)\big)\\
q\equiv 1 \!\! \pmod 4

\end{array}\qquad
  \begin{array}{c|cccccc}
       K                     & 1a      & 2a                    & la & lb & (x^r)^G  & (y^t)^G    \\
        |K|                  & 1 & \frac{1}{2}q(q\!+\!1) & \frac{1}{2}(q^2\!-\!1) & \frac{1}{2}(q^2\!-\!1) & q(q+1)  & q(q-1)  \\
       |\C_G(x_{\mbox{}_K})| & \frac{1}{2}q(q^2\!-\!1)&   q-1  & q  & q  & \frac{1}{2}(q-1) &  \frac{1}{2}(q+1)  \\
        x_{\mbox{}_K}        &  1      & z                     & u & v  &  x^r  & y^t \\
     \hline
    \x_1^{\vphantom{A^A}}    & 1      &  1    &  1  &  1 & 1 & 1 \\
    \x_2                     & q      & 1    &  .    & . & 1  &  \mathllap{-}1  \\
    \x_3                     & \frac{1}{2}(q+1)  & (-1)^{(q-1)/4}    & \frac{1}{2}(1+\sqrt{q}) & \frac{1}{2}(1-\sqrt{q}) & (-1)^r  &  .  \\
    \x_4                     & \frac{1}{2}(q+1)  & (-1)^{(q-1)/4}    & \frac{1}{2}(1-\sqrt{q}) & \frac{1}{2}(1+\sqrt{q}) & (-1)^r  &  .  \\
    \th_s                    & q+1    & 2\cdot(-1)^s &  1 & 1 & \z^{rs}\!+\!\z^{-rs}  &  .  \\
    \t_h                     & q-1    & .   &  \mathllap{-}1 &  \mathllap{-}1  &  .  & -\xi^{th}\!-\xi^{-th}
  \end{array}
$$
$$
  \begin{array}{c}
    \mbox{Обозначения} \\[3pt]
    q=l^m,\ \ l\ \mbox{простое} \\[10pt]

     \mbox{Параметры} \\[3pt]
    r,s\ = \ 1,\ld,\frac{q-5}{4};\\
    t,h\ = \ 1,\ld,\frac{q-1}{4}\\[10pt]
\end{array}\qquad
\begin{array}{c}
      \mbox{Элементы} \\[3pt]
    |u|=|v|=l,\\
  |x|=\frac{q-1}{2},\ |y|=\frac{q+1}{2},\\
z=x^{(q-1)/4}\\[10pt]

    \mbox{Иррациональные величины} \\[3pt]
    \z=e^{4\pi i/(q-1)},\ \xi=e^{4\pi i/(q+1)}
  \end{array}
$$

\medskip

$$
\begin{array}{l}
\Phi_p\big(PSL_2(q)\big)\\
q\equiv 1 \!\! \pmod 4,\\
2\ne p\mid(q-1)
\end{array}
\qquad
  \begin{array}{c|cccccc}
       K                     & 1a   & 2a  & la & lb  & (x_0^r)^G   & (y^t)^G \\
    |K|                  & 1 & \frac{1}{2}q(q\!+\!1) & \frac{1}{2}(q^2\!-\!1) & \frac{1}{2}(q^2\!-\!1) & q(q\!+\!1)  & q(q\!-\!1)  \\
    |\C_G(x_{\mbox{}_K})| & \frac{1}{2}q(q^2\!\!-\!1)&   q-1  & q  & q  & \frac{1}{2}(q\!-\!1) &  \frac{1}{2}(q\!+\!1)  \\
        x_{\mbox{}_K}        &  1      & z                     & u & v  &  x_0^r  & y^t \\
     \hline
    \vf_1^{\vphantom{A^A}}   & 1      &  1    &  1  &  1 & 1 & 1 \\
    \vf_2                     & q      & 1    &  .    & . & 1  &  \mathllap{-}1  \\
    \vf_3                     & \frac{1}{2}(q+1)  & (-1)^{(q-1)/4}    & \frac{1}{2}(1+\sqrt{q}) & \frac{1}{2}(1-\sqrt{q}) & (-1)^r  &  .  \\
    \vf_4                     & \frac{1}{2}(q+1)  & (-1)^{(q-1)/4}    & \frac{1}{2}(1-\sqrt{q}) & \frac{1}{2}(1+\sqrt{q}) & (-1)^r  &  .  \\
    \psi_j                      & q+1    & 2\cdot(-1)^j &  1 & 1 & \z_0^{rj}\!+\!\z_0^{-rj}  &  .  \\
    \eta_f                     & q-1    & .   &  \mathllap{-}1 &  \mathllap{-}1  &  .  & -\xi^{tf}\!-\xi^{-tf}
  \end{array}
$$

$$ \begin{array}{c|rrrrrr}
  D & \vf_1   & \vf_2 & \vf_3 & \vf_4 &\psi_j  &\eta_f \\
     \hline
    \x_1^{\vphantom{A^A}} & 1   & . & . & . & . & .\\
    \x_2 & .  & 1 & . & . & . & .\\
    \x_3 & . & . & 1 & . & . & . \\
    \x_4 & . & . & . & 1 & . & . \\
    \th_{s_+} & 1  & 1 & . & . & . & .\\
    \th_{s_-} & .  & . & 1 & 1 & . & .\\
    \th_{s_i} & .  & . & . & . & \delta_{ij} & . \\
    \t_h & . & . & . & . & . & \delta_{hf}
    \end{array} \qquad
  C=\left(
      \begin{array}{cccccc}
        \frac{p^d+1}{2} & \frac{p^d-1}{2} & . & . & . & .\\[2pt]
        \frac{p^d-1}{2} & \frac{p^d+1}{2} & . & . & . & .\\
         . & . & \frac{p^d+1}{2} & \frac{p^d-1}{2}  & . & .\\[2pt]
         . & . &  \frac{p^d-1}{2} & \frac{p^d+1}{2}  & . & .\\
         . & . & . & . & p^d \I_{(n/2)-1}& . \\
         . & . & . & . & .& \I_{(q-1)/4}
      \end{array}
    \right)
$$

$$
\begin{array}{c}
\mbox{Обозначения} \\[2pt]
 p^d=(q-1)_p,\ n=(\frac{q-1}{2})_{p'},\\
 x_0=x^{p^d},\ \z_0=\z^{p^d}=e^{2\pi i/n}\\
\end{array}\qquad
  \begin{array}{c}
    \mbox{Параметры} \\[3pt]
    r,i,j\ = \ 1,\ld,\frac{n}{2}-1;\\
    f,h,t\ = \ 1,\ld,\frac{q-1}{4};\\
    s_+=n,2n,\ld,\frac{p^d-1}{2}n;\\
    s_-=\frac{n}{2},n+\frac{n}{2},\ld,\frac{p^d-3}{2}n+\frac{n}{2};\\
    s_i=i,n\pm i,2n\pm i,\ld,\frac{p^d-1}{2}n\pm i
  \end{array}
$$
\medskip

$$
\begin{array}{l}
\Phi_p\big(PSL_2(q)\big)\\
q\equiv 1 \!\! \pmod 4,\\
2\ne p\mid(q+1)
\end{array}
\qquad
  \begin{array}{c|cccccc}
       K                     & 1a   & 2a  & la & lb  & (x^r)^G   & (y_0^t)^G \\
    |K|                  & 1 & \frac{1}{2}q(q\!+\!1) & \frac{1}{2}(q^2\!-\!1) & \frac{1}{2}(q^2\!-\!1) & q(q\!+\!1)  & q(q\!-\!1)  \\
    |\C_G(x_{\mbox{}_K})| & \frac{1}{2}q(q^2\!\!-\!1)&   q-1  & q  & q  & \frac{1}{2}(q\!-\!1) &  \frac{1}{2}(q\!+\!1)  \\
        x_{\mbox{}_K}        &  1      & z                     & u & v  &  x^r  & y_0^t \\
     \hline
    \vf_1^{\vphantom{A^A}}   & 1      &  1    &  1  &  1 & 1 & 1 \\
    \vf_2                     & q-1    & .   &  \mathllap{-}1    & \mathllap{-}1 & .  &  \mathllap{-}2  \\
    \vf_3                     & \frac{1}{2}(q+1)  & (-1)^{(q-1)/4}    & \frac{1}{2}(1+\sqrt{q}) & \frac{1}{2}(1-\sqrt{q}) & (-1)^r  &  .  \\
    \vf_4                     & \frac{1}{2}(q+1)  & (-1)^{(q-1)/4}    & \frac{1}{2}(1-\sqrt{q}) & \frac{1}{2}(1+\sqrt{q}) & (-1)^r  &  .  \\
    \psi_j                      & q+1    & 2\cdot(-1)^j &  1 & 1 & \z^{rj}\!+\!\z^{-rj}  &  .  \\
    \eta_f                     & q-1    & .   &  \mathllap{-}1 &  \mathllap{-}1  &  .  & -\xi_0^{tf}\!-\xi_0^{-tf}
  \end{array}
$$

$$ \begin{array}{c|rrrrrr}
  D & \vf_1   & \vf_2 & \vf_3 & \vf_4 &\psi_j  &\eta_f \\
     \hline
    \x_1^{\vphantom{A^A}} & 1   & . & . & . & . & .\\
    \x_2 & 1  & 1 & . & . & . & .\\
    \x_3 & . & . & 1 & . & . & . \\
    \x_4 & . & . & . & 1 & . & . \\
    \th_s & .  & . & . & . & \delta_{sj} & .\\
    \t_{h_0} & .  & 1 & . & . & . & .\\
    \t_{h_i} & . & . & . & . & . & \delta_{if}
    \end{array} \qquad
  C=\left(
      \begin{array}{cccccc}
        2 & 1 & . & . & . & .\\[2pt]
        1 & \frac{p^d+1}{2} & . & . & . & .\\
         . & . & 1 & .  & . & .\\[2pt]
         . & . &  . & 1 & . & .\\
         . & . & . & . & \I_{(q-5)/4}& . \\
         . & . & . & . & .& p^d \I_{(n-1)/2}
      \end{array}
    \right)
$$

$$
\begin{array}{c}
\mbox{Обозначения} \\[2pt]
 p^d=(q+1)_p,\ n=(\frac{q+1}{2})_{p'},\\
 y_0=y^{p^d},\ \xi_0=\xi^{p^d}=e^{2\pi i/n}\\
\end{array}\qquad
  \begin{array}{c}
    \mbox{Параметры} \\[3pt]
    r,s,j\ = \ 1,\ld,\frac{q-5}{4};\\
    i,f,t\ = \ 1,\ld,\frac{n-1}{2};\\
    h_0=n,2n,\ld,\frac{p^d-1}{2}n;\\
    h_i=i,n\pm i,2n\pm i,\ld,\frac{p^d-1}{2}n\pm i
  \end{array}
$$

\medskip

$$
\begin{array}{l}
\Phi_2\big(PSL_2(q)\big)\\
q\equiv 1 \!\! \pmod 4\\
\end{array}
\qquad
\begin{array}{c|ccccc}
     K                     & 1a   & la & lb  & (x_0^r)^G   & (y^t)^G \\
  |K|                  & 1 & \frac{1}{2}(q^2\!-\!1) & \frac{1}{2}(q^2\!-\!1) & q(q\!+\!1)  & q(q\!-\!1)  \\
  |\C_G(x_{\mbox{}_K})| & \frac{1}{2}q(q^2\!\!-\!1)  & q  & q  & \frac{1}{2}(q\!-\!1) &  \frac{1}{2}(q\!+\!1)  \\
      x_{\mbox{}_K}        &  1                         & u & v  &  x_0^r  & y^t \\
   \hline
  \vf_1^{\vphantom{A^A}}   & 1          &  1  &  1 & 1 & 1 \\
  \vf_2                     & \frac{1}{2}(q-1)      & \frac{1}{2}(-1+\sqrt{q}) & \frac{1}{2}(-1-\sqrt{q}) & . & \mathllap{-}1  \\
  \vf_3                     & \frac{1}{2}(q-1)      & \frac{1}{2}(-1-\sqrt{q}) & \frac{1}{2}(-1+\sqrt{q}) & . & \mathllap{-}1  \\
  \psi_j                      & q+1     &  1 & 1 & \z_0^{rj}\!+\!\z_0^{-rj}  &  .  \\
  \eta_f                     & q-1    &  \mathllap{-}1 &  \mathllap{-}1  &  .  & -\xi^{tf}\!-\xi^{-tf}
\end{array}
$$

$$ \begin{array}{c|rrrrr}
  D & \vf_1   & \vf_2 & \vf_3 &\psi_j  &\eta_f \\
     \hline
    \x_1^{\vphantom{A^A}} & 1   & . & . & . & .\\
    \x_2 & 1  & 1 & 1 & .  & .\\
    \x_3 & 1 & 1 & . & . & . \\
    \x_4 & 1 & . & 1 & . & . \\
    \th_{s_0} & 2  & 1 &  1 & . & .\\
    \th_{s_i} & .  & . &  . & \delta_{ij} & . \\
    \t_h & . & . & . & .  & \delta_{hf}
    \end{array} \qquad
  C=\left(
      \begin{array}{ccccc}
        2^d & 2^{d-1} & 2^{d-1} & . & . \\[2pt]
        2^{d-1} & 2^{d-2}+1 & 2^{d-2} & . & . \\
        2^{d-1} & 2^{d-2} & 2^{d-2}+1 & .  &  .\\[2pt]
         . & . &  . &  2^{d-1} \I_{(n-1)/2} & .\\
         . & . & . & . & \I_{(q-1)/4}
      \end{array}
    \right)
$$

$$
\begin{array}{c}
\mbox{Обозначения} \\[2pt]
 2^d=(q-1)_2,\ n=(q-1)_{2'},\\
 x_0=x^{2^{d-1}},\ \z_0=\z^{2^{d-1}}=e^{2\pi i/n}\\
\end{array}\qquad
  \begin{array}{c}
    \mbox{Параметры} \\[3pt]
    r,i,j\ = \ 1,\ld,\frac{n-1}{2};\\
    f,h,t\ = \ 1,\ld,\frac{q-1}{4};\\
    s_0=n,2n,\ld,(2^{d-2}-1)n;\\
    s_i=i,n\pm i,2n\pm i,\ld,(2^{d-2}-1)n\pm i,2^{d-2}n-i
  \end{array}
$$

\bigskip
$$
\begin{array}{l}
\XX\big(PSL_2(q)\big)\\
q\equiv -1 \!\! \pmod 4

\end{array}\quad
  \begin{array}{c|cccccc}
       K                     & 1a      & 2a                    & la & lb & (x^r)^G  & (y^t)^G    \\
        |K|                  & 1 & \frac{1}{2}q(q\!-\!1) & \frac{1}{2}(q^2\!-\!1) & \frac{1}{2}(q^2\!-\!1) & q(q\!+\!1)  & q(q\!-\!1)  \\
       |\C_G(x_{\mbox{}_K})| & \frac{1}{2}q(q^2\!\!-\!1)&   q+1  & q  & q  & \frac{1}{2}(q\!-\!1) &  \frac{1}{2}(q\!+\!1)  \\
        x_{\mbox{}_K}        &  1      & z                     & u & v  &  x^r  & y^t \\
     \hline
    \x_1^{\vphantom{A^A}}    & 1      &  1    &  1  &  1 & 1 & 1 \\
    \x_2                     & q      & \mathllap{-}1    &  .    & . & 1  &  \mathllap{-}1  \\
    \x_3                     & \frac{1}{2}(q-1)  & -(-1)^{(q+1)/4}    & \frac{1}{2}(-1\!+\!i\sqrt{q}) & \frac{1}{2}(-1\!-\!i\sqrt{q}) & . &-(-1)^t \\
    \x_4                     & \frac{1}{2}(q-1)  & -(-1)^{(q+1)/4}    & \frac{1}{2}(-1\!-\!i\sqrt{q}) & \frac{1}{2}(-1\!+\!i\sqrt{q}) & . &-(-1)^t \\
    \th_s                    & q+1    & . &  1 & 1 & \z^{rs}\!+\!\z^{-rs}  &  .  \\
    \t_h                     & q-1    & -2\cdot(-1)^h   &  \mathllap{-}1 &  \mathllap{-}1  &  .  & -\xi^{th}\!-\xi^{-th}
  \end{array}
$$
\medskip
$$
  \begin{array}{c}
    \mbox{Обозначения} \\[3pt]
    q=l^m,\ \ l\ \mbox{простое} \\[10pt]

      \mbox{Параметры} \\[3pt]
    r,s,t,h\ = \ 1,\ld,(q-3)/4\\[10pt]
\end{array}\qquad
\begin{array}{c}
      \mbox{Элементы} \\[3pt]
    |u|=|v|=l,\\
  |x|=\frac{1}{2}(q-1),\ |y|=\frac{1}{2}(q+1),\\
z=y^{(q+1)/4}\\[10pt]

    \mbox{Иррациональные величины} \\[3pt]
    \z=e^{4\pi i/(q-1)},\ \xi=e^{4\pi i/(q+1)}
  \end{array}
$$

\bigskip

$$
\begin{array}{l}
\Phi_p\big(PSL_2(q)\big)\\
q\equiv\! -1 \!\! \pmod 4,\\
2\ne p\mid(q-1)
\end{array}
\quad\!
  \begin{array}{c|cccccc}
       K                     & 1a      & 2a                    & la & lb & (x_0^r)^G  & (y^t)^G    \\
        |K|                  & 1 & \frac{1}{2}q(q\!-\!1) & \frac{1}{2}(q^2\!-\!1) & \frac{1}{2}(q^2\!-\!1) & q(q\!+\!1)  & q(q\!-\!1)  \\
       |\C_G(x_{\mbox{}_K})| & \frac{1}{2}q(q^2\!\!-\!1)&   q+1  & q  & q  & \frac{1}{2}(q\!-\!1) &  \frac{1}{2}(q\!+\!1)  \\
        x_{\mbox{}_K}        &  1      & z                     & u & v  &  x_0^r  & y^t \\
     \hline
    \vf_1^{\vphantom{A^A}}    & 1      &  1    &  1  &  1 & 1 & 1 \\
    \vf_2                     & q      & \mathllap{-}1    &  .    & . & 1  &  \mathllap{-}1  \\
    \vf_3                     & \frac{1}{2}(q-1)  & -(-1)^{(q+1)/4}    & \frac{1}{2}(-1\!+\!i\sqrt{q}) & \frac{1}{2}(-1\!-\!i\sqrt{q}) & . &-(-1)^t \\
    \vf_4                     & \frac{1}{2}(q-1)  & -(-1)^{(q+1)/4}    & \frac{1}{2}(-1\!-\!i\sqrt{q}) & \frac{1}{2}(-1\!+\!i\sqrt{q}) & . &-(-1)^t \\
    \psi_j                    & q+1    & . &  1 & 1 & \z_0^{rj}\!+\!\z_0^{-rj}  &  .  \\
    \eta_f                     & q-1    & -2\cdot(-1)^f   &  \mathllap{-}1 &  \mathllap{-}1  &  .  & -\xi^{tf}\!-\xi^{-tf}
  \end{array}
$$

\medskip

$$ \begin{array}{c|rrrrrr}
  D & \vf_1   & \vf_2 & \vf_3 & \vf_4 &\psi_j  &\eta_f \\
     \hline
    \x_1^{\vphantom{A^A}} & 1   & . & . & . & . & .\\
    \x_2 & .  & 1 & . & . & . & .\\
    \x_3 & . & . & 1 & . & . & . \\
    \x_4 & . & . & . & 1 & . & . \\
    \th_{s_0} & 1  & 1 & . & . & . & .\\
    \th_{s_i} & .  & . & . & . & \delta_{ij} & .\\
    \t_h & . & . & . & . & . & \delta_{hf}
    \end{array} \qquad
  C=\left(
      \begin{array}{cccccc}
        \frac{p^d+1}{2} & \frac{p^d-1}{2} & . & . & . & .\\[2pt]
        \frac{p^d-1}{2} & \frac{p^d+1}{2} & . & . & . & .\\
         . & . & 1 & .  & . & .\\[2pt]
         . & . &  . & 1 & . & .\\
         . & . & . & . & p^d \I_{(n-1)/2}& . \\
         . & . & . & . & . & \I_{(q-3)/4}
      \end{array}
    \right)
$$

\medskip

$$
\begin{array}{c}
\mbox{Обозначения} \\[2pt]
 p^d=(q-1)_p,\ n=(\frac{q-1}{2})_{p'},\\
 x_0=x^{p^d},\ \z_0=\z^{p^d}=e^{2\pi i/n}\\
\end{array}\qquad
  \begin{array}{c}
    \mbox{Параметры} \\[3pt]
    r,i,j\ = \ 1,\ld,\frac{n-1}{2};\\
    f,h,t\ = \ 1,\ld,\frac{q-3}{4};\\
    s_0=n,2n,\ld,\frac{p^d-1}{2}n;\\
    s_i=i,n\pm i,2n\pm i,\ld,\frac{p^d-1}{2}n\pm i
  \end{array}
$$

\newpage 

$$
\begin{array}{l}
\Phi_p\big(PSL_2(q)\big)\\
q\equiv\! -1 \!\! \pmod 4,\\
2\ne p\mid(q+1)
\end{array}
\quad\!
  \begin{array}{c|cccccc}
       K                     & 1a      & 2a                    & la & lb & (x^r)^G  & (y_0^t)^G    \\
        |K|                  & 1 & \frac{1}{2}q(q\!-\!1) & \frac{1}{2}(q^2\!-\!1) & \frac{1}{2}(q^2\!-\!1) & q(q\!+\!1)  & q(q\!-\!1)  \\
       |\C_G(x_{\mbox{}_K})| & \frac{1}{2}q(q^2\!\!-\!1)&   q+1  & q  & q  & \frac{1}{2}(q\!-\!1) &  \frac{1}{2}(q\!+\!1)  \\
        x_{\mbox{}_K}        &  1      & z                     & u & v  &  x^r  & y_0^t \\
     \hline
    \vf_1^{\vphantom{A^A}}    & 1      &  1    &  1  &  1 & 1 & 1 \\
    \vf_2                     & q-1      & \mathllap{-}2    &  \mathllap{-}1    & \mathllap{-}1 & .  &  \mathllap{-}2  \\
    \vf_3                     & \frac{1}{2}(q-1)  & -(-1)^{(q+1)/4}    & \frac{1}{2}(-1\!+\!i\sqrt{q}) & \frac{1}{2}(-1\!-\!i\sqrt{q}) & . &-(-1)^t \\
    \vf_4                     & \frac{1}{2}(q-1)  & -(-1)^{(q+1)/4}    & \frac{1}{2}(-1\!-\!i\sqrt{q}) & \frac{1}{2}(-1\!+\!i\sqrt{q}) & . &-(-1)^t \\
    \psi_j                    & q+1    & . &  1 & 1 & \z^{rj}\!+\!\z^{-rj}  &  .  \\
    \eta_f                     & q-1    & -2\cdot(-1)^f   &  \mathllap{-}1 &  \mathllap{-}1  &  .  & -\xi_0^{tf}\!-\xi_0^{-tf}
  \end{array}
$$

$$ \begin{array}{c|rrrrrr}
  D & \vf_1   & \vf_2 & \vf_3 & \vf_4 &\psi_j  &\eta_f \\
     \hline
    \x_1^{\vphantom{A^A}} & 1   & . & . & . & . & .\\
    \x_2 & 1  & 1 & . & . & . & .\\
    \x_3 & . & . & 1 & . & . & . \\
    \x_4 & . & . & . & 1 & . & . \\
    \th_s & . & . & . & . & \delta_{sj} & . \\
    \t_{h_+} & .  & 1 & . & . & . & .\\
    \t_{h_-} & .  & . & 1 & 1 & . & .\\
    \t_{h_i} & .  & . & . & . & . & \delta_{if} \\
    \end{array} \qquad
  C=\left(
      \begin{array}{cccccc}
        2 & 1 & . & . & . & .\\[2pt]
        1 & \frac{p^d+1}{2} & . & . & . & .\\
         . & . & \frac{p^d+1}{2} & \frac{p^d-1}{2}  & . & .\\[2pt]
         . & . &  \frac{p^d-1}{2} & \frac{p^d+1}{2}  & . & .\\
         . & . & . & . & \I_{(q-3)/4}& . \\
         . & . & . & . & . & p^d \I_{(n/2)-1}
      \end{array}
    \right)
$$

$$
\begin{array}{c}
\mbox{Обозначения} \\[2pt]
 p^d=(q+1)_p,\ n=(\frac{q+1}{2})_{p'},\\
 y_0=y^{p^d},\ \xi_0=\xi^{p^d}=e^{2\pi i/n}\\
\end{array}\qquad
  \begin{array}{c}
    \mbox{Параметры} \\[3pt]
    r,s,j\ = \ 1,\ld,\frac{q-3}{4};\\
    i,f,t\ = \ 1,\ld,\frac{n}{2}-1;\\
    h_+=n,2n,\ld,\frac{p^d-1}{2}n;\\
    h_-=\frac{n}{2},n+\frac{n}{2},\ld,\frac{p^d-3}{2}n+\frac{n}{2};\\
    h_i=i,n\pm i,2n\pm i,\ld,\frac{p^d-1}{2}n\pm i
  \end{array}
$$

\medskip

$$
\begin{array}{l}
\Phi_2\big(PSL_2(q)\big)\\
q\equiv -1 \!\! \pmod 4\\
\end{array}
\qquad
\begin{array}{c|ccccc}
     K                     & 1a   & la & lb  & (x^r)^G   & (y_0^t)^G \\
  |K|                  & 1 & \frac{1}{2}(q^2\!-\!1) & \frac{1}{2}(q^2\!-\!1) & q(q\!+\!1)  & q(q\!-\!1)  \\
  |\C_G(x_{\mbox{}_K})| & \frac{1}{2}q(q^2\!\!-\!1)  & q  & q  & \frac{1}{2}(q\!-\!1) &  \frac{1}{2}(q\!+\!1)  \\
      x_{\mbox{}_K}        &  1                         & u & v  &  x^r  & y_0^t \\
   \hline
  \vf_1^{\vphantom{A^A}}   & 1          &  1  &  1 & 1 & 1 \\
  \vf_2                     & \frac{1}{2}(q-1)      & \frac{1}{2}(-1+i\sqrt{q}) & \frac{1}{2}(-1-i\sqrt{q}) & . & \mathllap{-}1  \\
  \vf_3                     & \frac{1}{2}(q-1)      & \frac{1}{2}(-1-i\sqrt{q}) & \frac{1}{2}(-1+i\sqrt{q}) & . & \mathllap{-}1  \\
  \psi_j                      & q+1     &  1 & 1 & \z^{rj}\!+\!\z^{-rj}  &  .  \\
  \eta_f                     & q-1    &  \mathllap{-}1 &  \mathllap{-}1  &  .  & -\xi_0^{tf}\!-\xi_0^{-tf}
\end{array}
$$

$$ \begin{array}{c|rrrrr}
  D & \vf_1   & \vf_2 & \vf_3 &\psi_j  &\eta_f \\
     \hline
    \x_1^{\vphantom{A^A}} & 1   & . & . & . & .\\
    \x_2 & 1  & 1 & 1 & .  & .\\
    \x_3 & . & 1 & . & . & . \\
    \x_4 & . & . & 1 & . & . \\
    \th_s & .  & . &  . & \delta_{sj} & .\\
    \t_{h_0} & .  & 1 & 1 & . & . \\
    \t_{h_i} & . & . & . & .  & \delta_{if}
    \end{array} \qquad
  C=\left(
      \begin{array}{ccccc}
        2 & 1 & 1 & . & . \\[2pt]
        1 & 2^{d-2}+1 & 2^{d-2} & . & . \\
        1 & 2^{d-2} & 2^{d-2}+1 & .  &  .\\[2pt]
         . & . &  . &  \I_{(q-3)/4}  & .\\
         . & . & . & . & 2^{d-1} \I_{(n-1)/2}
      \end{array}
    \right)
$$

$$
\begin{array}{c}
\mbox{Обозначения} \\[2pt]
 2^d=(q+1)_2,\ n=(q+1)_{2'},\\
 y_0=y^{2^{d-1}},\ \xi_0=\xi^{2^{d-1}}=e^{2\pi i/n}\\
\end{array}\qquad
  \begin{array}{c}
    \mbox{Параметры} \\[3pt]
    r,s,j\ = \ 1,\ld,\frac{q-3}{4};\\
    i,f,t\ = \ 1,\ld,\frac{n-1}{2};\\
    h_0=n,2n,\ld,(2^{d-2}-1)n;\\
    h_i=i,n\pm i,2n\pm i,\ld,(2^{d-2}-1)n\pm i,2^{d-2}n-i
  \end{array}
$$

\renewcommand{\bibname}{Рекомендуемая литература} 

\printglossary

\printindex[
  style=mcolindexgroupfixed,
  title={Предметный указатель}
]

\glsaddallunused[symbols]
{\setlength{\glsdescwidth}{0.78\textwidth}%
\printsymbols[style=symbolslong,title={Список обозначений}]}

\end{document}

%% file: mth_indx.tex
%
\rusnewterm{iAut}         {автоморфизм}
\rusnewterm{iAlg}         {алгебра}
\rusnewterm{iAnn}         {аннулятор}
\rusnewterm{iAntAutRng}   {антиавтоморфизм кольца}
\rusnewterm{iAntHomRng}   {антигомоморфизм колец}
\rusnewterm{iAntIsoRng}   {антиизоморфизм колец}
\rusnewterm{iBas}         {базис}
\rusnewterm        [sort= {бла}]      {iPBlk}{$p$-блок}
\rusnewterm{iBlkAlg}      {блок групповой алгебры}
\rusnewterm        [sort= {вес}]      {iPWght}{$p$-вес}
\rusnewterm        [sort= {вкл}]      {iGIncl}{$G$-включение}
\rusnewterm        [sort= {вкл}]      {iHtIrrOrdChr} {$p$-высота неприводимого обыкновенного характера}
\rusnewterm{iCnj}         {гипотеза}
\rusnewterm{iHom}         {гомоморфизм}
\rusnewterm{iGrBr}        {граф Брауэра}
\rusnewterm{iGr}          {группа}
\rusnewterm{iDblCoset}    {двойной смежный класс}
\rusnewterm        [sort= {деа}]      {iPDef}{$p$-дефект}
\rusnewterm{iAct}         {действие}
\rusnewterm{iBrTr}        {дерево Брауэра}
\rusnewterm{iDef}         {дефект}
\rusnewterm{iCmplFr}      {дополнение Фробениуса}
\rusnewterm{iZakVzFr}     {закон взаимности Фробениуса}
\rusnewterm{iClAlgFld}    {замыкание алгебраическое поля}
\rusnewterm{iIdl}         {идеал}
\rusnewterm{iIdmp}        {идемпотент}
\rusnewterm{iIdmpsOrt}    {идемпотенты ортогональные}
\rusnewterm{iIso}         {изоморфизм}
\rusnewterm{iInvCart}     {инварианты Картана}
\rusnewterm{iInvl}        {инволюция}
\rusnewterm{iCls}         {класс}
\rusnewterm{iRng}         {кольцо}
\rusnewterm{iCmt}         {коммутант}
\rusnewterm{ixy}          {коммутатор}
\rusnewterm{iComp}        {компонента}
\rusnewterm{iMltIrrComp}  {кратность неприводимой компоненты}
\rusnewterm{iLem}         {лемма}
\rusnewterm{iMat}         {матрица}
\rusnewterm{iMnChU}       {множество частично упорядоченное}
\rusnewterm{iMod}         {модуль}
\rusnewterm{iMono}        {мономорфизм}
\rusnewterm{iNrmtr}       {нормализатор}
\rusnewterm{iNorPCmpl}    {нормальное $p$-дополнение}
\rusnewterm{iSuppElmGrAlg}{носитель элемента групповой алгебры}
\rusnewterm{iDmn}         {область}
\rusnewterm{iGOrb}        {орбита, $G$-}
\rusnewterm{iOrthBlk}     {ортогональность в блоке}
\rusnewterm{iMap}         {отображение}
\rusnewterm        [sort= {поа}]      {iASubMod}{$A$-подмодуль}
\rusnewterm{iSubAlg}      {подалгебра}
\rusnewterm{iSubGr}       {подгруппа}
\rusnewterm{iSubRng}      {подкольцо}
\rusnewterm{iSbBnd}       {подмножество, ограниченное сверху}
\rusnewterm{iSubMod}      {подмодуль}
\rusnewterm{iSubFld}      {подполе}
\rusnewterm{iFld}         {поле}
\rusnewterm{iGenAMod}     {порождающий $A$-модуля}
\rusnewterm{iRep}         {представление}
\rusnewterm{iRepsEqv}     {представления эквивалентные}
\rusnewterm{iPrmRt}       {примитивный корень}
\rusnewterm{iProd}        {произведение}
\rusnewterm{iRad}         {радикал}
\rusnewterm{iRadPSubGr}   {радикальная $p$-подгруппа}
\rusnewterm        [sort= {ран}]           {iZZRnk}{$\mathbb{Z}$-ранг}
\rusnewterm{iExtFld}      {расширение полей}
\rusnewterm{iSer}         {ряд}
\rusnewterm{iSersCompEqv} {ряды композиционные эквивалентные}
\rusnewterm        [sort= {сеч}]        {iPSecGr}{$p$-сечение группы}
\rusnewterm{idij}         {символ Кронекера}
\rusnewterm{iAtrA}        {симметрическая разность}
\rusnewterm{iImpSstMod}   {система импримитивности модуля}
\rusnewterm{iTrcMat}      {след матрицы}
\rusnewterm{iRelOrth}     {соотношение ортогональности}
\rusnewterm{iDeg}         {степень}
\rusnewterm{iStrConsZRG}  {структурные константы центра групповой алгебры}
\rusnewterm{iNlpCls}      {ступень нильпотентности}
\rusnewterm{iSum}         {сумма}
\rusnewterm{iTbl}         {таблица}
\rusnewterm{iSqrng}       {тело}
\rusnewterm{iThm}         {теорема}
\rusnewterm{iFct}         {фактор}
\rusnewterm{iQuoAlg}      {факторалгебра}
\rusnewterm{iQMod}        {фактормодуль}
\rusnewterm{iFormBurn}    {формула Бернсайда}
\rusnewterm{iFunc}        {функция}
\rusnewterm{iChr}         {характер}
\rusnewterm{iCen}         {центр}
\rusnewterm{iCntrl}       {централизатор}
\rusnewterm{iChn}         {цепь}
\rusnewterm{iPrt}         {часть}
\rusnewterm{iNumsDec}     {числа разложения}
\rusnewterm{iNum}         {число}
\rusnewterm{iExpGr}       {экспонента группы}
\rusnewterm{iElm}         {элемент}
\rusnewterm{iElms}        {элементы}
\rusnewterm{iEnd}         {эндоморфизм}
\rusnewterm{iEpi}         {эпиморфизм}
\rusnewterm{iKer}         {ядро}

\rusnewterm[           parent={iRep}]{iFRepGr}{группы, $F$-}

\rusnewterm[           parent={iAct}]{iActAlg}{алгебры}
\rusnewterm[           parent={iAct}]{iActGr}{группы}
\rusnewterm[           parent={iActGr}]{iActGrDblTrn}{дважды транзитивное}
\rusnewterm[           parent={iActGr}]{iActGrRReg}{правое регулярное}
\rusnewterm[           parent={iActGr}]{iActGrTr}{транзитивное}
\rusnewterm[           parent={iAlg}]{iAlgGr}{групповая}
\rusnewterm[           parent={iAlg}]{iAlgRng}{над кольцом $R$}
\rusnewterm[sort={зна},parent={iAlg}]{iAlgRvFuncs}{$R$-значных функций}
\rusnewterm[           parent={iAlg}]{iFAlgClsFncGr}{классовых функций на группе, $F$-}
\rusnewterm[           parent={iAlg}]{iFAlgSemSmp}{полупростая, $F$-}
\rusnewterm[           parent={iAlg}]{iFAlgSmp}{простая, $F$-}
\rusnewterm[sort={алб},parent={iAlg}]{iRAlg}{$R$-алгебра}
\rusnewterm[sort={мод},parent={iAnn}]{iAnnAmod}{$A$-модуля}
\rusnewterm[           parent={iAnn}]{iAnnElmAMod}{элемента $A$-модуля}
\rusnewterm[           parent={iAut}]{iAutRng}{кольца}
\rusnewterm[           parent={iBas}]{iBasDl}{двойственный}
\rusnewterm[sort={мод},parent={iBas}]{iBasMod}{модуля, $R$-}
\rusnewterm[           parent={iCen}]{iCenGr}{группы}
\rusnewterm[           parent={iCen}]{iCenRng}{кольца}
\rusnewterm[           parent={iChr}]{iChrAlgCnj}{алгебраически сопряжённый}
\rusnewterm[sort={мод},parent={iChr}]{iChrAMod}{$A$-модуля}
\rusnewterm[           parent={iChr}]{iChrBr}{брауэров}
\rusnewterm[           parent={iChr}]{iChrCen}{центральный, соответствующий}
\rusnewterm[           parent={iChr}]{iChrCompCnj}{комплексно сопряжённый}
\rusnewterm[           parent={iChr}]{iChrGnrl}{обобщённый}
\rusnewterm[           parent={iChr}]{iChrOrd}{обыкновенный}
\rusnewterm[sort={мое},parent={iChr}]{iChrPMod}{$p$-модулярный}
\rusnewterm[           parent={iChr}]{iChrRepAlg}{представления алгебры}
\rusnewterm[           parent={iChr}]{iPrincIndecChr}{главный неразложимый}
\rusnewterm[           parent={iChrBr}]{iChrBrIrr}{неприводимый}
\rusnewterm[           parent={iChrBr}]{iChrBrPrnc}{главный}
\rusnewterm[           parent={iChrCen}]{iChrCenCorBlk}{блоку}
\rusnewterm[           parent={iChrCen}]{iChrCenCorIrr}{неприводимому характеру}
\rusnewterm[           parent={iChrOrd}]{iChrOrdFth}{точный}
\rusnewterm[           parent={iChrOrd}]{iChrOrdHmg}{однородный}
\rusnewterm[           parent={iChrOrd}]{iChrOrdImpr}{импримитивный}
\rusnewterm[           parent={iChrOrd}]{iChrOrdIrr}{неприводимый}
\rusnewterm[           parent={iChrOrd}]{iChrOrdMnm}{мономиальный}
\rusnewterm[           parent={iChrOrd}]{iChrOrdPrm}{примитивный}
\rusnewterm[           parent={iCls}]{iClsCnj}{сопряжённости}
\rusnewterm[           parent={iCls}]{iClsDefPBlk}{дефектный $p$-блока}
\rusnewterm[sort={рег},parent={iClsCnj}]{iClsCnjPReg}{$p$-регулярный}
\rusnewterm[sort={й  },parent={iCmt}]{iCmtIth}{$i$-й}
\rusnewterm[           parent={iCnj}]{iCnjBrZerHt}{Брауэра о нулевой высоте}
\rusnewterm[           parent={iCnj}]{iCnjDade}{Дэйда}
\rusnewterm[           parent={iCnj}]{iCnjCharTr}{о тройках характеров}
\rusnewterm[           parent={iCnj}]{iConjAlpWt}{Альперина о весах}
\rusnewterm[           parent={iCntrl}]{iCntrlAlg}{в алгебре}
\rusnewterm[           parent={iCntrl}]{iCntrlGr}{в группе}
\rusnewterm[           parent={iComp}]{iCompIrr}{неприводимая}
\rusnewterm[           parent={iComp}]{iCompWedAMod}{Веддерберна $A$-модуля}
\rusnewterm[           parent={iComp}]{iHmgCompAMod}{однородная $A$-модуля}
\rusnewterm[           parent={iCompIrr}]{iCompIrrOrdChr}{обыкновенного характера}
\rusnewterm[           parent={iCompIrr}]{iCompIrrRep}{представления}
\rusnewterm[           parent={iCompIrrRep}]{iCompIrrRepLow}{нижняя}
\rusnewterm[           parent={iCompIrrRep}]{iCompIrrRepUpp}{верхняя}
\rusnewterm[           parent={iDef}]{iDefCls}{класса}
\rusnewterm[sort={бло},parent={iDef}]{iDefPBlk}{$p$-блока}
\rusnewterm[           parent={iDeg}]{iDegChr}{характера}
\rusnewterm[           parent={iDeg}]{iDegRep}{представления}
\rusnewterm[           parent={iDeg}]{iDegExt}{расширения полей}
\rusnewterm[           parent={iDegChr}]{iDegChrBr}{брауэрова}
\rusnewterm[           parent={iDegChr}]{iDegChrOrd}{обыкновенного}
\rusnewterm[           parent={iDmn}]{iDdkDmn}{дедекиндова}
\rusnewterm[           parent={iDmn}]{iIntDmn}{целостности}
\rusnewterm[           parent={iElm}]{iElmAlg}{алгебраический}
\rusnewterm[sort={p  },parent={iElm}]{iPPrElm}{$p'$-элемент}
\rusnewterm[sort={pi },parent={iElm}]{iPiElm}{$\pi$-элемент}
\rusnewterm[           parent={iElm}]{iElmKaz}{Казимира}
\rusnewterm[sort={рег},parent={iElm}]{iElmPReg}{$p$-регулярный}
\rusnewterm[           parent={iElm}]{iElmRng}{кольца}
\rusnewterm[           parent={iElm}]{iElmMax}{максимальный}
\rusnewterm[           parent={iElmRng}]{iInvElm}{обратимый (справа, слева)}
\rusnewterm[           parent={iElmRng}]{iNlpElm}{нильпотентный}
\rusnewterm[           parent={iElmRng}]{iObrElm}{обратный (правый, левый)}
\rusnewterm[           parent={iElms}]{iCnjElms}{сопряжённые}
\rusnewterm[           parent={iEnd}]{iAEnd}{$A$-эндоморфизм}
\rusnewterm[sort={мод},parent={iEnd}]{iEndRmod}{$R$-модулей}
\rusnewterm[           parent={iEpi}]{iEpiRng}{колец}
\rusnewterm[           parent={iExtFld}]{iExtFldAlg}{алгебраическое}
\rusnewterm[           parent={iExtFld}]{iExtFldFin}{конечное}
\rusnewterm[           parent={iExtFld}]{iExtFldNor}{нормальное}
\rusnewterm[           parent={iChr}]{iChrPerm}{подстановочный}
\rusnewterm[           parent={iChr}]{iFChrGen}{обобщённый, $F$-}
\rusnewterm[           parent={iChr}]{iFChrGr}{группы, $F$-}
\rusnewterm[           parent={iChr}]{iFChrInd}{индуцированный, $F$-}
\rusnewterm[           parent={iChr}]{iFChrLin}{линейный, $F$-}
\rusnewterm[           parent={iChr}]{iFChrReg}{регулярный, $F$-}
\rusnewterm[           parent={iChr}]{iChrRepGr}{представления группы}
\rusnewterm[           parent={iChrPerm}]{iChrPermFG}{групповой алгебры}
\rusnewterm[           parent={iChrPerm}]{iFChrAct}{действия группы}
\rusnewterm[           parent={iFChrGr}]{iFChrGrPrnc}{главный}
\rusnewterm[           parent={iFct}]{iFctCmpAMod}{композиционный $A$-модуля}
\rusnewterm[           parent={iFct}]{iQuoSerAMod}{ряда $A$-модулей}
\rusnewterm[           parent={iFld}]{iFldAlgCl}{алгебраически замкнутое}
\rusnewterm[           parent={iFld}]{iFldCcl}{круговое}
\rusnewterm[           parent={iFld}]{iFldDec}{разложения}
\rusnewterm[           parent={iFldDec}]{iFldDecGr}{группы}
\rusnewterm[           parent={iFldDec}]{iFldDecPol}{многочлена}
\rusnewterm[           parent={iFld}]{iFldPrim}{простое}
\rusnewterm[           parent={iFRepGr}]{iFRepGrReg}{регулярное}
\rusnewterm[           parent={iFunc}]{iFuncCls}{классовая}
\rusnewterm[           parent={iFunc}]{iFuncEul}{Эйлера}
\rusnewterm[           parent={iFuncCls}]{iFuncClsInd}{индуцированная}
\rusnewterm[sort={рег},parent={iFuncCls}]{iFuncClsPReg}{$p$-регулярная}
\rusnewterm[           parent={iFuncCls}]{iFuncClsConjElm}{сопряжённая элементом}
\rusnewterm[           parent={iFuncCls}]{iFuncClsConjAut}{сопряжённая автоморфизмом}
\rusnewterm[           parent={iGr}]{iGrGalExt}{Галуа расширения полей}
\rusnewterm[           parent={iGr}]{iGLnR}{полная линейная}
\rusnewterm[           parent={iGr}]{iGrDefPBlk}{дефектная $p$-блока}
\rusnewterm[           parent={iGr}]{iGrInr}{инерции}
\rusnewterm[           parent={iGr}]{iGrFr}{Фробениуса}
\rusnewterm[           parent={iGr}]{iGrInvElms}{обратимых элементов кольца}
\rusnewterm[           parent={iGr}]{iGrMltRng}{мультипликативная кольца}
\rusnewterm[           parent={iGr}]{iGrNlp}{нильпотентная}
\rusnewterm[sort={деф},parent={iGr}]{iGrPDefCls}{$p$-дефектная класса}
\rusnewterm[           parent={iGr}]{iGrElmBr}{элементарная (брауэрова)}
\rusnewterm[sort={эле},parent={iGr}]{iGrPElmBr}{$p$-элементарная (брауэрова)}
\rusnewterm[sort={ква},parent={iGr}]{iGrPQsElm}{$p$-квазиэлементарная}
\rusnewterm[sort={кв },parent={iGr}]{iGrpQsElm}{квазиэлементарная}
\rusnewterm[           parent={iGr}]{iGrpQsSimp}{квазипростая}
\rusnewterm[           parent={iGr}]{iGrPSol}{$p$-разрешимая}
\rusnewterm[           parent={iGr}]{iGrQuat}{кватернионов}
\rusnewterm[           parent={iGr}]{iGrSol}{разрешимая}
\rusnewterm[           parent={iGr}]{iGrTrv}{тривиальная}
\rusnewterm[           parent={iGr}]{iMGr}{$M$-группа}
\rusnewterm[           parent={iGrInr}]{iGrInrChr}{характера}
\rusnewterm[           parent={iGrInr}]{iGrInrMod}{модуля}
\rusnewterm[           parent={iGrInr}]{iGrInrRep}{представления}
\rusnewterm[           parent={iGrQuat}]{iGrQuatGen}{обобщённая}
\rusnewterm[           parent={iHom}]{iAHom}{$A$-гомоморфизм}
\rusnewterm[sort={алг},parent={iHom}]{iHomRAlg}{$R$-алгебр}
\rusnewterm[sort={мод},parent={iHom}]{iHomRMod}{$R$-модулей}
\rusnewterm[           parent={iHom}]{iHomRng}{колец}
\rusnewterm[           parent={iIdl}]{iIdlFrc}{дробный}
\rusnewterm[           parent={iIdl}]{iIdlGenX}{порождённый множеством}
\rusnewterm[           parent={iIdl}]{iIdlInt}{целый}
\rusnewterm[           parent={iIdl}]{iIdlNlp}{нильпотентный}
\rusnewterm[sort={алг},parent={iIdl}]{iIdlRAlg}{$R$-алгебры левый (правый, двусторонний)}
\rusnewterm[           parent={iIdl}]{iIdlRng}{кольца}
\rusnewterm[           parent={iIdl}]{iIdlSmp}{простой}
\rusnewterm[           parent={iIdlRng}]{iIdlRngL}{левый}
\rusnewterm[           parent={iIdlRng}]{iIdlRngR}{правый}
\rusnewterm[           parent={iIdlRng}]{iIdlRngTs}{двусторонний}
\rusnewterm[           parent={iIdmp}]{iIdmpCen}{центральный}
\rusnewterm[           parent={iIdmp}]{iIdmpOsm}{Осимы}
\rusnewterm[           parent={iIdmp}]{iIdmpPrim}{примитивный}
\rusnewterm[           parent={iIdmpCen}]{iIdmpCenCG}{алгебры $\mathbb{C}G$}
\rusnewterm[sort={алг},parent={iIdmpCen}]{iIdmpCenFAlg}{$F$-алгебры}
\rusnewterm[           parent={iIdmpCen}]{iIdmpCenFGCorr}{алгебры $FG$, соответствующий $p$-блоку}
\rusnewterm[           parent={iIntDmn}]{iIntClsdDmn}{целозамкнутая}
\rusnewterm[           parent={iInvl}]{iInvlIsl}{изолированная}
\rusnewterm[           parent={iIso}]{iAIso}{$A$-изоморфизм}
\rusnewterm[sort={мод},parent={iIso}]{iIsoRMod}{$R$-модулей}
\rusnewterm[           parent={iIso}]{iIsoRng}{колец}
\rusnewterm[           parent={iKer}]{iKerChr}{характера}
\rusnewterm[           parent={iKer}]{iKerHom}{гомоморфизма}
\rusnewterm[           parent={iKer}]{iKerFr}{Фробениуса}
\rusnewterm[           parent={iKerHom}]{iKerGrHom}{группового}
\rusnewterm[           parent={iKerHom}]{iKerRng}{кольцевого}
\rusnewterm[           parent={iKerChr}]{iKerChrOrd}{обыкновенного}
\rusnewterm[           parent={iLem}]{iLemNak}{Накаямы}
\rusnewterm[           parent={iLem}]{iLemBra}{Брауэра}
\rusnewterm[           parent={iLem}]{iLemPBur}{Бернсайда, $p^a$-}
\rusnewterm[           parent={iLem}]{iLemSchr}{Шура}
\rusnewterm[           parent={iLem}]{iLemZorn}{Цорна}
\rusnewterm[sort={лин},parent={iMap}]{iMapRlin}{$R$-линейное}
\rusnewterm[           parent={iMap}]{iMapRob}{Робинсона}
\rusnewterm[           parent={iMat}]{iMatCart}{Картана}
\rusnewterm[           parent={iMat}]{iMatEmp}{пустая}
\rusnewterm[           parent={iMat}]{iMatDec}{разложения}
\rusnewterm[           parent={iMat}]{iMatInvTr}{обратно-транспонированная}
\rusnewterm[           parent={iMat}]{iMatMnm}{мономиальная}
\rusnewterm[sort={бло},parent={iMatCart}]{iMatCartPBlk}{$p$-блока}
\rusnewterm[sort={бло},parent={iMatDec}]{iMatDecPBlk}{$p$-блока}
\rusnewterm[           parent={iMod}]{iAMod}{$A$-модуль}
\rusnewterm[           parent={iMod}]{iAModCompRed}{вполне приводимый, $A$-}
\rusnewterm[           parent={iMod}]{iAModIrr}{неприводимый, $A$-}
\rusnewterm[           parent={iMod}]{iAModRed}{приводимый, $A$-}
\rusnewterm[sort={ре1},parent={iMod}]{iAModReg}{регулярный, $A$-}
\rusnewterm[sort={ре2},parent={iMod}]{iRModReg}{регулярный, $R$-}
\rusnewterm[           parent={iMod}]{iAModSemSmp}{полупростой, $A$-}
\rusnewterm[           parent={iMod}]{iAmodSmp}{простой, $A$-}
\rusnewterm[           parent={iMod}]{iMidImpr}{импримитивный}
\rusnewterm[           parent={iMod}]{iModConj}{сопряжённый}
\rusnewterm[           parent={iMod}]{iModCorRep}{соответствующий представлению}
\rusnewterm[           parent={iMod}]{iModFinGen}{конечно порождённый}
\rusnewterm[           parent={iMod}]{iModInd}{индуцированный}
\rusnewterm[           parent={iMod}]{iModPrim}{примитивный}
\rusnewterm[           parent={iMod}]{iModRAlg}{над $R$-алгеброй}
\rusnewterm[           parent={iMod}]{iModRG}{групповой алгебры}
\rusnewterm[           parent={iMod}]{iRModDec}{разложимый, $R$-}
\rusnewterm[           parent={iMod}]{iRModDl}{двойственный, $R$-}
\rusnewterm[           parent={iMod}]{iRModFr}{свободный, $R$-}
\rusnewterm[           parent={iMod}]{iRModIndec}{неразложимый, $R$-}
\rusnewterm[           parent={iMod}]{iRModL}{левый, $R$-}
\rusnewterm[           parent={iMod}]{iRModR}{правый, $R$-}
\rusnewterm[           parent={iModConj}]{iModConjElm}{элементом}
\rusnewterm[           parent={iModConj}]{iModConjAut}{автоморфизмом}
\rusnewterm[           parent={iModRG}]{iModRGPerm}{подстановочный}
\rusnewterm[           parent={iModRG}]{iModRGPrinc}{главный}
\rusnewterm[           parent={iModRG}]{iModRGTriv}{тривиальный}
\rusnewterm[           parent={iModRG}]{iFGModCGrad}{контрагредиентный}
\rusnewterm[           parent={iModRGPerm}]{iModRGPermNat}{естественный}
\rusnewterm[           parent={iMono}]{iMonoRng}{колец}
\rusnewterm[           parent={iNum}]{iBiprNum}{бипримарное}
\rusnewterm[           parent={iNum}]{iNumAlg}{алгебраическое}
\rusnewterm[           parent={iNum}]{iNumIntAlg}{целое алгебраическое}
\rusnewterm[           parent={iNum}]{iNumIntRat}{целое рациональное}
\rusnewterm[           parent={iNum}]{iPrimNum}{примарное}
\rusnewterm[           parent={iNumsDec}]{iNumsDecHghr}{высшие}
\rusnewterm[sort={бло},parent={iNumsDec}]{iNumsDecPBlk}{$p$-блока}
\rusnewterm[           parent={iGOrb}]{iGOrbEtl}{элемента}
\rusnewterm[           parent={iPBlk}]{iPBlkBrChr}{брауэровых характеров}
\rusnewterm[           parent={iPBlk}]{iPBlkInd}{индуцированный}
\rusnewterm[           parent={iPBlk}]{iPBlkOrdChr}{обыкновенных характеров}
\rusnewterm[           parent={iPBlk}]{iPBlkPrnc}{главный}
\rusnewterm[           parent={iPDef}]{iPDefCls}{класса}
\rusnewterm[           parent={iPDef}]{iPDefOrdChr}{неприводимого обыкновенного характера}
\rusnewterm[           parent={iPrt}]{iBPrtClsFunc}{классовой функции, $B$-}
\rusnewterm[sort={чi },parent={iPrt}]{iPiPrtNum}{числа, $\pi$-}
\rusnewterm[sort={ч1 },parent={iPrt}]{iPPrPrtNum}{числа, $p'$-}
\rusnewterm[sort={ч0 },parent={iPrt}]{iPPrtNum}{числа, $p$-}
\rusnewterm[sort={э1 },parent={iPrt}]{iPPrPrtElm}{элемента, $p'$-}
\rusnewterm[sort={э0 },parent={iPrt}]{iPPrtElm}{элемента, $p$-}
\rusnewterm[           parent={iProd}]{iKronProdMat}{матриц кронекерово}
\rusnewterm[           parent={iProd}]{iProdDirClsFuncs}{прямое классовых функций}
\rusnewterm[           parent={iProd}]{iProdIds}{идеалов}
\rusnewterm[           parent={iProd}]{iProdTens}{тензорное}
\rusnewterm[           parent={iProd}]{iSclPrdCompClsFnc}{скалярное комплекснозначных классовых функций}
\rusnewterm[sort={мод},parent={iProdTens}]{iProdTensFMod}{$FG$-модулей}
\rusnewterm[           parent={iProdTens}]{iTensProdRep}{представлений}
\rusnewterm[sort={вес},parent={iPWght}]{iPWghtBlgBlk}{принадлежащий блоку}
\rusnewterm[           parent={iRad}]{iRadJac}{Джекобсона}
\rusnewterm[sort={алг},parent={iRad}]{iRadRAlg}{$R$-алгебры}
\rusnewterm[           parent={iRad}]{iRadRng}{кольца}
\rusnewterm[           parent={iRelOrth}]{iRelOrthFrs}{первое}
\rusnewterm[           parent={iRelOrth}]{iRelOrthGen}{обобщённое}
\rusnewterm[           parent={iRelOrth}]{iRelOrthScd}{второе}
\rusnewterm[           parent={iRep}]{iRepAlgMat}{алгебры матричное}
\rusnewterm[           parent={iRep}]{iRepCntGrd}{контрагредиентное}
\rusnewterm[           parent={iRep}]{iRepConj}{сопряжённое}
\rusnewterm[           parent={iRep}]{iRepNul}{нулевое}
\rusnewterm[           parent={iRep}]{iRepCorMod}{соответствующее модулю}
\rusnewterm[sort={алг},parent={iRep}]{iRepFAlg}{$F$-алгебры}
\rusnewterm[           parent={iRep}]{iRepGr}{группы}
\rusnewterm[           parent={iRep}]{iRepInd}{индуцированное}
\rusnewterm[           parent={iRep}]{iRepLin}{линейное}
\rusnewterm[           parent={iRep}]{iRepOrd}{обыкновенное}
\rusnewterm[           parent={iRepConj}]{iRepConjElm}{элементом}
\rusnewterm[           parent={iRepConj}]{iRepConjAut}{автоморфизмом}
\rusnewterm[           parent={iRepFAlg}]{iRepFAlgIrr}{неприводимое}
\rusnewterm[           parent={iRepFAlg}]{iRepFAlgRed}{приводимое}
\rusnewterm[           parent={iRepFAlg}]{iRepFAlgCrd}{вполне приводимое}
\rusnewterm[           parent={iRepFAlg}]{iRepFAlgInd}{неразложимое}
\rusnewterm[           parent={iRepFAlg}]{iRepFAlgReg}{регулярное}
\rusnewterm[           parent={iRepGr}]{iRepGrAbsIrr}{абсолютно неприводимое}
\rusnewterm[           parent={iRepGr}]{iRepGrZap}{записанное над}
\rusnewterm[           parent={iRepGr}]{iRepGrFth}{точное}
\rusnewterm[           parent={iRepGr}]{iRepGrPrnc}{главное}
\rusnewterm[           parent={iRepGr}]{iRepGrTrv}{тривиальное}
\rusnewterm[           parent={iRepGrZap}]{iRepGrZapR}{\ \ подкольцом}
\rusnewterm[           parent={iRepGrZap}]{iRepGrZapF}{\ \ подполем}
\rusnewterm[           parent={iRng}]{iRngGenFChr}{обобщённых $F$-характеров}
\rusnewterm[           parent={iRng}]{iRngIntFld}{целых величин поля}
\rusnewterm[           parent={iRng}]{iRngLoc}{локальное}
\rusnewterm[           parent={iRng}]{iRngOp}{противоположное}
\rusnewterm[           parent={iRng}]{iRngZer}{нулевое}
\rusnewterm[           parent={iSer}]{iSerCompAMod} {ряд композиционный для $A$-модуля}
\rusnewterm[           parent={iSer}]{iSerNor} {ряд нормальный}
\rusnewterm[           parent={iSubFld}]{iSubFldGenX}{порождённое подмножеством}
\rusnewterm[           parent={iSubGr}]{iSubGrCen}{центральная}
\rusnewterm[           parent={iSubGr}]{iSubGrFrt}{Фраттини}
\rusnewterm[           parent={iSubGr}]{iSubGrGenX}{порождённая подмножеством}
\rusnewterm[           parent={iSubGr}]{iSubGrMax}{максимальная}
\rusnewterm[           parent={iSubGr}]{iSubGrMin}{минимальная}
\rusnewterm[           parent={iSubMod}]{iSubModGenX}{порождённый множеством}
\rusnewterm[           parent={iSubMod}]{iSubModPermMod}{подстановочного модуля}
\rusnewterm[           parent={iSubModPermMod}]{iSubModPermModDiag}{диагональный}
\rusnewterm[           parent={iSubModPermMod}]{iSubModPermModDiff}{разностный}
\rusnewterm[           parent={iSubRng}]{iSubRngGenX}{порождённое подмножеством}
\rusnewterm[           parent={iSum}]{iClsSum}{классовая}
\rusnewterm[           parent={iSum}]{iSumClsK}{класса $K$}
\rusnewterm[           parent={iSum}]{iSumDir}{прямая}
\rusnewterm[           parent={iSum}]{iSumIds}{идеалов}
\rusnewterm[sort={под},parent={iSum}]{iSumRSubMods}{$R$-подмодулей}
\rusnewterm[           parent={iSumDir}]{iSumDirIds}{идеалов}
\rusnewterm[           parent={iSumDir}]{iSumDirRep}{представлений}
\rusnewterm[sort={мод},parent={iSumDir}]{iSumDirRMods}{$R$-модулей}
\rusnewterm[           parent={iSumDir}]{iSumDirRng}{колец}
\rusnewterm[sort={под},parent={iSumDir}]{iSumDirRSubMods}{$R$-подмодулей}
\rusnewterm[           parent={iTbl}]{iTblBrChr}{брауэровых характеров}
\rusnewterm[           parent={iTbl}]{iTblCntChr}{центральных характеров}
\rusnewterm[           parent={iTbl}]{iTblOrdChr}{обыкновенных характеров}
\rusnewterm[           parent={iThm}]{iThmBra}{Брауэра}
\rusnewterm[           parent={iThm}]{iThmBraSuz}{Брауэра--Сузуки}
\rusnewterm[           parent={iThm}]{iThmCliff}{Клиффорда}
\rusnewterm[           parent={iThm}]{iThmDade}{Дэйда}
\rusnewterm[           parent={iThm}]{iThmDblCntr}{о двойном централизаторе}
\rusnewterm[           parent={iThm}]{iThmGrn}{Грина}
\rusnewterm[           parent={iThm}]{iThmIto}{Ито}
\rusnewterm[           parent={iThm}]{iThmJorHol}{Жордана--Гёльдера}
\rusnewterm[           parent={iThm}]{iThmKnr}{Кнёрра}
\rusnewterm[           parent={iThm}]{iThmMacK}{Макки}
\rusnewterm[           parent={iThm}]{iThmMasch}{Машке}
\rusnewterm[           parent={iThm}]{iThmMinMax}{min--max}
\rusnewterm[           parent={iThm}]{iThmPQBur}{Бернсайда, $p^a q^b$-}
\rusnewterm[           parent={iThm}]{iThmRobn}{Робинсона}
\rusnewterm[           parent={iThm}]{iThmWedd}{Веддерберна}
\rusnewterm[           parent={iThm}]{iThmFrob}{Фробениуса}
\rusnewterm[           parent={iThm}]{iThmZGla}{Глаубермана, $\Z^*$-}
\rusnewterm[           parent={iThmBra}]{iThmBraFstMain}{первая основная}
\rusnewterm[           parent={iThmBra}]{iThmBraScdMain}{вторая основная}
\rusnewterm[           parent={iThmBra}]{iThmBraThdMain}{третья основная}
\rusnewterm[           parent={iThmCliff}]{iThmCliffChar}{для характеров}

%% file: mth_symb.tex

\newglossaryentry{aaa}{name={\ensuremath{\st}},                                      sort={00 },description={естественный эпиморфизм $\ov{\ZZ}\to \ov{\ZZ}/M$}}
\newglossaryentry{Ast}{name={\ensuremath{A^\st}},                                    sort={01 },description={матрица $(\a_{ij}^\st)$, где $A=(\a_{ij})$}}
\newglossaryentry{Xst}{name={\ensuremath{\X^\st}},                                   sort={02 },description={$F$-представление, построенное по обыкновенному представлению $\X$, записанному над $R$}}
\newglossaryentry{Mst}{name={\ensuremath{M^*}},                                      sort={03 },description={$R$-модуль, двойственный к $R$-модулю $M$}}
\newglossaryentry{VCGrad}{name={\ensuremath{V^*}},                                   sort={04 },description={модуль, контрагредиентный к $FG$-модулю $V$}}
\newglossaryentry{Xcon}{name={\ensuremath{\X^*}},                                    sort={05 },description={представление, контрагредиентное к представлению $\X$}}
\newglossaryentry{IdR}{name={\ensuremath{\nor}},                                     sort={06 },description={двусторонний идеал кольца}}
\newglossaryentry{IdRl}{name={\ensuremath{\nor_l}},                                  sort={07 },description={левый идеал кольца}}
\newglossaryentry{IdRr}{name={\ensuremath{\nor_r}},                                  sort={08 },description={правый идеал кольца}}
\newglossaryentry{XleGY}{name={\ensuremath{X\le_G Y}},                               sort={09 },description={сокращение для <<$X$\  $G$-содержится в $Y$>>}}
\newglossaryentry{AoB}{name={\ensuremath{A\ot B}},                                   sort={0a },description={кронекерово произведение матриц $A$ и $B$}}
\newglossaryentry{VotW}{name={\ensuremath{V\ot W}},                                  sort={0b },description={тензорное произведение $FG$-модулей $V$ и $W$}}
\newglossaryentry{XotY}{name={\ensuremath{\X\ot\Y}},                                 sort={0c },description={тензорное произведение представлений $\X$ и $\Y$}}
\newglossaryentry{fpsi}{name={\ensuremath{\vf\times \psi}},                          sort={0d },description={прямое произведение классовых функций}}                                                           %
\newglossaryentry{lvfpsrsG}{name={\ensuremath{(\vf,\psi)_{\mbox{}_G}}},              sort={0e },description={скалярное произведение классовых функций $\vf,\psi\in \cf(G)$}}
\newglossaryentry{lvfpsrsGpp}{name={\ensuremath{(\vf,\psi)_{\mbox{}_{G_{p'}}}}},     sort={0f },description={скалярное произведение ограничения на $G_{p'}$ классовых функций $\vf,\psi\in \cf(G)\cup\cf(G_{p'})$},text={\ensuremath{(\mu,\nu)_{\mbox{}_{G_{p'}}}}}}
\newglossaryentry{lmnr}{name={\ensuremath{(m,n)}},                                   sort={0g },description={наибольший общий делитель целых чисел $m$ и $n$}}
\newglossaryentry{laXra}{name={\ensuremath{\la X \ra}},                              sort={0h },description={подгруппа, порождённая множеством $X$}}
\newglossaryentry{AtrB}{name={\ensuremath{A\triangle B}},                            sort={0i },description={симметрическая разность множеств $A$ и $B$}}
\newglossaryentry{gmod}{name={\ensuremath{|g|}},                                     sort={0j },description={порядок элемента $g$ группы или модуль комплексного числа $g$}}
\newglossaryentry{xy}{name={\ensuremath{[x,y]}},                                     sort={0k },description={коммутатор $x^{-1}y^{-1}xy$}}
\newglossaryentry{Gpr}{name={\ensuremath{G'}},                                       sort={0l },description={коммутант группы $G$}}
\newglossaryentry{Glir}{name={\ensuremath{G^{(i)}}},                                 sort={0m },description={$i$-й коммутант группы $G$}}
\newglossaryentry{Sx}{name={\ensuremath{S[X]}},                                      sort={0n },description={кольцо, порождённое кольцом $S$ и множеством $X$}}
\newglossaryentry{Kx}{name={\ensuremath{K(X)}},                                      sort={0o },description={поле, порождённое полем $K$ и множеством $X$}}
\newglossaryentry{Rr0}{name={\ensuremath{R^\circ}},                                  sort={0p },description={регулярный правый $R$-модуль кольца $R$}}
\newglossaryentry{Rl0}{name={\ensuremath{\mbox{}^\circ\!R}},                         sort={0q },description={регулярный левый $R$-модуль кольца $R$}}
\newglossaryentry{A0}{name={\ensuremath{A^\circ}},                                   sort={0r },description={регулярный $A$-модуль $R$-алгебры $A$}}
\newglossaryentry{Rx}{name={\ensuremath{R^\times}},                                  sort={0s },description={группа обратимых элементов кольца $R$}}
\newglossaryentry{aov}{name={\ensuremath{\ov{A}}},                                   sort={0t },description={матрица $(\ov{\a_{ij}})$, где $(\a_{ij})=A$}}
\newglossaryentry{chiov}{name={\ensuremath{\ov{\x}}},                                sort={0u },description={комплексно сопряжённый характер}}
\newglossaryentry{Kwh}{name={\ensuremath{\wh{K}}},                                   sort={0v },description={сумма элементов класса сопряжённости $K$}}
\newglossaryentry{thh}{name={\ensuremath{\hat{\th}}},                                sort={0w },description={ограничение $\th$ на $G_{p'}$}}
\newglossaryentry{chkth}{name={\ensuremath{\chk\th}},                                sort={0x },description={отображение, задаваемое правилом $\chk\th(g)=\th(g_{p'})$}}

\newglossaryentry{0R}{name={\ensuremath{0_R}},                                       sort={10 },description={ноль кольца $R$}}
\newglossaryentry{On}{name={\ensuremath{\O_{n}}},                                    sort={11 },description={нулевая $n\times n$-матрица}}
\newglossaryentry{1G}{name={\ensuremath{1_G}},                                       sort={12 },description={главный $F$-характер}}
\newglossaryentry{1Gpp}{name={\ensuremath{1_{G_{p'}}}},                              sort={13 },description={главный брауэров характер}}
\newglossaryentry{1R}{name={\ensuremath{1_R}},                                       sort={14 },description={единица кольца $R$}}
\newglossaryentry{In}{name={\ensuremath{\I_{n}}},                                    sort={15 },description={единичная $n\times n$-матрица}}

\newglossaryentry{AlDr}{name={\ensuremath{A(D)}},                                    sort={A1 },description={\textit{см. с}}}
\newglossaryentry{aklm}{name={\ensuremath{a_{\mbox{}_{KLM}}}},                       sort={A2 },description={структурные константы $\Z(RG)$ в базисе из классовых сумм}}
\newglossaryentry{Akl}{name={\ensuremath{\A_{K,L}}},                                 sort={A3 },description={\textit{см. с}}}
\newglossaryentry{An}{name={\ensuremath{A_n}},                                       sort={A5 },description={знакопеременная группа степени $n$}}
\newglossaryentry{AnnlVr}{name={\ensuremath{\Ann(V)}},                               sort={A6 },description={аннулятор $A$-модуля $V$}}
\newglossaryentry{Annlvr}{name={\ensuremath{\Ann(v)}},                               sort={A7 },description={аннулятор $v$ из $A$-модуля $V$}}
\newglossaryentry{AutX}{name={\ensuremath{\Aut(X)}},                                 sort={A8 },description={группа автоморфизмов группы (поля, кольца, и т.\,п.) $X$}}
\newglossaryentry{AV}{name={\ensuremath{A_V}},                                       sort={A9 },description={образ $F$-алгебры $A$ при отображении $a\mapsto a_V$}}
\newglossaryentry{aV}{name={\ensuremath{a_V}},                                       sort={Aa },description={отображение $v\mapsto va$ из $A$-модуля $V$ в себя для элемента $a\in A$}}

\newglossaryentry{buG}{name={\ensuremath{b^G}},                                      sort={B0 },description={индуцированный $p$-блок}}
\newglossaryentry{BllGPr}{name={\ensuremath{\bl(G|P)}},                              sort={B1 },description={множество $p$-блоков с дефектной группой $P$}}
\newglossaryentry{BllGr}{name={\ensuremath{\bl(G)}},                                 sort={B2 },description={множество $p$-блоков группы $G$}}
\newglossaryentry{BlplGr}{name={\ensuremath{\bl_p(G)}},                              sort={B3 },description={множество $p$-блоков группы $G$}}

\newglossaryentry{C}{name={\ensuremath{\C}},                                         sort={C0 },description={матрица Картана}}
\newglossaryentry{CsB}{name={\ensuremath{\C_B}},                                     sort={C1 },description={матрица Картана $p$-блока $B$},text={\ensuremath{\C_{B_i}}}}
\newglossaryentry{CC}{name={\ensuremath{\CC}},                                       sort={C2 },description={поле комплексных чисел}}
\newglossaryentry{cfcirclGr}{name={\ensuremath{\cf^\circ(G)}},                       sort={C3 },description={множество функций из $\cf(G)$, тождественно равных нулю вне $G_{p'}$}}
\newglossaryentry{cfG}{name={\ensuremath{\cf(G)}},                                   sort={C4 },description={то же, что $\cf_\CC(G)$}}
\newglossaryentry{cfGpp}{name={\ensuremath{\cf(G_{p'})}},                            sort={C5 },description={пространство классовых функций на $p'$-элементах}}
\newglossaryentry{cfRlGr}{name={\ensuremath{\cf_R(G)}},                              sort={C6 },description={$R$-алгебра классовых функций на группе $G$ со значениями в коммутативном кольце $R$}}
\newglossaryentry{chFGr}{name={\ensuremath{\operatorname{ch}_F(G)}},                               sort={C7 },description={множество $F$-характеров группы $G$}}
\newglossaryentry{chGr}{name={\ensuremath{\operatorname{ch}(G)}},                                  sort={C8 },description={то же, что $\ch_\CC(G)$}}
\newglossaryentry{CGX}{name={\ensuremath{\C_G(X)}},                                  sort={C9 },description={централизатор подмножества $X$ в группе $G$}}
\newglossaryentry{CSX}{name={\ensuremath{\C_S(X)}},                                  sort={Ca },description={централизатор подмножества $X$ в алгебре $S$}}

\newglossaryentry{D}{name={\ensuremath{\DD}},                                        sort={D0 },description={матрица разложения}}
\newglossaryentry{DsB}{name={\ensuremath{\DD_B}},                                    sort={D1 },description={матрица разложения $p$-блока $B$},text={\ensuremath{\DD_{B_i}}}}
\newglossaryentry{degx}{name={\ensuremath{\deg \x}},                                 sort={D2 },description={степень обыкновенного характера $\x$}}
\newglossaryentry{degX}{name={\ensuremath{\deg\X}},                                  sort={D3 },description={степень представления $\X$}}
\newglossaryentry{dD}{name={\ensuremath{\df(B)}},                                    sort={D4 },description={дефект $p$-блока $B$}}
\newglossaryentry{dK}{name={\ensuremath{\df(K)}},                                    sort={D5 },description={дефект класса $K$}}
\newglossaryentry{dimV}{name={\ensuremath{\dim_F V}},                                sort={D6 },description={размерность векторного пространства $V$ над полем $F$}}
\newglossaryentry{dimZ}{name={\ensuremath{\dim_\ZZ A}},                              sort={D7 },description={ранг конечно порождённого свободного $\ZZ$-модуля $A$}}
\newglossaryentry{dxxvf}{name={\ensuremath{d^{\,x}_{\x\vf}}},                        sort={D8 },description={высшие числа разложения}}
\newglossaryentry{dxf}{name={\ensuremath{d_{\x\vf}}},                                sort={D9 },description={числа разложения}}

\newglossaryentry{esB}{name={\ensuremath{e_{\mbox{}_B}}},                            sort={E0 },description={центральный идемпотент алгебры $FG$, соответствующий $p$-блоку $B$}}
\newglossaryentry{EndFV}{name={\ensuremath{\End_F(V)}},                              sort={E1 },description={алгебра линейных преобразований векторного пространства $V$ над полем $F$}}
\newglossaryentry{esM}{name={\ensuremath{e_{\mbox{}_M}}},                            sort={E2 },description={центральный идемпотент $F$-алгебры, соответствующий её неприводимому модулю $M$}}
\newglossaryentry{EndR}{name={\ensuremath{\End_R(M)}},                               sort={E3 },description={то же, что $\Hom_R(M,M)$}}
\newglossaryentry{EndV}{name={\ensuremath{\End(V)}},                                 sort={E4 },description={кольцо эндоморфизмов абелевой группы $V$}}
\newglossaryentry{esx}{name={\ensuremath{e_\x}},                                     sort={E5 },description={центральные идемпотенты алгебры $\CC G$, соответствующий неприводимому характеру $\x$}}
\newglossaryentry{expG}{name={\ensuremath{\exp G }},                                 sort={E6 },description={экспонента группы $G$}}

\newglossaryentry{fsA}{name={\ensuremath{f_\A^{\protect\phantom{A}}}},               sort={F0 },description={то же, что $\sum_{\x\in\A}e_\x$}}
\newglossaryentry{fsB}{name={\ensuremath{f_{\mbox{}_B}}},                            sort={F1 },description={то же, что $f_{\irr(B)}$}}
\newglossaryentry{Fq}{name={\ensuremath{\FF_q}},                                     sort={F2 },description={конечное поле из $q$ элементов}}
\newglossaryentry{fRX}{name={\ensuremath{\f_R(X)}},                                  sort={F3 },description={алгебра $R$-значных функций на множестве $X$}}
\newglossaryentry{FrlGr}{name={\ensuremath{\Fr(G)}},                                 sort={F4 },description={подгруппа Фраттини группы $G$}}

\newglossaryentry{GalEF}{name={\ensuremath{\Gal(E,F)}},                              sort={G0 },description={группа Галуа расширения полей $F\le E$}}
\newglossaryentry{gchG}{name={\ensuremath{\gch(G)}},                                 sort={G1 },description={кольцо $\gch_\CC(G)$ обобщённых $\CC$-характеров}}
\newglossaryentry{gchFG}{name={\ensuremath{\gch_F(G)}},                              sort={G2 },description={кольцо обобщённых $F$-характеров}}
\newglossaryentry{gp}{name={\ensuremath{g_p}},                                       sort={G3 },description={$p$-часть элемента $g$}}
\newglossaryentry{GLnR}{name={\ensuremath{\GL_n(R)}},                                sort={G4 },description={полная линейная группа над коммутативным кольцом $R$}}
\newglossaryentry{gpp}{name={\ensuremath{g_{p'}}},                                   sort={G5 },description={$p'$-часть элемента $g$}}
\newglossaryentry{Gspi}{name={\ensuremath{G_{\pi}}},                                 sort={G6 },description={множество $\pi$-элементов группы $G$}}
\newglossaryentry{Gpp}{name={\ensuremath{G_{p'}}},                                   sort={G7 },description={множество всех $p$-регулярных элементов группы $G$}}

\newglossaryentry{HomR}{name={\ensuremath{\Hom_R(M,N)}},                             sort={H0 },description={множество всех гомоморфизмов $R$-модуля $M$ в $R$-модуль $N$}}

\newglossaryentry{IBrlBr}{name={\ensuremath{\iBr(B)}},                               sort={I0 },description={то же, что $B\cap\iBr(G)$, где $B$ --- $p$-блок группы $G$}}
\newglossaryentry{IBrlGr}{name={\ensuremath{\iBr(G)}},                               sort={I1 },description={то же, что $\iBr_p(G)$}}
\newglossaryentry{IBrMlGr}{name={\ensuremath{\iBr_M(G)}},                            sort={I2 },description={множество всех неприводимых $p$-модулярных характеров группы $G$ относительно идеала $M$}}
\newglossaryentry{IBrplGr}{name={\ensuremath{\iBr_p(G)}},                            sort={I3 },description={множество всех неприводимых $p$-модулярных характеров группы $G$}}
\newglossaryentry{IrrlBr}{name={\ensuremath{\irr(B)}},                               sort={I4 },description={то же, что $B\cap\irr(G)$, где $B$ --- $p$-блок группы $G$}}
\newglossaryentry{IrrlGr}{name={\ensuremath{\irr(G)}},                               sort={I5 },description={множество $\irr_\CC(G)$ всех неприводимых комплексных характеров группы $G$}}
\newglossaryentry{IrrFlGr}{name={\ensuremath{\irr_F(G)}},                            sort={I6 },description={множество характеров всех неприводимых $F$-представлений группы $G$}}
\newglossaryentry{IGth}{name={\ensuremath{\II_G(\th)}},                              sort={I7 },description={группа инерции характера $\th$}}
\newglossaryentry{IGU}{name={\ensuremath{\II_G(U)}},                                 sort={I8 },description={группа инерции модуля $U$}}
\newglossaryentry{IGX}{name={\ensuremath{\II_G(\X)}},                                sort={I9 },description={группа инерции представления $\X$}}

\newglossaryentry{JlAr}{name={\ensuremath{\J(A)}},                                   sort={J0 },description={радикал $R$-алгебры $A$}}
\newglossaryentry{JlSr}{name={\ensuremath{\J(S)}},                                   sort={J1 },description={радикал кольца $S$}}

\newglossaryentry{kervf}{name={\ensuremath{\ker\vf}},                                sort={K1 },description={ядро группового гомоморфизма $\vf$}}
\newglossaryentry{KlGPr}{name={\ensuremath{\K(G|P)}},                                sort={K2 },description={множество классов с $p$-дефектной группой $P$}}
\newglossaryentry{KGpp}{name={\ensuremath{\K(G_{p'})}},                              sort={K3 },description={множество всех $p$-регулярных классов группы $G$}}
\newglossaryentry{KlGpPr}{name={\ensuremath{\K(G_{p'}|P)}},                          sort={K4 },description={множество $p$-регулярных классов с $p$-дефектной группой $P$}}
\newglossaryentry{Ker}{name={\ensuremath{\Ker\vf}},                                  sort={K5 },description={ядро кольцевого гомоморфизма $\vf$}}
\newglossaryentry{KlGr}{name={\ensuremath{\K(G)}},                                   sort={K6 },description={множество классов сопряжённости группы $G$}}
\newglossaryentry{kerx}{name={\ensuremath{\ker \x}},                                 sort={K7 },description={ядро обыкновенного характера $\x$}}
\newglossaryentry{kerX}{name={\ensuremath{\ker \X}},                                 sort={K8 },description={ядро представления $\X$ группы}}

\newglossaryentry{LinG}{name={\ensuremath{\lin(G)}},                                 sort={L0 },description={множество всех линейных обыкновенных характеров группы $G$}}
\newglossaryentry{LinFG}{name={\ensuremath{\lin_F(G)}},                              sort={L1 },description={множество всех линейных $F$-характеров группы $G$}}

\newglossaryentry{Mtil}{name={\ensuremath{\protect\phantom{M}
          \mathllap{\widetilde{\protect\phantom{R}\mkern2mu}\mathllap{M}} }},        sort={M0 },description={максимальный идеал кольца $\wt{Z}$}}
\newglossaryentry{MlAr}{name={\ensuremath{\M(A)}},                                   sort={M1 },description={система представителей классов изоморфных неприводимых $A$-модулей}}
\newglossaryentry{MnR}{name={\ensuremath{\MM_n(R)}},                                 sort={M2 },description={алгебра $n\times n$\,-\,матриц над коммутативным кольцом $R$}}
\newglossaryentry{Mchi}{name={\ensuremath{M_\x}},                                    sort={M3 },description={неприводимый $\CC G$-модуль с характером $\x$}}

\newglossaryentry{NGg}{name={\ensuremath{\N_G(g)}},                                  sort={N0 },description={нормализатор подгруппы $\la g\ra$, где $g \in G$}}
\newglossaryentry{NGX}{name={\ensuremath{\N_G(X)}},                                  sort={N1 },description={нормализатор подмножества $X$ в группе $G$}}
\newglossaryentry{nMlVr}{name={\ensuremath{\n_M(V)}},                                sort={N2 },description={число неприводимых прямых слагаемых модуля $V$, изоморфных $M$}}
\newglossaryentry{np}{name={\ensuremath{n_p}},                                       sort={N3 },description={$p$-часть числа $n$}}
\newglossaryentry{npi}{name={\ensuremath{n_{\pi}}},                                  sort={N4 },description={$\pi$-часть числа $n$}}
\newglossaryentry{npp}{name={\ensuremath{n_{p'}}},                                   sort={N5 },description={$p'$-часть числа $n$}}

\newglossaryentry{OrbGm}{name={\ensuremath{\Orb_G(m)}},                              sort={O0 },description={орбита элемента $m$ относительно действия группы $G$}}
\newglossaryentry{opG}{name={\ensuremath{\OO_p(G)}},                                 sort={O1 },description={наибольшая нормальная $p$-подгруппа группы $G$}}
\newglossaryentry{opprG}{name={\ensuremath{\OO_{p'}(G)}},                            sort={O2 },description={наибольшая нормальная $p'$-подгруппа группы $G$}}
\newglossaryentry{Orbm}{name={\ensuremath{\Orb(m)}},                                 sort={O3 },description={орбита элемента $m$}}

\newglossaryentry{PX}{name={\ensuremath{\P(X)}},                                     sort={P0 },description={множество подмножеств множества $X$}}

\newglossaryentry{QQ}{name={\ensuremath{\QQ}},                                       sort={Q0 },description={поле рациональных чисел}}
\newglossaryentry{Q2n}{name={\ensuremath{Q_{2^n}}},                                  sort={Q1 },description={обобщённая группа кватернионов}}
\newglossaryentry{Q8}{name={\ensuremath{Q_8}},                                       sort={Q2 },description={группа кватернионов}}
\newglossaryentry{Qm}{name={\ensuremath{\QQ_m}},                                     sort={Q3 },description={то же, что $\QQ(\z)$, где $\z=e^{\frac{2\pi i}{m}}$}}

\newglossaryentry{RR}{name={\ensuremath{\RR}},                                       sort={R0 },description={поле вещественных чисел}}
\newglossaryentry{Rop}{name={\ensuremath{R^{op}}},                                   sort={R1 },description={кольцо, противоположное кольцу $R$}}
\newglossaryentry{RF}{name={\ensuremath{\R_F}},                                      sort={R2 },description={регулярное $F$-представление группы}}

\newglossaryentry{Sn}{name={\ensuremath{S_n}},                                       sort={S0 },description={симметрическая группа степени $n$}}
\newglossaryentry{SsGlgr}{name={\ensuremath{\SS_G(g)}},                              sort={S1 },description={$p$-сечение группы $G$, содержащее $g\in G$}}
\newglossaryentry{Stm}{name={\ensuremath{\St(m)}},                                   sort={S2 },description={стабилизатор элемента $m$}}
\newglossaryentry{StGm}{name={\ensuremath{\St_G(m)}},                                sort={S3 },description={стабилизатор элемента $m$ относительно действия группы $G$}}
\newglossaryentry{supp}{name={\ensuremath{\supp x}},                                 sort={S4 },description={носитель элемента $x\in RG$}}
\newglossaryentry{SylsplGr}{name={\ensuremath{\Syl_p(G)}},                           sort={S5 },description={множество $p$-силовских подгрупп группы $G$}}

\newglossaryentry{atop}{name={\ensuremath{M^\top}},                                  sort={T0 },description={транспонированная к матрице $M$}}
\newglossaryentry{trB}{name={\ensuremath{\tr B}},                                    sort={T1 },description={след матрицы $B$}}

\newglossaryentry{VupG}{name={\ensuremath{V^G}},                                     sort={V0 },description={индуцированный $FG$-модуль},text={\ensuremath{W^G}}}
\newglossaryentry{VH}{name={\ensuremath{V_H}},                                       sort={V1 },description={$RG$-модуль $V$, рассматриваемый как $RH$-модуль для подгруппы $H\le G$}}

\newglossaryentry{Wg}{name={\ensuremath{Wg}},                                        sort={W0 },description={модуль, сопряжённый с модулем $W$ элементом $g$}}
\newglossaryentry{WsB}{name={\ensuremath{\W_B}},                                     sort={W1 },description={то же, что $\W_B(FG)$}}
\newglossaryentry{WsBlFGr}{name={\ensuremath{\W_B(FG)}},                             sort={W2 },description={блок $e_{\mbox{}_B}FG$ алгебры $FG$, соответствующий $p$-блоку $B$}}
\newglossaryentry{WMlVr}{name={\ensuremath{\W_M(V)}},                                sort={W3 },description={$M$-однородная компонента модуля $V$}}
\newglossaryentry{Wal}{name={\ensuremath{W\a}},                                      sort={W4 },description={модуль, сопряжённый с модулем $W$ автоморфизмом $\a$}}

\newglossaryentry{XsB}{name={\ensuremath{\X_B}},                                     sort={X2 },description={ограничение представления $\X$ алгебры на её подалгебру $B$}}
\newglossaryentry{XE}{name={\ensuremath{\X^E}},                                      sort={X3 },description={представление $\X$ над расширением $E$ основного поля}}
\newglossaryentry{xG}{name={\ensuremath{x^G}},                                       sort={X4 },description={класс сопряжённости, содержащий элемент $x\in G$}}
\newglossaryentry{XsH}{name={\ensuremath{\X_H}},                                     sort={X5 },description={ограничение $F$-представления $\X$ группы на её подгруппу $H$}}
\newglossaryentry{xsK}{name={\ensuremath{x_{\mbox{}_K}}},                            sort={X6 },description={представитель класса сопряжённости $K$}}
\newglossaryentry{Xchi}{name={\ensuremath{\X_\x}},                                   sort={X7 },description={неприводимое $\CC$-представление с характером $\x$}}
\newglossaryentry{XuG}{name={\ensuremath{\X^G}},                                     sort={X8 },description={представление, индуцированное представлением $\X$}}
\newglossaryentry{Xg}{name={\ensuremath{\X^g}},                                      sort={X9 },description={представление, сопряжённое с представлением $\X$ элементом $g$}}
\newglossaryentry{xal}{name={\ensuremath{\X^\a}},                                    sort={Xa },description={представление, сопряжённое с представлением $\X$ автоморфизмом $\a$}}

\newglossaryentry{ZZ}{name={\ensuremath{\ZZ}},                                       sort={Z0 },description={кольцо целых чисел}}
\newglossaryentry{ZZZ}{name={\ensuremath{\ov{\ZZ}}},                                 sort={Z1 },description={кольцо всех целых алгебраических чисел}}
\newglossaryentry{ZZn}{name={\ensuremath{\ZZ_n}},                                    sort={Z2 },description={кольцо вычетов по модулю $n$}}
\newglossaryentry{Ztil}{name={\ensuremath{\wt{Z}}},                                  sort={Z3 },description={локализация кольца $\ov{\ZZ}$ относительно максимального идеала $M$}}
\newglossaryentry{ZG}{name={\ensuremath{\Z(G)}},                                     sort={Z4 },description={центр группы $G$}}
\newglossaryentry{ZilGr}{name={\ensuremath{\Z_i(G)}},                                sort={Z5 },description={подгруппа в $G$ такая, что $\Z_i(G)/\Z_{i-1}(G)=\Z(G/\Z_{i-1}(G))$}}
\newglossaryentry{ZlRr}{name={\ensuremath{\Z(R)}},                                   sort={Z6 },description={центр кольца $R$}}
\newglossaryentry{Zst}{name={\ensuremath{\Z^*(G)}},                                  sort={Z7 },description={прообраз в $G$ группы $\Z(G/\OO_{2'}(G))$}}
\newglossaryentry{Zpim1}{name={\ensuremath{\ZZ[\pi^{-1}]}},                          sort={Z8 },description={кольцо $\ZZ[\,\frac{1}{p}\mid p\in\pi]$ для некоторого множества простых целых чисел $\pi$}}
\newglossaryentry{Zchi}{name={\ensuremath{\Z(\x)}},                                  sort={Z9 },description={множество $\{g\in G\bigm|\  |\x(g)|=\x(1)\}$} для обыкновенного характера $\x$}

\newglossaryentry{asAlKr}{name={\ensuremath{\a_{\A}^{\protect\phantom{A}}(K)}},      sort={a0 },description={коэффициент при $\wh K$ в разложении $f_{\A}^{\phantom{A}}$}}

\newglossaryentry{bsP}{name={\ensuremath{\b_{\mbox{}_P}}},                           sort={b0 },description={гомоморфизм Брауэра}}

\newglossaryentry{DlKr}{name={\ensuremath{\D(K)}},                                   sort={d0 },description={множество $p$-дефектных групп класса $K$}}
\newglossaryentry{DplKr}{name={\ensuremath{\D_p(K)}},                                sort={d1 },description={множество $p$-дефектных групп класса $K$}}
\newglossaryentry{DlBr}{name={\ensuremath{\D(B)}},                                   sort={d2 },description={множество $p$-дефектных групп $p$-блока $B$}}
\newglossaryentry{dellBr}{name={\ensuremath{\d(B)}},                                 sort={d3 },description={дефектная группа $p$-блока $B$}}
\newglossaryentry{dij}{name={\ensuremath{\d_{ij}}},                                  sort={d4 },description={символ Кронекера}}
\newglossaryentry{dellKr}{name={\ensuremath{\d(K)}},                                 sort={d5 },description={$p$-дефектная группа класса $K$}}
\newglossaryentry{delplKr}{name={\ensuremath{\d_p(K)}},                              sort={d6 },description={$p$-дефектная группа класса $K$}}

\newglossaryentry{etasB}{name={\ensuremath{\eta_{\mbox{}_B}}},                       sort={e0 },description={$B$-часть классовой $\eta\in \cf(G)$}}

\newglossaryentry{Theta}{name={\ensuremath{\Theta}},                                 sort={f0 },description={таблица значений главных неразложимых характеров на $p$-регулярных элементах}}
\newglossaryentry{thsvf}{name={\ensuremath{\th_\vf}},                                sort={f2 },description={главный неразложимый характер, соответствующий брауэрову характеру $\vf$ }}

\newglossaryentry{Mi}{name={\ensuremath{M^\i}},                                      sort={i0 },description={обратно-транспонированная к невырожденной матрице $M$}}
\newglossaryentry{Ki}{name={\ensuremath{K^\i}},                                      sort={i1 },description={класс элементов, обратных к элементам класса $K$}}

\newglossaryentry{vklgr}{name={\ensuremath{\vk(g)}},                                 sort={k0 },description={число способов представить $g$ в виде коммутатора}}

\newglossaryentry{lsB}{name={\ensuremath{\l_{\mbox{}_B}}},                           sort={l0 },description={центральный характер, соответствующий блоку $B$}}
\newglossaryentry{lsvf}{name={\ensuremath{\l_\vf}},                                  sort={l1 },description={центральный характер $\Z(FG)\to F$, соответствующий брауэрову характеру $\vf\in \iBr(G)$}}
\newglossaryentry{lsx}{name={\ensuremath{\l_\x}},                                    sort={l2 },description={центральный характер $Z(FG)\to F$, соответствующий обыкновенному характеру $\x\in \irr(G)$}}

\newglossaryentry{rh}{name={\ensuremath{\r}},                                        sort={r0 },description={регулярный $\CC$-характер $\r_{\mbox{}_\CC}$}}
\newglossaryentry{rhF}{name={\ensuremath{\r_{\mbox{}_F}}},                           sort={r1 },description={регулярный $F$-характер}}
\newglossaryentry{Rb}{name={\ensuremath{\vr}},                                       sort={r2 },description={отображение Робинсона  $\Z(FG)\to \Z(FG)$}}

\newglossaryentry{PhiG}{name={\ensuremath{\Phi(G)}},                                 sort={s0 },description={таблица брауэровых характеров группы $G$}}
\newglossaryentry{PhiMG}{name={\ensuremath{\Phi_M(G)}},                              sort={s1 },description={таблица брауэровых характеров группы $G$ относительно максимального идеала $M$}}
\newglossaryentry{PhipG}{name={\ensuremath{\Phi_p(G)}},                              sort={s2 },description={таблица брауэровых характеров группы $G$ в характеристике $p$}}
\newglossaryentry{vfal}{name={\ensuremath{\vf^\a}},                                  sort={s3 },description={классовая функция, сопряжённая c $\vf$ автоморфизмом $\a$}}
\newglossaryentry{vfG}{name={\ensuremath{\vf^G}},                                    sort={s4 },description={индуцированная классовая функция}}
\newglossaryentry{vfg}{name={\ensuremath{\vf^g}},                                    sort={s5 },description={классовая функция, сопряжённая c $\vf$ элементом $g$}}
\newglossaryentry{vfsH}{name={\ensuremath{\vf_{H^{\protect\vphantom{A^a}}}}},        sort={s6 },description={ограничение классовой функции $\vf$ на подгруппу $H$}}
\newglossaryentry{phi}{name={\ensuremath{\phi}},                                     sort={s7 },description={функция Эйлера}}

\newglossaryentry{XiG}{name={\ensuremath{\XX(G)}},                                   sort={x0 },description={таблица обыкновенных характеров группы $G$}}
\newglossaryentry{chisi}{name={\ensuremath{\x^\s}},                                  sort={x1 },description={алгебраически сопряжённый характер}}
\newglossaryentry{chiuG}{name={\ensuremath{\x^G}},                                   sort={x3 },description={характер, индуцированный $F$-характером $\x$}}

\newglossaryentry{OmG}{name={\ensuremath{\Om(G)}},                                   sort={z0 },description={таблица центральных характеров $\om_\x$}}